\documentclass[english]{smfbook}
\usepackage[utf8]{inputenc}
\usepackage[english,francais]{babel}
\usepackage{smfthm}
\usepackage{mathrsfs}
\usepackage{graphics}
\usepackage{hyperref} 
\usepackage{amssymb, amsmath, mathabx}
\usepackage{lmodern}
\usepackage[all]{xy}
\usepackage{multicol}
\usepackage[shortlabels]{enumitem}
\usepackage{scalerel}
\makeindex
\DeclareUnicodeCharacter{00A0}{\relax}

\author{Antoine Ducros}
\address{Sorbonne Universités, UPMC Univ Paris 06, Institut de Mathématiques de Jussieu-Paris
Rive Gauche, UMR 7586, CNRS, Univ Paris Diderot, Sorbonne Paris Cité, F-75005, Paris, France
}
\email{antoine.ducros\at imj-prg.fr}
\urladdr{http://www.imj-prg.fr/$\sim$antoine.ducros}
\title{Families of Berkovich spaces}
\alttitle{Familles d'espaces de Berkovich}
\NumberTheoremsAs{subsection}
\SwapTheoremNumbers

\newcommand{\cf}{cf.\@~}
\newcommand{\eg}{e.\@g.\@}
\newcommand{\etc}{etc.\@}
\newcommand{\ie}{i.\@e.\@}
\newcommand{\ff}{ff.\@~}
\newcommand{\loccit}{loc.\@~cit.\@}
\newcommand{\opcit}{op.\@~cit.\@}
\newcommand{\resp}{resp.\@~}
\newcommand{\Chp}{Chapt.\@~}
\newcommand{\Cor}{Cor.\@~}
\newcommand{\Def}{Def.\@~}
\newcommand{\Prop}{Prop.\@~}
\newcommand{\Rem}{Rem.\@~}
\newcommand{\Th}{Thm.\@~}

\newcommand{\abs}[1]{\left|#1\right|}
\newcommand{\adht}[2]{\overline{ #1}^{#2}}
\newcommand{\adhz}[2]{\overline{ #1}^{#2_{\mathrm {Zar}}}}
\newcommand{\an}{^{\mathrm{an}}}
\newcommand{\al}{^{\mathrm{al}}}
\newcommand{\all}{^{\scaleobj{0.7}{\mathrm{al}}}}
\newcommand{\bij}[2]{{\mathrm {Bij}}(#1\to #2)}
\newcommand{\brac}[1]{\mathopen[#1\mathclose]}
\newcommand{\car}{{\rm char.}\@~}
\newcommand{\codim}{\mathrm{codim}}
\newcommand{\gpm}{^{\times}}
\newcommand{\grot}{_{\mathrm G}}
\newcommand{\hotimes}{\widehat \otimes}
\newcommand{\hr}[1]{\mathscr H(#1)}
\newcommand{\hrt}[1]{\widetilde{\mathscr H(#1)}}
\newcommand{\inv}{^{-1}}
\newcommand{\norm}[1]{\mathopen\|#1\mathclose \|}
\newcommand{\pfb}[2]{\mathsf Q_{\mathrm{fib}}(#1/#2)}
\newcommand{\rk}[1]{\mathrm{rk}_{#1}}
\newcommand{\spec}{\mathrm{Spec}\;}
\newcommand{\supp}[1]{\mathrm{Supp}(\mathscr #1)}
\newcommand{\uflat}[2]{\mathrm{Flat}(#1/#2)}
\newcommand{\zar}{_{\mathrm {Zar}}}

\newcommand{\field}{$(\mathsf F)$}
\newcommand{\gen}{$(\mathsf G)$}
\newcommand{\hci}{$(\mathsf H_{\mathrm {CI}})$}
\newcommand{\hreg}{$(\mathsf H_{\mathrm{reg}})$}
\newcommand{\hv}{$(\mathsf H)$}
\newcommand{\hwk}{$(\mathsf T_{\mathrm{weak}})$}
\newcommand{\hstr}{$(\mathsf T_{\mathrm{strong}})$}
\newcommand{\hwkk}{$(\mathsf T'_{\mathrm{weak}})$}
\newcommand{\hstrr}{$(\mathsf T'_{\mathrm{strong}})$}
\newcommand{\open}{$(\mathsf O)$}

\newcommand{\A}{\mathbf A}
\newcommand{\C}{\mathbf C}
\newcommand{\F}{\mathbf F}
\newcommand{\N}{\mathbf Z_{\geq 0}}
\renewcommand{\P}{\mathbf P}
\newcommand{\Q}{\mathbf Q}
\newcommand{\R}{\mathbf R}
\newcommand{\Z}{\mathbf Z}

\renewcommand{\H}{\mathrm H}
\renewcommand{\d}{\mathrm d}
\renewcommand{\phi}{\varphi}
\renewcommand{\epsilon}{\varepsilon}
\renewcommand{\leq}{\leqslant}
\renewcommand{\geq}{\geqslant}

\renewcommand{\theequation}{\alph{equation}}
\SetEnumerateShortLabel{1}{\textnormal{(\arabic{enumi})}}
\SetEnumerateShortLabel{i}{\textnormal{(\roman{enumi})}}
\SetEnumerateShortLabel{a}{\textnormal{(\alph{enumi})}}
\SetEnumerateShortLabel{A}{\textnormal{(\Alph{enumi})}}
\SetEnumerateShortLabel{2}{\textnormal{(\arabic{enumii})}}
\SetEnumerateShortLabel{j}{\textnormal{(\roman{enumii})}}
\SetEnumerateShortLabel{b}{\textnormal{(\arabic{enumi}\alph{enumii})}}
\SetEnumerateShortLabel{c}{\textnormal{(\alph{enumii})}}
\SetEnumerateShortLabel{B}{\textnormal{(\Alph{enumii})}}

\newcommand{\z}{\;\;\;\;\;}

\begin{document}
\newenvironment{step}{\begin{enonce}{Reduction Statement}}{\end{enonce}}

\frontmatter
\begin{abstract}
This memoir is devoted to the systematic study of relative properties
in the context of Berkovich analytic spaces. We first develop a theory of flatness in this setting. 
After showing through a counter-example that naive flatness cannot be the right notion because
it is not stable under base change, we
define flatness by {\em requiring}
invariance under base change, and we study a first important class of flat morphisms, that of quasi-smooth ones. 

We then show the existence of local {\em d\'evissages} (in the spirit of
Raynaud and Gruson's work) for coherent sheaves, which we use, together with a study of the local rings of ``generic fibers''
of morphisms, to prove that a {\em boundaryless}, naively flat morphism is flat. 

After that we prove that the image of a compact strictly analytic space by a 
flat morphism is an analytic domain of the target 
and that it admits, when the source is stricly analytic, a compact, flat multisection (i.e., a compact, flat cover of relative dimension zero over which there is a section).
This was first proved in the rigid-analytic context by Raynaud, but our proof is completely different: it is based upon Temkin's theory of the reduction of analytic germs and does not make any use of formal models.

In the last part of this work we study where various
interesting pointwise
relative properties are satisfied. We first prove
that
the flat
locus of a given morphism of analytic spaces is a Zariski-open subset of the source (we follow the method that was introduced by Kiehl 
for the complex analytic analogue of this statement). We then look at the
loci at which a point satisfies various commutative algebra properties {\em on its fiber}: being geometrically regular, 
geometrically $R_m$, complete intersection, or Gorenstein; being $S_m$ or Cohen-Macaulay. 
We prove that the results we could expect
actually hold: these loci are (locally) Zariski-constructible, and Zariski-open under suitable extra
assumptions (flatness, and also equidimensionality for~$S_m$ and~geometric $R_m$); for that purpose, 
we first study the general properties of the locally Zariski-constructible subsets of an analytic space. 

\end{abstract}

\begin{altabstract}
Ce mémoire est consacré à une étude systématique des propriétés relatives
dans le contexte des espaces de Berkovich. Nous commençons par développer
une théorie de la platitude dans ce cadre. L'acception naïve de cette notion est inadaptée : 
nous montrons en effet par un contre-exemple qu'elle n'est pas stable par changement de base, 
ce qui nous conduit à {\em imposer}
cette stabilité dans la définition. Nous étudions une première classe importante de morphismes
plats : celle des morphismes {\em quasi-lisses}.

Nous montrons ensuite l'existence de {\em dévissages} locaux (dans
l'esprit de Raynaud et Gruson) pour les faisceaux cohérents. Joint à une étude des anneaux
locaux des fibres «génériques» des morphismes, cela nous permet de montrer
qu'un morphisme {\em sans bord} qui est plat au sens naïf l'est encore au nôtre. 

Puis nous démontrons que l'image d'un espace
strictement analytique compact par un morphisme plat est un domaine analytique
du but, et qu'elle admet, lorsque la source est strictement analytique, une multisection compacte et plate
({\em i.e.} un revêtement plat, compact, et de dimension 
relative nulle sur lequel le morphisme considéré possède une section). 
Cela avait déjà été établi dans le contexte rigide-analytique par Raynaud, mais notre preuve est complètement différente : elle
repose sur la réduction à la Temkin des germes d'espaces analytiques et ne fait pas appel aux schémas formels.

Dans la dernière partie de ce travail, nous étudions les lieux de validité 
sur la source de certaines propriétés relatives. Nous y démontrons pour commencer
que le lieu de platitude d'un morphisme d'espaces analytiques est un ouvert de Zariski
de la source (nous suivons la méthode utilisée par Kiehl pour établir l'assertion correspondante
en géométrie analytique complexe). Cela nous permet de procéder à l'investigation 
systématique des ensemble de points satisfaisant {\em dans leur fibre} les propriétés classiques 
de l'algèbre commutative : être géométriquement
régulier, géométriquement $R_m$, d'intersection complète
ou de Gorenstein ; être
$S_n$ ou de Cohen-Macaulay. Nous prouvons que les énoncés auxquels
ont peut s'attendre sont effectivement vérifiés : ces lieux de validités sont (localement) Zariski-constructibles, et 
sont des ouverts de Zariski sous certaines hypothèses supplémentaires (platitude, ainsi qu'équidimensionalité pour 
les propriétés $R_m$ ou $S_m$) ; dans ce but, nous procédons tout d'abord à une étude générale des parties
localement Zariski-constructibles d'un espace analytique. 
\end{altabstract}

\altkeywords{Espaces de Berkovich, platitude, propriétés relatives, lieux constructibles}
\thanks{While writing this memoir the author was supported by the french {\em Agence nationale de la recherche} (ANR), through
the projects {\em Espaces de Berkovich} (ANR 07-JCJC-0004-CSD5), {\em Valuations, combinatoire et th\'eorie des mod\`eles} (ANR-13-BS01-0006), 
and {\em D\'efinissabilit\'e en g\'eom\'etrie non archim\'edienne}
(ANR-15-CE40-0008). He benefited from two {\sc Rosi and Max Varon}
professorships at Weizmann Institute (Rehovot, Israel), both with an invitation by Vladimir Berkovich: from January 15th till March 14th, 2009 and from April 19th till July 18th,  2015.
He was a junior member of the \emph{Institut universitaire de France}
from October 2012 till September 2017}. 

\maketitle
\tableofcontents
\mainmatter
\setcounter{chapter}{-1}
\chapter{Introduction}

This memoir is, roughly speaking, devoted to 
a systematic investigation
of
{\em families}
of objects in Berkovich's non-Archimedean analytic geometry (\cite{berkovich1990}, \cite{berkovich1993}).
More precisely, let us assume that we are given a morphism~$Y\to X$ between analytic spaces, 
and an object~$D$ on~$Y$ of a certain kind (think of the space~$Y$ itself, or a coherent sheaf, or a complex 
of coherent sheaves\ldots). Every point~$x$
of~$X$ gives then rise to an object~$D_x$, living on the fiber~$Y_x$,
and we thus get in some sense an analytic
family of objects parametrized by the space~$X$. The quite vague
problem we would like to address is the following: {\em how do the object~$D_x$ 
and its relevant properties vary?} Of course, such questions have been 
intensively studied for a long time in {\em algebraic} geometry, especially
by Grothendieck and his school, and 
our guideline has been to establish analytic 
avatars of their results every time it was possible. 
Let us now give a quick overview of our work.

\section{First step: flatness in the Berkovich setting}

\subsection{Motivation} 
In scheme theory, the key notion 
upon which the study of families
is based is {\em flatness}. 
This is a property of families of {\em coherent sheaves}, which encodes
more or less the intuitive idea of a {\em reasonable} variation
(this is why 
there is almost always a flatness assumption in 
the description of moduli problems). 

The point is that the 
study of general families is
often reduced 
(typically, through a suitable stratification of the base scheme) to the case where
some of the coherent sheaves involved are flat over the parameter space, 
which is easier to handle. 

But we would like to emphasize that flatness
is also a crucial technical tool for many other purposes. Let us mention for example:

\begin{itemize}
\item[$\bullet$] descent theory; 

\item[$\bullet$]  the first occurrence of flatness in algebraic geometry, in the celebrated paper GAGA by Serre \cite{serre1955}, where the following plays a major role: if~$X$ is a complex algebraic variety and if~$X\an$ denotes the corresponding analytic space, then for every~$x\in X(\C)$ the ring~$\mathscr  O_{X\an,x}$ is flat over~$\mathscr  O_{X,x}$. 
\end{itemize}

By analogy, our first step toward the understanding of analytic families has
been the development of a theory of flatness in Berkovich geometry, which is the core
of this memoir; 
and similarly to what happens for schemes, we hope that it will have
many applications beyond the study of families.

 In fact, flatness in non-Archimedean geometry had already been considered, but in the {\em rigid-analytic}
setting, following ideas of Raynaud; see the papers~\cite{frg1} and~\cite{frg2}
by Bosch and L\"utkebohmert, as well as
the more comprehensive recent study by Abbes in~\cite{abbes2010}. We now give
a more precise discussion of rigid-analytic flatness,  before saying some words
concerning our definition in the Berkovich framework. 

%\medskip
%What we do in this paper is starting a systematic investigation of flatness in {\em non-Archimedean analytic geometry}, and more precisely in the {\em Berkovich} setting (\cite{brk1}, \cite{brk2}).  We have chosen to give 'purely Berkovich' proofs, without using formal or algebraic models, and the related subtle, sophisticated results (flattening, semi-stable reduction, reduced fiber theorem...). Therefore we sometimes give new proofs of facts which were originally established using some of those highly non-trivial theorems about models (see below some examples of such facts in the summary of our results). The only 'reduction arguments' we will use are based upon Temkin's theory of the reduction of analytic germs, which involves general valuation theory, and especially Riemann-Zariski spaces. 

\subsection{Flatness in rigid geometry} 
The definition of flatness in the rigid setting is as simple as one may hope: if~$f: Y\to X$ is a morphism between rigid spaces and if~$\mathscr  F$
is a coherent sheaf on~$Y$, it is rig-flat
over~$X$ at a point~$y\in Y$ if it is
flat over~$X$ at~$y$ {\em in the sense of the theory of locally ringed spaces}; \ie,
the stalk~$\mathscr  F_y$ is a flat~$\mathscr  O_{X,f(y)}$-module (this is the definition by Abbes;
the original one by Bosch and L\"utkebohmert was slightly different, but both are easily seen to be equivalent
using the good properties of the completion of local rings as far as flatness is concerned, \cf
\cite{sga1}, Exposé  IV, \Cor 5.8 ).

Flatness in the above sense behaves well: it is stable under base change and ground field extension. But contrary to what happens in scheme theory, this is in
no way obvious, because base change and ground field extension are defined using {\em completed} tensor products. 
Roughly speaking, the proofs proceed as follows (see~\cite{frg1}, \cite{frg2} and~\cite{abbes2010}):

\begin{itemize}
\item[$\bullet$] the study of {\em rigid} flatness is reduced to that of {\em formal} flatness
through formal avatars of Raynaud-Gruson flattening techniques, which are used to build a flat formal model of any given rig-flat coherent sheaf; 

\item[$\bullet$] the study of {\em formal} flatness is reduced to that of {\em algebraic} flatness, in a more standard way, upon dividing by various ideals of definition and using flatness criteria in the spirit of \cite{sga1}, Exposé IV. 
\end{itemize}
Let us mention that this general strategy (formal flattening and reduction modulo an ideal of definition to replace an analytic problem with an algebraic one) was also used by Raynaud to prove the following fact: if~$\phi: Y\to X$ is a flat morphism between affinoid rigid spaces,
$\phi(Y)$ is a finite union of affinoid domains of~$X$ (\cf \cite{frg2},
\Cor  5.11).

\subsection{Flatness in Berkovich geometry}
We fix from now on and for the remaining part of this introduction 
a complete, non-Archimedean field~$k$. 
We shall
only consider {\em Berkovich} analytic spaces.
Any analytic space $X$
comes with a usual topology, but also with a set-theoretic Grothendieck topology 
which refines it, the so-called {\em {\rm G}-topology} --the corresponding site
is denoted by $X\grot$; the archetypal example
of a G-covering is the covering of an affinoid
space by finitely many affinoid domains.
The site $X\grot$ is equipped with a
coherent sheaf of rings $\mathscr O_{X\grot}$, whose
restriction $\mathscr O_X$ to the category of open subsets
makes $X$ a locally ringed space.
But $(X,\mathscr O_X)$ does not seem to be
easily tractable in general: for example, one does not know whether $\mathscr O_X$ is coherent, 
nor whether its stalks are noetherian (one only knows that they are henselian), 
and the category of coherent $\mathscr O_X$-modules is
poorly understood  --  one rather deals with coherent $\mathscr O_{X\grot}$-modules, which are well-behaved. 

Nevertheless, if $X$ is {\em good}, which means that
every point of $X$ has an affinoid neighborhood,
then $\mathscr O_X$ is coherent, 
its stalks are noetherian (and even excellent), and the category of coherent  $\mathscr O_{X\grot}$-modules
is equivalent through the restriction functor
to that 
of coherent $\mathscr O_X$-modules. This is the reason why many properties are first defined and studied for good analytic spaces, 
and thereafter extended to arbitrary analytic spaces by some ``G-gluing" process. Note however that the class
of good spaces is quite broad: it contains affinoid domains, analytification of schemes of finite type, generic fibers of both affine
and proper formal schemes\ldots ; and any open subset of a good space is still good. 

From now on we shall simply say ``coherent sheaf on $X$" for ``coherent $\mathscr O_{X\grot}$-module", and write 
$\mathscr O_X$ instead of $\mathscr O_{X\grot}$. But if $X$ is good, if $\mathscr F$ is a coherent sheaf
on $X$, and if $x\in X$, we shall use the notation
$\mathscr F_x$ to denote the stalk at $x$ of $\mathscr F$
\emph{viewed (by restriction to the
category of open subsets)
as a sheaf on the ordinary topological space underlying $X$}.

\subsubsection{Naive flatness}\label{sss-naive-flatness}
Let~$\phi \colon Y\to X$ be a morphism of
{\em good} $k$-analytic spaces
and let~$\mathscr  F$ be a coherent
sheaf on $Y$. Let us say that~$\mathscr  F$ is
{\em naively}
$X$-flat at a point~$y$ of~$Y$ is~$\mathscr  F_y$ is a flat~$\mathscr  O_{X,\phi(y)}$-module, exactly like in the rigid setting (if~$\mathscr  F=\mathscr  O_Y$
we simply say that~$Y$ is {\em naively~$X$-flat at~$y$}, or that~$\phi$ is naively flat at~$y$). 
We immediately face a big problem. Indeed, in this context, the use of completed tensor products
does not just make proofs of stability 
of naive flatness under base change or ground
field extension more complicated, as it does in the rigid case: naive flatness 
is actually {\em not}
stable under base change nor ground field extension.

Let us describe a counter-example, which had been suggested to us by Michael Temkin.
Roughly speaking, it is due to a boundary phenomenon: it consists
of the embedding into the affine plane of a curve which is drawn on some bi-disc and cannot be extended
to a larger disc; the problem occurs at the unique boundary point of the curve. 
We are now going to give some details; the reader will find proofs of what follows in
\ref{s-counter-example}. Choose~$r>0$
and ~$f=\sum a_iT^i\in k[\![T]\!]$ a power series whose radius of convergence is
exactly~$r$,
and let~$Y$ be the closed one-dimensional~$k$-disc of radius~$r$. The Shilov boundary of~$Y$ consists of
one point~$y$ (the one that corresponds to the semi-norm~$\sum b_i T^i\mapsto \max |b_i|r^i$).
Denote by~$\phi$ the morphism~$({\rm Id}, f):Y\to \A^{2,\rm an}_k$ and by~$X$ the closed analytic domain of~$\A^{2,\rm an}_k$ defined
by the inequality~$|T_1|\leq r$; note that~$\phi(Y)\subset X$; more precisely,~$\phi(Y)$ is the Zariski-closed subset of~$X$ defined by the equation~$T_2=f(T_1)$.
One can show that~$\mathscr  O_{\A^{2,\rm an}_k,\phi(y)}$ is a field: this is due to 
the fact that~$\phi(Y)$ cannot be extended to a curve defined {\em around}~$\phi(y)$, because the radius of convergence of~$f$ is exactly~$r$.
As a consequence,~$\phi$ is 
naively flat at the point~$y$. Now:

\begin{itemize}
\item[$\bullet$] ~$Y=\phi\inv (X)\to X$ is a closed immersion of a one-dimensional space in a purely two-dimensional space,  hence is not
naively flat at~$y$; 

\item[$\bullet$] if~$L$ is any complete extension of~$k$ such that~$Y_L$ has an~$L$-rational point~$y'$ lying
above~$y$, then~$\phi(y')$ belongs to the topological interior of~$X_L$ in~$\A^{2,\rm an}_L$ (because~$\phi(y')$ is a rigid point); therefore~$\phi_L : Y_L\to \A^{2,\rm an}_L$ is,
around~$y'$, a closed immersion of a one-dimensional space in a purely two-dimensional space, hence is not
naively flat at~$y'$. 
\end{itemize}

The above counter-example is, in some sense, archetypal, because boundary phenomena are actually the only obstruction for
naive flatness to be stable under base change: if~$Y\to X$ is a morphism between good~$k$-analytic spaces 
and if~$y$ is a point of~$Y$ at which~$Y\to X$ is inner, 
then every coherent sheaf~$\mathscr  F$ on~$Y$
which is naively~$X$-flat at~$y$
remains so after any
good
base change (Theorem \ref{thm-flat-naiveflat}). 
Our proof is based upon an analytic variant of Raynaud-Gruson's 
\emph{d\'evissages} \cite{raynaud-g1971},
which is developed in Chapter \ref{DEV}; let us mention
that this result had already been proved by Berkovich (in a completely different way) in some 
unpublished
work about flatness. 
It is
now clear why such problems cannot occur in the rigid setting: this is because boundary points are {\em never} rigid. 

\subsubsection{Our definition of flatness}
To overcome the problem we have just mentioned,
we define flatness as follows (\ref{def-flat-good},
\ref{def-flat-gen} 
\ff).
Let~$Y\to X$ be a morphism between~$k$-analytic spaces, let~$y$ be a
point of $Y$, and let $\mathscr  F$ be a coherent 
sheaf on $Y$. 

\begin{itemize}[label=$\bullet$]

\item If $Y$ and~$X$ are good, we say that~$\mathscr  F$ is~$X$-flat at~$y$ if
it 
is naively~$X$-flat at~$y$ and {\em if it remains so after any
good base change and any ground field extension}. 

\item In general, $\mathscr  F$ is~$X$-flat at~$y$ if there exists a good analytic domain~$V$ of~$Y$ containing~$y$
and a good analytic domain~$U$
of~$X$ containing the image of~$V$ such that~$\mathscr  F_{|V}$ is~$U$-flat at~$y$
(and if it is the case then this holds for every~$(V,U)$ as above). 

\end{itemize}
These definitions might seem somehow ad hoc, 
and not so easy to check nor to use. 
But we hope that this
memoir will convince the reader that
they 
provide a notion of flatness which behaves exactly as expected and is quite
convenient to work with. 

\section{Loci of validity}

\subsection{General presentation of the problem}\label{ss-loci-general}
We now turn back to the situation
we have described at the beginning of the Introduction. 
That is, we are given 
a morphism~$\phi\colon Y\to X$
between~$k$-analytic spaces, an object~$D$
on~$Y$, and we want
to understand the variation of properties of~$D_x$ for~$x$ going through~$X$. 
This may be investigated by considering
two kinds of problems. 

\begin{itemize}
\item[$\bullet$] {\em Problems on the source}. What can be said about the set of points $y\in Y$ such that $D$ 
satisfies a given (local or punctual) property $\mathsf P$ fiberwise at~$y$; \ie,
$D_{\phi(y)}$ satisfies~$\mathsf P$
at~$y$ as an object living on the~$\hr{\phi(y)}$-analytic space
~$Y_{\phi(y)}$ ? 

\item[$\bullet$] {\em Problems on the target}. What can be said about the set of points $x\in X$ such that $D_x$ 
satisfies a given global property $\mathsf Q$, as an object living on the $\hr x$-analytic space $Y_x$ ? 

\end{itemize}

In scheme theory, both kinds of problems are addressed. More precisely, 
one first proceeds to a study on the source, and one uses the latter 
\emph
{together with Chevalley's constructibility theorem}
in order to understand what happens on the target. Note that Chevalley's theorem 
itself can be described in a somehow pedantic way as an answer to the ``problem on the target"
for $\mathsf Q$ being the {\em non-emptiness}
property (of a scheme). 

There is kind of an avatar of Chevalley's theorem in analytic geometry, 
due to L. Lipshitz and Z. Robinson~\cite{lipshitz-rob2000}. 
It asserts that the image of a morphism
between strictly affinoid spaces is constructible
with respect to a 
broad class of functions, far larger than that
of global analytic ones. But unfortunately constructible subsets
in the sense of Lipshitz and Robinson are (up to now) poorly understood and not that tractable; 
for that reason, in this memoir we will
mainly {\em ignore}
the problems on the target, 
and focus on what happens {\em on the source}. 

\begin{rema}
The image of an {\em overconvergent} 
morphism between strictly affinoid spaces is constructible 
with respect to a class of functions which is smaller
than that of Lipshitz and Robinson, but still larger than that 
of global analytic ones; this had been proved by Schoutens 
in~\cite{schoutens1994}. Recently, F. Martin
\cite{martin2015} has given a purely 
geometric version of Schoutens's statements
and proofs (it involves some particular
finite sequences of blow-ups instead of
the class of recursively 
defined functions considered by Schoutens).
\end{rema}

\subsection{Loci of validity in
the source space}\label{intro-locivalid}

Let $Y\to X$ be a morphism
of $k$-analytic spaces, and let $\mathscr F$
be a coherent sheaf  
on $Y$. We prove the following (assertion (1)
is
Theorem ~\ref{thm-flatloc-zaropen},  
and the other ones are part of Theorem
\ref{thm-constloc-main} together with the fact
that quasi-smoothness is equivalent to flatness and
fiberwise regularity, see Th.\ref{thm-main-qsm} (1)):

\begin{enumerate}[1]

\item The set
of points of~$Y$ at which~$\mathscr  F$ is~$X$-flat is Zariski-open. 

\item The set of points of~$Y$ at which~$Y$ is fiberwise geometrically regular (\resp
geometrically reduced, \resp geometrically normal,  \resp Gorenstein, \resp  complete intersection) is
locally constructible. 

\item The set of points of~$Y$ at which~$Y$ is~$X$-flat
and fiberwise geometrically regular (\resp Gorenstein, \resp  complete intersection)
is Zariski-open; if~$Y$ is 
relatively equidimensional
over $X$ (\ref{ss-dim-nofield}; empty fibers are allowed, see \ref{ss-def-dimloc}),
the set of points of~$Y$ at which~$Y$ is~$X$-flat and fiberwise
geometrically reduced (\resp geometrically normal) is Zariski-open. 

\item The set of points of~$Y$ at 
which~$\mathscr  F$
is fiberwise Cohen-Macaulay (\resp $S_m$) is
locally constructible. 

 \item The set of points of~$Y$ at which~$\mathscr  F$ is~$X$-flat
and fiberwise Cohen-Macaulay
is Zariski-open. 

\end{enumerate}
Let us give some explanations about our terminology.

\subsubsection{}
A subset $E$ of $Y$ is called \emph{constructible}
if it is a finite boolean combination Zariski-open
subsets, and \emph{locally constructible} if every point of $Y$
has an open neighborhood $U$ such that $U\cap E$ is constructible
in $U$. 
But when
$Y$ is finite-dimensional, every locally constructible
subset of $Y$ is constructible (Proposition 
\ref{prop-cons-gloc}; see \ref{contrex-gcons}
for a counter-example in the infinite-dimensional 
case).

\subsubsection{}
Assume that $Y$ is good. It is said to be regular
(\resp reduced, \ldots) at $y$ if
the local ring $\mathscr O_{Y,y}$ is regular (\resp reduced,\ldots). The coherent
sheaf $\mathscr F$
is said to be Cohen-Macaulay at $y$
if the $\mathscr O_{Y,y}$-module $\mathscr F_y$ is Cohen-Macaulay. 

\subsubsection{}
The space $Y$ is no longer assumed to be good. 
Then $Y$ is said to be regular (\resp reduced, \ldots)
at $y$ if there exists a good analytic domain $V$ of $Y$ containing
$y$ which is regular (\resp reduced, \ldots)
at $y$, and it then holds for \emph{every} good analytic
domain containing $y$ 
(Lemma \ref{valid-at-concrete}); and
$\mathscr F$ is said to be Cohen-Macaulay
at $y$  if there exists a good analytic domain $V$ of $Y$ containing
$y$ such that the restriction of $\mathscr F$ to $V$ is Cohen-Macaulay
at $y$, and it then holds for \emph{every}  good analytic
domain containing $y$ 
(Lemma \ref{valid-at-concrete} again). 

\subsubsection{}
The properties of being Gorenstein or Complete
Intersection (for $Y$) or Cohen-Macaulay
(for $\mathscr F$) at a given point
are preserved by arbitrary ground field extensions. 
This is not the case in general
for the property of being regular, normal or reduced, and
when one of them
holds at
a point and 
remains valid after arbitrary ground field extensions it is said to
hold \emph{geometrically} (to ensure geometric validity, there is in fact no need to check
the property over
all possible extensions: it suffices that it holds after a single
\emph{perfect} extension). 
The behavior of
algebraic properties under ground field extension
is described in full detail in
\ref{s-ground-field-ext}.

\begin{rema}
The property for a morphism
to be flat and fiberwise geometrically regular at
a given point is of fundamental importance. It is 
called {\em quasi-smoothness}, and can also be defined
using an analytic avatar of the Jacobian criterion; 
this is the main topic of Chapter~\ref{QSM}. 
\end{rema}

\subsection{Loci of validity in the target space}
As we explained above, the general description of
the
set of points of the target whose fiber satisfies
some given property seems out of reach. But
we have been able to address that kind of question
in two very specific situations. 

\subsubsection{The case of a proper map}
Let $Y, X$ and $\mathscr F$ be as above, and assume moreover
that $Y\to X$ is \emph{proper}. Then under
this assumption 
we have kind of an analytic
Chevalley's theorem: indeed, by Kiehl's result on cohomological finiteness 
of proper morphisms, the image in $X$ of any Zariski-closed subset of
$Y$ is Zariski-closed (\ref{ss-kiehl-proper}); and one can deduce from this
that the image in $X$ of any locally constructible subset of $Y$ is locally
constructible (Theorem \ref{thm-chevalley-proper}).
The following then follow formally from
the aforementioned results on the source (this is part of Theorem
\ref{thm-locus-target}). 

\begin{enumerate}[1]

\item The set of points $x$ of $X$ such that 
$\mathscr F$ is $X$-flat at every point of $Y_x$ is a Zariski-open subset of $X$. 

\item The set
of points $x$ of~$X$
such that $Y_x$ is geometrically regular (\resp Gorenstein, \resp Complete Intersection)
at all of its points is a locally constructible subset of $X$. 

\item The set
of points $x$ of~$X$
such that $Y$ is $X$-flat at every point 
of $Y_x$ and $Y_x$ is geometrically regular (\resp Gorenstein, \resp Complete Intersection)
at all of its points is a Zariski-open subset of $X$. 

\item The set
of points $x$ of~$X$ such that 
the restriction
of $\mathscr F$ to
$Y_x$
is Cohen-Macaulay at each point of $Y_x$
is a locally constructible 
subset of $X$. 

\item The set
of points $x$ of~$X$ such that $\mathscr F$ is $X$-flat at each
point 
of $Y_x$ and
the restriction
of $\mathscr F$ to
$Y_x$
is Cohen-Macaulay at each point of $Y_x$
is a Zariski-open subset of $X$.

\item Assume moreover that $Y$ is relatively equidimensional over $X$.
Then the set
of points $x$ of~$X$  such that $Y$ is $X$-flat at every point 
of $Y_x$ and $Y_x$ is geometrically regular (\resp geometrically 
reduced, \resp geometrically normal)
at each of its points is a Zariski-open subset of $X$. 
\end{enumerate}

\subsubsection{The image of a morphism between affinoid spaces}
We prove
(Theorem~\ref{th-image-compact})
that
if $\phi\colon Y\to X$
if a flat morphism 
between affinoid spaces, then $\phi(Y)$
is a compact analytic domain (otherwise said, 
a finite union of affinoid domains) of $X$. 
In the strictly analytic case, 
this is nothing but the aforementioned
result of Raynaud (\cf \cite{frg2}, \Cor 5.11). But our methods
are completely different 
and provide a new proof
of his theorem, which does not involve any formal model. 

We first give a
direct proof  when the relative dimension of $\phi$
is zero (Proposition \ref{prop-image-cm}); it is based upon
Temkin's theory of reduction of analytic germs
and kind of
Chevalley theorem in the framework
of (graded) Riemann-Zariski spaces
(Theorem \ref{thm-elim-graded}). 
We then handle the general case (with no assumption
on the relative dimension of $\phi$)
by reducing through
a suitable ground field extension to the case where
$k$ is non-trivially valued and both $Y$ and $X$ are strict,
and then by showing the following, which
seems us to be of independent interest: there exist a strictly $k$-affinoid space $Z$, 
a flat morphism $\psi\colon  Z\to X$ of relative dimension $0$, and an $X$-map $Z\to Y$,
such that $\psi(Z)=\phi(Y)$ (Theorem \ref{thm-multisections-global}). Otherwise said, the image
of any
flat map between strictly $k$-affinoid spaces is covered by a flat, quasi-finite multisection
of the map. 
%
%\item
%\end{itemize}
%We are able to describe the validity 
%\end{rema*}

\subsection{From the target and the fibers to the source}
Let $\mathscr Y\to \mathscr X$ be a
morphism of noetherian schemes. If
$\mathscr X$ is regular (\resp reduced, \resp normal, \resp Gorenstein, \resp Complete Intersection), 
$\mathscr Y\to \mathscr X$ is flat, and its fibers are 
regular (\resp reduced, \ldots), then $\mathscr Y$ is regular (\resp reduced, \ldots). Analogously, 
let $\mathscr F$ be a coherent sheaf on $\mathscr Y$. If $\mathscr X$ is Cohen-Macaulay, $\mathscr F$
is flat over $\mathscr X$, and the restriction of $\mathscr F$ to every fiber of $\mathscr Y\to \mathscr X$ 
is Cohen-Macaulay, then $\mathscr F$ is Cohen-Macaulay. 

At the end of this memoir (Chapter \ref{c-TARG}),
we develop a general and systematic
method to transfer such theorems
from algebraic geometry to analytic geometry; in particular, 
the analytic counterparts of the above statements all hold
(see
Theorem \ref{thm-base-fibers-concrete}).

\section{About our proofs}

\subsection{Technical obstacles in analytic geometry}
Our general strategy is of course to adapt
to the Berkovich setting what had been done about relative
properties in algebraic 
geometry, and sometimes in complex
analytic geometry. For instance, our investigation
of quasi-smooth morphism is inspired by the study
of smooth maps by Bosch, L\"utkebohmert and Raynaud
in \cite{bosch-l-r1990}, our theory of d\'evissages
follows that of Raynaud and Gruson in \cite{raynaud-g1971}, and
our proof of the Zariski-openness of the flat locus
is mutatis mutandis the same as that of Kiehl in the complex
analytic setting \cite{kiehl1967}.

But most of the time (and that is the case for quasi-smoothness
as well as for the d\'evissages), techniques coming from the Grothendieck 
school
\emph{cannot}
be applied straightforwardly in our setting. Let us now give some
examples of obstacles we had to face, and then quickly explain
the way we overcame them -- we
hope that the techniques and methods developed for this purpose will
be useful in other circumstances. 

\begin{enumerate}[1]

\item If $Y$ is an analytic space, the Zariski topology of an analytic domain $V$ of $Y$
is in general strictly finer than the one inherited from the Zariski-topology of $Y$ -- even
if $V$ is a Zariski-open subset of $Y$ (note that
any infinite, discrete and closed subset of $\A^{1,\mathrm {an}}_k$ consisting
of rigid points is Zariski-closed in $\A^{1,\mathrm {an}}_k$, but it does not come from a Zariski-closed
subset of $\P^{1,\mathrm {an}}_k$). And analogously, if $Y\to X$ is a morphism of $k$-analytic spaces, the Zariski-topology
of a fiber $Y_x$ at a \emph{non-rigid}
point is also in general finer than the one inherited from the Zariski-topology of $Y$.

\item If $y$ and $z$ are two points of a good analytic space $Y$ and if $y$ belongs
to the Zariski closure of $z$ in $Y$, there is no way to relate in a simple way the local rings
$\mathscr O_{Y,y}$ and $\mathscr O_{Y,z}$ (on a scheme, the latter would be a localization
of the former). 

\item If $Y\to X$ is a morphism
of good analytic spaces and if $y$ is
a point of $Y$
whose image in $X$
is denoted by $x$,
the  
comparison between
$\mathscr O_{Y_x,y}$ and $\mathscr O_{Y,y}$ is
not that clear: for instance,
even if $\mathscr O_{X,x}$ is a field (which means that $Y_x$ should be thought of as kind 
of a generic fiber), the local ring $\mathscr O_{Y_x,y}$  is in general a huge $\mathscr O_{Y,y}$-algebra, 
though on
a scheme it would be equal to  $\mathscr O_{Y,y}$. 

\item If $\phi \colon Y\to X$ is a morphism between two
$k$-affinoid spaces, we have already explained
in \ref{ss-loci-general}
that $\phi(Y)$
is not easily understandable in general. 
In particular, when $Y$ is of dimension $d$
and $\phi$ of pure relative dimension $\delta$
for some $d$ and $\delta$, there is no reason why $\phi(Y)$ should
be contained in a $(d-\delta)$-dimensional Zariski-closed
subset of an analytic domain of $X$ (in scheme theory, one would
simply take $\overline{\phi(Y)}$). 
\end{enumerate}

\subsubsection{}Obstacle (1) is not so big a problem: one can more or less
overcome it because Zariski-closedness (or openness) is a G-local property 
and the irreducible components behave reasonably with respect 
to analytic domains. Moreover, some  work
had already been carried out by the author
in \cite{ducros2007} to remedy the fact that the Zariski-topology on 
a fiber is ``too fine"; see for instance
section 4 of \opcit, whose results are used
repeatedly in this memoir. 

\subsubsection{}Most of the time, we shall overcome obstacle (2) by working
with the respective images $\eta$ and $\zeta$
of $y$ and $z$ on the scheme $\mathscr Y:=\spec \mathscr O_Y(Y)$, 
which have the required property (\ie, $\mathscr O_{\mathscr Y,\zeta}$ is
a localization of $\mathscr O_{\mathscr Y,\eta}$), and then go back to 
our original space $Y$ by using some GAGA results. 

\subsubsection{}\label{ss-obstacle3}
Obstacle (3) is probably a priori
the most harmful 
for our purposes.
Indeed, the study of relative
properties in EGA,
rests, among other things, on the 
technique  of
``spreading out from generic fiber"
which we describe roughly.
One starts from a morphism of noetherian schemes $\mathscr Y\to \mathscr X$, and
with a point $\eta$ of $\mathscr Y$ at which some property $\mathsf P$ (of the scheme itself,
of a coherent sheaf on it, \ldots)
is satisfied fiberwise; let $\xi$ be the image of $\eta$. 
Then if one sets $\mathscr T=\mathscr Y\times_{\mathscr X}\overline{\{\xi\}}_{\mathrm{red}}$, the fiber $\mathscr Y_\xi$ is equal to the 
\emph{generic fiber}
of the map $\mathscr T\to \overline{\{\xi\}}_{\mathrm{red}}$; hence the local ring of this fiber at $\eta$ 
is equal to $\mathscr O_{\mathscr T,\eta}$; \emph{note that \textnormal{(3)}
tells precisely
that this step would fail in analytic geometry.}

The property $\mathsf P$
thus holds now absolutely (and not only fiberwise) on $\mathscr T$ at $\eta$; 
moreover the scheme $\mathscr T$ itself, as well as the restriction 
to $\mathscr T$ of any of the coherent sheaves possibly involved in our situation, are flat
over $\overline{\{\xi\}}_{\mathrm{red}}$ at $\eta$. This might help to ``spread"
fiberwise
validity of $\mathsf P$ from $\eta$ to a dense open subset of $\overline{\{\eta\}}$. 

\subsubsection{}
We remedy obstacle (3) as  follows. Let $Y\to X$ denote a morphism
of $k$-affinoid spaces, let $x$ be a point of $X$ whose local ring is a field, 
and let $y$ be a point of $Y_x$. We prove (Theorem \ref{thm-localring-generic})
that if $y$ does not belong
to the relative boundary of $Y$ over $X$, then the map $\spec \mathscr O_{Y_x,y}\to \spec \mathscr O_{Y,y}$
is flat (with complete intersection fibers). And flatness
suffices to ensure the transfer of all
familiar algebraic properties from $\mathscr O_{Y_x,y}$ to 
$\mathscr O_{Y,y}$, which is most of the time what is actually needed  (for instance in 
\ref{ss-obstacle3} above, the conclusion that $\mathsf P$ is satisfied on $\mathscr T$ at $\eta$
only requires the transfer of the property $\mathsf P$ from $\mathscr O_{\mathscr T_\xi,\eta}$
to $\mathscr O_{\mathscr T,\eta}$, and not the equality of the rings). 

\subsubsection{}
At the end of the memoir, we
remedy obstacle (4)
by showing the following: if $\phi \colon Y\to X$, $d$ and $\delta$ are as in (4),
then there exists a non-empty affinoid domain $V$
of $Y$ and a purely $(d-\delta)$-dimensional Zariski-closed subset $S$
of an analytic domain of $X$ such that $\phi(V)\subset S$
(\ref{thm-weak-chevalley}). This is slightly weaker than the 
``dream
property" discussed in (4) (one does not control
the whole of $\phi(Y)$, but only the image of an affinoid domain
which can possibly be very small), but this can nevertheless
be useful for some purposes, like carrying out an induction on the dimension
of the base space (this is the way we use it in the proof
of Theorem \ref{thm-base-fibers-concrete}).

\subsection{About the previous results by Kiehl}
In  \cite{kiehl1968}, Kiehl has established some
analogous
results for a morphism between two affinoid {\em rigid} spaces,
but they do not 
{\em a priori}
imply our theorems for the following reason.

In the rigid
analytic context, one only deals with fibers over {\em rigid} points.
To be sure, any 
point of an analytic
space can be made rigid after a suitable ground field extension, 
and it follows from the author's previous work
\cite{ducros2009}
that all properties involved can be checked after any ground field
extensions. 
But the combination of these remarks
and of Kiehl's theorems
does not yield directly 
the corresponding statements in the Berkovich setting,  {\em  unless
one knows that the formation
of the fiberwise validity
locus of a given property in rigid-analytic geometry commutes with scalar extension}. 
And this
has not been addressed by Kiehl, whose methods do not seem
to apply straightforwardly to such questions: indeed, these methods are very ``algebraic", 
and it is not clear -- at least to the author -- how one could use them to 
deal with scalar extension, which involves completion operations. 

Because of that, and also for
the reader's convenience, 
we have chosen to write self-contained
and purely ``Berkovich'' proofs; hence we
{\em recover}
and extend Kiehl's results. 

\section*{Acknowledgements}
I started to think about the topics this paper is devoted to when
I was asked some questions about flatness by Brian Conrad and Michael Temkin,
for their work \cite{conrad-tXXX}. I realized quite quickly that answering them
would take much more time than
I had originally expected... and it eventually gave rise to the present work.
I thus would like express my gratitude to
Conrad and Temkin for having given
to me the initial inspiration -- and also for a lot of fruitful discussions since then. 

I also would like to thank warmly the anonymous referee for invaluable work
when reading several versions of the manuscript extremely carefully. I greatly 
benefited of his/her (incredibly) numerous thorough comments, 
which helped me to significantly improve this work, and I am very grateful for this.

\chapter{Background material}

This chapter presents our general conventions in topology, algebra and algebraic geometry, and
then provides some basic reminders in analytic geometry. The reader who already knows well Berkovich's
theory
may skip most of it, but should nevertheless have a look at \ref{ss-space-nofield} \ff~for our conventions
about analytic spaces without mention of a ground field, and at
section \ref{ss-space-nofield}
and especially \ref{ss-stalks-fibers} for our conventions and notation about coherent sheaves and their stalks. 
Be aware that we depart from Berkovich's terminology in two cases: 

\begin{itemize}[label=$\bullet$]
\item We say ``boundaryless" instead of ``closed" (\ref{ss-boundary}). 

\item We say ``finite at a point" instead of ``quasi-finite at a point", ``locally finite" instead of
``quasi-finite" (\ref{ss-finite-maps}), 
and we have our own definition for ``quasi-finite" (\ref{ss-quasi-finite}, Remark \ref{rem-depart-qf}).
\end{itemize}

Let us also mention that we shall make much use
in this memoir of \emph{graded} commutative algebra,
after Temkin
\cite{temkin2004}. 
Since this theory essentially
consists of a
somehow tedious
transcription of classical 
commutative algebra with
the words ``graded" or ``homogeneous" added almost everywhere,
we have chosen to write the corresponding
reminders (almost without proofs, though some of them are sketched)
not in this chapter, but in Appendix \ref{c-graded}
at the end of the memoir. The reader may refer to it
if needed. Let us simply say here that for us, ``graded" will always mean
``$\R_+\gpm$-graded", with multiplicative graduation.

\section{Prerequisites and basic conventions}

\subsection{Prerequisites}
The understanding of this memoir requires
of course a robust general knowledge of commutative
algebra, algebraic geometry and non-archimedean analytic geometry; let us make this precise.

In commutative algebra, we will use freely the ``usual"
notions about commutative rings (and especially local noetherian rings) and
modules over the latter: flatness, Krull dimension, depth and codepth, regularity, properties
like $R_m$ and $S_m$ (due to Serre), Cohen-Macaulay, Gorenstein, or (local) complete intersection; 
and also Grothendieck's crucial
concept of excellence. Possible references on these topics are Matsumura's textbook 
\cite{matsumura1986}
and part of EGA  and SGA: \cite{ega41}, Chapter 0, \S15, \S16 and \S17; \cite{ega42}, \S 6 and \S7; and \cite{sga1}, Expos\'e V. 
We will also assume the reader to be familiar with elementary valuation theory; see for instance \cite{bourbaki1985}, Chapitre~VI.

In algebraic geometry, we will use scheme theory in the spirit
of EGA and of
the French school in algebraic geometry, and the reader will
be assumed to
master the corresponding language.
Parts of this memoir have
been more specifically inspired by the treatment of
smoothness in the book \cite{bosch-l-r1990}
on N\'eron models by Bosch, L\"utkebohmert and Raynaud, by the first part
of Raynaud and Gruson's seminal work on flatness \cite{raynaud-g1971}, 
and by \cite{ega43} \S12
which is devoted to relative properties. Familiarity with these
texts might therefore 
be helpful, but is not strictly needed. 

In
analytic geometry, 
the basics of Berkovich's approach as exposed in~\cite{berkovich1990}
and \cite{berkovich1993} \S1  will be considered known.
We will also make much use of the works by 
the author on dimension theory \cite{ducros2007}
and on ``commutative algebra properties" of analytic spaces~\cite{ducros2009}, as well 
as of Temkin's theory of the reduction of germs \cite{temkin2004}. Some acquaintance
with these topics is then recommended;  nonetheless, we will 
recall all related definitions
and statements
that are needed for our purposes.

\subsection{Conventions in algebra}\label{ss-conv-alg}
In this memoir, all rings (and algebras) are commutative
and have a unit element; 
morphisms of rings (and algebras) always respect the unit elements. 
The group of invertible elements of a ring $A$ will be denoted by $A^\times$. 
The acronyms CM and CI will stand respectively for ``Cohen-Macaulay" and ``complete intersection". 

If $M$ is a module over a ring $A$ and if $I$ denotes the annihilator of
$M$, the Krull dimension of 
$M$ will be by definition the Krull dimension of the ring $A/I$; in more geometric terms, this is
the dimension of the support of $M$ on $\spec A$. \index{dimension! of a module}

We shall
often use for short the {\em multi-index} notation, which consists
of the following.  
Let $\Lambda$ be a multiplicative
commutative monoid, let $n$ be a non-negative integer 
and let $\tau=(\tau_1, \ldots, \tau_n)$ be an $n$-uple
of elements of $\Lambda$ (the reader should have two
examples in mind: the case where $\Lambda=\R_+\gpm$ and the case where $\Lambda$
is the multiplicative monoid of an algebra of polynomials or power series
over some ground ring in the indeterminates $\tau_j$ ). 
Let $I=(i_1, \ldots, i_n)\in \Z^n$ be such that 
$i_j\geq 0$ as soon as $\tau_j$ is not invertible. Then the product
$\prod \tau_j^{i_j}$ will be simply denoted by $\tau^I$. 
If all $\tau_j$'s are invertible, we will set $\tau^{-1}=(\tau_1^{-1}, \ldots, \tau_n^{-1})$. 

\subsection{Conventions in topology}\label{ss-remind-topo}
We shall use the terminology of Bourbaki \cite{bourbaki1971}. A topological space $X$ will be
called \emph{quasi-compact}\index{topological space!quasi-compact}\index{quasi-compact (topological space)}
if every open covering of $X$ admits a finite sub-covering; it will be called
\emph{compact} if it is quasi-compact \emph{and Hausdorff}. \index{topological space!compact}\index{compact topological space}

Let $X$ be a topological space. It will be called:

\begin{itemize}[label=$\bullet$] 

\item \emph{locally compact}
if $X$ is Hausdorff and every point of $X$ has a compact neighborhood in \index{topological space!locally compact}\index{locally compact topological space}
$X$; 

\item \emph{countable at infinity}\index{topological space!countable at infinity}\index{countable at infinity (topological space)}
if $X$ is Hausdorff and $X$ is the union of (at most) countably many compact subsets; 

\item \emph{paracompact} if $X$ is Hausdorff and
every open covering of $X$ can be refined into a
locally finite open covering. \index{topological space!paracompact}\index{paracompact topological space}
\end{itemize}

If $X$ is locally compact, it is paracompact if and only if it is the disjoint union of
open subsets that are countable at infinity
(\cite{bourbaki1971}, Chapitre I, \S 9, \no 10, Th 5). If this is the case, then 
for every basis of neighborhoods $\mathscr B$ of $X$, every open covering of $X$
can be refined into a
locally finite covering consisting of elements of $\mathscr B$
(this is what the proof of \cite{bourbaki1971}
~actually shows). 

A continuous map $p\colon Y\to X$\index{continuous map!proper}\index{proper continuous map}
between two arbitrary topological spaces
is said to be \emph{proper}
if it is universally closed; \ie, for every continuous map $Z\to X$ and every
closed subset $F$ of $Y\times_X Z$, the image of $F$ in $Z$
(by the second projection) is closed.
(Be aware that 
$p$ is \emph{not}
required to be separated;  this is Bourbaki's convention and we have chosen to follow it despite the
inconsistency with the definition of properness in algebraic geometry, because some of our
results
actually hold for maps between analytic spaces
that are topologically proper in Bourbaki's sense
without any separatedness assumption; see Theorem \ref{th-image-compact}
and Theorem \ref{th-image-equi}). 
If $p$ is proper, 
$p\inv(K)$ is quasi-compact for every quasi-compact subset $K$ of $X$ (\cite{bourbaki1971}, Chapitre I, 
\S 10, \no 2, \Prop 6). 
If $Y$ is Hausdorff and $X$ is locally compact, $p$ is proper if and only if
$p\inv(K)$ is compact for every compact subset $K$ of $X$ (\cite{bourbaki1971}, Chapitre I, 
\S 10, \no 3, \Prop 7).

If $E$ is any subset of a topological space $X$, the closure
of $E$ inside $X$ will be denoted by $\adht EX$.

\subsection{Conventions in algebraic geometry}\label{ss-conv-scheme}
The word ``scheme" will be understood
here {\em without particular assumption}. We shall always
make
precise when the scheme we are working with is separated, noetherian, 
excellent, reduced, of finite type over some ground ring, etc. 

If $X$ is a scheme and if $x$ is a point of $X$, the maximal ideal of $\mathscr O_{X,x}$ will be denoted by $\mathfrak m_x$,\label{IN-sch-ox}\label{IN-sch-mx}
and its
residue field by $\kappa(x)$. If $\phi \colon Y\to X$\label{IN-sch-kx}
is a morphism of schemes, the {\em scheme-theoretic} fiber of $\phi$ over $x$ will be denoted by $\phi^{-1}(x) $ or $Y_x$\label{IN-sch-yx}
 ; this is
 a $\kappa(x)$-scheme. 

A morphism of schemes $Y\to X$ is called {\em regular} 
if it is flat
and if for every point $x$ of $X$ the fiber $Y_x$ is locally noetherian and
{\em geometrically}
regular; \ie, $Y_x\otimes F$ is regular for any finite, purely inseparable\index{morphism!of schemes!regular}\index{regular morphism (of schemes)}
extension $F$ of $\kappa(x)$.
A morphism of rings $A\to B$ is said to be regular if
 $\spec B\to \spec A$ is regular. 

If $\mathscr F$ is a coherent sheaf on a noetherian scheme $X$, and if  $x$ is a point of $X$, 
we shall denote by $\mathscr F_x$ the stalk of $\mathscr F$ at $x$, and by $\mathscr F_{\kappa(x)}$ the tensor\label{IN-sch-Fx}
product $\kappa(x)\otimes_{\mathscr O_{X,x}}\mathscr F$. If $Y$ is any $X$-scheme, the pull-back of $\mathscr F$
to $Y$ will be denoted by $\mathscr F_Y$. 

\section{Analytic geometry: basic definitions}

\begin{defi}
\label{def-anal-fields}\index{analytic field}\index{field!analytic}
An {\em analytic field}
is a field endowed with an $\R_+\gpm$-valuation for which it is {\em complete}; unless otherwise stated, 
the structure valuation of an analytic field will be denoted by $\abs{\cdot}$. 

\end{defi}

\begin{exem}
Any field 
endowed with the trivial 
valuation is an analytic field. 
\end{exem}

\subsection{}\label{ss-conv-tilde} If $k$
is an analytic field, we shall denote by $\widetilde k$
its \emph{graded}\index{graded!reduction}\label{IN-wkr}
reduction (\ref{ss-gradred-general} and
also \ref{ss-gradred-field}); \ie, 
\[\widetilde k=\bigoplus_{r>0}\{x\in k, \abs x\leq r\}/\{x\in k, \abs x <r\}.\]
For every positive
$r$, the $r$-th summand of $\widetilde k$ will be denoted by $\widetilde k^r$; note than
the usual residue field of $\widetilde k$ is nothing but $\widetilde k^1$. \label{IN-wk1}

If $x$ is any element of $k$ and $r$ any positive number $\geq \abs x$, we shall
denote by $\widetilde x^r$ the image of $x$ in $\widetilde k^r$; if $r=\abs x$
we shall
simply write $\widetilde x$. 

\subsection{}\label{sss-analytic-extension} 
An {\em analytic extension}\index{analytic extension}
of an analytic field
$k$ is an analytic field $L$ together with an isometric embedding $k\hookrightarrow L$; \label{IN-dkl}
for such an $L$ we shall denote by $d_k(L)$ the transcendance
degree of the graded extension $\widetilde k\hookrightarrow \widetilde L$
(see~\ref{ss-index-graded}
for a ``classical" interpretation of this invariant).

\begin{enonce}[remark]{Convention}\index{space!k@$k$-analytic}\index{k@$k$-analytic space}
The notion of a {\em $k$-analytic space}
will always be understood in the sense of Berkovich \cite{berkovich1993} \S 1.
\end{enonce}

\subsection{Topologies}\label{ss-topologies}
Let $k$
be an analytic field and
let $X$ be a $k$-analytic space. The space $X$ comes with a topology, which enjoys
very nice properties (this is one of the distinguished features of Berkovich's theory): every point of $X$
has a basis of open neighborhoods that are 
Hausdorff, locally compact, path-connected, and countable at infinity.

The space $X$
is also
equipped with (set-theoretic) \emph{Grothendieck} topology. Before describing it, let us introduce
some terminology. If $E$ is any subset of $X$ and if $(E_i)$ is a family of subsets of $E$, we shall say that $(E_i)$ 
is a \emph{\textnormal G-covering}
of $E$ if every point $x$ of $E$ has a neighborhood \emph{in E} of the form $\bigcup_{i\in I} E_i$
for some finite set $I$
such that $x\in \bigcap_{i\in I} E_i$.
A subset $V$ of $X$ is called an {\em analytic domain of $X$}\index{analytic domain (of an analytic space)}
if $V$ is G-covered by the affinoid domains of $X$ that are contained in $V$.
Affinoid domains and open subsets of $X$ are analytic domains; any 
analytic domain $V$ of $X$ inherits a canonical structure of a $k$-analytic space. 
If $x$ is a point of $X$, an \emph{analytic neighborhood}\index{neighborhood!analytic}\index{analytic neighborhood}
of $x$ in $X$ will be an analytic domain of $X$ which is a neighborhood of $x$ in $X$. 

The site $X_{\mathrm G}$
is then defined as follows:

\begin{itemize}[label=$\bullet$]

\item Its objects are the analytic domains of $X$. 

\item Its morphisms are the inclusion maps. 

\item Its topology is the so-called \emph{\textnormal G-topology}, whose covering families are
exactly the families $(V_i\subset V)_i$ for $(V_i)$ a G-covering of $V$. \index{G-topology}
\end{itemize}

Any locally finite covering of a Hausdorff analytic domain of $X$ by compact analytic domains is a G-covering.
Any open covering of an open subset of $X$ is a G-covering: the G-topology is finer than the usual topology.
This can be rephrased in a somehow pedantic way by saying that the forgertful functor induces a
morphism of sites from $X_{\mathrm G}$ to $X$.

If $X$ is quasi-compact, any G-covering of $X$ admits a finite sub-covering; in particular, $X$
admits a finite affinoid G-covering. Conversely, if $X$ admits a finite set-theoretic affinoid covering, 
it is quasi-compact because any affinoid domain is compact. 

If $X$ is Hausdorff, it is paracompact if and only if it admits a locally finite affinoid
covering. Indeed, the direct implication
is due to the comments following
the definition of paracompactness given in \ref{ss-remind-topo}
(by taking for $\mathscr B$ the set of finite unions of affinoid domains). 
For the converse implication, suppose we are given a locally finite affinoid
covering $(X_i)_{i\in I}$ of $X$; we may assume that every $X_i$ is non-empty. Let $\mathscr R$
be the finest equivalence relation on $I$ such that $i\mathscr R j$ as soon as $X_i\cap X_j\neq \emptyset$
(if the latter holds we shall say that $i$ and $j$ are elementary equivalent).
Fix $C\in I/\mathscr R$ and fix $i\in C$. For every $n\geq 1$, the set $C_n$ of indices $j\in C$
that are linked to $i$ by a chain of at most $n$ elementary equivalences is finite, because $(X_i)$ is locally finite.
Hence the union $X_{C_n}$ of all $X_j$'s for $j$ running through $C_n$ is compact; 
therefore $X_C:=\bigcup_{i\in C}X_i=\bigcup_n X_{C_n}$
is countable at infinity.
Since $X=\coprod_{C\in I/\mathscr R} X_C$ and since every $X_C$
is open in $X$ (again by local finiteness of the covering $(X_i)$), the space
$X$ is paracompact. 

Let $\phi \colon Y\to X$ be a morphism of $k$-analytic spaces. If $\phi$ is
\emph{topologically} proper (see \ref{ss-remind-topo}) 
then $\phi\inv(V)$ is quasi-compact for every quasi-compact subset $V$ of $X$, and in particular\index{topologically proper morphism}\index{morphism!of analytic spaces!topologically proper}
for every affinoid domain $V$ 
of $X$. Conversely, assume that $\phi\inv(V)$
is quasi-compact for every affinoid domain $V$
of $X$ ; the morphism $\phi$ is then topologically proper. Indeed, its fibers are quasi-compact, so it suffices to prove that it is closed. 
But this can be checked G-locally on $X$, which allows us to assume that $X$ is affinoid. The space $Y$ is then
quasi-compact, hence can be written as a finite union $\bigcup Y_i$ with $Y_i$ affinoid for all $i$. Now $Y_i\to X$ is closed
for every $i$, hence $\phi$ is closed. 

\subsection{}\label{ss-space-nofield}~If $X$ is a $k$-analytic space, and if~$L$ is an analytic extension of~$k$,
we shall denote by~$X_L$ the~$L$-analytic space deduced from~$X$ by extending the ground field to~$L$.
There is a natural map $X_L\to X$ which is surjective, \cf \cite{ducros2007}, \S 0.5. If $A$
is a $k$-affinoid algebra, we shall denote by $A_L$ the $L$-affinoid algebra $A\hotimes_k L$.

\subsection{}\label{ss-analytic-nofield}
An analytic space {\em without mention of any ground field} is a pair~$(k,X)$ where~$k$
is an analytic field and~$X$\index{space!analytic}\index{X@$X$-analytic space}
a~$k$-analytic space; a morphism between two analytic spaces~$(L,Y)$ and~$(k,X)$ consists of
an isometric embedding~$k\hookrightarrow L$ and a morphism~$Y\to X_L$ of~$L$-analytic spaces.
Similarly we define an affinoid algebra
(\resp space) as a pair~$(k,A)$ (resp.~$(k,X)$) where~$k$ is an analytic field and~$A$ a
$k$-affinoid algebra (\resp and~$X$ a~$k$-affinoid space).
While speaking about analytic spaces and morphisms between them, 
we shall of course most of the time omit mention of the fields and the isometric embeddings involved. 
Hence instead of saying ``let $(k,X)$ be an analytic space" we shall say ``let $X$ be an analytic space"
and we shall refer to $k$ as to the \emph{field of definition}
of $X$. 

\subsection{}\label{ss-hrx}
If $x$ is a point of an analytic space, its completed residue field will be denoted by $\hr x$.
We note that if $V$ is an analytic domain of $X$ containing $x$, then $\hr x$ does not depend
whether $x$ is viewed as belonging to $V$ or to $X$.

\subsection{}\label{ss-analytic-nullst}Let
$k$ be an analytic field. 
A point $x$ of a $k$-analytic space is called
\emph{rigid}\index{point!rigid}\index{rigid point}
if $\hr x$ is a finite extension of $k$. 

Assume that $\abs{k\gpm}\neq \{1\}$. 
By the analytic \emph{Nullstellensatz}
(\cite{bosch-g-r1984}, 6.1.2, \Cor 3), any
non-empty strictly $k$-affinoid space has a rigid point; 
this implies that every non-empty strictly $k$-analytic
space (hence in particular any non-empty boundaryless $k$-analytic space)
has a rigid point.

\subsection{}An analytic space $X$ is said to be {\em good} if every point of $X$ has \index{analytic space!good}\index{good analytic space}
an affinoid neighborhood, hence a basis of affinoid neighborhoods. 

\subsection{}\label{ss-def-xanalytic} Let $X$ be an analytic space, \index{space!X@$X$-analytic}
and let $k$ be its field of definition. 
An
{\em $X$-analytic space}
is a $k$-analytic space $Y$ together with a morphism $Y\to X$ of $k$-analytic
spaces; we emphasize that \emph{$Y$ and $X$ have the same field of definition}.
For such a space $Y$ and for $x$ a point of $X$, 
we shall denote by $Y_x$ the fiber of $Y$ over $x$; this is an $\hr x$-analytic space.
If the morphism $Y\to X$ has been given a name, say $\phi$, this fiber will
also be denoted by $\phi^{-1}(x)$.

\subsection{}\label{ss-fiber-products}
Be aware that the category of analytic spaces in the sense of \ref{ss-analytic-nofield}
does not
admit fiber products in general. For instance, $\mathscr M(\C_p)\times_{\mathscr M(\Q_p)}\mathscr M(\C_p)$ 
does not exist as an analytic space. But if $$\xymatrix{
Y\ar[rd]&&Z\ar[ld]\\
&X&
}$$ is a diagram in the category of analytic spaces and if $Y$ or $Z$ is $X$-analytic, then $Y\times_X Z$
does exist
in the category of analytic spaces (and is $Z$-analytic in the first case, and $Y$-analytic in the second case). 
Moreover, the natural continuous map from
$Y\times_XZ$ to the fiber product of $Y$ and $Z$ over $X$ \emph{in the category of topological spaces}
is surjective.
Indeed, let $x$ be a point of $X$, let $y$ be a pre-image of
$x$ on $Y$ and let $z$ be a pre-image of $x$ on $Z$. We want to prove that there exists
a point in $Y\times_XZ$ lying above both $y$ and $z$. For that purpose we may assume that $X, Y$ and $Z$ are
affinoid, say respectively equal to $\mathscr M(A),\mathscr M(B)$ and $\mathscr M(C)$. The completed tensor
product $\hr y\hotimes_{\hr x}\hr z$ is non-zero because it contains
$\hr y\otimes_{\hr x}\hr z$ by a result of Gruson (\cite{gruson1966}, \S 3.2 \Th 1 (4)); hence
its Banach spectrum $\mathscr M(\hr y\hotimes_{\hr x}\hr z)$ is non-empty (\cite{berkovich1990}, \Th 1.2.1).
Now if $t$ is a point of $\mathscr M(\hr y\hotimes_{\hr x}\hr z)$, its image on $\mathscr M(B\hotimes_AC)=Y\times_XZ$
lies above both $y$ and $z$, and we are done.

\begin{defi}\label{def-polyradius}\index{polyradius}
A {\em polyradius}
is a finite family of positive real numbers.
\end{defi}

\subsection{}\label{ss-intro-etar}
Let $r=(r_1,\ldots, r_n)$ be a polyradius and let $T=(T_1,\ldots, T_n)$
be a family of indeterminates. For every analytic field $k$ we denote
by $k_r$ the $k$-algebra of power series~$\sum_{I\in \Z^n} a_I T^I$ with coefficients in~$k$ such that~$\abs {a_I} r^I\to 0$ as~$\abs I\to \infty$. 
The map~$\sum a_IT^I\mapsto \max \abs {a_I} r^I$ is a multiplicative norm on $k_r$ (\cf \cite{ducros2007}, 1.2.1), 
which makes $k_r$ a $k$-affinoid algebra. Its analytic spectrum $\mathscr M(k_r)$ is the affinoid domain 
of the analytic
$n$-dimensional affine space $\A_k^{n,\mathrm {an}}$ (\cite{berkovich1993}, p.\@~25) defined by the conditions $\abs{T_i}=r_i$ for $i=1,\ldots, n$.

The norm of $k_r$ being multiplicative, it
defines a point $\eta_{k,r}$ on $\mathscr M(k_r)\subset \A_k^{n,\mathrm{an}}$, which is by its very definition\label{IN-etakr}
the unique point at which every function belonging to $k_r$ achieves its maximum; otherwise said, $\eta_{k,r}$
is the unique element of the Shilov
boundary of $\mathscr M(k_r)$ and we call it the {\em Shilov point}
of $\mathscr M(k_r)$. If there is no ambiguity with the gorund field involved, we shall often write simply 
$\eta_r$ instead of $\eta_{k,r}$. Note that by construction, the field $\hr {\eta_r}$ is nothing but the completion
of the valued field considered in Example \ref{exem-gradres-etar}; in particular, $\hrt {\eta_r}$ is isomorphic
to $\widetilde k(\tau/r)$ where $\tau=(\tau_1,\ldots, \tau_n)$ is a family of indeterminates, and $d_k(\eta_r)=n$. 

If $r$ is {\em $k$-free}, \ie, $r$ is free when viewed \index{polyradius!$k$-free}\index{k@$k$-free polyradius}
as a family of elements of the $\Q$-vector space $\R_+\gpm/\abs {k\gpm}^\Q$, then $k_r$ is an analytic field (\cf \cite{ducros2007}, 1.2.2). Hence
$\mathscr M(k_r)=\{\eta_r\}$ and $\hr {\eta_r} =k_r$ in such cases. 

Conversely if $\mathscr M(k_r)=\{\eta_r\}$ then $r$ is $k$-free. Indeed, suppose that $r$ is \emph{not}
$k$-free. Up to renumbering the $r_i$'s we can assume that there is some $j<n$ such that $(r_1,\ldots, r_j)$ is $k$-free 
and every $r_i$ with $i>j$ is torsion modulo $\abs{k\gpm}r_1^\Z\cdot \ldots \cdot r_j^\Z$. Set $L=k_{(r_1,\ldots, r_j)}$. 
Since $k_r=L_{(r_{j+1},\ldots, r_n)}$, the non-empty affinoid space $\mathscr M(k_r)$ is strictly $L$-affinoid, hence has an $L$-rigid point
$x$ 
(\ref{ss-analytic-nullst}; we can also see
this directly by choosing a finite extension $F$ of $L$
such that there exists 
$(a_i)_{j+1\leq i\leq n}$ in $F^{\{j+1,\ldots, n\}}$ satisfying the equality $\abs{a_i}=r_i$ for all $i\in \{j+1,\ldots, n\}$).
On the other hand, the equality $k_r=L_{(r_{j+1},\ldots, r_n)}$ implies that $\eta_r=\eta_{L,(r_{j+1},\ldots, r_n)}$. 
Hence $d_L(\eta_r)=n-j>0$, and $\eta_r$ is thus not $L$-rigid; in particular, $x\neq \eta_r$.

\subsection{}\label{ss-shilov-section}
Let $r$ be a polyradius. If $X$ is an analytic space
with field of definition $k$, we set $X_r=X\times_k \mathscr M (k_r)$; analogously, if $A$
is an affinoid algebra with field of definition $k$, we set $A_r=A\hotimes_k k_r$
The {\em Shilov section}
of the continuous arrow $X_r\to X$ is the map that sends a point $x$ to the Shilov point of its fiber $\mathscr M(\hr x_r)$.\index{Shilov section}
The Shilov section is continuous by 
 Lemma 3.3.2 (i) of \cite{berkovich1990}. 

\subsection{Finite morphisms (see \cite{berkovich1993}, \S 3.1)}\label{ss-finite-maps}\index{morphism!of analytic spaces!finite}\index{finite morphism}
Let $X$ be an analytic space and let $Y$ be an $X$-analytic space.
If $X$ is affinoid, say $X=\mathscr M(A)$, then $Y\to X$ is
said to be {\em finite} if $Y$ is $X$-isomorphic to $\mathscr M(B)$ for some finite Banach $A$-algebra $B$. 
In general, $Y\to X$ is said to be finite if $Y\times_X V\to V$ is finite for every affinoid domain 
$V$ of $X$, and this can be checked on a given affinoid G-covering of $X$.

Let $y$ be a point of $Y$. We shall say that $Y\to X$ is finite {\em at $y$}
if there exists an open neighborhood $V$ of $y$ in $Y$ and an open neighborhood
$U$ of the image of $y$ on $X$ such that $Y\to X$ induces a finite morphism $V\to U$. 
We shall sometimes say that $Y$ is {\em finite over $X$}, \resp {\em finite over $X$ at $y$}, instead
of saying that $Y\to X$ is finite, \resp finite at $y$. 

A morphism will be called {\em locally finite}
if it is finite at every point of the source space. \index{morphism!of analytic spaces!locally finite}\index{locally finite morphism}

\section{Coherent sheaves, Zariski topology and closed immersions}\label{s-coherent-sheaves}

\subsection{The structure sheaves}
Let $X$ be an analytic space. 
The site $X\grot$ inherits a sheaf of rings~$\mathscr O_{X\grot}$.
The restriction~$\mathscr O_{X}$ of~$\mathscr O_{X\grot}$ to the category
of open subsets of~$X$ makes~$X$ a locally ringed space. 
The sheaf of rings~$\mathscr O_{X\grot}$ is coherent, and so is~$\mathscr O_X$
when $X$ is good (for proofs\footnote{It was pointed out to the author by J\'er\^ome
Poineau that there is a gap in both proofs. Indeed, in each of them one starts with a {\em surjection}~$\mathscr O^n\to \mathscr O$ and proves that its kernel is locally finitely generated, though in order to get the coherence one should establish such a finiteness result for any, \ie,  not necessarily surjective,
map ~$\mathscr O^n\to \mathscr O$;
but the proofs do not make any use of the surjectivity assumptions.} see \cite{ducros2009}, Lemme 0.1).

\subsection{}\label{ss-tate-kiehl} If $X$ is an affinoid space, say $X=\mathscr M(A)$, the global section functor induces an equivalence between the category of coherent $\mathscr O_{X\grot}$-modules
and that of finitely generated~$A$-modules; if $M$ is a finitely generated $A$-module the corresponding coherent sheaf on $X$ assigns to any affinoid domain $V$ of $X$ the module $M\otimes_A \mathscr O_X(V)$ (this essentially follows from ``Tate acyclicity theorem" and a theorem by Kiehl; see~\cite{berkovich1993},\S 1.2). 

\subsection{} If~$X$ is good analytic space, the forgetful functor induces an equivalence between the category of coherent~$\mathscr O_{X\grot}$-modules and that of coherent~$ \mathscr O_{X}$-modules, which preserves the cohomology groups and maps locally free locally free $\mathscr O_{X\grot}$-modules to locally free $\mathscr O_{X}$-modules (\cite{berkovich1993},  \Prop 1.3.4 and  \Prop 1.3.6).

\begin{enonce}[remark]{Convention}\label{conv-ox}
Let $X$ be an analytic space. In the sequel, it will be sufficient for us to work with sheaves on~$X\grot$, and we shall not need
to pay special attention, in the good case, to the restriction of such a sheaf to the category of open subsets of~$X$ (except for stalks, see \ref{ss-stalks-fibers} below). For that reason, and for\index{coherent sheaf}
the sake of simplicity, a coherent~$\mathscr O_{X\grot}$-module will simply be called a {\em coherent sheaf on~$X$}, and we shall write~$\mathscr O_X$ instead of~$\mathscr O_{X\grot}$. 
\end{enonce}

\subsection{}\label{ss-cohsheaf-notations} If $\mathscr F$ is a coherent sheaf on an analytic space
$X$ and if $Y\to X$ is a morphism of analytic spaces we shall often denote
by $\mathscr F_Y$ the pullback of $\mathscr F$ on $Y$; in particular if\label{IN-FY}
$V$ is a analytic domain of $X$ we shall often write for short $\mathscr F_V$ instead of $\mathscr F|_V$; if~$L$ is analytic extension of $k$, we shall write
$\mathscr F_L$ instead of $\mathscr F_{X_L}$; \label{IN-FL}
if $r$ is a polyradius, we shall write for short $\mathscr F_r$ instead of $\mathscr F_{X_r}$.

\subsection{}\label{ss-fboxg}
Let $\phi \colon Z\to Y$ and $\psi \colon Z\to X$ be two morphisms of locally
noetherian schemes (\resp analytic spaces), let
$\mathscr F$
be a coherent sheaf on $Y$ and let $\mathscr G$ 
be a coherent sheaf on $X$. If $\phi$ and $\psi$ are clearly understood from the context,
we shall denote by $\mathscr F\boxtimes \mathscr G$ the tensor product $\phi^\ast\mathscr F\otimes_{\mathscr O_Z}\psi^\ast\mathscr G$. 

\subsection{Stalks}\label{ss-stalks-fibers}
Assume that $X$ is good, let $\mathscr F$ be a coherent sheaf on $X$,
and let $x$ be a point of $X$.  
We shall denote by~$\mathscr F_x$ the stalk at~$x$ of ~$\mathscr F$ {\em viewed as a sheaf on the underlying
ordinary topological space of $X$}.\label{IN-Fx}
In other words,~$$\mathscr F_x:=\lim_{\stackrel{\rightarrow}{U\;{\rm {\bf open}~neighborhood~of}\;x} }\mathscr F(U).$$
Note that this convention applies in particular
for $\mathscr F=\mathscr O_X$: in this text, $\mathscr O_{X,x}$ will always denote the stalk at $x$ of $\mathscr O_X$
viewed as a sheaf of rings
on the underlying ordinary topological space of $X$. We shall denote by $\mathfrak m_x$ its
maximal ideal, and by $\kappa(x)$ its residue field, which is a dense henselian
subfield of the valued field $\hr x$; 
see \cite{berkovich1993}, \Th 2.3.3 -- note that he uses the terminology ``quasi-complete" instead
of ``henselian". 
The tensor product $\kappa(x)\otimes_{\mathscr O_{X,x}}\mathscr F_x$
will be denoted by $\mathscr F_{\kappa(x)}$.

If $X$ is affinoid, it follows from \ref{ss-tate-kiehl}
that $\mathscr F_x=\mathscr O_{X,x}\otimes_{\mathscr O_X(X)}\mathscr F(X)$; 
we thus also have $\mathscr F_{\kappa(x)}=\kappa(x)\otimes_{\mathscr O_X(X)}\mathscr F(X)$.

\begin{rema}\label{rem-x-xv}
Be aware that 
the notation $\mathfrak m_x$ and $\kappa(x)$ might be ambiguous, because
they do not
mention the ambiant space $X$, though the object they denote actually depends on it
(contrary to $\hr x$).
Most of the time this will not cause any trouble; but we will sometimes 
write $x_V$ instead of $x$ in order to indicate that we think of $x$
as a point of a given good analytic domain $V$ of $X$: we shall
for instance write $\kappa(x_V)$ for the
residue field of $\mathscr O_{V,x}$, or $\mathscr F_{\kappa(x_V)}$ for
$\mathscr F_{V,x}\otimes_{\mathscr O_{V,x}}\kappa(x_V)$. 
Analogously, if $T$ is a good $X$-analytic space and $t$ denotes a point of $T_x$, 
then
we shall write $t_x$ whenever it is important to
make precise that $t$ is seen as a point 
of $T_x$; \eg, $\kappa(t_x)$ will denote the residue field of $\mathscr O_{T_x,t}$.

\end{rema}

\subsection{The Zariski topology}\index{Zariski topology}\label{ss-zariski-topology}
Let~$X$
be an analytic space. A {\em Zariski-closed}\index{subset!Zariski-closed}\index{Zariski-closed subset}
subset of $X$ is a subset of $X$ that is equal to the zero-locus of some coherent sheaf of ideals
$\mathscr I$ on $X$; \ie, it consists of points $x$ such that 
$f(x)=0$ for every analytic domain $V$ of $X$ containing $x$ and every function $f$ belonging to $\mathscr I(V)$. 
As suggested by the terminology, Zariski-closed subsets
of $X$ are exactly the closed subsets of a topology, the so-called {\em Zariski-topology}. 
If $E$ is any subset of $X$, we shall denote its closure
for the Zariski topology of $X$
by $\adhz EX$. 

\begin{rema}
Let $V$ be an analytic domain of $X$. 
The reader should be aware that the Zariski topology of $V$ is in general {\em strictly finer}
than the topology induced
by the Zariski topology of $X$.
Analogously, if $Y$ is an $X$-analytic
space and if $x$ is a point of $X$, the
Zariski-topology on the fiber $Y_x$ is in general
finer than the topology induced by the Zariski-topology of
$X$ (nonetheless, they coincide as soon as $x$
is rigid). 
This ``non-transitivity"
of Zariski topology
is one of the technical subtleties
that often
prevent from
transfering verbatim 
to the analytic framework what is done in scheme theory.

\end{rema}

\subsection{}\label{sss-zariski-affinoid}
Let $X$ be an affinoid space, say $X=\mathscr M(A)$. 
A subset $Y$ of $X$ is Zariski-closed if and only if there exists an ideal $I$ of $A$
such that $Y$ is equal to the set $\{x\in X, (\forall f\in I, f(x)=0)\}$; in other words, the Zariski topology on $X$ is nothing but the pre-image 
of the Zariski topology of $\spec A$ under the natural map $X\to \spec A$.

Let $Y$ be a Zariski-closed subset of $X$
and let $I$ be as above.
The morphism $\mathscr M(A/I)\to \mathscr M(A)$ induced by the quotient map
$A\to A/I$ establishes a homeomorphism $\mathscr M(A/I)\simeq Y$, which 
endows $Y$ with the structure of an affinoid space.

\subsection{}Let $X$ be an analytic space. Let $Y$ be a Zariski-closed subset of $X$, and let 
$\mathscr I$ be a coherent sheaf of ideals whose zero locus is $Y$.
By performing G-locally the construction described in \ref{sss-zariski-affinoid}
above, one gets a structure of an $X$-analytic space on $Y$, and the structure morphism $\iota \colon Y\to X$
satisfies the following property. 

\medskip
\begin{enumerate}
\item The map underlying $\iota$ is the inclusion, and $\iota_\ast\mathscr O_Y=\mathscr O_X/\mathscr I$. 

\item If $Z$ is an analytic space and if $f\colon Z\to X$ is a morphism such that $\mathscr I$
is contained in $\mathrm{Ker}\;(\mathscr O_X\to f_\ast\mathscr O_Z)$ then $f$ factors uniquely through $\iota$. 
\end{enumerate}
The Zariski-closed subset $Y$ together with this analytic structure is called the {\em closed analytic subspace of $X$ defined by $\mathscr I$}. 
Its Zariski topology only depends on the set $Y$, and not on $\mathscr I$:  this is simply the restriction to $Y$ of the Zariski-topology of $X$, 
and we shall call it the Zariski topology of $Y$. \index{closed analytic subspace}

\begin{rema}\label{rem-gerr-grau}
It follows immediately from Gerritzen-Grauert theorem (more precisely from Temkin's version of it for Berkovich spaces, see
\cite{temkin2005}, \Th 3.1) that
if~$Y$ is a closed analytic subspace of~$X$, and if $V$
is any analytic domain of $Y$, then~$V$ can be G-covered by affinoid domains of the
form~$U\cap Y$ with~$U$ an affinoid domain of~$X$.
\end{rema}

\subsection{} Let $k$ be an analytic field, let $X$ be a $k$-analytic
space and let $Y$ be a Zariski-closed subset of $X$. We shall denote by $Y_L$ the preimage of $Y$ on $X_L$; 
this is a Zariski-closed subset of $X_L$.

\subsection{}A morphism $Z\to X$ of analytic spaces is called
a {\em closed immersion}
if it induces an isomorphism between $Z$ and a closed analytic subspace of $X$. 

\subsection{}\label{ss-separated}
Let $X$ be an analytic space and let $Z$ be an $X$-analytic space. 
We say that the morphism $Z\to X$ 
\index{morphism!of analytic spaces!separated}\index{separated morphism}\index{morphism!of analytic spaces!locally separated}\index{locally separated morphism}
is {\em separated} (or that $Z$ is separated over $X$)
if the diagonal morphism $Z\to Z\times_X Z$ is a closed immersion; it is called
\emph{locally separated}
if every point of $Z$ has an open neighborhood $U$ such that $U\to X$
is separated. 

If $k$ is an analytic field, a $k$-analytic space $X$ is called separated
(\resp locally separated)
if the structure map $X\to \mathscr M(k)$ is\index{analytic space!separated}\index{separated analytic space}
separated (\resp locally separated).\index{analytic space!locally separated}\index{locally separated analytic space}

Every embedding of an analytic domain is separated. Every good $k$-analytic space is locally separated. Every morphism 
between locally separated $k$-analytic spaces 
is locally separated. 

A separated analytic space is Hausdorff
(\cite{berkovich1993}, \Prop 1.4.2)
but the converse is not true in 
general: for instance, the space $Y$ described in Example \ref{ex-temred}
below is Hausdorff (and even compact) and not separated.

\subsection{Reduced analytic spaces}\index{analytic space!reduced}\index{reduced analytic space}
\label{ss-defred-adhoc}
An analytic space $X$ is said to be \emph{reduced}
if $\mathscr O_X(V)$ is reduced for every analytic domain
$V$ of $X$; it suffices to check it on every \emph{affinoid}
domain of $X$.

\subsection{}\label{ss-red-gloc}
If $A$ is an affinoid algebra, then $\mathscr M(A)$ is reduced if
and only if $A$ is reduced (\cite{bosch-g-r1984}, 7.3.2 \Cor 10 in the strict case; one 
can easily reduce to this by
using Lemma 1.3 of \cite{ducros2007}).

Let us also mention that in general, an element $f$ of $A$ is
nilpotent if and only if $f(x)=0$ for all $x\in \mathscr M(A)$. 

\subsection{G-local nature of Zariski topology and the notion of reduced structure}\label{ss-red-structure}
Let $X$ be an anaytic space and let $Y$ be subset of $X$ which is G-locally
a Zariski-closed subset; \ie, there exists a G-covering $(X_i)$ of $X$
by analytic domains such that $Y\cap X_i$.
Let $\mathscr I$ be the sheaf of ideals defined by the assignment 
$$U\mapsto \{f\in \mathscr O_X(U), f(y)=0 \;\forall \; y\in Y\cap U\}.$$
The sheaf $\mathscr I$ is coherent (\cite{ducros2009}, \Prop 4.2 (i); the proof
rests in the crucial way on the result recalled in
\ref{ss-red-gloc}), and its zero-locus is equal to $Y$; hence
$Y$ is 
Zariski-closed. As a consequence, \emph{being Zariski-closed is of G-local nature}. 

It follows from the definition that $\mathscr I$ is the greatest coherent sheaf of ideals
on $X$ with zero-locus $Y$; if $\mathscr J$ another such sheaf, then $\mathscr I$
is the radical of $\mathscr J$; \ie, the sheaf of functions that are G-locally
nilpotent modulo $\mathscr J$ (if $Y=X$ we can take $\mathscr J=0$, 
hence $\mathscr I$ is the sheaf of G-locally nilpotent functions). 

The closed analytic subspace defined by $\mathscr I$ is denoted by $Y_{\mathrm{red}}$.
It is reduced
and it is the final object of the category of 
\emph{reduced}
 analytic spaces $Z$ equipped with a morphism $\phi \colon Z\to X$ such that $\phi(Z)\subset Y$. 
 
 Instead of writing $Y_{\mathrm{red}}$, we shall sometimes simply write $Y$
 after having
made precise
 that we endow it with its \emph{reduced} analytic structure. \index{reduced analytic structure}

\subsection{Boundary}\label{ss-boundary}\index{boundary}
Let $X$ be an analytic space and let $Y$ be an $X$-analytic space. There is a well-defined notion of {\em boundary}
of the morphism $Y\to X$, which is also called the {\em relative boundary of $Y$ over $X$} (\cite{berkovich1990}, \Def 2.5.7
for the affinoid case, and \cite{berkovich1993}, \Def 1.5.4 in general).
This is a closed subset of $Y$ which is denoted by $\partial (\phi)$ or $\partial (Y/X)$; its complement is called the {\em interior}
of $\phi$, or the {\em relative interior of $Y$ over $X$}, and is denoted by $\mathrm{Int}(\phi)$ \index{interior}
or $\mathrm{Int}(Y/X)$. The morphism $\phi$ is said to be {\em inner}
or {\em boundaryless} at some point $y$ of $Y$ if $y\in \mathrm{Int}(Y/X)$;
it is said to be inner or\index{morphism!of analytic spaces!boundaryless}\index{analytic space!boundaryless}\index{boundaryless analytic space}
boundaryless if it is so at every point of $Y$; \ie, if $\partial (Y/X)=\emptyset$ (Berkovich calls such a map a {\em closed}\index{morphism!of analytic spaces!inner}\index{inner morphism}
morphism, but we shall not use this terminology here: we shall reserve \emph{closed}
for denoting morphisms that are \emph{topologically}
closed). 

We shall write $\partial X$ and $\mathrm{Int}(X)$ instead of $\partial (X/\mathscr M(k))$ and
$\mathrm{Int}(X/\mathscr M(k))$, for $k$ the field of
definition of $X$; we shall call them respectively the boundary and the interior of $X$. 
The space $X$ will be called boundaryless if $\partial X=\emptyset$; any boundaryless space is good.

\subsection{}\label{ss-boundary-basics}
Let us list here some useful, basic properties of the boundary that will be useful. 

\begin{enumerate}[1]
\item If $Y\to X$ is finite, it is boundaryless. Conversely, if $Y\to X$ is boundaryless
\emph{and if both $Y$ and $X$ are affinoid}, then $Y\to X$ is finite (\cite{berkovich1990}, \Cor 2.5.13). 

\item If $Y$ is an analytic domain of $X$ then $\partial (Y/X)$ is nothing but the {\em topological}
boundary of $Y$ inside $X$ (\cite{berkovich1993}, \Prop 1.5.5 (i)).

\item Assume that $Y\to X$ is locally separated (\ref{ss-separated}; note that this assumption holds if both $Y$ and $X$ are good, or if $Y$ is an analytic
domain of $X$), let $Z$
be an $X$-analytic space
and let $\sigma \colon Z\to Y$ be an $X$-morphism. We then have the equality
$$\mathrm{Int}(Z/X)=\mathrm{Int}(\sigma)\cap (\sigma \inv (\mathrm{Int}(Y/X))$$
(\cite{temkin2004}, \Cor 5.7).  
In particular,
$\sigma(Z)\subset \mathrm{Int}(Y/X)$ as soon as $\sigma$ is 
boundaryless; \eg, $\sigma$ is finite, see (1).

\item The map $Y\to X$ is no longer assumed
to be locally separated. If $X'$ is an arbitrary analytic space (not necessarily $k$-analytic) and if $X'\to X$ is a morphism, 
then $\mathrm{Int}(Y/X)\times_X X'\subset \mathrm{Int}(Y\times_X X'\to X')$ (this rests on \Prop 3.1.3 of \cite{berkovich1990}). 

\item The property for a map of being boundaryless is G-local on the target (\cite{temkin2004}, \Cor 5.6; be aware that Temkin uses the word ``closed" instead of ``boundaryless").

\end{enumerate}

\subsection{Proper morphisms}\label{ss-proper-map}\index{morphism!of analytic spaces!proper}\index{proper morphism}
Let $k$ be an analytic field and let $\phi \colon Y\to X$ 
be a morphism of $k$-analytic spaces. We shall say that $\phi$
is \emph{proper} (or that $Y$ is proper over $X$, or proper $X$-analytic space) if it satisfies the two following
conditions: 

\begin{enumerate}
\item $\phi$ is boundaryless; 

\item $\phi$ is proper and separated
\emph{as a continuous map between topological spaces}, see \ref{ss-remind-topo}. 

\end{enumerate}

Properness can be checked G-locally on the target; any proper morphism is closed; any finite morphism is
proper.

\subsection{Kiehl's Theorem}\label{ss-kiehl-proper}
Let $\phi \colon Y\to X$ be a \emph{proper}
morphism of analytic spaces (\ref{ss-proper-map}). 
If $\mathscr F$ is any coherent sheaf on $Y$, then $\mathrm R^q\phi_\ast\mathscr F$
is a coherent sheaf on $X$ for every $q$; in particular, $\phi_\ast\mathscr F$
is coherent.
This is essentially due to Kiehl who proved it
in the rigid-analytic setting, \cite{kiehl1967b}; for the details about the transfer of this result
in Berkovich's theory, see \cite{ducros2015}, section 2. 

Let $Z$ be
a Zariski-closed
subset of $Y$ and let 
$\mathscr I$ be the sheaf of ideals on $Y$ that defines
the reduced structure 
on $Z$. 
By definition, a section of $\mathscr O_Y$ belongs to $\mathscr I$ if and only if
it vanishes pointwise on $Z$. Hence a section of $\mathscr O_X$ belongs to the kernel
$\mathscr J$
of $\mathscr O_X\to \pi_\ast (\mathscr O_Y/Z)$ if and only if it vanishes pointwise on $\phi(Z)$. 

By (the Berkovich avatar of) Kiehl's theorem, $\phi_\ast(\mathscr O_Y/\mathscr I)$
is a coherent sheaf on $X$, hence $\mathscr J$ is a coherent sheaf of ideals on $X$. By the above, $\phi(Z)$ is
contained in the zero-locus of $\mathscr J$.

On the other hand, since $\phi$ proper, it is in
particular  topologically proper and $\phi(Z)$ is thus a closed subset of $X$. Since $\mathscr J$ is the sheaf of functions
vanishing pointwise on $\phi(Z)$, 
we have the equality $\mathscr J_{X\setminus \phi(Z)}=\mathscr O_{X\setminus \phi(Z)}$. Hence the zero-locus
of $\mathscr J$ is contained in $\phi(Z)$, and thus equal to $\phi(Z)$. As a consequence,
\emph{$\phi(Z)$ is a Zariski-closed
subset of $X$}.

\section{Dimension theory}

\index{dimension}
We recall here some basic facts about the dimension theory of
analytic spaces. References are Berkovich's foundational work on the topic \cite{berkovich1990}, \S2 
and the author's paper \cite{ducros2007}. 

\subsection{Analytic dimension of an affinoid algebra}\label{ss-dim-affalg}
Let $k$ be an analytic field and let $A$ be a $k$-affinoid algebra. Let $L$ be an
analytic extension
of $k$ such that $A_L$ is strictly $L$-affinoid. The Krull dimension of $A_L$ 
does not depend on the choice of $L$, 
and is called the {\em $k$-analytic dimension}\index{dimension!k@$k$-analytic (of a $k$-affinoid algebra)}
of $A$; we shall denote it by $\dim_k A$ (it is finite unless $A=0$, in which case we have $\dim_k A=-\infty$).\label{IN-dimkA}

\subsection{}The analytic dimension may actually depend on $k$, and not only on the ring $A$. For instance, let $r=(r_1,\ldots, r_n)$ be a $k$-free
polyradius and let $L$ be any analytic extension of $k$ such that
$r_i\in \abs {L\gpm}$ for every $i$ (e.g., $L=k_r$). The $L$-algebra $L\hotimes_k k_r=L_r$ is then strictly $L$-affinoid, 
and of Krull dimension $n$; one has therefore $\dim_k k_r=n$. On the other hand, $k_r$ is an analytic field and $\dim_{k_r}k_r=0$. 

\subsection{}\label{ss-dim-dimkrull}
We have $\dim_{\mathrm{Krull}}A\leq \dim_k A$. Indeed, choose a 
$k$-free polyradius $r$ such that $\abs{k_r\gpm}\neq \{1\}$ and such that $A_r$ is strictly $k_r$-affinoid. 
If $B$ is any $k$-affinoid algebra that is a domain, then $B_r$ is a domain too (\cite{ducros2007}, Lemme 1.3). 
It follows that the pre-image of an irreducible Zariski-closed subset of $\spec A$ under the
natural map $\spec A_r\to \spec A$
is still irreducible. On the other hand, $\spec A_r\to \spec A$
is surjective (it is even faithfully flat, \cite{berkovich1993}, Lemma 2.1.2), which
implies that
two distinct subsets of $\spec A$ have distinct pre-images on $\spec A_r$. 
As a consequence, $\dim_{\mathrm{Krull}}A\leq \dim _{\mathrm{Krull}}A_r=\dim_k A$. 

\subsubsection{}\label{ss-dimanal-zero}
The $k$-analytic dimension of $A$ is zero if and only if $A$
is a non-zero finite $k$-algebra (\cite{ducros2007}, Lemme 1.7).  

\subsection{Dimension of an affinoid space}\label{sss-dim-affspace}
If $X$ is a $k$-affinoid space, its $k$-analytic dimension $\dim_k X$ is by definition the $k$-analytic dimension of \label{IN-dimkX}
the corresponding $k$-affinoid algebra
(hence when $X$ is strictly $k$-affinoid, it coincides with the Krull dimension of $X$ equipped with its Zariski topology). 
One has $\dim_k V\leq \dim_k X$ for every affinoid domain $V$ of $X$. It follows that $\dim_k X$ is equal to the supremum of
$\dim_k V$ for $V$ going through the set of all affinoid domains of $X$. 

\subsection{Dimension of an arbitrary space}\label{ss-dim-arbitrary-space}\index{dimension!k@$k$-analytic (of a $k$-analytic space)}
If $X$ is an arbitrary $k$-analytic
space, we {\em define}
$\dim_k X$ as the supremum of
$\dim_k V$ for $V$ going through the set of all affinoid domains of $X$ (this is compatible with the definition given
in \ref{sss-dim-affspace} in the affinoid case). If $V$ is any
analytic domain of such an $X$ then $\dim_k V\leq \dim_k X$. 

\subsection{}\label{ss-dimx-dkx} Let $X$ be a $k$-analytic space. If $x\in X$, we set $d_k(x)=d_k(\hr x)$ (see
~\ref{sss-analytic-extension}). This invariant plays a role similar to that
of the transcendence degree in classical dimension theory; indeed, one has
$$\dim_k X=\sup_{x\in X} d_k(x).$$

\subsection{}\label{ss-x-dimzero}
Let $X$ be a non-empty $k$-analytic space. If $X$ consists only of rigid points it follows
from \ref{ss-dimx-dkx}
above that $\dim_k X=0$. Conversely if $\dim_k X=0$ it follows from
\ref{ss-dimanal-zero} (applied to every non-empty affinoid domain of $X$)
that $X$ only consists of rigid points and is topologically discrete.

\subsection{Behavior of $d_k$ with respect to the Shilov section}\label{ss-dk-shilov}
Let $X$ be a $k$-analytic space and let $r=(r_1,\ldots, r_n)$ be a polyradius.
Let $\mathfrak s: X\to X_r$ denote the Shilov section (see \ref{ss-shilov-section}) and let 
$x$ be a point of $X$. It follows from \ref{ss-intro-etar} that $d_{\hr x}(\mathfrak s(x))=n$; we thus
have $d_k(\mathfrak s(x))=d_k(x)+n$.

Let us now assume that $r$ is $k$-free. We can then write
$$d_k(\mathfrak s(x))=d_{k_r}(\mathfrak s(x))+d_k(k_r)=d_{k_r}(\mathfrak s(x))+n$$
(because $d_k(k_r)=n$ by \ref{ss-intro-etar}), whence the equality
$$d_{k_r}(\mathfrak s(x))=d_k(x).$$

\subsection{Local dimension}\index{dimension!local}\label{ss-def-dimloc}
There is also a notion of {\em local}
$k$-analytic dimension: if $X$ is a $k$-analytic space, then for every $x\in X$ one defines
the $k$-analytic dimension of $X$ {\em at $x$} as the infimum of $\dim_k V$ for $V$ running through the set 
of analytic domains of $X$ containing $x$; we denote it by $\dim_{k,x} X$.  If $V$ is any analytic domain of $X$
containing $x$ then $\dim_{k,x}V=\dim_{k,x} X$. For $n\in \N$, we will say that $X$ is {\em of pure\label{IN-dimkxX}
dimension $n$}\index{dimension!pure (for an analytic space)}\index{pure dimension (for an analytic space)}
if $\dim_{k,x}X=n$ for every $x\in X$; 
if $X$ is of pure dimension $n$, so is any analytic domain of $X$. 

We shall say that $X$ is \emph{equidimensional}\index{analytic space!equidimensional}\index{equidimensional analytic space}
if it is of pure dimension $n$ for some $n$. Such an $n$
is necessarily equal to $\dim_k X$ if $X\neq \emptyset$;
but $\emptyset$ is equidimensional, of pure dimension $n$ for all $n$, 
and of dimension $-\infty$.

\subsection{Abhyankar points}\label{ss-abhyankar-points}\index{point!Abhyankar}\index{Abhyankar point}
Let $X$ be a $k$-analytic space and let $x$ be a point of $X$. In view of \ref{ss-dimx-dkx}
we have $d_k(x)\leq \dim_{k,x} X$. It can be seen as an analytic avatar of the classical
Abhyankar inequality stated in \ref{ss-index-graded} (2). For that reason, we shall say that
$x$ is an \emph{Abhyankar point}
if $d_k(x)=\dim_{k,x} X$.

\subsection{}\label{sss-dim-zarclos} Let $X$ be a $k$-analytic space and let $Y$
be a Zariski-closed subset of $X$. Let $\mathscr I$ be any coherent sheaf
of ideals on $X$ whose zero locus is equal to $Y$. The $k$-analytic dimension 
of the closed analytic subspace of $X$ defined by $\mathscr I$ only depends on $Y$, 
and not of the chosen ideal sheaf $\mathscr I$; the same holds for its $k$-analytic dimension 
at a given point $x$ of $Y$. 
We will then use the expressions ``$k$-analytic dimension of $Y$" and ``$k$-analytic dimension of $Y$ at $x$", 
and the notation $\dim_k Y$ and $\dim_{k,x} Y$, without fixing any $k$-analytic
structure on $Y$. One has $\dim_k Y\leq \dim_k X$ and $\dim_{k,x}Y\leq \dim_{k,x} X$.

Let $(X_i)$ be any {\em set-theoretic}
covering of $X$ by Zariski-closed subsets of analytic domains. One has 

$$\dim X=\sup_{x\in X} d_k(x)=\sup_i \sup_{x\in X_i} d_k(x)=\sup_i \dim X_i.$$

\subsection{} \label{dim-goundfield-extension}Analytic dimension behaves well under ground field extension: if $X$ is
a $k$-analytic space and $L$ is an analytic extension 
of $k$, then $\dim_L X_L=\dim_k X$, and $\dim_{L,y}X_L=\dim_{k,x} X$ for every $x\in X$
and every pre-image $y$ of $x$ on $X_L$. 

\subsection{}\label{ss-dim-nofield}\index{dimension}\index{local dimension}
When the ground field is clearly understood from the context, we shall omit it 
in the notation. For instance, if $X$ is an analytic space and if $x$\label{IN-dimX}
is a point of $X$, then $\dim X$\label{IN-dimxX}
and $\dim_x X$ will
denote analytic dimensions over the analytic field that is\index{dimension!local}
implicitly part of the definition of $X$.  If $Y$ is an $X$-analytic space and if $y$ is a point of $Y$, 
then $\dim Y_x$ and $\dim_y Y_x$ will
denote $\hr x$-analytic dimensions. The integer $\dim_y Y_x$ will be called
the {\em relative dimension of $Y$ over $X$ at $y$}; if the morphism $Y\to X$\index{dimension!relative}
has been given a name, say $\phi$, we shall also write $\dim_y\phi$ instead \label{IN-dimyphi}
of $\dim_y Y_x$, and call it the \emph{dimension of $\phi$ at $y$}\index{dimension!of a morphism}. 
The $\N$-valued function $y\mapsto \dim_y \phi$ is upper semi-continuous
for the Zariski topology of $T$ (\cite{ducros2007}, \Th 4.9). If $n$ is a non-negative
integer, we shall say that $Y$ is \emph{of pure 
relative dimension $n$ over $X$}, or that $\phi$ is \emph{of pure dimension $n$}, if $\dim_y \phi=n$\index{dimension!pure (for a morphism)}\index{pure dimension (for a morphism)}
for every $y\in Y$ or, which amounts to the same, if all fibers of $\phi$ are of pure dimension $n$. 
We shall say that $Y$ is \emph{relatively equidimensional over $X$}, or that \emph{$\phi$
is equidimensional}, if $\phi$ is of pure dimension $n$ for some $n\in \N$; the integer $n$ is then uniquely determined
as soon as $Y\neq \emptyset$: this is the common dimension of all non-empty fibers of $\phi$.

\subsection{}\label{ss-dim-sobafi}
The formula stated in \ref{ss-dimx-dkx}
enables to relate, to some extent, the dimensions
of the source, of the target, and of the fibers
of a given map. 
For instance, let 
let $\phi \colon Y\to X$ be a morphism of $k$-analytic spaces. 
For every $x\in X$, the equality
$\sup_{y\in \phi^{-1}(x)}d_{\hr x}(y)=\dim \phi\inv(x)$
implies that

\[
\dim Y\geq \sup_{y\in \phi^{-1}(x)}d_k(y)=\sup_{y\in \phi^{-1}(x)} [d_{\hr x}(y)+d_k(x)]=\dim \phi\inv(x)+d_k(x).\]
This has the following consequences:

\begin{enumerate}[1]

\item For every $x\in X$, one has $\dim \phi\inv(x)\leq \dim Y$. 

\item If $d$ is a non-negative integer such that $\dim \phi\inv(x)\leq d$ for every $x\in X$, then
\[
\dim Y=\sup_{y\in Y}d_k(y)\leq d+\sup_{x\in \phi(Y)} d_k(x)\leq d+\dim X.\]

\item If $d$ is a non-negative integer such that $\dim \phi\inv(x)=d$ for every $x\in \phi(Y)$, and 
$\phi(Y)$ is a Zariski-closed subset of an analytic domain of $X$, then 
\[\dim Y= d+\sup_{x\in \phi(Y)} d_k(x)=d+\dim \phi(Y)\]
(the assumption that $\phi(Y)$ is  Zariski-closed in an analytic domain of $X$
simply ensures that $\dim \phi(Y)$ \emph{makes sense}; of course, the above formula
remains valid without any assumption on $\phi(Y)$ if we \emph{define}
$\dim \phi(Y)$ as $\sup_{x\in \phi(Y)}d_k(x)$). 

\item If $y$ is a point of $Y$ such that $d_k(y)=\dim Y$
(such a point exists if and only if $Y$ is finite-dimensional), 
and if $x$ denotes
its image in $X$, then 
\[d_k(x)=\dim Y-\dim \phi\inv(x).\]
\end{enumerate}

\subsection{Quasi-finite morphisms}\label{ss-quasi-finite}\index{morphism!of analytic spaces!quasi-finite}\index{quasi-finite morphism}
Let $\phi\colon Y\to X$ be a morphism of $k$-analytic spaces. 
If $y$ is a point of $y$, we shall say that $\phi$
is \emph{quasi-finite}
at $y$ (or that $Y$ is quasi-finite at $y$ over $X$) if $\dim_y \phi=0$. 
The morphism $\phi$
is then finite at $y$ if and only if it is quasi-finite
and boundaryless at $y$ (\cite{berkovich1993}, \Cor 3.1.10). 

We shall say that $\phi$ is
\emph{quasi-finite}
(or that $Y$ is quasi-finite over $X$)
if $\phi$ is topologically proper and quasi-finite at every point 
of $y$; \ie, $\phi$ is topologically proper and of pure relative
dimension zero. 

A quasi-finite map is finite if and only if it is proper (which means that it is boundaryless and
topologically separated, since it is already topologically proper by definition).

\begin{rema}\label{rem-depart-qf}
The reader should be aware that our definition
of quasi-finiteness
 \emph{differs from Berkovich's}.
Indeed, Berkovich uses the expression ``quasi-finite at the point $y$" 
for ``finite at the point $y$" in our sense (\ref{ss-finite-maps}).
We have chosen to depart from Berkovich's definition
because we want
a quasi-\'etale map
(see Definition \ref{def-qet}) to be quasi-finite
at every point of the source space, and also for the sake of analogy with
scheme theory (see for instance Lemma \ref{lem-local-cm},
Theorem
\ref{thm-local-cm}, Theorem \ref{thm-multisections-local} or Theorem \ref{thm-multisections-global}).
 \end{rema}

\section{Irreducible components}\label{irr-comp}

The Zariski topology of an analytic space is far from being noetherian in general,
but there is nevertheless 
a reasonable theory of irreducible components in this settting. Brian Conrad developped
in it the rigid analytic framework in \cite{conrad-1999}. The author suggested another approach
in \cite{ducros2009}, \S4, which works for arbitrary analytic spaces and is perhaps more direct (it does not make any use of the normalization, contrary to \cite{conrad-1999}); 
we are now going to describe it. 

\subsection{Irreducible analytic spaces}\label{ss-irred-components}\index{analytic space!irreducible}\index{irreducible analytic space}
Let $X$ be an analytic space. A Zariski-closed subset of $X$ will be called
{\em irreducible}
if it
is irreducible for the Zariski topology of $X$, and $X$ will be called
{\em integral}\index{analytic space!integral}\index{integral analytic space}
if it is both irreducible and reduced. If $Y$ is an irreducible Zariski-closed subset
of $X$, it is purely
$d$-dimensional for some $d$
and if $Z$ is any Zariski-closed subset of $Y$
with $Z\neq Y$
then $\dim Z<d$.

\begin{defi}
Let $X$ be an analytic space. There exists a set $E$ of irreducible Zariski-closed subsets of $X$ having the following properties: 

\begin{itemize}
\item The set $E$ is G-locally finite; \ie, any affinoid domain of $X$ intersects only finitely many elements of $E$.

\item One has $X=\bigcup_{Z\in E}Z$. 

\item If $Y$ and $Z$ are two elements of $E$ with $Y\subset Z$ then $Y=Z$. 

\end{itemize}

The set $E$ is uniquely determined by those properties. Its elements are
exactly the maximal irreducible Zariski-closed subsets of $X$; moreover, 
every irreducible Zariski-closed subset of $X$ is contained in one element 
of $E$. The elements of $E$ are called the {\em irreducible components}\index{irreducible components}
of $X$. 
\end{defi}

\subsection{}If $X=\mathscr M(A)$ is an affinoid
space, an irreducible component of $X$ is nothing
but the pre-image of an irreducible component of $\spec A$. 
In general, i.e., if $X$ is no longer assumed
to be affinoid, for every Zariski-closed subset $Y$ of $X$ the following are equivalent: 

\begin{itemize}
\item[(i)] $Y$ is an irreducible component of $X$. 

\item[(ii)] There exist an affinoid domain $V$ of $X$ and an irreducible
component $Z$ of $V$ such that $Y=\overline Z^{X\zar}$. 
\end{itemize}

\subsection{}
Let $X$ be an analytic space and let $Y$
be a Zariski-closed subset of $X$; let $\mathscr I$ be a coherent sheaf of ideals with zero locus $Y$. 
The irreducible components of the closed analytic subspace of $X$ defined by $\mathscr I$ can be characterized purely in terms
of the Zariski topology of $Y$; therefore they only depend on the set $Y$, and not on $\mathscr I$; we shall call
them the irreducible
components of $Y$. 

For instance, if $E$ is a subset of the set of irreducible
components of $X$, then then $\bigcup_{Z\in E}Z$ is a Zariski-closed subset of $X$, whose irreducible
components are precisely the elements of $E$. 

\subsection{}\label{ss-dim-andom}
Let $X$ be an analytic space, let $V$
be an analytic domain of $X$, and let $d$ be a non-negative integer. 

\begin{itemize}[label=$\bullet$]

\item If $Y$ is an irreducible component of $X$ of dimension $d$, then $Y\cap V$ is a union (possibly
empty, possibly infinite) 
of irreducible components of $V$, each of which has dimension $d$. 

\item If $Z$ is an irreducible component of $V$ of dimension $d$, then $\overline Z^{X\zar}$ is an irreducible
component of $X$, of dimension $d$. 

\end{itemize}

\subsection{}\label{ss-dim-scalext}
Let $k$ be an analytic field, let $L$ be an analytic extension of $k$ and let $Y$ be an irreducible component of
a $k$-analytic space $X$. The Zariski-closed subset $Y_L$ of
$X_L$ has finitely many 
irreducible components. For every such component $Z$ we have $\dim Z=\dim Y$, the natural map $Z\to Y$ is surjective, and $Z$ is equal to $T\times_FL$
for some finite separable extention $F$ of $k$ inside $L$
and some irreducible component $T$ of $Y_F$; moreover $Z$
is an irreducible component of $X_L$. 

Conversely, if $Z$ is an irreducible component of $X_L$ then its image $Y$ in $X$ is an irreducible component of $X$ and $Z$ is an irreducible
component of $Y_L$. 

\subsection{}\label{ss-dimloc-open}
If $X$ is an analytic space, then $\dim X=\sup_Z \dim Z$ for $Z$ running through the set of irreducible components of $X$. For
every $x\in X$
the local dimension $\dim_x X$ is equal to $\max_Z\dim Z$ for $Z$ running through the (finite) set of irreducible components of $X$ containing $x$. 

\begin{rema}\label{rem-dimloc-open}
Let $x$ be a point of an analytic space $X$ and let $E$ (\resp $F$) be the set of irreducible components of $X$ that contain (\resp avoid) $x$; let $Y$ (\resp $Z$)
be the union of all components belonging to $E$ (\resp $F$). Both $Y$ and $Z$ are Zariski-closed subsets of $X$. The dimension of $Y$ being the supremum
of the dimensions of all components belonging to $E$, it is equal to $\dim_x X$. Now $U:=Y\setminus Z$ is a Zariski
open subset of $X$ that contains $x$, that is contained in $Z$, and that intersects
every irreducible component of $Z$. It follows that $\dim U=\dim Y=\dim_x X$. 
\end{rema}

\begin{rema}\label{rem-abh-point}
Let $X$ be an analytic space and met $x$
be a point of $X$. Assume that $x$
lies on a Zariski-closed
subset $Y$ of $X$ and is Abhyankar \emph{in $Y$} 
(\ref{ss-abhyankar-points}). Since $\dim \adhz {\{x\}}X\geq d_k(x)$ and since $d_k(x)=\dim_x Y$ is the maximum
of the dimensions of the irreducible components of $Y$ containing $x$, we see that 
$\adhz {\{x\}}X$ is an irreducible component of $Y$ of dimension $\dim_x Y=d_k(x)$; note that as $x$ is Zariski-dense in
$\adhz {\{x\}}X$, the latter is
even the \emph{only}
irreducible component of $Y$ containing $x$. 
\end{rema}

\subsection{}\label{ss-dimloc-morfin}
Let $\phi \colon Y\to X$
be a finite morphism
of analytic spaces.  
Let $(Y_i)$ be the family of irreducible components of $Y$. For every $i$, the image $\phi(Y_i)$ is an irreducible Zariski-closed subset of $X$, and 
$\dim \phi(Y_i)=\dim Y_i$ by \ref{ss-dim-sobafi}. Since $\phi$ is
topologically proper, the family $(\phi(Y_i))$ is G-locally finite; it follows that the irreducible
components of $\phi(Y)$ are precisely the maximal elements among the $\phi(Y_i)$'s. 

Now let $x$ be a point of $X$ having exactly one pre-image $y$ on $Y$; let $J$ be the set of indices $i$ such that $y\in Y_i$.
By the above, the irreducible components of $\phi(Y)$ that contain $x$ are exactly the maximal elements among the $\phi(Y_i)$'s for $i$
running through $J$. Since $\dim \phi(Y_i)=\dim Y_i$ for every $i$, this implies that
$$\dim_x \phi(Y)=\max_{i\in J}\dim Y_i=\dim_y Y.$$

\begin{lemm}\label{lem-abhyankar-dimension}
Let~$n$ and~$m$ be two integers,
and let~$Y\to X$ be a morphism of~$k$-analytic spaces,
with~$Y$ of pure dimension $m$
and $X$ of dimension $n$.
Let~$x$ be a point of $X$
such that~$d_k(x)=n$.
The fiber~$Y_x$ is then purely of dimension~$m-n$.
\end{lemm}

\begin{proof}
We may assume that~$Y$ and~$X$ are affinoid. 
Let~$T$ be an irreducible component of~$Y_x$, let~$d$
be its dimension, and let~$y$ be an Abhyankar point
of $T$. This condition 
implies that $T$ is the only
irreducible component of $Y_x$ that contains $y$ (Remark \ref{rem-abh-point}); 
as a consequence, $\dim_y Y_x=d$. 
One has~$d_k(y)=d_{\hr x}(y)+d=n+d$,
whence the inequality~$n+d\leq m$. It suffices to prove
the reverse inequality. 

Since $\dim_y Y_x=d$, 
the map~$Y\to X$ admits by \Th 4.6 of \cite{ducros2007}
a factorization~$Y\to \A^d_X\to X$
whose first step
is quasi-finite at~$y$. 
By \Th  4.9 (or, more simply, 
\Th 3.2) of \cite{ducros2007}, there exists an affinoid
neighborhood~$V$ of~$y$ in $Y$ such that~$V\to \A^d_X$ 
is quasi-finite. 

Since~$Y$ is purely~$m$-dimensional,
its non-empty affinoid domain~$V$ is~$m$-dimensional. Therefore, 
there exists~$v\in V$ with~$d_k(v)=m$. Let~$z$ be the image of~$v$ on~$\A^d_X$,
and let~$t$ be the image of~$z$ on~$X$.
The morphism~$V\to \A^d_X$ being quasi-finite, 
one has
\[m=d_k(v)=d_{\hr z}(v)+d_k(z)=d_k(z)=d_{\hr t}(z)+d_k(t)\leq d+n.\]
\end{proof}

\begin{lemm}\label{lem-closure-zartop}
Let $X$ be an analytic space, let~$U$ be a {\em Zariski}-open subset of~$X$ and let
$x$ be a point of $X$ ; set~$F=X\setminus U$. The following are equivalent:

\begin{enumerate}[i]

\item The point $x$ belongs to $\adht U X$.

\item The point $x$ belongs to $\adhz U X$. 

\item The Zariski-open subset $U$ intersects at least one of the irreducible components of~$X$ that contain~$x$;

\item There exists an irreducible component~$Z$ of~$X$ that contains~$x$ and satisfies~$\dim (F\cap Z)<\dim Z$.
\end{enumerate}
\end{lemm}

\begin{proof}
It is clear that (i)$\Rightarrow$(ii). Assume that~(ii) is true. Let~$X'$ be the union of all irreducible components of~$X$ that do not contain~$x$. It is a Zariski-closed 
subset of~$X$. Its complement $X\setminus X'$ is then a Zariski-open neighborhood of $x$, hence it instersects $U$ by assumption (ii), whence (iii).

Suppose that~(iii) is true. Let~$Z$ be an irreducible component of~$X$ containing~$x$ and intersecting~$U$. The intersection~$F\cap Z$
is then a proper Zariski-closed subspace of the irreducible analytic space~$Z$; it follows that~$\dim {} (F\cap Z)<\dim {} Z$.

Assume that~(iv) is true and let~$V$ be an open neighborhood of~$x$. The intersection~$Z\cap V$ cannot be contained
in~$F\cap Z$ because it is of dimension~$\dim {} Z$. Therefore~$V\cap U\neq \emptyset$, whence~(i).\end{proof}

\begin{coro}\label{cor-zarclosure-dom}
Let $X$ be an analytic space, 
let $x$ be a point of $X$, let $V$ be
an analytic domain of~$X$ containing~$x$, and
let $U$ be a Zariski-open subset of $X$. The point $x$
belongs to
 $\adht U X$ if and only if~$x\in \adht {(U\cap V)}V$.
\end{coro}

\begin{proof}
If~$x\in \adht {(U\cap V)}V$ it is obvious that~$x\in \adht U X$.
Assume that~$x\in \adht U X$, and let~$Z$ be as in
assertion~(iv) of Lemma~\ref{lem-closure-zartop} above. 
Let~$T$ be an irreducible component of~$Z\cap V$ containing~$x$. It is of dimension~$\dim {} Z$. Therefore~$\dim {}(F\cap T)<\dim {} T$, 
and~$x\in \adht {U\cap V}V$ by Lemma~\ref{lem-closure-zartop}. 
\end{proof}

\begin{coro}\label{cor-zarclosure-xl}
Let $X$ be a $k$-analytic space, let $L$ be an analytic extension of~$k$, let $x$
be a point of $X$ and let $y$ be a point of~$X_{L}$ lying above~$x$. Let $U$ be a Zariski-open subset of $X$.
One has the equivalence 
$$x\in \adht U X\iff y\in \adht {U_{L}}{X_L}.$$
\end{coro}

\begin{proof}
If~$y$ belongs to $\adht {U_L}{X_L}$, then $x$ obviously belongs to $\adht U X$. Assume conversely that~$x\in \adht U X$, and let~$Z$ be as in
assertion~(iv) of Lemma \ref{lem-closure-zartop}
above. 
Let~$T$ be an irreducible component of~$Z_L$ containing~$y$. Since~$Z$ is equidimensional, $T$ is of dimension~$\dim Z$. Therefore~$\dim (F_L\cap T)<\dim  T$, 
and~$y\in \adht {U_L}{X_L}$ by Lemma~\ref{lem-closure-zartop}. 
\end{proof}

\subsection{Codimension}\label{ss-codim-def}
Let~$X$ be a~$k$-analytic space and let~$Y$ be a Zariski-closed subset of~$X$.
The \emph{codimension}~$\codim( Y,X)$ of~$Y$ in~$X$ is defined as follows.\label{IN-codim}\index{codimension}

\begin{itemize}[label=$\bullet$]

\item If both~$Y$ and~$X$ are irreducible, $\codim(Y,X)=\dim X-\dim Y.$

\item If~$Y$ is irreducible, $\codim(Y,X)=\sup_Z\codim (Y,Z)$, where~$Z$ varies
through the set of irreducible components of~$X$ that
contain~$Y$. 

\item In the general case,~$\codim(Y,X) =\inf_Z \codim (Z,X)$ where~$Z$ varies through the set of irreducible components of~$Y$.
\end{itemize}

It
is trivially checked that these definitions are consistent with each other. If
$X$ and $Y$ are non-empty and equidimensional, then 
$\codim (Y,X)=\dim X-\dim Y$, generalizing the formula which defined codimension when $X$ and $Y$ are each irreducible.

If~$x\in X$, we define the codimension of~$Y$ in~$X$ \emph{at}~$x$ as being equal to~$\inf_Z \codim (Z,X)$ where~$Z$
varies through the set of irreducible components of~$Y$ \emph{that contain~$x$}; it is denoted by~$\codim_x (Y,X)$. Note that~$\codim_x(Y, X)$ does make sense even if~$x\notin Y$, and that one has in
this case~$\codim_x (Y,X)=+\infty$ by the definition. \label{IN-codimx}

\subsection{}\label{ss-codim-basics}
Let us now list some basic properties of the codimension. 

\begin{enumerate}[1]
\item Let $X$ be an analytic space, let $Y$ be a Zariski-closed subset of $X$
and let $x$ be a point of $X$. If~$V$ is any analytic domain of~$X$ containing~$x$, 
it follows from \ref{ss-dim-andom} that
$\codim_x (V\cap Y,V)=\codim_x (Y,X).$

\item Let $X=\mathscr M(A)$ be an affinoid space and let $Y$ be a Zariski-closed subset of $X$
defined by an ideal $I$ of $A$. 
Let $x$ be a point of $X$ and let $\xi$ denote its image in $\spec A$. 
By \cite{ducros2007}, \Prop 1.11 we have the following equalities: 

\begin{enumerate}[b]
\item $\codim (Y,X)=\codim (\spec (A/I), \spec A).$

\item $\codim_x(Y,X)=\codim_\xi(\spec (A/I), \spec A)$
\end{enumerate}

\end{enumerate}

\chapter{Algebraic properties in analytic geometry}

This chapter is devoted to a general study of algebraic properties
(like being regular, Gorenstein, Cohen-Macaulay...) in analytic geometry. 
Section \ref{s-local-rings}
provides some reminders about the analytification of a scheme
of finite type over an affinoid agebra, 
and about algebraic properties of (local and global) rings of analytic functions. As a first application, 
it describes an elementary procedure (\ref{ss-approx-finite})
which we shall use quite often 
to reduce the algebraic study of local analytic rings
to that of
affinoid spaces.  

Section \ref{s-category-framework}--\ref{s-ground-field-ext}
may appear slightly unattractive. Their motivation is the following: since we
will have to deal with several kind of objects ``living on an analytic space"
(like the analytic space itself, coherent sheaves, diagrams in the category
of coherent sheaves, \etc) and with several properties, we have chosen 
to introduce a rather abstract framework, consisting of objects and properties
satisfying some axioms. From our viewpoint, this offers three advantages:

\begin{itemize}[label=$\bullet$]

\item This allows
us to write proofs once for all, and not to repeat them for every kind of object and/or
property of interest. 

\item This emphasizes which arguments are actually needed for proofs. 

\item This could be potentially applied to other objects and properties. 

\end{itemize}
But it might of course be unpleasant or boring to read. For that reason,
every important definition and statement has been given a concrete counterpart
involving only explicit objects and properties,
to which readers can directly refer if they prefer to avoid considering our 
dry formalism. 

More precisely, sections \ref{s-category-framework}
and \ref{s-alg-properties} are essentially devoted to the presentation 
of the abstract framework alluded to above. Then in
section \ref{s-appendix-analytic}, we explain what it means
for one of the properties we consider to hold 
at a point of an analytic space;
\eg, see the ``concrete" Lemma-Definition \ref{valid-at-concrete}; 
the point is that one cannot use local rings of arbitrary
(\ie, not necessarily good)
analytic spaces, hence the definition has to be given a G-local flavor, which 
requires checking some compatibilities of restriction to analytic domains. 
Thereafter we establish GAGA results about those properties;
see the ``concrete" Lemma \ref{gaga-concrete}. And we show that some 
of them have a Zariski-open locus of validity, 
sometimes automatically non-empty whenever the ambient space is reduced; 
see the ``concrete" Lemma \ref{locus-concrete}.

In section \ref{s-fibers-coherent},
we investigate
the validity 
at a point of some usual properties of a
morphism of coherent sheaves (like injectivity, surjectivity, and
bijectivity), and it almost
does not involve our general abstract setting. We prove that surjectivity
can be checked at the level of fibers (\ref{ss-surj-nakayama}; this is a
straightforward consequence of Nakayama's Lemma), and that our 
notions of injectivity, surjectivity and bijectivity are compatible 
with the usual ones in sheaf theory (\ref{ss-pointwise-sheaf}). 

In section \ref{s-ground-field-ext}
we go back to our general formalism in order to study
the behavior of algebraic properties
under ground field extension; some of them are preserved
by arbitrary such extensions, some of them (essentially, those
that involve regularity)
only
by \emph{analytically separable}
extensions; see the ``concrete" Proposition
\ref{prop-concrete-extension}, and Definition \ref{def-an-sep}
for the notion of an analytically separable extension. 

The final section of this chapter (\ref{s-complement-xan})
is essentially independent of the preceding ones. It aims at
extending some GAGA results which are known for affinoid spaces,
but not for arbitrary finitely generated scheme over an affinoid
algebra (or at least, they are not available in the literature in such generality). 
For instance, we get GAGA principles for local dimension (Lemma \ref{lem-gaga-dimension}; note
that it only works in the strict case), for codimension 
(Lemma \ref{lem-codim-kr} (3)), for normalization (Lemma \ref{lem-gaga-rednorm} (2)),
and for irreducible components (Proposition \ref{prop-gaga-irrcomp}). 

\section{Analytification of schemes, algebraic properties of analytic rings}
\markright{\thesection.~ANALYTIFICATION, ANALYTIC RINGS}
\label{s-local-rings}

\subsection{}\label{ss-def-analytify}
Let $X$ be an affinoid space, say $X=\mathscr M(A)$. We shall denote
the scheme $\spec A$ by $X\al$ (here ``al" stands for {\em algebraic}). Let $\mathscr X$ be an $X\al$-scheme
locally of finite type.
The category of good $X$-analytic spaces $Y$  endowed with a morphism {\em of locally ringed spaces} $Y\to \mathscr X$ making the diagram
$$\xymatrix{
Y\ar[r]\ar[d]&{\mathscr X}\ar[d]\\
X\ar[r]&{X\al}
}$$ commute admits a final object, which is denoted by $\mathscr X\an$ and is called the \index{analytification}
{\em analytification}
of $\mathscr X$; the  canonical map $\mathscr X\an \to \mathscr X$
is surjective, and the analytic space $\mathscr X\an$ is
relatively boundaryless over $X$ (\cite{berkovich1993}, Prop. 2.6.2). 
The space $\mathscr X\an$ is
Hausdorff, \resp compact, 
if and only if the $X\al$-scheme $\mathscr X$ is separated, \resp proper
(\cite{berkovich1993}, \Cor 2.6.7 and \Prop 2.6.9). 
The assignment 
$\mathscr X\mapsto \mathscr X\an$ is functorial in $\mathscr X$. Note that $(X\al)\an=X$. 
This construction commutes with affinoid base change: if
$Z$ is an affinoid space and if $Z\to X$ is a morphism, then $\mathscr X\an\times_XZ=(\mathscr X\times_{X\al}Z\al)\an$. 

If $x$ is a point of $\mathscr X\an$, its image
in $\mathscr X$ will be denoted by $x\al$; 
if $F$
is a subset of $\mathscr X\an$, its image in $\mathscr X$ will be denoted
by $F\al$.  
If $\mathscr Y$ is a closed (\resp open) subscheme of $\mathscr X$, then $\mathscr Y\an$
is closed analytic subspace (\resp an open subspace) of $\mathscr X\an$, which is
set-theoretically
equal to the pre-image of $\mathscr Y$ on $\mathscr X\an$. For that reason, we shall
more generally denote by $E\an$ the pre-image on $\mathscr X\an$ of any 
subset $E$ of $\mathscr X$. By surjectivity of $\mathscr X\an \to \mathscr X$,
we have $(E\an)\al=E$ for every subset $E$ of $\mathscr X$, so the assignment $E\mapsto E\an$
is injective. 
If $\mathscr F$ is a coherent sheaf on $\mathscr X$, its pull-back on $\mathscr X\an$ will be denoted by $\mathscr F\an$.

If $\mathscr X$ is proper over $X\al$, 
then non-Archimedean GAGA holds (\cf for instance \cite{poineau2010}, Annexe A;
the case where $\mathscr X=\spec A$ is essentially due to Tate and Kiehl, see~\ref{ss-tate-kiehl}): the functor 
$\mathscr F\mapsto \mathscr F\an$ induces an equivalence, which moreover preserves
cohomology, between 
the category of coherent sheaves on $\mathscr X$ and that of coherent sheaves on $\mathscr X\an$; we shall denote
by $\mathscr G\to \mathscr G\al$ a quasi-inverse of the latter. \label{IN-Gal}
Therefore $\mathscr Y\mapsto \mathscr Y\an$ induces a bijection between the set of closed
closed subschemes of $\mathscr X$ and the set of closed
analytic subspaces of $\mathscr X\an$. The inverse bijection
will be denoted by $Y\mapsto Y\al$. This implies that any Zariski-closed subset of $\mathscr X$
is of the form $E\an$ for some Zariski-closed subset $E$ of $\mathscr X$. By injectivity of $E\mapsto E\an$, 
it follows that $E\mapsto E\an$ induces a bijection between the set of Zariski-closed subsets of $\mathscr X$
and that of Zariski-closed subsets of $\mathscr X\an$; the converse bijection is induced by the assignment
$F\mapsto F\al$.

\subsection{}
\label{ss-conventions-gaga}
Let $X$ be an analytic space, let $x$
be a point of $X$ and let $\mathscr F$ be a coherent sheaf on $X$. 
Let $T\to X$ be a morphism of analytic
spaces, with $T$ affinoid. 

\begin{enumerate}[1]

\item We shall write $\mathscr F_T\al$ instead of $(\mathscr F_T)\al$; this
should not cause any confusion. \label{IN-FTal}

\item If $V$ is an affinoid domain of $X$ containing
$x$ we shall denote by $x_V\al$ the image of $x$ on $V\al$ (this is consistent with
the notation $x_V$ introduced in Remark \ref{rem-x-xv}
to indicate
that 
$x$ is viewed as belonging to $V$). 

\item If $X$ is affinoid and the affinoid space
$T$ is $X$-analytic, 
we shall write $T_x\al$ instead of $(T_x)\al$; this should
not cause any confusion. 
We shall denote by $t_x\al$ the image of $t$ on $T_x\al$ (this is consistent with
the notation $t_x$ introduced in Remark \ref{rem-x-xv}
to indicate
that $t$ is viewed as belonging to the fiber
$T_x$). We shall use
(in acccordance with our general conventions
in scheme theory)
the notation $T\al_{x\al}$
for the scheme-theoretic fiber of $T\al \to X\al$
at $x\al$. 

\end{enumerate}

\subsection{Algebraic properties of analytic rings}\label{ss-algprop-analytic}
Let $k$ be an analytic field and let $L$ be an analytic extension of $k$. Let $A$ be a $k$-affinoid algebra, let 
$B$ be the algebra of analytic functions on some affinoid domain  of $\mathscr M(A)$. 
Let $X$ be a good $k$-analytic space, let $V$ be a good analytic domain of $X$, and let
$x$ be a point of $V$.

\begin{enumerate}[1]

\item The ring $A$ is excellent (\cite{ducros2009}, \Th 2.13; the strictly affinoid case is due to Kiehl \cite{kiehl1969}).

\item\label{aa} The $A$-algebra $B$ is regular (\cite{ducros2009}, \Th 3.3; flatness follows from \cite{berkovich1990}, \Prop 2.2.4 and from \cite{bosch-g-r1984}, \Cor 6 of \S 7.3.2 in the strict case). 

\item The $A$-algebra $A_L$ is faithfully flat  (\cite{berkovich1993}, Lemma 2.1.2).

\item The local ring~$\mathscr O_{X,x}$ is noetherian, henselian (\cite{berkovich1993}, \Th 2.1.4 and \Th 2.1.5), and excellent (\cite{ducros2009}, \Th 2.13).

\item The morphism $\mathscr O_{X,x}\to \mathscr O_{V,x}$ is regular (\cite{ducros2009}, \Th 3.3; flatness is a straightforward consequence
of statement (2) above). 

\item The morphism $X_L\to X$ is flat when viewed as a morphism of locally ringed spaces (\cite{berkovich1993}, \Cor 2.1.3). 

\end{enumerate}

\subsection{Regularity of analytification}\label{ss-gaga-alg}
Let $\mathscr X$ be a scheme 
locally of finite type over an affinoid algebra. For every $x\in \mathscr X\an$, the 
morphism $\mathscr O_{\mathscr X,x\al}\to \mathscr O_{\mathscr X\an, x}$ is regular (\cite{ducros2009}, \Th 3.3; flatness is due
to Berkovich, \cite{berkovich1993}, \Prop 2.6.2).

\subsection{}\label{ss-gaga-domain}
Let $X$ be an affinoid space and let $x$ be a point of $X$. Faithful flatness
of $\mathscr O_{X\al,x\al}\to \mathscr O_{X,x}$
has the following consequences. 

\begin{enumerate}[1]

\item Assume that $\mathscr O_{X,x}$ is a domain. Being a subring 
of $\mathscr O_{X,x}$, the local ring $\mathscr O_{X\al, x\al}$
is a domain too. This implies that $x\al$ lies on a unique irreducible
component $\mathscr X$ of $X\al$; hence $x$ lies on a unique irreducible component of $X$,
namely $\mathscr X\an$
(for another proof of this fact, see~\cite{ducros2009}, Lemme 0.11).

\item
Assume moreover that the domain $\mathscr O_{X,x}$ is a field. By surjectivity
of the map $\spec \mathscr O_{X,x}\to \spec \mathscr O_{X\al, x\al}$, 
the scheme $\spec \mathscr O_{X\al, x\al}$
consists of one point; therefore the domain $\mathscr O_{X\al, x\al}$ is a field and $x\al$ is the generic point of $\mathscr X$. 
\end{enumerate}

\subsection{Approximation of finite algebras}\label{ss-approx-finite}
Let $X$ be a good analytic space
and let $x$ be a point of $X$. 
Let $B$ be a finite $\mathscr O_{X,x}$-algebra. Since $\mathscr O_{X,x}$ is noetherian, $B$ is finitely presented. 
Therefore there exists an affinoid neighborhood $X'$ of $x$ in $X$ and a finite $\mathscr O_X(X')$-algebra $R$ such that
$B=\mathscr O_{X,x}\otimes_{\mathscr O_X(X')}R$. The affinoid space $Z:=\mathscr M(R)$ is finite over $X'$. We shall say
that $B$ is {\em induced} by the finite $X'$-analytic space $Z$; if $B$ is a quotient of $\mathscr O_{X,x}$, we can chose
$X'$ and $R$ so that $R$ is a quotient of $\mathscr O_X(X')$, and $Z$ is then a
closed analytic subspace of $X'$. 
If $z_1, \ldots, z_r$ denote the pre-images of $x$ on $Z$, one has $B=\prod_i \mathscr O_{Z,z_i}$ because
$Z\to X'$ is topologically
closed (for details see \cite{berkovich1993}, Lemma 2.1.6, which is the key point for proving that $\mathscr O_{X,x}$ is henselian).

Assume moreover that $B$ is a domain.
The point $x$ has thus only one pre-image $z$ on $Z$, 
and since $\mathscr O_{Z,z}=B$ is a domain, 
it follows from \ref{ss-gaga-domain}
(1)
that $z$ lies on a unique irreducible component
$T$ of $Z$. Let $\mathscr I$ be the sheaf of ideals on $Z$ that defines $T_{\mathrm {red}}$. 
Since $T$ is the unique irreducible component of $Z$ containing $z$, this is a neighborhood of 
$z$ in $Z$; therefore, any section of $\mathscr I$ vanishes pointwise in a neighborhood of $z$, 
hence is nilpotent in $\mathscr O_{Z,z}$, which is a domain. As a consequence, $\mathscr I_z=0$
and the closed analytic subspace $T_{\mathrm {red}}$ of $Z$ induces 
the identity map $\mathrm{Id}_{\mathscr O_{Z,z}}$. Hence the finite $X'$-space
$T_{\mathrm{red}}$ still induces the finite $\mathscr O_{X,x}$-algebra
$B$. 

Therefore any finite $\mathscr O_{X,x}$-algebra
which is a domain (\resp any quotient of $\mathscr O_{X,x}$ by a prime ideal)
can be induced by an \emph{integral}
finite $X'$-space (\resp an \emph{integral}
closed analytic subspace of $X'$) for some
affinoid neighborhood $X'$ of $X$.

\section{A rather abstract categorical framework}\label{s-category-framework}

In this memoir
we deal with with various kinds of
geometric objects: analytic spaces equipped with their G-topology and the corresponding structure sheaf, 
schemes of finite type 
over an affinoid algebra, (spectra of) local rings of good analytic spaces, etc., 
and the goal of this section is to define a category $\mathfrak T$ that encompass
all of them. 

Technically speaking, all
aforementioned objects
are locally ringed toposes, and all relevant morphisms between them are morphisms
of locally ringed toposes (the reader should not be afraid about that: one does not need to know
what a topos is, let alone a morphism of toposes; we will only use these concepts here as 
unifying terminology). This leads to the following definition of $\mathfrak T$. 

\begin{defi}
\label{def-category-framework}
We denote by $\mathfrak T$ the smallest subcategory of the category of locally ringed toposes such that
the following hold (we allow ourselves to write for short that a given object, \resp arrow, {\em belongs}
to $\mathfrak T$):

\begin{enumerate}[1]
\item If $X$ is an analytic space, then $X$ (viewed as equipped with its G-topology and the corresponding structure sheaf)
belongs to $\mathfrak T$. 

\item Let $A$ be a ring which is either an affinoid algebra or of the form $\mathscr O_{X,x}$ with $X$ a good analytic space
and $x$ a point of $X$. If $B$ is any $A$-algebra essentially of finite type, the affine scheme $\spec B$ belongs to $\mathfrak T$. 

\item Let $\mathscr X$ be a scheme. If it admits a
covering by open subschemes that belong to $\mathfrak T$, then $\mathscr X$ belongs to $\mathfrak T$. 

\item Any morphism $Y\to X$ between analytic spaces belongs to $\mathfrak T$.

\item Let $\mathscr Y\to \mathscr X$ be a morphism of schemes. If both $\mathscr Y$ and $\mathscr X$ belong  to $\mathfrak T$, then
$\mathscr Y\to \mathscr X$ belongs to $\mathfrak T$.

\item Let $\mathscr X$ be a scheme of finite type over an affinoid algebra. The morphism of locally ringed toposes
$\mathscr X\an \to \mathscr X$ belongs to $\mathfrak T$.

\end{enumerate}
\end{defi}

\begin{rema}
By construction, all objects of $\mathfrak T$ have a coherent structure sheaf, hence
admit a nice theory of coherent sheaves. 
\end{rema}

\begin{rema}\label{rem-scheme-t}
The only schemes that are required
to belong to $\mathfrak T$ are those mentioned in (2) and (3). 
Therefore a given scheme belongs to $\mathfrak T$ if and only if it admits
an open covering by affine schemes of the form described in (2).

Note that any field can be seen as an affinoid algebra (once equipped with the trivial 
absolute value); hence any scheme of finite type over a field belongs to $\mathfrak T$. 

Let $A$ be a complete equicharacteristic local notherian ring. It is
a quotient of a formal power series ring $k[\![T_1,\ldots, T_n]\!]$ (\cf
\cite{matsumura1986}, Thm. 28.3 and
the discussion at the beginning of \S 29).
Since the latter can be seen 
as an affinoid algebra over the trivially valued field $k$ (it is isomorphic to $k\{T_1/r_1,\ldots, T_n/r_n\}$ 
as soon as all $r_i$'s are smaller than $1$), the affine scheme $\spec A$ belongs to $\mathfrak T$. 
\end{rema}

\subsection{}\label{ss-embed-regular}
Let us list some consequences of Remark~\ref{rem-scheme-t}
above. 

\begin{enumerate}[1]

\item Any scheme belonging to $\mathfrak T$
is excellent, and even locally
embeddable into a {\em regular} excellent scheme. Indeed if $k$
is an analytic field and $X$ a $k$-affinoid space, then
$X$ admits a closed immersion into some compact polydisc $Y$
over $k$. Since the local rings of $Y$ as well as those of $Y\al$ are regular 
(for an elementary proof, see \cite{ducros2009}, Lemme 2.1), the scheme
$X\al$ admits a closed immersion into a regular excellent scheme, and every
local ring of $X$ is the quotient of an excellent regular local ring, whence our claim. 

\item Let $\mathscr X$ be a scheme. If $\mathscr X$ belongs to $\mathfrak T$, every $\mathscr X$-scheme
locally of finite type belongs to $\mathfrak T$, and  $\spec \mathscr O_{\mathscr X,x}$ belongs to $\mathfrak T$
for every point $x$ of $\mathscr X$. 

\item Let $\mathscr Y\to \mathscr X$ be a morphism of schemes.  If both $\mathscr Y$ and $\mathscr X$
belong to $\mathfrak T$, then $\mathscr Y_x$
and $\mathscr Y\times_{\mathscr X}\spec \mathscr O_{\mathscr X,x}$
belong to $\mathfrak T$ for every $x\in \mathscr X$.  

\end{enumerate}

\subsection{}\label{ss-fiber-category}
We fix a fibered category~$\mathfrak F$
over~$\mathfrak T$; if $ Y\to X$ is an arrow of $\mathfrak T$, the corresponding pull-back functor
$\mathfrak F_X\to \mathfrak F_Y$ will be denoted by $D\mapsto D_Y$. If $A$
is a ring whose spectrum belongs to $\mathfrak T$, we shall sometimes write for short
$\mathfrak F_A$ instead of $\mathfrak F_{\spec A}$, 
and $D_A$ instead of $D_{\spec A}$. If $\mathscr X$ is any scheme belonging
to $\mathfrak T$,
$x$ is a point of $\mathscr X$, and 
$D$ is an object of $\mathfrak F_{\mathscr X}$, we shall write $D_x$ instead
of $D_{\mathscr O_{\mathscr X,x}}$. If $\mathscr X$ is a scheme of finite type over an affinoid algebra, 
the pull-back functor $\mathfrak F_{\mathscr X}\to \mathfrak F_{\mathscr X\an}$ will be denoted by $D\mapsto D\an$. 
If $k$ is an analytic field, $X$ is a $k$-analytic
space, and $L$ is a complete extension of $k$, the pull-back functor
$\mathfrak F_X\to \mathfrak F_{X_L}$ will be denoted by $D\mapsto D_L$. 
We assume moreover
that $\mathfrak F$ satisfies the following property:

\begin{enonce*}{GAGA axiom for $\mathfrak F$}
For every affinoid space, the pull-back functor
$D\mapsto D\an$ from $\mathfrak F_{X\al}$ to $\mathfrak F_X$ is an equivalence; we denote
by $D\mapsto D\al$ a quasi-inverse of the latter.
\end{enonce*}

\subsection{}\label{ss-notation-dx}
Let $X$ be an analytic space and let $x$ 
be a point of $X$.
Let us first assume that $X$ is good. Let $V$ be an affinoid neigborhood of $x$ in $X$. The composition 
functor 
$$\xymatrix{
{\mathfrak F_X}\ar[rr]&&{\mathfrak F_V}\ar[rr]_{D\mapsto D\al}&&
{\mathfrak F_{V\al}}\ar[rr]&&\mathfrak F_{\mathscr O_{X,x}}
}
$$
only depends on $x$, and not on $V$. It will be denoted by $D\mapsto D_x$. 
If $\mathscr X$ is any $\mathscr O_{X,x}$-scheme belonging to $\mathfrak T$ and if $D$
is an object of $\mathfrak F_X$ we shall often write
$D_{\mathscr X}$ instead of $(D_x)_{\mathscr X}$, if there is no ambiguity;  for example, we shall use the notation
$D_{\kappa(x)}$ and $D_{\hr x}$. 

\subsection{}\label{ss-mathfrak-l}
Let $\mathfrak L$ be the category whose objects
are local noetherian rings $A$ such that $\spec A$ belongs to $\mathfrak T$, and
whose arrows are {\em local}
maps. 
The contravariant functor $\spec$
anti-identifies $\mathfrak L$ with a subcategory
$\spec \mathfrak L$ of $\mathfrak T$. The fibered category $\mathfrak F$ gives rise by restriction 
to a fibered category over $\spec \mathfrak L$; we shall denote it by $\mathfrak F_{\mathfrak L}$. 

We are now going to give three examples
the reader should keep in mind while working with this general, abstract fibered category $\mathfrak F$. For each of them, we shall
only describe the corresponding fiber categories;
the definition of
pull-back functors is straightforward and left to the reader.

\begin{exem}\label{ex-fiber-trivial}
We may take for $\mathfrak  F$ the category $\mathfrak T$, viewed as fibered category over itself in the obvious way: for every object $X$ of $\mathfrak T$, 
the fiber category $\mathfrak T_ X$ 
simply consists of the single object $X$ with $\mathrm{Id}_X$ as unique endomorphism. 

The fibered category $\mathfrak T_{\mathfrak L}$ can then be anti-identified with $\mathfrak L$. 
This allows us to view objects of $\mathfrak T_{\mathfrak L}$ as objects of $\mathfrak L$. 

\end{exem}

\begin{exem}\label{ex-fiber-coh}
We may take for $\mathfrak F$ the category $\mathfrak{Coh}\to \mathfrak T$ of coherent sheaves, defined as follows: for every object $X$ of $\mathfrak T$, the fiber category
$\mathfrak{Coh}_X$ is the category of coherent sheaves on $X$. 

In particular, objects of $\mathfrak{Coh}_{\mathfrak L}$ can be seen as pairs $(A,M)$ with $A\in \mathfrak L$
and $M$ a finitely generated $A$-module. 
\end{exem}

\begin{exem}\label{ex-fiber-diag}
Let $\mathfrak I$ be
a small category. We may take for $\mathfrak F$ the
fibered category~$\mathfrak{Coh}^{\mathfrak I}\to \mathfrak T$, defined as follows: for every object $X$ of $\mathfrak T$, the fiber category $\mathfrak{Coh}^{\mathfrak I}_X$
is the category of $\mathfrak I$-diagrams of coherent sheaves on $X$; i.e., of covariant functors from~$\mathfrak I$ to~$\mathfrak{Coh}_X$ (morphisms are natural transformations of functors). 

In particular, objects of~$\mathfrak{Coh}^{\mathfrak I}_{\mathfrak L}$ can be seen as pairs $(A,D)$ with $A\in \mathfrak L$
and $D$ an $\mathfrak I$-diagram of finitely generated $A$-modules. 
\end{exem}

\begin{rema}\label{rem-diag-includes}
Examples~\ref{ex-fiber-trivial}
and~\ref{ex-fiber-coh}
can actually been interpreted as
particular cases of Example~\ref{ex-fiber-diag}.

Indeed, let us begin with Example~\ref{ex-fiber-trivial}. For any category $\mathfrak C$, 
there is a unique functor from $\emptyset$ to $\mathfrak C$. It simply does {\em nothing},
because the source category has no object (it is analogous to the empty map from $\emptyset$ to an arbitrary set); and it has
a unique endomorphism. Hence we have for every object $X$ of $\mathfrak T$ an equivalence $\mathfrak{Coh}^\emptyset_X\simeq \mathfrak T_X$, 
both categories involved having a single object and a single morphism.  Those equivalences essentially commute
with pull-back
functors, 
whence
an equivalence $\mathfrak{Coh}^\emptyset\simeq \mathfrak T.$

Let us now deal with Example~\ref{ex-fiber-coh}. Let $\{\ast \}$ be the category with one single element and one single
morphism. For any category $\mathfrak C$, the assignment $F\mapsto F(\ast)$ induces
an equivalence between the category of covariant functors from $\{\ast \}$ to $\mathfrak C$
and  $\mathfrak C$ itself. This yields for every object $X$ of $\mathfrak T$ an equivalence
$\mathfrak {Coh}^{\{\ast  \}}_X\simeq \mathfrak {Coh}_X$.
Those equivalences essentially commute
with pull-back functors, 
whence
eventually an equivalence $\mathfrak{Coh}^{\{\ast\}}\simeq \mathfrak{Coh}.$
 \end{rema}

 \begin{enonce}[remark]{Convention}\label{ss-conf-mathfrakl}
If $\mathfrak F$ is one of fibered categories considered
in Examples \ref{ex-fiber-trivial}, \ref{ex-fiber-coh}
and \ref{ex-fiber-diag}, objects
of $\mathfrak F_{\mathfrak L}$ can be interpreted as objects from
commutative algebra. We shall always use this viewpoint here. In other words, 
we shall consider that objects of $\mathfrak T_{\mathfrak L}$
 (\resp $\mathfrak{Coh}_{\mathfrak L}$, \resp $\mathfrak{Coh}^{\mathfrak I}_{\mathfrak L}$) actually {\em are}
 objects of $\mathfrak L$ (\resp pairs $(A,M)$ as in
 Example~\ref{ex-fiber-coh},  \resp pairs $(A,D)$ as in Example~\ref{ex-fiber-diag}).

 \end{enonce}

\section{Formalisation of algebraic properties}\label{s-alg-properties}

\subsection{}\label{ss-alg-properties}
We fix from now on
until the end of section \ref{s-ground-field-ext}
a property~$\mathsf P$ whose validity makes sense for every object of~$\mathfrak F_{\mathfrak L}$.

\begin{exem}\label{ex-algprop-space}
If~$\mathfrak F=\mathfrak T$ then objects of $\mathfrak F_{\mathfrak L}$ are local noetherian rings
(belonging to $\mathfrak L$). We can therefore
take for~$\mathsf P$ the property of being
regular, Gorenstein, CI, or $R_m$ for some specified $m$ (for a definition of the $R_m$
property, see~\cite{ega42}, Def. 5.8.2).

\end{exem}

\begin{exem}\label{ex-algprop-coh} If~$\mathfrak F=\mathfrak{Coh}$ then objects of $\mathfrak F_{\mathfrak L}$ are pairs $(A,M)$ with
$A$ a local noetherian ring (belonging to $\mathfrak L$) and $M$ a finitely generated $A$-module.
We can therefore take for~$\mathsf P$ the property of being CM, or free of given rank, or of given residue rank, or of given depth or codepth, 
or $S_m$ for some specified $m$ (for a definition of the $S_m$
property, see \cite{ega42}, 5.7.2; note that the zero module is $S_m$ for every $m\geq 0$).\index{depth of a module}\index{codepth of a module}

The notions of depth and codepth play a crucial role in this memoir for the construction
of d\'evissages. Recall that they are related to each other by the formula
\[\mathrm{codepth}_A(M)=\dim_{\text{Krull}}M-\mathrm{depth}_A(M) \]
if $M\neq 0$ (the Krull dimension of $M$ is defined in \ref{ss-conv-alg}); if $M=0$ we have $\dim_{\text{Krull}}M=-\infty$
and $\mathrm{depth}_A(M)=+\infty$, but $\mathrm{codepth}_A(M)=0$
\emph{by convention}. 
\end{exem}

\begin{exem}\label{ex-algprop-complex}
If~$\mathfrak I$ is an interval of~$\Z$ (viewed as a category through its natural ordering) and if $\mathfrak F=\mathfrak{Coh}^{\mathfrak I}$,
then objects of $\mathfrak F_{\mathfrak L}$ are pairs $(A,D)$ with $A$ a local
noetherian ring (belonging to $\mathfrak L$) and $D$ an $\mathfrak I$-indexed sequence
$$\ldots \to M_{i-1}\to M_i\to M_{i+1}\to \ldots$$ of~$A$-linear maps between finitely generated
$A$-modules. 
We can therefore take for~$\mathsf P$ the property of being a complex, or of being exact, or of being exact
at some specified position $i\in \mathfrak I$.  
\end{exem}

\begin{rema}
The properties considered in Example \ref{ex-algprop-coh}
also make sense for $\mathfrak F=\mathfrak T$, by viewing any local noetherian ring (belonging to $\mathfrak L$)
as a module over itself. 

\end{rema}

\begin{rema}\label{rem-norm-r1s2}
One might also of course consider some relevant 
combinations of the aforementioned properties. The most important
examples the reader should have in mind
are: the property of being reduced, which amounts
satisfying both $R_0$ and $S_1$; and that of being normal, which
amounts satisfying both $R_1$ and $S_2$
(\cf \cite{ega42}, \Prop 5.8.5 and \Th 5.8.6). 

\end{rema}

\begin{defi}\label{def-valid-good}
Let $X$ be either a good analytic space or a scheme
belonging to $\mathfrak T$. Let $x$ be a point of $X$ and let $D$ be an object
of $\mathfrak F_X$. We shall say that $D$ {\em satisfies $\mathsf P$ at $x$}
if $D_x$ satisfies $\mathsf P$.

\end{defi}

\begin{rema}\label{rem-valid-localization}
If $D$ satisfies $\mathsf P$ at every point of~$X$, one could
be tempted to say for short that $D$ satisfies $\mathsf P$, but the reader should be aware that such
terminology might be ambiguous. Indeed, if $X$ belongs to $\spec \mathfrak L$,
it might happen that $D$ satisfies $\mathsf P$ in the original meaning
(validity of $\mathsf P$ makes sense for any object of $\mathfrak F_{\mathfrak L}$), but
does not satisfy $\mathsf P$ everywhere at $X$. 

By definition, this problem cannot occur as soon as $\mathsf P$ satisfies condition \gen~below.
Hence under this assumption, we shall actually use the expression ``$D$ satisfies $\mathsf P$"
instead of ``$D$ satisfies $\mathsf P$ at every point of $X$". 
\end{rema}

\subsubsection*{Condition \gen}
For every local noetherian ring $A$ belonging to $\mathfrak L$, every prime ideal $\mathfrak p$
of $A$ and every object $D$ of $\mathfrak F_A$, the following implication holds: 

$$(D\;\mathrm{satisfies~}\mathsf P)\Rightarrow(D_{A_{\mathfrak p}}\;\mathrm{satisfies~}\mathsf P).$$

This amounts requiring
that for
every scheme $\mathscr X$ belonging to $\mathfrak T$ and every object 
$D$ of $\mathfrak F_{\mathscr X}$, 
the set of points of $\mathscr X$ at which $D$ satisfies $\mathsf P$ is stable under generization. 

\begin{exem}

The following properties satisfy \gen : if $\mathfrak F=\mathfrak T$, 
the property of being CI \cite{avramov1975}, of being Gorenstein (\Th 18.2 of \cite{matsumura1986}), 
regular (Serre, see \Th 19.3 of \cite{matsumura1986}), or $R_m$ by its very definition; if $\mathfrak F=\mathfrak {Coh}$, the property
of being CM (\cite{matsumura1986}, \Th 17.3), of being $S_m$ (by its very definition) or of being of codepth bounded above by $m$
for some specified $m$
(see \cite{ega42}, \Prop 6.11.5); if $\mathfrak F=\mathfrak {Coh}^{\mathfrak I}$ for some
interval $\mathfrak I$ of $\Z$, the property of being exact at some specified
position $i\in \mathfrak I$.

As far as the CI property is concerned, one may give an alternative, simpler proof, using the facts that any local ring
belonging to $\mathfrak L$ is a quotient of a regular local ring; see for instance the discussion at the end of~\cite{matsumura1986}, \S 21. 

\end{exem}

\begin{rema}
Let $n$ be a non-negative integer. The property of being
free of rank $n$ satisfies \gen, but
we shall of course use ``locally free of rank $n$"
for ``free of rank $n$ at every point", and not ``free of rank $n$"
which we reserve for globally free sheaves as usual. 
\end{rema}

\subsection{Geometric validity}
Let $k$ be a field and let $\mathscr X$ be a $k$-scheme belonging to $\mathfrak T$. Let $x$
be a point of $\mathscr X$ and let $D$ be an object of $\mathfrak F_{\mathscr X}$. We shall say
that $D$ satisfies $\mathsf P$ {\em geometrically}
at $x$, or for short that $D$ satisfies $\mathsf P_{\mathrm{geo}}$
at $x$, if for every finite extension $F$ of $k$,
the object $D_{\mathscr X\times_k F}$ satisfies $\mathsf P$ at every point of $\mathscr X\times_k F$. 
If $\mathsf P$ satisfies condition \gen, we shall say that $D$ satisfies $\mathsf P_{\mathrm{geo}}$ if it does
so at every point of $\mathscr X$. 

Note that if $\mathsf P$ satisfies \gen, the set of points of $\mathscr X$ at which $D$ satisfies 
$\mathsf P_{\mathrm{geo}}$ is stable under generization: this  is a formal consequence
of the fact that finite morphisms are closed. 

\begin{defi}
Let $\mathscr Y\to \mathscr X$
be a morphism of schemes belonging to $\mathfrak T$, and let
$D$ be an object of $\mathfrak F_{\mathscr X}$. Let $y$ be a point of $\mathscr Y$
and let $x$ be its image on $\mathscr X$. We shall say that $D$ satisfies $\mathsf P$
(\resp $\mathsf P_{\mathrm{geo}}$)
{\em fiberwise}
at $y$ if $D_{\mathscr Y_x}$ satisfies $\mathsf P$ at $y$ (\resp satisfies $\mathsf P_{\mathrm{geo}}$ at $y$ with respect 
to the canonical morphism $\mathscr Y_x\to \spec \kappa(x)$). 
If $\mathsf P$
satisfies \gen, we shall say that $D$ satisfies $\mathsf P$ (\resp $\mathsf P_{\mathrm{geo}}$)
fiberwise
if it does so at every point. 
\end{defi}

The study of fiberwise validity of $\mathsf P$ or $\mathsf P_{\mathrm{geo}}$ can often be reduced to 
that of
the usual validity of $\mathsf P$ on the source space through 
a standard trick, which is frequently used in SGA.

\subsection{The standard trick}\label{ss-trick-fibval}
Suppose that we are given a morphism $\mathscr Y\to \mathscr X$
of schemes belonging to $\mathfrak T$, a point $x$ on $\mathscr X$
and an object $D$ of $\mathfrak F_{\mathscr Y}$. Let $F$ be a finite extension of
$\kappa(x)$ (if we are interested in $\mathsf P$ we shall only consider the case where $F=\kappa(x)$)
and let $y$ be a point of $\mathscr Y_x\times_{\kappa(x)}F$. 
Let us endow $\overline{\{x\}}$ with its reduced structure, and let $\mathscr Z$ be any integral 
finite $\overline{\{x\}}$-scheme whose function field
is equal to $F$; \eg, we can
take for $\mathscr Z$ the normalization of the Japanese scheme $\overline{\{x\}}$ inside $F$, or
the scheme $\overline{\{x\}}$ itself if $F=\kappa(x)$. 
Let $\xi$ be the generic point of $\mathscr Z$, and set $\mathscr T=\mathscr Y\times_{\mathscr X}\mathscr Z$.
The scheme $\mathscr Y_x\times_{\kappa(x)}F$ is nothing but the generic fiber $\mathscr T_\xi$. Therefore
$\mathscr O_{\mathscr Y_x\times_{\kappa(x)}F,y}=\mathscr O_{\mathscr T_\xi,y}=\mathscr O_{\mathscr T, y}$. 
Hence $D_{\mathscr Y_x\times_{\kappa(x)}F}$ satisfies $\mathsf P$ at $y$ if and only if so does $D_{\mathscr T}$. 

\subsection{Analytic version of the standard trick} \label{ss-trick-fibvalan}
Suppose that we are given a morphism $Y\to X$
of good analytic spaces, a point $y$ on $Y$
whose image in
$X$ is denoted by $x$, 
and an object $D$ of $\mathfrak F_Y$. Let $\xi$ be a point
of $\spec \mathscr O_{X,x}$, let $F$ be a finite extension of
$\kappa(\xi)$ (if we are interested in $\mathsf P$ we shall only consider the case where $F=\kappa(\xi)$), 
and let $\eta$ be a point of $(\spec \mathscr O_{Y,y})_\xi \times _{\kappa(\xi)}F$. Let $B$ 
be the quotient of $\mathscr O_{X,x}$ by its prime ideal that corresponds to $\xi$. 

Let $C$ be any finite $B$-subalgebra of $F$
with fraction field $F$; \eg, we may take for $C$ the integral closure
of $B$ in $F$ (the ring $\mathscr O_{X,x}$ is universally japanese), or the ring $B$ itself
when $F=\kappa(\xi)$. Let us shrink $X$ so that $X$ is affinoid and so that
$C$ is induced by a finite integral $X$-analytic space $Z$; \cf \ref{ss-approx-finite}.
Let $z$ be the unique pre-image of $x$ in $Z$. 
Set $T=Y\times_X Z$ and let $t_1,\ldots, t_r$ be the pre-images of $y$ on $T$.
Let $\zeta$ be the generic point of $\spec \mathscr O_{Z,z}$; note that $\zeta$ lies over the generic
point of $Z\al$ by flatness of $\mathscr O_{Z\al,z\al}\to \mathscr O_{Z,z}$. 
Let $p$ the natural map $\coprod \spec \mathscr O_{T,t_i}\to \spec \mathscr O_{Z,z}$. 
One has
\begin{eqnarray*}
\prod_i \mathscr O_{T,t_i}&=&\mathscr O_{T}(T)\otimes_{\mathscr O_Y(Y)}\mathscr O_{Y,y}\\
&=&\mathscr O_{Z}(Z)\otimes_{\mathscr O_X(X)}\mathscr O_{Y,y}\\
&=&\mathscr O_{Z,z}\otimes_{\mathscr O_{X,x}}\mathscr O_{Y,y}.
\end{eqnarray*}
Therefore $(\spec \mathscr O_{Y,y})_\xi \times _{\kappa(\xi)}F$ is nothing
but the generic fiber $p^{-1}(\zeta)$. 
The point $\eta$ lies on $\spec \mathscr O_{T,t_i}$ for some $i$, and the local rings
at $\eta$ of
$\spec \mathscr O_{T,t_i}$ and of $p^{-1}(\zeta)$ coincide.
It follows that $D_{(\spec \mathscr O_{Y,y})_\xi \times _{\kappa(\xi)}F}$
satisfies $\mathsf P$ at $\eta$ if and only if
$D_{T,t_i}$ satisfies $\mathsf P$ at $\eta$. 

\subsection{}\label{ss-list-hreg}
We are now going to introduce various technical conditions that make sense for $\mathsf P$.

\subsubsection*{Condition \hreg}For any flat morphism~$A\to B$ of $\mathfrak L$
and any~$D\in \mathfrak F_A$ the following implications hold: 

\begin{itemize}[label=$\bullet$]
\item If $D_B$ satisfies $\mathsf P$ then $D$ satisfies $\mathsf P$. 

\item If $D$ satisfies $\mathsf P$ and if moreover the fibers of $\spec B\to \spec A$ are regular then $D_B$ satisfies $\mathsf P$. 
\end{itemize}

\subsubsection*{Condition \hci}
For any flat morphism~$A\to B$  of $\mathfrak L$
and any~$D\in \mathfrak F_A$ the following implications hold: 

\begin{itemize}[label=$\bullet$]

\item If $D_B$ satisfies $\mathsf P$ then $D$ satisfies $\mathsf P$.

\item If $D$ satisfies $\mathsf P$ and if moreover the fibers of $\spec B\to \spec A$ are CI then $D_B$ satisfies $\mathsf P$. 
\end{itemize}

\subsubsection*{Condition \hv}
For any flat morphism~$A\to B$ of $\mathfrak L$
and any~$D\in \mathfrak F_A$,  the object $D$ satisfies $\mathsf P$ if and only if $D_B$ satisfies $\mathsf P$.

\subsubsection*{Condition \field}
For every field $k$, every object $D$ of $\mathfrak F_k$ satisfies $\mathsf P$.

\subsubsection*{Condition \open}
For any scheme $X$ belonging to $\mathfrak T$ and any object $D$ of $\mathfrak F_X$, the subset of $X$
consisting of points at which $D$
satisfies $\mathsf P$ is open.

\begin{rema}
It follows from the definitions that~\hv$\Rightarrow$\hci$\Rightarrow$\hreg~and that \open$\Rightarrow$\gen. 
\end{rema}

\begin{exem}\label{ex-hreg-rings}
Assume that~$\mathfrak F=\mathfrak T$. Let $m$ be an integer. 

The following properties satisfy \hreg~: being regular, and being~$R_m$; see \cite{ega42},
\Prop 6.5.1 \Prop 6.5.3 (ii). The following properties satisfy~\hci:
being Gorenstein, see \cite{matsumura1986}, \Th 23.4;  being CI, see \cite{avramov1975}.

The properties of being regular, $R_m$, Gorenstein, and CI obviously satisfy \field.
They also satisfy \open: see \cite{ega42}, Scholie 7.8.3 (iv)
for regularity and the $R_m$ property; see \cite{greco-m1978}, \Cor 1.5 
for the Gorenstein property; and see \cite{greco-m1978}, \Cor 3.3 for the CI property. 

As far as condition \open~for Gorenstein and CI properties is concerned, one may find simpler, 
alternative proofs using the fact that every scheme belonging to $\mathfrak T$ is locally embeddable
into a regular scheme; see \cite{matsumura1986}, Exercise 24.3 and \cite{ega44}, 
\Cor 19.3.3. 

\end{exem}

\begin{exem}\label{ex-hreg-modules}
Assume that~$\mathfrak F=\mathfrak{Coh}$. Let $m$ be an integer. 

The properties of being CM, of being~$S_m$
and of being of codepth~$m$
satisfy \hci; see~
\cite{ega42}, \Cor 6.3.2, \Prop 6.4.1 (i) and (ii). 
The properties of being CM, of being $S_m$
and of being of codepth $\leq m$ obviously satisfy \field. They
also satisfy \open; see \cite{ega42},  Scholie 7.8.3 (iv). 

The property of being free satisfies \hv, see
\cite{ega42}, \Prop 6.2.1 (ii); the property of being of given residue rank
obviously satisfies
\hv, hence being free of specified rank also satisfies \hv. 

The property of being free satisfies \field~and  \open. The property of being free {\em of specified rank}
still satisfies \open~but it does not satisfy \field.

\end{exem}

\begin{exem}\label{ex-hreg-diag}
Assume that~$\mathfrak F=\mathfrak{Coh}^{\mathfrak I}$ 
for~$\mathfrak I$ an interval of~$\Z$.  We will only
list here some properties satisfying~\hv: the properties of being a complex, of being exact, of being a complex
having its~$i$-th homology of a given residue rank for some specified~$i\in \mathfrak I$.
In the particular case where~$\mathfrak I=\{0,1\}$ (in which case~$\mathfrak{Coh}^{\mathfrak I}_X$
is the category of maps between two coherent sheaves on~$X$),
let us mention the properties of being an isomorphism, an injection, a surjection. 

\end{exem}

\subsection{}\label{ss-fiberwise-chiant}
Assume that $\mathsf P$ satisfies \hreg, and let 
$\mathscr C$ be a class of morphisms between analytic spaces that is
stable under base change by inclusions of analytic domains
and
finite morphisms. Let us consider the following four assertions. 

\begin{itemize}
\item[(A)] Let $Y\to X$ be an arrow belonging to $\mathscr C$ with $Y$ and $X$ affinoid and with
$X$ integral. Let $D$ be an object of $\mathfrak F_{Y\al}$. The object $D$ satisfies $\mathsf P$ at every point of $Y\al$ lying
above the generic point of $X\al$.

\item[($\text A^\ast$)] Let $Y\to X$ be an arrow belonging to $\mathscr C$, with $Y$ and $X$ affinoid and with $X$ integral. 
Let $y$ be a point of $Y$
 and let $x$ be its image in $X$. Let $D$ be an object of $\mathfrak F_Y$.
 The object $D_y$ satisfies $\mathsf P$
at every point of $\spec \mathscr O_{Y,y}$ lying
above the generic point of $X\al$.

\item[(B)] Let $Y\to X$ be an arrow belonging to $\mathscr C$, with $Y$ and $X$ affinoid. 
Let $D$ be an object of $\mathfrak F_{Y\al}$. The object $D$ satisfies $\mathsf P_{\mathrm{geo}}$ fiberwise at
every point of $Y\al$, with respect to the morphism $Y\al \to X\al$. 

\item[($\text B^\ast $)] Let $Y\to X$ be an arrow belonging to $\mathscr C$, with $Y$ and $X$ good. 
Let $y$ be a point of $Y$
 and let $x$ be its image in $X$. Let $D$ be an object of $\mathfrak F_Y$.
 The object $D_y$ satisfies $\mathsf P_{\mathrm{geo}}$
fiberwise at every point of $\spec \mathscr O_{Y,y}$, with respect to the map $\spec \mathscr O_{Y,y}\to \spec \mathscr O_{X,x}$. 

\end{itemize}

Then if (A) is true, so are ($\text A^\ast$), (B) and ($\text B^\ast$). 
Indeed, we may perform the ``standard trick" described in 
\ref{ss-trick-fibval} in order to reduce assertion (B) to assertion (A), and
the one described in \ref{ss-trick-fibvalan}
in order to reduce assertion assertion ($\text B^\ast$) to assertion ($\text A^\ast$) (note that both tricks involve base change
by inclusions of affinoid domains and finite maps, which preserve $\mathscr C$ by assumption).

It thus remains to ensure that (A)$\Rightarrow$($\text A^\ast$). So let us assume that (A) holds, and let us prove ($\text A^\ast$). 
Let $\eta$ be a point of $\spec \mathscr O_{Y,y}$ lying over
the generic point of $X\al$. It follows from (A) that $D\al$ satisfies $\mathsf P$ at the image 
of $\eta$ on $Y\al$. Since $\mathsf P$
satisfies \hreg~
and since $\mathscr O_{Y\al, y\al}\to \mathscr O_{Y,y}$ is regular,
this implies that $D_y$ 
satisfies $\mathsf P$ at $\eta$, whence ($\text A^\ast$).

\section{Validity in analytic geometry}\label{s-appendix-analytic}

If $X$ is an
analytic space and if $D$ is an object of $\mathfrak F_X$, 
we have explained in Definition~\ref{def-valid-good} what it means for $D$ to satisfy 
$\mathsf P$ at a given point of $X$, under the assumption that 
$X$ is {\em good}. The purpose of what follows is to extend this definition to arbitrary analytic spaces, 
provided that $\mathsf P$ satisfies \hreg.

\begin{enonce}[lemm]{Lemma-Definition}\label{lem-equiv-valid}
Let $X$ be an analytic space, let $x$ be a point of $X$
and let $D$ be an object of $\mathfrak F_X$. Assume that $\mathsf P$
satisfies \hreg. The following are then equivalent: 

\begin{enumerate}[i]

\item For every good analytic domain $U$ of $X$ containing $x$, the object $D_U$ satisfies $\mathsf P$ at $x$; 

\item There exists a good analytic domain $U$ of $X$ containing $x$
such that $D_U$ satisfies $\mathsf P$ at $x$. 
\end{enumerate}
We say that $D$ {\em satisfies $\mathsf P$ at $x$}
if equivalent conditions {\rm(i)} and {\rm (ii)} are fulfilled. If $D$ satisfies
$\mathsf P$ at every point of $X$ and  if $\mathsf P$ satisfies \gen, we shall
say that $D$ {\em satisfies $\mathsf P$}
(\cf
Remark \ref{rem-valid-localization}).

\end{enonce}

\begin{proof}
Implication (i)$\Rightarrow$(ii) follows from the fact that there exists a good
analytic domain of $X$ containing $x$; \eg, an affinoid domain. Now assume that 
there exists $U$ as in (ii). Let $V$ be a good analytic domain of $X$ contaning $x$. 
Choose a good analytic domain $W$ of $U\cap V$ that contains $x$.
The morphisms $\mathscr O_{V,x}\to \mathscr O_{W,x}$ and $\mathscr O_{U,x}\to \mathscr O_{W,x}$ are
regular by~\ref{ss-algprop-analytic} (2). By choice of $U$, the object $D_{U,x}$ satisfies $\mathsf P$. Using twice the fact that
$\mathsf P$ satisfies \hreg, we see that $D_{W,x}$ satisfies $\mathsf P$ and that $D_{V,x}$ satisfies $\mathsf P$. 
\end{proof}

\begin{rema}\label{rem-equiv-valid}
We keep the notation of Lemma-Definition~\ref{lem-equiv-valid}.
If~$V$ is any analytic domain of~$X$ containing $x$, 
it follows from the definition that $D$ satisfies $\mathsf P$
at $x$ if and only if so does $D_V$. 
\end{rema}

In view of Examples \ref{ex-hreg-rings}--\ref{ex-hreg-diag}
(and Remark \ref{rem-norm-r1s2}), 
Lemma-Definition \ref{lem-equiv-valid}
above leads to the following more concrete statement. 

\begin{enonce}{Lemma-Definition (concrete
version of Lemma-Definition \ref{lem-equiv-valid})}\label{valid-at-concrete}
 Let $X$ be a $k$-analytic space,
let $\mathscr F$ be a coherent sheaf on $X$, let $\mathscr G\to \mathscr H$
be a morphism of
coherent sheaves on $X$,  and let
$\mathsf S$ be a complex of
coherent sheaves on $X$. Let $x$ be a point of $X$ and let $n$ be a non-negative integer. 

\begin{enumerate}[1]\index{analytic space!Gorenstein}\index{analytic space!regular}\index{analytic space!reduced}\index{Gorenstein analytic space}\index{regular analytic space}
\index{reduced analytic space}
\item The space $X$ is
said to be regular
(resp $R_n$, \resp Gorenstein, \resp CI, \resp normal, \resp reduced)\index{analytic space!CI (Complete Intersection)}\index{analytic space!normal}
at\index{analytic space!R@$R_n$}\index{normal analytic space}\index{CI (Complete Intersection) analytic space}
$x$ if 
there exists a good analytic domain $V$
of $X$
containing $x$ such that
the local ring
$\mathscr O_{V,x}$ is regular
(resp $R_n$, \resp Gorenstein, \resp CI, \resp normal, \resp reduced);
and if this is the case,
this holds for \emph{every} 
good analytic domain of $X$
containing $x$. 

\item The coherent sheaf
$\mathscr F$ is said to be\index{coherent sheaf!CM (Cohen-Macaulay)}\index{coherent sheaf!S@$S_n$}\index{CM (Cohen-Macaulay) coherent sheaf}
$S_n$ (\resp  CM, \resp of codepth $n$, \resp free of
rank $n$)
at $x$ if 
there exists a good analytic domain $V$
of $X$
containing $x$ such that the $\mathscr O_{V,x}$-module $\mathscr F_{V,x}$
is $S_n$ (\resp  CM, \resp of codepth $n$, \resp free of
rank $n$); and if this is the case,
this holds for \emph{every} 
good analytic domain of $X$
containing $x$. 

\item The morphism $\mathscr G\to \mathscr H$
is said to be
injective (\resp surjective, \resp bijective)
at $x$ if 
there exists a good analytic domain $V$
of $X$
containing
$x$ such that the $\mathscr O_{V,x}$-linear map
$\mathscr G_{V,x}\to \mathscr H_{V,x}$
is injective (\resp surjective, \resp bijective); 
and if this is the case,
this holds for \emph{every} 
good analytic domain of $X$
containing $x$. 

\item The complex $\mathsf S$ is
said to be exact at $x$ if 
there exists a good analytic domain $V$
of $X$
containing
$x$ such that 
the complex $\mathsf S_{V,x}$ of $\mathscr O_{V,x}$-modules
is exact; and if this is the case,
this holds for \emph{every} 
good analytic domain of $X$
containing $x$.

\end{enumerate}

\end{enonce}

\begin{rema}
According to our conventions, a space
$X$ will be called reduced it is reduced at every of its points, 
in the sense of Example \ref{valid-at-concrete}
above. One checks that this is consistent with the former
definition of reducedness (\ref{ss-defred-adhoc}). 
\end{rema}

\begin{lemm}[GAGA Principles]\label{ss-gaga-principles}
Let~$\mathscr X$
be scheme locally
of finite type over an affinoid algebra and let~$x$ be a point of~$\mathscr X\an$. 
Assume that $\mathsf P$ satisfies \hreg.
The following are equivalent for every object $D$ of $\mathfrak F_{\mathscr X}$:

\begin{enumerate}[i]
\item The object $D$ satisfies $\mathsf P$ at~$x\al$. 

\item The object $D\an$ satisfies $\mathsf P$ at $x$.

\end{enumerate}
\end{lemm}

\begin{proof}
Since $\mathsf P$ satisfies \hreg, this is an immediate consequence of the
regularity of the map $\mathscr O_{\mathscr X,x\al}\to \mathscr O_{\mathscr X\an,x}$. 
\end{proof}

In view of Examples \ref{ex-hreg-rings}--\ref{ex-hreg-diag}
(and Remark \ref{rem-norm-r1s2}), the GAGA principles stated above 
lead to the the following
more concrete statement. 

\begin{enonce}{Lemma (GAGA principles, concrete version of Lemma
\ref{ss-gaga-principles})}\label{gaga-concrete}
Let~$\mathscr X$
be scheme
of finite type over an affinoid algebra and let~$x$ be a point of~$\mathscr X\an$. 
Let $\mathscr F$ be a coherent sheaf on $\mathscr X$ , let $\mathscr G\to \mathscr H$
be a morphism
of coherent sheaves on $\mathscr X$, and let 
$\mathsf S$ be a complex of
coherent sheaves on $\mathscr X$. 

\begin{enumerate}[1]
\item The analytic space $\mathscr X\an$ is regular
(resp $R_n$, \resp Gorenstein, \resp CI, \resp normal, \resp reduced)
at
$x$ if and only 
the scheme $\mathscr X$ is
regular
(resp $R_n$, \resp Gorenstein, \resp CI, \resp normal, \resp reduced)
at
$x\al$. 

\item The coherent sheaf $\mathscr F\an $
is $S_n$ (\resp  CM, \resp of codepth $n$, \resp free of
rank $n$) at $x$ if and only if
the coherent sheaf $\mathscr F$
is $S_n$ (\resp  CM, \resp of codepth $n$, \resp free of
rank $n$) 
$x\al$. 

\item The morphism of coherent sheaves 
$\mathscr G\an\to \mathscr H\an$ is injective (\resp surjective, \resp bijective)
at $x$ if and only if
$\mathscr G\to \mathscr H$
is injective (\resp surjective, \resp bijective)
at $x\al$.  

\item The complex $\mathsf S\an$ is exact at $x$ if and only if 
the complex $\mathsf S$ is exact at $x\al$. 

\end{enumerate}

\end{enonce}

\begin{rema}
Both proofs of Lemma-Definition 
\ref{lem-equiv-valid}
and Lemma \ref{ss-gaga-principles}
rest in a crucial way on the regularity of some
local maps ($\mathscr O_{V,x}\to \mathscr O_{W,x}$ and $\mathscr O_{U,x}\to \mathscr O_{W,x}$
in the first one, $\mathscr O_{\mathscr X,x\al}\to \mathscr O_{\mathscr X\an,x}$ in the second one),
which is required to apply the axiomatic condition \hreg. 
But for some of the explicit
properties considered in Lemmas \ref{valid-at-concrete} 
and \ref{gaga-concrete}, it would be enough to know that the local maps involved
are flat: this is the case for the property of being free of rank $n$ (for a coherent sheaf),
of being injective or bijective (for a map of coherent sheaves), 
or of being exact (for a complex of coherent sheaves); note that for surjectivity
of a map of coherent sheaves, flatness is even not necessary: 
we would be done
via Nakayama's Lemma.
\end{rema}

\begin{lemm}\label{lem-loci-abstract}
Let $X$ be a $k$-analytic space and let $D$
be an object of $\mathscr F_X$. Assume
that $\mathsf P$ satisfies \open~and \hreg. 
Let $U$ be the set of points of $X$ at which $D$
satisfies $\mathsf P$. 

\begin{enumerate}[1]

\item The set $U$ is Zariski-open in $X$.

\item Assume moreover that $\mathsf P$ satisfies \field~and
the space $X$ is reduced and non-empty (\eg, $X$ is integral). 
Then $U\neq \emptyset$.

\end{enumerate}
\end{lemm}

\begin{proof}
Both assertions
are G-local, hence we can assume
that $X$ is affinoid. 
Since $\mathsf P$ satisfies \open, the subset $\mathscr U$
of $X\al$ consisting of points at which $D\al$ satisfies $\mathsf P$ is Zariski-open. 
Since $\mathsf P$ satisfies \hreg, Lemma \ref{ss-gaga-principles}
ensures that $U=\mathscr U\an$, whence (1). 

Assume now that $\mathsf P$ satisfies
\field~
and $X$ is reduced and non-empty. In this case $X\al$ is reduced too, 
hence there exists a point $\xi$ on $X\al$
such that $\mathscr O_{X\al,\xi}$ is a field
(take the generic point of any irreducible component). 
As a consequence $\mathscr U\neq \emptyset$, 
so $U=\mathscr U\an$ is non-empty. 
\end{proof}

In view of Examples
\ref{ex-algprop-space}--\ref{ex-algprop-complex}
(and Remark \ref{rem-norm-r1s2}), Lemma \ref{lem-loci-abstract}
leads to the more concrete statement. 

\begin{enonce}{Lemma (concrete version of Lemma \ref{lem-loci-abstract})}\label{locus-concrete}
Let $X$ be a $k$-analytic space,
let $\mathscr F$ be a coherent sheaf on $X$, let $\mathscr G\to \mathscr H$
be a morphism of
coherent sheaves on $X$,  and let
$\mathsf S$ be a complex of
coherent sheaves on $X$. Let $x$ be a point of $X$ and let $n$ be a non-negative integer. 

\begin{enumerate}[1]
\item The subset of $X$ consisting of points at which 
$X$ is regular
(resp $R_n$, \resp Gorenstein, \resp CI, \resp normal, \resp reduced)
is Zariski-open, and non-empty if $X$ is non-empty and reduced. 

\item The subset of $X$ consisting of points
at which the coherent sheaf
$\mathscr F$ is $S_n$ (\resp  CM, \resp of codepth $\leq n$, \resp free) is  
Zariski-open, and non-empty if $X$ is non-empty and reduced.

\item The subset of $X$ consisting of points
at which the morphism $\mathscr G\to \mathscr H$
is injective (\resp surjective, \resp bijective)
is Zariski-open. 

\item The subset of $X$ consisting of points
at which the complex $\mathsf S$ is exact
is Zariski-open. 
\end{enumerate}
\end{enonce}

\section{Fibers of coherent sheaves}\label{s-fibers-coherent}\index{fiber (of a coherent sheaf)}

\begin{defi}\label{defi-fhx-rkx}
Let $X$ be an analytic space and let $\mathscr F$ be a coherent sheaf on $X$.
Let $x$ be a point of $X$ and let $V$ be a good analytic domain containing $x$. 
The tensor product $\hr x\otimes_{\kappa(x_V)}\mathscr F_{V,x}$
does not depend on $V$, and will be denoted by $\mathscr F_{\hr x}$.
This is a finite dimensional $\hr x$-vector space which is called
the \emph{fiber} of $\mathscr F$ at $x$. Its dimension is called the\label{IN-Fhrx}
\emph{fiber rank}
of $\mathscr F$ at $x$ and is denoted by $\rk x(\mathscr F)$. \label{IN-rkxF}
\end{defi}

\subsection{Basic properties}\label{ss-pointwise-rk}
Let $X$ be an analytic space
and let $\mathscr F$ be a coherent sheaf on $X$.
The following properties 
show that the fiber rank behaves
like the usual fiber rank
on a noetherian scheme (or on a locally ringed space with coherent structure sheaf).

\begin{enumerate}[1]
\item For every good analytic domain $V$ of $X$ and every point $x$
of $V$, we have the equality 
$\mathscr F_{\hr x}=\hr x\otimes_{\kappa(x_V)} \mathscr F_{\kappa(x_V)}$. 

\item For every affinoid domain $V$ of $X$ and every point $x$ of $V$, we have the equalities

\begin{eqnarray*}
\mathscr F_{\hr x}
&=&\hr x\otimes_{\mathscr O_{V,x}}\mathscr F_{x_V}\\
&=&\hr x \otimes_{\mathscr O_X(V)}\mathscr F(V)\\
&=&\hr x\otimes_{\mathscr O_{V\al}(V\al)}\mathscr F_V\al(V\al)\\
&=&\hr x\otimes_{\kappa(x_V\al)}(\mathscr F_V\al)_{\kappa(x_V\al)}.\\
\end{eqnarray*} 

\item For every point $x$ of $X$, the functor $\mathscr F\mapsto \mathscr F_{\hr x}$ 
is right-exact. 

\item For every coherent sheaf $\mathscr F$, the function $x\mapsto \rk x \mathscr F $ is upper semi-continuous
for the {\em Zariski}
topology of $X$. 

\item For every point $x$ of $X$ and every coherent sheaf $\mathscr F$ on $X$, the following are equivalent:

\begin{itemize}
\item [(i)]The vector space $\mathscr F_{\hr x}$ is zero. 

\item [(ii)]  The coherent sheaf $\mathscr F$ is zero at $x$
(in the sense of Lemma-Definition \ref{lem-equiv-valid}; it means that
$\mathscr F$ is ``free of rank $0$" at $x$ in the sense of the more concrete
Lemma \ref{valid-at-concrete}). 

\item [(iii)] There exists an open neighborhood $U$ of $x$ such that $\mathscr F_U=0$.
\end{itemize}

\end{enumerate}

Indeed, (1) and (2) follow immediately from the definitions
and the fact that $\mathscr F_{x_V}=\mathscr O_{V,x}\otimes_{\mathscr O_X(V)}\mathscr F(V)$
(\ref{ss-stalks-fibers}). 
Assertion (3) is a consequence of (1),  and (4) is a consequence of (2)
and of upper semi-continuity of the pointwise rank in the 
locally ringed space context (which itself comes from Nakayama's Lemma).
Concerning (5), let us first assume that $\mathscr F_{\hr x}=0$.
Then for every good analytic domain $V$ of $X$ containing $x$ one has $\mathscr F_{\kappa(x_V)}=0$ in view of (1), 
and hence $\mathscr F_{V,x}=0$ by Nakayama's Lemma,
whence (ii).
Assume now that (ii) holds, and let $V$ be an affinoid domain of $X$ containing $x$.
By assumption, one has $\mathscr F_{V,x}=0$, which implies that $\mathscr F_U=0$
for some open neighborhhod $U$ of $x$ {\em in $V$}; now (iii) follows from the fact that $x$ has a neighborhood {\em in $X$}
that is the union of finitely many affinoid domains containing $x$. The
implication (iii)$\Rightarrow$(i) is obvious. 

\subsection{The support of a coherent sheaf}\label{ss-interpret-support}\index{support (of a coherent sheaf)}
Let $\mathscr F$ be a coherent sheaf on $X$ and let $\mathscr I$ be the
(coherent) annihilator ideal of $\mathscr F$ (on the site $X\grot$).
The {\em support} of $\mathscr F$ is the closed analytic subspace of $X$ defined
by $\mathscr I$; it is denoted by $\supp F$. If $i\colon \supp F\hookrightarrow X$
is the canonical closed immersion, then $\mathscr F=i_\ast i^\ast \mathscr F=i_\ast i\inv \mathscr F$. 

Let $W$ be the subset of $X$ consisting of points $x$ such that $\mathscr F_{\hr x}=0$. By assertion (4) of \ref{ss-pointwise-rk}
above, $W$ is a Zariski-open subset of $X$, and $\mathscr F_W=0$ by assertion (5). 
The complement
$X\setminus W$ is then equal to (the Zariski-closed subset underlying) $\supp  F$.
Indeed, by arguing G-locally on $X$ we may assume that it is good. Let $x$ be a point of $X$. 
By the general theory
of locally ringed spaces with coherent structure sheaf
$$(\mathscr F_x\neq 0)\iff (\mathscr I_x\subset \mathfrak m_x)\iff(f(x)=0\;\text{for all}\;f\in \mathscr I_x),$$
and the latter condition exactly means that $x$ lies in $\supp F$. 

In accordance with the usual convention in commutative algebra
(\ref{ss-conv-alg}), the \emph{dimension}\index{dimension!of a coherent sheaf}
of $\mathscr F$ (\resp the \emph{dimension}
of $\mathscr F$ \emph{at $x$} for $x$ a given point of $X$) is
by definition the dimension of $\supp F$ (\resp the dimension of $\supp F$ at $x$
if $x$ lies in $\supp F$, and $-\infty$ otherwise). We denote it by $\dim \mathscr F$
(\resp $\dim_x \mathscr F$).

\subsection{Surjectivity can be checked fiberwise}\label{ss-surj-nakayama}
Let $\mathscr F\to \mathscr G$ be a morphism between two coherent sheaves on $X$; let $x$ be a point of $X$.
Let $\mathscr Q$ be the cokernel
of $\mathscr F\to \mathscr G$. Since formation of fibers is a right-exact functor by
\ref{ss-pointwise-rk} (3), the sequence
$$\mathscr F_{\hr x}\to \mathscr G_{\hr x}\to \mathscr Q_{\hr x}\to 0$$ is exact. By assertion (5) of
\loccit, 
the vector space $\mathscr Q_{\hr x}$ is zero if and only if $\mathscr Q$ is zero at $x$, that is, if and only if
$\mathscr F\to \mathscr G$ is surjective at $x$. In other words, $\mathscr F\to \mathscr G$ is surjective at $x$
if and only if $\mathscr F_{\hr x}\to \mathscr G_{\hr x}$ is surjective.

\subsection{Consistency with the usual sheaf-theoretic notions}\label{ss-pointwise-sheaf}
According to our general conventions, a morphism of coherent sheaves on $X$ should be called
injective, \resp surjective, \resp bijective
if and only if satisfies this property at every point of $X$. But these notions also make sense
in the general sheaf-theoretic context; our purpose is now to ensure
that both terminologies are compatible. 

Let $\mathscr F\to \mathscr G$ be a morphism between two coherent sheaves on $X$. Let us denote its kernel by $\mathscr K$, and its cokernel 
by $\mathscr Q$. Let $U$ be the complement of $\supp K$ in $X$, and let $V$ be that of $\supp Q$. 
It follows from \ref{ss-interpret-support} that $U$ is the set of points of $X$ at which 
$\mathscr K=0$, that $V$ is the set of points of $X$ at which 
$\mathscr Q=0$, and that $\mathscr K_U=0$ and $\mathscr Q_V=0$. 
This implies that 
$U$ is the set of points at which $\mathscr F\to \mathscr G$ is 
injective,
that $V$ is the set of points at which $\mathscr F\to \mathscr Q$ is surjective,
that $\mathscr F_U\to \mathscr G_U$
is sheaf-theoretically
injective and that $\mathscr F_V\to \mathscr G_V$ is sheaf-theoretically
surjective. 
It follows that $U\cap V$ is the set of points at which $\mathscr F\to \mathscr G$
is bijective, and that $\mathscr F_{U\cap V}\to \mathscr G_{U\cap V}$ is sheaf-theoretically
bijective. As a consequence, our terminology is compatible
with that from sheaf theory.

\section{Ground field extension}\label{s-ground-field-ext}

Our aim is now to state some
properties of the ground field extension functor in analytic geometry, beyond
faithful flatness.
In {\em algebraic} geometry, some properties
of schemes of finite type over a field (like being CM, CI, Gorenstein, $S_m$ for some specified $m$) behave well
under {\em any} ground field extension; but some others (like being regular, or $R_m$ for some specified $m$)
are only preserved in full generality by {\em separable} extensions of the ground field. 

An analogous
phenomenon holds in analytic geometry. In order to describe it, we first need 
to explain what the analytic analogue
of a separable extension is.  

\begin{defi}[after \cite{ducros2009}]\label{def-an-sep}
Let $k$ be an analytic field. An analytic extension $L$ of $k$ is called {\em analytically separable} if char.\@~$k=0$\index{extension!analytically separable}\index{analytically separable
extension}
or if char.\@~$k=p>0$ and the semi-norm
$$L\hotimes_k k^{1/p}\to \R_+,\;\sum \ell_i \otimes x_i\mapsto\underbrace{\abs{\sum \ell_i x_i}}_{\in L^{1/p}}$$
is a norm equivalent to the tensor norm of $L\hotimes_k k^{1/p}$. 
\end{defi}

Let us give some examples (for detailed proofs, see \cite{ducros2009}, \S 1.2, Lemme 1.8 and Exemple 1.9). 

\begin{exem} If $k$ is a perfect analytic field, every analytic extension of $k$ is analytically separable.
\end{exem}

\begin{exem} A finite extension of an analytic field is analytically separable if and only if it is separable. 
\end{exem}

\begin{exem}\label{ex-etar-ansep}
Let $k$ be an analytic field and let $r$ be a polyradius. 
The analytic extension $k\hookrightarrow \hr {\eta_r}$ (\ref{ss-intro-etar}) is analytically separable (recall that $\hr {\eta_r}=k_r$ when
$r$ is $k$-free). 
\end{exem}

\begin{theo}\label{th-fiber-xl}
Let~$X=(k,X)$ be a good analytic space
and let~$L$ be an analytic extension of
$k$.
Let $x$ be a point of $X$ and let $y$ be a point of $X_L$
lying over $x$. 

\begin{enumerate}[1]

\item  The fibers of the faithfully flat morphism
$\spec \mathscr O_{X_L,y}\to \spec \mathscr O_{X,x}$ 
are CI, and geometrically regular when $L$ is analytically separable 
over $k$. 

\item
If $X$ is affinoid, the fibers of the faithfully flat morphism 
$(X_L)\al\to X\al$ are CI, and geometrically regular when
$L$ is analytically separable 
over $k$.

\end{enumerate}
\end{theo}

\begin{proof}
Assume that $X$ is affinoid and integral. 
By \Prop
2.2 (a) of \cite{ducros2009},
there exists a non-empty open subset~$\mathscr U$ of $X\al$ such that~$\mathscr U\times_{X\al}(X_L)\al$ is CI, and regular
when $L$ is analytically separable over $k$. 
The theorem then follows by applying~\ref{ss-fiberwise-chiant}
with $\mathfrak F=\mathfrak T$, with $\mathsf P$ being the CI property
in general,  and the regularity
property when $L$ is analytically separable,
and with $\mathscr C$ being
the class of morphisms of the form $Y_L\to Y$ for $Y$ a $k$-analytic space.
\end{proof}

By the definition of \hci~and \hreg, the above theorem implies the following proposition. 

\begin{prop}\label{prop-behavior-extension}
Let~$k$ be an analytic field, let $X$
be a $k$-analytic space, let $D$ be an object of $\mathfrak F_X$,
and let $L$ be an analytic extension
of $k$. Let $x$ be a point of $X$ and let $y$ be a pre-image
of $x$ in $Y$. 
Assume that $\mathsf P$ satisfies \hreg. 

\begin{enumerate}[1]
\item If~$D_L$ satisfies~$\mathsf P$ at~$y$ then~$D$ satisfies~$\mathsf P$ at~$x$; 

\item If~$D$ satisfies~$\mathsf P$ at~$x$, then  
$D_L$
satisfies~$\mathsf P$ at~$y$
if either $L$
is analytically
separable over~$k$ or $\mathsf P$ satisfies~\hci. 
\end{enumerate}

\end{prop}

In view of Examples
\ref{ex-algprop-space}--\ref{ex-algprop-complex}
(and Remark \ref{rem-norm-r1s2}),
Proposition \ref{prop-behavior-extension}
leads to a more concrete statement:

\begin{prop}[Concrete version of Proposition \ref{prop-behavior-extension}]\label{prop-concrete-extension}
Let~$k$ be an analytic field, let $X$ be a $k$-analytic space, 
let $\mathscr F$ be a coherent sheaf on $X$, let $\mathscr G\to \mathscr H$
be a morphism of
coherent sheaves on $X$,  and let
$\mathsf S$ be a complex of
coherent sheaves on $X$. Let $x$ be a point of $X$ and let $n$ be a non-negative integer. 
Let $L$ be an analytic extension of $k$ and let $y$ be a pre-image
of $x$ on $X_L$. 

\begin{enumerate}[1]
\item The space $X_L$ is Gorenstein (\resp CI) at $y$
if and only if $X$ is Gorenstein (\resp CI) at $x$. 

\item If $X_L$ is regular (\resp $R_n$, \resp normal, \resp reduced)
at $y$, then $X$ is regular (\resp $R_n$, \resp normal, \resp reduced)
at $x$. 

\item If $X$ is regular (\resp $R_n$, \resp normal, \resp reduced) at $x$ and $L$ is analytically 
separable over $k$, then $X_L$  is regular (\resp $R_n$, \resp normal, \resp reduced)
at $y$. 

\item The coherent sheaf
$\mathscr F_L$ is 
$S_n$ (\resp  CM, \resp of codepth $n$, \resp free of
rank $n$) at $y$ if and only if $\mathscr F$ is 
$S_n$ (\resp  CM, \resp of codepth $n$, \resp free of
rank $n$) at $x$.

\item The morphism
of coherent sheaves $\mathscr G_L\to \mathscr H_L$
is
injective (\resp surjective, \resp bijective)
at $y$ if and only if
$\mathscr G\to \mathscr H$
is
injective (\resp surjective, \resp bijective)
at $x$.

\item The complex $\mathsf S_L$ is
exact at $y$ if and only if 
$\mathsf S$ is exact at $x$.
\end{enumerate}
\end{prop}

\begin{rema}
Proposition \ref{prop-behavior-extension}
rests on Theorem
\ref{th-fiber-xl}
which is needed for using
the axiomatic conditions \hreg~or
\hci. 
But for some of the explicit
properties considered in Proposition \ref{prop-concrete-extension}, 
it would be enough to know that the local maps involved
are flat: this is the case for the property of being free of rank $n$ (for a coherent sheaf),
of being injective or bijective (for a map of coherent sheaves), 
or of being exact (for a complex of coherent sheaves); note that for surjectivity
of a map of coherent sheaves, flatness is even not necessary: 
we would be done via Nakayama's Lemma.
\end{rema}

\subsection{Geometric validity}\label{ss-geom-valid}
 Let~$k$ be an analytic field, let $X$
 be a $k$-analytic space,
and let~$x$ be a point of $X$.
Let~$D$ be an object of $\mathfrak F_X$. We shall
say that~$D$ satisfies~$\mathsf P$ {\em geometrically}
at~$x$ if for every analytic extension~$L$ of~$k$ and every point $y$
of $X_L$ lying above $X$, the object~$D_L$
satisfies~$\mathsf P$ at~$x_L$.

Assume
that $\mathsf P$
satisfies \hreg.
In order for~$D$ to satisfy geometrically~$\mathsf P$ at~$x$, it is sufficient that there exists a {\em perfect}
analytic extension~$L$ of~$k$ and a  pre-image~$y$ of~$x$ on~$X_L$
such that the object~$D_L$
satisfies~$\mathsf P$ at~$y$. Indeed, assume that it is the
case, let~$F$ be an analytic extension of~$k$ and let~$z$ be a pre-image of~$x$ on~$X_F$. 
By Lemma~\ref{lem-pointlying}
below, there exists an analytic extension~$K$ of~$k$, equipped
with two isometric~$k$-embeddings $L\hookrightarrow K$ and~$F\hookrightarrow K$, 
and a point~$t\in X_K$ lying above both $y$ and~$z$. Since~$L$ is perfect,~$K$
is analytically separable over~$L$. It thus
follows from
Proposition~\ref{prop-behavior-extension}
above that~$D_K$ satisfies~$\mathsf P$ at~$t$, and 
then that~$D_F$ satisfies~$\mathsf P$ at~$z$.

\begin{lemm}\label{lem-pointlying}
Let ~$k$ be an analytic field  and let~$L$ and~$F$ be two analytic extensions of~$k$. Let~$Y\to X$ be a morphism of~$k$-analytic spaces.
Let~$X'$ be an~$F$-analytic space, let~$X'\to X$ be  a morphism of analytic spaces and set
$$Y'=Y\times_XX'=Y_F\times_{X_F}X'.$$ Let~$y$ be a point of~$Y$. Let~$u$ (resp.~$y'$)
be a point of~$Y_L$ (resp.~$Y'$) lying above~$y$. There exist a complete extension~$K$ of~$k$, equipped with
two isometric~$k$-embeddings~$F\hookrightarrow K$ and 
$L\hookrightarrow K$, and a point~$\omega$ of~$$Y'_K:= Y'\times_FK\simeq Y\times_XX'_K\simeq Y_K\times_{X_K}X'_K\simeq  Y_L\times_{X_L}X'_K\simeq Y_F\times_{X_F}X'_K$$
lying above both~$y'$ and~$u$.
\end{lemm}

\begin{proof}
We immediately reduce to the case where~$Y$ and
$X$ are $k$-affinoid and~$X'$
is~$F$-affinoid. Let~$A, A',B$
and~$B'$ be the respective algebras of analytic functions on~$X,X', Y$ and $Y'$.
The points~$u$ and~$y'$ furnish a
pair of characters~$$(B_L\to \hr u, B'\to \hr {y'}),$$ the restriction of
which to~$B$ go through~$B\to \hr y$. The Banach algebra~$\hr u \hotimes_{\hr y}\hr {y'}$ is non-zero: a result by Gruson ensures that it contains~$\hr u\otimes_{\hr y}\hr {y'}$, \cite{gruson1966}, \S 3.2, \Th 1 (4). There exists therefore an analytic field~$K$ and a bounded homomorphism~$\hr u\hotimes_{\hr y}\hr {y'}\to K$ (\cite{berkovich1990}, \Th 1.2.1). 
This makes~$K$ an analytic extension of both~$\hr u$ and~$\hr {y'}$ over $\hr y$. 
One thus gets a new pair of characters~$$(B_L\to K, B'\to K)$$ whose restriction to~$B$ coincide; that
pair induces tautologically a character from $B_L\hotimes_B B'=B_L\hotimes_A A'$ to  $K$, which extends canonically to a character
$$B_K\hotimes_{A_K}(K\hotimes_FA')\to K.$$ The latter defines a point~$\omega$ on~$Y_K\times_{X_K}X'_K$ lying by construction above both~$y'$ and~$u$.
\end{proof}

\begin{rema}\label{rem-hci-geo}
By
Proposition
\ref{prop-behavior-extension}
above, if~$\mathsf P$ satisfies~\hci~then
its validity at
a point is equivalent to
its geometric validity at the latter. Therefore in practice, the notion
of geometric validity will be of specific interest only when~$\mathfrak F=\mathfrak T$
and when the property we are considering involves 
regularity or~$R_m$ for some specified~$m$; \eg, geometric regularity, geometric
reducedness, and geometric normality. 
\end{rema}

\section{Complements on analytifications}\label{s-complement-xan}

Our next goal is to extend some results that hold for
(spectra of) affinoid algebras to schemes locally of finite type over such spectra. We
fix an analytic field $k$ and a $k$-affinoid algebra $A$; we
first consider the strict situation.

\begin{lemm}\label{lem-closed-point}
Assume that $\abs{k\gpm}\neq \{1\}$ and that $A$ is strict. Let $\mathscr X$ be an $A$-scheme
locally of finite type and let $x$ be a point of $\mathscr X$. The point $x$
is closed if and only if $\kappa(x)$ is a finite extension of $k$.
\end{lemm}

\begin{rema}
The result is well-known for $\mathscr X=\spec A$: this is nothing but the classical analytic {\em Nullstellensatz}; our proof
essentially consists if reducing to this situation. 

\end{rema}

\begin{proof}[Proof of Lemma \ref{lem-closed-point}]
If $\kappa(x)$ is finite over $k$, then $\spec \kappa(x)\to \mathscr X$ is finite, and $x$ is closed. Conversely,
let us assume that $x$ is closed, and let $\xi$ be its image
in $\spec A$. Since $x$ is closed in its fiber $\mathscr X_\xi$ which is locally of finite
type over $\kappa(\xi)$, 
it suffices to prove that $\kappa(\xi)$ is a finite extension of $k$.
By Chevalley's theorem, $\{\xi\}$ is a constructible subset of $\spec A$, which means that $\xi$
is open in its Zariski closure $\mathscr Y$. Set $\mathscr Z=\mathscr Y\setminus \{\xi\}$; this is 
a proper Zariski-closed subset of $\mathscr Y$. 

Let us endow $\mathscr Y$ with its reduced structure, and let $f$ be a non-zero function on $\mathscr Y$ 
that vanishes pointwise on $\mathscr Z$; let $r$
be the spectral norm of $f$, seen as a function on the strict affinoid space $\mathscr Y\an$. The strict affinoid domain 
of $\mathscr Y\an$ defined by the condition $\abs f=r$ is non-empty, hence admits a rigid point $z$. By construction, $f(z\al)\neq 0$, which implies that 
$z\al=\xi$. Therefore $\kappa(\xi)$ is a finite extension of $k$. 
\end{proof}

\begin{lemm}\label{lem-catenarity}
Assume that $\abs{k\gpm}\neq \{1\}$ and that $A$ is strict. Let $\mathscr X$ be an irreducible $A$-scheme
locally of finite type 
and let $\mathscr Y$ be an irreducible closed
subset of $\mathscr X$. 

\begin{enumerate}[1]

\item The scheme $\mathscr X$ is finite-dimensional, and we have the equality
\[\dim \mathscr Y+\codim (\mathscr Y,\mathscr X)=\dim \mathscr X.\]

\item Every closed point of $\mathscr X$ is of codimension $\dim \mathscr X$ in $\mathscr X$.

\end{enumerate}

\end{lemm}

\begin{rema}\label{rem-catenarity}
The result is already known for $\mathscr X=\spec A$: see \cite {ducros2007}, \Prop 1.11, and more precisely \S 1.11.1 in its proof, in which the strict case is handled. Our proof
essentially consists in reducing to this situation. 

\end{rema}
\begin{proof}
Note that assertion (2) is a particular case of (1); we have written it down because we wanted to
emphasize it. But since $A$ is excellent, hence universally catenary, both statements are actually equivalent, and we shall in fact prove (2). 

Let us first assume that $\mathscr X$ is affine. By the analytic version of Noether's normalization lemma there exist $n\in \N$ and a finite, dominant morphism from $\spec A$
to the scheme $\mathscr Z:=\spec k\{T_1,\ldots, T_n\}$. Let us choose a factorization of the composite
morphism $\mathscr X\to \mathscr Z$ through a closed immersion $\mathscr X\hookrightarrow \A^m_{\mathscr Z}$ (for some $m$). 
Let $x$ be a closed point of $\mathscr X$.
By Lemma \ref{lem-closed-point}
above, $\kappa(x)$ is a finite extension of $k$; if $z$ denotes the image of $x$ on $\mathscr Z$, then $\kappa(z)$
is also a finite extension of $k$, hence $z$ is closed in $\mathscr Z$. By Remark \ref{rem-catenarity}
above, $\{z\}$ is of codimension $n$ in $\mathscr Z$. Being a closed point of the  fiber $\A^m_z$, the point $x$ is
of codimension $m$ in it. By flatness of $\A^m_{\mathscr Z}\to \mathscr Z$, the point $x$ is then of codimension $m+n$ in $\A^m_{\mathscr Z}$. By catenarity of the latter
scheme, $x$ is of codimension $d:=n+m-\codim(\mathscr X,\A^m_{\mathscr Z})$ in $\mathscr X$. Since this holds
for every closed point of $\mathscr X$, the integer $d$ coincides with $\dim \mathscr X$, which ends the proof when $\mathscr X$ is affine. 

For the general case, let us chose a covering $(\mathscr X_i)$ of $\mathscr X$ by non-empty affine open subsets.
If
$i$ and $j$ are two indices, the intersection $\mathscr X_i\cap \mathscr X_j$ is non-empty, hence has a closed point $x$. By Lemma \ref{lem-closed-point}, $\kappa(x)$ is finite
over $k$, and $x$ is therefore closed in both $\mathscr X_i$ and $\mathscr X_j$. By the affine case already proven
we have $$\dim \mathscr X_i=\dim \mathscr X_j=\dim_{\mathrm{Krull}}(\mathscr O_{\mathscr X,x}).$$
Hence the $\mathscr X_i$'s all have the same dimension $d$. We then have $\dim \mathscr X=d$. If $x$ is a closed point of $\mathscr X$, it belongs to $\mathscr X_i$ for some $i$, whence (again by the affine case)
we get the equality
$\dim_{\mathrm{Krull}}(\mathscr O_{\mathscr X,x})=\dim \mathscr X_i=d$. 
\end{proof}

\begin{rema}
We still assume that $k$ is non-trivially valued
and $A$ is strict. Let $\mathscr X$ be an $A$-scheme
locally of finite type, let $\mathscr U$ be a locally closed subscheme of $\mathscr X$ and let $x$ be a point of $\mathscr U$.
It follows from Lemma \ref{lem-closed-point} that $x$ is closed in $\mathscr U$ if and only if it is closed in $\mathscr X$, because both assertions are equivalent to the finiteness of $\kappa(x)$ over $k$. 

Assume moreover that $\mathscr X$ is irreducible, and let 
$\mathscr U$ be a non-empty open subscheme of $\mathscr X$. 
We then have $\dim \mathscr U=\dim \mathscr X$. Indeed, choose a closed point $x\in \mathscr U$. 
By the above and by Lemma \ref{lem-catenarity}, we then have 
\[\dim \mathscr U=\dim_{\mathrm{Krull}}\mathscr O_{\mathscr X,x}=\dim \mathscr X.\] 
\end{rema} 

\begin{lemm}\label{lem-gaga-dimension}
Assume that $\abs{k\gpm}\neq \{1\}$ and that $A$ is strict. Let $\mathscr X$ be
an $A$-scheme locally of finite type 
and let $x$ be a point of $\mathscr X\an$. 
One has the equality
$$\dim_x \mathscr X\an=\dim_{x\al}\mathscr X.$$
\end{lemm}

\begin{proof}
Since both sides of the equality are local on $\mathscr X$, we can assume (by replacing $\mathscr X$ with a projective
compactification of any affine neighborhood of $x$) that $\mathscr X$ is projective. Let $(\mathscr X_i)$ be the family of irreducible
components of $\mathscr X$ that contain $x$. By GAGA
(see \ref{ss-def-analytify}), the $\mathscr X_i\an$ are the irreducible components of $\mathscr X\an$ that contain $x$ (this is in fact
true without the projectivity assumption, but that is more involved; see \Prop \ref{prop-gaga-irrcomp}
below). Therefore it suffices to prove the required equality for every $\mathscr X_i$; we thus can assume that $\mathscr X$
is irreducible, and in particular purely of dimension $d$ for some $d$. It suffices now to prove that $\mathscr X\an$ is purely $d$-dimensional too. 

Let $V$ be a non-empty strict affinoid
domain of $\mathscr X\an$. Let $y$ be a rigid point of $V$. We have the equalities 
$$\widehat{\mathscr O_{V\al, y_V\al}}=\widehat{\mathscr O_{V,y}}=\widehat{\mathscr O_{\mathscr X\an,y}}=\widehat{\mathscr O_{\mathscr X, y\al}}$$ (the first one and the third one come from Lemma 6.3 of \cite{berkovich1993}, and the middle one from the fact that $y$ belongs to $\mathrm{Int}(V/\mathscr X\an)$ since $y$ is rigid). Since a local noetherian ring and its completion have the same Krull dimensions, one has $$\dim_{\mathrm{Krull}}\mathscr O_{V\al, y_V\al}=\dim_{\mathrm{Krull}}\mathscr O_{\mathscr X, y\al}=d,$$ where the last equality comes from Lemma~\ref{lem-catenarity} (2). Since the closed points of the scheme $V\al$ are exactly the points of the form $y_V\al$ for $y$ a rigid point of $V$, the dimension of the scheme $V\al$ is equal to $d$. Hence $\dim V=d$ and $\mathscr X\an$ is purely $d$-dimensional. 
\end{proof}

\begin{prop}\label{prop-gaga-dim}
(The valuation on $k$ is no longer assumed to be non-trivial, nor the algebra $A$ to be strictly affinoid). 

\begin{enumerate}[1]
\item Let $\mathscr Y\to \mathscr X$ be a morphism of $A$-schemes locally of finite type, 
and let $y$ be a point of $\mathscr Y\an$. 
The relative dimension of $\mathscr Y\an \to \mathscr X\an$ at $y$ is equal to the relative dimension of $\mathscr Y\to \mathscr X$ at
$y\al$.

\item Let $\mathscr Z$ be a $k$-scheme locally of finite type. For every point $z$ of $\mathscr Z\an$ one has
$\dim_z \mathscr Z\an=\dim_{z\al}\mathscr Z.$
\end{enumerate}
\end{prop}

\begin{proof}Note that assertion (2) is a particular case of (1); we have written it down because we wanted to emphasize it.
But we shall in fact first prove (2), 
and deduce (1). 

Due to Lemma \ref{lem-gaga-dimension}, assertion (2) holds as soon as $\abs {k\gpm}\neq \{1\}$. But since both sides of the equality are 
preserved by ground field extension (\ref{dim-goundfield-extension}) we can always reduce to this case. 

Let us now prove (1). If $x$ denotes the image of $y$ on $\mathscr X\an$, the fiber $\mathscr Y\an_x$ is naturally isomorphic to 
the analytification $(\mathscr Y_{x\al}\times_{\kappa(x\al)}\hr x)\an$; hence (1) follows from (2) (and from the invariance of the local dimension 
on a scheme locally of finite type over a field under ground field extension). 
\end{proof}

Our purpose is now to extend Lemma \ref{lem-catenarity} to the non-strict case. In order to do that, we 
must understand what happens when one performs a base change from $k$ to $k_r$ for some $k$-free polyradius $r=(r_1,\ldots, r_n)$. 
The key point will be the following lemma. 

\begin{lemm}\label{lem-xxr-scheme}
Let $\mathscr X$ is a reduced (\resp irreducible)
$A$-scheme locally
of finite type and let $r$ be a $k$-free polyradius. 
The scheme $\mathscr X\times_A A_r$ is also reduced (\resp irreducible).
\end{lemm}

\begin{rema}
When $\mathscr X=\spec A$ this is nothing but 
Lemma 1.3 of \cite{ducros2007}, 
and we have simply adapted its proof to our more general setting. 
Let us also mention that we shall not need the result about reducedness
for extending Lemma  \ref{lem-catenarity}, and that it also follows from 
Lemma \ref{gaga-concrete} and Proposition \ref{prop-concrete-extension}; 
but we have chosen to include it in this lemma because our elementary 
method for dealing
with irreducibility provides it almost for free.
\end{rema}

\begin{proof}[Proof of Lemma \ref{lem-xxr-scheme}]
Let us first suppose that
$\mathscr X$ is affine, say $\mathscr X=\spec B$ for some $A$-algebra $B$ of finite type.
Set $T=(T_1,\ldots, T_n)$. Any function belonging to $B\otimes_A A_r$ can be written
in a unique way as an infinite sum $\sum_{I\in \Z^n} b_I T^I$, with $b_I\in B$ for all $I$,
that converges on every affinoid domain of $(\mathscr X\an)_r=(\mathscr X\times_A A_r)\an$. For such a function $f$, 
we denote by $\mathscr Z(f)$ the Zariski-closed subset of $\mathscr X$ defined
by the ideal $(b_I)_I$. 

Let us now make a general remark. 
Let $f=\sum b_I T^I$ and $g=\sum \beta_I T^I$ be two elements of $B\otimes_A A_r$
with $fg=0$. 
Let $x$ be a point of $\mathscr X\an$. The fiber $Y$ of the base-change map
$(\mathscr X \an)_r \to \mathscr X\an$
at $x$ is isomorphic to $ \mathscr M(\hr x)_r$. Since $fg=0$, we have $(f|_Y)(g|_Y)=0$. 
As $\hr x_r$ is a domain (its norm is multiplicative, \cf \cite{ducros2007}, 1.2.1), $f|_Y=0$ or $g|_Y=0$; 
otherwise said, we have $\sum b_I(x) T^I=0$ or $\sum \beta_I(x)T^I=0$ in $\hr x_r$,
which means that every $b_I$ vanishes at $x$ (hence at $x\al$) or every $\beta_I$ vanishes
at $x$ (hence at $x\al$). As a consequence, $\mathscr X=\mathscr Z(f)\bigcup \mathscr Z(g)$.

Assume that $B$ is reduced. Let $f=\sum b_IT^I$ be a nilpotent element of $B\otimes_AA_r$. 
By the above remark, $\mathscr X=\mathscr Z(f)$.
This means that for every $I$, the function $b_I$ vanishes
at every point of $\mathscr X$, hence is nilpotent, hence is zero because $B$ is reduced. 
Therefore $f=0$ and $B\otimes_A A_r$ is reduced.

Assume that $\spec B$ is irreducible (but $B$ can now have non-trivial nilpotents), 
and let us prove that $\spec B\otimes_A A_r$ is irreducible. 
By quotienting by the nilradical of $B$ (which does not modify the topology of the schemes
involved), we can assume that $B$ is reduced, hence a domain. Note that $B\to B\otimes_AA_r$
is injective
by faithful flatness of $A_r$ over $A$, hence $B\otimes_A A_r\neq 0$. 
Let $f=\sum b_I T^I$ and $g=\sum \beta_I T^I$ be two elements of $B\otimes_A A_r$
with $fg=0$. We have then $\mathscr X=\mathscr Z(f)\cup \mathscr Z(g)$.
By irreducibility of $\mathscr X$, we have $\mathscr X=\mathscr Z(f)$ or $\mathscr X=\mathscr Z(g)$. 
Suppose that $\mathscr X=\mathscr Z(f)$. For every $I$, the function $b_I$
vanishes
at every point of $\mathscr X$, hence is nilpotent, hence is zero because $B$ is a domain. 
We thus have $f=0$. Analogously, $g=0$ if $\mathscr X=\mathscr Z(g)$. 
This ends the proof when 
$\mathscr X$ is affine.

The assertion about reducedness
is local, hence holds in fact for arbitrary $\mathscr X$. 
Let us now assume that $\mathscr X$ is irreducible, but not
necessarily affine. Choose a non-empty covering $(\mathscr X_i)$ of $\mathscr X$ by non-empty open affine subsets
(as $\mathscr X$ is irreducible, it is non-empty). 
By the affine case
already proven, $\mathscr X_i\times_A A_r$ is
irreducible for every $i$. And for every $(i,j)$ the intersection $(\mathscr X_i\times_AA_r)\cap (\mathscr X_j\times_A A_r)$ is non-empty,
because it surjects onto $\mathscr X_i\cap \mathscr X_j$
by faithful flatness of $A\to A_r$; it follows that $\mathscr X\times_A A_r$ is irreducible. 
\end{proof}

%\begin{rema}
%Let $\mathscr X$ be an irreducible $A$-scheme
%locally of finite type, and let $r$ be a $k$-free
%polyradius such that $\abs{k_r\gpm}\neq \{1\}$
%and $A_r$ is strictly $k_r$-affinoid.
%By Lemma \ref{lem-xxr-scheme}
%above, $\mathscr X\times_A A_r$ is irreducible. It is thus finite-dimensional
%by Lemma \ref{lem-catenarity}; since
%\[\dim_k \mathscr X\an=\dim_{k_r} (\mathscr X\an)_r=\dim_k( \mathscr X \times_A A_r)\an=\dim \mathscr X \times_A A_r\]
%(the last equality
%follows from \ref{lem-gaga-dimension}), the analytic space $\mathscr X\an$
%is finite-dimensional. 
%\end{rema}

\begin{lemm}\label{lem-codim-kr}
Let $\mathscr X$ be an irreducible
$A$-scheme locally
of finite type and let $\mathscr Y$
be an irreducible closed subset of $\mathscr X$.

\begin{enumerate}[1]
\item The scheme $\mathscr X$ is finite-dimensional.

\item The analytic space $\mathscr X\an$ is finite-dimensional
and
\[ \dim \mathscr X\an=\dim \mathscr Y\an +\codim (\mathscr Y,\mathscr X)\]

\item Let $r$ be a $k$-free polyradius. We have the following
equalities and inequalities (recall that both 
$ \dim
\mathscr X\times_A A_r$ and $\mathscr Y\times_A A_r$ are irreducible
by Lemma \ref{lem-xxr-scheme}
above):

\begin{enumerate}[b]
\item $\dim \mathscr X \leq  \dim \mathscr X\times_A A_r$.
\item $\mathrm{codim}( \mathscr Y\times_A A_r, \mathscr X\times_A A_r)=
\mathrm{codim}(\mathscr Y, \mathscr X).$
\end{enumerate}

\end{enumerate}
\end{lemm}

\begin{proof}
We begin with (3).
Lemma \ref{lem-xxr-scheme} 
ensures that the pre-image
in $\mathscr X\times_A A_r$ of any irreducible closed subset of $\mathscr X$
is still irreducible. Moreover, since 
$\mathscr X\times_A A_r\to \mathscr X$
is surjective
(by faithful flatness of $A\to A_r$), the pre-images in $\mathscr X\times_A A_r$
of two distinct subsets of $\mathscr X$ are still distinct. These facts
imply (3a) as well as the inequality
\[\codim (\mathscr Y,\mathscr X)\leq \codim (\mathscr Y\times_A A_r,\mathscr X\times_A A_r).\]
We want to prove the reverse inequality. 
By catenarity of the scheme $\mathscr X\times_A A_r$
(which is excellent),
we can argue by induction on codimension; hence
we only have to consider the case where $\codim (\mathscr Y,\mathscr X)=1$.
Choose a parameter $f$ of the one-dimensional
local ring $\mathscr O_{\mathscr X,\eta}$; \ie, $f$ is an element
of $\mathscr O_{\mathscr X,\eta}$ whose zero locus on $\spec \mathscr O_{\mathscr X,\eta}$
is set-theoretically equal to the closed point (this existence of such an $f$ comes from \Th 13.4
of \cite{matsumura1986}). Then there exists an open neighborhood
$\mathscr U$ of the generic point $\eta$ of $\mathscr Y$ in $\mathscr X$
on which $f$ is defined and admits $\mathscr Y\cap \mathscr U$ as set-theoretical
zero locus. 

Now $(\mathscr Y\cap \mathscr U)\times_A A_r$ is the zero-locus of $f$ in the irreducible scheme $\mathscr U\times_A A_r$,
hence the codimension of the proper
irreducible closed subset $\mathscr Y\times_A A_r$ of 
$\mathscr X\times_A A_r$ is equal to 1 by the
\emph{Hauptidealsatz}, whence (3b). 

Now take $r$ such that $\abs{k_r\gpm}\neq \{1\}$
and $A_r$ is strictly $k_r$-affinoid. The irreducible
scheme $\mathscr X\times_A A_r$ is finite-dimensional
by Lemma \ref{lem-catenarity}, and (1) follows
then from (3a). 

Since one has the equalities
$\dim_k \mathscr X\an=\dim_{k_r} (\mathscr X\an)_r=\dim_{k_r}( \mathscr X \times_A A_r)\an$
and the same for $\mathscr Y$, assertion (2) follows from (3b),  Lemma \ref{lem-catenarity}.
and Lemma \ref{lem-gaga-dimension}.
\end{proof}

\begin{rema}
Assertion (2) was already known in the affinoid case, \cf \cite{ducros2007}, \Prop 1.11;
we have simply adapted the latter's proof to our more general setting. 
\end{rema}

\begin{rema}
Let $\mathscr X$ be a scheme locally of finite type over an affinoid algebra $A$. By Lemma \ref{lem-codim-kr}
above,
every irreducible component of $\mathscr X$ is finite-dimensional. The scheme
$\mathscr X$ itself is then finite-dimensional
if and only if the dimensions of its irreducible components
are uniformly bounded above. 

Assume that this is the case. The dimensions of the fibers of the map $\mathscr X\to \spec A$ are
then bounded above  by some integer $d$. It follows then from Proposition \ref{prop-gaga-dim} (1)
that the dimensions of the fibers of the map $\mathscr X\an \to \mathscr M(A)$
are bounded above by $d$. This implies in view of \ref{ss-dim-sobafi} (2)
that $\dim \mathscr X\an\leq \dim \mathscr M(A)+d$; in particular, $\mathscr X\an$ is finite-dimensional. 

\end{rema}

\begin{coro}
\label{cor-codim-schft}
Let $\mathscr X$ be a finite-dimensional
scheme locally of finite type over an affinoid algebra
and let $\mathscr Y$ be a Zariski-closed subset of $\mathscr X$. Let $x$ be a point
of $\mathscr X\an$. We have the following equalities: 

\begin{enumerate}[1]
\item $\codim (\mathscr Y\an, \mathscr X\an)=\codim (\mathscr Y,\mathscr X).$
\item $\codim_x (\mathscr Y\an, \mathscr X\an)=\codim_{x\al} (\mathscr Y,\mathscr X).$
\end{enumerate}
\end{coro}

\subsection{}
In the proof of Lemma \ref{lem-gaga-dimension}, 
we have used the fact (due to GAGA over an affinoid algebra, see \ref{ss-def-analytify})
that in the proper case,
the irreducible components of the analytification
are the analytifications of the irreducible components. We are now going to explain why
this holds without any properness assumption. For that purpose, we first need to establish
the compatibility between analytfication and normalization (\cite{ducros2009}, \S 5);
this is achieved by the following lemma (which follows directly from the constructions involved when $\mathscr X=\spec A$). 

\begin{lemm}
\label{lem-gaga-rednorm}
Let $\mathscr X$ be an $A$-scheme locally of finite type, and let $\mathscr Y$ denote
its normalization. 

\begin{enumerate}[1]
\item The closed immersion $(\mathscr X_{\mathrm{red}})\an\hookrightarrow \mathscr X\an$ identifies $(\mathscr X_{\mathrm{red}})\an$
with $(\mathscr X\an)_{\mathrm{red}}$. 

\item The finite morphism $\mathscr Y\an \to \mathscr X\an$ identifies $\mathscr Y\an$ with the normalization of $\mathscr X\an$.
\end{enumerate}
\end{lemm}

\begin{proof} 
By construction, the closed immersion $(\mathscr X_{\mathrm{red}})\an\hookrightarrow \mathscr X\an$ 
is defined by a nilpotent ideal, and its source is reduced (Lemma \ref{gaga-concrete}), whence (1). 

Let us now prove (2). In view of (1), we can replace $\mathscr X$ with $\mathscr X_{\mathrm{red}}$, hence assume
that $\mathscr X$ is reduced. The formation of the normalization of an analytic space
commutes to restriction to 
analytic domains (\cite{ducros2009}, Lemme 5.1.11) and is G-local on the target (\cite{ducros2009}, proof of \Th 5.13).
Hence by using compactification of affine charts on $\mathscr X$ we reduce to the case where the latter is proper over $A$.
Now let $Z$ be the normalization of $\mathscr X\an$. The morphism $\pi \colon Z\to \mathscr X\an$ is finite,
the image of every irreducible component of $Z$ is
an irreducible component of $\mathscr X\an$, and the set of point
of $\mathscr X\an$ at which $\mathscr O_{\mathscr X\an}\to \pi_\ast \mathscr O_Z$ is {\em not} 
an isomorphism does not contain any irreducible component of $\mathscr X\an$
(indeed, the latter property is G-local,
and is fulfilled by construction of the normalization on any affinoid chart). 

By GAGA over affinoid algebras (\ref{ss-def-analytify})
applied to the coherent $\mathscr O_{\mathscr X\an}$-algebra $\pi_*\mathscr O_Z$, 
the finite morphism $Z\to \mathscr X\an$ arises from a (unique)
finite morphism $\mathscr Z\to \mathscr X$ (which we still denote by $\pi$). 
The scheme $\mathscr Z$ is normal because $\mathscr Z\an$
is normal  (Lemma \ref{gaga-concrete}),
the image of every irreducible component of $\mathscr Z$ is an irrreducible component of $\mathscr X$,
and the set of points of $\mathscr X$ at which $\mathscr O_{\mathscr X}\to \pi_\ast \mathscr O_{\mathscr Z}$ is {\em not}
an isomorphism does not contain any irreducible component of $\mathscr X$. If $(\mathscr X_i)$ denotes the family of irreducible
of $\mathscr X$, we thus may write $\mathscr Z=\coprod \mathscr Z_i$ where each $\mathscr Z_i$ is normal
and maps birationally onto $\mathscr X_i$ (equipped with its reduced structure). In other words, $\mathscr Z_i$
is the normalization of $\mathscr X_i$ for all $i$, and $\mathscr Z$
is therefore the normalization of $\mathscr Y$ of $\mathscr X$,
whence the equality $Z=\mathscr Y\an$. 
\end{proof}

\begin{prop}\label{prop-gaga-irrcomp}
Let $\mathscr X$ be an $A$-scheme
locally of finite type, and let $(\mathscr X_i)$ be the family of irreducible components
of $\mathscr X$. The $\mathscr X_i\an$'s are the irreducible components of $\mathscr X\an$.
\end{prop}

\begin{proof}
We first note that $(\mathscr X_i\an)$ is a locally finite family of Zariski-closed
subsets of $\mathscr X\an$, which are pairwise not comparable with respect to the inclusion relation (by surjectivity of $\mathscr X\an \to \mathscr X$). 
It is therefore sufficient to prove that $\mathscr X_i\an$ is irreducible for every $i$.
Otherwise said, we have reduced to the case where $\mathscr X$ is irreducible. 
Let $\mathscr Y$ be its normalization. Since $\mathscr Y\an$ is the normalization of $\mathscr X\an$ by
Lemma \ref{lem-gaga-rednorm}, 
it is sufficient by \Th 5.17
of \cite{ducros2009}
to prove that $\mathscr Y\an$ is connected. Let us choose a covering $(\mathscr Y_i)$ of the
irreducible, normal scheme $\mathscr Y$ by non-empty affine open subschemes. Since
$\mathscr Y_i\cap \mathscr Y_j\neq \emptyset$ for every$ (i,j)$, the intersection
$(\mathscr Y_i\an)\cap (\mathscr Y_j\an)$ is non-empty; hence it suffices to prove that $\mathscr Y_i\an$ is connected for every $i$. 

Fix $i$ and choose a normal, projective compactification $\mathscr Z$ of $\mathscr Y_i\an$. 
The analytic space $\mathscr Z\an$ is normal
(Lemma \ref{gaga-concrete})
and irreducible by GAGA over an affinoid algebra (\ref{ss-def-analytify}),
and $\mathscr Y_i\an$ is a non-empty Zariski-open subset of $\mathscr Z\an$.
It is then connected by the non-Archimedean avatar of Riemann's extension theorem (\cite{berkovich1990}, \Prop 3.3.14;
it is based upon the rigid-analytic version due to Lütkebohmert, \cite{lutkebohmert1974} \Th 1.6). 
\end{proof}

\chapter{Germs, Temkin's reduction and $\Gamma$-strictness}
\label{c-gstrict}

\begin{enonce*}{Convention}
We fix from now on and until the end of the memoir
an analytic field $k$.
\emph{We do not make any assumption on it}: it is not assumed to be algebraically closed, 
it can be of any characteristic and residue characteristic and is not necessarily perfect, the value group
$\abs {k\gpm}$ can be any subgroup of $\R\gpm_+$, such as $\{1\}$ or
$\R\gpm_+$, \etc
\end{enonce*}

Berkovich's theory makes a distinction between strictly
$k$-affinoid spaces, which are defined by a finite system of equations
on a \emph{unit}
compact polydisc, and general affinoid spaces, whose definition allows arbitrary radii. In Section \ref{s-gstrict}, 
we introduce an intermediate class of $k$-affinoid spaces, namely those whose
definition allows radii belonging to a given
subgroup $\Gamma$ of $\R_+\gpm$, which are called \emph{$\Gamma$-strict}; these are
the building blocks of the category of \emph{$\Gamma$-strict $k$-analytic spaces}.
The motivation for introducing such a notion 
is to keep under control the real parameters that are needed to define our spaces,
especially as far as the
description of the image
of a map is involved (Section \ref{s-images-maps}); the reader can ignore it at first reading.

If $k$ is non-trivially valued, 
the category of strictly $k$-analytic spaces is a full subcategory of the category of all analytic spaces, but
this result is by no way obvious; it was shown by Temkin in \cite{temkin2004}, using the theory of 
\emph{graded} reduction of (punctual)
analytic germs which he introduced for this purpose (and which is based upon graded Riemann-Zariski
spaces; \ie, spaces of
graded valuations). 
Using the same method, 
we
shall prove that the category of $\Gamma$-strict $k$-analytic spaces is a full subcategory of the categories
of all analytic spaces (\ref{ss-gstrict-fullfaith}). But we first give a detailed account
of Temkin's theory in Sections \ref{s-RZar} and \ref{s-temkin} -- we shall need it throughout the whole memoir, 
not only for questions related to $\Gamma$-strictness, but also because it is a powerful tool
for the local study of analytic spaces, and often a very efficient substitute for the theory of formal 
models (which is technically more involved, and moreover not available in the non-strict case).

Before doing this, we of course have to say a few words
about the notion of a (punctual) analytic germ. 
This is done in Section
\ref{s-an-germ}, in which 
we also introduce
the \emph{central dimension}
of such a germ (Definition \ref{def-centdim}).
This turns out to measure the difference
between the Krull dimension 
of a local ring of a good analytic space and the local dimension 
of the space at the corresponding point (Corollary \ref{cor-interp-centdim})
and will play an important role in this work.

Let us end this introduction 
by mentioning that in the strict
context, there
is no need for considering Temkin's graded reduction: the non-graded reduction
(based upon usual Riemann-Zariski spaces) which he had 
developped in \cite{temkin2000}
and is in some sense the ``degree 1" part of
the graded reduction, 
is sufficient. This phenomenon extends to the $\Gamma$-strict context: if one
works with a $\Gamma$-strict $k$-analytic germ, there is no need
to consider its whole graded reduction; it is sufficient to consider
its ``$\Gamma$-graded part", as explained in \ref{ss-graded-reduction}.

\section{$\Gamma$-strictness}\label{s-gstrict}

\begin{enonce}[remark]
{Notation}
We fix for the whole chapter a subgroup $\Gamma$ of $\R_+\gpm$ such that $\Gamma\cdot \abs {k\gpm}\neq \{1\}$; otherwise said,
$\Gamma$ is non-trivial whenever $k$ is trivially valued. 
\end{enonce}

\subsection{$\Gamma$-strict affinoid algebras}
\label{ss-gamstr-affalg}
Let~$A$ be a~$k$-affinoid algebra.
We shall say that a~$A$ is {\em ~$\Gamma$-strict} if it is a quotient of~$k\{r_1\inv T_1,\ldots, r_n\inv T_n\}$ for some $r_j$ {\em belonging to~$\Gamma$.} 

If~$A$ is a quotient of~$k\{r_1\inv T_1,\ldots, r_n\inv T_n\}$ for some $r_j$ {\em belonging to~$(\abs {k^\times}\cdot\Gamma)^{\Q}$}, then~$A$ is~$\Gamma$-strict. Indeed,
we can choose a~$k$-free polyradius~$s=(s_1,\ldots,s_m)$ such that $\abs{k_s^\times}\neq \{1\}$, every~$s_i$ belongs to~$\Gamma$,
and every~$r_j$ belongs to~$\abs {k_s^\times}^{\Q}$. This implies
that~$A_s$ is strictly~$k_s$-affinoid; the proof of \Cor
2.1.8 of \cite{berkovich1990} together with an easy induction on~$m$ shows then that~$A$ is~$\Gamma$-strict.

\subsection{$\Gamma$-strict $k$-affinoid spaces}\label{ss-gstr-ansp}
A~$k$-affinoid space will be said to be~$\Gamma$-strict if its algebra of analytic functions is~$\Gamma$-strict. 
If~$X$ is such a space, the spectral semi-norm of any analytic function~$f$ on~$X$
belongs to~$(\abs {k\gpm}\cdot \Gamma)^\Q\cup\{0\}$. Indeed, as we saw in~\ref{ss-gamstr-affalg} above, there exists a~$k$-free
polyradius~$s=(s_1,\ldots,s_m)$ such that $\abs{k_s\gpm}\neq \{1\}$, every $s_i$ belongs to~$\Gamma$,
and~$X_s$ is strictly~$k_s$-affinoid.
The spectral semi-norm of~$f$ can be computed on~$X_s$; hence it follows from \cite{bosch-g-r1984}, 6.2.1/4 that
it belongs to~$\abs {k_s\gpm}^\Q\cup \{0\}\subset (\abs {k^\times}\cdot \Gamma)^\Q\cup\{0\}$. 

Conversely, let~$X$ be a~$k$-affinoid space such that the spectral semi-norm of every analytic function on~$X$
belongs to~$(\abs {k^\times}\cdot \Gamma)^\Q\cup\{0\}$; then~$X$ is~$\Gamma$-strict.  Indeed, let~$A$ be the algebra
of analytic functions on~$X$, and let us fix an admissible epimorphism~$k\{T_1/r_1,\ldots, T_n/r_n\}\to A$. For every~$i$, let~$s_i$ be the
spectral radius of the image of~$T_i$ in~$A$. By assumption,~$s_i\in  (\abs {k^\times}\cdot \Gamma)^\Q\cup\{0\}$. Now set~$t_i=s_i$ if~$s_i\neq 0$
and take for~$t_i$ any element of~$(\abs {k^\times}\cdot \Gamma)^\Q\cap [0,r_i]$ if~$s_i=0$. 
The admissible epimorphism~$k\{T_1/r_1,\ldots, T_n/r_n\}\to A$
then factors through an admissible epimorphism~$k\{T_1/t_1,\ldots, T_n/t_n\}\to A$,
whence the~$\Gamma$-strictness of~$A$. 

\subsection{$\Gamma$-strict $k$-analytic spaces}\index{analytic space!g@$\Gamma$-strict}\index{g@$\Gamma$-strict analytic space}
The class of~$\Gamma$-strict affinoid spaces is a {\em dense} class in
the sense of \cite{berkovich1993}, \S 1 (the assumption that~$\abs {k^\times}\cdot \Gamma\neq\{1\}$ has
benne made precisely
to ensure this property). It thus gives rise to a corresponding category of analytic spaces, which is called
the category of {\em $\Gamma$-strict} $k$-analytic spaces. 

\begin{rema}\label{rem-gstrict-careful}
If~$U$ and~$V$ are two~$\Gamma$-strict affinoid domains of a {\em separated}~$k$-analytic space~$X$,
then~$U\cap V$ is a~$\Gamma$-strict affinoid domain of~$X$ by the standard argument:
it
is a closed analytic subspace of $U\times_k V$. As a consequence, a separated~$k$-analytic space admits a $\Gamma$-strict $k$-analytic
structure if and only if it admits a G-covering by~$\Gamma$-strict affinoid domains, and this structure is then unique.
This is for instance the case for separated
boundaryless $k$-analytic spaces: indeed,
every point of such a space has a $\Gamma$-strict affinoid neighborhood. 

But be aware that a general (\ie, non necessarily separated) $k$-analytic space could {\em a priori}
admit several $\Gamma$-strict $k$-analytic structures. We
shall see in
\ref{ss-gstrict-germ-can}
that this is actually impossible, but to avoid
circular reasoning, we
shall for the moment say that a $k$-analytic space $X$
``is" $\Gamma$-strict if $X$ admits {\em some}
$\Gamma$-strict $k$-analytic structure. 
A $\Gamma$-strict analytic domain of such a space
will
mean an analytic domain admitting {\em some}
$\Gamma$-strict structure, but not necessarily a $\Gamma$-strict analytic domain for the given 
$\Gamma$-strict $k$-analytic structure on the ambient space.

\end{rema}

\begin{rema}
If $\abs {k\gpm}\neq \{1\}$ then what what we call $\{1\}$-strictness is nothing but the usual
notion of strictness, and we shall of course say strict instead of $\{1\}$-strict. 
\end{rema}

\begin{rema}\label{rem-gamma-strict}
There is a small difference between our definition of $\Gamma$-strictness and that of Conrad and Temkin in \cite{conrad-tXXX}: they require the group $\Gamma$ to contain $|k^\times|$, which we do not. Nevertheless, this is more or less irrelevant: indeed, an analytic space is $\Gamma$-strict in our sense if and only if it is $\Gamma\cdot |k^\times|$-strict in Conrad and Temkin's sense. 

The reason why we have chosen to relax the assumption on $\Gamma$ is the following. Using our convention, if $X$ is a $\Gamma$-strict $k$-analytic space, then $X_L$ is a $\Gamma$-strict $L$-analytic space for any analytic extension $L$ of $k$; but such a statement simply {\em does not make any sense in general} if one uses Conrad and Temkin's convention, because even if $|k\gpm|\subset \Gamma$ it may happen that $|L\gpm|$ is not contained in $\Gamma$. 

\end{rema}

\begin{rema}
One can also define the notion of a $\Gamma$-strict analytic space, without any mention of the ground field: this is a pair $(L,X)$ where $L$ is 
an analytic field such that $|L^\times|\cdot \Gamma \neq 1$ and where $X$ is a $\Gamma$-strict $L$-analytic space; morphisms between such spaces are defined in the obvious way.
\end{rema}

\section{Analytic germs}\label{s-an-germ}

\subsection{}\label{ss-recall-germs} Let us recall briefly the definition
of the category of (punctual) germs\index{germ of analytic space}
of $k$-analytic spaces (which we shall simply call for short
{\em $k$-analytic germs}) given by Berkovich in \cite{berkovich1993}, 3.4. This is the localization 
of the category of pointed $k$-analytic spaces by the class of morphisms $\phi \colon (Y,y)\to (X,x)$ having the following
property: {\em $\phi$ induces an isomorphism between an open neighborhood of $y$ and an open neighborhood
of $x$}. A germ $(X,x)$ is said to have a given property (preserved by restriction to open subsets)
if $x$ admits a neighborhood in $X$ having this property. 

By replacing $k$-analytic spaces with analytic spaces in the above construction,
we get the notion of an analytic germ without mention of the ground field. 

If $(X,x)$ is an analytic germ, a germ of the form $(V,x)$ for
$V$ an analytic domain of $X$ containing $x$ will simply be called an analytic domain of $(X,x)$. \index{analytic domain (of an analytic germ)}

We
can define in a similar way the 
category of $\Gamma$-strict $k$-analytic germs
(as a localization of the category of pointed $\Gamma$-strict
analytic spaces). Be aware that Remark \ref{rem-gstrict-careful}
applies mutatis mutandis in this context, and we 
shall say that a germ $(X,x)$ ``is" $\Gamma$-strict if $x$ admits a $\Gamma$-strict
analytic neighborhood in $X$; \ie, $(X,x)$ admits  {\em some}
$\Gamma$-strict $k$-analytic structure, which is not a priori compatible
with the one on $X$ nor
unique
unless $(X,x)$ is separated (but it will be unique a posteriori).

\begin{defi}\label{def-centdim}
Let $X$ be an analytic
space and let $x$ be a point of $X$. 
The infimum of the integers $\dim_k \overline{\{x\}}^{V\zar}$ for $V$ running through the set of analytic {\em neighborhoods}
of $x$ in $X$ only depends on the germ $(X,x)$;  it will be called the {\em~$k$-analytic central dimension}\index{dimension!central}\index{central dimension}
of the germ~$(X,x)$ and will be denoted by~$\mathrm{centdim}_k(X,x)$, or
usually simply by~$\mathrm{centdim}(X,x)$ if there is no ambiguity about the ground field. 
\end{defi}

\begin{rema}
We obviously have $\mathrm{centdim}(X,x)\leq \dim \adhz {\{x\}}X$, and this inequality can be strict; 
see Remark \ref{rem-centdim-counter}.

\end{rema}

\subsection{Basic properties of central dimension}\label{ss-basicprop-centdim}
Let $X$
be an analytic space
and
let $x$ be a point of $X$. There exists an analytic neighborhood of $x$
in $X$ of dimension $\dim_xX$
(Remark \ref{rem-dimloc-open}); we thus have
$\mathrm{centdim}(X,x)\leq \dim_x X$. 

For every analytic neighborhood 
$V$ of $x$ in $X$, one has $\dim \overline{\{x\}}^{V\zar}\geq d_k(x)$; therefore $\mathrm{centdim}(X,x)\geq d_k(x)$.

If $Y$ is a closed analytic subspace of $X$ contaning $x$, it follows
from the definition that $\mathrm{centdim}(Y,x)=\mathrm{centdim}(X,x)$.
More generally, if $Y$ is a finite $X$-analytic space and if $y$ is a pre-image of $x$ on $Y$, then 
$\mathrm{centdim}(Y,y)=\mathrm{centdim}(X,x)$. Indeed, by topological
properness and topological
separatedness of finite morphisms
we may choose an open neighborhood $V$ of $x$ in $X$ such that
$\overline{\{x\}}^{V\zar}=\mathrm{centdim}(X,x)$ and such that $\overline{\{y\}}^{W\zar}=\mathrm{centdim}(Y,y)$,
where $W$ is the connected component of $y$ inside $Y\times_X V$. The image of~$\overline{\{y\}}^{W\zar}$ on~$V$ is a Zariski-closed subset of $V$ in which $x$ is Zariski-dense, hence it
coincides
with~$\overline{\{x\}}^{V\zar}$. One has then $\dim  \overline{\{x\}}^{V\zar}=\dim\overline{\{y\}}^{W\zar}$
by \ref{ss-dim-sobafi}, whence our claim.

\begin{exem} Since an irreducible $k$-analytic space is zero-dimensional if and only if it consists of one rigid point
(\ref{ss-x-dimzero}), the central dimension of a germ $(X,x)$ 
is zero if and only
if $x$ is rigid, in which case $d_k(x)=0$. 
\end{exem}

\begin{exem}\label{ex-centdim-abhyankar}
Let $X$ be an analytic space, let $Y$ be a Zariski-closed subset of $X$
and let $x$ be
a point of $Y$. Assume that 
$\dim_x Y=d_k(x)$; \ie, $x$ is an Abhyankar point of $Y$. By \ref{ss-basicprop-centdim}
we have
\[\dim_x Y\geq \mathrm{centdim}(Y,x)=\mathrm{centdim}(X,x)\geq d_k(x).\]
It follows that $\mathrm{centdim}(X,x)=d_k(x)$. 
\end{exem}

\begin{exem}Let $X$ be
a  curve and let $x$
be a point of type 4 of $X$
(according to Berkovich's classification described in \cite{berkovich1990}, Chapter 4). 
Then $d_k(x)=0$ but $x$ is not rigid; therefore $\mathrm{centdim}(X,x)=1$. 
\end{exem}

\begin{lemm}\label{lem-mox-max}
Let $X$ be an affinoid space and let $x$ be a point of $X$.
The following are equivalent:

\begin{enumerate}[label=\textnormal{(\roman{enumi})}]

\item $\mathfrak m_{x\al} \mathscr O_{X,x}=\mathfrak m_x$. 

\item $\dim_{\mathrm {Krull}} \mathscr O_{X,x}= \dim_{\mathrm{ Krull}} \mathscr O_{X\al,x\al}$. 

\item $\mathrm {centdim}(X,x)=\dim \overline{\{x\}}^{X\zar}$. 

\end{enumerate}
\end{lemm}

\begin{proof}
Let $\omega$ be the closed point of $\spec  \mathscr O_{X,x}$ and let $p\colon \spec  \mathscr O_{X,x}
\to \spec \mathscr O_{X\al, \al}$ be the canonical map. Let us consider the following assertions:

\begin{enumerate}[a]
\item $p^{-1}(x\al)=\{\omega\}_{\mathrm{red}}$  {\em scheme-theoretically}. 

\item $p^{-1}(x\al)=\{\omega\}$ {\em set-theoretically}. 

\item $\dim p^{-1}(x\al)=0$. 

\end{enumerate}

Assertion (i) is tautologically equivalent to (a). Since the fibers of~$p$ are (geometrically) regular, and in particular reduced, one 
has (a)$\iff$(b), and it is clear that (b)$\iff$(c). As 
(c)$\iff$(ii) by flatness of $p$, we eventually get the equivalence
(i)$\iff$(ii). 

We are now going to prove the equivalence (i)$\iff$(iii). For that purpose, we may replace $X$
with
any of its closed analytic subspaces containing $x$, and in particular with $\overline{\{x\}}^{X\zar}_{\mathrm{red}}$.  
Hence we may assume that $X$ is integral and that $x\al$ is the generic point of $X\al$; note that under this assumption
$\mathscr O_{X\al, x\al}$ is a field and $\mathfrak m_{x\al}=0$. Let $d$
be the dimension of $X$. Since $X$ is irreducible, it is purely $d$-dimensional, and so are all
of its analytic domains. 

Assume that~(i) is true, \ie, $\mathscr O_{X,x}$ is a field. Let~$V$ be an affinoid neighborhood of~$x$ in~$X$ and let~$Z$
be a Zariski-closed subset of~$V$ containing~$x$. Let~$f_1,\ldots,f_n$ be analytic
functions on~$V$ that generate the ideal of functions vanishing pointwise on~$Z$.
For any~$i$, we have~$f_i(x)=0$; the image of~$f_i$ in~$\mathscr O_{X,x}$ is then not invertible,
hence is zero. Therefore~$Z$ contains a neighborhood~$U$ of~$x$ in~$V$. 
The dimension of~$U$ being equal to~$d$, the dimension of~$Z$ (which is bounded by~$d$) is equal to~$d$ too; therefore, (iii) is proved.

Assume that~(iii) is true; \ie, $\mathrm{centdim}(X,x)=d$. Let~$f$ be a non-zero element of~$\mathscr O_{X,x}$,
and let~$V$ be an affinoid neighborhood of~$x$ on which~$f$ is defined. Let~$Y$ and~$Z$ be two
irreducible components of~$V$ containing~$x$. Both are of dimension~$d$; if~$Y\neq Z$,
their intersection is a Zariski-closed subset of~$V$ containing~$x$ and of dimension strictly smaller than~$d$, which contradicts~(iii). 
Hence there is only one irreducible component~$Y$ of~$V$ containing~$x$.
The zero-locus of~$f$ on~$Y$ is a Zariski-closed subset of~$Y$ which is not equal to the whole
of $Y$,
because otherwise
$f$ would vanish pointwise on a neighborhood of~$x$, hence would vanish
in the reduced local ring ~$\mathscr O_{X,x}$; but by assumption, this is not the case.
Therefore the dimension of the zero-locus of~$f$ is  strictly smaller than~$d$, hence this locus can not contain~$x$ because of (iii).
As a consequence,~$f$ is invertible in~$\mathscr O_{X,x}$. Thus the latter is a field, and (i) is proved. 
\end{proof}

\begin{coro}\label{cor-interp-centdim}
Let $X$ be a good analytic space and let $x$
be a point of $X$. One has the equality
$$\mathrm{centdim}(X,x)+\dim_{\mathrm{Krull}}\mathscr O_{X,x}=\dim_x X.$$
In particular if $x$ is rigid then $\dim_{\mathrm{Krull}}\mathscr O_{X,x}=\dim_x X.$
\end{coro}

\begin{proof}
Set~$d=\mathrm{centdim} (X,x)$. Up to shrinking $X$ we may assume
that it is~affinoid and that~$\dim \overline{\{x\}}^{X\zar}=d$. Let~$X_1,\ldots,X_n$ be the irreducible components of~$X$ that
contain~$x$. For every~$i$, set~$d_i=\dim X_i$ and~$\delta_i=\codim_{\mathrm{Krull}}(\overline{\{x\}}^{X\zar},X_i)$. 

One has $\dim_x X=\max d_i$, and since $\dim \overline{\{x\}}^{X\zar}=d$, Lemma~\ref{lem-mox-max}
ensures that  $\dim_{\mathrm {Krull}}\mathscr O_{X,x}$ coincides
with $\dim_{\mathrm{Krull}}\mathscr O_{X\al, x\al}=\codim_{\mathrm{Krull}}(\overline{\{x\}}^{X\zar},X)=\max \delta_i$.
But $d_i=\delta_i+d$ for every $i$ (see \cite{ducros2007}, \Prop 1.11), whence we
get the equality
$\dim_x X=\dim_{\mathrm {Krull}}\mathscr O_{X,x}+d$.
\end{proof}

\begin{exem}\label{ex-centdim-recap}
Let $X$ be a good analytic space and let $Y$
be a Zariski-closed subset of $X$. Let $x$ be a point of $Y$, 
and assume that $d_k(x)=\dim_x Y$; \ie, $x$ is Abhyankar
in $Y$. We then have
$\mathrm{centdim}(X,x)=d_k(x)=\dim_x Y$ (Example \ref{ex-centdim-abhyankar})
and $\dim \adhz {\{x\}}X=d_k(x)$ (Remark \ref{rem-abh-point}). 
It follows then from Corollary \ref{cor-interp-centdim}
that 
\[\dim _{\mathrm{Krull}}\mathscr O_{X,x}=\dim_x X-\dim_x Y\]
(in particular, $\mathscr O_{X,x}$ is artinian as soon as $x$ is Abhyankar in $X$; \eg, $Y=X$). 

If
$X$ is moreover assumed to be affinoid, 
we deduce from Lemma \ref{lem-mox-max}
that $\dim _{\mathrm{Krull}}\mathscr O_{X\al,x\al}$ is also equal to $\dim_xX-\dim_x Y$
and that $\mathfrak m_{x\al}\mathscr O_{X,x}=\mathfrak m_x$.

\end{exem}

\section{Around graded Riemann-Zariski spaces}\label{s-RZar}

Our purpose is now to give
a short account of Temkin's theory of (graded) reduction of analytic germs; our reference is the foundational article \cite{temkin2004}. In this 
section, we shall introduce all the required notions
about graded Riemann-Zariski spaces; the application to analytic germs
will be explained in Section \ref{s-temkin}.

\subsection{}\label{ss-def-plkzr}\index{graded!Riemann-Zariski space}\index{Riemann-Zariski space!graded}
Let~$K$ be a graded field. If~$L$ is any graded extension of~$K$, we shall denote by~$\P_{L/K}$ the ``graded Riemann-Zariski space of~$L$ over~$K$"; 
\ie, the set of equivalence classes of
graded valuations on~$L$ whose restriction to~$K$ is trivial (or, in other words, whose graded ring contains~$K$).  For any set~$S$ of homogeneous elements of~$L$, we denote by~$\P_{L/K}\{S\}$ the subset of~$\P_{L/K}$ that consists of graded valuations~$\abs {\cdot}$ such that~$\abs f\leq 1$ for every~$f\in S$; note that $\P_{L/K}\{S\}=\P_{L/K}\{S\setminus (S\cap K)\}$ (in particular if $0\in S$ then it can be removed without modifying $\P_{L/K}\{S\}$).
We endow~$\P_{L/K}$ with the topology generated by the sets of the form ~$\P_{L/K}\{S\}$
for $S$ a {\em finite} set of homogeneous elements of $L$, which are called {\em affine}
open subsets of $\P_{L/K}$. Note that $\P_{L/K}=\P_{L/K}\{\emptyset\}$ is itself affine, and that the intersection of two affine open subsets is affine. Note also that every affine open subset of $\P_{L/K}$ contains the trivial
graded valuation; the latter is thus a generic point of $\P_{L/K}$, and is easily seen to be the only one.
Any affine open subset of~$\P_{L/K}$ (especially,~$\P_{L/K}$ itself) is quasi-compact (\cite{temkin2004}, 5.3.6); since the intersection of two
affine open subsets of $\P_{L/K}$ is affine, $\P_{L/K}$ is quasi-separated.

\subsection{Functoriality}\label{ss-rz-funct}
If~$F$ is any graded extension of~$L$, and if~$E$ is a graded subfield of~$F$
 such that~$E\cap L\supset K$, restriction of graded valuations defines a map~$r\colon \P_{F/E}\to \P_{L/K}$. For every set $S$ of homogeneous elements of $L$
we have $r^{-1}(\P_{L/K}\{S\})=\P_{F/E}\{S\}$; by applying this when $S$ is finite
we see that $r$ is continuous and quasi-compact. Morever, $r$ is surjective as soon as $E=K$.

\subsection{$\Gamma$-strictness}\label{ss-delta-strict}
If~$\Gamma$ is a subgroup of~$\R\gpm_+$, a quasi-compact open subset~$\mathsf U$ of~$\P_{L/K}$
is said to be {\em
$\Gamma$-strict} if~$\mathsf U$ is the pre-image of some
(possibly empty)
quasi-compact open subset of~$\P_{L ^\Gamma/K ^\Gamma}$; equivalently, $\mathsf U$ is~$\Gamma$-strict if and only if it is a finite union of affine open subsets whose definition only involve homogeneous elements of $L^\Gamma$. We shall simply say {\em strict} instead of~$\{1\}$-strict. 

Let $S$ be a finite set of homogeneous elements of $L$. If $S\subset L^{(\mathfrak D(K)\cdot \Gamma)^\Q}$, then there exists a finite set
$S'$ of homogeneous elements of $L^\Gamma$ such that $\P_{L/K}\{S\}=\P_{L/K}\{S'\}$, and $\P_{L/K}\{S\}$ is therefore $\Gamma$-strict. Indeed, for every $f\in S$ there exists a non-zero homogeneous element $a_f$ of $K$ and a positive integer $n_f$ such that $a_ff^{n_f}\in L^\Gamma$, and we can set $S'=\{a_f f^{n_f}\}_{f\in S}$.  Conversely, if $\P_{L/K}\{S\}$ is $\Gamma$-strict, then $S\subset L^{(\mathfrak D(K)\cdot \Gamma)^\Q}$ by \Prop 2.5~(i) of \cite{temkin2004}.

\begin{rema}
Our definition of $\Gamma$-strictness is not exactly that of Temkin in \cite{temkin2004}.
Indeed, Temkin requires $\Gamma$ to contain $\mathfrak D(K)$, 
which we do not. More precisely, $\Gamma$-strictness in our sense is equivalent to $\Gamma \cdot \mathfrak D(K)$-strictness
in Temkin's sense. We have
made this choice for the sake of consistency
with our definition of strictness in analytic geometry, see Remark \ref{rem-gamma-strict}. 
\end{rema}

\begin{rema}\label{rem-temkin-transfer}
Because of this difference between our definition and Temkin's, there are some results of \cite{temkin2004}
that
cannot be applied directly in our setting.
But we shall
remedy this by using the following fact: {\em the natural continuous
map $\P_{L^{(\mathfrak D(K)\cdot \Gamma)^\Q}/K}\to \P_{L^\Gamma/K^\Gamma}$ is a homeomorphism}.

To see
this, we first note that in view of~\ref{ss-delta-strict}, it is sufficient to prove that this map 
is bijective. Now let $\abs {\cdot}$ be an element of
$\P_{L^\Gamma/K^\Gamma}$ and let $f$ be a homogeneous element of $L^{(\mathfrak D(K)\cdot \Gamma)^\Q}$. By definition, 
there exists a non-zero homogeneous element $a$ of $K$, a homogeneous element $g$ of $L^\Gamma$ and an integer $n$ such that $f^n=ag$. 
One checks straightforwardly that the element $\abs g^{1/n}$ only depends on $f$, and not on $(a,n,g)$, and that
the assignment $f\mapsto \abs g ^{1/n}$ is 
is the unique pre-image of $\abs {\cdot}$ on 
$\P_{L^{(\mathfrak D(K)\cdot \Gamma)^\Q}/K}$. 

Let us mention an important consequence of the above: $\P_{L/K}\to \P_{L^\Gamma/K^\Gamma}$
is the composition of the surjection 
$\P_{L/K}\to \P_{L^{(\mathfrak D(K)\cdot \Gamma)^\Q}/K}$ and of the homeomorphism
$ \P_{L^{(\mathfrak D(K)\cdot \Gamma)^\Q}/K}\to \P_{L^\Gamma/K^\Gamma}$, so
it is surjective. 

\end{rema}

\subsection{The category $\mathscr S_{L/K}$}\label{ss-def-slk}
We denote by $\mathscr S_{L/K}$ the full subcategory of the category of topological spaces over
$\P_{L/K}$ consisting of objects $\mathsf X\to \P_{L/K}$ satisfying the following conditions:

\begin{enumerate}[1]
\item The space $\mathsf X$ is quasi-compact and quasi-separated. 

\item The morphism $\mathsf X\to \P_{L/K}$ is a local homeomorphism.

\end{enumerate}

For instance, any quasi-compact open subset of  $\P_{L/K}$
is an object of $\mathscr S_{L/K}$ through its inclusion into $\P_{L/K}$. 

Let $\mathsf X$ be an object of $\mathscr S_{L/K}$. A {\em chart}
of $\mathsf X$ is a quasi-compact open subset $\mathsf U$ of $\mathsf X$
such that $\mathsf U\to \P_{L/K}$ is an open immersion. An {\em atlas}
of $\mathsf X$ is a finite covering of $\mathsf X$ by charts. 

Let $\eta$
be the unique generic point of $\P_{L/K}$
and let $\mathsf X_\eta$ be the set of pre-images of $\eta$ on $\mathsf X$. Every non-empty chart $\mathsf U$ of $\mathsf X$
has a unique intersection point with $\mathsf X_\eta$, which is the unique generic point of $\mathsf U$. 
For $\xi \in \mathsf X_\eta$, let us denote by $\mathsf U_\xi$ the union of all charts of $\mathsf X$
that contain $\xi$. By the above, $\xi$ is the unique generic point of $\mathsf U_\xi$,
and
$\mathsf X$ is the disjoint union of the $\mathsf U_\xi$'s for $\xi$ running through $\mathsf X_\eta$; 
it follows by quasi-compactness that $\mathsf X_\eta$ is finite
and each $\mathsf U_\xi$ is quasi-compact. 
Note that the $\mathsf U_\xi$'s are the connected 
components of $\mathsf X$. 

Let~$F$ be a graded extension of~$L$, let~$E$ be a graded subfield of~$F$
containing $K$, 
let $r\colon \P_{F/E}\to \P_{L/K}$ be the natural map,
and let
$$\xymatrix{
{\mathsf Y}\ar[r]\ar[d]&{\P_{F/E}}\ar[d]^r\\
{\mathsf X}\ar[r]&{\P_{L/K}}
}$$ be a commutative diagram of topological spaces in which $\mathsf Y\to \P_{F/E}$ and
$\mathsf X\to \P_{L/K}$
belong
respectively to 
$\mathscr S_{F/E}$ and
$\mathscr S_{F/K}$. Since $r$ is quasi-compact by \ref{ss-rz-funct}, 
the continuous map $\mathsf Y\to \mathsf X$ is quasi-compact too. 
If $\mathsf V$ is a connected component of $\mathsf Y$, its generic point lies above
the generic point of a connected component of $\mathsf X$: this comes from
the fact that $r$ sends the trivial graded
valuation on $F$ to the trivial graded valuation on $L$. 

\subsection{$\Gamma$-strict objects of $\mathscr S_{L/K}$}\label{ss-def-dstrict}
Let $\Gamma$ be a subgroup of $\R_+\gpm$ and let $\mathsf X$ be an object of $\mathscr S_{L/K}$. 
A chart of $\mathsf X$ is said to be $\Gamma$-strict
if its image on $\P_{L/K}$ is $\Gamma$-strict. An atlas of $\mathsf X$ is called $\Gamma$-strict if
it consists of $\Gamma$-strict charts with pairwise $\Gamma$-strict intersections. 
The object $\mathsf X$ is called $\Gamma$-strict if it admits a $\Gamma$-strict atlas; this amounts to require that $\mathsf X$ is isomorphic
to $\mathsf Y\times_{\P_{L^\Gamma/K^\Gamma}}\P_{L/K}$ for some 
object $\mathsf Y$ of $\mathscr S_{L^\Gamma/K^\Gamma}$.
A non-empty quasi-compact open subset of $\mathsf X$ will be said to be $\Gamma$-strict
if it is $\Gamma$-strict as an objet of $\mathscr S_{L/K}$ (this is consistent with the previous definition when
$\mathsf X=\P_{L/K}$). 

\subsection{}\label{ss-dstrict-ff}
We are now going to list some properties
related to the notion of a $\Gamma$-strict object, with references
to Temkin's seminal paper \cite{temkin2004}. Note that Remark \ref{rem-temkin-transfer}
allows us to use Temkin's result in our setting (let us also mention that Temkin only deals with connected
non-empty spaces, but
this assumption is not seriously needed for what follows because one can argue componentwise
due to \ref{ss-def-slk}). 
Let $\mathsf X\to \P_{L/K}$ be a $\Gamma$-strict object of $\mathscr S_{L/K}$ and let
$F$ be a graded extension of $L$.

\begin{enumerate}[1]

\item {\em Canonicity}. The object of $\mathscr S_{L^\Gamma/K^\Gamma}$ from which $\mathsf X$ comes
is unique up to a unique homeomorphism by \cite{temkin2004}
\Prop 2.6. We shall denote it by $\mathsf X^\Gamma$. 
Note that in view of Remark~\ref{rem-temkin-transfer},
the natural continuous map $\mathsf X\to \mathsf X^\Gamma$ is surjective. 

\item {\em Functoriality}. Let
$$\xymatrix{
{\mathsf Y}\ar[r]\ar[d]&{\P_{F/K}}\ar[d]\\
{\mathsf X}\ar[r]&{\P_{L/K}}
}$$ be a commutative diagram
of topological spaces
with $\mathsf Y\to  \P_{F/K}$ a $\Gamma$-strict object of $\mathscr S_{F/K}$.
There exists a 
unique continuous map $\mathsf Y^\Gamma\to \mathsf X^\Gamma$ making the diagram 
$$\xymatrix{
&{\mathsf Y}\ar[rr]\ar'[d][dd]\ar[ld]&&{\P_{F/K}}\ar[dd]\ar[ld]\\
{\mathsf Y^\Gamma}\ar[rr]\ar[dd]&&{\P_{F^\Gamma/K^\Gamma}}\ar[dd]&\\
&{\mathsf X}\ar'[r][rr]\ar[ld]&& {\P_{L/K}}\ar[ld]\\
{\mathsf X^\Gamma}\ar[rr]&& {\P_{L^\Gamma/K^\Gamma}}
}$$ commute; this follows again from  \Prop 2.6 of \cite{temkin2004}~Note
that the uniqueness of the map
is an obvious consequence of the surjectivity of $\mathsf Y\to \mathsf Y^\Gamma$.

\item Let $\mathsf U$ be a $\Gamma$-strict, quasi-compact open subset of $\mathsf X$. By (2) the open immersion $\mathsf U\hookrightarrow \mathsf X$ is obtained from a continuous map $\mathsf U^\Gamma \to \mathsf X^\Gamma$ in the category $\mathscr S_{L^\Gamma/K^\Gamma}$, through the base-change functor by the map $\P_{L/K}\to \P_{L^\Gamma/K^\Gamma}$.
Since  the map $\P_{L/K}\to \P_{L^\Gamma/K^\Gamma}$ is surjective and since $\mathsf U\to \mathsf X$ is injective, $\mathsf U^\Gamma\to \mathsf X^\Gamma$ is injective, hence is 
an open immersion. Taking into account the uniqueness part in the assertion of (2),
or more directly the surjectivity of $\mathsf X\to \mathsf X^\Gamma$, we get the following: {\em the map $\mathsf V\mapsto \mathsf V\times_{\mathsf X^\Gamma}\mathsf X$ establishes a bijection
(inclusion preserving in both directions) from the set of quasi-compact open subsets of $\mathsf X^\Gamma$ 
to that of $\Gamma$-strict quasi-compact open subsets of $\mathsf X$}. 

\item Let $\Delta$ be a subgroup of $\R_+\gpm$ containing $\Gamma$. 
It follows immediately from the definitions that $\mathsf X$
 is $\Delta$-strict and that $\mathsf X^\Delta=\mathsf X^\Gamma\times_{\P_{L^\Gamma/K^\Gamma}}\P_{L^\Delta/K^\Delta}$
(and hence $X^\Delta$ is $\Gamma$-strict and $(X^\Delta)^\Gamma=X^\Gamma$). 

\item It follows immediately from the definition that $\mathsf Z:=\mathsf X\times_{\P_{L/K}}\P_{F/K}\to \P_{F/K}$ is a $\Gamma$-strict
object of $\mathscr S_{F/K}$ and that $\mathsf Z^\Gamma=\mathsf X^\Gamma\times _{\P_{L^\Gamma/K^\Gamma}}\P_{F^\Gamma/K^\Gamma}.$

\end{enumerate}

\section{Temkin's construction}\label{s-temkin}
To every $k$-analytic germ
$(X,x)$, Temkin associates a non-empty, connected object of $\mathscr S_{\hrt x/\widetilde k}$, which is denoted by $\widetilde{(X,x)}$
and called the {\em (graded) reduction}
of $(X,x)$. Let us first explain how it is defined, and then list some of its basic properties; proofs can be found in Section 4 of \cite{temkin2004}.\index{reduction (of an analytic germ)}

\subsection{Definition of $\widetilde{(X,x)}$: the good case}\label{ss-defredgrad-good}
Assume that $(X,x)$ is good; \ie, $x$ has an affinoid neighborhood in $X$. Let $V$ be such a neighborhood, say $V=\mathscr M(A)$. 
We endow $A$ with its spectral semi-norm, which allows to define 
a residue graded $\widetilde k$-algebra $\widetilde A$. It is finitely generated (this follows from \cite{temkin2004}, 
\Prop 3.1 (iii), applied to any presentation of $A$ as an admissible quotient of a Tate algebra). 
The map $f\mapsto f(x)$ induces a morphism of graded $\widetilde k$-algebras $\widetilde A\to \hrt x$; let $B$ denotes its image. The graded $\widetilde k$-algebra $B$ is finitely generated. The subset $\P_{\hrt x/\widetilde k}\{B\}$ of $\P_{\hrt x/\widetilde k}$ being equal to $\P_{\hrt x/\widetilde k}\{S\}$ for any finite set $S$ of homogeneous generators of $B$ over $\widetilde k$, it is an affine open subset of $\P_{\hrt x/\widetilde k}$ and in particular an object of $\mathscr S_{\hrt x/\widetilde k}$.  It only depends on $X$, and not on $V$ ; it is denoted by $\widetilde{(X,x)}$. 

Let $(f_1,\ldots, f_n)$ be invertible functions on $(X,x)$; for every $i$, set $r_i=\abs{f_i(x)}$. Let $(Y,x)$ be the analytic domain of $(X,x)$ defined by the conjunction of inequalities 
$\abs{f_i}\leq r_i$. Then the germ $(Y,x)$ is good and
$$\widetilde{(Y,x)}=\widetilde{(X,x)}\cap \P_{\hrt x/\widetilde k}\{\widetilde{f_1(x)}, \ldots, \widetilde{f_n(x)}\}\subset \P_{\hrt x/\widetilde k}. $$
Moreover, every good analytic domain of $(X,x)$ is of the above form (this is a consequence of Gerritzen-Grauert theorem, see
\cite{bosch-g-r1984}, \S 7.3.5 \Th 1 \Cor 3
in the strict case, and \cite{ducros2003}, Lemme 2.4 or  \cite {temkin2005}, \Prop 3.5
for the general case). 

\subsection{Definition of $\widetilde{(X,x)}$: the general case}\label{ss-defredgrad-general}
We do not suppose anymore that $(X,x)$ is good. The graded reduction $\widetilde{(X,x)}$ is then defined as
the colimit of the $\P_{\hrt x/\widetilde k}$-spaces $\widetilde{(Y,x)}$ for $(Y,x)$ running through the set of good analytic domains of $(X,x)$. 
This has the following concrete meaning: 

\begin{itemize}[label=$\bullet$]
\item For every good analytic domain $(Y,x)$ of $(X,x)$ the space
$\widetilde{(X,x)}$ 
is endowed with an open immersion $\iota_{(Y,x)}\colon
\widetilde{(Y,x)}\hookrightarrow \widetilde{(X,x)}$ of $\P_{\hrt x\widetilde k}$-spaces. 
\item For every good analytic domain $(Y,x)$ of $(X,x)$ and every good analytic domain $(Z,x)$ of $(Y,x)$, the open immersion 
$\iota_{(Z,x)}$ is equal to the restriction of $\iota_{(Y,x)}$ to $\widetilde{(Z,x)}$
(which is in a natural way
an open subset 
of $\widetilde{(Y,x)}$ as explained
in \ref{ss-defredgrad-good}).  
\item The space $\widetilde{(X,x)}$ is equal to $\bigcup_{(Y,x)}\iota_{(Y,x)}(\widetilde{(Y,x)})$ for $(Y,x)$ running through the set of good analytic
domains of $(X,x)$. 
\end{itemize}
%Let $(Y,x)$ and $(Z,x)$ be two good analytic domains of $(X,x)$. Both spaces
%$\widetilde{(Y,x)}$ and $\widetilde{(Z,x)}$ embed into $\widetilde{(X,x)}$. If
%$(Y\cap Z,x)$
%is good (\eg, $(X,x)$ is separated) then $\widetilde{(Y,x)}$ 
%and $\widetilde{(Z,x)}$ are simply glued along $\widetilde{(Y\cap Z,x)}$ 
%inside $\widetilde{(X,x)}$. The situation is slightly more complicated in general
%(one has to write $(Y\cap Z,x)$ as a finite
%union of good analytic domains). 

\subsection{Properties of $(X,x)$ than can be
seen on $\widetilde{(X,x)}$}
\label{ss-first-prop-reduction}
The germ $(X,x)$ is separated, \resp good, \resp boundaryless if and only if the $\P_{\hrt x/\widetilde k}$-space $\widetilde{(X,x)}$ is an open subset of $\P_{\hrt x/\widetilde k}$, \resp an affine open subset of $\P_{\hrt x/\widetilde k}$, \resp the whole of $\P_{\hrt x/\widetilde k}$. 

\begin{exem}\label{ex-temred}
Let us assume that $X=\mathscr M(k\{T\})$ and that $x$ is its Shilov points (in other words, $x=\eta_1$).
There is a $\widetilde k$-isomorphism $\widetilde k(\tau)\simeq \hrt x$ that sends $\tau$ to $\widetilde{T(x)}$.
Therefore 
$$\P_{\hrt x/\widetilde k}\simeq \P_{\widetilde k(\tau)/\widetilde k}=\P_{\widetilde k_1(\tau)/\widetilde k_1}=\P^1_{\widetilde k_1}$$
(the middle equality comes from Remark
\ref{rem-temkin-transfer} applied with $\Gamma=\{1\}$, together with
the fact that $\mathfrak D(\widetilde k(\tau))=\mathfrak D(\widetilde k)$). 
The reduction $\widetilde{(X,x)}$ is equal to the quasi-compact open subset $\P_{\hrt x/\widetilde k}\{\widetilde{T(x)}\}$ of $\P_{\hrt x/\widetilde k}$,
which is identified with $\A^1_{\widetilde k_1}$ through the above homeomorphism\footnote{
Be aware that this homeomorphism does not preserve the notion of an affine open subset: indeed,
as remarked above, the whole space $\P_{\hrt x/\widetilde k}$ is affine, though $\P^1_{\widetilde k^1}$
is not affine as a scheme! But this local terminology inconsistency should not cause any trouble in practice.} 
between $\P_{\hrt x/\widetilde k}$ and
$\P^1_{\widetilde k_1}$. 

Now let $Y$ be the $k$-analytic space obtained by gluing $\mathscr M(k\{T\})$ and $\mathscr M(k\{S\})$ along the isomorphism
$\mathscr M(k\{T,T\inv\})\simeq \mathscr M(k\{S,S\inv\})$ given by $S\mapsto T$. The Shilov points of $\mathscr M(k\{T\})$ and $\mathscr M(k\{S\})$ are identified and give rise to a single point $y$ on $Y$. By the above, there is a homeomorphism between the Zariski-Riemann space 
$\P_{\hrt y/\widetilde k}$ 
and $\P^1_{\widetilde k_1}$, modulo which Temkin's reduction $\widetilde{(Y,y)}$ is the {\em affine line with double origin}, viewed as a  $\P^1_{\widetilde k_1}$-space through the open immersion of each of the two copies of $\A^1_{\widetilde k_1}$ it contains
(by design). Hence $\widetilde{(Y,y)}\to \P_{\hrt y/\widetilde k}$ is not one-to-one, which witnesses the fact that $(Y,y)$ is not separated. 
\end{exem}

\subsection{Functoriality}\label{ss-temred-properties}
Let $L$ be an analytic extension of $k$, let $(X,x)$ be a $k$-analytic germ, and let $(Y,y)$ be an $L$-analytic germ. 
Any morphism $(Y,y)\to (X,x)$ gives rise in a natural way to a commutative diagram of topological spaces

$$
\xymatrix{
{\widetilde{(Y,y)}}\ar[r]\ar[d]&{\widetilde{(X,x)}}\ar[d]\\
{\P_{\hrt y/\widetilde k}}\ar[r]& {\P_{\hrt x/\widetilde k}}
}$$
in which the bottom arrow is the one induced by the extension $\hrt x\hookrightarrow \hrt y$.
Note that $\widetilde{(Y,y)}\to {\widetilde{(X,x)}}$ is quasi-compact
by \ref{ss-def-slk}. 
Let us now list some
very useful properties of this construction. 

\begin{enumerate}[1]
\item If $(X,x)$ is any analytic germ, then $(Y,x)\mapsto \widetilde{(Y,x)}$ induces a bijection between the set of analytic domains of $(X,x)$ and the set of quasi-compact, non-empty open subsets of $\widetilde{(X,x)}$; moreover, this bijection commutes
with finite unions and intersections.

\item If $k$ is an analytic field, a morphism $(Y,y)\to (X,x)$ of $k$-analytic germs is boundaryless if and only if the
local homeomorphism
\[\widetilde{(Y,y)}\to \P_{\hrt y/\widetilde k}\times_{\P_{\hrt x/\widetilde k}}\widetilde {(X,x)}\] is bijective (hence a homeomorphism). 

\item If $(X,x)$ is any analytic germ and if $(Y,x)$ is a closed analytic subgerm
of $(X,x)$, then $\widetilde{(Y,x)}\to \widetilde{(X,x)}$ is a homeomorphism. 

\item Let $X$ be a $k$-analytic space and let $Y$ be an $X$-analytic space. Let $L$ be an analytic extension of $k$, let $Z$ be an $L$-analytic space,
and let $Z\to X$ be a morphism of analytic spaces. Set $T=Y\times_X Z$, let $t$ be a point of $T$ and let $x, y$ and $z$ denote the images of $t$ in $X$, $Y$ and $Z$ respectively. Let us set for short

\begin{eqnarray*}
\mathsf X&=&\widetilde{(X,x)}\times_{\P_{\hrt x/\widetilde k}}\P_{\hrt t/\widetilde L}\\
\mathsf Y&=&\widetilde{(Y,y)}\times_{\P_{\hrt y/\widetilde k}}\P_{\hrt t/\widetilde L}\\
\mathsf Z&=&\widetilde{(Z,z)}\times_{\P_{\hrt z/\widetilde L}}\P_{\hrt t/\widetilde L}\end{eqnarray*}
The natural continuous $\P_{\hrt t/\widetilde L}$-map 
$\widetilde{(T,t)}\to \mathsf Y\times_{\mathsf X}\mathsf Z$ is then a homeomorphism.

\item Let $(Y,y)\to (X,x)$ be a morphism of analytic germs, let $(V,x)$ be an analytic domain of $(X,x)$,
and set $(W,y)=(Y,y)\times_{(X,x)}(V,x)$. The reduction $\widetilde{(W,y)}$ is equal 
to the pre-image of $\widetilde{(V,x)}$ in $\widetilde{(Y,y)}$.   

\end{enumerate}

\begin{enonce}[remark]{Remarks} Assertion (4) is stated by Temkin only when $L=k$: this is \Prop 4.6 of \cite{temkin2004}~But
in view of \Prop 3.1 of \opcit, its proof can be straightforwardly adapted to work for arbitrary $L$. 
Assertion (3) can be seen as a particular case of (2), but it can be checked directly from the definition after reduction to the affinoid case. Assertion (5) is a particular case of (4). 
\end{enonce}

\section{Temkin's reduction and $\Gamma$-strictness}\label{s-temkin-gamma}

\begin{lemm}\label{lem-gstrict-good}
Let~$(X,x)$ be a~$k$-analytic germ. The following are equivalent: 

\begin{enumerate}[i]
\item The point $x$ has a~$\Gamma$-strict~$k$-affinoid neighborhood in~$X$. 

\item The germ $\widetilde{(X,x)}$ is a $\Gamma$-strict affine open subset of $\P_{\hrt x/\widetilde k}$. 

\end{enumerate}

\end{lemm}

\begin{proof}
Let us assume that (i) holds, and let us choose a $\Gamma$-strict affinoid neighborhood $V$ of $x$, say $V=\mathscr M(A)$. By \ref{ss-gstr-ansp}, the spectral semi-norm on $A$ takes values in $(\abs{k\gpm}^\Q\cdot\Gamma)_0$; the image $B$ of the natural morphism $\widetilde A\to \hrt x$ is thus contained in $\hrt x^{\abs{k\gpm}^\Q\cdot\Gamma}$. This implies, in view of the equality
$\widetilde{(X,x)}=\P_{\hrt x/\widetilde k}\{S\}$ for any finite set $S$ of homogeneous generators of $B$, that $\widetilde{(X,x)}$ is affine and $\Gamma$-strict. 

Let us now assume that (ii) holds. We can then write $\widetilde{(X,x)}=\P_{\hrt x/\widetilde k}\{f_1,\ldots,f_n\}$ where each $f_i$ is a non-zero element
of $\hrt x^{r_i}$ for some $r_i$ in $\Gamma$.
Since $\widetilde{(X,x)}$ is an affine open subset of $\P_{\hrt x/\widetilde k}$, 
the germ $(X,x)$ is good (\ref{ss-first-prop-reduction});  in other words,~$x$ has an affinoid neighborhood~$V$ in~$X$. By shrinking~$V$ if needed,
we may and do
assume that there exist invertible analytic functions~$h_1,\ldots,h_n$ on~$V$ satisfying for every~$i$ the equalities~$|h_i(x)|=r_i$ and~$\widetilde{h_i(x)}=f_i$. Let $h \colon V\to \A^{n,\mathrm {an}}_k$ be the morphism induced by the $h_i$'s; set ~$t=h(x)$. Let~$W$ be the affinoid domain of~$\A^{n, \mathrm {an}}_k$ defined by the inequalities~$|T_i|\leq r_i$ for~$i=1,\ldots,n$ (where the~$T_i$'s are the coordinate functions on the affine space). 
Since the quasi-compact open subset~$\widetilde{(X,x)}=\P_{\hrt x/\widetilde k}\{f_1,\ldots, f_n\}$ of~$\P_{\hrt x/\widetilde k}$ is by construction the pre-image of~$\P_{\hrt t/\widetilde k}\{\widetilde{T_1(t)},\ldots,\widetilde {T_n(t)}\}=\widetilde{(W,t)}$,
it follows from \ref{ss-temred-properties} (5)
that $(X,x)\to (\A^{n, \rm an}_k,t)$ goes through~$(W,t)$. Hence we can shrink~$V$ so that there exist an affinoid neighborhood~$W'$ of~$t$ in~$\A^{n, \mathrm {an}}_k$ such that $h(V)$ is contained in $W\cap W'$. Since~$\A^{n,\mathrm {an}}_k$ has no boundary, we may and do assume that~$W'$ is~$\Gamma$-strict. As~$\widetilde{(V,x)}=\widetilde{(X,x)}$ is the pre-image of~$\widetilde{(W\cap W',t)}=\widetilde{(W,t)}$ inside~$\P_{\hrt x/\widetilde k}$, the morphism~$V\to W\cap W'$ is inner at~$x$ by \ref{ss-temred-properties} (2).
As~$W\cap W'$ is~$\Gamma$-strict, Lemma 2.5.11 of \cite{berkovich1990} immediately implies that~$x$ has a~$\Gamma$-strict affinoid neighborhood in~$V$, hence in~$X$.
\end{proof}

\begin{lemm}\label{lem-gstrict-general}
Let~$(X,x)$ be a~$k$-analytic germ. The following are equivalent: 

\begin{enumerate}[i]
\item The germ $(X,x)$ is $\Gamma$-strict. 

\item The reduction $\widetilde{(X,x)}$ is $\Gamma$-strict
 (see \ref{ss-def-dstrict}).

\end{enumerate}
\end{lemm}

\begin{proof}
The implication (i)$\Rightarrow$(ii) follows directly from Lemma~\ref{lem-gstrict-good}. Now let
us assume that $\widetilde{(X,x)}$ is $\Gamma$-strict, and let $(\mathsf U_i)$ 
be a $\Gamma$-strict atlas of $\widetilde{(X,x)}$.
For every $i$, let $(X_i,x)$ be the analytic domain of $(X,x)$ that corresponds to $\mathsf U_i$; for every $(i,j)$, the analytic
domain of $(X,x)$ that corresponds to $\mathsf U_i\cap \mathsf U_j$ is $(X_i\cap X_j,x)$. In order to prove that $(X,x)$ is $\Gamma$-strict, it is sufficient to prove that $(X_i,x)$ and $(X_i\cap X_j,x)$ are
$\Gamma$-strict for all $i,j$; hence we reduce to the case where $\widetilde{(X,x)}$ is a $\Gamma$-strict, non-empty, quasi-compact open subset of $\P_{\hrt x/\widetilde k}$. 

Under this assumption  $\widetilde{(X,x)}$ admits a finite covering $(\mathsf V_j)$ by $\Gamma$-strict
{\em affine}
open subsets. For every $j$, let $(V_j,x)$ denote the analytic domain of $(X,x)$ that corresponds to $\mathsf V_j$. 
Since the intersection of two $\Gamma$-strict affine open subsets of $\P_{\hrt x/\widetilde k}$ is still affine and $\Gamma$-strict, 
it follows from Lemma \ref{lem-gstrict-good}
that $(V_j,x)$ and $(V_j\cap V_\ell,x)$ 
are $\Gamma$-strict (and good) for all $j,\ell$. This implies that $(X,x)$ is $\Gamma$-strict, which ends the proof.
\end{proof}

\subsection{The $\Gamma$-graded reduction}\label{ss-graded-reduction}

Let $(X,x)$ be a $\Gamma$-strict $k$-analytic germ. Its reduction $\widetilde{(X,x)}$ is $\Gamma$-strict 
by Lemma~\ref{lem-gstrict-general} above; recall
that $\widetilde{(X,x)}^\Gamma$ then denotes the unique object of $\mathscr S_{\hrt x^\Gamma/\widetilde k^\Gamma}$ from which $\widetilde{(X,x)}$ arises.
If $\Delta$ is any subgroup of $\R_+\gpm$ containing $\Gamma$, then
$(X,x)$ is $\Delta$-strict as well and
$$\widetilde{(X,x)}^\Delta=\widetilde{(X,x)}^\Gamma\times_{\P_{\hrt x^\Gamma/\widetilde k^\Gamma}}\P_{\hrt x^\Delta/\widetilde k^\Delta}$$
(see \ref{ss-dstrict-ff} (4)). 

Let $L$ be an analytic extension of $k$, let $(Y,y)$ be a $\Gamma$-strict $L$-analytic germ, and let $(Y,y)\to (X,x)$ be a morphism of analytic germs. 
There is a unique continuous
map 
$\widetilde{(Y,y)}^\Gamma\to \widetilde{(X,x)}^\Gamma$ making the diagram 
$$\xymatrix{
&{\widetilde{(Y,y)}}\ar[rr]\ar'[d][dd]\ar[ld]&&{\widetilde{(X,x)}}\ar[dd]\ar[ld]\\
{\widetilde{(Y,y)}^\Gamma}\ar[rr]\ar[dd]&&{\widetilde{(X,x)}^\Gamma}\ar[dd]&\\
&{\P_{\hrt y/\widetilde k}}\ar'[r][rr]\ar[ld]&& {\P_{\hrt x/\widetilde k}}\ar[ld]\\
{\P_{\hrt y^\Gamma/\widetilde k^\Gamma}}\ar[rr]&& {\P_{\hrt x^\Gamma/\widetilde k^\Gamma}}
}$$
commute. 
Indeed, if $L=k$ this is a direct application of \ref{ss-dstrict-ff} (2); and if $Y=X_L$, this is a consequence of the
equalities $$\widetilde{(X_L,y)}=\widetilde {(X,x)}\times_{\P_{\hrt x/\widetilde k}}\P_{\hrt y/\widetilde L}$$ and 
$$\widetilde{(X_L,y)}^\Gamma=\widetilde {(X,x)}^\Gamma\times_{\P_{\hrt x^\Gamma/\widetilde k^\Gamma}}\P_{\hrt y^\Gamma/\widetilde L^\Gamma}$$
(the first one  comes from \ref{ss-temred-properties} (4), and it implies the second one).
The general case now follows formally by combining those two particular cases.

It follows from \ref{ss-dstrict-ff} (3) that a
quasi-compact open subset of $\widetilde{(X,x)}$ is $\Gamma$-strict if and only if
it is the pre-image of a quasi-compact open subset of $\widetilde{(X,x)}^\Gamma$. This implies
that a finite union or a finite intersection of $\Gamma$-strict quasi-compact
open subsets of $\widetilde{(X,x)}$ is a $\Gamma$-strict quasi-compact open subset, and that
the pre-image in $\widetilde{(Y,y)}$ of any $\Gamma$-strict
quasi-compact open subset of $\widetilde{(X,x)}$ is a $\Gamma$-strict
quasi-compact open subset
of $\widetilde{(Y,y)}$.

\subsection{Canonicity of the $\Gamma$-strict structure}\label{ss-gstrict-germ-can}
Let $(X,x)$ be a $\Gamma$-strict $k$-analytic germ. 
Let $(V_i,x)_i$ be a $\Gamma$-strict affinoid atlas defining a $\Gamma$-strict 
analytic structure on $(X,x)$, and let $(V,x)$ be a $\Gamma$-strict
analytic domain of $(X,x)$. By Lemma~\ref{lem-gstrict-general}, $\widetilde{(V,x)}$ and the $\widetilde{(V_i,x)}$'s
are $\Gamma$-strict. Therefore for every index $i$,
the intersection $\widetilde{(V,x)}\cap \widetilde{(V_i,x)}$ is $\Gamma$-strict by~\ref{ss-graded-reduction}
above, which implies again 
by Lemma~\ref{lem-gstrict-general}
that $(V,x)\cap (V_i,x)$ is $\Gamma$-strict. Since $(V_i,x)$ is separated, 
its $\Gamma$-strict structure is unique, hence $(V,x)\cap (V_i,x)$ is a $\Gamma$-strict analytic domain
of $(X,x)$ {\em for the given $\Gamma$-strict structure on $(X,x)$}. As this holds for every $i$, 
the analytic domain $(V,x)$ is a $\Gamma$-strict analytic domain
{\em for the given $\Gamma$-strict structure on $(X,x)$}. 

As a consequence, 
there exists a {\em unique}
$\Gamma$-strict $k$-analytic structure on $(X,x)$; the corresponding $\Gamma$-strict
analytic domains are simply the analytic domains $(V,x)$ that are $\Gamma$-strict in the sense
that was given up to now to this notion (\ie, analytic domains that admit {\em some}
$\Gamma$-strict analytic structure); 
by Lemma~\ref{lem-gstrict-general}, these are precisely
the analytic domains $(V,x)$ whose reduction $\widetilde{(V,x)}$ is $\Gamma$-strict. 

These facts immediately extend to the category of $k$-analytic spaces: if a $k$-analytic space
$Z$ admits a $\Gamma$-strict analytic structure, the latter is unique
and the corresponding $\Gamma$-strict
analytic domains are simply the analytic domains that admit a $\Gamma$-strict analytic structure. 

Therefore all (possible) ambiguities mentioned in Remark \ref{rem-gstrict-careful}
vanish, so now we can use the notion of $\Gamma$-strictness without worrying
about such unpleasant subtleties. 

\subsection{$\Gamma$-strictness is a local notion} \label{ss-gstrict-local}
Let us emphasize an important consequence of \ref{ss-gstrict-germ-can}: a $k$-analytic
space $X$ is $\Gamma$-strict if and only if the germ $(X,x)$ is $\Gamma$-strict
for every point $x$ of $X$. Indeed, the direct implication is obvious. Assume now that 
$(X,x)$ is $\Gamma$-strict for all $x\in X$. Then every point of $X$ admits a $\Gamma$-strict
analytic neighborhood, which can be chosen to be an {\em open} subset of $X$, because any open
subset of a $\Gamma$-strict analytic space is $\Gamma$-strict. Therefore $X$ can be covered
by $\Gamma$-strict open subsets. Since
the intersection of two $\Gamma$-strict
open subsets of $X$ is $\Gamma$-strict (again, because $\Gamma$-strictness is inherited
by open subsets of a $\Gamma$-strict space), $X$ is $\Gamma$-strict. 

\subsection{Fullness of the $\Gamma$-strict subcategory}\label{ss-gstrict-fullfaith}
Let $(X,x)$ be a $\Gamma$-strict $k$-analytic germ, let $L$
be an analytic extension, and let $(Y,y)$ be a $\Gamma$-strict $L$-analytic germ. 
Suppose that we are given a morphism $(Y,y)\to (X,x)$ 
of analytic germs. 

Let $(V,x)$ be a $\Gamma$-strict analytic domain 
of $(X,x)$ and let $(W,y)$ be the fiber product $(Y,y)\times_{(X,x)}(V,x)$.
By Lemma~\ref{lem-gstrict-general}
and~\ref{ss-graded-reduction}
above, the pre-image of $\widetilde{(V,x)}$ on $\widetilde{(Y,x)}$ is open, quasi-compact 
and $\Gamma$-strict. 
But this pre-image can be identified with $\widetilde{(W,y)}$ by \ref{ss-temred-properties} (5); hence $(W,y)$
is a $\Gamma$-strict analytic domain of $(Y,y)$ by Lemma ~\ref{lem-gstrict-general},
and $(Y,y)\to (X,x)$ is thus a morphism 
of $\Gamma$-strict analytic germs (note that here we have used implicitly \ref{ss-gstrict-germ-can}). 

We have thus proved that the category of $\Gamma$-strict analytic germs is a {\em full} subcategory of the category 
of analytic germs. As a consequence, the category of $\Gamma$-strict $k$-analytic germs (\resp analytic spaces, \resp $k$-analytic spaces)
is a full subcategory of the category of $k$-analytic germs (\resp analytic spaces, \resp $k$-analytic spaces). 

\subsection{Preservation of $\Gamma$-strictness under boundaryless pullback}\label{ss-pullback-gstrict}
Let $(Y,y)\to (X,x)$ be a boundaryless morphism of $k$-analytic germs. Assume that $(X,x)$ is $\Gamma$-strict. 
Its reduction $\widetilde{(X,x)}$ 
is then $\Gamma$-strict by Lemma \ref{lem-gstrict-general}. 
Since $(Y,y)\to (X,x)$ is boundaryless, 
$$\widetilde{(Y,y)}=\widetilde{(X,x)}\times_{\P_{\hrt x/\widetilde k}}\P_{\hrt y/\widetilde k}$$
by \ref{ss-temred-properties} (2), 
hence $\widetilde{(Y,y)}$ is $\Gamma$-strict by \ref{ss-dstrict-ff} (5). Using again 
Lemma \ref{lem-gstrict-general}, we see that $(Y,y)$ is $\Gamma$-strict. 

In view of \ref{ss-gstrict-local}, this immediately implies that if $Y\to X$ is a boundaryless
morphism
of $k$-analytic spaces, then $Y$ is $\Gamma$-strict as soon as $X$ is $\Gamma$-strict. 

\begin{rema}\label{rem-gstrict-fini-transfer}
By \ref{ss-pullback-gstrict}
above, if $Y\to X$ is a finite morphism between $k$-analytic spaces
and if $X$ is $\Gamma$-strict, then $Y$ is $\Gamma$-strict. This can also be seen
without using Temkin's reduction, as follows. First of all, it is sufficient to ensure that the pull-back
of a given $\Gamma$-strict affinoid atlas on $X$ is a $\Gamma$-strict affinoid atlas
on $Y$. Hence we reduce to the case where both $Y$ and $X$
are affinoid, say $Y=\mathscr M(B)$ and $X=\mathscr M(A)$. 

Now since $A$ is $\Gamma$-strict, there exists a polyradius $r=(r_1,\ldots, r_n)$
consisting of elements of $\Gamma$ such that $\abs{k_r\gpm}\neq \{1\}$ and 
$A_r$ is $k_r$-strict. The finite $A_r$-algebra $B_r$ is then $k_r$-strict as well. 
Therefore by \cite{bosch-g-r1984} 6.2.1/4, the spectral radius of every element of $B_r$ belongs to
$$\abs{k_r\gpm}^\Q\cup\{0\}\subset \cup (\abs {k\gpm}\cdot \Gamma)^\Q\cup\{0\}.$$
This holds in particular for every element of $B$, whence the $\Gamma$-strictness
of $B$ by \ref{ss-gstr-ansp}.
\end{rema}

\subsection{}\label{ss-prop-gammared}
It follows from \ref{ss-graded-reduction}
that the assignment $(X,x)\mapsto \widetilde{(X,x)}^{\Gamma}$ is functorial in
the $\Gamma$-strict analytic germ $(X,x)$. Using straightforward descent arguments
(based upon the surjectivity of the continuous map $\P_{\hrt x/\widetilde k}\to \P_{\hrt x^\Gamma/\widetilde k^\Gamma}$)
we deduce from
\ref{ss-first-prop-reduction}, \ref{ss-temred-properties} and Lemma \ref{lem-gstrict-good}
that it enjoys the following properties: 

\begin{enumerate}[1]\setcounter{enumi}{-1}
\item A $\Gamma$-strict $k$-analytic germ is separated, \resp good, \resp boundaryless if and only if the
$\P_{\hrt x^\Gamma/\widetilde k^\Gamma}$-space $\widetilde{(X,x)}^\Gamma$ is an open subset of $\P_{\hrt x^\Gamma/\widetilde k^\Gamma}$, \resp an affine open subset of $\P_{\hrt x^\Gamma/\widetilde k^\Gamma}$, \resp the whole of $\P_{\hrt x^\Gamma/\widetilde k^\Gamma}$.

\item If $(X,x)$ is any $\Gamma$-strict analytic germ, then $(V,x)\mapsto \widetilde{(V,x)}^\Gamma$ induces a bijection between the set of $\Gamma$-strict
analytic domains of $(X,x)$ and the set of quasi-compact, non-empty open subsets of $\widetilde{(X,x)}^\Gamma$; moreover, this bijection commutes
with finite unions and intersections.

\item If $k$ is an analytic field, a morphism $(Y,y)\to (X,x)$ of $\Gamma$-strict $k$-analytic germs is boundaryless if and only if the
continuous local homeomorphism
$$\widetilde{(Y,y)}^\Gamma\to \P_{\hrt y^\Gamma/\widetilde k^\Gamma}\times_{\P_{\hrt x^\Gamma/\widetilde k^\Gamma}}\widetilde {(X,x)}^\Gamma$$ is bijective (hence a homeomorphism). 

\item If $(X,x)$ is any $\Gamma$-strict
analytic germ and if $(Y,x)$ is a closed analytic subspace of $(X,x)$, then $\widetilde{(Y,x)}^\Gamma\to \widetilde{(X,x)}^\Gamma$ is a homeomorphism. 

\item Let $X$ be a $\Gamma$-strict $k$-analytic space, and let $Y$ be a $\Gamma$-strict $X$-analytic space.
Let $L$ be an analytic extension of $k$, let $Z$ be a $\Gamma$-strict $L$-analytic space and let $Z\to X$ be a morphism of analytic spaces.
Set $T=Y\times_X Z$, let $t$ be a point of $T$, and let $x, y$ and $z$ denote the images of $t$ in $X$, $Y$, and $Z$ respectively. Let us set for short

\begin{eqnarray*}
\mathsf X&=&\widetilde{(X,x)}^\Gamma\times_{\P_{\hrt x^\Gamma/\widetilde k^\Gamma}}\P_{\hrt t^\Gamma/\widetilde L^\Gamma}\\
\mathsf Y&=&\widetilde{(Y,y)}^\Gamma\times_{\P_{\hrt y^\Gamma/\widetilde k^\Gamma}}\P_{\hrt t^\Gamma/\widetilde L^\Gamma}\\
\mathsf Z&=&\widetilde{(Z,z)}^\Gamma\times_{\P_{\hrt z^\Gamma/\widetilde L^\Gamma}}\P_{\hrt t^\Gamma/\widetilde L^\Gamma}
\end{eqnarray*}
The natural continuous $\P_{\hrt t^\Gamma/\widetilde L^\Gamma}$-map 
$\widetilde{(T,t)}^\Gamma\to \mathsf Y\times_{\mathsf X}\mathsf Z$ is then a homeomorphism.

\item Let $(Y,y)\to (X,x)$ be a morphism of $\Gamma$-strict
analytic germs, let $(V,x)$ be a $\Gamma$-strict analytic domain of $(X,x)$,
and set $(W,y)=(Y,y)\times_{(X,x)}(V,x)$.
The reduction $\widetilde{(W,y)}^\Gamma$ is equal 
to the pre-image of $\widetilde{(V,x)}^\Gamma$ in $\widetilde{(Y,y)}^\Gamma$.

\end{enumerate}

\begin{rema}
Assume that $|k^\times|\neq \{1\}$. If $(X,x)$ is a strictly $k$-analytic space, we shall write
$\widetilde{(X,x)}^1$ instead of $\widetilde{(X,x)}^{\{1\}}$. The $\P_{\hrt x^1/\widetilde k^1}$-space
$\widetilde{(X,x)}^1$ is nothing but the {\em non-graded}
Temkin reduction of $(X,x)$ defined in \cite{temkin2000}, and all properties of $\Gamma$-strict analytic
spaces and $\Gamma$-graded reduction that we have just established were already known by Temkin's
works \cite{temkin2000} and \cite{temkin2004} in the case
$\Gamma=\{1\}$. What we have done here simply consists in explaining how Temkin's methods
actually apply to an arbitrary subgroup of $\R_+\gpm$. 
\end{rema}

\begin{rema}\label{rem-directdef-gred}
Let $(X,x)$ be a $\Gamma$-strict good germ, and let $V$ be a $\Gamma$-strict
affinoid neighborhood of $x$ in $X$, say $V=\mathscr M(A)$. Let $B$ denote the image
of $\widetilde A$ in $\hrt x$ through the evaluation map at $x$.

Let $f$ be any non-nilpotent
function belonging to $A$. Its spectral semi-norm belongs to $(|k^\times|\cdot \Gamma)^\Q$; hence there exists
$\lambda\in k^\times$ and an integer $n$ such that the spectral norm of $\lambda f^n$ belongs to $\Gamma$. 
If $\abs \cdot$ belongs to $\P_{\hrt x/\widetilde k}$, then $\abs{\widetilde{f(x)}}\leq 1$
if and only if $\abs{\widetilde{\lambda f^n(x)}}\leq 1$. As a consequence,
$$\widetilde{(X,x)}=\P_{\hrt x/\widetilde k}\{B\}=\P_{\hrt x/\widetilde k}\{B^\Gamma\}.$$ Since $B^\Gamma\subset \hrt x^\Gamma$, it follows
that $\widetilde{(X,x)}^\Gamma=\P_{\hrt x^\Gamma/\widetilde k^\Gamma}\{B^\Gamma\}.$ Note that $B^\Gamma$ can be described directly as the image 
of the finitely generated $\widetilde k^\Gamma$-graded algebra $\widetilde A^\Gamma$ inside $\hrt x^\Gamma$ through the evaluation map at $x$. 

Let $(Y,x)$ be a good $\Gamma$-strict analytic domain of $(X,x)$. It can be described by a finite conjunction of inequalities
$$\abs{f_1}\leq r_1\;{\rm and}\;\ldots\;{\rm and}\;\abs{f_n}\leq r_n$$ where the $f_i$'s are invertible function on $(X,x)$ and $r_i=\abs{f_i(x)}$
for every $i$. Now by $\Gamma$-strictness of $(Y,y)$, every $r_i$ appears as the spectral semi-norm
of $f_i$ on some $\Gamma$-strict affinoid neighborhood of $x$ in $Y$, hence belongs to $(|k^\times|\cdot \Gamma)^\Q$. 
Therefore by replacing every $f_i$ with $\lambda_i f_i^{n_i}$ for suitable $(\lambda_i,n_i)$ in $k\gpm\times \Z_{\geq 0}$ we may assume
that $r_i\in \Gamma$ for all $i$, and the equality $\widetilde{(Y,x)}=\widetilde{(X,x)}\cap \P_{\hrt x/\widetilde k}\{\widetilde{f_1(x)},\ldots, \widetilde{f_i(x)}\}$
then implies
that $\widetilde{(Y,x)}^\Gamma=\widetilde{(X,x)}^\Gamma\cap \P_{\hrt x^\Gamma/\widetilde k}\{\widetilde{f_1(x)},\ldots, \widetilde{f_n(x)}\}$. 

One could use the above to describe directly the $\Gamma$-graded reduction of a general $\Gamma$-strict germ, analogously
to what was done in~\ref{ss-defredgrad-good} and \ref{ss-defredgrad-general} for general graded reductions (\ie, $\R_+\gpm$-graded reductions),
by considering good charts and gluing. When $\Gamma=\{1\}$, this is essentially the way Temkin's ungraded reductions were built in \cite{temkin2000}.

\end{rema}

\begin{rema}\label{rem-topo-gstrict}
Let $X$ be a $\Gamma$-strict analytic space. If $X$ is quasi-compact, then it admits a finite G-covering
by $\Gamma$-strict affinoid domains. If $X$ is paracompact, it admits a locally finite covering by $\Gamma$-strict 
affinoid domains: the arguments are \emph{mutatis mutandi}
the same as in \ref{ss-topologies}.

\end{rema}

\chapter{Flatness in analytic geometry}\label{FLA}

In this chapter, we introduce one of the key notions
of this memoir, namely flatness. The definition
is given in Section \ref{s-def-flatness} (Definition \ref{def-flat-gen})
and as explained in the Introduction, it differs
from naive flatness even in the good case: one \emph{requires}
stability under arbitrary base change (including ground field extension). 
The basic properties are then stated, and some simple
examples
are given: for instance, 
any $k$-analytic space 
is flat over $k$ (Lemma \ref{lem-field-flat}; note that the stability
under base change requires some work). 

Section \ref{s-flatan-flatalg}
is devoted to GAGA results about flatness.
For instance, let $\mathscr Y\to \mathscr X$
be a morphism of schemes of finite
type over a given affinoid algebra, and let $y$
be a point of $\mathscr Y\an$. Then $\mathscr Y\an$
is flat over $\mathscr X\an$ at $y$
if and only if $\mathscr Y$ is flat over $\mathscr X$
at the image of $y$ (the ``only if" part
is the easiest one, see Lemma \ref{lem-gagaflat-easy}; but
the ``if" is still quite straightforward, see Proposition \ref{prop-gagaflat-ft}).
Let us now consider 
a morphism $Y\to X$ between $k$-affinoid spaces, and let $y$ denote
a point of $Y$. If $Y$ is flat over $X$
at $y$, then $\spec \mathscr O_Y(Y)$ is flat
over $\spec \mathscr O_X(X)$ at the image of $y$
(this is also covered by Lemma \ref{lem-gagaflat-easy}), 
but the converse is false in general: we give a counter-example
in 
\ref{ss-counterex-gagaflat}. Nevertheless, if  $\spec \mathscr O_Y(Y)$ is flat
over $\spec \mathscr O_X(X)$ at the image of $y$
and if
$y$ lies on a closed analytic subspace of $Y$ which is finite over $X$,
then $Y\to X$ is naively flat at $y$: this is Theorem \ref{th-gagaflat-hard}, whose
proof
rests on \emph{Crit\`eres locaux de platitude}, \cite{sga1} Expos\'e IV (we shall see
later that under these assumptions, $Y\to X$ is even flat at $y$; see Theorem
\ref{thm-flat-improve}). 

In Section \ref{s-flat-finite}, we
investigate finite flat
morphisms. The results we present there are
essentially due
to Berkovich
(see \cite{berkovich1993}, 3.2), but 
we include proofs (sometimes different from Berkovich's)
for the reader's convenience. Let us simply
mention here one of them:  finite flat
maps are open (Corollary \ref{cor-finflat-open}). 

In Section \ref{s-counter-example} we
present a counter-example
(suggested by Temkin, and described in \ref{ss-counterflat-intro}) showing that 
naive flatness is not preserved by base-change; note that the detailed
study of this counter-example
uses some results of section \ref{s-flat-finite}, which is the reason
why we have not carried it out immediately after having given the definition
of naive flatness. 

We end this Chapter by showing in Section \ref{s-flatness-expected}
that our notion of flatness behaves similarly
to flatness in algebraic 
geometry. For instance: usual algebraic properties
descend under flat maps (Lemma \ref{lem-flat-desc2}); flatness
can be checked after flat base change or arbitrary ground
field extension (\Prop \ref{prop-flatbc} and \ref{prop-flat-xl});
flatness has the expected properties as far as exactness
of complexes of coherent sheaves is concerned (Proposition
\ref{prop-flat-univexact}); it ensures that some properties
that hold inside a fiber can be spread out to the whole space
(Lemma \ref{lem-iso-fiber}); and it
implies the usual 
formula relating the local dimensions of the source space, 
of the target space, and of the fibers (Lemma \ref{lem-flatness-dim}).

\section{Naive and non-naive flatness}\label{s-def-flatness}

\subsection{} Let~$Y\to X$ be a morphism of {\em good}
analytic spaces, let~$y$ be a point of $Y$
and let~$x$ be its image on~$X$. Let~$\mathscr F$ be a coherent sheaf on~$Y$. One could
be tempted 
to say that~$\mathscr F$ is~$X$-flat at~$y$ if it is the case
in the framework
of locally ringed spaces; \ie, 
if~$\mathscr F_y$ is a flat~$\mathscr O_{X,x}$-module.

But we have chosen
to call the latter property {\em naive}\index{morphism!of analytic spaces!naively flat}\index{coherent sheaf!naively flat}\index{naively flat morphism}\index{naively flat coherent sheaf}
$X$-flatness of~$\mathscr F$ at~$y$, because it turns out that it
is not a reasonable candidate to be the analytic avatar of scheme-theoretic flatness.
Indeed, as we are going to see below
in Section \ref{s-counter-example}
through an explicit example,
{\em it is not stable under base change}. 
Let us now give the ``right" definition of flatness. 

\begin{defi}[Analytic flatness: the good case]\label{def-flat-good}
Let~$Y\to X$ be a morphism between good~$k$-analytic spaces, and let~${\mathscr F}$
be a coherent sheaf on~$Y$. Let~$y$ be a point of $Y$.
We shall say that~${\mathscr F}$ is {\em~$X$-flat at~$y$} if for any any good analytic space~$X'$,
for any morphism~$X'\to X$, and for any point~$y'$ lying above~$y$ on~$Y':=Y\times_{X}X'$,
the pull-back of~$\mathscr F$ on~$Y'$ is naively~$X'$-flat at~$y'$.
\end{defi}

\begin{rema}
We emphasize that in Definition~\ref{def-flat-good} above, the space $X'$
can be any good analytic space in the sense of \ref{ss-analytic-nofield}:  naive flatness
has to be checked after base change by a space defined over
an arbitrary analytic extension of $k$; in particular, it is required to hold after arbitrary ground field extension. 
\end{rema}

\begin{rema}
Theoretically, checking $X$-flatness of $\mathscr F$
at $y$ requires to consider all possible base-changes. But we shall see in fact later (Theorem \ref{thm-check-flatness}) that it if there exists an analytic extension $L$ of $k$
and an $L$-rigid point on $Y_L$ over $y$ at which $\mathscr F_L$ is naively $X_L$-flat, then $\mathscr F$ is $X$-flat 
at $y$. 

\end{rema}

\begin{exem}\label{ex-flat-dom}
Let $X$ be a good $k$-analytic space and let $V$ be a good analytic domain of $X$. 
For every $x\in V$, the coherent sheaf $\mathscr O_V$ is $X$-flat at $x$ by \ref{ss-algprop-analytic} (2). 

\end{exem}

For further investigation about flatness, we shall need the following
technical (but very easy) lemma. 

\begin{lemm}\label{lem-flat-formal}
Let
$$\xymatrix{
D&C\ar[l]\\B\ar[u]&A\ar[l]\ar[u]
}$$ be a commutative diagram of commutative rings such that~$C$, \resp $D$,
is flat over $A$, \resp faithfully flat over~$B$. 
If~$M$ is a~$B$-module such that~$D\otimes_BM$ is~$C$-flat, then~$M$ is~$A$-flat.
\end{lemm}

\begin{proof} Let~$N\hookrightarrow N'$ be an injective linear map between two~$A$-modules.
As~$C$ is~$A$-flat,~$C\otimes_AN\hookrightarrow C\otimes_AN'$. As~$M\otimes_B D$ is~$C$-flat,
$$\underbrace{(M\otimes_B D) \otimes_C (C\otimes_A N)}_{(M\otimes_B D) \otimes_A N}\hookrightarrow
\underbrace{(M\otimes_B D) \otimes_C (C\otimes_A N')}_{(M\otimes_B D) \otimes_A N'}.$$ In other words,
$(N\otimes_A M)\otimes_BD\hookrightarrow  (N'\otimes_A M)\otimes_BD$.
Faithful flatness of the~$B$-algebra~$D$ now implies that~$N\otimes_A M\hookrightarrow  N'\otimes_A M$.
\end{proof}

\subsection{}\label{flat-good-basics} Let $\phi \colon Y\to X$ be a morphism
of good $k$-analytic spaces, and let~$\mathscr F$ be a coherent sheaf on $Y$. 
Let $y$ be a point of $Y$. 

\begin{enumerate}[1]

\item If $\mathscr F$ is $X$-flat at $y$, it is in particular naively $X$-flat at $y$ (by definition).

\item Let~$V$  be a good analytic domain of~$Y$ containing $y$, 
and let~$U$ be a good analytic domain of~$X$ containing $\phi(Y)$. Let us consider the four following
properties: 

\begin{enumerate}[j]
\item $\mathscr F$ is naively $X$-flat at~$y$;

\item ${\mathscr F}_{|V}$ is naively $U$-flat at~$y$.

\item $\mathscr F$ is $X$-flat at~$y$;

\item ${\mathscr F}_{|V}$ $U$-flat at~$y$.
\end{enumerate}
It follows straightforwardly from
Example~\ref{ex-flat-dom} and from Lemma~\ref{lem-flat-formal}
that (i)$\iff$(ii). Applying this after an arbitrary good base change we see that
(iii)$\iff$(iv) as well. 
\end{enumerate}

\begin{defi}[Analytic flatness: the general case]\label{def-flat-gen}\index{morphism!of analytic spaces!flat}\index{coherent sheaf!flat}\index{flat morphism}\index{flat coherent sheaf}
Let~$\phi \colon Y\to X$ be a morphism of {\em non-necessarily good}~$k$-analytic spaces and let~$y\in Y$. Let~$\mathscr F$ be a coherent sheaf on~$Y$. 
It follows from~\ref{flat-good-basics} (2) that the following are equivalent: 

\begin{enumerate}[i]

\item For all pairs~$(V,U)$, where~$V$  is a good analytic domain of~$Y$ containing~$y$ and where~$U$ is a good analytic domain of~$X$ containing~$\phi(V)$, the coherent sheaf~${\mathscr F}_{|V}$ is~$U$-flat at~$y$. 

\item There exist a good analytic domain~$V$ of~$Y$ containing~$y$ and a good analytic domain~$U$ of~$X$ containing~$\phi(V)$ such that the coherent sheaf~${\mathscr F}_{|V}$ is~$U$-flat at~$y$. 

\item There exist an affinoid domain~$V$ of~$Y$ containing~$y$ and an affinoid domain~$U$ of~$X$ containing~$\phi(V)$ such that the coherent sheaf~${\mathscr F}_{|V}$ is~$U$-flat at~$y$. 

\end{enumerate}

We shall say that~$\mathscr F$ is {\em $X$-flat at~$y$} if it satisfies the equivalent assertions (i), (ii) and (iii) above (this definition is obviously equivalent to the previous one when~$Y$ and~$X$
are good). We shall say that~$\mathscr F$ is~$X$-flat if it is~$X$-flat at every point of~$Y$. We shall say that~$Y$ is~$X$-flat as~$y$ or~$X$-flat
if so is~$\mathscr O_Y$.
We shall say that $\mathscr F$ is flat at $y$, \resp flat, if it is $Y$-flat at $y$, \resp $Y$-flat (with respect to $\mathrm {Id}\colon Y\to Y$) 

\end{defi}

We are now going to state some basic facts (\ref{ss-flat-bc}--\ref{ss-flat-andom}), each of which can be proved by reduction to the good case
(which is allowed by the very definition of flatness). After such a reduction, the first two follow straightforwadly from the definition, 
the third from Example \ref{ex-flat-dom}, and the fourth from \ref{flat-good-basics} (2). We shall then investigate some situations 
in which flatness is expected, and actually holds -- but some work is needed to prove it.

\subsection{Stability under base-change}\label{ss-flat-bc} Let~$Y\to X$ be a morphism between~$k$-analytic spaces and let~$y$ be a point of $Y$. 
Let~$X'$ be an analytic space,  let~$X'\to X$ be a morphism, and let~$y'$ be a point of~$Y':=Y\times_XX'$ lying over~$y$.
If~$\mathscr F$ is a coherent sheaf on~$Y$ that is~$X$-flat at~$y$, its pull-back to $Y'$ is~$X'$-flat at~$y'$: this follows
directly from the definition, which {\em requires} that flatness be preserved under base change. 

\subsection{Stability under composition}\label{ss-flat-compo}  Let~$Z\to Y$ and~$Y\to X$ be morphisms between~$k$-analytic spaces, 
let~$z$ be a point of $Z$ and let~$y$ be its image on~$Y$.

\begin{enumerate}[1]

\item If~$\mathscr F$ is a coherent sheaf on~$Z$
that is~$Y$-flat at~$z$, and if~$Y$ is~$X$-flat at~$y$, then~$\mathscr F$ is~$X$-flat at~$Z$.

\item If $\mathscr G$
is a coherent sheaf on $Y$ that is $X$-flat at $y$, and if $Z$ is $Y$-flat at $z$, then $\mathscr G_Z$ is $X$-flat at $z$. 
\end{enumerate}

Indeed, both assertions are G-local, which allows to reduce to the good case and then to the corresponding
naive statements, which are
are obvious. 

\subsection{Flatness of analytic domains}\label{ss-flat-resdom} The inclusion of an analytic domain is flat: we see this again by reducing to the good case, 
and then to the corresponding naive statement, which is Example \ref{ex-flat-dom}.

\subsection{Good behavior by restriction to analytic domains}\label{ss-flat-andom}
Let~$Y\to X$ be a morphism of~$k$-analytic spaces, 
let~$V$ be an analytic domain of~$Y$ and let~$U$ be an analytic domain of~$X$ which contains the image of~$V$. 
Let~$\mathscr F$ be a coherent sheaf on~$Y$ and let~$y$ be a point of $V$. The coherent sheaf~$\mathscr F$ is~$X$-flat at~$y$ if and only if~$\mathscr F_V$ is~$U$-flat at~$y$: this follows once again by reducing to the good case and then to the corresponding naive statement, which is \ref{flat-good-basics} (2). 

\begin{lemm}\label{lem-field-flat}
Let~$Y$ be a~$k$-analytic space. The structure map~$Y\to \mathscr M(k)$ is flat.
\end{lemm}

\begin{proof}
We may and do assume that~$Y$ is~$k$-affinoid. 
Let~$X$ be an affinoid space, let~$y$ be a point of $Y\times_kX$, and let~$x$ be its image on~$X$. Let~$U$ be an affinoid neighborhood of~$x$ in~$X$, and let~$V$ be an affinoid neighborhood of~$y$ in~$Y\times_X U$. 
Let~$A$, $B$ and $C$ be the respective algebras of analytic functions on~$U$, $Y$ and $V$.
The~$A\hotimes_k B$-algebra~$C$ is flat (\ref{ss-algprop-analytic} (2)); the~$A$-algebra~$A\hotimes_k B$ is flat (\cite{berkovich1993},
Lemma 2.1.2; its statement
involves an
analytic extension $K$ of $k$, but its proof works for $K$ any Banach space over $k$, and we apply
it with $K=B$); hence~$C$ is~$A$-flat. By a straightforward limit argument,~$\mathscr O_{Y\times_k X,y}$ is a flat~$\mathscr O_{X,x}$-algebra, 
whence the lemma.
\end{proof}

\begin{lemm}\label{lem-cohsheaf-flat}
Let $X$ be a good analytic space and let $\mathscr F$ ba a coherent sheaf on $X$. Let $x$ be a point of $X$
at which $\mathscr F$ is naively flat; \ie, $\mathscr F_x$ is a flat $\mathscr O_{X,x}$-module. There exists 
an open neighborhood $U$ of $x$ in $X$ such that $\mathscr F_U$ is a free $\mathscr O_U$-module, 
and 
$\mathscr F$ is flat at $x$. 
\end{lemm}

\begin{proof}
Since $\mathscr F_x$ is flat over $\mathscr O_{X,x}$, it is a free- $\mathscr O_{X,x}$-module. This implies the existence of an open neighborhood
$U$ such that $\mathscr F_U$ is a free $\mathscr O_U$-module. Now for every good analytic space $Z$
and every morphism
$Z\to X$, the pre-image of $\mathscr F$ on $Z\times_X U$ is free, and thus naively flat
at every point of $Z\times_X U$, and in particular at every pre-image of $x$ on $Z$. The coherent sheaf
$\mathscr F$ is then flat at $x$. 
\end{proof}

\begin{lemm}\label{lem-cohsheaf-finflat}
\label{lem-stalk-finite}
Let~$Y\stackrel\pi{\longrightarrow}T\to X$
be a diagram of analytic spaces, with $\pi$ finite. Let $t$ be a point of $T$ and let $y_1,\ldots,y_r$ be the pre-images
of $t$ in $Y$. 
Let~$\mathscr F$ be a coherent sheaf on~$Y$.

\begin{enumerate}[1]
\item Assume that $Y$ and $T$ are good. The $\mathscr O_{T,t}$-module~$(\pi_\ast \mathscr F)_t$ is
then naturally isomorphic to~$\prod \mathscr F_{y_i}$.

\item Assume that $Y, T$ and $X$ are good, and let us consider the following 
assertions. 

\begin{enumerate}[j]
\item The coherent sheaf $\pi_\ast \mathscr F$ is naively $X$-flat at $t$. 
\item The coherent sheaf $\mathscr F$ naively $X$-flat at every $y_i$.
\item The coherent sheaf $\pi_\ast \mathscr F$ is $X$-flat at $t$. 
\item The coherent sheaf $\mathscr F$ is $X$-flat at every $y_i$. 
\end{enumerate}

We then have the equivalences
$$\text{\rm (i)}\iff\text{\rm (ii)}\;\;\;\text{and}\;\;\;\text{\rm(iii)}\iff\text{\rm(iv)}.$$ 

\item If we drop the goodness assumption, the equivalence
$$\mathrm{(iii)}\iff\mathrm{(iv)}$$
still holds. 
\end{enumerate}
\end{lemm}

\begin{proof}
Let us assume that $Y$ and $T$ are good.
The finite morphism $Y\to T$ is in particular
closed. Hence
for every neighborhood $V$ of~$\{y_1,\ldots,y_r\}$, there exist an affinoid neighborhood of~$t$ in~$T$
whose pre-image is included in~$V$ and is a {\em disjoint} union~$\coprod V_i$, where~$V_i$ is for every~$i$ an affinoid
neighborhood of~$y_i$ in~$Y$. This implies
that $(\pi_\ast \mathscr F)_t=\prod \mathscr F_{y_i}$, and (1) holds. 

Now let us prove (2). Let
let $x$ denote the image of $t$ in $X$. By the above,the $\mathscr O_{X,x}$-module
$(\pi_\ast \mathscr F)_t$ is flat if and only if $\mathscr F_{y_i}$ is $\mathscr O_{X,x}$-flat for every $i$, whence
the equivalence (i)$\iff$(ii). 

Assume that (iii) holds. Let $X'$ be a good analytic space and let $X'\to X$ be a morphism. Set $T'=T\times_X X'$ and $Y'=Y\times_X X'$.
Let $i$ be an element
of $\{1,\ldots, r\}$ and let $z$ be a pre-image of $y_i$ on $Y'$; let $t'$ and $x'$ denote the images of $z$ on $T'$ and $X'$,
and let $\pi'$ be the natural finite map $Y'\to T'$. 
Since $\pi_\ast \mathscr F$ is $X$-flat at
$t$, the coherent sheaf $(\pi_\ast \mathscr F)_{T'}$ is naively $X'$-flat at $t'$.
But $(\pi_\ast \mathscr F)_{T'}=\pi'_\ast (\mathscr F_{Y'})$ (to see it, one can assume that $T$ and $Y$ are affinoid, in which
case it is obvious by viewing coherent sheaves as modules). By the implication (i)$\iff$(ii) already proven, $\mathscr F_{Y'}$
is naively $X'$-flat at every primage of $t'$, and in particular at $z$; hence
the coherent sheaf $\mathscr F$ is $X$-flat at $y_i$, and (iv) holds. 

Assume conversely that (iv) holds. Let $X'$ be an affinoid space and let $X'\to X$ be a morphism. Set $T'=T\times_X X'$ and $Y'=Y\times_X X'$, 
and let $\pi'$ be the natural finite map $Y'\to T'$. Let $t'$ be a pre-image of $t$ on $X'$. If $z$ is a
pre-image of $t'$ on $Y'$, then the image of $z$ in $Y$ is equal 
to $y_i$ for some $i$. By assumption, $\mathscr F$ is universally $X$-flat at $y$. Hence $\mathscr F_{Y'}$ is naively $X'$-flat at $z$.
Since this holds for any such $z$, 
the implication (ii)$\iff$(i) already proved ensures that $\pi'_\ast (\mathscr F_{Y'})=(\pi_\ast \mathscr F)_{T'}$ is naively $X'$-flat at $t'$.
Hence $\pi_\ast \mathscr F$ is $X$-flat at $t$,
and (iii) holds.

Since flatness in the general case
can {\em by definition}
be checked on good analytic domains, (3) follows from (2). 
\end{proof}

\section{Algebraic flatness {\em versus} analytic flatness}\label{s-flatan-flatalg}

We shall first prove some GAGA principles for flatness ``in the easy direction" under very weak assumptions; we shall
then prove the converse implication in some particular (but nonetheless signifiant) cases. 

\begin{lemm}\label{lem-gagaflat-easy}
Let~$A\to B$ be a morphism between~$k$-affinoid algebras; let~$\mathscr Y$  (\resp $\mathscr X$)
be a
$B$-scheme of finite type (\resp an $A$-scheme of finite type). Let $\mathscr F$ be a cohrent sheaf on $\mathscr Y$, 
and let~$\mathscr Y\to \mathscr X$ be an~$A$-morphism. Let~$y$ be a point of $\mathscr Y\an$ at which 
$\mathscr F\an$ is naively~$\mathscr X\an$-flat. The coherent sheaf $\mathscr F$ is then~$\mathscr X$-flat at~$y\al$. 
\end{lemm}

\begin{proof}
Let $x$ be the image of $y$ on $\mathscr X\an$.
In the commutative diagram
$$\xymatrix{
{\mathscr O_{\mathscr Y\an, y}}&{\mathscr O_{\mathscr X\an,x}}\ar[l]\\
{\mathscr O_{\mathscr Y, y\al}}\ar[u]&{\mathscr O_{\mathscr X,x\al}}\ar[u]\ar[l]
}$$
the vertical arrows are faithfully flat (\ref{ss-gaga-alg}),
and naive flatness of $\mathscr F\an$ at $y$ means that the $\mathscr O_{X\an,x}$-module $\mathscr F\an _y$ is flat. 
Lemma \ref{lem-flat-formal}
then ensures that $\mathscr F_{y\al}$ is flat over $\mathscr O_{\mathscr X,x\al}$; \ie, $\mathscr F$
is flat at $y\al$. 
\end{proof}

We are now interested in the converse
of Lemma~\ref{lem-gagaflat-easy}. We are first going to mention in
\ref{ss-gagaflat-finite}
and \ref{ss-gagaflat-rigid} two special cases 
in which it is more or less well-known. Then we shall generalize  \ref{ss-gagaflat-finite}
to any morphism between schemes of finite type over a {\em given}
affinoid algebra (Proposition \ref{prop-gagaflat-ft}), and hence prove (Theorem \ref{th-gagaflat-hard})
a GAGA principle
for a morphism between affinoid spaces that extends
both \ref{ss-gagaflat-finite}
and \ref{ss-gagaflat-rigid}. 

\subsection{}\label{ss-gagaflat-finite} Let~$Y\to X$ be a {\em finite} morphism between affinoid spaces,
let $y$ be a point  of $Y$ and let $\mathscr F$ is a coherent sheaf on~$Y$. 
The coherent sheaf $\mathscr F$ is
naively~$X$-flat at~$y$ if and only if $\mathscr F\al$ is $X\al$-flat at $y\al$: this is
essentially \Prop 3.2.1 of \cite{berkovich1993} -- the latter is written only for~$\mathscr F=\mathscr O_Y$,
but using Lemma \ref{lem-cohsheaf-finflat} one can easily easily adapt it so that it works for any coherent sheaf. 

\subsection{}\label{ss-gagaflat-rigid} Let~$Y\to X$ be a
morphism between affinoid spaces, let~$y$ be a rigid
point of $Y$, and let $x$ be its image on $X$; let $\mathscr F$ be a coherent sheaf on $Y$.
We claim that $\mathscr F$ is
naively~$X$-flat at~$y$ if and only if $\mathscr F\al$ is~$X\al$-flat
at~$y\al$. If~$\abs{k\gpm}\neq\{1\}$ and if~$Y$ and~$X$ are strictly~$k$-affinoid, this is a classical assertion
of rigid-analytic geometry, but its proof is very simple and immediately extends to our situation: indeed, one
knows from \cite{sga1}, Expos\'e IV \Cor 5.8 that~$\mathscr F_y$ is flat over~$\mathscr O_{X,x}$ if and only if
its completion $\widehat{\mathscr F_y}=\mathscr F_y\otimes_{\mathscr O_{Y,y}}\widehat{\mathscr O_{Y,y}}$ is
flat over~$\widehat{\mathscr O_{X,x}}$, and that~$\mathscr F\al_{y\al}$ is flat over~$\mathscr O_{X\al,x\al}$
if and only if~$\widehat{\mathscr F\al_{y\al}}=\mathscr F_{y\al}\al\otimes_{\mathscr O_{Y\al ,y\al}}\widehat{\mathscr O_{Y\al,y\al}}$
is flat over~$\widehat{\mathscr O_{X,x\al}}$; but as~$x$ and~$y$ are
rigid,~$ \widehat {\mathscr O_{X,x}}=\widehat {\mathscr O_{X\al,x\al}}$ and~$ \widehat {\mathscr O_{Y,y}}=\widehat {\mathscr O_{Y\al,y\al}}$ (\cite{berkovich1993}, 
Lemma 2.6.3), whence our claim.

\begin{prop}\label{prop-gagaflat-ft}
Let~$\mathscr Y\to \mathscr X$ be a morphism between schemes of finite type over a given affinoid algebra. Let~$\mathscr F$ be a coherent sheaf on~$\mathscr Y$
and let~$y$ be a point of~$\mathscr Y\an$. Assume that $\mathscr F$ is~$\mathscr X$-flat at~$y\al$. The coherent sheaf $\mathscr F\an$ is
then~$\mathscr X\an$-flat at~$y$.
\end{prop}

\begin{proof}
We can assume that~$\mathscr X$ is affine.
Let us first prove that the coherent sheaf~$\mathscr F\an$ is {\em naively}
flat
at~$y$.
Let~$x$ be the image of~$y$ on~$\mathscr X\an$ and let~$U$ be an
affinoid neighborhood of~$x$ in~$\mathscr X\an$.  There  is a natural map from~$\mathscr O_{\mathscr X}(\mathscr X)$ to~$\mathscr O_U(U)$
which induces a morphism~$U\al\to \mathscr X$, and the space~$\mathscr Y\an\times_{\mathscr X\an}U~$ can be identified with the analytification
$(\mathscr Y\times_{\mathscr X}U\al)\an$ of the $U\al$-scheme of finite type $\mathscr Y\times_{\mathscr X}U\al$.
Since flatness is preserved by any {\em scheme-theoretic} base change,
the pull-back of $\mathscr F$ on $(\mathscr Y\times_{\mathscr X}U\al)$ is $U\al$-flat
at the image of $y$. This implies, in view of the fact that $\mathscr Y\an \times_{\mathscr X\an}U\to \mathscr Y\times_{\mathscr X}U\al$
is flat as a morphism
of locally ringed spaces (\ref{ss-gaga-alg}), that $\mathscr F\an_y$ is a flat $\mathscr O_{U\al, x_U\al}$-module. Since $\mathscr O_{X,x}$
is the direct limit of the local rings $\mathscr O_{U\al, x\al_U}$ for $U$ running through the set of affinoid neighborhoods of $x$ in $\mathscr X\an$,
we conclude that  $\mathscr F\an_y$ is a
flat $\mathscr O_{\mathscr X\an,x}$-module; otherwise said, $\mathscr F\an$ is naively $\mathscr X\an$-flat at $y$. 

Let us now prove that $\mathscr F\an$ is flat at~$y$.
Let~$V$ be an affinoid space, let~$V\to \mathscr X\an$ be a morphism, and let~$z$ be a point of~$\mathscr Y\an\times_{\mathscr X\an}V$
lying above~$y$.
There  is a natural map~$\mathscr O_{\mathscr X}(\mathscr X)\to \mathscr O_V(V)$ which induces a morphism
$V\al\to \mathscr X$, and~$\mathscr Y\an\times_{\mathscr X\an}V~$ can be identified with
the analytification $(\mathscr Y\times_{\mathscr X}V\al)\an$ of the $V\al$-scheme of finite type $\mathscr Y\times_{\mathscr X}V\al$.
As flatness behaves well under {\em scheme-theoretic} base change, the pull-back of~$\mathscr F$ on~$\mathscr Y\times_{\mathscr X}V\al$
is~$V$-flat at the image of $z$. It follows then from the naive version of the proposition (which we have proved above)
that the pull-back of~$\mathscr F\an$ on~$\mathscr Y\an\times_{\mathscr X\an}V$ is naively~$V$-flat at~$z$, which ends the proof.
\end{proof}

\begin{theo}\label{th-gagaflat-hard}
Let $Y\to X$ be a morphism between~$k$-affinoid spaces, let $\mathscr F$ be a coherent sheaf on $Y$ and let $y$ be a point of $Y$.
Assume that there exists a closed analytic
subspace~$Z$ of~$Y$ containing~$y$ such that~$Z\to X$ is finite,
and that $\mathscr F\al$ is $X\al$-flat
at $y\al$. The coherent sheaf $\mathscr F$ is then naively $X$-flat at $y$. 
\end{theo}

\begin{rema}
In Theorem \ref{thm-flat-improve}
we shall in fact prove that $\mathscr F$ is actually $X$-flat at $y$. 
\end{rema}

\begin{proof}[Proof of Theorem \ref{th-gagaflat-hard}]
We denote by $A$, \resp $B$, the algebra of analytic functions on $X$, \resp $Y$; and we 
set $M=\mathscr F(Y)$. Let $x$ be the image of $y$ in $X$. Our purpose is now
to prove that the $\mathscr O_{X,x}$-module $M\otimes_B\mathscr O_{Y,y}$ is flat. 
Let~$V$ be an affinoid neighborhood of~$x$ in~$X$ and set~$W=Y\times_X V$. We denote by~$A_V$, \resp $B_W$, the algebra
of analytic functions on~$V$, \resp $W$.
Let~$\mathfrak p$ the prime ideal of~$ A$ that corresponds to~$x\al$, and let~$I$ be the ideal of~$B$ that corresponds to~$Z$.

The core of the proof consists in showing that ${\mathscr O}_{V\al, x_V\al}$-module $M\otimes_B \mathscr O_{W\al, y_W\al}$ is flat.
By \cite{sga1}, Expos\'e IV, \Th 5.6 the latter is true if and only if the two following conditions are satisfied:

\begin{enumerate}[A]
\item $M\otimes_B \mathscr O_{W\al,y_W\al}/\mathfrak{p}$ is~${\mathscr O}_{V\al,x\al_V}/\mathfrak{p}$-flat.

\item For every~$d>0$, the natural map~$$M\otimes_B(\mathfrak p^d \mathscr O_{W\al,y\al_W}/\mathfrak{p}^{d+1}\mathscr O_{W\al, y_W\al})\to
\mathfrak p^d (M\otimes_B \mathscr O_{W\al,y_W\al})/\mathfrak{p}^{d+1}(M\otimes_B \mathscr O_{W\al,y_W\al})$$ is an isomorphism.

\end{enumerate}

We first prove (A). For that purpose, let us begin with a general remark. Let $\Lambda$ be an arbitrary $B/I$-module.
As~$\mathscr O_{X\al,x\al}/\mathfrak p$ is a field, $\Lambda\otimes_A\mathscr O_{X\al,x\al}/\mathfrak p$ is a
flat $\mathscr O_{X\al,x\al}/\mathfrak p$-module. It follows
that~$\Lambda\otimes_A\mathscr O_{V\al,x_V\al}/\mathfrak p$ is a flat~${\mathscr O}_{V\al,x_V\al}/\mathfrak p$-module, 
which on the other hand can be rewritten as

$$\Lambda \otimes_{B/I}(B/I)\otimes_A A_V\otimes_{A_V} {\mathscr O}_{V\al,x\al}/\mathfrak p
=\Lambda\otimes_{B/I}(B_W/I)\otimes_{A_V} {\mathscr O}_{V\al,x_V\al}/\mathfrak p, 
$$
where the equality comes from the fact that
$B_W/ I=(B/I)\hotimes_A A_V=(B/I)\otimes_A A_V$
by finiteness of $B/I$ over $A$. Since ${\mathscr O}_{W\al,y_W\al}/(\mathfrak p+I)$ is a localization
of~$(B_W/I)\otimes_{A_V} \mathscr O_{V\al,x_V\al}$, 
the  ${\mathscr O}_{V\al,x_V\al}/\mathfrak p$-module
$$\Lambda\otimes_B {\mathscr O}_{W\al,y_W\al}/\mathfrak p=\Lambda\otimes_{B/I} {\mathscr O}_{W\al,y_W\al}/(\mathfrak p+I)$$
is flat as well.

Set $N=M/\mathfrak pM$, and let~$n$ be an element of $\Z_{>0}$.  By applying
the above to the
$B/I$-module $\Lambda=I^nN/I^{n+1}N$, we see that
$$(I^nN/I^{n+1}N)
\otimes_B {\mathscr O}_{W\al,y_W\al}=(I^nN/I^{n+1}N)
\otimes_B {\mathscr O}_{W\al,y_W\al}/\mathfrak p$$ is a flat  ${\mathscr O}_{V\al,x_V\al}/\mathfrak p$-module. 
As~$ {\mathscr O}_{W\al, y_W\al}$ is a flat~$B$-algebra (indeed, it is a localization of $B_W$ which is
itself $B$-flat;
see \ref{ss-algprop-analytic} (2)),  
we have the equality $$(I ^nN/I^{n+1} N)\otimes_B\mathscr O_{W\al, y_W\al}=
 I^n(N\otimes_B  \mathscr O_{W\al, y_w\al})/I^{n+1}(N\otimes_B  \mathscr O_{W\al,y_W\al}).$$
The~${\mathscr O}_{V\al,x_V\al}/\mathfrak p$-module
$I^n(N\otimes_B  \mathscr O_{W\al,y_W\al})/I^{n+1}(N\otimes_B  {\mathscr O}_{W\al,y_W\al})$
is thus flat for any non-negative~$n$. It obviously implies that for any
such~$n$, the~${\mathscr O}_{V\al,x_V\al}/\mathfrak p$-module
$$(N\otimes_B\mathscr O_{W\al,y_W\al})/I^{n+1}(N\otimes_ B\mathscr O_{W\al, y_W\al})$$ is flat. By \cite{ega31},
Chapter 0, \S 10.2.6,~$M\otimes_B \mathscr O_{W\al,y_W\al}/\mathfrak p=N\otimes_B\mathscr O_{W\al,y_W\al}$ is
then~$\mathscr O_{V\al,x_V\al}/\mathfrak p$-flat; hence~(A) is true.

Let us now prove (B). Let~$d$ be a positive integer. By assumption, $\mathscr F\al$ is $X\al$-flat at $y\al$,
which means that~$M\otimes_B {\mathscr O}_{Y\al, y\al}$ is~${\mathscr O}_{X\al, x\al}$-flat. 
Therefore the natural map 
$$(M\otimes_B {\mathscr O}_{Y\al, y\al})\otimes_{\mathscr O_{X\al, x\al}} \mathfrak p^d\mathscr O_{X\al, x\al}
\to \mathfrak p^d M\otimes_B {\mathscr O}_{Y\al, y\al}$$ is an isomorphism. But it can be written as the composition 
$$(M\otimes_B {\mathscr O}_{Y\al, y\al})\otimes_{\mathscr O_{X\al, x\al}} \mathfrak p^d\mathscr O_{X\al, x\al}\to
M\otimes_B \mathfrak p^d \mathscr O_{Y\al, y\al}\to \mathfrak p^d M\otimes_B {\mathscr O}_{Y\al, y\al}.$$
The left arrow being surjective and the composition of the two arrows being injective, the right arrow is injective.
As it is also
surjective, it is an isomorphism. It follows that
$$M\otimes_B(\mathfrak p^d \mathscr O_{Y\al, y\al}/\mathfrak{p}^{d+1}\mathscr O_{Y\al,y\al})
\to \mathfrak p^d (M\otimes_B \mathscr O_{Y\al, y\al})/\mathfrak{p}^{d+1}(M\otimes_B \mathscr O_{Y\al, y\al})$$
is an isomorphism as well. As~${\mathscr O}_{W\al, y_W\al}$ is a flat~${\mathscr O}_{Y\al, y\al}$-algebra
by \ref{ss-algprop-analytic} (2), assertion (B) follows by tensoring with $\mathscr O_{W\al,y_W\al}$ over $\mathscr O_{Y\al, y\al}$.
Therefore $M\otimes_B\mathscr O_{W\al, y_W\al}$ is flat over $\mathscr O_{V\al, x_V\al}$, as announced. 

Now let~$T$ be any affinoid neighborhood of~$y$ in~$W$. As~${\mathscr O}_{T\al, y\al_T}$ is a
flat~${\mathscr O}_{W\al, y_W\al}$-algebra by \ref{ss-algprop-analytic} (2),
the~${\mathscr O}_{V\al,v\al}$-module~$M\otimes_B {\mathscr O}_{T\al, y_T\al}$ is flat. 

We have thus shown the following: if~$V$ is any neighborhood of~$x$ in~$X$
and if~$T$ is any affinoid neighborhood of~$y$ in the pre-image of~$V$ inside~$Y$,
then~$M\otimes_B {\mathscr O}_{T,y_T\al}$ is~${\mathscr O}_{V\al,x_V\al}$-flat.
A straightforward limit argument then ensures that~$M\otimes_B {\mathscr O}_{Y,y}$ is~${\mathscr O}_{X,x}$-flat. 
\end{proof}

\begin{rema}
In the above proof, the existence of a closed analytic
subspace of $Y$ containing $y$
and finite over $X$ was used only while proving assertion (A), and $X\al$-flatness of $\mathscr F\al$ 
at $y\al$ was used only while proving assertion (B). 
\end{rema}

\section{The flat, (locally) finite morphisms}\label{s-flat-finite}
We shall use what we have just done to show some results which were already proven in \cite{berkovich1993},
\S 3.2 when~$\mathscr F=\mathscr O_Y$; we include the proofs (which 
partially differ of those of \cite{berkovich1993}) for the convenience of the reader.

\begin{prop}\label{prop-finite-flat}
Let~$Y\to X$ be a morphism between good~$k$-analytic spaces and let~$y$ be a point of~$Y$
at which this morphism is finite; let~$x$ be the image of~$y$ on~$X$. Let~$\mathscr F$ be a coherent sheaf on~$Y$. The following are equivalent:

\begin{enumerate}[i]
\item  The coherent sheaf $\mathscr F$ is naively $X$-flat at~$y$.

\item There exist an affinoid neighborhood~$T$ of~$y$ in~$Y$ and an affinoid neighborhood
$S$ of~$x$ in~$X$ such that~$T\to X$ goes through a finite map~$\pi : T\to S$ and such that~$\pi_\ast \mathscr F_{|T}$ is flat at~$x$. 

\item There exist an affinoid neighborhood~$T$ of~$y$ in~$Y$ and an affinoid
neighborhood~$S$ of~$x$ in~$X$ such that~$T\to X$ goes through a finite map~$\pi : T\to S$
and such that~$\pi_\ast (\mathscr F_T)$ is a free~$\mathscr O_S$-module. 

\item There exist an affinoid neighborhood~$T$ of~$y$ in~$Y$ and an affinoid neighborhood~$S$
of~$x$ in~$X$ such that~$T\to X$ goes through a finite map~$\pi : T\to S$ and such
that~$\mathscr F(T)$ is a flat~$\mathscr O_S(S)$-module. 

\item There exist an affinoid neighborhood~$T$ of~$y$ in~$Y$ and an affinoid
neighborhood~$S$ of~$x$ in~$X$ such that~$T\to X$ goes through a finite map~$\pi : T\to S$ and such that
$(\mathscr F_T)\al$ is $S\al$-flat at $y_T\al$. 

\item The coherent sheaf $\mathscr F$ is $X$-flat at $y$. 

\end{enumerate}
\end{prop}

\begin{proof}
Suppose that (i) is true. As~$Y\to X$ is finite at~$y$, there exist an affinoid neighborhood~$T$
of~$y$ in~$Y$ and an affinoid neighborhood~$S$ of~$x$ in~$X$ such that~$T\to X$ goes
through a finite map~$\pi : T\to S$ for which~$y$ is the only pre-image of~$x$. As~$\mathscr F$
is naively $X$-flat at~$y$, Lemma~\ref{lem-cohsheaf-finflat} applied to the diagram $T\to S\stackrel{\mathrm{Id}}{\longrightarrow}S$
yields the naive flatness of~$\pi_\ast (\mathscr F_T)$ at~$x$; by Lemma \ref{lem-cohsheaf-flat}, $\pi_\ast (\mathscr F_T)$
is then flat at $x$, whence (ii).
If~(ii) is true, then by Lemma~\ref{lem-cohsheaf-flat}
we can shrink $S$ (and $T$) so that $\pi_\ast (\mathscr F_T)$ is a free $\mathscr O_S$-module,
whence (iii). If (iii) holds then $\mathscr F(T)$ is a free, thus flat, $\mathscr O(S)$-module, whence (iv).
The
implication (iv)$\Rightarrow$(v) is obvious. Now if (v) holds, Proposition \ref{prop-gagaflat-ft} ensures that $\mathscr F_T$ is~$S$-flat $y$,
whence (vi). Implication (vi)$\Rightarrow$(i) is tautological. 
\end{proof}

\begin{coro}\label{cor-finflat-open}
Let~$Y\to X$ be a morphism of $k$-analytic spaces,
and let $\mathscr F$ be a coherent sheaf on $Y$. 
Let $y$ be a point of $\supp F$ at which $\supp F\to X$ is finite and at which $\mathscr F$
is $X$-flat, and let $x$ be the image of $y$
in $X$. The image of~$\supp F$ on~$X$ is then a neighborhood of~$x$
in $X$.
\end{coro}

\begin{proof} By arguing G-locally on $X$, shrinking $Y$ around $y$, and replacing $Y$
with $\mathrm{Supp}(\mathscr F)$ we may assume that both $Y$
and $X$ are good and $Y=\mathrm{Supp}(\mathscr F)$, so $Y\to X$ is finite. 
Let us now choose~$T$ and~$S$ as in 
(iii) above. Since $\mathscr F_y\neq 0$, it follows from Lemma \ref{lem-cohsheaf-finflat} (1)
that $(\pi_\ast (\mathscr F_T))_x\neq 0$; hence the free~$\mathscr O_S$-module~$\pi_\ast (\mathscr F_T)$
has positive rank. It follows that $(\pi_\ast (\mathscr F_T))_z\neq 0$ for every point $z$ of $S$; using again
Lemma \ref{lem-cohsheaf-finflat} (1),
we conclude that $S$ is contained in the image of $\supp F$. 
\end{proof}

\begin{coro}
\label{cor-finflat-dim}
Let~$Y\to X$ be a morphism of $k$-analytic spaces, let~$y\in Y$ be
such that~$Y\to X$ is finite at~$y$ and let~$x$ be the image of~$y$ in $Y$.
Let~$\mathscr F$ be a coherent sheaf on~$Y$. If~$y$ belongs to $\supp F$ and
$\mathscr F$ is $X$-flat at~$y$, then~$\dim_x X=\dim_y \mathscr F$
(\ref{ss-interpret-support}). 
\end{coro}

\begin{proof}
By arguing $G$-locally on $X$ and shrinking $Y$ around $y$ we may assume that both $Y$ and $X$
are affinoid, that
$Y\to X$ is finite, and that $y$ is the only pre-image
of~$x$ in~$Y$. The image of~$\supp F$  in~$X$ is then a Zariski-closed subset~$T$ of~$X$,
and one has~$\dim_xT=\dim_y \supp F=\dim_y \mathscr F$
(the first equality comes from \ref{ss-dimloc-morfin}
and the second one holds by definition
of $\dim_y \mathscr F$); on the other hand,~$T$ is a neighborhood of~$x$ in~$X$ by Corollary~\ref{cor-finflat-open}
above, hence~$\dim_xT=\dim_x X$.
\end{proof}
 
\section{Naive flatness is not preserved by base change}\label{s-counter-example}

As we have already mentioned, the reason why we have introduced (in the good case)
a sophisticated definition 
for flatness instead of dealing with naive flatness is the fact that the latter is not preserved by base change.
The purpose
of this section is to discuss in full detail a counter-example  that was initially suggested by Michael Temkin.

\subsection{}\label{ss-typic-nonflat}
Before introducing the counter-example, 
let us describe a situation that we shall encouter several times, in which we can conclude
that a given morphism between good spaces is {\em not}
naively flat at a given point; the key argument will be the violation of the dimension 
equality provided by Corollary \ref{cor-finflat-dim}. 

Let $\phi \colon Y\to \Omega$ be a morphism between good $k$-analytic spaces; we assume that $\phi$
factorizes through a closed immersion $Y\hookrightarrow X$ for $X$ a good analytic domain 
of $\Omega$ and that $Y$ and $\Omega$ are respectively of pure dimension $d$ and $d'$ with $d<d'$
(in all specific examples considered below, $d=1$ and $d'=2$). Let $y$ be a point of $Y$ at which 
$\phi$ is inner; this can be the case for instance if $X=\Omega$ (because then $\phi$ is boundaryless) or
if $y$ is rigid (because $y$ is then an inner point of $Y$). Since $Y\hookrightarrow X$ is boundaryless, 
$\mathrm{Int}(\phi)$ is equal to $\phi\inv(\mathrm{Int}(X/\Omega))$; as a consequence $\phi$ induces a 
closed immersion $\mathrm{Int}(\phi)\hookrightarrow \mathrm{Int}(X/\Omega)$, and
since $\mathrm{Int}(X/\Omega)$
is an open subset of $\Omega$ (this is the {\em topological} interior
of $X$ inside $\Omega$), the morphism $\phi$
is finite at every point of $\mathrm{Int}(\phi)$. It is particular finite at $y$, and
since $\dim_y Y=d$ and $\dim_{\phi(y)}\Omega=d'>d$, 
it follows from Corollary \ref{cor-finflat-dim} that $\phi$ is not naively flat at $y$ (strictly speaking, it
only follows from Corollary \ref{cor-finflat-dim}
that
$\phi$ is not flat at $y$; but as $\phi$ is finite at $y$, it is flat at $y$ if and only if it is
naively flat at $y$, by \Prop \ref{prop-finite-flat}). 

\subsection{Presentation of the counter-example}\label{ss-counterflat-intro}
Let~$r$ be a positive real number and let~$f=\sum \alpha_iT^i$ be a power series with coefficients in~$k$
whose radius
of convergence
is exactly $r$; \ie, 
$$\xymatrix{
{\abs{\alpha_i}r^i}\ar[r]_{i\to+\infty}&0}$$ and
$(\abs{\alpha_i}s^i)_i$ is non-bounded as soon as~$s>r$. We
denote by~$p : \A^{2, \mathrm {an}}_k\to \A^{1, \mathrm {an}}_k$
the first projection. Let~$X$ be the analytic domain of~$\A^{2, \mathrm{an}}_k$ defined by the
inequality~$|T_1|\leq  r$, and let~$Y$ be
the one-dimensional closed disc of radius~$r$; note that~$X=p\inv (Y)$; \ie, $X$ can
be identified with~$Y\times_k \A^{1, \mathrm {an}}_k$. The map~$\phi:=({\rm Id},f)$
from~$Y$ to~$Y\times_k \A^{1,\mathrm{an}}_k$ induces 
a closed immersion $Y\hookrightarrow X$, and more precisely an isomorphism
between~$Y$ and the Zariski-closed subspace~$Z$ of~$X$ defined by the sheaf
of ideals $(T_2-f(T_1))\mathscr O_X$; the inverse isomorphism is nothing but~$p|_Z$.  Set~$x=\phi(\eta_r)$.
We are going to show that $\phi$ is naively flat at $\eta_r$, but not flat; \ie, naive flatness of $\phi$
at $\eta_r$ does not hold universally. 

\subsection{}\label{ss-notflat-easy} The easiest part of our study is the negative one; \ie, the fact that $\phi$
is {\em not}
flat at $y$. Let us give two examples of base-change functors that witness it. 

\begin{enumerate}[1]
\item The base-change of the morphism $\phi$ by the inclusion $X\hookrightarrow \A^{2,\mathrm{an}}_k$ 
is the closed immersion $Y\hookrightarrow X$, and it follows from the general discussion in \ref{ss-typic-nonflat}
that $Y\hookrightarrow X$ is not naively flat at $\eta_r$. 

\item Let $L$ be any analytic extension of $k$ such that $\eta_r$ has an $L$-rational
pre-image $y$ in $X_L$; \eg, $L=\hr {\eta_r}$. The morphism $\phi_L\colon Y_L\to \A^{2,\mathrm{an}}_L$
is then inner at $y$, and it follows again from the general discussion in \ref{ss-typic-nonflat}
that $\phi_L$ is not naively flat at $y$. 
\end{enumerate}

\subsection{} We are now going to prove that $\phi$ is naively flat at $\eta_r$. We shall
in fact prove the following stronger result: {\em the local ring
$\mathscr O_{\A^{2,\mathrm{an}}_k,x}$ is a field}.

Before giving the rigorous proof (see Proposition \ref{prop-oanx-field}
below), let us roughly explain
what is going on. 
We want to prove that any analytic
function defined in a neighborhood of $x$ and that vanishes at $x$ vanishes
around $x$; or otherwise said, that 
any Zariski-closed subgerm $(W,x)$ of $(\A^{2,\mathrm{an}}_k,x)$ is equal to the whole
of $(\A^{2,\mathrm{an}}_k,x)$. 
So let us consider such a $(W,x)$. The point $x$ is not rigid, so the
dimension of $(W,x)$ is $\geq 1$. Assume that it is equal to 1. 
Using again the fact that $x$ is not rigid, we see that the one-dimensional
Zariski-closed subgerms $(X\cap W,x)$ and $(Z,x)$ of $(X,x)$ 
coincide; hence by gluing the germ $(W,x)$ to the curve $Z$ (whose boundary is $\{x\}$), we
can in some sense extend $Z$ beyond $x$ in $\A^{2,\mathrm{an}}$, and we get a
contradiction with
the fact that the radius of convergence
of $f$ is exactly $r$. Therefore $\dim (W,x)=2$ and  $(W,x)=(\A^{2,\mathrm{an}}_k,x)$.

\begin{lemm}\label{lem-onedim-field}
Let~$T$ be
a
reduced one-dimensional good analytic space and let~$t$ be a non-rigid point of $T$. 
The local ring~$\mathscr O_{T,t}$ is a field.
\end{lemm}

\begin{proof}
By Corollary \ref{cor-interp-centdim}, 
$\mathrm{centdim} (T,t)+\dim_{\mathrm{Krull}}\mathscr O_{T,t}=\dim_t T$. As~$t$ is not
rigid,~$\mathrm{centdim} (T,t)>0$; as~$T$ is one-dimensional,~$\dim_ t T\leq 1$. 
Therefore~$\dim_{\mathrm{Krull}}\mathscr O_{T,t}=0$; being reduced,~$\mathscr O_{T,t}$ is thus a field.
\end{proof}

\begin{prop}\label{prop-oanx-field}
The local ring~$\mathscr O_{\A^{2, \mathrm {an}}_k, x}$ is a field.
\end{prop}

\begin{proof}
As the analytic space~$\A^{2, \mathrm {an}}_k$ is reduced, it
is sufficient to prove that
$\dim_{\mathrm{Krull}}\mathscr O_{\A^{2,\mathrm{an}}_k, x}=0$, which is equivalent, in view of Corollary
\ref{cor-interp-centdim}, to the fact that~$\mathrm{centdim}(\A^{2,\mathrm{an}}_k, x)=2$.
Since~$x$ is not a rigid point (because one has by construction~$\hr x=\hr {\eta_r}$, and because
$d_k(\eta_r)=1$ by \ref{ss-intro-etar}, 
$\mathrm{\rm centdim}(\A^{2,\mathrm{an}}_k, x)>0$. 

{\em Assume that~$\mathrm{centdim}(\A^{2,\mathrm{an}}_k, x)=1$}. Then there exists an affinoid neighborhood~$V$ of~$x$ 
in~$\A^{2,\mathrm {an}}_k$
and an irreducible one-dimensional Zariski-closed subset~$W$ of~$V$ that contains~$x$. 
Both~$W\cap X=W\cap (V\cap X)$ and~$Z\cap V=Z\cap (V\cap X)$ are purely
one-dimensional Zariski-closed subsets of~$V\cap X$ containing~$x$.
As~$x$ is not a rigid point, it belongs to a unique irreducible component of
$(W\cap X)\cup(Z\cap V)$. Therefore it belongs to a unique irreducible
component of~$W\cap X$ and to a unique irreducible component of~$Z\cap V$, so those two components coincide.
One can hence shrink~$V$ so that~$W\cap X=Z\cap V$; by endowing~$W$ with its reduced structure, this equality becomes
an equality of closed analytic subspaces of~$X\cap V$.

If~$w$ is any point of~$W$ such that~$p(w)=\eta_r$, then in view of the equality
$d_k(\eta_r)=1$, 
the inequality~$d_k(w)\leq 1$ (due to the fact that~$W$ is one-dimensional)
forces~$d_{\hr {\eta_r}}(w)$ to be equal to zero.
Therefore~$(p|_W)\inv (\eta_r)$ is zero-dimensional. In particular $p|_W$ is
quasi-finite at $x$; moreover, since~$x$ belongs to the topological interior of~$V$ in~$\A^{2,\mathrm{an}}_k$,
the map~$p|_W$ is inner at~$x$; \Prop 3.1.4 of \cite{berkovich1993} then ensures that~$p|_W$
is finite at~$x$; as~$\mathscr O_{\A^{1,\mathrm{an}}_k,\eta_r}$ is a field by Lemma \ref{lem-onedim-field},~$p|_W$ is 
naively flat at~$x$. It follows then from Proposition~\ref{prop-finite-flat}
that there exists an affinoid neighborhood $W_0$
of $x$ in $W$ and an affinoid neighborhood~$U$ of~$\eta_r$ in~$\A^{1,\mathrm{an}}_k$ such that~$p(W_0)\subset U$, 
and such that~$p|_{W_0}: W_0\to U$ is finite and makes~$\mathscr O_{W_0}(W_0)$ a free~$\mathscr O_U(U)$-module
of finite positive rank, say~$r$ (note that we thus have~$p(W_0)=U$). By restricting to~$Y$
and using the fact that $X=p\inv (Y)$, one sees that~$p|_{W_0\cap X}:W_0\cap X\to U\cap Y$
is finite and makes~$\mathscr O_{W_0\cap X}(W_0\cap X)$ a free~$\mathscr O_{U\cap Y}(U\cap Y)$-module of rank~$r$;
we thus have~$p(W_0\cap X)=U\cap Y$.

It follows from the inclusion $W_0\cap X\subset W\cap X=Z\cap V$ that $W_0\cap X$ is an analytic domain of $Z$.
Since $p|_Z$ induces an isomorphism~$Z\simeq Y$ whose
inverse isomorphism is induced by $\phi$, the morphism $p|_{W_0\cap X}$ induces an isomorphism from $W_0\cap X$
 to $p(W_0\cap X)=U\cap Y$ whose inverse is~$\phi_{|U\cap Y}$. 
 As~$\mathscr O_{W_0\cap X}(W_0\cap X)$ is a free~$\mathscr O_{U\cap Y}(U\cap Y)$-module
 of rank~$r$, we have~$r=1$, which means that $p$ induces an isomorphism~$W_0\simeq U$.
 The inverse isomorphism defines a section~$\sigma$ of the first projection~$U\times_k\A^{1,\mathrm{an}}_k\to U$;
 we have ~$\sigma_{|U\cap Y}=\phi_{|U\cap Y}$. We can thus glue~$\sigma$ and~$\phi$ to obtain
 a section of the first projection~$(U\cup Y)\times_k\A^{1,\mathrm{an}}_k\to (U\cup Y)$ that
 coincides with~$\phi$ on~$Y$\; \ie, an analytic function~$g$ on~$U\cup Y$ that coincides
 with~$f$ on~$Y$. As~$U$ is a neighborhood of~$\eta_r$ in~$\A^{1,\mathrm{an}}_k$,
the analytic domain~$U\cup Y$ of~$\A^{1,\mathrm{an}}_k$ contains a closed disc centered
at the origin and whose radius is $>r$; but on the other hand the radius of
convergence of $f$
is exactly $r$, hence $f$ does not extend to any disc of radius $>r$, contradiction.
As a consequence,~$\mathrm{centdim}(\A^{2,\mathrm{an}}_k, x)\neq 1$ and
thus $\mathrm{centdim}(\A^{2,\mathrm{an}}_k, x)=2.$
 \end{proof}

\begin{rema}
The above counter-example rests on a boundary phenomenon. In fact, it turns out that
such phenomena are the only obstructions for naive flatness to be preserved by base-change; indeed, we
shall see later that naive flatness implies flatness in the boundaryless
case (Theorem \ref{thm-flat-naiveflat}). 
\end{rema}

\begin{rema}\label{rem-centdim-counter}
Let $\rho $ be an element of $(0,r)$ and set $z=\phi(\eta(\rho))$. We have $\hr z=\hr {\eta_{\rho}}$,
and $d_k(\eta_\rho)=1$ by \ref{ss-shilov-section}. Therefore
$\d_k(z)=1$; since $z$ lies on the one-dimensional irreducible Zariski-closed subset $Z$ of $X$, it follows that
$\adhz {\{z\}}X=Z$ (\ref{rem-abh-point}). Hence $\adht {\{z\}}{\A^{2,\mathrm{an}}_{k,\mathrm{Zar}}}$ 
contains $x$, and is thus of dimension 2
since $\mathrm{centdim}(\A^{2,\mathrm{an}}_k, x)=2$, as seen in the proof of Proposition \ref{prop-oanx-field}.
This implies that $\adht {\{z\}}{\A^{2,\mathrm{an}}_{k,\mathrm{Zar}}}$ is the whole of $\A^{2,\mathrm{an}}_k$ 
because the latter is irreducible; see
Proposition \ref{prop-gaga-irrcomp}
or simply note that $\A^{2,\mathrm{an}}_k$ is non-empty, connected, and normal. Note that the same reasoning
would more generally show that $\adhz {\{z\}}D=D$ for any irreducible analytic domain $D$ of $\A^{2,\mathrm{an}}_k$
containing $z$ and such that $\mathrm{centdim}(D,x)=2$; \eg, $D$ is a neighborhood of $x$. 

Choose $s\in (\rho, r)$. Let $X'$ be the open subset of $\A^{2,\mathrm{an}}_k$ defined by the inequality $\abs{T_1}<s$.
The intersection $Z\cap X'$ is a
closed analytic subspace of $X'$ containing $z$ which is isomorphic through $p_1$ to the open disc of radius $s$,
and is thus one-dimensional. 
Since $d_k(z)=1$, it follows from Example \ref{ex-centdim-abhyankar}
that $\mathrm{centdim}(X', z)=1$. We thus have 

\[\mathrm{centdim}(\A^{2,\mathrm{an}}_k, z)=\mathrm{centdim}(X', z)=1<2=\dim \adht {\{z\}}{\A^{2,\mathrm{an}}_{k,\mathrm{Zar}}}.\]

\end{rema}

We are now going to explain how the above construction also provides a counter-example to 
general GAGA-principle for naive flatness, and another one to stability of scheme-theoretic flatness under analytic
ground field extension.

Let 
$D$ be a closed two-dimensional polydisc centered at the origin such that $Z$ is contained in the corresponding open polydisc. 
Let $\psi$ denote the morphism $Y\to D$ induced by $\phi$. Note that
$\mathscr O_{D,x}=\mathscr O_{\A^{2,\mathrm{an}},x}$; since the latter is a field, $\psi$ is naively flat at $\eta_r$. Note also that $\psi$ factorizes
by construction through a closed immersion
$Y\hookrightarrow D\cap X$. 

\subsection{}\label{ss-counterex-gagaflat} Since $\eta_r$ is a norm, it lies above the generic point $\xi$
of $Y\al$. As~$\psi$ is naively flat at~$\eta_r$, it follows from Lemma \ref{lem-gagaflat-easy}
that the induced map~$Y\al\to D\al$ is flat at~$\xi$. 
Now let $y$ be a non-rigid point of $Y\setminus \{\eta_r\}$ (\eg, $y=\eta_{r'}$ for some $r'\in (0,r)$). The point $y$ does not
lie on any proper Zariski-closed subset of $Y$, which means that  $y\al =\xi$. But since~$y\in {\rm Int}(Y/k)$,
it follows from the general discussion in \ref{ss-typic-nonflat}
that $\psi$ is not naively flat at $y$.

\subsection{}\label{ss-counterex-algflat}
Let us now assume that~$r\notin \abs{k\gpm}^\Q$. In this case,~$\{y\}$ is an affinoid domain of~$Y$ (defined by the equality~$|T|=r$)
whose corresponding~$k$-affinoid algebra is nothing but~$k_r$. In view of Remark \ref{rem-centdim-counter}, 
the Zariski-closure $\adhz{\{\psi(y)\}}D$ is equal to the whole of $D$, which implies
that $\psi(y)$ lies above $\xi$. As a consequence, the morphism $\spec k_r \to D\al$
induced by $\psi|_{\{y\}}$ is flat. 

Let~$L$ be any analytic extension of~$k$ such that~$r$ belongs to
the group $\abs{L\gpm}^\Q$ (e.g.,~$L=k_r$). The space~$\mathscr M(L\hotimes_kk_r)$ is strictly~$L$-affinoid
and non-empty; it has thus an~$L$-rigid point, say~$t$. Since $t$ is rigid, it belongs to $\mathrm{Int}(\mathscr M(L\hotimes_kk_r))$, 
and since $\psi|_{\{y\}}$ factorizes through a closed immersion from $\{y\}$
to the affinoid domain of $D$ defined by the equality $|T_1|=r$, it follows
again from the general discussion in \ref{ss-typic-nonflat}
that $\mathscr M(L\hotimes_kk_r)\to D_L$ is not naively flat at $t$. We deduce then
from Theorem \ref{th-gagaflat-hard} (for a direct and simpler proof,
see \ref{ss-gagaflat-rigid})
that~$\spec (L\hotimes_kk_r)\to D_L\al$ is not flat at $t\al$.

\section{Analytic flatness has the expected
properties}\label{s-flatness-expected}

In this
section, we show that flatness in our sense behaves reasonably; \ie, the analogues
of classical results from algebraic geometry hold in our setting. We begin with the descent
of algebraic properties. We shall first write a statement that holds 
in the abstract settings
of \ref{s-category-framework} and \ref{s-alg-properties},
where we deal with general objects and properties of the latter satisfying various axioms; and then
we shall then write down what it means for some {\em explicit} properties of interest. 
For the notion of validity of a property at a point, the reader may refer
to Lemma-Definition \ref{lem-equiv-valid} in our general 
abstract setting and to Lemma-Definition
\ref{valid-at-concrete}
for a more
concrete version.

\begin{lemm}\label{lem-flat-desc1}
Let $Y\to X$ be a morphism of analytic spaces, let $y$ be a point of $Y$ at which $Y$ is $X$-flat, 
and let $x$ be the image of $y$ in $X$. Let $\mathfrak F$ be a fibered category as in \ref{s-category-framework}
and let $\mathsf P$ be a property as in \ref{ss-alg-properties}. Let $D$ be an object of $\mathfrak F_X$.

\begin{enumerate}[1]

\item If $\mathsf P$ satisfies condition \hv~ of \ref{ss-list-hreg} and if $D$ satisfies $\mathsf P$ at $x$, then $D_Y$ 
satisfies $\mathsf P$ at $y$. 

\item  If $\mathsf P$ satisfies condition \hreg~ of \ref{ss-list-hreg} and if $D_Y$ satisfies $\mathsf P$ at $y$, then $D$ 
satisfies $\mathsf P$ at $x$. 
\end{enumerate}

\end{lemm} 

\begin{proof}
We may and do assume that $Y$ and $X$ are good. Being $X$-flat at $y$, the space $Y$ is in particular naively $X$-flat
at $y$, which means that $\mathscr O_{Y,y}$ is a flat $\mathscr O_{X,x}$-algebra. The lemma follows then immediately from 
the definitions of conditions \hv~and \hreg. 
\end{proof}

\begin{lemm}[A concrete version of Lemma \ref{lem-flat-desc1}]
\label{lem-flat-desc2}
Let~$Y\to X$ be a morphism between~$k$-analytic spaces, let~$y$ be a point of $X$
at which $Y$ is $X$-flat, and let
$x$ be its image on~$X$. Let $\mathscr F$
be a coherent sheaf on~$X$, and let $m$ be an element of $\N$. Let $\mathsf S=(\mathscr E'\to \mathscr E\to \mathscr E'')$
be a short
complex of coherent sheaves on $X$.

\begin{enumerate}[1]
\item If~$Y$ is regular, \resp $R_m$, \resp Gorenstein, \resp CI, at~$y$, 
so is~$X$ at~$x$. If~$\mathscr F_Y$ is CM, \resp $S_m$, \resp free of rank $m$, at~$y$, so is~$\mathscr F$ at~$x$. 
If $\mathsf S_Y$ is exact at $y$, so is $\mathsf S$ at $x$. 

\item If~$\mathscr F$ is free of rank $m$ at $x$, so is $\mathscr F_Y$ at $y$. 
If $\mathsf S$ is exact at $x$, so is $\mathsf S_Y$ at $y$. 
\end{enumerate}
\end{lemm}

We are now going to prove that flatness can be checked after 
flat base change (Proposition \ref{prop-flatbc})
and ground field extension (Proposition \ref{prop-flat-xl}). 

\begin{lemm}\label{lem-flatbc}
Let~$$\xymatrix{
Z\ar[r]\ar[d]&T\ar[d]\\
Y\ar[r]&X
}$$ be a commutative diagram of good analytic spaces, let~$z$ be a point of $Z$ and let~$t$, \resp $y$,
be its image in~$T$, \resp $Y$. Let~$\mathscr F$ be a coherent sheaf on~$Y$
and let~$\mathscr G$ be its pull-back on~$Z$. Suppose that~$T$ is~
naively~$X$-flat at~$t$ and that~$Z$ is
naively~$Y$-flat at~$z$. If~$\mathscr G$ is
naively~$T$-flat at~$z$ then~$\mathscr F$ is
naively~$X$-flat at~$y$.
\end{lemm}

\begin{proof}
This follows straightforwardly from Lemma \ref{lem-flat-formal}.
\end{proof}

\begin{rema}
We emphasize that in Lemma \ref{lem-flatbc}
above, the spaces involved are {\em not}
assumed to be $k$-analytic; indeed, we want typically to apply it to diagrams arising
from ground field extension.
\end{rema}

\begin{prop}\label{prop-flatbc}
Let
$$\xymatrix{
Z\ar[r]\ar[d]&T\ar[d]\\Y\ar[r]&X
}$$ be a commutative diagram of~$k$-analytic spaces, let~$z$ be a point of $Z$ and let~$t$, \resp $y$,
be its image in~$T$, \resp $Y$. Let~$\mathscr F$ be a coherent sheaf on~$Y$.
Suppose that~$T$ is~$X$-flat at~$t$, that~$Z$ is ~$Y$-flat at~$z$ and
that~$\mathscr F_Z$ is~$T$-flat at~$z$;  under those assumptions~$\mathscr F$ is~$X$-flat at~$y$.
\end{prop}

\begin{proof}
One immediately reduces to the case where all spaces are affinoid. Let~$X'$
be a good analytic space and let~$X'\to X$ be a morphism. We set~$Y'=Y\times_X X'$
and so on. Let~$y'$ be a point on~$Y'$ lying above~$y$, and let~$z'$ be a point on~$Z'$
lying above both~$z$ and~$y'$ (such a point always exists by Lemma \ref{lem-pointlying}); 
denote by~$t'$ the image of~$z'$ on~$T'$. Since~$\mathscr F_Z$ is~$T$-flat at~$z$, the sheaf~$\mathscr F_{Z'}$ is
naively~$T'$-flat at~$z'$. Since~$T$, \resp $Z$, is~$X$-flat at~$t$, \resp $Y$-flat at~$z$,
the space $T'$, \resp $Z'$,  is
naively~$X'$-flat at~$t'$, \resp naively~$Y'$-flat at~$z'$. Lemma \ref{lem-flatbc}  above now implies that~$\mathscr F_{Y'}$ is
naively~$X'$-flat at~$y'$.
\end{proof}

\begin{prop}\label{prop-flat-xl}
Let~$Y\to X$ be a morphism of~$k$-analytic spaces and let~$L$ be an analytic extension of~$k$. Let~$y\in Y$
and let~$\mathscr F$ be a coherent sheaf on~$Y$. Let~$u$ be a point of~$Y_L$ lying above~$y$.
Suppose that $\mathscr F_L$ is~$X_L$-flat at~$u$;  the coherent sheaf~$\mathscr F$ is then~$X$-flat at~$y$.
\end{prop}

\begin{proof}
One can assume that both~$Y$ and~$X$ are affinoid.
Let~$X'$ be a good~$F$-analytic space for some analytic extension~$F$ of~$k$. We set~$Y'=Y\times_X X'$.
Let~$y'$ be a point on~$Y'$ lying above~$y$; we
are going to show that $\mathscr F_{Y'}$ is
naively~$X'$-flat at~$y$; by shrinking~$X'$, one can assume that it is~$F$-affinoid. 

\medskip
By Lemma \ref{lem-pointlying} there exists an analytic extension~$K$ of both~$F$
and~$L$ and a point~$\omega$ on~$Y'_K:=Y_K\times_{X_K}X'_K$ lying above both~$u$ and~$y'$.
By
$X_L$-flatness of~$\mathscr F_L$ at~$u$ the coherent sheaf~$\mathscr F_{Y'_K}$ is
naively~$X'_K$-flat at~$\omega$. Applying  Lemma \ref{lem-flatbc}
above to the diagram
$$\xymatrix{
{Y'_K}\ar[r]\ar[d] &{X'_K}\ar[d] \\{Y'}\ar[r]&X'
}$$ (which is possible in view of \ref{ss-algprop-analytic} (3)), 
one immediately gets the
naive~$X'$-flatness of~$\mathscr F_{Y'}$ at~$y'$.
\end{proof}

Proposition \ref{prop-flat-univexact}
below describes some consequences of flatness
on the homology of complexes of coherent
sheaves
(for the notion of exactness, injectivity, bijectivity, \etc\@~at
a given point, see 
Lemma-Definition \ref{valid-at-concrete}). Assertions (1)
and (2) are analogues of well-known
results in scheme theory.
Assertions (3) and (4)
are stated and proved 
for further use in the study of the loci of validity
of various properties (Chapter \ref{LOC}).

\begin{prop}
\label{prop-flat-univexact}
Let~$Y\to X$ be a morphism between~$k$-analytic spaces, let~$y$ be a point of $Y$ and let~$x$ be its image in~$X$.
Let $L$ be an analytic extension of $k$, let $X'$ be an $L$-analytic space and let $X'\to X$
be a morphism of analytic spaces.
Let~$y'$ be a pre-image of~$y$ on~$Y':=Y\times_XX'$ and let~$x'$ denote
the image of~$y'$ in~$X'$. Let~$\mathscr H$ be a coherent sheaf on~$X'$.

\begin{enumerate}[1]

\item Let~$\mathsf S=(\mathscr E\to  \mathscr E'\to \mathscr E'')$ be a 
sequence of coherent sheaves on~$Y$ which is exact at~$y$. 
If $\mathscr H$ is~$X_L$-flat at~$x'$, the sequence
$\mathsf S\boxtimes \mathscr H$
is exact at $y'$.

\item
Let~$0\to \mathscr G\to  \mathscr F\to \mathscr E\to 0$ be a sequence of coherent
sheaves on~$Y$ which is exact at~$y$; moreover, assume
that~$\mathscr E$ is~$X$-flat at~$y$. The sequence~$$0\to \mathscr G\boxtimes \mathscr H\to  \mathscr F\boxtimes \mathscr H
\to \mathscr E\boxtimes \mathscr H\to 0$$ of coherent sheaves on $Y'$ is then exact at~$y'$.

\item Let~$\mathsf S=(\mathscr F_n\to \mathscr F_{n-1}\to \ldots \to \mathscr F_0\to 0)$ be a sequence of coherent sheaves on~$Y$. 
Assume that all the~$\mathscr F_i$'s are~$X$-flat
at $y$, 
and that~$\mathsf S$ is exact at~$y$. The 
sequence~$\mathsf S\boxtimes \mathscr H$ on $Y'$
is then exact at~$y'$, and if~$n\geq 1$
the kernel of~$\mathscr F_n\to \mathscr F_{n-1}$ is~$X$-flat at~$y$. 

\item Let~$n$ be a positive integer
and let~$\mathsf S=(\mathscr F_n\to \mathscr F_{n-1}\to \ldots \to \mathscr F_0\to 0)$ 
be a complex of coherent sheaves on~$Y$.
Asssume that~$\mathscr F_i$ is ~$X$-flat at~$y$
for every~$i\leq n-1$, and that~$\mathsf S$ is exact 
at~$y$ except possibly at $\mathscr F_{n-1}$. 
The natural map
$$\mathrm H_{n-1}(\mathsf S)\boxtimes \mathscr H \to \mathrm H_{n-1}(\mathsf S \boxtimes \mathscr H)$$ of coherent
sheaves on $Y'$ is an isomorphism
at~$y'$.
\end{enumerate}
\end{prop}

\begin{proof}
For all assertions one can assume that~$X,Y,X'$, and~$Y'$ are affinoid.
Let us first prove (1). As~$\mathscr H$ is~$X_L$-flat at~$x'$, the coherent sheaf $\mathscr H_{Y'}$
is
naively~$Y_L$-flat at~$y'$, hence is naively $Y$-flat
at $y'$ since~$Y_L$ is naively
flat over~$Y$. In other words, $\mathscr H_{Y',y'}$ is a flat~$\mathscr O_{Y,y}$-module. By assumption,~$\mathsf S_y$ is
exact. Tensoring with the flat~$\mathscr O_{Y,y}$-module $\mathscr H_{Y',y'}$
then yields
the exactness of the sequence
$\mathsf S_{y'},$ whence (1). 

Let us prove (2). As~$X'$ is affinoid, it can be identified
with a
Zariski-closed subspace of~$X_L\times_L D$ where~$D$ is some closed polydisc over the field~$L$.
Right-exactness of the tensor product ensures that 
~$$\mathscr G_{Y_L\times_L D,y'}\to \mathscr F_{Y_L\times_L D,y'}\to \mathscr E_{Y_L\times_L D,y'}\to 0$$
is exact. Since~$X_L\times_L D\to X_L$ is flat
by Lemma~\ref{lem-field-flat}, it follows from assertion (1)
already proven (and applied with $X'=X_L\times_L D$ and $Y'=Y_L\times_L D$) that
the arrow $\mathscr G_{Y_L\times_L D,y'}\to \mathscr F_{Y_L\times_L D,y'}$ is
injective; hence~$$0\to \mathscr G_{Y_L\times_L D,y'}\to \mathscr F\otimes _{Y_L\times_L D,y'}\to \mathscr E_{Y_L\times_L D,y'}\to 0$$ is exact. 
As~$X'$ is a Zariski-closed subspace of~$X_L\times_L D$, the local ring~$\mathscr O_{Y',y'}$ is
naturally isomorphic to
$\mathscr O_{Y_L\times_L D,y'}\otimes_{\mathscr O_{X_L\times_L  D,x'}} \mathscr O_{X',x'}$.
Therefore the sequence
$$0\to (\mathscr G \boxtimes \mathscr H)_{y'}\to  (\mathscr F\boxtimes \mathscr H)_{y'}
\to (\mathscr E\boxtimes \mathscr H)_{y'}\to 0$$ is simply deduced from the exact sequence
$$0\to \mathscr G_{Y_L\times_L D,y'}\to \mathscr F _{Y_L\times_L D,y'}\to \mathscr E_{Y_L\times_L D,y'}\to 0$$ 
by applying the functor~$\otimes_{\mathscr O_{X_L\times_L  D,x'}} \mathscr H_{x'}$. As
the coherent sheaf~$\mathscr E$ is~$X$-flat at~$y$, the~$\mathscr O_{X_L\times_L D,x'}$-module
$\mathscr E_{Y_L\times_L D,y'}$ is flat; it follows then immediately from the~$\mathrm{Tor}_\bullet$
exact sequence
that
$$0\to (\mathscr G \boxtimes \mathscr H)_{y'}\to  (\mathscr F\boxtimes \mathscr H)_{y'}\to (\mathscr E\boxtimes \mathscr H)_{y'}\to 0$$
is exact, whence (2). 

Let us prove (3). We argue by induction on~$n$. For~$n=0$ there is nothing
to prove. 
Let us now assume that~$n\geq 1$ and that the required assertion is true for
all integers~$<n$. Let~$\mathscr N$ be the kernel of~$\mathscr F_1\to \mathscr F_0$. 
The two sequences 
$$\mathscr F_n\to \mathscr F_{n-1}\to \ldots \to \mathscr F_2\to \mathscr N\to 0\;{\rm and}\; \;0
\to \mathscr N\to \mathscr F_1\to \mathscr F_0\to 0$$
are exact at~$y$. 

Let us prove that  $\mathscr N$ is~$X$-flat at~$y$. Let~$Z\to X$
be any morphism with~$Z$ affinoid, let~$t\in T:=Y\times_X Z$ be a point lying above~$y$,
and let~$z$ be the image of~$t$ on~$Z$. 
Since the coherent sheaf~$\mathscr F_0$ is~$X$-flat at~$y$,  assertion
(2) applied with $X'=Z$, $Y'=T$, $x'=z$, $y'=t$, and $\mathscr H=\mathscr O_Z$ ensures that
the sequence
$$0\to \mathscr N_{T,t}\to \mathscr F_{1,T,t}\to \mathscr F_{0,T,t}\to 0$$
is
exact. As~$\mathscr F_1$ and~$\mathscr F_0$ are~$X$-flat at~$y$,
 the~$\mathscr O_{Z,z}$-modules~$\mathscr F_{1,T,t}$
and~$\mathscr F_{0,T,t}$ are flat. By a~$\mathrm{Tor}$ computation,
 it follows that~$\mathscr N_{T,t}$ is also flat over~$\mathscr O_{Z,z}$; therefore,~$\mathscr N$
 is~$X$-flat at~$y$. 
 
 It follows from
(2) 
that
the sequence $$0\to \mathscr N\boxtimes \mathscr H \to \mathscr F_1\boxtimes \mathscr H\to \mathscr F_0 \boxtimes \mathscr H\to 0$$
 is exact at~$y'$. Since we have just seen that~$\mathscr N$ is~$X$-flat at~$y$, the case $n=1$ is settled 
and it follows from the induction hypothesis
 that
the kernel of $\mathscr F_n\to \mathscr F_{n-1}$ is~$X$-flat at~$y$ if~$n\geq 2$
and that~$$\mathscr F_n\boxtimes \mathscr H\
\to \mathscr F_{n-1}\boxtimes \mathscr H \to \ldots \to \mathscr F_2\boxtimes \mathscr H\to \mathscr N\boxtimes \mathscr H\to 0$$
is exact at~$y'$, which yields the exactness of 
~$$\mathscr F_n\boxtimes \mathscr H\to \mathscr F_{n-1}\boxtimes \mathscr H\to \ldots \to \mathscr F_0\boxtimes \mathscr H\to 0$$
at $y'$ and ends the proof of (3). 

Let us prove (4). If $n=1$ then (4) simply means that 
$$\mathrm{Coker}(\mathscr F_1\to \mathscr F_0)\boxtimes \mathscr H
\to \mathrm{Coker}(\mathscr F_1\boxtimes \mathscr H\to \mathscr F_0\boxtimes \mathscr H)$$
is an isomorphism at $y'$, which is true; indeed, it is an isomorphism
at every point of $Y'$ by right-exactness of the tensor product. 
Now assume that $n\geq 2$, and let~$\mathscr N$ be the kernel of~$\mathscr F_{n-1}\to \mathscr F_{n-2}$. 
By assertion (3) the coherent
sheaf~$\mathscr N$ is~$X$-flat at~$y$. Applying (3)
to the complex 
$$0\to \mathscr N\to \mathscr F_{n-1}\to \mathscr F_{n-2}\to \ldots \to \mathscr F_0\to 0$$
(which is exact at $y$),
we see that
$$\mathscr N\boxtimes \mathscr H\to \mathrm{Ker}(\mathscr F_{n-1}\boxtimes \mathscr H \to \mathscr F_{n-2}\boxtimes \mathscr H)$$
is an isomorphism at~$y'$. 
On the other hand,~$\mathrm H_{n-1}(\mathsf S)$ is the cokernel of~$\mathscr F_n\to \mathscr N$.
By right-exactness of~$\bullet \boxtimes \mathscr H$, 
the coherent sheaf~$\mathrm H_{n-1}(\mathsf S)\boxtimes \mathscr H$ is then the cokernel of
$\mathscr F_n\boxtimes \mathscr H \to \mathscr N\boxtimes \mathscr H$.
Hence we see that~$\mathrm H_{n-1}(\mathsf S)\boxtimes \mathscr H
\to \mathrm H_{n-1}(\mathsf S\boxtimes \mathscr H)$ is an isomorphism at~$y'$.
\end{proof}

In scheme theory, flatness is often useful to spread out some properties from 
a given fiber across the ambiant space. As an application of the preceding proposition, 
we give a first example of such a phenomenon in analytic geometry.

\begin{lemm}\label{lem-iso-fiber}
Let~$Y\to X$ be a morphism between~$k$-analytic spaces, let~$y$ be a point of $Y$ and let~$x$ be its image on~$X$. 
Let~$\mathscr G\to \mathscr F$ be a morphism between coherent sheaves on~$Y$. If~$\mathscr G_{Y_x}\to \mathscr F_{Y_x}$
is an isomorphism at~$y$
and~$\mathscr F$ is~$X$-flat at~$y$, then~$\mathscr G\to \mathscr F$ is an isomorphism at~$y$.
\end{lemm}

\begin{proof}
We may
and do assume that both~$Y$ and~$X$ are~$k$-affinoid. Let~$\mathscr N$ be the kernel of~$\mathscr G\to \mathscr F$. 
Since~$\mathscr G_{Y_x}\to \mathscr F_{Y_x}$ is surjective at~$y$, the map
$\mathscr G_{\hr y}\to \mathscr F_{\hr y}$ is surjective, which implies 
that~$\mathscr G\to \mathscr F$ is surjective at~$y$ by \ref{ss-surj-nakayama}; hence the sequence
$$0\to \mathscr N\to \mathscr G \to \mathscr F\to 0$$ is exact
at~$y$. The coherent sheaf~$\mathscr F$ being~$X$-flat at~$y$, \Prop \ref{prop-flat-univexact}
(2) applied to the
morphism~$\mathscr M(\hr x)\to X$ yields the exactness of 
$$0\to \mathscr N_{Y_x}\to \mathscr G_{Y_x} \to \mathscr F_{Y_x}\to 0$$ at~$y$. Since~$ \mathscr G_{Y_x} \to \mathscr F_{Y_x}$
is by assumption an isomorphism at~$y$, this implies that $\mathscr N_{Y_x,y}=0$; therefore $\mathscr N_{\hr y}=0$ and
$\mathscr N_y=0$ by \ref{ss-surj-nakayama}.
As a consequence,~$\mathscr G\to \mathscr F$ is an isomorphism at~$y$.
\end{proof}

We are now going to give two flatness
criteria. The first one describes the behavior of flatness
with respect to extensions. The second one might
look somehow specific, but it will be crucial 
for the study of quasi-smooth morphisms in the next chapter 
and of fiberwise regular sequences in \ref{s-reg-seq}. 

\begin{lemm}\label{lem-flat-extension}
Let $Y\to X$ be a morphism of $k$-analytic spaces, let
$y$ be a point of $Y$, and let $0\to \mathscr E\to \mathscr E'\to \mathscr E''\to 0$ be a
sequence of coherent sheaves on $Y$. Assume that this sequence is exact at $y$.

\begin{enumerate}[1]
\item If $\mathscr E$ and $\mathscr E''$ are $X$-flat at $y$, so is $\mathscr E'$. 
\item If $\mathscr E'$ and $\mathscr E''$ are $X$-flat at $y$, so is $\mathscr E$.
\end{enumerate}
\end{lemm}

\begin{proof}
Let us prove (1) (\resp 2). So we assume that
$\mathscr E''$ and $\mathscr E$ (\resp $\mathscr E'$) 
is $X$-flat.
We immediately reduce to the case where 
$X$ and $Y$ are good. Let $Z\to X$ be any morphism
of analytic spaces with $Z$ good, let $t$ be a pre-image 
of $y$ on $T:=Y\times_X Z$, and let $z$ be the image of $t$ in $Z$. We
have to show that 
$\mathscr E'_{T,t}$
(\resp $\mathscr E_{T,t}$)
is flat over $\mathscr O_{Z,z}$. 

By Proposition 
\ref{prop-flat-univexact}
(2), the sequence
\[0\to \mathscr E'_{T,t}\to \mathscr E_{T,t}\to \mathscr E''_{T,t}\to 0\]
is still exact.
It follows from our flatness assumptions that $\mathscr E''_{T,t}$ and $\mathscr E_{T,t}$
(\resp $\mathscr E'_{T,t}$) are $\mathscr O_{Z,z}$-flat, hence so is $\mathscr E'_{T,t}$
(\resp $\mathscr E_{T,t}$)
by a straightforward $\mathrm{Tor}$
computation. 
\end{proof}

\begin{lemm}\label{lem-flattrick-reg}
Let~$Y\to X$ be a morphism of~$k$-analytic spaces, let~$y$
be a point of  $Y$ and let~$x$ be its image in~$X$. Let~$\mathscr G\to \mathscr F$
be a morphism of coherent sheaves on~$Y$, and
let~$\mathscr E$ be its cokernel.
Assume that~$\mathscr G$ and~$\mathscr F$ are~$X$-flat at~$y$, and that $\mathscr G_{Y_x}\to \mathscr F_{Y_x}$
is injective at~$y$.
Under this assumption,~$\mathscr E$
is~$X$-flat at~$y$.
\end{lemm}

\begin{proof}
We can assume that~$Y$ and~$X$ are~$k$-affinoid. Let~$L$ be an analytic extension of~$k$, 
let~$X'$ be an~$L$-affinoid space, and let~$X'\to X$ be a morphism. 
Let~$y'$ be a pre-image of~$y$ on~$Y':=Y\times_X X'$ and let~$x'$ be its image on~$X'$. We have to prove
that $\mathscr E_{Y'}$ is naively
flat over~$X'$
at~$y'$. By \Prop~\ref{prop-flat-xl}, 
it suffices to do this after
some ground field extension. 
Hence we may assume that~$x'$ is~$L$-rational.

Since~$\mathscr F$ and~$\mathscr G$ are~$X$-flat by assumption,
the sequence
$$\mathscr G_{Y',y'}\to \mathscr F_{Y',y'}\to \mathscr E_{Y',y'}\to 0$$ is the truncation of a flat resolution 
of the~$\mathscr O_{X',x'}$-module~$\mathscr E_{Y',y'}$. 
The map~$\mathscr G_{Y_x}\to \mathscr F_{Y_x}$ is injective at~$y$; by 
naive flatness of ground field extension, 
$\mathscr G_{Y'_{x'}, y'}\to \mathscr F_{Y'_{x'}, y'}$ is injective too. But since~$x'$
is~$L$-rational, the local ring~$\mathscr O_{Y'_{x'}, y'}$ is equal to
the quotient $\mathscr O_{Y',y'}/\mathfrak m_{x'}\mathscr O_{Y',y'}$. 
Hence by passing to the quotient of the truncated flat resolution
$$\mathscr G_{Y',y'}\to \mathscr F_{Y',y'}\to \mathscr E_{Y',y'}\to 0$$
modulo~$\mathfrak m_{x'}$, one gets an exact sequence whose first arrow is injective. This immediately implies that
~$\mathrm{Tor}_1^{\mathscr O_{X',x'}}(\mathscr E_{Y',y'}, \mathscr O_{X',x'}/\mathfrak m_{x'})=0$. 
As a consequence,~$\mathscr E_{Y',y'}$
is a flat~$\mathscr O_{X',x'}$-module (\cite{sga1}, Expos\'e IV, \Th 5.6). 
\end{proof}

We end this section by showing that flatness
ensures a reasonable behavior of local
dimension (this is a generalization of \Cor \ref{cor-finflat-dim}).
We use the notion of dimension for modules
and coherent sheaves (\ref{ss-conv-alg},\ref{ss-interpret-support}).

\begin{lemm}\label{lem-flatness-dim}
Let~$Y\to X$ be a morphism between~$k$-analytic spaces, let~$\mathscr F$ be a coherent sheaf on~$Y$,
let~$y$ be a point of $Y$ and let $x$ be its image in $X$. Assume that $y$ belongs to $\supp F$
(this is the case if and only it belongs to $\mathrm{Supp}(\mathscr F_{Y_x})$,
because by
\ref{ss-pointwise-rk} both properties are equivalent to the fact that $\mathscr F_{\hr y}\neq 0$)
and that $\mathscr F$ is~$X$-flat at~$y$. 
One has
then the equality~$\dim_y \mathscr F=\dim_y \mathscr F_{Y_x}+\dim_x X.$
\end{lemm}

\begin{proof}
We immediately reduce to the case where~$Y$ and~$X$ are affinoid.
By extending the ground field, we may assume that~$y$ is a~$k$-point (hence, so is~$x$). 
By assumption,~$\mathscr F_y$ is a flat~$\mathscr O_{X,x}$-module;
moreover, the quotient~$\mathscr O_{Y,y}/\mathfrak m_x\mathscr O_{Y,y}$ is equal
to $\mathscr O_{Y_x,y}$ because~$x$ is a rigid point.
By \Cor 6.1.2 of \cite{ega42} one has therefore the equality
$$\dim_{\mathrm{Krull}} \mathscr F_{Y,y}=\dim_{\mathrm {Krull}}\mathscr F_{Y_x,y}+\dim_{\mathrm {Krull}}\mathscr O_{X,x}.$$
Since $y$ and $x$ are rigid, Corollary \ref{cor-interp-centdim} and GAGA
for the support of a coherent sheaf allow
us to rewrite 
it as $\dim_y \supp F= \mathscr F=\dim_y \mathscr F_{Y_x}+\dim_x X.$
\end{proof}

\chapter{Quasi-smooth morphisms}\label{QSM} 

One of the most important
example of flat morphism in algebraic geometry, 
both conceptually and technically, is that of a smooth morphism:
this is a locally finitely presented flat morphism with geometrically
regular fibers. The purpose of this chapter is to introduce and study
the corresponding class of maps in analytic geometry, which are said to be
\emph{quasi-smooth}. 

In fact, the definition of a scheme-theoretic
smooth morphism we have just given is not the usual one:
one classically defines smooth morphisms using
the sheaf of relative K\"ahler differentials, and then proves
that they can be characterized by the aforementioned property. 
This is what we shall do here
(we have more precisely been inspired by the approach
on smoothness
of Bosch, L\"utkebohmert and Raynaud
in \cite{bosch-l-r1990}). Thus after having recalled the definition
and the basic properties of
the coherent sheaf $\Omega_{Y/X}$
attached to a morphism of $k$-analytic spaces $Y\to X$
(Section \ref{s-remind-omega}), we use kind of 
a Jacobian criterion to define what it means for $Y\to X$
to be quasi-smooth at a given point $y$ of $Y$
(Definition \ref{def-qsm}); we say that $Y\to X$
is \emph{quasi-\'etale} at $y$ if it is quasi-smooth
and quasi-finite at $y$ (a former definition of quasi-\'etaleness
had been given by Berkovich in \cite{berkovich1994}; it is consistent with ours
by
\ref{lem-qetad-qetvb}). 

Then in Section \ref{s-qsm-flat}
we prove the expected characterization 
of quasi-smooth morphism (Theorem \ref{thm-main-qsm}):
if $y$ is a point of $Y$ lying above a point $x$ of $X$, then $Y\to X$
is quasi-smooth at $y$ if and only if $Y$ is $X$-flat at $y$
and $Y_x$ is geometrically regular at $y$. 

Section \ref{s-qsm-sm}
explains the links between
quasi-smoothness 
and quasi-\'etaleness on one hand, and smoothness
and \'etaleness in the sense of Berkovich (\cite{berkovich1993} 3.3 and 3.5)
on the other hand. 
We prove more precisely the
following: $Y\to X$ is 
\'etale at $y$ if and only it it is quasi-\'etale and boundaryless at $y$; 
\emph{if morevoer $Y$ and $X$ are good}, then $Y\to X$ is
smooth at $y$ if and only if it is quasi-smooth and boundaryless at $y$
(see
Corollary \ref{cor-qsm-sminner} for smoothness, and Remark \ref{rem-qsm-goodness}
for some comments about the goodness assumption and the \'etale case). 

The chapter ends with Section \ref{s-transfer-qsm}, in which we prove
that all usual algebraic
properties
are preserved by quasi-smooth maps (Proposition \ref{prop-tranferqsm-concrete}).
This rests on the following fact: if $Y$
and $X$ are good and $y$ is a point of $Y$ lying over a point $x$ of $X$ and at which
$Y\to X$ is quasi-smooth, then $\spec \mathscr O_{Y,y}\to \spec \mathscr O_{X,x}$
is regular; \ie, flat with geometrically regular fibers (Theorem \ref{thm-qsm-schemereg}). 
 
\section{Reminders about the sheaf of relative differentials}\label{s-remind-omega}
We begin with some reminders about the 
sheaf of K\"ahler differentials in analytic geometry; a 
general reference for the results of this section is \cite{berkovich1993}, \S 3.3. 

\subsection{}\label{ss-remind-omega}
Let~$Y\to X$ be a morphism between~$k$-analytic spaces.
The diagonal map~$\delta : Y\to Y\times_XY$ is G-locally a closed immersion and it has therefore
a {\em conormal sheaf} (see the paragraph before Remark 1.3.8 in \cite{berkovich1993});
this is a coherent sheaf on~$Y$ which is denoted\footnote{Berkovich denotes it by
$\Omega_{Y_{\mathrm G}/X_{\mathrm G}}$; for the sake of simplicity, and according to our general
conventions, we have decided to simply denote it by~$\Omega_{Y/X}$.} by~$\Omega_{Y/X}$
and is called the {\em sheaf of relative (K\"ahler) differentials of~$Y$ over~$X$}. 

\subsection{Contravariant functoriality}\label{ss-contravariant-omega} Let

$$\xymatrix{
{Y'}\ar[r]\ar[d]&{X'}\ar[d]\\
Y\ar[r]&X
}$$
be a commutative diagram of 	analytic spaces. It gives rise to a natural morphism
$(\Omega_{Y/X})_{Y'}\to \Omega_{Y'/X'}$, which is an isomorphism if (at least) one of the two following conditions is fullfilled: 

\begin{itemize}
\item The above diagram is cartesian (see \cite{berkovich1993}, \Prop 3.3.2). 

\item Both maps $Y'\to X'$ and $Y\to X$ are inclusion of analytic domains. 
\end{itemize}

In particular if $V$ is an analytic domain of $Y$, then $\Omega_{V/X}=(\Omega_{Y/X})_V$; if $U$ is an analytic domain of $X$, then $\Omega_{U/X}=0$. 
If $x$ is a point of $X$, then $(\Omega_{Y/X})_{Y_x}=\Omega_{Y_x/\hr x}$ and $(\Omega_{Y/X})_{\hr y}=(\Omega_{Y_x/\hr x})_{\hr y}$ for every 
$y\in Y_x$.

\subsection{The universal differential}
Let $V$ be an analytic domain of $Y$ and let $p_1$ and $p_2$ be the two projections from $V\times_XV$ to $X$. If $f$ is an
analytic function on $V$, then $p_1^\ast  f-p_2^\ast  f$ vanishes on the diagonal $V\hookrightarrow V\times_X V$, hence defines
an element $\d f$ of $\Omega_{V/X}(V)=\Omega_{Y/X}(V)$. The map $\d \colon \mathscr O_Y\to \Omega_{Y/X}$ is an $X$-derivation, 
and $(\Omega_{Y/X},\d)$ is the initial object of the category of coherent sheaves $\mathscr F$ on $Y$ equipped with an $X$-derivation 
$\mathscr O_Y\to \mathscr F$. In the situation of \ref{ss-contravariant-omega} above, the
natural morphism
$(\Omega_{Y/X})_{Y'}\to \Omega_{Y'/X'}$ commutes with the derivations.

%\subsection{}
%We are now going to describe the basic rules that one can use in order to make computations on Kähler differentials. 

\subsection{GAGA principle}
\label{ss-omega-gaga}
If~$A$ is a~$k$-affinoid algebra and if~$\mathscr Y\to \mathscr X$ is a morphism between~$A$-schemes of finite type, there is a natural
isomorphism
$$(\Omega_{\mathscr Y/\mathscr  X})\an \simeq  \Omega_{\mathscr  Y\an/\mathscr  X\an},$$
which commutes with the derivations. 

\subsection{}\label{ss-omega-exactseq}
Let~$$\xymatrix{
Z\ar[r]& Y\ar[r]&X}$$ be a diagram in the category of~$k$-analytic spaces. The natural sequence 
$$(\Omega_{Y/ X})_Z\to \Omega_{Z/X}\to \Omega_{Z/Y}\to 0$$ is then exact (\cite{berkovich1993}, \Prop 3.3.2 (i) ).

\subsection{}\label{ss-omega-affspace}
Let~$X$ be a~$k$-analytic space. For every non-negative integer $n$, the coherent
 $\mathscr O_{\A^n_X}$-module
~$\Omega_{\A^n_X/X}$ is free 
with basis~$\d T_1,\ldots,\d T_n$.  

\subsection{}\label{ss-omega-zarclos}
Let~$Y\to X$ be a morphism of~$k$-analytic spaces, let~$(f_i)$ be a
family of analytic functions on~$Y$, and let~$Z$ be the closed
analytic subspace
of~$Y$ defined by the sheaf of ideals~$(f_i)_i$.
The morphism $(\Omega_{Y/X})_Z\to \Omega_{Z/X}$ is a surjection
whose kernel is generated by the pullbacks of the $\d f_i$'s (\cite{berkovich1993}, \Prop 3.3.2 (ii) ). 

%\subsection{} Let~$X$ be a~$k$-analytic space and let~$x$
%be a point of $X$. 

\subsection{} Let~$X$ be a
good $k$-analytic space, let~$x$
be a point of $X$ such that $\hr x=k$, and let
$\mathfrak m$ be the maximal ideal of
$\mathscr O_{X,x}$. The map $f\mapsto f-f(x)$ from $\mathscr O_{X,x}$ to $\mathfrak m$ is then a $k$-derivation which
induces an 
isomorphism
$$(\Omega_{X/k})_{\kappa(x)}=(\Omega_{X/k})_x/\mathfrak m(\Omega_{X/k})_x\simeq \mathfrak m/\mathfrak m^2, $$
whose
inverse
isomorphism is induced by the derivation $\d$ (this can be checked by direct computation, as in algebraic geometry). 

\subsection{}\label{ss-omega-inequality}
Let $X$ be an analytic space and let $x$ be a point of $X$. 
If follows respectively
from Lemma 6.2 and \Prop 6.3 of \cite{ducros2009}
that $\rk x (\Omega_{X/k})\geq \dim_x X$ and
that the following are equivalent: 

\begin{enumerate}[i]

\item One has the equality~$\rk x (\Omega_{X/k})=\dim_x X$. 

\item The space $X$ is geometrically regular at~$x$. 
\end{enumerate}

Moreover if~$\hr x=k$ or if~$k$ is perfect,
then (i) and (ii)
hold if and only if~$X$ is regular at~$x$.
Indeed, this also follows from Lemma 6.2 and \Prop 6.3 of \cite{ducros2009}
except possibly when $\hr x=k$ and $\abs{k\gpm}=1$. But in this situation we can first perform
a scalar extension to $k_r$ for some arbitrary positive $r$; since $k_r$ is analytically separable over $r$, this operation
has no effect on regularity, hence allows us
to reduce
to the non-trivially valued case.

\subsection{}
Assume that (i) and (ii) above hold, and let~$d$ be the dimension of~$X$ at~$x$. There exists a purely~$d$-dimensional open neighborhood~$U$ of~$x$ in $X$
such that the coherent sheaf~$\Omega_{U/k}$ is free of rank~$d$ (\cite{ducros2009}, \Prop 6.6; note that this implies that (i) and (ii) hold at {\em every} point of~$U$). 

\begin{rema}\label{rem-georeg-fiber}
Let $Y\to X$ be a morphism of $k$-analytic spaces, let $y$
be a point of $Y$ and let $x$
be its image in $X$. It follows from \ref{ss-omega-inequality}
that
$$\rk y (\Omega_{Y/X})\leq \dim_y Y_x$$
with equality if and only if the fiber $Y_x$ is geometrically regular at $y$. 

\end{rema}

\section{Quasi-smoothness: definition and first properties}

We begin with a technical lemma,
which will be needed for dealing with the Jacobian criterion we have in mind. 

\begin{lemm}
\label{lem-pre-qsm}
Let~$X$ be a~$k$-analytic space, let~$n$ be a
non-negative integer, and let~$V$ be an affinoid domain of~$\A^n_X$. 
Let~$Y$ be a Zariski-closed subset of~$V$, let~$I$ be the corresponding ideal of the ring of analytic functions
on~$V$, and let~$g_1,\ldots, g_r$ be a generating family of~$I$. For every~$z\in Y$,
let~$s(z)$ denote the rank of the family~$((\d g_1)(z), \ldots, (\d g_r)(z))$ in the vector space~$(\Omega_{V/X})_{\hr z}$.

\begin{enumerate}[1]
\item For every~$z\in Y$ the rank $\rk z(\Omega_{Y/X})$ is equal to~$n-s(z)$; in particular,~$$\rk z (\Omega_{Y/X})\geq n-r.$$

\item Assume that there exists~$y\in Y$ with~$s(y)=r$ (which is equivalent, by {\rm (1)}, to the fact that~$\rk y (\Omega_{Y/X})=n-r$). Then:

\begin{itemize}[label=$\bullet$]
\item Every generating family of~$I$ has cardinality at least~$r$.  

\item There exists an affinoid neighborhood~$U$ of~$y$ in~$V$ such that
the morphism $U\cap Y\to X$ is purely of relative dimension~$n-r$, and such that~$s(z)=r$ for every~$z$
belonging to $U\cap Y$; in particular, every fiber of~$U\cap Y\to X$ is geometrically regular.
\end{itemize}
\end{enumerate}

\end{lemm}

\begin{proof}
Let us first prove (1). By \ref{ss-omega-affspace}, the~$\hr z$-vector space~$(\Omega_{V/X})_z$ is~$n$-dimensional;
and it follows from \ref{ss-omega-zarclos} that ~$(\Omega_{Y/X})_{\hr z}$ is naturally isomorphic to~$$(\Omega_{V/X})_{\hr z}/((\d g_1)(z),\ldots, (\d g_r)(z)),$$ whence (1).

Now let us come to assertion (2). If~$(h_1,\ldots, h_t)$ is a a generating family of~$I$, applying (1) to it yields the inequality~$n-r\geq n-t$; \ie, $t\geq r$, as required. 
Let~$x$ be a point of $X$. Being an affinoid domain of~$\A^{n, \rm an}_{\hr x}$, the fiber~$V_x$ is purely~$n$-dimensional. As the ideal of~$Y_x$ in~$V_x$ is generated by~$r$ functions, it follows from the {\em Hauptidealsatz} applied on the noetherian scheme $(V_x)\al$
that the Krull codimension of any irreducible component of~$Y_x$ in~$V_x$ is at most~$r$.
Therefore, the dimension of such a component is at least~$n-r$, and it follows that~$\dim_z (Y\to X)\geq n-r$ for every~$z\in Y$.

By upper-semi-continuity of the rank of the stalks of a given coherent sheaf
(\ref{ss-pointwise-rk}), 
there exists an affinoid neighborhood~$U$ of~$y$ in~$V$ such that~$\rk z (\Omega_{Y/X})$ is
bounded by~$n-r$ for every~$z\in U\cap Y$; note that this rank is then
actually {\em equal} to~$n-r$ in view of (1). Let~$z$ be a point of $U\cap Y$ and let~$x$ be its image on~$X$.
The~$\hr x$-analytic space~$Y_x$ if of dimension {\em at least}~$n-r$ at~$z$; and~$\rk z \Omega_{Y_x/\hr x}= \rk z (\Omega_{Y/X})=n-r$. We thus deduce from \ref{ss-omega-inequality} that ~$Y_x$ is
of dimension~$n-r$ at~$z$, which ends the proof (the claim about geometric regularity comes from \ref{ss-omega-inequality}).
\end{proof}

\begin{defi}\label{def-jacob-pres}
Let~$Y\to X$ be a morphism of~$k$-analytic spaces and let~$y$ be a point of $Y$. Let~$W$ be an affinoid
domain of~$Y$ containing~$y$, let~$n$ be an element of $\N$, and let~$V$ be an affinoid domain of~$\A^n_X$ such that~$W\to X$
goes through a closed immersion~$W\hookrightarrow V$; let us denote by~$I$ the ideal
defining the latter (in the ring of analytic functions on~$V$), and set~$r=n-\rk y(\Omega_{Y/X})$.
We say that the diagram~$W\hookrightarrow V\subset \A^n_X$ is a {\em Jacobian 
presentation of~$Y\to X$ at~$y$} if~$I$ can be generated by~$r$ elements. \index{Jacobian presentation}
\end{defi}

\subsection{}\label{ss-jacob-pres}
We use the notation of Definition \ref{def-jacob-pres} above. 
If $W\hookrightarrow V\subset \A^n_X$ is a  Jacobian  presentation of~$Y\to X$ at~$y$, it follows from
Lemma~\ref{lem-pre-qsm}
that~$r$ is the minimal cardinality of a generating family of~$I$, that~$Y\to X$ is of dimension~$n-r$ at~$y$, and that the fiber of~$Y\to X$ containing~$y$ is geometrically regular at~$y$. 
Lemma \ref{lem-pre-qsm}
also ensures that there exists an affinoid neighborhood~$V'$ of~$y$ inside~$V$ such that~$W\times_VV'\hookrightarrow V'$ is a Jacobian presentation of~$Y\to X$ at {\em each of its points}, and such that $W\times_V V'$ is purely of relative dimension~$n-r$ over~$X$.

\begin{defi}\label{def-qsm}
Let~$Y\to X$ be a morphism
of~$k$-analytic spaces, and let~$y$
be a point of $Y$. We say that~$Y\to X$ is {\em quasi-smooth} at~$y$ if there exists  a Jacobian presentation of~$Y\to X$ at~$y$.\index{morphism!of analytic spaces!quasi-smooth}\index{quasi-smooth morphism}
We say that it is {\em quasi-smooth}
if it is quasi-smooth at every point of $Y$.
\end{defi}

\subsection{}\label{ss-qsm} Let $Y\to X$ be a morphism of $k$-analytic spaces, let $y$ be a point of $Y$, and let $x$ be its image in $X$. 
If $Y\to X$ is quasi-smooth at $y$, it follows from \ref{ss-jacob-pres} that $Y_x$ is geometrically regular at $y$. 

\begin{defi}\label{def-qet}
A morphism  of~$k$-analytic spaces~$Y\to X$ is {\em quasi-\'etale} at a point~$y$ of~$Y$ if it is\index{morphism!of analytic spaces!quasi-etale@quasi-\'etale}
quasi-smooth and quasi-finite\index{quasi-étale morphism}
at~$y$; and it is {\em quasi-\'etale} if it is quasi-\'etale at every point of~$Y$. 
\end{defi}

\begin{rema}
An analytic space~$(k,X)$ is called quasi-smooth, \resp quasi-\'etale, at a given point~$x$ of $X$ if~$X\to \mathscr M(k)$ is.  
\end{rema}

\subsection{Comments on the terminology}
Berkovich has defined (\cite{berkovich1993}, \S3) the notions of  \'etale and smooth maps.\index{morphism!of analytic spaces!smooth}\index{morphism!of analytic spaces!etale@\'etale}\index{smooth morphism}\index{e@étale morphism}
We shall see below (Theorem
\ref{thm-qsm-embeddsm}
and Remark \ref{rem-qsm-goodness}) that a map is \'etale at
a given point if and only if it is quasi-\'etale and
boundaryless
at that point; and that a map between {\em good}~$k$-analytic spaces is smooth at a given point if and only if it is quasi-smooth and
boundaryless
at that point. (For some comments about the need for a goodness assumption, see Remark \ref{rem-qsm-goodness}).

There is already a notion of quasi-\'etale morphism, which was defined by Berkovich (\cite{berkovich1994}, \S3); we
shall see below that his definition is equivalent to ours (Lemma \ref{lem-qetad-qetvb}).  

In \cite{ducros2009}, \S 6, an analytic space was said to be quasi-smooth ({\em quasi-lisse} in French) at~$x$ if it is
geometrically regular at~$x$; this turns out to be consistent with our current definition of quasi-smoothness
(Corollary \ref{cor-qsmfobs-qsmexc}.)

If~$\abs {k\gpm} \neq \{1\}$, if~$Y$ and~$X$ are strictly~$k$-analytic spaces and if~$y$ is a rigid point of~$Y$, quasi-smoothness of~$Y\to X$ at~$y$ is nothing but {\em rig-smoothness} of~$Y\to X$ at~$y$; we nevertheless have chosen to use ``quasi-smooth" instead of ``rig-smooth" to be consistent with the terminology ``quasi-\'etale".

%\subsection{} We are going to list some basic properties of quasi-smooth morphisms, and to give some examples. But before doing it, let us make a technical remark we shall
%use several times in the sequel.
%
%\begin{rema}\label{rem-gerr-grau}
%It follows immediately from Gerritzen-Grauert theorem (more precisely from Temkin's version of it for Berkovich spaces, \cite{temkin2005}, \Th 3.1) that if~$X$ is an analytic space
%(not necessarily good) and if~$Y$ is a Zariski-closed subspace of~$X$, and if $V$
%is any analytic domain of $Y$, then~$V$ can be G-covered by affinoid domains of the
%form~$U\cap Y$ with~$U$ an affinoid domain of~$X$.
%\end{rema}

\subsection{} \label{ss-aff-qsm}
Let $X$ be an analytic space and let $n$ be a non-negative integer. The space $\A^n_X$ is then quasi-smooth over $X$ of relative dimension~$n$: indeed, for every~$y\in \A^n_X$ and every affinoid domain~$W$ of~$\A^n_X$ containing~the point $y$, the diagram~$W\simeq W\subset \A^n_X$ is a Jacobian presentation of~$Y\to X$ at~$y$. 

\subsection{} \label{ss-andom-qsm}
Let~$Y\to X$ be a morphism of $k$-analytic spaces, let $y$ be a point of $Y$, let~$V$ be an analytic domain of~$Y$ containing~$y$, and let~$U$ be an analytic domain of~$X$ containing the image of~$V$. Then~$Y\to X$ is quasi-smooth at~$y$ if and only if $V\to U$ is quasi-smooth at $y$.

Indeed, let us first assume that~$V\to U$ is quasi-smooth at~$y$. Then if~$Z\hookrightarrow T\subset \A^n_U$ is a Jacobian presentation of~$V\to U$ at~$y$, it follows immediately from the definition that~$Z\hookrightarrow T\subset \A^n_X$ is a Jacobian presentation of~$Y\to X$ at~$y$.

Conversely, let us assume that~$Y\to X$ is quasi-smooth at~$y$, and let us choose a Jacobian presentation~$Z\hookrightarrow T\subset \A^n_X$ of~$Y\to X$ at~$y$; the image of~$y$ in~$\A^n_X$ belongs to~$\A^n_U$. Let~$T'$ be an affinoid domain of~$\A^n_U$ that
contains the image of~$y$. The fiber product~$Z':=Z\times_TT'$ is an affinoid domain of~$Y$ which contains~$y$, and~$Z'\to T'$ is a closed immersion. 
By Remark \ref{rem-gerr-grau}
there exists an affinoid domain~$T''$ of~$T'$ containing the image of~$y$ such that~$Z'':=Z'\times_{T'}T''$ is included in~$V\cap Z'$. It follows from the construction that~$Z''\hookrightarrow T''\subset \A^n_U$ is a Jacobian presentation of~$V\to U$ at~$y$.

\subsection{}\label{ss-id-qsm}
If $X$ is an analytic space, the morphism~${\rm Id}_X$ is quasi-mooth (\ref{ss-aff-qsm} with~$n=0$); it follows by \ref{ss-andom-qsm} that if~$Y$ is an analytic domain of~$X$ then~$Y\to X$ is quasi-smooth.

\subsection {Behavior with respect to base change}
\label{ss-qsm-bc} Let ~$X'$ be an analytic space and let~$X'\to X$ be a morphism. If~$y$ is a point of $Y$, $y'$ is a point of~$Y':=Y\times_X X'$ lying over~$y$, and $Y\to X$ is quasi-smooth at~$y$ then~$Y'\to X'$ is quasi-smooth at~$y'$. Indeed, let~$W\hookrightarrow V\subset \A^n_X$ be a Jacobian presentation of~$Y\to X$ at~$y$, and let~$V'$ (resp.~$W'$) denote the fiber product~$V\times_X X'$ (resp.~$W\times_X X'$). Let~$V''$ be any affinoid domain of~$V'$ which contains the image of~$y'$ by the closed immersion~$V'\hookrightarrow W'$; the fiber product~$W'':=W'\times_{V'}V''$ is an affinoid domain of~$Y'$, and it is easily seen that~$W''\hookrightarrow V''\subset \A^n_{X'}$ is a Jacobian presentation of~$Y'\to X'$ at~$y'$.

\subsection{Behavior with respect to composition}
\label{ss-qsm-composition}
Let~$Z$ be an analytic space, let~$Z\to Y$ be a morphism, let~$z$ be a point of $Z$, and let~$y$ be its image in~$Y$. If~$Z\to Y$ is quasi-smooth at~$z$ and if~$Y\to X$ is quasi-smooth at~$y$ then the composite map~$Z\to Y\to X$ is quasi-smooth at~$z$. 

Indeed, let~$W\hookrightarrow V\subset \A^n_X$ be a Jacobian presentation of~$Y\to X$ at~$y$. As~$Z\to Y$ is quasi-smooth at~$z$, the map~$Z\times_YW\to W$ is quasi-smooth at~$z$ too by \ref{ss-andom-qsm} or \ref{ss-qsm-bc}. Let~$T\hookrightarrow S\subset \A^m_W$ be a Jacobian presentation of~$Z\times_YW\to W$ at~$z$. As~$S$ is an affinoid domain of
the Zariski-closed subspace $\A^m_W$ of
$\A^m_V$, it follows from Remark \ref{rem-gerr-grau}
that there exists an affinoid domain~$S'$ of~$\A^m_V$ such that~$S'\cap \A^m_W$ is contained in~$S$ and contains the image of~$z$; set~$T'=T\times_{S}(S'\cap \A^m_W)$. The morphism~$T'\hookrightarrow S'$ is equal to the composition of~$T'\hookrightarrow S'\cap \A^m_W$ and~$S'\cap \A^m_W\hookrightarrow S'$, hence is a closed immersion. Being  an affinoid domain of~$\A^{m}_V$, which is itself an analytic domain of~$\A^{n+m}_X$, the space~$S'$ is an affinoid domain of~$\A^{n+m}_X$. 
Set~$d=\rk y (\Omega_{Y/X})$ and $\delta=\rk z (\Omega_{Z/Y})$. 
It follows from \ref{ss-omega-zarclos} that $\rk z(\Omega_{Z/X})\leq d+\delta$.

Now, as~$W\hookrightarrow V\subset \A^n_X$ is a	Jacobian presentation of~$Y\to X$ at~$y$, the Zariski-closed subspace~$W$ of~$V$ can be defined by~$n-d$ equations; hence the Zariski-closed subspace~$S'\cap \A^m_W$ of~$S'$ can also be defined by~$n-d$ equations. 
And as~$T\hookrightarrow S\subset \A^m_W$ is a Jacobian presentation of~$Z\times_YW\to W$ at~$z$, the closed analytic subspace~$T$ of~$S$ can be defined by~$m-\delta$ equations; hence the closed analytic subspace~$T'$ of~$S'\cap \A^m_W$ can be defined by~$m-\delta$ equations. It follows that the closed
analytic subspace~$T'$ of~$S'$ can be defined using~$m+n-d-\delta$ equations.
We deduce
then from Lemma~\ref{lem-pre-qsm}~(1) that $\rk z (\Omega_{Z/X})\geq d+\delta.$ On the other hand, we have proven above that~$\rk z(\Omega_{Z/X})\leq d+\delta$, whence the equality
$$\rk z (\Omega_{Z/X})= d+\delta.$$
Therefore~$T'\hookrightarrow S'\subset \A^{n+m}_X$ is a Jacobian presentation of~$Z\to X$ at~$z$,
so~$Z\to X$ is quasi-smooth at~$z$.

\subsection{}
\label{ss-qsm-gaga}
Let~$A$ be a~$k$-affinoid algebra and let~$\mathscr Y\to \mathscr X$ be a morphism
between~$A$-schemes of finite type. Let~$y$
be a point of $\mathscr Y\an$. If~$\mathscr Y\to \mathscr X$ is smooth at~$y\al $ then
$\mathscr Y\an \to \mathscr X\an$ is quasi-smooth at~$y$. Indeed, as~$\mathscr Y\to \mathscr X$ is
smooth at~$y\al$, there exists an integer~$n$, an affine open neighborhood~$\mathscr V$ of~$y\al$,
and an affine  open subset~$\mathscr U$ of~$\A^n_{\mathscr X}$ so that~$\mathscr V\to \mathscr X$ goes
through a closed immersion~$\mathscr V\hookrightarrow \mathscr U$ whose ideal can be generated by~$r$
elements, where~$r=n-\rk {y\al }( \Omega_{\mathscr Y/\mathscr X})$. Now if~$U$ is any affinoid domain of~$\mathscr U\an$
containing the image of~$y$ and if we set~$V=\mathscr V\an\times_{\mathscr U\an}U$ then~$V\hookrightarrow U\subset \A^n_{\mathscr X\an}$
is a Jacobian presentation of~$\mathscr Y\an\to \mathscr X\an$
at~$y$ (due to \ref{ss-omega-gaga}). 

\subsection{}\label{ss-qet-obvious}
By obvious relative dimension arguments, the conclusions
in \ref{ss-andom-qsm}-\ref{ss-qsm-gaga} remain true with ``quasi-smooth" replaced by ``quasi-\'etale".

\subsection{}
\label{ss-et-qet}
Let~$Y\to X$ be a map between~$k$-analytic spaces, let~$y$
be a point of $Y$ and let~$x$ be its image on~$X$. Assume that~$Y\to X$ is \'etale at~$y$. Under this assumption, there exists an affinoid domain~$U$ of~$X$ containing~$x$ and an affinoid domain~$V$ of~$Y$ containing~$y$ such that~$V\to X$ goes through a finite \'etale map~$V\to U$. But saying that~$V\to U$ is finite \'etale simply means that~$V\al\to U\al$ is finite \'etale, which implies that~$V\to U$ is quasi-\'etale (\ref{ss-qsm-gaga} and \ref{ss-qet-obvious}); as a consequence,~$Y\to X$ is quasi-\'etale at~$y$ by \ref{ss-andom-qsm}.

\subsection{}
\label{ss-sm-qsm}
Let~$Y\to X$ be a map between~$k$-analytic spaces and let~$y$
be a point of $Y$. Assume that~$Y\to X$ is smooth at~$y$. By definition, there exists an open neighborhood~$V$ of~$y$ in~$Y$ such that~$V\to X$ goes through an \'etale map~$V\to \A^n_X$ for some~$n$. It follows from \ref{ss-et-qet} above that~$V\to \A^n_X$ is quasi-\'etale.
Since~$\A^n_X\to X$ is quasi-smooth by \ref{ss-aff-qsm}, one deduces from \ref{ss-qsm-composition} and \ref{ss-andom-qsm} that~$Y\to X$ is quasi-smooth at~$y$.  

\section{Quasi-smoothness, flatness and fiberwise geometric regularity}\label{s-qsm-flat}
\markright{\thesection.~Quasi-smoothness and flatness}

Our goal is to establish some expected properties of quasi-smoothness: the fact that a morphism $\phi \colon Y\to X$
is quasi-smooth at a point $y$ of $Y$ if and only if $Y\to X$ is flat at $y$
and $Y_{\phi(y)}$ is geometrically regular at $y$; and the fact that
if this is the case, then $\Omega_{Y/X}$ is free of rank $\dim_y \phi$ at $y$. We begin with a slightly technical lemma
that
has no interest by itself -- especially in view of the announced results -- but which we shall need to
argue by induction on the number of equations in a given Jacobian presentation. 

\begin{lemm}\label{lem-qsm-technical}
Let $d$ be a non-negative integer, let~$Y\to X$ be a morphism of~$k$-analytic spaces, and
let~$y$ be a point of $Y$ at which
$Y\to X$ is quasi-smooth of relative dimension $d$. Assume that
both $\mathscr O_Y$ and $\Omega_{Y/X}$ are $X$-flat at $y$.

\begin{enumerate}[1]
\item
The coherent sheaf $\Omega_{Y/X}$ is free of rank~$d$ at~$y$.

\item Let $f$ be an analytic function on $Y$ such that the element $(\d f)(y)$ of $(\Omega_{Y/X})_{\hr y}$
is non-zero, and let $Z$ be the closed analytic subspace of $Y$ defined by the ideal $(f)$. 

\begin{enumerate}[b]

\item The morphism~$Z\to X$ is quasi-smooth of relative dimension~$d-1$ at~$y$. 

\item The coherent sheaves
$\mathscr O_Z$ and~$\Omega_{Z/X}$ are~$X$-flat at~$y$.
\end{enumerate}
\end{enumerate}
\end{lemm}

\begin{proof}
All properties involved can be checked G-locally,
hence we may and do assume that both $Y$ and $X$ are $k$-affinoid. Let
$x$ denote the image of $y$ in $X$.

As~$Y\to X$ is quasi-smooth of relative
dimension~$d$ at~$y$, the dimension of~$(\Omega_{Y/X})_{\hr y}$ is equal to~$d$. Let us choose global forms~$\omega_1,\ldots,\omega_d$ belonging to~$\Omega_{Y/X}(Y)$ such that~$(\omega_i(y))_i$ is a basis of ~$(\Omega_{Y/X})_{\hr y}$, 
which is possible
in view of~\ref{ss-pointwise-rk} (2). 
The~$\omega_i$'s define a morphism~$\mathscr O_Y^d\to \Omega_{Y/X}$. Since~$Y_x$ is geometrically regular at~$y$, the sheaf~$\Omega_{Y_x/\hr x}$ is free of rank~$d$ at~$y$; it follows therefore
from Nakayama's Lemma that $\mathscr O_{Y_x}^d\to  \Omega_{Y_x/\hr x}$ is an isomorphism at~$y$. Since $\Omega_{Y/X}$ is~$X$-flat at~$y$ by assumption, 
Lemma \ref{lem-iso-fiber}
implies that~$\mathscr O_Y^d\to \Omega_{Y/X}$ is an isomorphism at~$y$, whence (1).

The morphism $Y\to X$ being quasi-smooth at $y$, it admits a Jacobian presentation $W\hookrightarrow V\subset \A^n_X$ at~$y$. There exists a finite family~$(g_1,\ldots,g_{n-d})$ of analytic functions  on~$V$ such that the ideal~$(g_1,\ldots, g_{n-d})$ defines the closed immersion~$W\hookrightarrow V$, and such that the family~$(\d g_i(y))_i$ of elements of the vector space~$(\Omega_{V/X})_{\hr y}$ is free. As~$$(\Omega_{Y/X})_{\hr y}\simeq (\Omega_{V/X})_{\hr y}/((\d g_1)(y),\ldots, (\d g_{n-d})(y)),$$ the fact that~$(\d f)(y)$ is non-zero in ~$(\Omega_{Y/X})_{\hr y}$ simply means that the family
$((\d g_1)(y),\ldots, (\d g_{n-d})(y), (\d f)(y))$
is free in $(\Omega_{V/X})_{\hr y}$.
As the ideal~$(g_1,\ldots, g_{n-d}, f)$ defines precisely the closed immersion~$W\cap Z\hookrightarrow V$, the diagram~$(W\cap Z)\hookrightarrow V\subset \A^n_X$ is a Jacobian presentation of~$Z\to X$ at~$y$, and~$Z\to X$ is quasi-smooth at~$y$ of relative dimension~$d-1$, whence (2a).

Denote by $\iota$ the closed immersion~$Z\hookrightarrow Y$. We have an exact sequence

$$\xymatrix{ \mathscr O_Y\ar[rr]^{\times\; f}&&\mathscr O_Y\ar[r] &\mathscr \iota_\ast  \mathscr O_Z\ar[r]&0}.$$

By assumption,~$\mathscr O_Y$ is~$X$-flat at~$y$. As~$Y$ is quasi-smooth over~$X$ at~$y$, the 
fiber~$Y_x$ is (geometrically) regular at~$y$. Since we have made the hypothesis
that~$(\d f)(y)\neq 0$ in  $(\Omega_{Y/X})_{\hr y}$, 
the element~$f$ of the regular local ring~$\mathscr O_{Y_x,y}$ is non-zero.
As a regular local ring is a domain, 
the multiplication by~$f$ from~$\mathscr O_{Y_x}$ to itself is injective at~$y$. It now follows from
Lemma~\ref{lem-flattrick-reg}
that~$\iota_\ast  \mathscr O_Z$ is~$X$-flat at~$y$ ; in other words,~$\mathscr O_Z$ is~$X$-flat
at~$y$. 

Consider now the exact sequence
$$\xymatrix{
{\mathscr O_Z }\ar[rr]^(.4)u&& {(\Omega_{Y/X})_Z}\ar[r] & {\Omega_{Z/X}}\ar[r]&0},$$
where $u$ is the multiplication by the pullback of $\d f$.  We have just proven that~$\mathscr O_Z$
is~$X$-flat at~$y$. Since~$\Omega_{Y/X}$ has been seen to be free at~$y$, its 
pull-back~$(\Omega_{Y/X})_Z$ is free at~$y$ as well, hence
is $X$-flat at~$y$ because
so is $\mathscr O_Z$. 
Since~$(\d f)(y)\neq 0$ in  $(\Omega_{Y/X})_{\hr y}$, 
the image of~$\d f$ in~$(\Omega_{Y/X})_{Z_x, y}$ is non-zero. As~$Z$
is quasi-smooth over~$X$ at~$y$, the fiber~$Z_x$ is (geometrically) regular at~$y$; in particular,~$\mathscr O_{Z_x,y}$ is a domain, 
and the free $\mathscr O_{Z_x,y}$ -module
$(\Omega_{Y/X})_{Z_x, y}$ is thus
torsion-free. It
follows that $u$ induces an {\em injection}
from~$\mathscr O_{Z_x,y}$
into~$(\Omega_{Y/X})_{Z_x,y}$.  We then deduce
from Lemma~\ref{lem-flattrick-reg}
that~$\Omega_{Z/X}$ is~$X$-flat at~$y$, whence (2b). 
\end{proof}

\begin{coro}\label{cor-qsm-flat}
Let~$Y\to X$ be a morphism of~$k$-analytic spaces, let~$d$
be an element of $\N$,
and let~$y$ be a point of $Y$
at which~$Y\to X$ is quasi-smooth of relative dimension~$d$. The sheaf~$\mathscr O_Y$ is~$X$-flat at~$y$, 
and~$\Omega_{Y/X}$ is free of rank~$d$ at~$y$. 
\end{coro}

\begin{proof}
It is sufficient to prove that both coherent sheaves
$\mathscr O_Y$ 
and $\Omega_Y$ are $X$-flat at $y$; it will then follow from Lemma~\ref{lem-qsm-technical} (1)
that
$\Omega_{Y/X}$ is free of rank~$d$ at~$y$. 
Let us choose a Jacobian presentation~$W\hookrightarrow V\subset \A^n_X$ of~$Y\to X$ at~$y$, and set~$r=n-\dim_{\hr y}(\Omega_{Y/X})_{\hr y}$. By definition
of such a presentation,
there exists a family~$(g_1,\ldots, g_r)$ of analytic functions on~$V$ such that the ideal~$(g_1,\ldots, g_r)$ defines the closed immersion~$W\hookrightarrow V$, and such that the family
$(\d g_1)(y),\ldots, (\d g_r)(y))$ of elements of~$(\Omega_{V/X})_{\hr y}$ is free. 
For every~$i\in \{0,\ldots, r\}$, denote by~$V_i$ the closed analytic subspace of~$V$ defined by the ideal~$(g_1,\ldots, g_i)$; note that~$V_0=V$ and that~$V_r=W$. 
The map~$V\to X$ is quasi-smooth at~$y$; moreover, ~$\mathscr O_V$ and~$\Omega_{V/X}$ are~$X$-flat. Indeed, $\A^n_X$ is flat over $X$: this comes from Lemma \ref{lem-field-flat}, or from the fact that $\A^n_Z$ is flat over~$Z$ for every affinoid space $Z$ in view of Proposition~\ref{prop-gagaflat-ft}, because~$\A^n_{Z\al}$ is flat over~$Z\al$; as a
consequence, $V$ is flat over $X$. And since~$\Omega_{V/X}$ is a free~$\mathscr O_V$-module (with basis~$({\rm d}T_i)_i$), it is flat over~$X$ too. 

Now Lemma~\ref{lem-qsm-technical} together with a straightforward induction on $i$ shows that for every $i\in \{0,\ldots, r\}$, the space $V_i$ is quasi-smooth at $y$ and the coherent sheaves
$\mathscr O_{V_i}$ and~$\Omega_{V_i/X}$ are~$X$-flat at~$y$. By taking $i=r$ we get the expected statements. 
\end{proof}

\subsection{}\label{ss-qsm-equations}
Let~$Y\to X$ be a morphism of~$k$-{\em affinoid} spaces, let~$y$ be a point of $Y$ and let~$x$ be its image on~$X$. Let us assume that~$Y_x$ is geometrically regular, and that~$Y\to X$ is
flat at~$y$.
There exists~$n\in \N$ so that the morphism~$Y\to X$ goes through a closed immersion~$Y\hookrightarrow D\times_kX$,
where~$D$ is a closed~$n$-dimensional polydisc.
Set~$r=n- \rk y (\Omega_{Y/X})$ (note that $Y\to X$ is then of relative dimension $n-r$ at $y$), and let~$\mathscr I$ be the sheaf of ideals of $\mathscr O_{D\times_k X}$ that defines the closed immersion~$Y\hookrightarrow D\times_kX$. 
Let us choose
global sections $g_1,\ldots, g_r$ of $\mathscr I$
such that~$$(\Omega_{Y/X})_{\hr y}=(\Omega_{D\times_k X/X})_{\hr y}/((\d g_1)(y),\ldots, (\d g_r)(y)),$$
which is possible by \ref{ss-omega-zarclos}
and \ref{ss-pointwise-rk} (2); note that $((\d g_1)(y),\ldots, (\d g_r)(y))$ is then a free family of elements of the vector space $(\Omega_{D\times_k X/X})_{\hr y}$.
Let~$Z$ be the Zariski-closed subspace of~$D\times_kX$ defined by the ideal sheaf~$\mathscr J:=(g_1,\ldots,g_r)$; by construction,~$Y$ is a Zariski-closed subspace of~$Z$, and~$Z\to X$ is quasi-smooth at~$y$ of relative dimension~$n-r$. 

As~$Z\to X$ is quasi-smooth at~$y$ of relative dimension~$n-r$, the~$\hr x$-space~$Z_x$ is geometrically regular at~$y$ of relative dimension~$n-r$. In particular, there exists a connected affinoid neighborhood~$U$ of~$y$ in~$Z_x$ that is normal, connected and of dimension $n-r$.  
As the fiber $Y_x$ is geometrically regular at~$y$ and as we have the equalities $\rk y (\Omega_{Y_x/\hr x})=\rk y(\Omega_{Y/X})=n-r$, we can shrink~$U$ so that~$U\cap Y_x$ is of dimension $n-r$. The intersection $U\cap Y_x$ being a closed analytic subspace of the reduced, irreducible,~$(n-r)$-dimensional space~$U$, it coincides with~$U$.

The  natural surjection~$\mathscr O_{Z,y}\to \mathscr O_{Y,y}$ is bijective. Indeed, we have proven above that~$U=U\cap Y_x$, which implies that~$\mathscr O_{Z_x,y}\to \mathscr O_{Y_x,y}$ is a bijection;  and since~$\mathscr O_Y$ is
$X$-flat at~$y$ by assumption, it follows then from Lemma \ref{lem-iso-fiber} that~$\mathscr O_{Z,y}\to \mathscr O_{Y,y}$ is bijective. 

The bijectivity of $\mathscr O_{Z,y}\to \mathscr O_{Y,y}$ is equivalent to the bijectivity of $\mathscr J_y\to \mathscr I_y$. As a consequence, 
there exists an affinoid neighborhood~$V$ of~$y$ in~$D\times_kX$ such that the closed immersion
$V\cap Y\hookrightarrow V\cap Z$ is an isomorphism; note that~$V\cap Y$ is an affinoid neighborhood of~$y$ in~$Y$, and that~$V\cap Y\hookrightarrow V\subset \A^n_X$ is a Jacobian presentation of~$Y\to X$ at~$y$.

\begin{theo}
\label{thm-main-qsm}
Let~$Y\to X$ be a morphism of~$k$-analytic spaces, let~$y$ be a point of $Y$ and let~$x$ be its image on~$X$.

\begin{enumerate}[1]

\item The following are equivalent: 

\begin{enumerate}[j]
\item The morphism $Y\to X$ is quasi-smooth at~$y$.

\item The morphism $Y\to X$ is flat at~$y$, and the fiber $Y_x$ is geometrically regular at~$y$.

\item The morphism $Y\to X$ is flat at $y$, and $\rk y (\Omega_{Y/X})=\dim_y Y_x$. 

\end{enumerate}

\item If moreover~$Y$ and~$X$ are good, then those properties hold if and only if there exists a Jacobian presentation~$W\hookrightarrow V\subset \A^n_X$ of~$Y\to X$
 at~$y$ with~$W$ being an affinoid {\em neighborhood} of~$y$ in~$Y$.
\end{enumerate}
\end{theo}

\begin{proof}
We first remark that (ii)$\iff$(iii) by Remark \ref{rem-georeg-fiber}.
If~$Y\to X$ is quasi-smooth at~$x$, we already know that~$Y_x$ is geometrically regular at~$y$, and flatness of~$Y\to X$ at~$y$ is part of Corollary \ref{cor-qsm-flat}; so (i)$\Rightarrow$(ii).

Assume now that~$Y\to X$ is flat at~$y$, and that~$Y_x$ is geometrically regular at~$y$. In order to prove that~$Y\to X$ is quasi-smooth at~$y$, we may assume that both~$Y$ and~$X$ are~$k$-affinoid. But under that assumption, we have seen in \ref{ss-qsm-equations} that there exists a Jacobian presentation~$W\hookrightarrow V\subset \A^n_V$ of~$Y\to X$ at~$y$ with~$W$ being an affinoid neighborhood of~$y$ in~$Y$, which  at the same time ends the proof of (ii)$\Rightarrow$(i) and proves (2).
\end{proof}

\begin{coro}\label{cor-qsmfobs-qsmexc}
If~$X$ is a~$k$-analytic space and if~$x$ is a point of $X$, then~$X$ is quasi-smooth at~$x$ if and only if it is geometrically regular at~$x$.
\end{coro}

\begin{proof}
This is an immediate consequence of Theorem
\ref{thm-main-qsm} above, together with the fact that~$X\to \mathscr M(k)$ is flat (Lemma \ref{lem-field-flat}). 
\end{proof}

\begin{coro}\label{cor-gaga-qsm}
Let $\mathscr Y\to \mathscr X$ be a morphism of schemes of finite type over a given affinoid algebra, and let $y$
be  point of $\mathscr Y\an $. The map $\mathscr Y\an \to \mathscr X\an$ is quasi-smooth at $y$ if and only 
if $\mathscr Y\to \mathscr X$ is smooth at $y\al$.
\end{coro}

\begin{proof}
Let $x$ be the image of $y$ on $\mathscr X\an$. By Theorem~\ref{thm-main-qsm}, 
the
natural morphism $\mathscr Y\an \to \mathscr X\an$ is quasi-smooth at $y$ if and only if 
it is flat at $y$ and the dimension of the $\hr y$-vector space $(\Omega_{\mathscr Y\an/\mathscr X\an})_{\hr y}$ 
is equal to the relative dimension of $\mathscr Y\an$ over $\mathscr X\an$ at $y$.

On
the other hand, $\mathscr Y \to \mathscr X$ is smooth at $y\al $ if and only if 
it is flat at $y\al $ and the dimension of the $\kappa(y\al)$-vector space $(\Omega_{\mathscr Y/\mathscr X})_{\kappa(y\al)}$ 
is equal to the relative dimension of $\mathscr Y$ over $\mathscr X$ at $y\al$.

The corollary now follows from GAGA principles about flatness (Lemma \ref{lem-gagaflat-easy} and Proposition \ref{prop-gagaflat-ft}),
about the sheaf of differential forms (\ref{ss-omega-gaga}), and about relative dimension
(\Prop \ref{prop-gaga-dim}). 
 \end{proof}

\begin{coro}
\label{coro-qsm-open}
Let $Y\to X$ be a morphism between $k$-analytic spaces and let $d$ be an integer. The set of points of $Y$ at which 
$Y\to X$ is quasi-smooth of relative dimension $d$ is an open subset of $Y$.
\end{coro}

\begin{proof}
We immediately reduce to the case where both $Y$ and $X$ are affinoid. Let $y$ be a point of $Y$ at which $Y$
is quasi-smooth of relative dimension $d$ over $X$. By Theorem~\ref{thm-main-qsm}, there exists an affinoid neighborhood $W$
of $y$ in $Y$ and a Jacobian presentation $W\hookrightarrow V\subset \A^m_X$ of $Y\to X$ at $y$. Now by \ref{ss-jacob-pres}
there exists an affinoid neighborhood $V'$ of $y$
in $V$ such that $W\times_V V'\hookrightarrow V'$ is a Jacobian presentation of $Y\to X$
at each point of $W\times_V V'$, and such that $Y\to X$ is of relative dimension $d$ at each point of $W\times_V V'$. Hence
$Y\to X$ is quasi-smooth of relative dimension $d$ at each point of the affinoid neighborhood $W'\times_V V'$ of $y$. 
\end{proof}

\begin{rema}
We shall see later that the set of points of $Y$ at which 
$Y\to X$ is quasi-smooth of relative dimension $d$ is 
even \emph{Zariski}-open in $Y$ (Theorem \ref{thm-constloc-main}). 
\end{rema}

\section{Links with \'etale and smooth morphisms}\label{s-qsm-sm}

Our purpose is now to investigate the link between quasi-smooth morphisms, 
and smooth morphisms in the sense of Berkovich \cite {berkovich1993}, 3.5. As we shall see, 
{\em as far as the spaces involved are good}, the situation is very pleasant: a morphism is smooth
(at a given point of the source)
if and only if it is quasi-smooth and boundaryless; and it is quasi-smooth if and only if 
it is (locally) the composition of an inclusion of an analytic domain and of a smooth map. We shall
discuss thereafter this goodness assumption, and see that it is not needed in the 
case of relative dimension zero; \ie, for the comparison between quasi-étaleness and
étaleness. 

\begin{defi}\label{def-unramified}
Let $Y\to X$ be a morphism of $k$-analytic spaces and let $y$ be a point of $Y$. The morphism $Y\to X$ is
said to be
{\em unramified at $y$} (we shall also say that $Y$ is
{\em unramified over $X$} at $y$) if $(\Omega_{Y/X})_{\hr y}=0$; \ie, if $y$ does not belong to the support of $\Omega_{Y/X}$ (\ref{ss-pointwise-rk}). The morphism $Y\to X$
is called
{\em unramified}\index{morphism!of analytic spaces!unramified}\index{unramified morphism}
if it is unramified at every point of $Y$; \ie, if $\Omega_{Y/X}=0$
(we shall also say ``$Y$ is {\em unramified over $X$}"). 
\end{defi}

\subsection{}
Let $Y\to X$ be a morphism of $k$-analytic spaces and let $y$ be a point of $X$. If $Y\to X$ is unramified at $y$, then 
the inequality $\rk y (\Omega_{Y/X})\geq \dim_y(Y\to X)$ implies that $Y$ is of relative dimension zero over $X$ at $y$. 
As a consequence, it follows from Theorem \ref{thm-main-qsm} that $Y\to X$ is quasi-\'etale at $y$ if and only if it is 
flat and unramified at $y$.

\begin{lemm}\label{lem-unramified-diagonal}
Let $Y\to X$ be a morphism of $k$-analytic spaces. It is unramified if and only if 
the diagonal $Y\to Y\times_X Y$ is G-locally on $Y$ an open immersion
with closed (and open) image; \ie, there exists a G-covering $(Y_i)$ of $Y$ such that
for every $i$, the diagonal $Y_i\to Y_i\times_X Y_i$ is an open immersion with closed image. 
\end{lemm}

\begin{proof}
If $Y\to Y\times_X Y$ is G-locally on $Y$ an
open immersion, its conormal sheaf is trivial
and $\Omega_{Y/X}=0$. Conversely, assume that $\Omega_{Y/X}=0$, and let us prove that 
the diagonal $Y\to Y\times_X Y$
is G-locally on $Y$ an open immersion
with closed image. By arguing G-locally we reduce to the case were both $Y$ and $X$ are affinoid.
The diagonal map $Y\to Y\times_X Y$ is then a closed immersion, 
inducing a closed immersion $Y\al \hookrightarrow (Y\times_XY)\al$ at the scheme-theoretic level.
Since $\Omega_{Y/X}=0$ the conormal sheaf
of the closed immersion $Y\hookrightarrow Y\times_X Y$ is trivial, hence the conormal sheaf of the closed
immersion $Y\al \hookrightarrow ( Y\times_XY)\al$ is trivial
as well, which means that this closed immersion is
also an open immersion. Therefore $Y\hookrightarrow Y\times_X Y$
is an open immersion with closed image. 
\end{proof}

\begin{coro}\label{cor-qet-qet}
Let $$\xymatrix{
Z\ar[rr]\ar[rd]&&Y\ar[ld]\\
&X&
}$$be a commutative diagram in the category of $k$-analytic spaces. Let $z$ be a point of $Z$ and let $y$ be its image 
on $Y$. Assume that $Z$ is quasi-étale over $X$ at $z$, and that $Y$ is unramified over $X$ at $y$.
Then $Z$ is quasi-étale over $Y$
at $z$. 
\end{coro}

\begin{proof}
The quasi-étaleness locus of a morphism is an open subset of the source by Corollary \ref{coro-qsm-open}, and 
its unramified locus also (indeed, this is the complement of the support of a coherent sheaf). Hence we may shrink $Z$ and $Y$ around $z$
and $y$ respectively so that $Z$ becomes quasi-étale over $X$ and $Y$ becomes unramified over $X$. 

The morphism $Z\to Y$ is the composition of its graph $Z\hookrightarrow Z\times_X Y$ and the second
projection $Z\times_X Y\to Y$. Since $Z\to X$ is quasi-étale, $Z\times_X Y\to Y$ is quasi-étale too. And the graph 
$Z\hookrightarrow Z\times_X Y$ arises from the diagonal $Y\hookrightarrow Y\times_X Y$ through the base-change functor by $Z\to Y$
(the product $Y\times_X Y$ being seen as a $Y$-space through the first projection). Since $Y\to X$ is unramified, the diagonal
$Y\to Y\times_X Y$ is G-locally on $Y$ an open immersion by Lemma \ref{lem-unramified-diagonal}
above; in particular, it is quasi-étale. Therefore the graph $Z\to Z\times_X Y$ is quasi-étale too, and $Z\to X$ is quasi-étale
as the composition of two quasi-étale maps. 
\end{proof}

\begin{lemm}\label{lem-partial-qsm}
Let~$Y\to X$ be a morphism between~$k$-analytic spaces and let~$y$
be a point of $Y$. Assume that~$Y\to X$ is quasi-smooth at~$y$ of relative dimension~$d$,
and let ~$f_1,\ldots,f_\ell$ be analytic functions on~$Y$ such that~$((\d f_i)(y))_i$ is
a free family of elements of $(\Omega_{Y/X})_{\hr y}$
(note that one then has $\ell \leq d$). 
The map $Y\to \A^\ell_X$ defined by the~$f_i$'s is quasi-smooth of relative dimension~$d-\ell$
at $y$. 
\end{lemm}

\begin{proof}
One immediately reduces to the case where~$Y$ and~$X$ are affinoid.
Under that assumption, the
map~$Y \to \A^\ell_X$ goes through~$D\times_k X$ for some~$\ell$-dimensional compact polydisc~$D$, and it is sufficient to prove that~$Y\to D\times_kX$ is quasi-smooth of relative dimension~$d-\ell$ on an analytic neighborhood of $y$.

Both spaces~$Y$ and~$D\times_kX$ are~$k$-affinoid, hence~$Y\to D\times_kX$ goes through a closed immersion~$Y\hookrightarrow \Delta\times_kD\times_kX$ where~$\Delta$ is a closed polydisc. Let $\mathscr I$ be the
corresponding ideal sheaf on $\Delta\times_k D\times_k X$  and let $\delta$ be the dimension of $\Delta$.

By the choice of the~$f_i$'s, and in view of \ref{ss-omega-exactseq} applied to the diagram
$$Y\to D\times_k X\to X,$$ the rank $\rk y (\Omega_{Y/D\times_kX})$ is equal to $d-\ell$. As
a consequence and in view of \ref{ss-stalks-fibers}
and \ref{ss-omega-zarclos},
we can find $\delta-d+\ell$ global sections~$g_1,\ldots,g_{\delta-d+\ell}$
of $\mathscr I$ such that~$((\d g_j)(y))_j$ is a free family of elements of~$(\Omega_{\Delta\times_kD\times_kX/D\times_k X})_{\hr y}$;
it remains free when viewed as a family of vectors of~$(\Omega_{\Delta\times_kD\times_kX/X})_{\hr y}$, because the former vector space is a quotient of the latter
by \ref{ss-omega-exactseq}.~The quotient of $(\Omega_{\Delta\times_kD\times_kX/X})_{\hr y}$ by the $(\d g_j)(y)$'s has then dimension $d$; the natural surjection 
$$(\Omega_{\Delta\times_kD\times_kX/X})_{\hr y}/((\d g_j)(y))_j\to (\Omega_{Y/X})_{\hr y}$$ is thus an isomorphism. 

On the other hand,~$Y\to X$ is by assumption quasi-smooth at~$y$;
by Corollary \ref{cor-qsm-flat}, this implies that $Y\to X$ is
flat at~$y$, and
that $Y_x$ is geometrically regular at~$y$.

By the above and in view of \ref{ss-qsm-equations},
there exists an affinoid neighborhood~$V$ of~$y$ in~$\Delta\times_kD\times_kX$ such that the
closed analytic subspace~$Y\cap V$ of~$V$ is defined by the ideal~$(g_j)_j$; since ~$((\d g_j(y))_j$ is a free family of elements
$(\Omega_{\Delta\times_kD\times_kX/D\times_k X})_{\hr y}$, the map~$Y\to D\times_k X$ is quasi-smooth of relative dimension~$d-\ell$ at~$y$.
\end{proof}

\begin{theo}\label{thm-qsm-embeddsm}
Let~$Y\to X$ be a morphism of {\em good}~$k$-analytic spaces, and let~$y$
be a point of $Y$.
The following are equivalent : 

\begin{enumerate}[i]
\item There exists an affinoid neighborhood~$Y_0$ of~$y$ in~$Y$ and a smooth~$X$-space~$Z$ such that~$Y_0$ is~$X$-isomorphic to an affinoid domain of~$Z$. 

\item The morphism $Y\to X$ is quasi-smooth at~$y$.
\end{enumerate}
\end{theo}

\begin{proof}
As embeddings of analytic domains and smooth morphisms are quasi-smooth
(\ref{ss-id-qsm},
\ref{ss-sm-qsm}),  assertion (i) implies assertion (ii).
Let us now assume that (ii) is true. In order to prove (i), we may and do assume
that~$Y$ is affinoid. Let~$d$ be the dimension of~$(\Omega_{Y/X})_{\hr y}$ and let~$f_1,\ldots, f_d$ be analytic functions on~$Y$ such that the family~$((\d f_i)(y))_i$ is a basis of~$(\Omega_{Y/X})_{\hr y}$; let~$\phi\colon Y\to \A^d_X$ be the morphism induced by the~$f_i$'s. By Lemma~\ref{lem-partial-qsm}, the morphism~$\phi$ is quasi-\'etale at~$y$. 
Set $\xi=\phi(y)$. As~$\phi$ is
quasi-finite
at~$y$, the
analytic Zariski's Main Theorem
(\cite{ducros2007}, \Th 3.2)
ensures that~$Y$ can be shrunken so that~$\phi$ admits a factorization~$Y\to T_0\to T\to \A^d_X$ where~$T$ is finite \'etale over an open neighborhood~$U$ of~$\xi$,~$T_0$ is an affinoid domain of~$T$, and~$Y\to T_0$ is finite. 

{\em The finite morphism~$Y\to T_0$ is \'etale at~$y$.} Indeed, $T_0$ is quasi-étale over $\A^d_X$ because both maps $T_0\to T$ and $T\to \A^d_X$ are quasi-étale (the first one is the embedding of an analytic domain, and the second one is étale). Since $Y\to \A^d_X$ is also quasi-étale at $y$, it follows from 
Corollary \ref{cor-qet-qet} that $Y\to T_0$ is quasi-étale at $y$.  
This implies that~$Y\to T_0$ is flat at~$y$ and that~$(\Omega_{Y/T_0})_{\hr y}=0$; the map~$Y\to T_0$ being finite, those conditions exactly mean that it is \'etale at~$y$.

Let~$t$ be the image of~$y$ in~$T$. The categories of finite \'etale covers of the germ~$(T_0,t)$ and~$(T,t)$ are naturally equivalent (both are equivalent to the category of finite \'etale~$\hr t$-algebras, \cite{berkovich1993} \Th 3.4.1). Therefore there exists: 

\begin{itemize}[label=$\bullet$]

\item an open neighborhood~$T_1$ of~$t$ in~$T$; 

\item a finite \'etale map~$Z\to T_1$; 

\item an isomorphism between~$Z\times_{T_1}(T_1\cap T_0)$ and an open neighborhood~$Y_1$ of~$y$ in~$Y$. 

\end{itemize}
All morphisms in the diagram
$$Z\to T_1\to U\to \A^d_X\to X$$ are smooth; hence their composition $Z\to X$ is smooth. Now one can take~$Y_0$ equal to any affinoid neighborhood of~$y$ inside ~$Y_1$.
\end{proof}

\begin{rema}
In the strictly~$k$-analytic case, such a result has already been proved by Berkovich (\cite{berkovich1999}, Remark 9.7). 
\end{rema}

\begin{coro}\label{cor-qsm-sminner}
Let~$Y\to X$ be a morphism between {\em good}~$k$-analytic spaces and let~$y$ be a point of $Y$. The following are equivalent:

\begin{enumerate}[1]
\item The morphism ~$Y\to X$ is smooth at~$y$.  

\item The morphism $Y\to X$ is quasi-smooth and boundaryless at~$y$.
\end{enumerate}
\end{coro}

\begin{proof}
If $Y\to X$ is smooth at $y$,  it is quasi-smooth and boundaryless at $y$ (without a goodness assumption). 
Assume conversely that $Y\to X$ is quasi-smooth and boundaryless at $y$. Since $Y\to X$ is quasi-smooth at $y$ and since
$Y$ and $X$ are good, it follows from Theorem~\ref{thm-qsm-embeddsm}
that there exists an affinoid neighborhood $Y_0$ of $y$ in $Y$ such that $Y_0$ can be identified with an affinoid domain of a smooth
$X$-analytic space $Z$. Since $Y\to X$ is boundaryless at $Y$, the morphism $Y_0\to X$ is also boundaryless at $y$.
This imply that $Y_0\hookrightarrow Z$ is boundaryless as $y$; since $Y_0$ is an affinoid domain of $Z$, this means that $y$
belongs to the topological interior $U$ of $Y_0$ in $Z$. But $U$ is a neighborhood of $y$ in $Y$ and is smooth
over $X$ (as an open subset of $Z$); hence $Y\to X$ is smooth at $y$. 
\end{proof}

\begin{rema}\label{rem-qsm-goodness}
The author does not know if
Corollary \ref{cor-qsm-sminner} above is true without the goodness assumption. By looking carefully at what happens, the reader should be convinced that the main problem to face in the non-good case is the following: if~$Y\to X$ be a morphism of analytic spaces and if~$y\in Y$, there is no reason why there should exist analytic functions~$f_1,\ldots, f_r$ defined in a {\em neighborhood} of~$y$ such that~$((\d f_1)(y)),\ldots, (\d f_r)(y))$ generate~$(\Omega_{Y/X})_{\hr y}$. 

However, in the case where the relative dimension is zero (where the aforementioned problem vanishes), Corollary \ref{cor-qsm-sminner} is true without any goodness assumption.
Indeed, let us assume that~$Y\to X$ is quasi-\'etale and boundaryless at~$y$. Being zero-dimensional and boundaryless at~$y$, it is finite at~$y$
(\cite{berkovich1993}, \Prop 3.14),
hence we can shrink~$Y$ and~$X$ so that~$Y\to X$ is finite, and so that~$y$ is the only pre-image of its image $x$ in $X$.
Now choose a compact analytic neighborhood of~$x$ in $X$ of the form~$V_1\cup \ldots\cup V_m$ where the~$V_i$'s are affinoid domains of~$X$ containing $x$. For every
$i$ the map~$Y\times_X V_i\to V_i$ is finite; being quasi-\'etale at~$y$, it is in particular flat and unramified at~$y$, hence \'etale at~$y$, which is the only pre-image of~$x$. As a consequence, there exists an affinoid neighborhood~$W_i$ of~$x$ in~$V_i$ such that~$Y\times_X W_i\to W_i$ is \'etale. If one sets~$W=\bigcup W_i$,
then~$W$ is a compact analytic neighborhood of~$x$ and~$Y\times_XW\to W$ is \'etale, whence our claim. 

\end{rema}

In his work \cite{berkovich1994}
on vanishing cycles for formal schemes, Berkovich has defined
quasi-étale morphisms. We are going to check that our notion of quasi-étaleness
coincides with Berkovich's. 

\subsection{Berkovich's definition (\cite{berkovich1994}, \S 3)} Let $Y\to X$ be a morphism of $k$-analytic spaces and let $y$ be a point of $Y$. The morphism
$Y\to X$ is quasi-étale at $y$ in the sense of Berkovich if $y$ has a compact analytic neighborhood in $Y$ of the form $V_1\cup \ldots\cup V_n$ where
every $V_i$ is an affinoid domain of $Y$ that is $X$-isomorphic to an affinoid domain of an étale $X$-analytic space. 

\begin{lemm}\label{lem-qetad-qetvb}
Let $Y\to X$ be a morphism of $k$-analytic spaces and let $y$ be a point of $Y$. 
The following are equivalent: 

\begin{enumerate}[i]
\item  The morphism
$Y\to X$ is quasi-\'etale at~$y$ in the sense of Berkovich. 

\item The morphism~$Y\to X$ is quasi-\'etale at~$y$ in our sense. 
\end{enumerate}
\end{lemm}

\begin{proof}
In what follows, ``quasi-\'etale" will mean ``quasi-\'etale in our sense", and we will write ``quasi-\'etale in the sense of Berkovich" when needed. 
Let us assume (i). There exists in particular an affinoid domain~$V$ of~$Y$ containing~$y$ such that~$V$ can be identified with an affinoid domain of an analytic space~$X'$ which is \'etale over~$X$. As~$V\hookrightarrow X'$ and~$X'\to X$ are quasi-\'etale,~$V\to X$ is quasi-\'etale; in particular,~$V\to X$ is quasi-\'etale at~$y$, and~$Y\to X$
is therefore quasi-\'etale at~$y$.

Let us now assume (ii), and let~$x$ denote the image of~$y$ in~$X$.
Let us choose a compact analytic neighborhood of~$Y$ which is a finite union~$\bigcup V_i$ of affinoid domains of~$Y$ containing~$y$; we may assume that there exists for every~$i$ an affinoid domain~$U_i$ of~$X$ such that~$V_i\to X$ goes through~$U_i$. Fix~$i$. As~$Y\to X$ is quasi-\'etale at~$y$, the morphism~$V_i\to U_i$ is quasi-\'etale at~$y$ too (\ref{ss-andom-qsm}). Hence it follows from
Theorem
\ref{thm-qsm-embeddsm}
that there exists an affinoid neighborhood~$V'_i$ of~$y$ in~$V_i$ and an \'etale~$U_i$-space~$U'_i$ such that~$V'_i$ is $U_i$-isomorphic to an affinoid domain of~$U'_i$. The categories of finite \'etale covers of the germs~$(X,x)$ and~$(U_i,x)$ are naturally equivalent (both are equivalent to the category of finite \'etale~$\hr x$-algebras by \cite{berkovich1993}, \Th 3.4.1). Therefore there exists an open neighborhood~$X_i$ of~$x$ in~$X$ and a finite \'etale morphism~$X'_i\to X_i$ such that~$X'_i\times_X U_i$ can be identified with an open neighborhood of~$y$ in~$U'_i$; let us choose an affinoid neighborhood~$V''_i$ of~$y$ in~$V'_i$ such that~$V''_i\subset X'_i\times_X U_i\subset X'_i$. 
The union of the~$V''_i$'s is a neighborhood of~$y$, and for every~$i$ one can identify~$V''_i$ with an affinoid domain of the~$X$-\'etale space~$X'_i$; therefore~$Y\to X$ is quasi-\'etale at~$y$ in the sense of Berkovich. 
\end{proof}

\section{Transfer of algebraic properties}\label{s-transfer-qsm}

Let $Y\to X$ be a morphism of $k$-affinoid
spaces, let $y$ be a point of $Y$ and let $x$ be its image in $X$.
Assume that $Y$ is quasi-smooth over $X$
at $y$. 
By Theorem \ref{thm-main-qsm}, 
$Y$ is $X$-flat at $y$, and $Y_x$ is geometrically regular at $y$. The purpose of what follows is to establish 
algebraic avatars of this result. Namely, we shall prove that both morphisms
$\mathscr O_{X\al, x\al} \to \mathscr O_{Y\al, y\al}$ 
and $\mathscr O_{X,x}\to \mathscr O_{Y,y}$ are regular.
In fact, the general results established in
\ref{ss-fiberwise-chiant}
will enable us to deduce both statements from a particular case of the first one, 
which is the object of Lemma \ref{lem-transfer-sm}
below. 

\begin{lemm}\label{lem-transfer-sm}
Let $Y\to X$ be a quasi-smooth morphism between $k$-affinoid spaces, let $y$ be a point of $Y$ and let $x$
be its image in $X$. Assume that $X$ is integral
and that 
$x\al$ is the generic point of $X$. The scheme $Y\al$ is then regular at $y\al$.
\end{lemm}

\begin{proof}
By flatness of the map of locally ringed spaces $Y\to Y\al$ it is enough to prove that 
$Y$ is regular at $y$. For that purpose we may shrink $Y$ around $y$; hence using
Theorem
\ref{thm-qsm-embeddsm} we reduce to the case where $Y$ is $X$-isomorphic to an 
affinoid domain of some $X$-smooth space $X'$, and it suffices then
to prove that $X'$ is regular at
every point lying above $x$. 

Let $x'$ be such a point. There exists an open neighborhood $U$ of $x'$ in $X'$ and an étale 
$X$-morphism $U\to \A^n_X$ for some $n$; let $z$ be the image of $x'$ in $\A^n_X$. 
Since $x\al$ is the generic point of the integral scheme $X\al$, the local ring
$\mathscr O_{\A^n_{X\al}, z\al}$ coincides with $\mathscr O_{\A^n_{\kappa(x\al)}, z\al}$, hence is regular. 
In view of Lemma
\ref{gaga-concrete}, this implies that $\mathscr O_{\A^n_X,z}$ is regular.
Since 
$\mathscr O_{X',x'}$ is a finite étale $\mathscr O_{\A^n_X,z}$-algebra, it is regular too. 
\end{proof}

\begin{lemm}\label{lem-reg-sorital}
Let
$$\xymatrix{
{\mathscr Z}\ar[d]_f\ar[rd]^g&\\
{\mathscr Y}\ar[r]_h&{\mathscr  X}
}$$ be a commutative diagram of locally noetherian schemes. 
If $g$ is regular
and if 
$f$ is faithfully flat then~$h$
is regular.
\end{lemm}

\begin{proof}
The faithful flatness
of~$f$
and the flatness of $g$ imply the flatness of $h$ (one can check it directly, or see it as a particular case
of Lemma \ref{lem-flat-formal}). 
It remains to show that the fibers of~$h$ are geometrically regular.

Let~$x$ be a point of $\mathscr X$
and let~$L$ be a finite extension of~$\kappa(x)$. Since the map~$g$ is regular, the scheme
$\mathscr Z_{x,L}$ is regular. The morphism $\mathscr Z_{x,L}\to\mathscr Y_{x,L}$ is faithfully flat
because it is deduced from $f$ by base change along the map $\mathscr Y_{x,L}\to \mathscr Y$; as a consequence, $\mathscr Y_{x,L}$ is regular
(\cite{ega42}, \Prop 6.5.3 (i)). 
\end{proof}

\begin{theo}\label{thm-qsm-schemereg}
Let $Y\to X$ be a morphism between good $k$-analytic 
spaces, let $y$ be a point of $Y$ and let $x$ be its image in $X$. 
Assume that $Y\to X$ is quasi-smooth at $y$. 

\begin{enumerate}[1]
\item The morphism $\mathscr O_{X,x}\to \mathscr O_{Y,y}$ is regular. 
\item If moreover $Y$ and $X$ are affinoid, then $\mathscr O_{X\al, x\al}\to \mathscr O_{Y\al, y\al}$ is regular. 
\end{enumerate}

\end{theo}

\begin{proof}
We shall use the general abstract results of \ref{ss-fiberwise-chiant} by taking for $\mathfrak F$ the categeory $\mathfrak T$ itself (see \ref{s-category-framework}, and especially Example
\ref{ex-fiber-trivial}), for $\mathsf P$ the property of being regular, and for $\mathscr C$ the class of quasi-smooth morphisms. With this convention, 
Lemma \ref{lem-transfer-sm}
is nothing but assertion (A) of \ref{ss-fiberwise-chiant}, which
implies assertions ($\text A^\ast$), (B) and ($\text B^\ast$) of \loccit, as explained there. 
Now ($\text B^\ast$ ) implies that the fibers of
the map $\spec \mathscr O_{Y,y}\to \spec \mathscr O_{X,x}$ are geometrically regular; on the other hand, 
$\mathscr O_{Y,y}$ is a flat $\mathscr O_{X,x}$-algebra because $Y\to X$ is quasi-smooth at $y$, hence flat at $y$ by Corollary
\ref{cor-qsm-flat}; this ends the proof of (1).

Let us now come to (2). Since the quasi-smooth locus of $Y\to X$ is open by Corollary \ref{coro-qsm-open}, there exists an affinoid neighborhood
$V$ of $y$ in $Y$ such that the arrow
$V\to Y$ is quasi-smooth. We have mentioned above that assertion (B) of \ref{ss-fiberwise-chiant}
holds (with our choices of $\mathfrak F,\mathsf P$, and $\mathscr C$); this implies that he morphism 
$V\al \to X\al$ has geometrically regular fibers, and on the other hand it is flat because $V\to X$ is flat as a quasi-smooth morphism, and in view of Lemma
\ref{lem-gagaflat-easy}: moreover the morphism $V\al \to Y\al$ is regular by \ref{ss-algprop-analytic}
(2). We can thus apply Lemma \ref{lem-reg-sorital}
to the commutative diagram 
$$\xymatrix{
{\spec \mathscr O_{V\al, y_V\al}}\ar[d]\ar[rd]&\\
{\spec \mathscr O_{Y\al, y\al}}\ar[r]&{\spec \mathscr O_{X\al, x\al}}
}$$
and it yields the regularity of $\mathscr O_{X\al, x\al}\to \mathscr O_{Y\al, y\al}$. 
\end{proof}

Lemmas \ref{lem-flat-desc1}
and \ref{lem-flat-desc2} provide some descent and transfer results
for flat morphisms. Due to Theorem  \ref{thm-qsm-schemereg},
the transfer results alluded to above can be strengthened in the case
of a quasi-smooth map.
Again, we shall first write a statement that holds 
in the abstract settings
of \ref{s-category-framework} and \ref{s-alg-properties},
where we deal with general objects and properties of the latter satisfying various axioms; and then
we shall write down what it means for some {\em explicit} properties of interest. 
For the notion of validity of a property at a point, the reader may refer
to Lemma-Definition \ref{lem-equiv-valid} in our general 
abstract setting and to Lemma-Definition
\ref{valid-at-concrete}
for a more
concrete version.

\begin{prop}
\label{prop-transferqsm-general}
Let~$\mathfrak F$ be a fibered category as in \ref{s-category-framework}, and let $\mathsf  P$ be
a property as in \ref{ss-alg-properties}. Let us assume that $\mathsf P$ satisfies \hreg
(\ref{ss-list-hreg}). 
Let~$Y\to X$ be a morphism between~$k$-analytic spaces, let~$y$ be a point of $Y$, let~$x$ be its image on~$X$, and let~$D$
be an object of the fiber category $\mathfrak F_X$. Assume that $Y$
is quasi-smooth over $X$ at $x$. 
If~$D$ satisfies~$\mathsf P$ at~$x$, then
$D_Y$
satisfies $\mathsf P$ at~$y$.

\end{prop}

\begin{proof}

For both assertions, we can assume that
the spaces $Y$ and~$X$ are affinoid; now~$\spec \mathscr O_{Y,y}\to \spec \mathscr O_{X,x}$ is flat with (geometrically) regular fibers
in view of Theorem \ref{thm-qsm-schemereg}.The proposition then follows immediately from the fact that the property~$\mathsf P$ satisfies~\hreg~by
assumption.
\end{proof}

\begin{prop}[A concrete version of Proposition \ref{prop-transferqsm-general}]
\label{prop-tranferqsm-concrete}
Let~$Y\to X$ be a morphism between~$k$-analytic spaces, let~$y$ be a point of $X$ and let
$x$ be its image on~$X$. Let $\mathscr F$
be a coherent sheaf on~$X$, and let $m$ be an element of $\N$. Let $\mathsf S=\mathscr E\to \mathscr E'\to \mathscr E''$
be a complex of coherent sheaves on $X$. Assume that $Y$ is quasi-smooth
over $X$ at $y$.
If~$X$ is regular, \resp $R_m$, \resp Gorenstein, \resp CI at~$x$, 
so is~$Y$ at~$y$. If~$\mathscr F$ is CM, \resp $S_m$, \resp free of rank $m$ 
at~$x$, so is
$\mathscr F_Y$ at~$y$. Is $\mathsf S$ is exact at $x$, so is $\mathsf S_Y$ at $y$. 

\end{prop}

\chapter{Generic fibers in analytic geometry}\label{GEN}

In the study of relative properties in scheme
theory, a key role is
played by the technique of ``spreading out 
from the generic fiber", 
which is based upon the following obvious remark. Let $\mathscr Y\to \mathscr X$
be a morphism of schemes
and let $\xi$ be a point of $\mathscr X$ such that $\mathscr O_{\mathscr X,\xi}$ is a field (if $\mathscr X$
is noetherian and reduced, this amounts
to require that $\xi$ is the generic point of an irreducible component of $\mathscr X$).
Then for every point $y$ of the fiber $\mathscr Y_\xi$, 
the natural map $\mathscr O_{\mathscr Y,y}\to \mathscr O_{\mathscr Y_\xi,y}$ is an isomorphism. 

Let us now consider an analytic analogue of our situation: namely, $Y\to X$
is a morphism between good $k$-analytic spaces, $x$ 
is a point of $X$ such that $\mathscr O_{X,x}$ is a field, and $y$ is a point of $Y_x$.
Except in some very particular cases (\eg, if $x$ is rigid and $X=\{x\}$), the natural map
$\mathscr O_{Y,y}\to \mathscr O_{Y_x,y}$ cannot be expected to be an isomorphism,
because the formation of $Y_x$ involves several completion 
operations: indeed, if $X$ and $Y$ are affinoid, say $X=\mathscr M(A)$ and $Y=\mathscr M(B)$,
then $Y_x=\mathscr M(B\hotimes_A \hr x)$, and $\hr x$ is itself the completion 
of $\kappa(x)$.

But nevertheless, it does not
prevent us
from implementing a technique of ``spreading out from the generic fiber"
in analytic geometry. The point is the following: 
in the scheme-theoretic situation
we have described above, what actually matters
for carrying out this technique (as far as one is only
interested in the locus
of validity of the usual properties) is not the equality
$\mathscr O_{\mathscr Y,y}\to \mathscr O_{\mathscr Y_\xi,y}$
by itself, but the
slightly weaker
fact that every algebraic property of interest that holds over $\mathscr O_{\mathscr Y_\xi,y}$
descends to $\mathscr O_{\mathscr Y,y}$. 
And in this section, we prove the following (Theorem \ref{thm-localring-generic}): if $Y, X, y$, and $x$ are 
as above, \emph{and if moreover $y$ belongs
to $\mathrm{Int}(Y/X)$}, then $\spec \mathscr O_{Y_x,y}\to \spec \mathscr O_{Y,y}$
is flat (with CI fibers,
and even with regular fibers whenever \car $k=0$); this will be sufficient
to ensure descent 
of algebraic properties that descent along flat local maps
(\eg, all properties that satisfy condition \hreg~of
\ref{ss-list-hreg} and more concretely, all properties mentioned in Definition-Lemma \ref{valid-at-concrete}),
hence
to perform spreading out from the generic fiber -- and
we shall use it repeatedly
throughout the rest of this
memoir. 
Let us make some additional comments about this
result. 

\begin{itemize}[label=$\bullet$] 
\item The assumption that $y$ belongs to $\mathrm{Int}(Y/X)$ cannot
be dropped: we give a counter-example in \ref{ss-optimal-flat}
where $y$ belongs to $\partial (Y/X)$ and the morphism $\spec \mathscr O_{Y_x,y}\to \spec \mathscr O_{Y,y}$
is not flat. This means that for spreading out properties from generic fibers, one always will have
to reduce to the inner case. 

\item The fact that the fibers of $\spec \mathscr O_{Y_x,y}\to \spec \mathscr O_{Y,y}$
are CI, and regular if moreover \car $k=0$ (as soon as $y\in \mathrm{Int}(Y/X)$), is not used in this memoir; but
we get it almost for free by the method we use for proving flatness, and it seems us to be of independent interest. 
Let us mentions that the characteristic assumption on $k$ cannot be dropped as far as regularity is concerned, 
as witnessed by the counter-example \ref{ss-optimal-reduced} (which was communicated
to the author by Temkin). 

\item If $x$ is moreover an Abhyankar point (\ref{ss-abhyankar-points}), the map
 $\spec \mathscr O_{Y_x,y}\to \spec \mathscr O_{Y,y}$ is regular without having
 to assume that $y\in \mathrm{Int}(Y/X)$
 or  \car $k=0$ (Theorem \ref{thm-loc-gen-abhy}). This ensures that all properties satisfying \hreg~satisfy
descent and transfer between $\mathscr O_{Y,y}$ and $\mathscr O_{Y_x,y}$, and illustrates
 a general phenomenon: the best analogue of a generic fiber in analytic geometry (for a map between
 good spaces) is a fiber
 over an \emph{Abhyankar} point whose local ring is a field; \ie, an Abhyankar point at which the 
 target space is reduced (because the local ring of any Abhyankar point is artinian, \cf
 Example \ref{ex-centdim-recap}). 
 
 \end{itemize}

Before investigating local rings of generic fibers, we shall need some preparatory work. Section \ref{s-spread-prelim}
collects slightly technical (but easy) lemmas involving the completed residue field of a
point of an analytic space.
Section \ref{s-spread-polydiscs}
is devoted to a local study of smooth morphisms, which can be of independent interest. Let us say a few words about it. 

In complex analytic geometry, a morphism if smooth if and only if it is, locally on the source
and on the target, of the form $D\times X\to X$ for some open polydisc $D$. There is no hope for so
nice a description in the non-archimedean setting, even in the absolute context. Indeed, if
$k$ is algebraically closed and if $X$ is an irreducible, smooth,
projective curve over $k$ of positive genus, there always exists a point $x$ of $X\an$ that has {\em no}
neighborhood isomorphic to an open disc; nevertheless if $k$ is moreover non-trivially valued, any non-empty
smooth $k$-analytic space has a $k$-point, hence contains an open polydisc (choose a suitable open neighborhood of this point).
In Section \ref{s-spread-polydiscs}, we extend the latter absolute result to a relative situation. 
We prove the following (Proposition \ref{prop-nghbr-smooth}), which holds over an arbitrary analytic field $k$. 
Let $Y\to X$ be a smooth morphism 
of good $k$-analytic spaces. Let $x$ be a point of $X$
such that $\abs{\hr x\gpm}\neq \{1\}$
and $Y_x\neq \varnothing$; there exists an \'etale map $X'\to X$ whose image contains $x$
and an open subset of $Y':=Y\times_X X'$ that is $X'$-isomorphic to $D\times_k X'$ for some open polydisc $D$ (Proposition  \ref{prop-nghbr-smooth}
also gives a similar, but slightly more complicated result, when $\hr x$ is trivially valued). 

This almost immediately implies the openness of quasi-smooth boundaryless morphisms (Corollary \ref{cor-smooth-open}), which had already been proved by Berkovich with more 
sophisticated tools (see detailed comments at Remark \ref{rem-qsmopen-berk}), and also (at least when the ground field is non-trivially valued), the fact that 
a surjective smooth morphisms between good $k$-analytic spaces admits sections locally for the \'etale topology on the target (Corollary \ref{cor-qsm-sections}).

\section{Preliminary lemmas}\label{s-spread-prelim}

We begin with some general results about extensions of valued fields.

\begin{lemm}\label{lem-ext-ktt}
Let ~$F$ be a field and let~$L$ be a finite, separable extension of~$F(\!(t)\!)$. There exists a finite extension~$K$ of~$F$ such that~$L\otimes_F K$ admits a quotient isomorphic to~$K(\!(\tau)\!)$.
\end{lemm}

\begin{proof}
Let us consider~$F(\!(t)\!)$ as the completion of the function field of~$\P_F^1$ at the origin. Krasner's lemma ensures that there exists a projective, normal, irreducible~$F$-curve~$Y$ equipped with a finite, generically \'etale map to~$\P^1_F$, such that~$L$ can be identified with the completion of~$F(Y)$ at a closed point~$P$ lying above the origin. There exists a finite extension~$K$ of~$F$ such that the normalization~$Z$ of~$Y\times_F K$ is smooth and admits a~$K$-rational point~$Q$ over~$P$. Since the
completion of~$K(Z)$ at~$Q$ is isomorphic to~$K(\!(\tau)\!)$, the extension~$K$ satisfies the required property.
\end{proof}

\begin{lemm}\label{lem-valfield-dense}
Let~$K$ be a real valued field and let~$V$ be a subgroup of~$K$. Assume that there exists~$\rho\in (0,1)$ such that for every~$\lambda\in K$ there exists~$\mu\in V$ with 
$\abs{\lambda-\mu}\leq \rho \abs \lambda$. The group~$V$ is dense in~$K$.
\end{lemm}

\begin{proof}
For every~$\lambda \in K$
we choose an element~$\phi(\lambda)$ in~$V$ such that~$\abs{\lambda-\phi(\lambda)}\leq \rho \abs \lambda$. Now let~$\lambda$ be an element of $K$. Define inductively the sequence~$(\lambda_i)_i$ by setting~$\lambda_0=0$ and~$\lambda_{i+1}=\lambda_i+\phi(\lambda-\lambda_i)$. By induction, one sees that~$\lambda_i\in V$ and that~$\abs {\lambda-\lambda_i}\leq \rho^i\abs \lambda$ for every~$i$; hence~$\lambda_i\to \lambda$.
\end{proof}

\begin{lemm}
\label{lem-discreteval-finite}
Let~$K$ be a real valued field
such that $\abs {K\gpm}$ is free of rank 1, and let~$F$ be a complete subfield of~$K$ such that $\abs{F\gpm}\neq \{1\}$
and
the \emph{classical}
residue extension~$\widetilde F^1\hookrightarrow \widetilde K^1$ is
finite. The field~$K$ is a finite extension of~$F$.
\end{lemm}

\begin{proof}
The assumptions on the value groups ensures that~$\abs{K\gpm}/\abs{F\gpm}$ is finite; hence the {\em graded}
extension~$\widetilde F\hookrightarrow \widetilde K$ is finite too. Let~$\lambda_1,\ldots, \lambda_n$ be elements
of~$K\gpm$ such that~$(\widetilde {\lambda_i})_i$ is a basis of~$\widetilde K$ over~$\widetilde F$. Let us call~$V$
the~$F$-vector subspace of~$K$ generated by the~$\lambda_i$'s. Let~$\lambda$ be an element of $K$;
there exist~$a_1,\ldots, a_n\in F$ such that~$\widetilde\lambda =\sum \widetilde {a_i}\widetilde{\lambda_i}$, which
exactly means, if~$\lambda\neq 0$, that~$\abs{\lambda-\sum a_i\lambda_i}<\abs \lambda$. As~$\abs{K\gpm}$ is free of rank one,~$\abs{K\gpm}\cap (0,1)$ has a maximal element~$\rho$.
By the above, for every~$\lambda \in K$ there exists~$\mu\in V$
with~$\abs{\lambda-\mu}\leq \rho\abs \lambda$. By Lemma \ref{lem-valfield-dense}, the group $V$ is dense in~$K$. Since~$V$ is a finite dimensional~$F$-vector space, it is complete; hence~$V=K$.
\end{proof}

We are now going to use the above lemmas
to establish some
results related to the completed residue field of a point of
an analytic space.

\begin{lemm}
\label{lem-nullst-trivval}
Let~$F$ be a trivially valued field and let~$X$ be a non-empty, boundaryless~$F$-analytic space. There exists~$x\in X$ such that~$\hr x$ is either a finite extension of~$F$ or a finite extension of~$F_r$ for some~$r\in ]0;1[$.
\end{lemm}

\begin{proof}
Choose an arbitrary~$s\in ]0;1[$. As~$X_s$ is a non-empty, boundaryless space over the non-trivially valued field~$F_s$,
it 
has an $F_s$-rigid point $y$ (\ref{ss-analytic-nullst}); let~$x$ be the image of~$y$ on~$X$. Note that~$\hrt x^1$ is a subfield of~$\hrt y^1$, which is itself finite over~$\widetilde{F_s}^1=F$; hence~$\hrt x^1$ is finite over~$F$.

If~$\abs{\hr x\gpm}=\{1\}$, then as~$\hr x=\hrt x^1$, it is finite over~$F$ and we are done.

If~$\abs{\hr x\gpm}\neq \{1\}$, choose ~$r\in \abs{\hr x\gpm}\cap (0,1)$, and choose $\lambda\in \hr x \gpm$ such that~$\abs \lambda=r$. The complete subfield~$E$ generated by~$\lambda$ over~$F$ in~$\hr x$ is isomorphic to~$F_r$. As~$\hr x$ is a subfield of~$\hr y$, the non-trivial group~$\abs{\hr x\gpm}$ is free of rank 1; together with the fact that~$\widetilde E^1=\widetilde{F_r}^1=F$ this implies, in view of Lemma \ref{lem-discreteval-finite}, that~$\hr x$ is a finite extension of~$E\simeq F_r$. 
\end{proof}

\begin{lemm}
\label{lem-hx-kr}
Let~$r=(r_1,\ldots, r_n)$ be a~$k$-free polyradius and let~$S_1,\ldots, S_n$ be elements of~$k_r$ such that~$\abs {S_i}=r_i$ for every~$i$. The complete subfield~$L$ of~$k_r$ generated by
$k$ and the~$S_i$'s over~$k$ is equal to~$k_r$; in other words,~$S_1,\ldots, S_n$ are coordinate functions of~$k_r$.
\end{lemm}

\begin{proof}
Let~$T_1,\ldots, T_n$ be coordinate functions of~$k_r$; note that there is a well-defined $k$-isometry~$\phi: \sum a_I T^I\mapsto \sum a_I S^I$ between~$k_r$ and
$L$.

For every~$i$ one can write~$S_i=\alpha_iT_i+u_i$ where~$\alpha_i\in k$ and~$u_i\in k_r$, and where~$\abs{\alpha_i}=1$ and~$\abs{u_i}<r_i$. By replacing~$S_i$ with~$\alpha_i\inv S_i$, we may assume that~$\alpha_i=1$ for all~$i$. Therefore, there exists~$\rho\in (0,1)$ such that~$\abs{T_i-S_i}\leq \rho r_i$
and~$\abs{T_i^{-1}-S_i^{-1}}\leq \rho r_i^{-1}$
for every~$i$; it follows
that~$\abs{\lambda-\phi(\lambda)}\leq \rho \abs \lambda$ for every~$\lambda\in k_r$. Lemma \ref{lem-valfield-dense} then ensures that~$L$ is dense in~$k_r$; as $L$ is complete
we get $L=k_r$, as required.
\end{proof}

\begin{lemm}\label{lem-kr-omega}
Let~$Y$ be a quasi-smooth~$k$-analytic space and let~$y$ be a point of $Y$ such that~$\hr y\simeq k_r$  for some~$k$-free polyradius~$r=(r_1,\ldots,r_m)$
(as an analytic extension of $k$).
Let~$(g_1,\ldots,g_m)$ be analytic functions on~$Y$ such that~$\abs {g_i(y)}=r_i$ for every~$i$. The~$(\d g_i)(y)$'s are $\hr y$-linearly independent elements of~$(\Omega^1_{Y/k})_{\hr y}$. 
\end{lemm}

\begin{proof}
One can assume that~$Y$ is~$k$-affinoid and of pure dimension, say,~$n$. Let~$V$ be the affinoid domain of~$Y$ defined as the
locus of simultaneous validity of the equalities~$\abs{g_i}=r_i$. Its~$k$-affinoid structure factorizes through a~$k_r$-analytic
structure provided by the~$g_i$'s, for which~$y$ is~$k_r$-rational
by Lemma~\ref{lem-hx-kr}
above. For every $x\in V$ we have $d_{k_r}(x)=d_k(x)-d_k(k_r)=d_k(x)-m$
(\ref{ss-intro-etar}). Therefore
$$\dim_{k_r}V=\sup_{x\in V}d_{k_r}(x)=\sup_{x\in V}d_k(x)-m=\dim_k V-m=n-m.$$
As~$V$ is quasi-smooth,~$\mathscr O_{V,y}$ is regular; since~$y$ is~$k_r$-rational,
this implies that $(\Omega^1_{V/k_r})_{\hr y}$ is of dimension~$n-m$ (\ref{ss-omega-inequality}); otherwise said, 
$V$ is quasi-smooth {\em over~$k_r$} at~$y$. 
On the other hand, by quasi-smoothness of $Y$ (hence of $V$) over $k$, the~$\hr y$-vector space~$(\Omega^1_{V/k})_{\hr y}$ is of dimension~$n$. 

By \ref{ss-omega-exactseq}, we have an exact sequence
$(\Omega^1_{k_r/k})_V\to \Omega_{V/k}\to \Omega_{V/k_r}\to 0.$
Let $(T_i)$ 
denotes the family of coordinates functions of $k_r$. By \ref{ss-omega-affspace} (and \ref{ss-contravariant-omega})
the family $(\d T_i)_i$ is a basis of the $k_r$-vector space $\Omega_{k_r}/k$. Since the $k_r$-analytic
structure on $V$ arises from the morphism $(T_i\mapsto g_i)_i$, it follows from the above 
that 
$$(\Omega_{V/k_r})_{\hr y}=(\Omega_{V/k})_{\hr y}/((\d g_1)(y),\ldots, (\d g_m)(y)).$$ By considering the dimensions
of both sides, we see that $((\d g_1)(y),\ldots, (\d g_m)(y))$ is free
as a family of elements of $(\Omega_{V/k_r})_{\hr y}$, as announced. 
\end{proof}

\section{Relative polydiscs inside relative smooth spaces} 
\label{s-spread-polydiscs}

\begin{lemm}\label{lem-nghbr-affine}
Let~$X$ be a good~$k$-analytic space, let~$x$
be a point of $X$ and let~$n$ be an element of $\N$. Let~$m$ be a
non-negative integer $\leq n$, 
let~$r=(r_1,\ldots,r_m)$ be an~$\hr x$-free polyradius, and set~$r_i=0$ for~$m<i\leq n$. Let~$\xi$ be the point of~$\A^n_X$ lying above~$x$ and defined by the semi-norm~$$\sum a_I T^I\mapsto \max \abs{a_I} r^I$$ on the ring~$\hr x[T]$, and let~$V$ be an open neighborhood of~$\xi$ in~$\A^n_X$. The open subset $V$ of $\A^n_X$
contains an open neighborhood of~$\xi$ of the form~$$U\times_kD_1\times_k\ldots\times_k D_n$$ where~$U$ is an open neighborhood of~$x$ in~$X$ and $D_i$ is for every
integer $i\leq m$ (\resp $i>m$) a one-dimensional open annulus (\resp disc) with coordinate function~$T_i$.
\end{lemm}

\begin{proof}
Through a straightforward induction argument one immediately reduces to the case where~$n=1$; in that situation~$r$ is either zero
or
an~$\hr x$-free positive number. Let~$X_0$ be an affinoid neighborhood of~$x$ in~$X$ and set~$A=\mathscr O_X(X_0)$; let~$X'_0$ be the topological interior of~$X_0$ in~$X$. By the explicit description of the topology of the analytification of an~$A$-scheme of finite type
(\cf for instance
\cite{ducros2007b}, \S 1.4), there exist a finite family~$P_1,\ldots, P_\ell$ of elements of~$A[T]$, and a finite family~$I_1,\ldots, I_\ell$ of (relatively) open intervals of~$\R_+$ such that the open subset of~$\A^1_{X_0}$ defined by the conditions~$\abs {P_j}\in I_j$ (for~$j=1,\ldots, \ell$) contains~$\xi$ and is included in~$V$; write $P_j=\sum a_{i,j}T^i$. 

Assume that~$r=0$. In that case one has~$\abs{P_j(\xi)}=\abs{a_{0,j}(x)}$ for every~$j$. There exists for every~$j$ an open neighborhood~$I'_j$ of~$\abs {a_{0,j}(x)}$ in~$I_j$
and a positive number~$R_j$ such that~$\abs{P_j(\eta)}\in I_j$ as soon as~$\abs{a_{0,j}(\eta)}\in I'_j$ and~$\abs{T(\eta)}<R_j$. Let us denote by~$U$ the set of points~$y\in X'_0$ such that~$\abs{a_{0,j}(y)}\in I'_j$ for every~$j$, and let~$R$ be any positive number smaller than all~$R_j$'s.
The product of~$U$ and of the open disc centered at the origin with radius~$R$ is then included in~$V$ and contains~$\xi$, which ends the proof when~$r=0$. 

Assume that $r$ is an~$\hr x$-free positive number. In that case there exists for every~$j$ an index~$i_j$ such that~$\abs{a_{i_j,j}(x)}r^{i_j}>\abs{a_{i,j}(x)}r^i$ for all~$i\neq i_j$.
One can find for every~$j$ two positive numbers~$S_j$ and~$R_j$ with~$S_j<r<R_j$ and a family $(I'_{ij})_{i\leq \deg P_j}$ of (relatively) open intervals of $\R_+$,
each of which contains $\abs{a_{ij}(x)}$, such that~$\abs{P_j(\eta)}$ is equal to $\abs{a_{i_j,j}(\eta)}\cdot \abs{T(\eta)}^{i_j}$ and belongs to $I_j$  as soon as~$\abs{a_{i,j}(\eta)}\in I'_{ij}$ for all $i\leq \deg P_j$
and~$S_j<\abs{T(\eta)}<R_j$. Let us denote by~$U$ the set of points~$y\in X'_0$ such that~$\abs{a_{i,j}(y)}\in I'_{i,j}$ for every~$j$
and every $i\leq \deg P_j$, and let~$R$ and~$S$ be two positive number such that~$S<r<R$ and such that~$S_j<S$ and~$R<R_j$ for every~$j$. 
The product of~$U$ and of the open annulus with bi-radius~$(S,R)$ is then included in~$V$ and contains~$\xi$, which ends the proof.
\end{proof}

\begin{lemm}\label{lem-nghbr-smooth}
Let~$X$ be a good~$k$-analytic space, let~$x$ be a point of $X$, and let~$n$
be an element of $\N$; let~$m$ be 
a non-negative integer $\leq n$ and let~$r=(r_1,\ldots, r_m)$ be an~$\hr x$-free polyradius.
Let~$Y\to X$ be a smooth morphism of relative dimension~$n$, and assume that~$Y_x$ contains a point~$y$ with~$\hr y\simeq \hr x_r$
as analytic extensions of $\hr x$.
There exists an open subset
$V$
of~$Y$ which is~$X$-isomorphic to~$U\times_k D\times_k \Delta,$ where~$U$ is an open neighborhood of~$x$ in~$X$, $D$ is an~$m$-dimensional open poly-annulus,
and~$\Delta$ is an~$(n-m)$-dimensional open polydisc.
\end{lemm}

\begin{rema}\label{rem-polydisc-nocenter}
We emphasize the fact that (contrary to Lemma~\ref{lem-nghbr-affine}
with $\xi$), we do {\em not}
require the open subset $V$ to contain $y$, and our proof actually
does not enable us to achieve it: as the reader
will see, we need at some point to replace $y$ with a suitable 
approximation $z$, and the neighborhood we seek will be built
around $z$, and it might avoid $y$. (We shall encounter a similar restriction while dealing with \'etale multisections; see Remark 
\ref{rem-sections-around}
below.)

But this should not cause any trouble. Indeed, as explained
in the introduction of this chapter, we are interested in exhibiting nice open subsets inside relative
smooth spaces, but there is no hope for having a {\em basis}
of such open subsets.

\end{rema}

\begin{proof}[Proof of Lemma \ref{lem-nghbr-smooth}]
Let us choose analytic functions~$f_1,\ldots, f_m$ defined in an open
neighborhood of~$y$ such that~$\abs{f_i(y)}=r_i$ for every~$i$. According to Lemma \ref{lem-kr-omega}
(which one applies to the~$\hr x$-analytic space~$Y_x$), the elements~$(\d f_1)(y),\ldots, (\d f_m)(y)$ are linearly independent in~$(\Omega_{Y/X})_{\hr y}$ ; one can hence choose~$f_{m+1},\ldots, f_n$ in~$\mathscr O_{Y,y}$ so that ~$((\d f_1)(y),\ldots, (\d f_m)(y),(\d f_{m+1})(y),\ldots, (\d f_n)(y))$
 is a {\em basis} of~$(\Omega_{Y/X})_{ \hr y}$.

The~$X$-morphism~$Y\to \A^n_X$ induced by the~$f_i$'s is quasi-\'etale at~$y$ by Lemma \ref{lem-partial-qsm}; as~$Y\to X$ is boundaryless (it is smooth),~$Y\to \A^n_X$ is boundaryless at~$y$, and thus \'etale (Remark \ref{rem-qsm-goodness}). Lemma \ref{lem-hx-kr} ensures that the complete subfield of~$\hr y$ generated
over $\hr x$ by the~$f_i(y)$'s for~$i=1,\ldots,m$ is equal to~$\hr y$ itself. Therefore, if~$y'$ denotes the image of~$y$ on~$\A^n_X$, one has~$\hr {y'}=\hr y$; as~$Y\to \A^n_X$ is \'etale at~$y$, this implies that~$Y\to X$ induces an isomorphism between an open neighborhood of~$y$ in~$Y$ and an open neighborhood of~$y'$ in~$\A^n_X$
(\cite{berkovich1993}, \Th 3.4.1).
We thus may reduce to the case where~$Y$ is an open subset of~$\A^n_X$, and where the following holds: for the projection~$p\colon \A^n_X\to \A^m_X$ defined by~$(T_1,\ldots, T_m)$, the point~$p(y)$ of $\A^{m,\mathrm {an}}_{\hr x}$ is the one that corresponds
to the
Gauß norm~$\sum a_IT^I\mapsto \max \abs{a_I}r^I$, and~$y$ is an~$\hr {p(y)}$-rational point of the fiber~$p\inv (p(y))$.

Let~$\kappa$ be the residue field of~$\mathscr O_{\A^m_X,p(y)}$; by density of~$\kappa$ inside~$\hr {p(y)}$, the fiber~$(p_{|Y})\inv (p(y))$ possesses an~$\hr {p(y)}$-point~$z$ such that~$T_i(z)\in \kappa$ for every integer~$i\in \{m+1,\ldots, n\}$. Let~$V$ be an open neighborhood of~$p(y)$ in~$\A^m_X$ on which the~$T_i(z)$'s are defined; translation by~$(0,\ldots,0, T_{m+1}(z),\ldots, T_n(z))$ identifies over~$V$ the space~$Y\times_{\A^m_X}V$ with an open subset of~$\A^n_X\times_{\A^m_X}V$ whose fiber over~$p(y)$ contains the origin of~$\A^{m,\mathrm an}_{\hr {p(y)}}$. It follows then from Lemma \ref{lem-nghbr-affine} that there exists an open subset~$W$ of~$Y$ which is~$V$-isomorphic to~$V'\times \Delta$ where~$V'$ is an open neighborhood of~$p(y)$ in~$V$ and where~$\Delta$ is a~$(n-m)$-dimensional open polydisc.

Now by applying once again Lemma \ref{lem-nghbr-affine}, but this time to the map~$\A^m_X\to X$ and at the point~$p(y)$, one sees that there exists an open neighborhood~$V''$ of~$p(y)$ in~$V'$ that is~$X$-isomorphic to~$U\times_kD$ for some open neighborhood~$U$ of~$x$ in~$X$ and some~$m$-dimensional open poly-annulus~$D$. Now~$W\times_{V'}V''$ is an open subset of~$Y$ that is~$X$-isomorphic to 
$U\times_kD\times_k\Delta$, as required.
\end{proof}

We are now ready to investigate the local structure of arbitrary smooth morphisms
between good analytic spaces. 

\begin{prop}\label{prop-nghbr-smooth}
Let~$n$ be an element of $\N$ and let~$Y\to X$ be a smooth map of pure
relative
dimension~$n$ between good~$k$-analytic spaces. Let~$x$
be a point of $X$, and let~$W$ be a non-empty open subset of~$Y_x$. There exist: 

\begin{itemize}[label=$\bullet$]
\item a flat, locally finite morphism~$X'\to X$  which can be chosen to be \'etale if~$\abs{\hr x\gpm}\neq\{1\}$; 

\item a pre-image~$x'$ of~$x$ on~$X'$, such that the morphism~$\spec \mathscr O_{X',x'}\to \spec \mathscr O_{X,x}$
has a reduced closed fiber
(note that this is kind of a weak substitute to \'etaleness in the case where $\hr x$ is trivially valued);

\item an open subset~$V$ of~$Y':=Y\times_XX'$ whose intersection with 
$Y'_{x'}$ is contained in~$W_{\hr {x'}}$, and which is
$X'$-isomorphic to~$X'\times_k D$ where~$D$ is:

\begin{itemize}[label=$\diamond$]
\item an~$n$-dimensional open polydisc if~$W$ has an~$\hr x$-rigid point, which is always the case if~$\abs{\hr x\gpm}\neq\{1\}$ or if~$n=0$; 

\item the product of a 1-dimensional open annulus
and an~$(n-1)$-dimensional open polydisc if~$W$ has a no~$\hr x$-rigid point, which can occur only if~$\hr x$ is trivially valued and~$n>0$.
\end{itemize}
\end{itemize}
\end{prop}

\begin{proof} By replacing~$Y$ with an open subset of~$Y$ whose intersection with~$Y_x$ is equal to~$W$, we may assume that~$W=Y_x$. By assumption,~$Y_x\neq \emptyset$.

Let us assume that~$\abs{\hr x\gpm}\neq \{1\}$. Let us choose~$y\in Y_x$. As~$Y\to X$ is smooth, there exists a neighborhood~$Z$ of~$y$ in~$Y$ such that~$Z\to X$ goes through an \'etale map~$Z\to \A^n_X$; the image of~$Z$ on~$\A^n_X$ is an open subset~$U$ of the latter, and~$U_x$ is non-empty.  Let~$K$ be the completion of an algebraic closure of~$\hr x$. Since~$\hr x$ is not trivially valued, the analytic {\em Nullstellensatz} ensures that~$U_x(K)\neq \emptyset$(\ref{ss-analytic-nullst}); as the separable closure of~$\hr x$ in~$K$ is dense (again because~$\abs{\hr x\gpm}\neq \{1\}$), there exists~$u\in U_x$ with~$\hr u$ a finite, separable extension of~$\hr x$. Now let us choose a pre-image~$z$ of~$u$ on~$Z$; as~$Z\to U$ is \'etale,~$\hr z$ is a finite, separable extension of~$\hr u$, and hence a finite, separable extension of~$\hr x$ too. 
The categories of finite \'etale covers of the germ~$(X,x)$ and of finite \'etale~$\hr x$-algebras being naturally equivalent (\cite{berkovich1993}, \Th 3.4.1), there exists an \'etale morphism~$X'\to X$ and a pre-image~$x'$ of~$x$ on~$X$ such that~$z$ has a pre-image~$z'$ on~$Y':=Y\times_X X'$ with~$\hr {z'}=\hr {x'}$. This implies, in view of Lemma \ref{lem-nghbr-smooth}, that one can shrink~$X'$ around~$x'$ so that~$Y'$ possesses an open subset which is~$X'$-isomorphic to the product of~$X'$
and an~$n$-dimensional open polydisc. This ends the proof in the case where~$\abs{\hr x\gpm}\neq \{1\}$.

Let us assume that $\abs{\hr x\gpm}= \{1\}$ and that $Y_x$ has an~$\hr x$-rigid point~$y$. As~$\hr x$ is trivially valued, it coincides with the residue field~$\kappa$ of~$\mathscr O_{X,x}$. Therefore, there exists a finite, flat, local~$\mathscr O_{X,x}$-algebra~$A$ with~$A\otimes_{\mathscr O_{X,x}}\kappa\simeq\hr y$ (over $\kappa$). Since $X$ is good, one
can find a locally finite, flat map~$X'\to X$ and a pre-image~$x'$ of~$x$ on~$X'$ such that~$\mathscr  O_{X',x'}\simeq A$; note that~$\hr {x'}\simeq \hr y$ and that the closed fiber of~$\spec \mathscr O_{X',x'}\to \spec \mathscr O_{X,x}$ is reduced. By construction,~$y$ has a pre-image~$y'$ on~$Y':=Y\times_XX'$ lying above~$x'$ and satisfying~$\hr {y'}=\hr {x'}$. This implies, in view of Lemma \ref{lem-nghbr-smooth}, that one can shrink~$X'$ around~$x'$ so that~$Y'$ possesses an open subset which is~$X'$-isomorphic to the product of~$X'$ and an~$n$-dimensional open polydisc. This ends the proof in the case where~$\abs{\hr x\gpm}= \{1\}$ and where~$Y_x$ has a rigid point.

Let us assume that $\abs{\hr x\gpm}= \{1\}$ and that $Y_x$ has no~$\hr x$-rigid point. In this case, there exists~$t\in Y_x$ and~$r\in (0,1)$ such that~$\hr t$ is a finite extension of~$\hr x_r$ (Lemma \ref{lem-nullst-trivval}). Due to Lemma \ref{lem-discreteval-finite}, there exists a finite~$\hr x$-extension~$F$ such that~$F\otimes_{\hr x}\hr t$ admits a quotient $\hr x$-isomorphic (as a valued extension
of $F$)
to~$F_s$ for some~$s\in (0,1)$. As~$\hr x$ is trivially valued, it coincides with the residue field $\kappa(x)$
of~$\mathscr O_{X,x}$. Therefore, there exists a finite, flat, local~$\mathscr O_{X,x}$-algebra~$A$ with~$A\otimes_{\mathscr O_{X,x}}\kappa(x)\simeq F$. One can find a locally finite, flat map~$X'\to X$ and a pre-image~$x'$ of~$x$ on~$X'$ such that~$\mathscr  O_{X',x'}\simeq A$; note that~$\hr {x'}\simeq F$ and that the closed fiber of~$\spec \mathscr O_{X',x'}\to \spec \mathscr O_{X,x}$ is reduced. By construction,~$y$ has a pre-image~$y'$ on~$Y':=Y\times_XX'$ lying above~$x'$ and such that~$\hr {y'}$ is $\hr {x'}$-isomorphic (as a valued field)
to~$\hr {x'}_s$ for some~$s\in (0,1)$. This implies, in view of Lemma \ref{lem-nghbr-smooth} that one can shrink~$X'$ around~$x'$ so that~$Y'$ possesses an open subset which is~$X'$-isomorphic to the product of~$X'$
and an~$(n-1)$-dimensional open polydisc and 
a
1-dimensional open annulus, which ends the proof.
\end{proof}

\begin{coro}\label{cor-smooth-open}
Any quasi-smooth, boundaryless map is open. 
\end{coro}

\begin{proof}
Let~$Y\to X$ be a quasi-smooth, boundaryless map. To prove that it is open, one may argue G-locally on~$X$, hence assume that~$X$ is good. Now~$Y$ is also good because~$Y\to X$ is boundaryless, and it follows from Corollary \ref{cor-qsm-sminner} that~$Y\to X$ is smooth. Its openness then follows immediately from Proposition \ref{prop-nghbr-smooth} above together with openness of flat, locally finite morphisms (Corollary \ref{cor-finflat-open}).  
\end{proof}

\begin{rema}\label{rem-qsmopen-berk}
Corollary \ref{cor-smooth-open} had already been proven by Berkovich (\cite{berkovich1993}, \Cor 3.7.4).
Strictly speaking, Berkovich
proved the openness of smooth maps, and not of all quasi-smooth, boundaryless maps -- and the latter class is likely broader than the former
(see Remark
\ref{rem-qsm-goodness}); but openness can be checked G-locally on the target, hence in the good case where both notions coincide (Corollary \ref{cor-qsm-sminner}).

Our proof is different from Berkovich's. The latter used the d\'evissage of a smooth morphism in ``elementary"
curve fibrations, whose existence comes from the semi-stable reduction theorem. Ours involves less sophisticated tools: it 
is essentially based upon easy explicit computations (Lemma \ref{lem-nghbr-affine}) and the fact that
if $(Y,y)\to (X,x)$ is an \'etale morphism  of analytic germs such that $\hr x\simeq \hr y$, it is an isomorphism
(\cite{berkovich1993}, \Th 3.4.1).
\end{rema}

\begin{coro}\label{cor-qsm-sections}
Let~$Y\to X$ be a smooth morphism between good~$k$-analytic spaces and let~$x$
be a point of $X$ such that~$\abs{\hr x\gpm}\neq\{1\}$. If~$Y_x\neq \emptyset$, there exists an \'etale morphism~$ X'\to X$ whose image contains~$x$ and an~$X$-map~$X'\to Y$.
\end{coro}

%\begin{coro}\label{cor-qsm-dimension}
%Let~$n$ and~$d$ be two integers and let~$\phi: Y\to X$ be a quasi-smooth morphism between~$k$-analytic spaces. Assume that~$X$ is purely~$d$-dimensional, and that~$\phi$ is of pure relative dimension~$n$. The space~$Y$ is purely~$n+d$-dimensional.
%\end{coro}

%\begin{proof}
%One can argue G-locally on~$Y$ and~$X$ ; therefore, one can assume first, that~$X$ and~$Y$ are good, and then (Theorem \ref{thm-qsm-embeddsm}) that~$Y$ is an analytic domain of a smooth~$X$-space of pure relative dimension~$n$. Eventually  (by replacing~$Y$ with the latter), we reduce to the case where~$\phi$ itself is smooth. Now if~$V$ is any non-empty open subset of~$Y$, corollary \ref{smoothopen} above ensures that~$\phi(V)$ is an open subset of~$X$. It is therefore of dimension~$d$. By \ref{imanddim},~$\dim{}V=n+d$; by \ref{remdimopen},~$Y$ is purely of dimension~$n+d$.~$\Box$ 

\begin{rema}\label{rem-sections-around}
Let us keep using the notation
of Corollary \ref{cor-qsm-sections}
above. Let $y$ be any rigid point of the fiber $Y_x$ such that $\hr y$ is separable over $\hr x$. 
We do {\em not}
claim
that we can build an \'etale multisection of $Y\to X$ around $x$ that goes through $y$. And in fact, 
such a multisection {\em does not exist in general}. Indeed, let $X'$ be as in Corollary \ref{cor-qsm-sections},
let $x'$ be a pre-image of $x$ in $X'$, and let $y$ be the image of $x'$ in $Y$. The morphism $X'\to Y$ is
quasi-finite and boundaryless at $x'$, 
hence finite at $x'$
(\cite{berkovich1993}, \Cor 3.1.10).
Therefore
$$\mathrm{centdim}(Y,y)=\mathrm{centdim}(X',x')=\mathrm{centdim}(X,x)$$
(both equalities follow from \ref{ss-basicprop-centdim}).

Now let $f$ be any power series belonging to $k[\![t]\!]$ whose radius
of convergence
is a positive real number $r$, let $D$
be the closed $k$-analytic disc of radius $r$, 
and let $\phi$ be the morphism $(\mathrm{Id},f)$ from $D$ to $\A^{2,\mathrm{an}}_k$. Let $x$ be the Shilov point of $D$, and set $y=\phi(x)$. If $p$
denotes the first projection $\A^{2,\mathrm{an}}\to  \A^{1,\mathrm{an}}$, then $p(y)=x$ and $\hr y=\hr x$. 
By the proof of Proposition \ref{prop-oanx-field} (be aware that our current point $y$ is denoted by $x$ in \loccit), $\mathrm{centdim}(\A^{2,\mathrm{an}},y)=2$. 
Therefore $\mathrm{centdim}(\A^{2,\mathrm{an}},y)\neq \mathrm{centdim}(\A^{1,\mathrm{an}},x)$ because the latter is $\leq 1$ (it is in fact equal to $1$ because $x$
is not rigid, but we do not need that). Hence by the above, $p$ does not admit any \'etale
multisection around $x$ and going through $y$.

\end{rema}

\section{Local rings of generic fibers}
\label{s-spread-general}

We begin this section
with a lemma that ensures kind of a ``spreading out
from the generic fiber" of a topological property, in the particular case
 of a smooth morphism; it will play a key role in our study of local rings 
 of generic fibers of an arbitrary map at an inner point
 (Theorem \ref{thm-localring-generic}). 

\begin{prop}
\label{prop-nghbr-generic}
Let~$Y\to X$ be a map between good~$k$-analytic spaces. Let~$y$ be a point of $Y$
and~let~$x$ be its image in~$X$. Assume that ~$Y\to X$ is smooth at~$y$ and
the local ring ${\mathscr O}_{X,x}$ is artinian. Let
$Z$ be a Zariski-closed subset of~$Y$
that contains a neighborhood of~$y$ in~$Y_x$. Under those assumptions,
$Z$ is a neighborhood of~$y$ in~$Y$.
\end{prop}

\begin{proof}
We can shrink $Y$
around $y$ so that $Y\to X$
is smooth.
The required property being purely topological, one may assume that~$X$ is reduced; in this
case,~$\mathscr O_{X,x}$ is a field, and is in particular normal.
The normality locus of $X$ is (Zariski-)open
by Lemma \ref{locus-concrete} (1). We thus may shrink~$X$ so that it is itself normal. Now
the~$X$-quasi-smooth space~$Y$ is normal too, in view of Proposition \ref{prop-tranferqsm-concrete};
by shrinking~$Y$ (and~$Z$, accordingly) we eventually reduce to the case where~$Y$ is connected, hence irreducible
(and equi-dimensional) and where~$Z$ is the zero-locus of a finite family~$(f_1,\ldots, f_n)$ of analytic functions on~$Y$. We
shall prove that~$Z$ contains a non-empty open subset of~$Y$, which will force it to coincide with~$Y$, and end the proof.

By Proposition~\ref{prop-nghbr-smooth}, there exists a flat, locally finite map~$X'\to X$, a point~$x'$ on~$X'$ lying above~$x$, and a~$k$-analytic space~$D$ such that:

\begin{itemize}[label=$\bullet$] 

\item $\mathscr O_{X',x'}$ is a field ; 
\item $D$ is an open polydisc, or the product of an open polydisc and a 1-dimensional open annulus;
\item if one sets~$Y'=Y\times_X X'$ and~$Z'=Z\times_X X'$, there exists an open subset~$V$ of~$Y'$ which is~$X'$-isomorphic to~$D\times_kX'$ and such that~$V_{x'}\subset Z'_{x'}$. 

\end{itemize}

As~$\mathscr O_{X',x'}$ is a field, we may shrink~$X'$ so that it is connected and normal, by the same reasoning as above. 
Let us still denote by~$f_1,\ldots, f_n$ the pull-backs of the~$f_i$'s on~$Y'$. Analytic functions on~$V\simeq D\times_kX'$ consist
of
power series~$\sum a_I T^I$ where the~$a_I$'s are analytic functions on~$X'$. For any~$j$, let us write~$f_{j|V}=\sum a_{I,j} T^I$. By construction,~$V_{x'}\subset Z'_{x'}.$ Therefore~$a_{I,j}(x')=0$ for every~$(I,j)$. As~$\mathscr O_{X',x'}$ is a field,~$a_{I,j}$ vanishes
in a neighborhood of~$x'$ in the normal
connected space~$X'$ for every~$(I,j)$; therefore~$a_{I,j}=0$ for every~$(I,j)$. This implies that~$V\subset Z'$; hence~$Z'$ contains a non-empty open subset of~$Y'$. As~$Y'\to Y$ is flat
and locally finite, it is open by Corollary  \ref{cor-finflat-open}.
Therefore~$Z$ contains a non-empty open subset of~$Y$, which we have seen is sufficient. 
\end{proof}

\subsection{}\label{ss-property-r}
Let~$\mathscr X$ be a locally noetherian scheme over a field~$\kappa$. For technical purposes,
we shall have to consider the following property~$\mathfrak R$: {\em there exist a subfield~$\kappa_0$
of~$\kappa$ with~$[\kappa:\kappa_0]<+\infty$ and a {\em regular}~$\kappa_0$-scheme~$\mathscr X_0$
such that~$\mathscr X\simeq \mathscr X_0\otimes_{\kappa_0}\kappa$.} Let us mention some basic facts
about this property.

\begin{itemize}[label=$\bullet$]

\item It follows from the 
definition that if the~$\kappa$-scheme~$\mathscr X$ satisfies~$\mathfrak  R$, then so does 
the~$\lambda$-scheme $\mathscr X\otimes_\kappa \lambda$ for every finite extension~$\lambda$ of~$\kappa$. 

\item If~$\kappa$ is of char.~$0$, then
it is immediate that~$\mathscr X$ satisfies~$\mathfrak R$ if and only if it is regular. 

\item Assume that the~$\kappa$-scheme $\mathscr X$
satisfies~$\mathfrak R$. Let us choose a closed immersion
of~$\spec \kappa$ into~$\A^m_{\kappa_0}$ for some~$m$. Since~$\kappa$ is CI, this closed immersion
is a regular embedding. 
By taking the fiber product with the flat~$\kappa_0$-scheme~$\mathscr X_0$, we get a regular embedding of~$\mathscr X$
into the regular scheme~$\A^m_{\mathscr X_0}$. 

\end{itemize}

\begin{theo}\label{thm-localring-generic}
Let~$Y\to X$ be a morphism between good~$k$-analytic spaces. Let~$y$ be a point of $Y$ and let~$x$ be its image
in~$X$.
Assume that $y$ belongs to $\mathrm{Int}(Y/X)$
and that~${\mathscr O}_{X,x}$ is a field. 
The morphism~$\spec{\mathscr O}_{Y_{x},y}\to \spec{\mathscr O}_{Y,y}$ is then flat, and its fibers satisfy property~$\mathfrak R$ of~\ref{ss-property-r}, In particular, its
fibers are CI, and are
regular if \car$k=0$. 
\end{theo}

\begin{proof}
Set~$n=\dim_y Y_x$. We
argue by induction on~$n$; let us begin with some preparation.

According to \Th 4.6 of \cite{ducros2007}, one can shrink~$Y$ around~$y$ such that~$Y\to X$ goes through a map~$Y\to \A^{n}_{X}$ that is quasi-finite at~$y$.
By assumption, $y$ belongs to  $\mathrm{Int}(Y/X)$, so it belongs to $\mathrm{Int}(Y\to \A^{n}_{X})$; hence $Y\to X$ it is finite at~$y$
by \cite{berkovich1993}, \Prop 3.1.4. Denote by~$t$ the image of~$y$ in~$\A^{n}_{X}$. One can shrink~$Y$ around~$y$ so that it is finite
over an affinoid neighborhood~$V$ of~$t$ in~$\A^{n}_{X}$ (note that $V\to X$ is then boundaryless at $t$) and so that~$y$ is the only preimage of~$t$ in~$Y$.  Let~$A$ (\resp $B$) be the algebra of analytic functions on~$V$ (\resp $Y$). Then~${\mathscr O}_{Y,y}=B\otimes_A{\mathscr O}_{V,t}$ and~${\mathscr O}_{Y_{x},y}=B\otimes_A{\mathscr O}_{V_{x},t}$; hence~$ {\mathscr O}_{Y_{x},y}={\mathscr O}_{V_{x},t}\otimes_{{\mathscr O}_{V,t}}{\mathscr O}_{Y,y}.$
It is thus
sufficient
(by finiteness of $\mathscr O_{Y,y}$ over $\mathscr O_{V,t}$)
to prove that~$\spec{\mathscr O}_{V_{x},t}\to \spec{\mathscr O}_{V,t}$ is flat and that its fibers satisfy~$\mathfrak R$. Let us list some facts that will be useful for the proof. 

\begin{enumerate}[A] 

\item As~${\mathscr O}_{X,x}$ is a field, it is regular; hence~${\mathscr O}_{V,t}$ is regular by Proposition \ref{prop-tranferqsm-concrete}.  

\item The local ring~${\mathscr O}_{V_{x},t}$ is regular by \ref{ss-omega-inequality}. 

\item The ring~${\mathscr O}_{V,t}$ being regular by (A), it is in particular reduced; this implies that if~$f$ is a non-zero element of it, its zero-locus (which is a Zariski-closed subset of a suitable neighborhood of~$t$) contains no neighborhood of~$t$ in~$V$. By Proposition~\ref{prop-nghbr-generic}, it follows that this zero-locus contains no neighborhood of~$t$ in~$V_{x}$. 

\end{enumerate}

Let us now go back to the proof by induction on~$n$. If~$n=0$ then~$t=x$, and both~${\mathscr O}_{V,t}={\mathscr O}_{X,x}$ and~${\mathscr O}_{V_{x},t}=\hr x$ are fields, hence we are done. 
Assume that~$n>0$ and that the theorem has been proved for any integer $<n$. 

Let us first prove that $\mathscr O_{V_x,t}$ is a flat $\mathscr O_{V,t}$-algebra.
By \cite{sga1}, Exposé IV, Th. 5.6, it is sufficient to prove that~${\mathscr O}_{V_{x},t}/\mathfrak m_t^d{\mathscr O}_{V_{x},t}$ is a flat~${\mathscr O}_{V,t}/\mathfrak m_t^d$-algebra for any~$d>0$. If~$\mathfrak m_t=0$, then~${\mathscr O}_{V,t}$ is a field and we are done. 
Suppose now that $\mathfrak m_t\neq 0$ and let~$d$ be a positive integer. Due to assertion (C) above, the closed analytic
subspace~$Z$ defined in a neighborhood of~$t$ by the finitely generated ideal~$\mathfrak m_t^d$ contains no neighborhood of~$t$ in~$V_{x}$ ; therefore~$Z\to X$ is of dimension $<n$ at~$t$ (and also boundaryless at $t$, because $Z\hookrightarrow V$ is a closed immersion). By induction,~${\mathscr O}_{Z_{x},t}$ is a flat~${\mathscr O}_{Z,t}$-algebra. But
${\mathscr O}_{Z_{x},t}$ (resp.~${\mathscr O}_{Z,t}$) is nothing but~${\mathscr O}_{V_{x},t}/\mathfrak m_t^d\mathscr O_{V_x,t}$  (resp.~${\mathscr O}_{V,t}/\mathfrak m_t^d$), whence the
desired claim is proved.

Now, let us prove that any fiber of~$\spec{\mathscr O}_{V_{x},t}\to \spec{\mathscr O}_{V,t}$ satisfies~$\mathfrak R$.
Let~$\mathfrak p$ be a prime ideal of~${\mathscr O}_{V,t}$. If~$\mathfrak p=0$, the fiber of~$\spec{\mathscr O}_{V_{x},t}$ over~$\mathfrak p$ is the spectrum of a localization of~${\mathscr O}_{V_{x},t}$ ; but the latter is regular by assertion (B) above, hence we are done. 
Suppose now that~$\mathfrak p \neq 0$. By assertion (C) above, the closed analytic subspace~$Z$ defined in a neighborhood of~$t$ by the finitely generated ideal~$\mathfrak p$ contains no neighborhood of~$t$ in~$V_{x}$ ; therefore~$Z\to X$ is of dimension $<n$ at $t$ (and also boundaryless at $t$, because $Z\hookrightarrow V$ is a closed immersion).
The fiber of~$\spec{\mathscr O}_{V_{x},t}$ over~$\mathfrak p$ is nothing but the generic fiber of
the map $\spec{\mathscr O}_{Z_{x},t}\to \spec{\mathscr O}_{Z,t}$. By the induction hypothesis, the latter satisfies~$\mathfrak R$, 
which ends the proof.
\end{proof}

We are now going to show that the assumptions of Theorem \ref{thm-localring-generic}
are
probably not far from being optimal; we still denote by $Y\to X$ a morphism of good $k$-analytic spaces, by $x$ a
point of $X$ whose local ring is a field, and by $y$ a point of $Y_x$. 

\subsection{}\label{ss-optimal-flat}
One cannot expect in general flatness of $\spec{\mathscr O}_{Y_{x},y}$ over $\spec{\mathscr O}_{Y,y}$ if~$y$ belongs to $\partial (Y/X)$. 
Indeed, let~$r$ be a positive real number let and let~$f$ be a power series in one variable
with coefficients in~$k$ whose radius
of convergence is equal to $r$. Let~$V$ be the analytic domain of~$\A^{2,\mathrm{an}}_k$ defined by the condition~$\abs{T_1}=r$.
There is a natural closed immersion~$\phi : D\to V$ given by~$(\mathrm{Id}, f)$, where~$D$ is the closed disc of radius~$r$; let~$x$ denote the image under~$\phi$ of the Shilov point
of~$D$. 
Proposition~\ref{prop-oanx-field} ensures that~$\mathscr O_{\A^{2,\rm an}_k,x}$ is a field.
The fiber of~$V\hookrightarrow \A^{2,\rm an}_k~$ at~$x$ is nothing but~$\mathscr M(\hr x)$, and~$\mathscr O_{V_x,x}$ is then simply the field~$\hr x$. As~$x$ lies on a one-dimensional Zariski-closed subset
(namely,~$\phi(D)$) of the purely 2-dimensional space~$V$,  the local ring~$\mathscr O_{V,x}$ {\em cannot} be a field (Corollary \ref{cor-interp-centdim}). As a consequence, ~$\spec \hr x\to \spec \mathscr O_{V,x}$ is not flat. 
\subsection{}\label{ss-optimal-reduced}
One cannot expect in general regularity of the fibers of the morphism~$\spec{\mathscr O}_{Y_{x},y}\to \spec{\mathscr O}_{Y,y}$ if~$k$ is of positive characteristic. Indeed, let us give the following counter-example, which was communicated to the author by M. Temkin. Assume that~$k$ is a non-algebraically closed field of char.~$p>0$ such that~$\abs {k\gpm}\neq   \{1\}$.  Let~$k^a$ be an algebraic closure of~$k$, and let~$k^s$ be the separable closure of~$k$ inside~$k^a$; assume moreover
that $k^s$ is of countable dimension over $k$ (for instance, we can take for $k$ the field $\F_p(\!(t)\!)$ equipped with an arbitrary $t$-adic absolute value). By assumption, there exists an increasing sequence ~$(k_n)_{n\in \N}$ of subfields of~$k^s$ that are finite over~$k$ and whose union is equal to~$k^s$. For any~$n$, the complement of a finite union of proper~$k$-vector subspaces of the~$k$-Banach space~$k_n$ is a dense subset of it. Therefore there exists a sequence~$(\lambda_n)$ of elements of~$k^s$ and a decreasing sequence~$(r_n)$ of positive real numbers
such that the following hold:

\begin{enumerate}[a]
\item for any~$n$, one has~$k[\lambda_n]=k_n$;
\item for any~$n$ and any non-trivial $k$-conjugate $\mu$ of~$\lambda_n$ in~$k^a$, one has
$|\mu-\lambda _n|>r_n$;
\item for any~$m>n$, one has~$|\lambda_m-\lambda_n|<r_n$;
\item one has  $r_n\to 0$ as~$n\to \infty$.
\end{enumerate} 

For any~$n$ let us denote by~$D_n$ the affinoid domain of~$\A^{1, \rm an}_{k_n}$ defined by the inequality~$\abs{T-\lambda_n}\leq r_n$. It follows from~(a) and~(b) that the natural map~$\A^{1, \mathrm{an}}_{k_n}\to \A^{1, \mathrm{an}}_k$ induces an isomorphism beween~$D_n$ and an affinoid domain~$\Delta_n$ of~$\A^{1, \mathrm{an}}_k$. It follows from~(c) and~(d) that~$(D_n(\widehat{k^a}))$ is a decreasing sequence of (naive) closed discs of $\widehat{k^a}$ whose intersection consists of a single element~$\lambda\in \widehat{k^a}$. Let~$x\in \A^{1, \mathrm{an}}_k$ be the point that corresponds to~$\lambda$. We have~$x\in \bigcap \Delta_n$. Therefore,~$k_n$ embeds into~$\hr x\subset \widehat{k^a}$ for every~$n$. Hence~$\hr x$ is a closed subfield of~$\widehat{k^a}$ containing~$k^s$; the latter being dense in~$\widehat{k^a}$, 
one has
$\hr x=\widehat{k^a}$; in particular,~$x$ is not a rigid point. 

Let~$\phi$ be the (finite, flat) morphism~$\A^{1, \mathrm {an}}_k\to \A^{1, \mathrm {an}}_k$ induced
by the morphism from~$k[T]$ to itself that sends~$T$ to~$T^p$, and let~$y$ be the unique preimage of~$x$ by~$\phi$. As~$x$ in non-rigid,~$y$ is non-rigid
and~$\mathscr O_{\A^{1, \mathrm {an}}_k,y}$ is a field due to Lemma \ref{lem-onedim-field}.
Now~$\mathscr O_{\phi\inv (x),y}=\hr x[\tau]/(\tau^p-T(x))$. Since ~$\hr x=\widehat{k^a}$, the local ring ~$\mathscr O_{\phi\inv (x),y}$ is non-reduced, and in particular, non-regular. Hence the map
$$\spec \mathscr O_{\phi\inv (x),y}\to \spec \mathscr O_{\A^{1,\mathrm{an}}_k,y}$$ is not regular.

\subsection{}
Let~$X$ be a reduced good $k$-analytic space, and let~$x$ be an Abhyankar 
point of $X$ (\ref{ss-abhyankar-points}).
The local ring
$\mathscr O_{X,x}$ is then artinian (Example \ref{ex-centdim-recap}), hence it a field
by reducedness. If $Y$ is any good
$X$-analytic space, we can thus
apply Theorem \ref{thm-localring-generic} to
the map $Y\to X$
over the point $x$. But in fact, due to the Abhyankar property of $x$, 
we have the following slightly stronger result, whose proof is totally
different.

\begin{theo}
\label{thm-loc-gen-abhy}

Let~$Y\to X$ be a morphism between good $k$-analytic spaces, with $X$ reduced. 
Let $x$ be an  Abhyankar point of~$X$ (\ref{ss-abhyankar-points}),
and let
$y$ be a point of the fiber $Y_x$. 
The morphism
$$\spec \mathscr O_{Y_x,y}\to \spec \mathscr O_{Y,y}$$ is regular
(note that we do not assume $y\in \mathrm{Int}(Y/X)$ nor \car$k=0$).
\end{theo}

\begin{proof} Set~$n=\dim_x X$. By shrinking~$X$ and~$Y$  we immediately reduce to the 
case where both are affinoid and where~$\dim X=n$ (we can proceed to such a reduction
without modifying $\mathscr O_{Y,y}$ because $n$ is the infimum of the dimensions of all
analytic {\em neighborhoods}
of $x$ in $X$, by \ref{rem-dimloc-open}); we denote by~$A$
(\resp $B$) the algebra of analytic functions on~$X$
(\resp $Y$).

Let us first assume that~$X$ is a closed $n$-dimensional polydisc centered at the origin of $\A^{n,\mathrm{an}}_k$
and $x=\eta_r$ for some polyradius $r=(r_1,\ldots, r_n)$.
The point~$x$ has a canonical~$\hr x$-rational pre-image
$t$ on the~$\hr x$-analytic disc~$X_{\hr x}$. Choose $\epsilon \in \R_+\gpm$ such that~$\epsilon <r_i$ for every~$i$. Let
$V$ be the 
affinoid domain of~$X_{\hr x}$ defined by the inequalities~$\abs{T_i-T_i(t)}\leq \epsilon$ for~$i=1,\ldots, n$ (where the $T_i$'s are the coordinate functions on $X$)
and let~$C$ be the algebra of analytic functions on~$V$.
For every point $v$ of $V$ and every $i$, we have $\widetilde{T_i(v)}=\widetilde{T_i(t)}$ (the latter equality makes sense because $\hr t=\hr x\subset \hr v$),
and $\abs{T_i(v)}=r_i$. Since $t$ lies above $x$, the elements $\widetilde{T_1(t)},\ldots, \widetilde{T_n(t)}$ 
of $\hrt t=\hrt x$ are algebraically independent over the graded field $\widetilde k$; therefore $\widetilde{T_1(v)},\ldots, \widetilde{T_n(v)}$ are algebraically independent over $\widetilde k$, which implies
(together with the fact that $\abs{T_i(v)}=r_i$ for every $i$) that the image of $v$
in $X$ is equal to $\eta_r=x$; hence $V$ is contained in the fiber of $X_{\hr x}$ over $x$. 
As a consequence, for every $a\in A$  the image of~$a$ under the composite map
$A\to A_{\hr x}\to C$ is equal to the element~$a(x)$ of~$\hr x$ (indeed, since~$C$ is reduced, it is sufficient to check this equality pointwise on~$V$). In other words, both maps
$$A \to A_{\hr x}\to C\;\;{\rm and}\;\xymatrix{
A\ar[rr]^{a\mapsto a(x)}&&{\hr x}\ar[r]&C
}$$ coincide.
Let~$W$ be the preimage of~$V$ in~$Y_{\hr x}$. This is an affinoid domain of~$Y_{\hr x}$,
and
$$\mathscr O_{Y_x}(W)=B\hotimes_A C=(B\hotimes_A\hr x)\hotimes_{\hr x}C.$$
Hence the morphism $W\to Y$ goes through~$Y_x$, and the~$Y_x$-analytic space~$W$ 
is isomorphic to~$Y_x\times_{\hr x} V$. 
As $\hr t=\hr x$, the natural map from $Y_{\hr x, t}$
to $Y_x$ is an isomorphism; let $y'$ be the unique point of~$Y_{\hr x, t}$
lying above~$y$. Since $V$ is a neighborhood of $t$ in $X_{\hr x}$, the affinoid
domain $W$ is a neighborhood of $y'$
in $Y_{\hr x}$, and we thus have $\mathscr O_{W,y'}=\mathscr O_{Y_{\hr x}, y'}$. As $\hr x$ is analytically
separable over~$k$, the morphism from~$\spec \mathscr O_{W,y'}=\spec \mathscr O_{Y_{\hr x}, y'}$ to~$\spec \mathscr O_{Y,y}$ is regular
(Theorem \ref{th-fiber-xl}). 
And by flatness of $V$ over~$\hr x$ (Lemma \ref{lem-field-flat}) and in view of the $Y_x$-isomorphism
$W\simeq Y_x\times_{\hr x}V$, the map
$\spec \mathscr O_{W,y'}\to \spec \mathscr O_{Y_x,y}$
is (faithfully) flat (in fact, by quasi-smoothness of~$V$
and in view of Theorem
\ref{thm-main-qsm}, it is even regular
but we shall not need that). 
Lemma~\ref{lem-reg-sorital}
then applies to the commutative diagram
$$\xymatrix{
{\spec \mathscr O_{W,y'}}\ar[rd]\ar[d]&\\
{\spec \mathscr O_{Y_x,y}}\ar[r]&
{\spec \mathscr O_{Y,y}}
}$$ and yields the regularity
of the map~$\spec \mathscr O_{Y_x,y}\to \spec \mathscr O_{Y,y}$, which ends the proof
in the special case considered.

Let us now prove the theorem for arbitrary $X$. Since~$d_k(x)=n$, there exist analytic functions~$f_1,\ldots, f_n$ on~$X$, invertible at~$x$
and such that the~$\widetilde{f_i(x)}$'s are algebraically independent over~$\widetilde k$. Let~$R$ be the supremum of the
spectral norms
of the~$f_i$'s. Those functions induce a morphism from~$X$
to
the polydisc~$Z:=\mathscr M(k\{T_1/R,\ldots, T_n/R\})$; let~$z$ denote the image of~$x$ in~$Z$. 
If we set~$r_i=|f_i(x)|$ 
for every~$i$, then our assumption on the~$f_i$'s means that $z=\eta_r$ for $r=(r_1,\ldots, r_n)$.
Let~$t$ be a point of $X_z$. By assumption, $\dim X=n$ and we know that $d_k(\hr z)=n$ (\cf Example \ref{exem-gradres-etar}).
As a consequence, $n\geq d_k(t)=d_{\hr z}(t)+n$, whence the equality $d_{\hr z}(t)=0$. It follows that the fiber~$X_z$
is zero-dimensional. It thus consists of finitely many points~$x_1=x,x_2,\ldots, x_m$. 
By assumption, the space~$X$ is reduced. By the particular case proven above, 
$\spec \mathscr O_{X_z,x_j}\to \spec  \mathscr O_{X,x_j}$ is regular for every~$j$, which implies that the~$\mathscr O_{X_z,x_j}$'s
are reduced (this follows for instance from \cite{ega42}, \Prop 6.4.1 (ii) and \Prop 6.5.3 (ii)). 
Therefore the fiber~$X_z$ is isomorphic to $\coprod \mathscr M(\hr {x_j})$; it follows that
$$Y_z\simeq \coprod Y_{x_j}.$$ In particular, $Y_x$ is an open subspace of~$Y_z$, and~$\mathscr O_{Y_x,y}$ is thus
equal to~$\mathscr O_{Y_z,y}$. Using again the particular case proven above (now applied to $Y\to Z$) we
get the regularity of the map~$\spec \mathscr O_{Y_z,y}\to \spec \mathscr O_{Y,y}$; hence
$\spec \mathscr O_{Y_x,y}\to \spec \mathscr O_{Y,y}$ is regular. 
\end{proof}

\chapter{Images of morphisms: local results}\label{GR}

The general purpose of this chapter is the study of the ``image"
of a morphism between analytic germs. Most of the work is in fact
carried out in the realm of (graded)
Riemann-Zariski spaces, and is transferred thereafter to the analytic
world through Temkin's theory. 

Our main result about Riemann-Zariski spaces is kind of
a ``Chevalley theorem" in this setting
(Theorem \ref{thm-elim-graded}). It tells the following
(we use the definitions and notation introduced in \ref{ss-def-plkzr}). 
Let $K$ be a graded field, let $F$ be a graded extension of $K$, and
let $L$ be a graded extension of $F$. Let $\mathsf U$
be a quasi-compact open subset of
$\P_{L/K}$. The image of $\mathsf U$
in $\P_{F/K}$ is a quasi-compact open subset of $\P_{F/K}$. 

The proof goes as follows. Using scalar extension to $K(T/r)$ for suitable
polyradius $r$ (this notation 
is introduced in \ref{sss-notation-ktr}), one reduces to the non-graded case; 
one can
clearly
assume that $\mathsf U=\P_{L/K}\{f_1,\ldots, f_n\}$ 
for some elements $f_1,\ldots, f_n$ of $L$. One first considers the case where $L$
is algebraic over $F$,
for which one
proves the required assertion by
performing explicit computations based on Newton polygons (Proposition \ref{pro-elim-simple}). 
For handling the general case, one sets $A=F[f_1,\ldots, f_n]\subset L$, and shows
(Theorem \ref{theo-elim-sections}) that there exist finitely many closed point
$\mathsf y_1,\ldots, \mathsf y_m$ of $\spec A$ such that the image
of $\mathsf U$ in $\P_{F/K}$ is the union of the images of
the sets $\P_{\kappa(\mathsf y_j)/K}\{f_1(\mathsf y_j),\ldots, f_n(\mathsf y_j)\}$
in $\P_{F/K}$
for $j$ varying between $1$ and $m$; it is thus a quasi-compact
open subset of $\P_{F/K}$ by the algebraic case already investigated. 
Beside some
quasi-compactness arguments, the proof of Theorem \ref{theo-elim-sections} rests
on
the so-called \emph{quantifier elimination in the theory of non-trivially
valued algebraically closed fields}, 
which can be seen as a valuative
avatar
of Chevalley's theorem. For the reader's convenience, 
Section \ref{appendix-acvf} at the end of the chapter
explains what
it consists of, and also provides 
a way to replace it (for our purposes) by a scheme-theoretic 
argument. 

Applications of the above to analytic
geometry are carried out in Section \ref{s-smallest-image}. 
The main result is the following: if $(Y,y)\to (X,x)$ is a morphism
of analytic germs, then there exists a smallest analytic domain
$(Z,x)$ of $(X,x)$ through which
the map $(Y,y)\to (X,x)$ factorizes
(Theorem \ref{thm-smallest-germ} (1)). Moreover, if $\Gamma$
is a subgroup of $\R_+\gpm$ such that $\Gamma \cdot \abs{k\gpm}\neq \{1\}$,
then $(Z,x)$ is $\Gamma$-strict if $(Y,y)$ is $\Gamma$-strict and $(X,x)$
is separated (\loccit, (3)). The latter fact has a global consequence: if $\phi\colon Y\to X$
is a morphism between affinoid spaces and $Y$ is $\Gamma$-strict, 
then $\phi(Y)$ is contained in a $\Gamma$-strict analytic domain of $X$
(Proposition \ref{prop-image-subgstrict}). 

\section{Maps between Riemann-Zariski Spaces: the non-graded case}
\markright{\thesection.~MAPS BETWEEN R-Z SPACES: THE NON-GRADED CASE}

The purpose of what follows is to prove an avatar of Chevalley's constructibility theorem
in the setting of (non-graded) Riemann-Zariski spaces.

\begin{lemm}\label{lem-extend-gauss}
Let~$K$ be a field
and let~$\abs {\cdot}$ be a valuation on it. Let~$L$ be an algebraic extension of~$K$ ;  let~$\abs {\cdot}'$ be a valuation on~$L$ extending~$\abs {\cdot}$.
Let~$S_1,\ldots,S_n$ be indeterminates, and let~$\abs {\cdot}''$ be an extension of~$\abs {\cdot}'$ to~$L(S_1,\ldots,S_n)$ whose restriction to~$K(S_1,\ldots,S_n)$ is equal to~$\abs {\cdot}_{\mathrm{Gau\ss}}$. Then ~$\abs.''=\abs.'_{\mathrm{Gau\ss}}$.

\end{lemm}

\begin{proof} We denote by~$k$, \resp~$\ell$,  the residue field of~$\abs {\cdot}$, \resp~$\abs {\cdot}'$. For every $i$ one has $\abs{S_i}_{\mathrm{Gau\ss}}=1$, and the images of the~$S_i$'s in the residue field of~$\abs {\cdot}_{\mathrm{Gau\ss}}$ are algebraically independent over~$k$. As~$L$ is algebraic over~$F$, the field~$\ell$ is algebraic over~$k$. Therefore, the images of the~$S_i$'s in the residue field of~$\abs {\cdot}''$ are algebraically independent over~$\ell$, whence the required equality $\abs {\cdot}''=\abs.'_{\mathrm{Gau\ss}}$.
\end{proof}

\begin{lemm}\label{lem-size-root}
Let~$F$ be a field and let~$P=T^{n}+ a_{n-1}T^{n-1}+\ldots+a_{0}$ be a monic polynomial belonging to~$F[T]$; set~$a_n=1$. Assume that~$P$ is totally split.
Let~$\abs {\cdot}$ be a valuation on~$F$. The following are equivalent: 

\begin{enumerate}[i]
\item $\abs \lambda>1$ for every root~$\lambda$ of~$P$ ;

\item $\abs {a_0}>\abs{a_j}$ for every~$j>0$.
\end{enumerate}
\end{lemm}

\begin{proof}
If (i) is true then (ii) follows immediately from the usual relations between the coefficients and the roots of~$P$. Suppose that (ii)
is true, and let~$\lambda$ be a root of~$P$. As~$P(\lambda)=0$ there exists~$j>0$ such that~$\abs{a_0}\leq \abs{a_j \lambda^{j}}$. Since~$\abs{a_0}>\abs{a_j}$, this implies that~$\abs \lambda>1$.
\end{proof}

\begin{prop}\label{pro-elim-simple}
Let~$k$ be a field, let~$K$ be an extension of~$k$, and let~$L$ be an algebraic extension of~$K$. Let~$\mathsf U$ be a quasi-compact open subset of ~$\P_{L/k}$. Its image~$\mathsf V$ on~$\P_{K/k}$ is a quasi-compact open subset of the latter.
\end{prop}

\begin{proof}
We can assume that~$\mathsf U$ is equal to~$\P_{L/k}\{f_1,\ldots,f_\ell\}$ for suitable elements~$f_1,\ldots,f_\ell$ of~$L$. Let $S=(S_1,\ldots, S_\ell)$ be a family of indeterminates, and let $f$ be the element $f_1S_1+\ldots f_\ell S_\ell$ of $L(S)$. Let~$P=T^{n}+a_{n-1}T^{n-1}+\ldots+a_0$ be the minimal polynomial of~$f$ over~$K(S)$. Let~$\Lambda$ be a finite extension of~$K(S)$ over which~$P$ splits. Let~$\abs \cdot$ be a valuation on~$K$ whose restriction to~$k$ is trivial. We fix an extension~$\abs \cdot_0$ of~$\abs\cdot_{\mathrm{Gau\ss}}$ to~$\Lambda$.

We are going to prove that the valuation~$\abs \cdot$ belongs to ~$\mathsf V$ if and only if~$\abs\cdot_{\mathrm{Gau\ss}}$
extends to a valuation~$\abs \cdot''$ on~$L(S)$  such that~$\abs  f''\leq 1$. Let us first assume
that~$\abs \cdot  \in \mathsf V$. This means that it extends to a valuation~$\abs \cdot '$ on~$L$ such that~$\abs{f_i}'\leq 1$ for every ~$i$, so~$\abs \cdot '_{\mathrm{Gau\ss}}$ is then an extension of~$\abs\cdot_{\mathrm{Gau\ss}}$ to~$L(S)$ satisfying the required properties. 
Conversely, assume that~$\abs \cdot_{\mathrm{Gau\ss}}$ extends to a valuation~$\abs \cdot ''$ on~$L(S)$
such that~$\abs f''\leq 1$ and let~$\abs \cdot '$ denotes the restriction of~$\abs \cdot ''$ to~$L$. By Lemma~\ref{lem-extend-gauss}, one has~$\abs \cdot''=\abs \cdot '_{\mathrm{Gau\ss}}$; therefore, the inequality~$\abs f''\leq 1$ simply means that~$\abs {f_i}'\leq 1$ for every ~$i$, and we are done.

On the other hand, the valuation~$\abs\cdot_{\mathrm{Gau\ss}}$ admits an extension~$\abs \cdot ''$ to~$L(S)$ such that ~$\abs f''\leq 1$ if and only if there exists a root~$\lambda$ of~$P$ in~$\Lambda$ such that~$\abs \lambda _0 \leq 1$; according to Lemma \ref{lem-size-root}, the latter condition is equivalent to the existence of $j>0$ such that~$\abs{a_0}_{\mathrm{Gau\ss}}\leq \abs{a_j}
_{\mathrm{Gau\ss}}$ (note that if there exists such $j$, it can always be chosen so that $a_j\neq 0$:
this is obvious if $a_0\neq 0$, and if $a_0=0$ we have $n=1$ and we take $j=1$). 

We thus have proved that $\mathsf V$ is equal to the preimage under~$\abs \cdot \mapsto \abs \cdot_{\mathrm{Gau\ss}}$ of ~$\bigcup\limits_{a_j\neq 0}  \P_{K(S)/k}\{a_0/a_j\}$; and this preimage is easily seen, by the very definition of~$\abs \cdot_{\mathrm{Gau\ss}}.$ for a given~$\abs \cdot$,
to be a quasi-compact open subset of~$\P_{K/k}$.
\end{proof}

\begin{theo}\label{theo-elim-sections}

Let~$K$ be a field, let~$F$ be an extension of~$K$, and let~$L$ be an extension of~$F$. Let~$f_{1},\ldots,f_{n}$ be finitely many elements of~$L$ and let~$A$ be the~$F$-subalgebra of~$L$ generated by the~$f_{i}$'s. For any~$y\in \spec A$, let $p_y$ denote the map~$\P_{\kappa(y)/K}\to\P_{F/K}$.

\begin{enumerate}[1]
\item The image~$\mathsf V$ of~$\P_{L/K}\{f_1,\ldots,f_n\}$ in~$\P_{F/K}$ is a quasi-compact open subset of the latter.

\item There exist finitely many {\em closed} points~$y_1,\ldots,y_m$ of~$\spec A$ such that ~$$\mathsf V=\bigcup_j p_{y_j}(\P_{\kappa(y_j)/K}\{f_1(y_j),\ldots,f_n(y_j)\}\;).$$
\end{enumerate}

\end{theo}

\begin{proof}
If~$y$ is any closed point of~$\spec A$, its residue field~$\kappa(y)$ is finite over~$F$. Proposition~\ref{pro-elim-simple}
above then ensures that~$p_{y}(\;\P_{\kappa(y)/K}\{f_1(y),\ldots,f_n(y)\}\;)$ is a quasi-compact open subset of~$\P_{F/K}$. As~$\mathsf V$ is a quasi-compact topological space, it is therefore enough, to establish both~(1) and (2),  to prove that~$$   \mathsf V=\bigcup_{
y\in \mathsf C} p_{y}(\P_{\kappa(y)/K}\{f_1(y),\ldots,f_n(y)\}\;)$$ where~$\mathsf C$ is the set of {\em all} closed points of~$\spec A$. 

Let us first prove that~$$\bigcup_{
y\in \mathsf C} p_{y}(\;\P_{\kappa(y)/K}\{f_1(y),\ldots,f_n(y)\}\;)\subset  \mathsf V.$$ We shall
even show that~$$\bigcup_{
y\in \tiny \spec A} p_{y}(\;\P_{\kappa(y)/K}\{f_1(y),\ldots,f_n(y)\}\;)\subset  \mathsf V.$$ Let~$y$ be any point of~$\spec A$ and let~$\abs \cdot$ be a valuation on~$F$ that is trivial on~$K$ and that belongs to~$p_{y}(\;\P_{\kappa(y)/K}\{f_1(y),\ldots,f_n(y)\}\;)$; \ie, $\abs \cdot$ extends to a valuation $\abs \cdot '$ on~$\kappa(y)$ which satisfies the inequality~$\abs{f_j(y)}'\leq 1$ for every~$j$. Let~$\abs \cdot ''$ be a valuation on~$L$ whose ring dominates~$\mathscr O_{\tiny \spec A,y}$. The residue field~$\Lambda$ of~$\abs {\cdot}''$ is an extension of~$\kappa(y)$; we choose an extension~$\abs {\cdot}'''$ of~$\abs {\cdot}'$ to~$\Lambda$. The composition of~$\abs {\cdot}''$ and~$\abs \cdot'''$ is a valuation on~$L$ whose restriction to~$F$ is equal to~$\abs {\cdot}$ and whose ring contains the~$f_j$'s. Hence~$\abs {\cdot}\in \mathsf V$.

Let us now prove that~$$\mathsf V\subset \bigcup_{
y\in \mathsf C} p_{y}(\;\P_{\kappa(y)/K}\{f_1(y),\ldots,f_n(y)\}\;).$$ Let~$\abs {\cdot}$ be a valuation on~$F$ that is trivial on~$K$ and that  belongs to~$\mathsf V$; \ie, it extends to a valuation~$\abs {\cdot}'$ on~$L$ that satisfies the inequalities~$\abs{f_i}'\leq 1$ for all~$i$.

Let~$L^{\mathrm a}$  be an algebraic closure of~$L$ and let~$F^{\mathrm a}$ be the algebraic closure of~$F$ inside~$L^{\mathrm a}$. We choose an extension~$\abs {\cdot}''$ of~$\abs {\cdot}'$ to~$L^{\mathrm a}$. Let~$(P_1,\ldots,P_m)$ be polynomials that generate the ideal of relations between the~$f_i$'s over the field~$F$, 
so $A=F[T_1,\ldots, T_n]/(P_1,\ldots, P_m)$. 
The system of equations and inequalities (in variables~$x_1,\ldots,x_n$)~$$\{P_j(x_1,\ldots,x_n)=0\}_{j=1,\ldots,m}\;{\rm and}\;\{\abs{x_i}''\leq 1\}_{i=1,\ldots,n}$$ has a solution in~$L^{\mathrm a}$, provided by the~$f_i$'s. This implies that it has a solution $(g_1,\ldots, g_n)$ in $F^{\mathrm a}$. Indeed, if $\abs \cdot$ is trivial
this comes from the {\em Nullstellensatz}, because the inequality $\abs x\leq 1$ is then satisfied by every $x\in F$. And if $\abs .$ is non-trivial, this is a particular case
of the so-called {\em model completeness} of the theory of algebraically closed, non-trivially valued fields, which itself follows from {\em quantifier elimination}; but this can also be 
given a direct, algebro-geometric proof. Reminders on model-completeness and quantifier elimination
have been postponed to
Section \ref{appendix-acvf} (see Theorems \ref{thm-qe-acvf}
and \ref{thm-mc-acvf}), as well as the direct proof alluded to (Proposition~\ref{prop-mc-acvf}
and Theorem~\ref{thm-MC-acvf}). 

Now evaluation at $(g_1,\ldots, g_n)$ defines a map from $A=F[T_1,\ldots, T_n]/(P_1,\ldots, P_m)$ to $F^{\mathrm a}$, which sends $f_i$ to $g_i$ for every $i$. Its kernel is a closed point $y$ of $\spec A$, and we have an $F$-embedding $\iota$
of $\kappa(y)$ into $F^{\mathrm a}$ mapping $f_i(y)$
to $g_i$ for every $i$. The composition $\abs \cdot ''\circ \iota$ is then a valuation on $\kappa(y)$
that extends $\abs .$ and whose ring contains the $f_i(y)$'s, and we are done. 
\end{proof}

\section{Maps between Riemann-Zariski spaces: the general case}\markright{\thesection.~MAPS BETWEEN R-Z SPACES: THE GENERAL
CASE}
We are now going to give a graded version of Theorem~\ref{theo-elim-sections} (1). We refer the reader to \ref{s-grad-alg} 
for our general conventions in graded commutative algebra, to \ref{s-graded-linear} for graded linear algebra
(and especially to Definition \ref{def-grad-tens} \ff
for graded tensor products
and flatness in the
graded context), and to \ref{s-grad-val}
for graded valuations. 

\begin{lemm}\label{lem-extension-bivaluation}
Let~$K$ be a graded field and let $E$ and $L$ be two graded field extensions of $K$.
Let~$A$, \resp $B$, \resp $C$ be a graded valuation ring of~$K$, 
\resp $E$, \resp $L$; assume that~$$B\cap K=C\cap K=A.$$ 

\begin{enumerate}[1]
\item Let $F$ be a graded field equipped with an injective morphism
$E\otimes_K L\hookrightarrow F$, making $F$ a graded extension of both $E$ and $L$. There exists 
a graded valuation ring~$D$ of~$F$ such that~$D\cap E=B$ and~$D\cap L=C$.

\item Assume that $E$ and $L$ are algebraic over $K$ (\ref{ss-graded-algebraic}). There exists
a common graded extension $\Lambda$ of $E$ and $L$ over $K$
and a
graded valuation ring $\Delta$ of $\Lambda$
such that $\Delta\cap E=B$ and~$\Delta\cap L=C$.

\end{enumerate}

\end{lemm}

\begin{proof}
We denote by~$a, b$ and~$c$ the respective residue graded fields of~$A,B$ and~$C$.
We choose a maximal homogeneous ideal of the non-zero graded ring~$b\otimes_ac$
and denote by $d$ be the corresponding quotient.
Note that since $B$ and $C$ have no
non-zero homogeneous $A$-torsion, both are flat as graded $A$-modules:
this follows from the flatness criterion given at the end of  \ref{ss-grad-flat} and from the fact that every ideal $I$ of $A$ generated by a finite set $S$ of homogeneous elements is principal (if $S=\emptyset$ one has $I=0$, and if $S\neq \emptyset$ one has $I=(s)$ for any element $s$ of $S$ of maximal valuation).

Let us show (1). By $A$-flatness of $C$, 
the natural map $$u\colon B\otimes_A C \to E\otimes_A C=E\otimes_K (K\otimes_A C)$$
is injective. As~$K\otimes_A C$ is simply a graded localization of~$C$ by a homogeneous multiplicative subset which does not contain zero,~it embeds into $L$; the natural map
$v\colon E\otimes_A C\to E\otimes_K L\hookrightarrow F$ is thus injective. As a consequence,
the natural map $v\circ u \colon B\otimes_A C\to F$ induces an isomorphism 
between $B\otimes_A C$ and the graded subring $B\cdot C$ of $F$ generated by $B$ and $C$. Hence there exists a (unique) map~$B\cdot C \to d$ extending both~$B\to b\to d$ and~$C \to c\to d$. The kernel of this map is a homogeneous prime ideal (because its target is a graded field); by Zorn's Lemma the corresponding graded localization of~$B\cdot C$ is dominated by a graded valuation ring~$D$ of~$F$. By construction,~$D\cap E=B$ and~$D\cap L=C$, whence (1). 

Let us now prove (2).
Let $\mathfrak p$
be the kernel of the morphism $B\otimes_AC\to b\otimes_ac\to d$ and let $R$ be
the graded localization $(B\otimes_AC)_{\mathfrak p}$. 
Since $B$ and $C$ are flat over $A$, the graded tensor product $B\otimes_AC$ is flat over $A$, hence has no non-zero 
homogeneous $A$-torsion.
The graded $A$-algebra $R$
is then also torsion-free, and $R\otimes_A K$ is consequently
non-zero. In particular, it has a homogeneous prime ideal, which corresponds (as in classical commutative algebra) to a homogeneous prime ideal $\mathfrak q$ of $B\otimes_A C$ such that $\mathfrak q\subset \mathfrak p$ and $\mathfrak q\cap A=0$. Let $\Lambda$ denote the graded fraction field of $(B\otimes_AC)/\mathfrak q$. 

Let $r$ be a positive real number and let $x$ be an element of $E^r$.
By assumption, $x$ is algebraic over $K$; therefore there exists a unitary homogeneous element $P$ of $K[T/r]$ such that $P(x)=0$. Choose a non-zero homogeneous element $a$ of $A$ such that 
$ab\in A$ for every coefficient $b$ of $P$. The element $ax$ of $E$ is then easily seen to be the root of a unitary homogeneous element
of $K[T/(r\cdot\deg a)] $ {\em with coefficients in $A$}. By a straightforward valuation computation, this implies that $ax\in B$. As a consequence, the localization $B\otimes_A K$ of $B$ is equal to $E$, and $C\otimes_AK=L$ analogously, whence the equality
$$(B\otimes_AC)\otimes_A K=E\otimes_KL.$$

The latter implies, together with the fact that $\mathfrak q\cap A=0$, 
that $B\otimes_AC\to \Lambda$ admits a factorization 
$$B\otimes_AC\to E\otimes_KL\to \Lambda.$$ Hence $\Lambda$ can be seen as common
graded
extension of $E$ and $L$ over $K$
so that the image of $B\otimes_A C\to \Lambda$ 
is the graded subring $B\cdot C$ of $\Lambda$ generated by $B$ and $C$. Since $\mathfrak q\subset \mathfrak p$, 
the map $B\otimes_AC\to b\otimes_ac\to d$ induces a morphism $B\cdot C\to d$, 
extending both~$B\to b\to d$ and~$C \to c\to d$. The kernel of this latter map is the
homogeneous prime ideal $\mathfrak p/\mathfrak q$; by Zorn's Lemma the corresponding graded localization of~$B\cdot C$
is dominated by a graded valuation ring~$\Delta$ of~$\Lambda$. By construction,~$\Delta\cap E=B$ and~$\Delta\cap L=C$, whence (2). 
\end{proof}

\begin{coro}\label{coro-extension-bivaluation}
Let
\[\xymatrix{
E\ar[r]&F\\
K\ar[u]\ar[r]&L\ar[u]
}
\]
be a commutative diagrams of graded fields
such that $E\otimes_KL\to F$ is injective. Let~$\ell$ be a graded subfield of~$L$; set~$k=K\cap \ell$. Define~$\pi,\rho,\phi$ and~$\psi$ by the commutative diagram
$$\xymatrix{
{\P_{F/\ell}}\ar[r]^\pi\ar[d]_\phi&{\P_{E/k}}\ar[d] ^\psi\\
{\P_{L/\ell}}\ar[r]^\rho&{\P_{K/k}}
}$$ and let $\mathsf U$ be any subset of~$\P_{E/k}$. One has $\phi(\pi\inv (\mathsf U))=\rho\inv (\psi(\mathsf U))$.
In particular, $\P_{F/\ell}\to \P_{E/k}\times_{\P_{K/k}}\P_{L/\ell}$ is surjective. 
\end{coro}

\begin{proof}
The inclusion~$\phi(\pi\inv (\mathsf U))\subset \rho\inv (\psi(\mathsf U))$ follows formally from the commutativity of the diagram. Now, let~$\abs {\cdot}$ be a valuation belonging to $\rho\inv (\psi(\mathsf U))$. This means that~$\abs {\cdot}$ is a graded valuation on~$L$, trivial over~$\ell$, and that there exists a graded valuation~$\abs {\cdot}'\in \mathsf U$ such that~$\abs {\cdot}'_{|K}=\abs {\cdot} _{|K}$. By Lemma \ref{lem-extension-bivaluation} (1), there exists a graded valuation~$\abs {\cdot}''$ on~$F$ whose restriction to~$E$ is equal to~$\abs {\cdot}'$, and whose restriction to~$L$ is equal to~$\abs {\cdot}$. The latter fact implies that the restriction of~$\abs {\cdot}''$ to~$\ell$ is trivial; therefore~$\abs {\cdot}''\in \pi\inv (\mathsf U)$ and~$\abs {\cdot}\in \phi(\pi\inv (\mathsf U))$; as a consequence $\phi(\pi\inv (\mathsf U))=\rho\inv (\psi(\mathsf U))$, as required. 
The last assertion follows by applying this equality for $\mathsf U$ a singleton. 
\end{proof}

Corollary \ref{coro-extension-bivaluation}
above can be used every time one has a commutative diagram of graded fields~$$\xymatrix{
 E\ar[r]& F\\
K\ar[u]\ar[r]&L\ar[u]}$$ such that~$E\otimes_K L\to F$ is injective. Let us give two examples of such a diagram, which will play a role in the sequel.

\begin{exem}\label{ex-tensor1}
Let~$K$ be a graded field, let~$L$ be a graded extension of~$K$ and let $\Gamma$ be  a subgroup of $\R_+\gpm$. {\em The natural map~$L^{\Gamma}\otimes_{K^\Gamma}K\to L$ is injective}. Indeed, for every class $c$ of $\R_+\gpm$ modulo $\Gamma$, set $K^c=\bigoplus_{r\in c}K^r$, and define $L^c$ analogously. Let $\mathscr C$ be the susbet of $\R_+\gpm$ consisting of classes $c$ such that $K^c\neq 0$.
One has~$K=\bigoplus\limits_{c\in \mathscr C} K^c$. For every~$c\in \mathscr C$, the summand $K^c$ is a one-dimensional graded vector space over $K^\Gamma$, and $L^c$ is a one-dimensional graded vector space over  $L^\Gamma$. Therefore~$L^\Gamma\otimes_{K^\Gamma}K^c\simeq L^c$
for every $c\in \mathscr C$ 
and $$L^\Gamma\otimes_{K^\Gamma}K\simeq \bigoplus_{c\in \mathscr C} L^c\subset L,$$ whence the claim. 
\end{exem}

\begin{exem}\label{ex-tensor2}
Let~$K$ be a graded field, let~$s=(s_1,\ldots,s_n)$ be a polyradius, and let $S=(S_1,\ldots, S_n)$ be a finite family of indeterminates. Let~$L$ be a graded extension of~$K$.
{\em The natural map~$L\otimes_K K(S/s)\to L(S/s)$ is injective}. Indeed, it follows directly from the definition that~$$L\otimes_K K[S/s]\simeq  L[S/s].$$ Therefore~$L\otimes_K K(S/s)$ appears as a graded localization of the graded domain ~$L[S/s]$ by a homogeneous multiplicative system which does not contain zero; hence it embeds in the graded fraction field~$L(S/s)$ of~$L[S/s]$.
\end{exem}

\begin{theo}\label{thm-elim-graded}
Let~$K$ be a graded field, let~$F$ be a graded extension of~$K$ and let~$L$ be a graded extension of~$F$. Let~$\Gamma$ be a subgroup of~$\R\gpm$, and let~$\mathsf V$ be a~$\Gamma$-strict quasi-compact open subset of~$\P_{L/K}$. Its image in~$\P_{F/K}$ is a~$\Gamma$-strict quasi-compact open subset of the latter.
\end{theo}

\begin{proof}
Let us first assume that $\Gamma=\{1\}$. Consider the commutative diagram
$$\xymatrix{
{\P_{L/K}}\ar[r]^\pi\ar[d]_\phi&{\P_{L^1/K^1}}\ar[d] ^\psi\\
{\P_{F/K}}\ar[r]^\rho&{\P_{F^1/K^1}}
}$$ in which $\rho, \pi, \phi$ and $\psi$ are the obvious maps. By assumption, there exists a quasi-compact open subset~$\mathsf U$ of~$\P_{L^1/K^1}$ such that~$\mathsf V=\pi\inv (\mathsf U)$. By Corollary~\ref{coro-extension-bivaluation} and Example \ref{ex-tensor1}, one has~$\phi(\mathsf V)=\rho\inv (\psi(\mathsf U))$. By Theorem~\ref{theo-elim-sections} (1), the image $\psi(\mathsf U)$ is a quasi-compact open subset of~$\P_{F^1/K^1}$, whence the result.

Let us treat now the general case. We may
and do assume that $\mathsf V$
is equal to $\P_{L/K}\{f_1,\ldots, f_n\}$ for some non-zero homogeneous elements
$f_i$ of $L^\Gamma$. For every $i$, we denote by $s_i$ the degree of $f_i$; set $s=(s_1,\ldots, s_n)$.  Let us consider the commutative diagram
$$\xymatrix{
{\P_{L(s\inv S)/K(s\inv S)}}\ar[rr]^\mu\ar[d]_\theta&&{\P_{L/K}}\ar[d] ^\phi\\
{\P_{F(s\inv S)/K(s\inv S)}}\ar[rr]^\nu&&{\P_{F/K}}
\ar@/^/ [ll]^\sigma
}$$ in which $\mu,\theta,\nu$ and $\phi$ are the obvious maps
and in which $\sigma$ is the map $\abs \cdot \mapsto \abs \cdot_{\mathrm{Gau\ss}}$.
The quasi-compact open subset~$\mu\inv (\mathsf V)$ of~$\P_{L(s\inv S)/K(s\inv S)}$ is equal to
$$\P_{L(s\inv S)/K(s\inv S)}\{f_1,\ldots, f_n\}=\P_{L(s\inv S)/K(s\inv S)}\{f_1/S_1,\ldots, f_n/S_n\},$$ hence is strict.
It follows therefore from the case $\Gamma=\{1\}$ 	already proven that $\theta(\mu\inv(\mathsf V))$ is a strict
quasi-compact open subset of~$\P_{F(s\inv S)/K(s\inv S}$.
By Corollary~\ref{coro-extension-bivaluation} and Example \ref{ex-tensor2}, one has
$\theta(\mu\inv (\mathsf V))=\nu\inv (\phi(\mathsf V))$, which
formally implies that $\phi(\mathsf V)=\sigma\inv(\theta(\mu\inv(\mathsf V))$. Therefore in order to conclude, it suffices to show that $\sigma\inv(\mathsf W)$ is a $\Gamma$-strict quasi-compact open subset of $\P_{F/K}$ for every {\em strict}
quasi-compact open subset $\mathsf W$ of $\P_{F(s\inv S)/K(s\inv S)}$. 

So, let us consider such a $\mathsf W$. We may and do assume that $\mathsf W=\P_{L(s\inv S)/K(s\inv S)}\{g\}$ for some homogeneous element $g$ of degree 1. Let us write 
$g=\sum a_I S^I/\sum b_J S^J$ (with $\sum b_I S^J\neq 0$)
and let $\abs \cdot$ be a
graded valuation belonging to $\P_{F/K}$. One has the equivalence
$$\abs \cdot \in \sigma\inv(\mathsf W)\iff \frac{\max_I \abs{a_I}}{\max_J \abs{b_J}}\leq 1.$$

As a consequence,
$$\sigma\inv(\mathsf W)=\bigcap_I \bigcup_{J,b_J\neq 0} \P_{F/K}\{a_I/b_J\}.$$
Fix $I$ and $J$ with $b_J\neq 0$. Since $g$ is homogeneous of degree 1, the elements $a_Is^I$ and $b_JS^J$ of $L(s\inv S)$ are homogeneous of the same degree. It follows
that $a_I/b_J$ is homogeneous of degree $s^J/s^I=s^{J-I}$, which belongs to $\Gamma$. Hence $\sigma\inv(\mathsf W)$ is a $\Gamma$-strict quasi-compact open
subset of $\P_{F/K}$, which ends the proof. 
\end{proof}

\subsection{The case of a normal graded extension}\label{ss-zr-normal}
Let $K$ be a graded field, let $F$ be a graded extension of $K$ and let $L$ be an algebraic
graded algebraic
extension of $F$ (\ref{ss-graded-algebraic}). 
We assume moreover
that $L$ is normal over $F$ and we set $G=\mathrm{Gal}(L/F)$ (\ref{ss-graded-galois}). 
The group $G$ acts on $\P_{L/K}$ and $\P_{L/K}\to \P_{F/K}$ is $G$-equivariant. We
are going to prove that $\P_{L/K}/G\to \P_{F/K}$ is a homeomorphism.

This map is surjective, 
and it is open by Theorem \ref{thm-elim-graded}. It suffices
therefore to prove that it is injective. Let $\abs \cdot _1$ and $\abs \cdot_2$
be two graded
valuations on $L$ that have the same restriction to $F$, and let us
show that they are conjugate under $G$. By Lemma \ref{lem-extension-bivaluation}, there 
exists a graded extension $M$ of $F$, a graded
valuation $\abs \cdot$ on $M$, and two $F$-embeddings $j_1$ and $j_2$
from $L$ to $M$ such that $\abs \cdot_1=\abs \cdot \circ j_1$ and
$\abs \cdot _2=\abs \cdot \circ j_2$. Let $r$ be a positive real number. 
Since $F\hookrightarrow L$ is normal, $j_1(L)^r$ is the subset of $M^r$ 
consisting of elements $x$ such that there exists an element of
$L^r$ whose
minimal polynomial over $F$ vanishes at $x$, and the same holds for $j_2(L)^r$. 
As a consequence, $j_1(L)=j_2(L)$
and
$\abs \cdot _1=\abs \cdot _2 \circ j_2\inv \circ j_1$; otherwise said, 
$\abs \cdot _1$ and $\abs \cdot_2$ are conjugate by the element 
$j_2\inv \circ j_1$ of $G$. 

\begin{rema}\label{rem-zr-radicial}
Assume that $L$ is radicial over $F$ (\ref{ss-graded-galois}). The group $G$ is then trivial, and
\ref{ss-zr-normal} thus ensures that $\P_{L/K}\to \P_{F/K}$ is a homeomorphism.
But this can easily be seen directly. Indeed, let $\abs \cdot $ be a
graded valuation on $L$, 
and let $a$ be a homogeneous element of $L$. As $L$ is radicial 
over $K$, there exists $N\in \N$ such that $a^N\in F$; since
$\abs a\leq 1\iff \abs{a^N}\leq 1$, we see that $\abs \cdot $ is uniquely determined
by its restriction to $F$, whence our claim. 
\end{rema}

\section{Applications to analytic geometry}
\label{s-smallest-image}

We fix an analytic field $k$ and a subgroup $\Gamma$ of $\R_+\gpm$ such that $\abs{k\gpm}\cdot \Gamma\neq \{1\}$. 
The purpose of what follows is to give some consequences of Theorem \ref{thm-elim-graded}
concerning morphisms between $k$-analytic germs. We shall use freely Temkin's reduction of analytic germs,
as well as its $\Gamma$-graded avatar; 
see Sections \ref{s-temkin} and \ref{s-temkin-gamma}. 

\begin{theo}\label{thm-smallest-germ}
Let $(Y,y)\to (X,x)$ be a morphism of $k$-analytic germs. 

\begin{enumerate}[1]
\item There exists a smallest analytic domain $(Z,x)$ of $(X,x)$ through which the morphism
$(Y,y)\to (X,x)$ factorizes; the reduction $\widetilde{(Z,x)}\subset \widetilde{(X,x)}$ is equal to
the image of the map $\widetilde{(Y,y)}\to \widetilde{(X,x)}$. 

\item If both germs $(Y,y)$ and $(X,x)$ are $\Gamma$-strict, then $(Z,x)$ is $\Gamma$-strict too, and $\widetilde{(Z,x)}^\Gamma\subset \widetilde{(X,x)}^\Gamma$ is equal to
the image of the map $\widetilde{(Y,y)}^\Gamma\to \widetilde{(X,x)}^\Gamma$. 

\item If $(Y,y)$ is $\Gamma$-strict and if $(X,x)$ is separated (but not necessarily $\Gamma$-strict), then $(Z,x)$ is $\Gamma$-strict.

\end{enumerate}

\end{theo}

\begin{proof}

Let $p$ denote the natural map $\widetilde {(Y,y)}\to \widetilde{(X,x)}$. 
Let us
choose an atlas $\mathscr U$ of $\widetilde{(X,x)}$
(\ref{ss-def-slk})
and an atlas $\mathscr V$ of $\widetilde{(Y,y)}$ such that the covering $\mathscr V$ is a refinement of $p\inv(\mathscr U)$
(this is possible because $p$ is quasi-compact, see~\ref{ss-rz-funct}); if moreover $(Y,y)$ and $(X,x)$ are $\Gamma$-strict, we require that $\mathscr U$ and $\mathscr V$ be
$\Gamma$-strict.  Let $\mathsf V$ be an element of
$\mathscr V$, and let $\mathsf U$ be an element of $\mathscr U$ containing $p(\mathsf V)$. 
By considering the diagram
$$\xymatrix{
{\mathsf V}\ar@{^{(}->}[r]\ar[d]_p&{\P_{\hrt y/\widetilde k}}\ar[d]\\
{\mathsf U}\ar@{^{(}->}[r]&{\P_{\hrt x/\widetilde k}}
}$$
and applying
Theorem~\ref{thm-elim-graded}
we see that $p(\mathsf V)$ is a quasi-compact open subset of $\widetilde{(X,x)}$, which is $\Gamma$-strict whenever
$(Y,y)$ 
and $(X,x)$ are $\Gamma$-strict (because then so is $\mathsf V$ by construction). In view of the equality
$$p(\widetilde{(Y,y)})=\bigcup_{\mathsf V\in \mathscr V} p(\mathsf V),$$ this implies that
$p(\widetilde{(Y,y)})$ is a (necessarily non-empty) quasi-compact open subset of
$\widetilde{(X,x)}$, which is $\Gamma$-strict as soon as $(Y,y)$ and $(X,x)$ are. 

This non-empty quasi-compact open subset of $\widetilde{(X,x)}$ is equal to $\widetilde{(Z,x)}$ for a uniquely determined
analytic domain $(Z,x)$ of $(X,x)$, which is $\Gamma$-strict as soon as $(Y,y)$ and $(X,x)$ are
(\ref{ss-temred-properties} (1) and Lemma \ref{lem-gstrict-general}).
In view of assertion (5) of \ref{ss-temred-properties}, $(Z,x)$ is the smallest analytic domain of~$(X,x)$ through which~$p$ factorizes, which ends the proof of  (1), and of the first statement of (2).
If both $(Y,y)$ and $(X,x)$ are $\Gamma$-strict, we have a commutative diagram 
$$\xymatrix{
{\widetilde{(Y,y)}}\ar[d]\ar[r]&{\widetilde{(Y,y)}^\Gamma}\ar@/^2pc/[dd]\\
{\widetilde{(Z,x)}}\ar@{^{(}->}[d]\ar[r]&{\widetilde{(Z,x)}^\Gamma}\ar@{^{(}->}[d]\\
{\widetilde{(X,x)}}\ar[r]&{\widetilde{(X,x)}^\Gamma}
}$$
in which $\widetilde{(Y,y)}\to \widetilde{(Z,x)}$ is surjective by definition of $(Z,x)$, and in which all horizontal arrows are surjective. 
It follows that $\widetilde{(Y,y)}^\Gamma\to \widetilde{(X,x)}^\Gamma$ factorizes through a surjection
$\widetilde{(Y,y)}^\Gamma\to \widetilde{(Z,z)}^\Gamma$, which ends the proof of (2). 

Now we assume that $(Y,y)$ is $\Gamma$-strict and that $(X,x)$ is separated. Let 
$\mathscr V$ be any $\Gamma$-strict atlas of $\widetilde{(Y,y)}$. Let $\mathsf V$ be
a chart belonging to $\mathscr V$. By considering
the commutative diagram $$\xymatrix{
{\mathsf V}\ar@{^{(}->}[r]\ar[d]_p&{\P_{\hrt y/\widetilde k}}\ar[d]\\
{\widetilde{(X,x)}}\ar@{^{(}->}[r]&{\P_{\hrt x/\widetilde k}}
}$$
and applying
Theorem~\ref{thm-elim-graded}
we see that $p(\mathsf V)$ is a $\Gamma$-strict quasi-compact open subset 
of $\P_{\hrt x/\widetilde k}$. In view of the equality
$$p(\widetilde{(Y,y)})=\bigcup_{\mathsf V\in \mathscr V} p(\mathsf V),$$ this implies that
$\widetilde{(Z,x)}=p(\widetilde{(Y,y)})$ is a $\Gamma$-strict quasi-compact open subset of $\P_{\hrt x/\widetilde k}$; therefore $(Z,x)$ is $\Gamma$-strict, which proves (3). 
\end{proof}

\subsection{}\label{ss-smallest-basics}
The following facts follow straightforwardly from the characterization of $(Z,x)$ by its
reduction $\widetilde{(Z,x)}$:

\begin{enumerate}[1]
\item If~$(Y_1,y),\ldots, (Y_n,y)$ are analytic domains of~$(Y,y)$ such that~$(Y,y)=\bigcup (Y_i,y)$ and
if~$(Z_i,x)$ denotes (for every~$i$) the smallest analytic domain of~$(X,x)$
through which~$(Y_i,y)\to(X,x)$ factorizes, then~$(Z,x)=\bigcup (Z_i,x)$. 

\item If $(X,x)\to (T,t)$ is another morphism of germs and if $(W,t)$ denotes the smallest analytic domain of $(T,t)$ through which $(Z,x)\to (T,t)$ factorizes, then
$(W,t)$ is also the smallest analytic domain of $(T,t)$ through which the composite map
$$(Y,y)\to (X,x)\to (T,t)$$ factorizes. 
\end{enumerate}

\begin{exem}\label{ex-smallest-trivial}
Let $X$ be a $k$-analytic space and let $V$ be an analytic domain of $X$. If $x$ is
a point of $V$, it follows immediately from the definition that $(V,x)$ is the smallest analytic domain 
of $(X,x)$ through which the
map $(V,x)\hookrightarrow (X,x)$ factorizes. 
\end{exem}

\begin{exem}\label{ex-smallest-boundaryless}
Let $X$ and $Y$ be two $k$-analytic spaces and let $\phi \colon Y\to X$
be a boundaryless morphism (\eg, a closed immersion).
Let $y$
be a point of $Y$ and let $x$ be its image in $X$. 
Since $\widetilde{(Y,y)}=\widetilde{(X,x)}\times_{\P_{\hrt x/\widetilde k}}\P_{\hrt y/\widetilde k}$, 
the continuous map $\widetilde{(Y,y)}\to \widetilde{(X,x)}$ is surjective. Hence $(X,x)$ is the smallest
analytic domain of $(X,x)$ through which $(Y,y)\to (X,x)$ factorizes. 

This has the following concrete meaning: if $V$ is any analytic domain of $X$ containing $x$ 
such that $\phi\inv(V)$ is a neighborhood of $y$ in $Y$, then $V$ is a neighborhood of $x$ in $X$.

\end{exem}

\begin{exem} In view of \ref{ss-smallest-basics}
one may combine Examples \ref{ex-smallest-trivial}
and \ref{ex-smallest-boundaryless} and get the following. 
Let $X$ and $Y$ be two $k$-analytic spaces, let $V$ be an analytic domain of $X$,
and let $Y\to X$ be the composition of a boundaryless
morphism $Y\to V$ and of the inclusion
$V\hookrightarrow X$. 
Then $(V,x)$ is the smallest analytic domain of $(X,x)$ through which 
$(Y,y)\to (X,x)$ factorizes. 
\end{exem}

We are now going to give
the first {\em global}
consequence of Theorem \ref{thm-smallest-germ}.
Note that the source is not assumed to be Hausdorff, 
but that the target is assumed to be separated, which is stronger than Hausdorff.

\begin{prop}\label{prop-image-subgstrict}
Let $\phi \colon Y\to X$ be morphism between $k$-analytic spaces. Assume that $Y$ is quasi-compact and $\Gamma$-strict
and $X$ is separated. The image $\phi(Y)$ is contained in a compact $\Gamma$-strict analytic domain of $X$.
\end{prop}

\begin{proof}
Let $y$ be a point of $Y$, and set $x=\phi(y)$. It follows from Theorem \ref{thm-smallest-germ}
that there exists a smallest analytic domain $(Z,x)$ of $(X,x)$ through which the map $(Y,y)\to (X,x)$ factorizes, and that it is
$\Gamma$-strict. Using simply the fact that the map $(Y,y)\to (X,x)$ factorizes through $(Z,x)$, we get the existence of an analytic neighborhood
$V_y$ of $y$ in $Y$ and a compact $\Gamma$-strict analytic neighborhood $U_x$ of $x \in X$ such that $\phi(V_y)\subset U_x$. 
By quasi-compactness of $Y$ there is a finite subset $E$ of $U$ such that $Y=\bigcup_{y\in E}V_y$. We then have $\phi(Y)\subset \bigcup_{y\in E}U_{\phi(y)}$. 
Since $X$ is separated and
every $U_{\phi(y)}$ for $y\in E$ is a compact, $\Gamma$-strict analytic domain of $X$, 
$\bigcup_{y\in E}U_{\phi(y)}$ is a compact $\Gamma$-strict analytic domain of $X$
(Remark \ref{rem-gstrict-careful}). This ends the proof. 
\end{proof}

\section{Complement: around quantifier elimination in the theory {\sc acvf}}\label{appendix-acvf}
\markright{\thesection. QUANTIFIER ELIMINATION IN THE THEORY ACVF}

In
the proof of Theorem \ref{theo-elim-sections}
we have referred at one
point to \emph{quantifier elimination} in the theory of non-trivially valued
algebraically closed fields. The aim of this section is to quickly recall what it consists of, to explain
why it implies the statement that was needed for showing Theorem \ref{theo-elim-sections},
and to give another proof 
of this statement, based upon classical arguments of algebraic geometry. 

\subsection{}Let us first define recursively what a {\em formula}\footnote{
The formulas we define here are formulas {\em in the language of valued fields};
but since we shall not use any other language in this section, we do not need such
precise terminology.}
is;
for the moment, it should be seen as a purely formal
syntactic object, without any meaning.
We fix a countable set of symbols which are called
the {\em variables}. Every formula $\Phi$ will involve 
two disjoint\footnote{The general definition of a formula is much more complicated: in a given formula
the same variable
may occur at some
places  with the ``free" status and at some other
places
with the ``bound" status. 
But any formula in this sense is, up to renaming some bound
variables, equivalent (at
least as far as the ``concrete" interpretation is concerned)
to a formula in our sense; that is why we have chosen this simpler definition.} 
finite
sets of variables: the set $\mathscr V_{\mathrm b}(\Phi)$ of
{\em bound}
variables, and the set $\mathscr V_{\mathrm f}(\Phi)$ of {\em free}
variables (the status -- free or bound -- of a variable is not an absolute notion: a variable
is free or bound {\em in a given formula}). 
 
\begin{enumerate}[1]

\item Let $\Phi$ be an inequality 
of the form $\abs P \Join \abs Q$, where $P$ and $Q$ belong to $\Z[x_1,\ldots, x_n]$ 
for some variables $x_1,\ldots, x_n$ 
and where $\Join$
is a symbol belonging to $\{<,>,\leq,\geq \}$. Then $\Phi$ is a formula, 
$\mathscr V_{\mathrm b}(\Phi)=\emptyset$ and $\mathscr V_{\mathrm f}(\Phi)$ is the set of $x_i$'s 
that actually occur in $P$ or $Q$. 

\item If $\Phi$ is a formula, its negation $\Psi$ is a formula too, and
one has
$$\mathscr V_{\mathrm b}(\Psi)=\mathscr V_{\mathrm b}(\Psi)\;{\rm and}\;
\mathscr V_{\mathrm f}(\Psi)=\mathscr V_{\mathrm f}(\Psi).$$

\item If $\Phi$ is a formula and if $x$ belongs to $\mathscr V_{\mathrm f}(\Phi)$,
then $\Psi:=(\forall x,\Phi)$ is a formula and one has
$$\mathscr V_{\mathrm f}(\Psi)=\mathscr V_{\mathrm f}(\Phi)\setminus \{x\}\;{\rm and}\;
\mathscr V_{\mathrm b}(\Psi)=\mathscr V_{\mathrm b}(\Phi)\cup \{x\}.$$

\item If $\Phi$ and $\Psi$ are two formulas with $(\mathscr V_b(\Phi)\cup \mathscr V_b(\Psi))\cap (\mathscr V_f(\Phi)\cup \mathscr V_f(\Psi)=\emptyset$
then $\Theta:=(\Phi\;\mathrm{and}\;\Psi)$ is a formula
and $$\mathscr V_{\mathrm b}(\Theta)=\mathscr V_{\mathrm b}(\Phi)\cup \mathscr V_{\mathrm b}(\Psi)\;{\rm and}\;
\mathscr V_{\mathrm f}(\Theta)=\mathscr V_{\mathrm f}(\Phi))\cup \mathscr V_{\mathrm f}(\Psi).$$

\end{enumerate}

For instance, 
$$\exists x, \forall y, (\abs{3x^2-8y+7}\leq \abs{-x+z+t})\;\mathrm{or}\;(\abs {5xy-y^2}>\abs{w+x+2y^3-zt})$$ is a formula
with free variables $z,t$ and $w$.

\subsection{}\label{ss-interp-formula}
Let $\Phi$ be a formula and let $x_1,\ldots, x_r$ be the free variables of $\Phi$. Let $K$ be a valued field and let $a_1,\ldots, a_r$
be elements of $K$. By replacing $x_i$ with $a_i$ for every $i$, one gets a statement $\Phi(a_1,\ldots, a_r)$ whose truth value in $K$, and more generally in any valued
extension of $K$, makes sense. 

\begin{rema}\label{rema-qfree}
Assume that $\Phi$ is quantifier-free (equivalently,
it only involves free variables).
Then for every valued extension $L$ of $K$ and every
$r$-uple $(a_1,\ldots, a_r)$ of elements of $K$, the statement
$\Phi(a_1,\ldots, a_r)$ holds in $L$ if and only if it holds in $K$. 

\end{rema}

We are now going to state the {\em quantifier elimination}
theorem for the theory of algebraically closed (non-trivially) valued field, {\sc acvf} for short. 
It is usually attributed to Robinson, though the result 
proved in his book \cite{robinson1977}
is weaker -- but the main
ideas are there; for a complete proof, \cf
\cite{prestel1986}
or~\cite{weispfenning1984}. 

\begin{theo}[quantifier elimination in {\sc acvf}]\label{thm-qe-acvf}
Let $\Phi$ be any formula with free variables $(x_1,\ldots, x_r)$.
There exists a {\em quantifier-free}
formula $\Psi$ with free variables $(x_1,\ldots, x_r)$ such that for every algebraically closed, non-trivially valued field $K$\index{quantifier elimination (for {\sc ACVF})}
and any $r$-tuple $(a_1,\ldots, a_r)$ of elements on $K$, the statement $\Phi(a_1,\ldots, a_r)$ holds in $K$
if and only if $\Psi(a_1,\ldots, a_r)$ holds in $K$. 

\end{theo}

In view of Remark~\ref{rema-qfree}
above, quantifier elimination in {\sc acvf}
implies the so-called {\em model-completeness}
of {\sc acvf}: 

\begin{theo}[model-completeness of {\sc acvf}]\label{thm-mc-acvf}
Let $\Phi$ be any formula with free variables $(x_1,\ldots, x_r)$, let $K$ be an algebraically closed, 
non-trivially valued field and let $L$ be an algebraically closed valued extension of $K$. 
For every $r$-uple $(a_1,\ldots, a_r)$ of elements of $K$, 
the statement $\Phi(a_1,\ldots, a_r)$ holds in $L$ if and only if it holds in $K$.
\end{theo}

\begin{rema}
The assumption that $K$ is non-trivially valued cannot be removed. Indeed, let $K$ be a trivially valued algebraically closed field; choose 
an algebraically closed, non-trivially valued extension $L$ of $K$. The statement
$$\exists x, (x\neq 0)\;{\rm and}\;\abs x <1$$ then 
holds in $L$, but
not in $K$.

\end{rema}

The {\em Nullstellensatz}
in the trivially valued case and model-completeness of {\sc acvf} (Theorem \ref{thm-mc-acvf} above)
in the non-trivially valued case immediately imply the following proposition, which was
the key ingredient
in the proof of Theorem \ref{theo-elim-sections}.

\begin{prop}\label{prop-mc-acvf}
Let $K$ be an algebraically closed valued field,
and let $L$ be a valued extension of $K$; let $(P_1,\ldots,P_m)$ be elements
of $K[T_1,\ldots, T_n]$ for some $n$. If the system of equations $\{P_i=0\}_{1\leq i\leq m}$ has a solution $(x_1,\ldots, x_n)$ in $L^n$ such that $\abs {x_i}\leq 1$ 
for all $i$, then it already has a solution $(y_1,\ldots, y_n)$ in $K^n$ such that $\abs{y_i}\leq 1$ for all $i$.
\end{prop} 

We now aim at giving a direct, algebro-geometric proof of Proposition~\ref{prop-mc-acvf}
above. We shall in fact prove the following slightly stronger result.

\begin{theo}\label{thm-MC-acvf}
Let $K$ be an algebraically closed
valued field and let $A$ be its valuation ring. Let $\mathscr X$ be a finitely presented $A$-scheme
and let $B$ be a faithfully flat $A$-algebra. If $\mathscr X(B)\neq \emptyset$, then
$\mathscr X(A)\neq \emptyset$.

\end{theo}

Before proving it, let us explain 
why this implies Proposition
\ref{prop-mc-acvf}. 
We use the notation of \loccit~and call $A$, \resp $B$,  the valuation ring of $K$, \resp $L$.
By multiplying every $P_i$ 
by a suitable element of $K\gpm$, we may and do assume that $P_i$
belongs to $A[T_1,\ldots, T_n]$ for all $i$. Let $\mathscr X$ be the 
$A$-scheme $\spec A[T_1,\ldots, T_n]/(P_1,\ldots, P_m)$. The $A$-algebra $B$ is faithfully flat, 
and $\mathscr X(B)\neq \emptyset$ by assumption.	
Therefore if we assume that Theorem \ref{thm-MC-acvf}
holds, 
we can
conclude that $\mathscr X(A)\neq \emptyset$, which is exactly the required assertion. 

\begin{proof}[Proof of Theorem \ref{thm-MC-acvf}]
By assumption, $\mathscr X$ has a $B$-point. Since $\mathscr X$ is finitely presented, this $B$-point is induced by a $B'$-point of $\mathscr X$ for some
finitely generated subalgebra $B'$ of $B$. Since being $A$-flat simply means having no
non-zero $A$-torsion, the $A$-algebra $B'$ is flat; being finitely generated, it is finitely presented
by a result of Nagata \cite{nagata1966} (this holds more generally for $A$ any {\em domain}; see \cite{raynaud-g1971}, \Cor 3.4.7); moreover, since $\spec B\to \spec A$ goes through $\spec B'$, the map
$\spec B'\to \spec A$ is surjective, which means that $B'$ is faithfully flat over $A$.
Hence by replacing $B$ with $B'$, we may and do assume that $B$ is finitely presented over $A$. It is now sufficient to prove that $\spec B\to \spec A$ has a section, because by composing the latter with any $B$-point $\spec B\to \mathscr X$ we will get an $A$-point
$\spec A\to \mathscr X$. 

As $B$ is finitely presented over $A$, there exists a subring $A_0$ of $A$ finitely generated over $\Z$, and a flat $A_0$-algebra of finite type $B_0$ 
such that $B\simeq A\otimes_{A_0}B_0$. Let $x$ be the closed point of $\spec A$, and let $x_0$ be the image of $x$ on $\spec A_0$. Since the fiber
$(\spec B)_x$ is non-empty, the finitely generated $\kappa(x_0)$-scheme $(\spec B_0)_{x_0}$ is non-empty; its
CM locus is thus non-empty (it contains the maximal points), and therefore has a closed point $y$. Let $(f_1,\ldots, f_r)$ be a finite family of elements of $B_0$ that 
lifts a maximal regular sequence of $(\spec B_0)_{x_0}$ at $y$, and let $Z$ be the closed subscheme of $\spec B_0$ defined by the $f_i$'s. At the point $y$, the $A_0$-scheme $Z$
is quasi-finite by construction, and flat by \cite{ega43}, \Th 11.3.8 (c). 
By Zariski's Main Theorem, there is an open neighborhood $U$ of $y$ in $Z$ and an open immersion $U\hookrightarrow V$ over $A_0$ for some finite $A_0$-scheme $V$. 

Set $U'=U\times_{A_0}A$ and $V'=V\times_{A_0}A$.  By construction, $U'$ is a locally closed subscheme of $\spec B$, flat over $A$ at some point $y'$ lying over $x$; the scheme $V'$ is finite
over $A$, and there is an open immersion of $A$-schemes $U'\hookrightarrow V'$. 

Since $K$ is algebraically closed, the valuation ring $A$ is henselian. The finite $A$-scheme $V'$ is therefore a disjoint union of finitely many finite, local $A$-schemes. This implies that the connected component $W$ of $y'$ in $U'$ is a finite, local $A$-scheme. Since $U'$ is flat over $\spec A$ at $y'$, so is $W$; as a consequence, the image of $W\to \spec A$ contains the generic point of $\spec A$.
Let us choose a point $w$ on the generic fiber of $W$, and let
endow the Zariski-closed subset $\overline{\{w\}}$ of $\spec B$ with its reduced structure. One has $\overline{\{w\}}=\spec C$ for some finite $A$-algebra $C$; the ring $C$ is a domain and the map $A\to C$ is injective. 
The tensor product $C\otimes_A K$ is a localization of the domain $C$ by a mutliplicative subset which does not
contain $0$, so $C\otimes_A K$ is a domain. Since it is finite over $K$, this is a field
and thus a finite extension of $K$. 
As $K$ is algebraically closed, it follows that $C\otimes_A K=K$. Since $A$ is normal, 
this implies that $C=A$, and the closed embedding $\overline{\{w\}}\hookrightarrow \spec B$ then defines a section of $\spec B\to \spec A$. 
\end{proof}

\chapter{D\'evissages \`a la Raynaud-Gruson}\label{DEV} 

Most of this chapter 
is inspired by the celebrated work of Raynaud and Gruson
on flatness \cite{raynaud-g1971}. In that paper,
they consider a finitely presented morphism of schemes $X\to S$
and a quasi-coherent $\mathscr O_X$-module of finite type $\mathscr M$. The
key notion they introduce
is that of a \emph{d\'evissage}
of $\mathscr M$ over $S$ at a given point of $\supp M$
(\cite{raynaud-g1971}, \Def 1.2.2); they prove that
it always exists after some 
Nisnevich localization on the target and on the source; see
\Prop 1.2.3 of \cite{raynaud-g1971}~for the precise statement. 
They use d\'evissages for studying the 
$S$-flat locus of $\mathscr M$, describing the local structure
of $\mathscr M$ at a point at which it is $S$-flat, or flattening it through a 
blow-up in general\ldots

Now, let $Y\to X$ be a morphism of good
$k$-analytic 
spaces, and let $\mathscr F$ be a coherent sheaf on $Y$.
We define the notion of an $X$-d\'evissage of $\mathscr F$
at a given point of $\supp F$
(Definition \ref{def-devissages}), and prove that such a d\'evissage
always exists
(see Theorem 
\ref{thm-devissages}); note that there is no need
for Nisnevich localization here:
it suffices to work on a small enough
affinoid neighborhood of $y$ in $Y$ -- this
ultimately relates
to the henselian property
of local rings of good analytic spaces. 

We then give two applications of d\'evissages. We first use them to prove
that if $y\in \mathrm{Int}(Y/X)$ and $\mathscr F$ is naively $X$-flat at $y$, then 
$\mathscr F$ is $X$-flat at $y$ (Theorem \ref{thm-flat-naiveflat}).
Otherwise said, naive flatness at a relatively inner point
is automatically universal;
\ie, it remains true after arbitrary good base change (including ground field extension). 
More precisely, in the spirit of
\Cor 2.3 of \cite{raynaud-g1971}, Theorem \ref{thm-flat-naiveflat} characterizes
naive $X$-flatness and $X$-flatness of $\mathscr F$ at $y$ in terms of an $X$-d\'evissage
of $\mathscr F$ at $y$ (provided $y\in \supp F$; if not, $\mathscr F$
is obviously $X$-flat at
$y$). It turns out that both characterizations are equivalent when $y\in \mathrm{Int}(Y/X)$, 
which ultimately rests on the nice
properties of local rings of generic fibers at inner points (Theorem \ref{thm-localring-generic}). 

The second application of d\'evissages concerns the local structure of relatively CM (\ie, flat and fiberwise CM)
coherent sheaves. There are two basic examples of such sheaves:

\begin{enumerate}[a]

\item If $Y\to X$ is finite and $\mathscr F$ is
$X$-flat at $y$, then $\mathscr F$
is relatively CM at $y$. 

\item If $Y\to X$ is quasi-smooth at $y$, then $\mathscr O_Y$
is relatively CM at $y$ (this follows from Theorem \ref{thm-main-qsm}).

\end{enumerate}
Now Theorem \ref{thm-local-cm}
essentially states that every coherent sheaf relatively CM
at a point arises around this
point as a ``combination"
of 
a relatively CM coherent sheaf of the kind described in (a), and of another one of the kind
described in (b). 

Though the two aforementioned applications of
d\'evissages are the only ones in this memoir, we hope that they
will be useful in the future for other purposes, such as for
development of
flattening techniques in the non-archimedean setting.

\section{Universal injectivity and flatness} 
The purpose of this section is to introduce a technical notion, namely that of
\emph{universal injectivity},
and to study how it interacts with (naive and non-naive)
flatness. 
We shall use freely and repeatedly the notions of exactness (of a complex of coherent sheaves), and injectivity, surjectitvity 
and bijectivity (of a morphism between coherent sheaves) {\em at a given point of an analytic space}, and the related affinoid GAGA principles; 
see Lemma-Definition \ref{valid-at-concrete} (3) and (4), 
and Lemma \ref{gaga-concrete} (3) and (4). 

We shall also use freely  basic results about the fibers of coherent sheaves (\ref{s-fibers-coherent}), and especially
the fact that the property of being zero at a point (for a coherent sheaf) or of being surjective at a point (for a map of coherent sheaves)
can be checked fiberwise; this follows basically from Nakayama's Lemma, see \ref{ss-pointwise-rk} and \ref{ss-surj-nakayama} for details.

\begin{defi}\label{def-univ-inj}\index{morphism!of coherent sheaves!universally injective}\index{universally injective morphism (of coherent sheaves)}
Let~$Y\to X$ be a morphism between~$k$-analytic spaces,
and let~${\mathscr G}\to {\mathscr F}$ be a linear map between coherent sheaves on~$Y$. Let~$y$ be a point of $Y$.
We say that~${\mathscr G}\to {\mathscr F}$ is {\em~$X$-universally injective at~$y$}
if for every analytic space~$X'$, for every morphism~$X'\to X$, and for every point~$y'$ lying above~$y$ on~$Y':=Y\times_{X}X'$,  the map~$\mathscr G_{Y'}\to \mathscr F_{Y'}$ is injective at~$y'$. 

\end{defi}

\begin{rema}
Definition \ref{def-univ-inj} of universal injectivity is equivalent to the same
in which one (apparently) weakens the condition by taking for~$X'$ a good space, or even an affinoid one. 
\end{rema}

\subsection{Basic properties of universal injectivity}\label{ss-prop-univinj}
It follows from its definition that universal injectivity is preserved by any
base change (including ground field extensions). 

Let~$Y\to X$ be a morphism of $k$-analytic spaces,
and let ${\mathscr G}\to {\mathscr F}$ be a morphism between coherent sheaves on~$Y$. Let~$U$ be an analytic domain of~$X$ and 
let $V$ be an analytic domain of~$Y\times_XU$. For every $y\in V$, the map~$\mathscr G\to \mathscr F$ is~$X$-universally injective at~$y$ if and only if~$\mathscr G_V\to \mathscr F_V$ is~$U$-universally injective at~$y$: this comes from the fact that the validity of injectivity at a given point is insensitive to the restriction to an analytic domain.

\begin{rema}
We could also have defined universal surjectivity and universal bijectivity at a point in the same way, and obtained analogous properties. But there is no need for such notions, 
because surjectivity and bijectivity (of a morphism of coherent sheaves) at a point are {\em automatically universal}: indeed, after reduction to the good case this simply follows from
the fact that surjectivity and bijectivity (of a morphism of modules) are preserved by tensor product. 
\end{rema}

\subsection{About the bijectivity locus}\label{ss-bij-locus}
If~$X$ is an analytic space and if~$\mathscr F\to \mathscr G$ is a morphism
of coherent sheaves on~$X$, we shall denote by $\bij{\mathscr F}{\mathscr G}$ the set of points\label{IN-bijFG}
of $X$ at which $\mathscr F\to \mathscr G$ is bijective. By \ref{ss-pointwise-sheaf}, $\bij {\mathscr F}
{\mathscr G}$ is a Zariski-open subset of~$X$
(hence $\adht {\bij {\mathscr F}
{\mathscr G}}X=\adhz  {\bij {\mathscr F}
{\mathscr G}}X$ by Lemma \ref{lem-closure-zartop})
and $\mathscr F_{\bij{\mathscr F}{\mathscr G}}\to \mathscr G_{\bij{\mathscr F}{\mathscr G}}$ is an isomorphism. 
 
Let $Y\to X$ be a morphism of analytic spaces. It follows from the above that
$$(\bij {\mathscr F}{\mathscr G})_Y\subset \bij {\mathscr F_Y}{\mathscr G_Y}.$$
Moreover, if $y$
is a point of $Y$ at which $Y$ is $X$-flat (\eg, $Y$ is an analytic domain of $X$, or 
the space $X_L$ for some analytic extension $L$ of $k$),
then $y$ belongs to $\bij {\mathscr F_Y}{\mathscr G_Y}$
{\em if and only if}
it belongs to  $(\bij {\mathscr F}{\mathscr G})_Y$:
this is a particular case of Lemma \ref{lem-flat-desc1} (2) or of Lemma \ref{lem-flat-desc2} (2).

\subsection{}\label{ss-stand-not}
In this chapter, we shall often encounter
the following situation. We are given a quasi-smooth morphism $T\to X$
of $k$-affinoid~spaces
and a point $t$ of $T$
whose image in $X$ is
denoted by $x$,
and our investigation requires to choose
a point $z$ whose image $z\al_x$ on~$T_x\al$ is the generic point of
the connected component of~$t_x\al$ 
(we use the notation described in \ref{ss-conventions-gaga}).
Let us now make some basic remarks.

\begin{enumerate}[1]

\item Since $T_x$ is a quasi-smooth $\hr x$-analytic space, it is geometrically
regular by Lemma \ref{lem-pre-qsm} (2); therefore $T_x\al$
is a regular scheme by affinoid GAGA, see Lemma \ref{gaga-concrete} (1).
Hence
${\mathscr O}_{T_x\al, z_x\al}$ is the fraction field of the regular local ring~${\mathscr O}_{T_x\al, t_x\al}$.

\item The point $z$ can be chosen
in any given open subset~$U$ of~$T_{x}$
which intersects the connected component of~$t$ in $T_x$. Indeed, let~$n$ be the
dimension of~$T_{x}$ at~$t$ and let~$V$ be the intersection of~$U$ and
the
connected component of~$t$ in~$T_{x}$.
The dimension of~$V$ is also~$n$ and as a consequence,  there exists~$z$ in~$V$ such that $d_{\hr x}(z)=n$, which
fulfills our requirement (Remark \ref{rem-abh-point}). 

\item If $U$ is a Zariski-open subset of~$T_x$ such that $t$ belongs to $\adht U {T_x}$,
then $U$ has a non-empty intersection with the connected component of~$t$ in~$T_{x}$, hence contains~$z$. 

\item If $V$ is a Zariski-open subset of $T$ containing $t$, then it contains $z$: apply (3) above with $U=V\cap T_x$. 
As a consequence, $t\al$ belongs to $\adht {\{{z\al}\}}{T\al}$, so~${\mathscr O}_{T\al ,z\al}$ is
a localization of~${\mathscr O}_{T\al ,t\al}$.

\end{enumerate}

\subsection{}
We are now going to prove two technical results. 
The first one is
Proposition \ref{prop-equiv-univinj}
below,
which is the analogue of  \cite{raynaud-g1971}, Lemme 2.2; its statement
is not very enlightening, and it will only be used as an intermediate step for the
second technical result, namely
Proposition \ref{prop-pre-devisssage}.
The latter is
the analogue 
of \Th 2.1 of \cite{raynaud-g1971}, and its statement is designed for our purposes: it provides
a general flatness criterion, and also a criterion for naive flatness in the inner case, which
we shall apply to some of
the coherent sheaves that appear in a d\'evissage
(of course, the distinction between flatness and naive flatness has no counterpart
in the work of Raynaud and Gruson). 

\begin{prop}\label{prop-equiv-univinj}
Let~$T\to X$ be a quasi-smooth morphism  between~$k$-analytic spaces. Let~$\mathscr L$
be a free~${\mathscr O}_T$-module of finite rank and let~${\mathscr N}$ be a coherent sheaf on~$T$. Let~$t$ be
a point of $T$, and let~$x$ be its image on~$X$. Let~${\mathscr L}\to {\mathscr N}$ be a map such that~$t$ belongs to 
$\adht { \bij{{\mathscr L}_{T_{x}}}{{\mathscr N}_{T_{x}}}}{T_x}$. The following are equivalent:

\begin{enumerate}[i]
\item  The map~${\mathscr L}\to {\mathscr N}$ is~$X$-universally injective at~$t$.

\item The map~${\mathscr L}\to {\mathscr N}$ is injective at~$t$.

\item The point $t$ belongs to $\adht {\bij{\mathscr L}{\mathscr N}_x}{T_x}$. 
\end{enumerate}
\end{prop}

\begin{rema}\label{rem-equiv-univinj}
Let $X'$ be an analytic space, and let $X'\to X$ be a morphism; set $T'=T\times_X X'$. Let $t'$ be a pre-image of $t$ on $T'$, 
and let $x'$ be the image of $t'$ on $X'$. The coherent sheaf $\mathscr L_{T'}$ is then a free $\mathscr O_{T'}$-module. 
Moreover, $t'$ belongs to $\adht {\bij{\mathscr L_{T'_{x'}}}{\mathscr N_{T'_{x'}}}}{T'_{x'}}$. Indeed, 
$$\bij{\mathscr L_{T'_{x'}}}{\mathscr N_{T'_{x'}}}=(\bij{\mathscr L_{T_{x}}}{\mathscr N_{T_{x}}})_{\hr {x'}}$$ by
\ref{ss-bij-locus}, and our claim thus follows from Corollary \ref{cor-zarclosure-xl}. Hence 
the data $(T',X',t',x',\mathscr L_{T'}\to \mathscr N_{T'})$ also fulfills the assumptions of the proposition (possibly over a ground
field larger than $k$). 

Now assume that assertion (iii) is satisfied.
Corollary \ref{cor-zarclosure-xl} then 
ensures that $t'$ belongs to
$$\adht {(\bij{\mathscr L}{\mathscr N}_x)_{\hr {x'}}}{T'_{x'}}=\adht{(\bij{\mathscr L}{\mathscr N}_{T'})_{x'}}{T'_{x'}},$$
and since $\bij{\mathscr L}{\mathscr N}_{T'}\subset \bij{\mathscr L_{T'}}{\mathscr N_{T'}}$ by \ref{ss-bij-locus},
the point $t'$ belongs to $\adht {\bij{\mathscr L_{T'}}{\mathscr N_{T'}}_{x'}}{T'_{x'}}$; otherwise said, the analogue of
assertion (iii) with respect to the data $(T',X',t',x',\mathscr L_{T'}\to \mathscr N_{T'})$ also holds. 
\end{rema}

\begin{proof}[Proof of Proposition \ref{prop-equiv-univinj}]
We begin with the equivalence (ii)$\iff$(iii). We can assume that~$T$ and~$X$ are affinoid. 
We chose $z$ in $T$ such that $z_x\al$ is the generic point of the connected
component of $T_x\al$ containing $t_x\al$.

Assume that (ii) holds. This means that the arrow~$\mathscr L_t\to \mathscr N_t$ is injective. Therefore:

\begin{enumerate}[a]

\item The map $\mathscr L\al_{t\al} \to \mathscr N\al_{t\al}$ is injective.

\item The map $\mathscr L\al_{z\al}\to \mathscr N\al_{z\al}$ is injective by~(a) and \ref{ss-stand-not} (4). 

\end{enumerate}
By assumption,~$t$ belongs to $\adht{\bij{{\mathscr L}_{T_{x}}}{{\mathscr N}_{T_{x}}}} {T_x}$ ; therefore $\bij{{\mathscr L}_{T_{x}}}{{\mathscr N}_{T_{x}}}$ contains 
$z$ by \ref{ss-stand-not} (3) ; as a consequence, $\mathscr L_{\hr z}\to  \mathscr N_{\hr z}$ is an isomorphism, and is in particular
surjective. Therefore: 

\begin{enumerate}[a]\setcounter{enumi}{2}

\item The map $\mathscr L\to  \mathscr N$ is surjective at $z$. 

\item The map $\mathscr L\al_{z\al}\to \mathscr N\al_{z\al}$ is surjective by (c).

\item The map
$\mathscr L\al_{z\al} \to\mathscr N\al_{z\al}$ is bijective by (b) and (d). 

\end{enumerate}
By~(e), the point $z$ belongs to $\bij{\mathscr L}{\mathscr N}$;  since $t$
belongs the Zariski-closure of $z$ in $T_x$, it belongs to  $\adht{\bij{\mathscr L}{\mathscr N}_x}{T_x}$; hence (iii) holds.

Assume conversely that~(iii) holds, and let us prove (ii).
We can do it after an arbitrary scalar
extension: Remark \ref{rem-equiv-univinj}
ensures that our hypotheses (including (iii))
will remain the same; 
and Proposition \ref{prop-concrete-extension} (5)
ensures that (ii) will descent to the original ground field. We can thus assume that $x$
is a rational point.  

By assumption, $t$
belongs to~$\adht {\bij{\mathscr L}{\mathscr N}_x}{T_x}$; hence $\bij{\mathscr L}{\mathscr N}$
 contains $z$ by \ref{ss-stand-not} (3); 
therefore $\mathscr L\al_{z\al}\simeq  \mathscr N\al_{z\al}$.

We want to prove that
the top horizontal arrow of the following commutative diagram
$$\xymatrix{
{\mathscr L\al_{t\al}}\ar[r]\ar[d]&{\mathscr N\al_{t\al}}\ar[d]\\
{\mathscr L\al_{z\al}}\ar[r]&{\mathscr N\al_{z\al}}
}$$ is an injection; the bottom horizontal arrow being an isomorphism, it is enough to establish the injectivity of the left vertical arrow.
The~${\mathscr O}_T$-module~$\mathscr L$ is free of finite rank; therefore, it suffices to prove
that the map~$\mathscr O_{T\al, t\al}\to \mathscr O_{T\al, z\al}$ is injective.
We denote by $S$ be the multiplicative subset of ~${\mathscr O}_{T\al, t\al}$ 
that
consists of all elements~$a$ such that~$a(z\al)\neq 0$; by \ref{ss-stand-not} (4), we have $\mathscr O_{T\al, z\al}=S\inv \mathscr O_{T\al,t\al}$ . 

Let~$a$ be an element of $S$. Since~$a(z\al)$ is non-zero, the image of~$a$ in~$\kappa(z\al_x)$ is non-zero. But~$\kappa(z\al_x)$ coincides with~${\mathscr O}_{T\al_x, z\al_x}$; \ie, with~$\mathrm{Frac}(\mathscr O_{T\al_x,t\al_x})$. Therefore the image of~$a$ in the domain~$\mathscr O_{T\al_x,t\al_x}$ is non-zero, and hence is not a zero divisor.
On the other hand, ${\mathscr O}_{T_x\al, t_x\al}=\mathscr O_{T\al, t\al}/\mathfrak m_{x\al}{\mathscr O}_{T\al, t\al}$ because $x$ is a $k$-point;
and since $T$ is quasi-smooth over $X$, the ring $\mathscr O_{T\al, t\al}$ is flat over $\mathscr O_{X\al,x\al}$ by Theorem \ref{thm-qsm-schemereg} (2)
(or more directly by Corollary \ref{cor-qsm-flat} and Lemma \ref{lem-gagaflat-easy}). The above properties imply that the multiplication by~$a$ in~${\mathscr O}_{T\al,t\al}$ is injective
(\cite{matsumura1986}, \Th 22.5). The
set $S$ thus only consists of elements that are not zero divisors. As a consequence, 
the localization map
${\mathscr O}_{T\al, t\al}\to S\inv{\mathscr O}_{T\al, t\al}={\mathscr O}_{T\al, z\al}$ is
injective, and (ii) holds, whence the equivalence (ii)$\iff$(iii). 

Let us now prove that (i)$\iff$(ii), the spaces $T$ and $X$ being no longer assumed to be affinoid. 
The direct implication is tautological; it thus remain to show that (ii)$\Rightarrow$(i). So let us
assume that $\mathscr L\to \mathscr N$ is injective at $t$. Let $X'\to X$ be an arbitrary morphism, 
and set $T'=T\times_X X'$. Let $t'$ be a pre-image of $t$ in $T'$ and let $x'$ be the image of $t'$ in $X'$. By the implication (ii)$\Rightarrow$(iii) already proven, 
the point $t$ belongs to $\adht {\bij{\mathscr L}{\mathscr N}_x}{T_x}$. Remark \ref{rem-equiv-univinj}
then ensures that $t'$ belongs to $\adht {\bij{\mathscr L_{T'}}{\mathscr N_{T'}}_{x'}}{T'_{x'}}$; since it also ensures
that $(T',X',t',x',\mathscr L_{T'}\to \mathscr N_{T'})$ fulfills the assumptions of the proposition, the implication (iii)$\Rightarrow$(ii)
already proven yields injectivity of $\mathscr L_{T'}\to \mathscr N_{T'}$ at $t'$. Therefore $\mathscr L\to \mathscr N$ is $X$-universally injective at $t$. 
\end{proof}

\begin{prop}\label{prop-pre-devisssage}
Let~$T\to X$ be a quasi-smooth morphism  between~$k$-analytic spaces.
Let~$\mathscr L$ be a free~${\mathscr O}_T$-module of finite rank and let~${\mathscr N}$ be a coherent sheaf on~$T$.
Let~$t$ be a point of $T$, and let~$x$ be its image on~$X$. Let~${\mathscr L}\to {\mathscr N}$ be a map such that~$t$ belongs to 
$\adht { \bij{{\mathscr L}_{T_{x}}}{{\mathscr N}_{T_{x}}}}{T_x}$. Let~$\mathscr P$ be the cokernel of~${\mathscr L}\to {\mathscr N}$.

\begin{enumerate}[1]
\item The following are equivalent:

\begin{enumerate}[j]

\item  The coherent sheaf $\mathscr N$ is~$X$-flat at~$t$. 

\item The map~${\mathscr L}\to {\mathscr N}$ is injective at~$t$ and ~$\mathscr P$ is~$X$-flat at~$t$.

\end{enumerate}

\item Assume that $X$ and $T$ are good and that
$t$ belongs to $\mathrm{Int}(T/X)$. The following are equivalent:

\begin{enumerate}[j]\setcounter{enumii}{2}
\item The coherent sheaf $\mathscr N$ is
naively~$X$-flat 
at~$t$.

\item The map ${\mathscr L}\to {\mathscr N}$ is injective at~$t$
and ~$\mathscr P$ is
naively~$X$-flat at~$t$.
\end{enumerate}
\end{enumerate}
\end{prop}

\begin{proof}
We begin with (2). 
By shrinking $X$ and $T$, we may assume that both are affinoid, and that the maximal ideal $\mathfrak m_x$ is generated by
an ideal of $\mathscr O_X(X)$; we denote by $Y$ the corresponding closed analytic subspace of $X$, and by $S$ the fiber product $T\times_X Y$. 
We have by construction
$\mathscr O_{Y,x}=\mathscr O_{X,x}/\mathfrak m_x=\kappa(x)$; hence $\mathscr O_{Y,x}$ is a field, 
and $\mathscr O_{Y\al, x\al}$ is thus a field as well by \ref{ss-gaga-domain}; 
we therefore have $\mathscr O_{Y\al, x\al}=\mathscr O_{X\al, x\al}/\mathfrak m_{x\al}$, and thus
$\mathscr O_{S\al,\tau}=\mathscr O_{T\al,\tau}/\mathfrak m_{x\al} \mathscr O_{T\al, \tau}$
for every point $\tau$ of $T\al$ lying above $x\al$. 
We choose a point $z$ in the open neighborhood $\mathrm{Int}(T/X)_x$ of $t$ in $T_x$ such that $z_x\al$
is the generic point of the connected component of $T_x\al$
containing $t_x\al$, which is possible by \ref{ss-stand-not} (2). 
~Note that $z$ then also belongs to $\mathrm{Int}(S/Y)$, by base change. 

Assume that (iii) holds. 
Let~$\mathscr R$ be the kernel of
the map ${\mathscr L}\to {\mathscr N}$. The
point $t$ lies on $\adht  {\bij{{\mathscr L}_{T_{x}}}{{\mathscr N}_{T_{x}}}}{T_x}$. By \ref{ss-stand-not} (3), this implies that 
$z$ belongs to $\bij{{\mathscr L}_{T_{x}}}{{\mathscr N}_{T_{x}}}$. 
In particular, $\mathscr P_{\hr z}=0$; hence $\mathscr P_z=0$ and $\mathscr P\al_{z\al}=0$. 
The sequence
$$0\to \mathscr R\al_{z\al}\to  \mathscr L\al_{z\al}\to \mathscr N\al_{z\al}\to 0$$
is thus exact.

By assumption,~$\mathscr N_t$ is~$\mathscr O_{X,x}$-flat. By  Lemma \ref{lem-gagaflat-easy},  the
$\mathscr O_{X\al, x\al}$-module~$\mathscr N\al_{t\al}$ is therefore flat; in view of \ref{ss-stand-not} (4),
the  ~$\mathscr O_{X\al, x\al}$-module~$\mathscr N\al_{z\al}$ is also flat.  Hence the sequence
$$0\to \mathscr R\al_{z\al}/(\mathfrak m_{x\al} \mathscr R\al_{z\al})\to  \mathscr L\al_{z\al}/(\mathfrak m_{x\al} \mathscr L\al_{z\al})\to \mathscr N\al_{z\al}/(\mathfrak m_{x\al}\mathscr N_{z\al})\to 0$$ is exact; note that it can be rewritten 
$$0\to   \mathscr R\al_{S\al,z\al}\to  \mathscr L\al_{S\al,z\al}\to \mathscr N\al_{S\al, z\al}\to 0.$$
As a consequence, the sequence
$$0\to \mathscr R_{S,z}\to \mathscr L_{S,z}\to \mathscr N_{S,z}\to 0$$ is exact too. 
Since $\mathscr O_{Y,x}$ is a field
and $z$ belongs to $\mathrm{Int}(S/Y)$,
Theorem \ref{thm-localring-generic}
tells us that~${\mathscr O}_{S_x,z}$ is a flat ~${\mathscr O}_{S,z}$-algebra; this yields the exactness 
of the sequence 
$$ 0\to   \mathscr R_{S_x,z}\to  \mathscr L_{S_{x}, z}\to \mathscr N_{S_{x},z  }\to 0,$$ 
which can also be written $$ 0\to   \mathscr R_{T_{x}, z}\to  \mathscr L_{T_{x}, z}\to \mathscr N_{T_{x},z  }\to 0.$$
As $z$ belongs to $\bij{{\mathscr L}_{|T_{x}}}{{\mathscr N}_{|T_{x}}}$, we have $\mathscr R_{T_{x}, z}=0$; hence~$\mathscr R_{\hr z}=0$
and $\mathscr R_z=0$. Since $\mathscr P_z=0$
too, the point
$z$ belongs to the Zariski-open subset $\bij {\mathscr L}{\mathscr N}$; being a Zariski-specialization of $z$ inside $T_x$, the point 
$t$ belongs to $\adht{ \bij{ \mathscr L}{{\mathscr N}}_{x}}{T_{x,\mathrm{Zar}}}=\adht {\bij{ \mathscr L}{{\mathscr N}}_{x}}{T_x}$ (the equality comes
from Lemma \ref{lem-closure-zartop}).  By Proposition~\ref{prop-equiv-univinj}, this implies that $\mathscr L\to \mathscr N$ is injective at $t$, and even
\emph{$X$-universally}
injective
at $t$.
In particular, it remains injective after base-change by the map $\mathscr M(\hr x)\to X$ induced by $x$; this means that $\mathscr L_{T_x}\to \mathscr N_{T_x}$
is injective at $t$. The map
$\mathscr L_{T_x,t}\to \mathscr N_{T_x,t}$ is thus injective, and it can also be written
$\mathscr L_{S_x,t}\to \mathscr N_{S_x,t}$.
Since $\mathscr O_{X,x}$ is a field
and $t$ lies in
$\mathrm{Int}(S/Y)$, Theorem \ref{thm-localring-generic} ensures that $\mathscr O_{S_x,t}$ is a flat $\mathscr O_{S,t}$-algebra. 
As a consequence, the map $\mathscr L_{S,t}\to \mathscr N_{S,t}$ is injective; it can be rewritten 
$\mathscr L_t/(\mathfrak m_x \mathscr L_t)\to  \mathscr N_t/(\mathfrak m_x \mathscr N_t)$. 

Now since $\mathscr N_t$ is a flat $\mathscr O_{X,x}$-module by assumption, the
exact sequence $$0\to \mathscr L_t\to \mathscr N_t\to \mathscr P_t\to 0$$ induces a
long exact sequence
$$0\to \mathrm{Tor}_1^{\mathscr O_{X,x}}(\mathscr P_t,\kappa(x))\to \mathscr L_t/(\mathfrak m_x \mathscr L_t)\to  \mathscr N_t/(\mathfrak m_x \mathscr N_t)\to
\mathscr P_t/(\mathfrak m_x \mathscr P_t)\to 0$$ 
The injectivity of $\mathscr L_t/(\mathfrak m_x \mathscr L_t)\to  \mathscr N_t/(\mathfrak m_x \mathscr N_t)$
thus yields the equality
$$\mathrm{Tor}_1^{\mathscr O_{X,x}}(\mathscr P_t,\kappa(x))=0,$$ which implies that $\mathscr P_t$ is a flat
$\mathscr O_{X,x}$-module (\cite{sga1}, Expos\'e IV, \Th 5.6). Otherwise said, $\mathscr P$ is
naively $X$-flat at $t$, and (iv) holds.

Assume conversely that (iv) holds. The $\mathscr O_{X,x}$-module $\mathscr P_t$ is then flat. Since $T\to X$ is quasi-smooth, it is flat (Corollary \ref{cor-qsm-flat}) and 
$\mathscr O_{T,t}$ is thus flat over $\mathscr O_{X,x}$. As $\mathscr L$ is a free $\mathscr O_T$-module, $\mathscr L_t$ is flat
over $\mathscr O_{X,x}$ too. 
The map $\mathscr L\to \mathscr N$ being injective at $t$, the sequence
$0\to \mathscr L_t\to \mathscr N_t\to \mathscr P_t\to 0$ is exact. Since both $\mathscr P_t$ and $\mathscr L_t$ 
are flat over $\mathscr O_{X,x}$, an easy Tor computation shows that $\mathscr N_t$ is flat over $\mathscr O_{X,x}$; \ie, 
$\mathscr N$ is naively $X$-flat at $t$, and (iii) holds.

We can now prove (1). 
We may assume that $T$ and $X$ are affinoid. Let us suppose that (i) holds. 
Let $L$ be a complete extension of $k$, let $X'$ be an $L$-affinoid space, let $t'$ be a pre-image of $t$ on $T':=X'\times_XT$,
and let $x'$ be the image of $t'$ in $X'$. It suffices to show that the
map $\mathscr L_{T'}\to \mathscr N_{T'}$ is injective at $t'$ (this will yield the injectivity  of $\mathscr L\to \mathscr N$ at $t$ by taking $X'=X$)
and that $\mathscr P_{T'}$ is naively $X'$-flat at $t'$.
It suffices
to show
both properties after enlarging the field $L$
(Proposition \ref{prop-flat-xl} and Proposition \ref{prop-concrete-extension} (5)),
which allows
us
to assume that $t'$ is $L$-rational; in particular, $t'$ belongs to $\mathrm{Int}(T'/X')$. Since $\mathscr N$ is universally $X$-flat at $t$, the coherent sheaf
$\mathscr N_{T'}$ is $X'$-flat at $t'$. Moreover, by Remark \ref{rem-equiv-univinj}, the data $(T',X',t',x',\mathscr L_{T'}\to \mathscr N_{T'})$ fulfill the assumptions
of our proposition (which are the same as those
of Proposition \ref{prop-equiv-univinj}). Hence we can apply the assertion (2) already proven; it ensures that $\mathscr L_{T'}\to \mathscr N_{T'}$ is injective
at $t'$, and that $\mathscr P_{T'}$ is naively $X'$-flat at $t'$, whence (ii). 

Let us now suppose that (ii) holds and prove (i). Let $L$ be a complete extension of $k$, let $X'$ be an $L$-affinoid space, let $t'$ be a pre-image of $t$ on $T':=X'\times_XT$
and let $x'$ be the image of $t'$ in $X'$. It suffices to show that the
map $\mathscr N_{T'}$ is naively $X'$-flat at $t'$.
It is harmless to enlarge the field $L$
(Lemma \ref{lem-flatbc}), which allows
to assume that $t'$ is $L$-rational; in particular, $t'$ belongs to $\mathrm{Int}(T'/X')$. Since $\mathscr L\to \mathscr N$ is injective at $t$, it follows from Proposition \ref{prop-equiv-univinj}
that it is in fact {\em $X$-universally}
injective at $t$; hence $\mathscr L_{T'}\to \mathscr N_{T'}$ is injective at $t'$. 
Since $\mathscr P$ is universally $X$-flat at $t$, the coherent sheaf $\mathscr P_{T'}$ is naively $X'$-flat at $t'$. Moreover, 
by Remark \ref{rem-equiv-univinj}, the data $(T',X',t',x',\mathscr L_{T'}\to \mathscr N_{T'})$ fulfill the assumptions
of our proposition (which are the same as those
of Proposition \ref{prop-equiv-univinj}).
Hence we can apply the assertion (2) already proven; it ensures that $\mathscr N_{T'}$ is naively $X'$-flat at $t'$, whence (i). 
\end{proof}

\section{D\'evissages: definition and existence}

In this section we will make much use of the notions
of dimension, depth and codepth of a finitely generated
module over a local noetherian ring
(\ref{ss-conv-alg}, \ref{ex-algprop-coh}). We shall also need
the notions of dimension and codepth of a coherent sheaf at a given
point of the ambient space; dimension in this setting is defined
in \ref{ss-interpret-support}, 
and we are now going to define codepth. 

\begin{defi}\label{def-codepth-cohsheaf}
Let $Y$ be an analytic space, 
let $y$ be a point of $Y$, and let $\mathscr F$ be a coherent sheaf on $Y$. 
In view of Lemma-Definition \ref{lem-equiv-valid}
and Example \ref{ex-hreg-modules}, 
for a good analytic domain $V$
of $Y$ containing $y$,
the codepth of the $\mathscr O_{V,y}$-module $\mathscr F_{V,y}$ only depends
on $y$, and not on $V$ (Lemma-Definition \ref{lem-equiv-valid}
together with Example \ref{ex-hreg-modules}). It is called
the \emph{codepth} of $\mathscr F$ at $y$,\index{codepth of a coherent sheaf at a point}
and we denote it by $\mathrm{codepth}_y \;\mathscr F$. 
\end{defi}

\subsection{}\label{ss-codepth-basics}
Let $Y$ be an analytic space, let $y$ be a point of $Y$
and let $\mathscr F$
be a coherent sheaf on $Y$.

\begin{enumerate}[1]

\item Since the codepth of the zero-module is equal to zero by convention, 
we have $\mathrm{codepth}_y\;\mathscr F=0$ as soon as $y\notin \supp F$.

\item If $y$ belongs to $\supp F$, then $\mathrm{codepth}_y\;\dim_y \mathscr F$. 
Indeed, one immediately reduces to the good case, and one can then write
\[\mathrm{codepth}_y \mathscr F=\mathrm{codepth}_{\mathscr O_{Y,y}}\mathscr F_y
\leq \dim_{\mathrm{Krull}}\mathscr F_y\leq \dim_y \mathscr F,\]
where the last inequality comes from Corollary  \ref{cor-interp-centdim}
In particular we see that
$\mathrm{codepth}_y \;\mathscr F=0$ if $\dim_y Y=0$.
\item If $Y$ is regular at $y$, then $\mathrm{codepth}_y\;\mathscr O_Y=0$: this follows
from the fact that any regular local ring has codepth zero. 

\end{enumerate}

\begin{enonce}[remark]{Convention}\label{not-before-deviss}
Let $Y$ be an analytic space and let~$\mathscr F$ be a coherent sheaf on~$Y$. 
The unique coherent sheaf on~$\supp F$ that induces~$\mathscr F$ will be also denoted by~$\mathscr F$, if there is no risk of confusion.
\end{enonce}

\begin{defi}\label{def-devissages}\index{d\'evissage}
Let~$Y\to X$ be a morphism between
good~$k$-analytic spaces. Let~$\mathscr F$ be a coherent sheaf on~$Y$, let~$y$ be a point of  $\supp F$
and let~$x$ be its image on~$X$.  Let~$r$ be a positive integer, and let~$n_{1}>n_{2}>\ldots >n_{r}$ be a decreasing sequence of non-negative
integers.  A {\em~$\Gamma$-strict~$X$-d\'evissage of~$\mathscr F$ at~$y$ in dimensions~$n_{1},\ldots,n_{r}$ }
is a list of data~$(V,\{T_{i},\pi_{i},t_{i},{\mathscr L}_{i},{\mathscr P}_{i}\}_{i\in\{1,\ldots,r\}})$, where:

\begin{itemize}[label=$\bullet$]

\item~$V$ is a~$\Gamma$-strict affinoid neighborhood of~$y$ in~$Y$;

\item~$T_{i}$ is for every
$i$ a~$\Gamma$-strict~$k$-affinoid domain of a smooth~$X$-space  of pure relative dimension~$n_{i}$ and~$t_{i}$ is a point of~$T_{i}$ lying over~$x$;

\item for every~$i$,~${\mathscr L}_i$ and~${\mathscr P}_i$ are coherent~${\mathscr O}_{T_{i}}$-modules and~${\mathscr L}_i$ is free;

\item $t_i\in \mathrm{Supp}(\mathscr P_i)$ if~$i<r$, and~$\mathscr P_r=0$;

\item $\pi_1$ is a 
finite~$X$-map from $\mathrm{Supp}(\mathscr F_V)$ to $T_1$ such that we have~$\pi_1\inv (t_1)=\{y\}$ {\em set-theoretically};

\item $\pi_{i}$ for any~$i\in\{2,\ldots,r\}$ is
a finite~$X$-map from~$\mathrm{Supp}(\mathscr P_{i-1})$
to $T_{i}$
such that we have~$\pi_{i}\inv (t_{i}) =\{t_{i-1}\}$ {\em set-theoretically};

\item $\mathscr L_1$ is endowed with a morphism~${\mathscr L}_{1}\to \pi_{1\ast}{\mathscr F_{V}}$ whose cokernel is~${\mathscr P}_{1}$ and such that~$t_1\in \adht{\bij{ {(\mathscr L_{1})}_{T_{1,x}}}{(\pi_{1\ast} {\mathscr F}_{V})_{T_{1,x}}\;}}{T_{1,x}}$;

\item for any~$i\in\{2,\ldots,r\}$,~$\mathscr L_i$ is endowed with a morphism~${\mathscr L}_{i}\to \pi_{i\ast}{\mathscr P}_{i-1}$ whose cokernel is~${\mathscr P}_{i}$ and such that~$t_i\in \adht {\bij{ ({\mathscr L_{i}})_{T_{i,x}}}{(\pi_{i\ast}{\mathscr P}_{i-1})_{T_{i,x}}}}{T_{i,x}}\;$.
\end{itemize}
If we do not care about $\Gamma$ (\ie, if we take $\Gamma=\R_+\gpm$), we shall
simply say \emph{$X$-d\'evissage}. 
\end{defi}

The following commutative diagram of pointed spaces will hopefully make things easier to understand; at the beginning of every line, we have put the corresponding exact sequence of coherent sheaves (they live on the space~$T_i$ that lies on the line). 

\newpage {\small~$$\rotatebox{270}{\xymatrix{&&&&&(\mathrm{Supp}(\mathscr F_V),y)\ar[d]^{\pi_1}\ar@{^{(}->}[r]&(V,y)\ar[ddddddd]&\\  \mathscr L_1\to \pi_{1\ast}\mathscr F_V\to \mathscr P_1\to 0&&&&(\supp {P_1},t_1)\ar@{^{(}->}[r]\ar[d]^{\pi_2}&(T_1,t_1)\ar[rdddddd]&&\\\mathscr L_2
\to \pi_{2\ast}\mathscr P_1\to \mathscr P_2\to 0\ar@{.}[dddd]&&&(\supp {P_2},t_2)\ar@{^{(}->}[r]\ar@{.}[dd]&(T_2,t_2)\ar[rrddddd]&&&\\&&&&&&&\\&\ar@{.}[dd]&&&&&\\&&&&&&&\\ \mathscr L_{r-1}\to \pi_{r-1,\ast}\mathscr P_{r-2}\to \mathscr P_{r-1}\to 0&(\supp {P_{r-1}},t_{r-1})\ar[d]^{\pi_r}\ar@{^{(}->}[r]&(T_{r-1},t_{r-1})\ar[rrrrd]&&&&&\\ \mathscr L_r\to \pi_{r\ast}\mathscr P_{r-1}\to \mathscr P_r=0&(T_r,t_r)\ar[rrrrr]&&&&&(X,x)}}$$}\newpage

\begin{theo}\label{thm-devissages}
Let~$Y\to X$ be a map between good~$k$-analytic spaces, let~$\mathscr F$ be a coherent sheaf on~$Y$, let~$y$ be a point of $\supp F$,
and let~$x$ be its image in~$X$. Assume that the germ~$(Y,y)$ is~$\Gamma$-strict. Let~$c=\mathrm{codepth}_y \; \mathscr F_{Y_x}$ and let~$n=\dim_y \mathscr F_{Y_x}$.
There exists a~$\Gamma$-strict~$X$-d\'evissage of~$\mathscr F$ at~$y$ in dimensions belonging to~$[n-c\; ;\;n]$ (recall that $c\leq n$ by \ref{ss-codepth-basics} (2)). 
\end{theo}

\begin{proof}
According to \Cor 4.7 of \cite{ducros2007}, there exist an affinoid neighborhood~$Z$ of~$y$ in~$\supp F$, an affinoid domain~$T$ of a smooth~$X$-analytic
space of pure relative dimension~$n$, and a finite map~$\pi \colon Z\to T$ through which~$Z\to X$ factorizes. Let us set~$t=\pi(y)$. We can first assume, by shrinking~$T$, that~$y$ is the only pre-image of~$t$.

By Example \ref{ex-smallest-boundaryless}, the smallest analytic domain of $(T,t)$ through which the map $(Z,y)\to (T,t)$
factorizes is the whole of $(T,t)$; since $(Y,y)$
is $\Gamma$-strict, so is $(Z,y)$
and Theorem \ref{thm-smallest-germ} (3) then ensures that $(T,t)$ is $\Gamma$-strict as well.
By shrinking again $T$ and $Z$, we can thus assume that $T$ is $\Gamma$-strict; the $k$-affinoid space $Z$ is then $\Gamma$-strict too
by \ref{ss-pullback-gstrict}
(see also Remark \ref{rem-gstrict-fini-transfer}). 

Let $V$ be any $\Gamma$-strict affinoid neighborhood of $y$ in $Y$ such that $V\cap \supp F\subset Z$.
By the very definition of the topology on $T$, every point
of $T$ has a basis of open neighborhoods that are
finite intersections of subsets of $T$
defined by inequalities of the form $\abs g\in I$, for $g$ an analytic function on $T$
and $I$ an interval of $\R_+$ open in $\R_+$. It follows that every point 
also has a basis of \emph{compact}
neighborhoods that are finite intersections of subsets of $T$
defined by inequalities of the form $\abs g\in I$, for $g$ an analytic function on $T$
and $I$ an compact interval of $\R_+$ with non-empty interior in $\R_+$; by density of $(\Gamma\cdot\abs{k\gpm})^\Q$
in $\R_+\gpm$,
we can even restrict to
such intervals with endpoints in $(\Gamma\cdot\abs{k\gpm})^\Q\cup\{0\}$.
By the above and by topological properness (and topological separatedness) of $T\to V$, there exists analytic functions $g_1,\ldots, g_m$ on $T$, and elements $r_1,\ldots, r_m,s_1,\ldots, s_m$
of $(\Gamma\cdot \abs {k\gpm})^\Q\cup\{0\}$ with $r_i<s_i$ for every $i$, such that the $\Gamma$-strict affinoid
domain $T'$ of $T$ defined by the inequalities 
$r_i\leq \abs {g_i}\leq s_i$ for $i=1,\ldots, m$ is a neighborhood of $t$ and satisfies the inclusion $\pi\inv(T')\subset V\cap \supp F$. For every $i$, choose a lifting $h_i$ in $\mathscr O_Y(V)$ of the pull-back of $g_i$ in $\mathscr O_Z(V\cap Z)=\mathscr O_Z(V\cap \supp F)$; let $V'$ be the $\Gamma$-strict affinoid domain of $V$ defined by the inequalities 
$r_i\leq \abs {h_i}\leq s_i$ for $i=1,\ldots, m$. We then have $\pi\inv(T')=V' \cap Z$.
Hence by replacing $T$ with $T'$, $Z$ with $V'\cap Z$
and $V$ with $V'$, we may assume that $Z=V\cap \supp F$. (Note that in the construction above,
we may also require that $T'$ be contained in any given neighborhood of $t$; this allows us if needed to modify $T$, $Z$ and $V$ so that $T$ becomes arbitrary small.) 
We set $\mathscr G=\pi_\ast(\mathscr F_V)$, and $\delta=\dim_{\mathrm{Krull}} \mathscr O_{Z_x,y}=\dim_{\mathrm {Krull}}\mathscr F_{Y_x,y}$. 

By finiteness of $Z_x\to T_x$, one has $\mathrm{centdim}(T_x,t)=\mathrm{centdim}(Z_x,y)$ (\ref{ss-basicprop-centdim}).  
Since~$\dim_y Z_x=\dim_t T_x=n$, Corollary~\ref{cor-interp-centdim}
then ensures that
\begin{equation}\setcounter{equation}{1}
\dim_{\mathrm{Krull}} \mathscr O_{T_x,t}=\dim_{\mathrm{Krull}} \mathscr O_{Z_x,y}=\delta.
\end{equation}

The support of $\mathscr G$ is equal to $\pi(Z)$. It follows thus from
\ref{ss-dimloc-morfin}
that
\begin{equation}\label{eq-diman-gtx}
\dim_t \mathscr G_{T_x}=\dim_t \pi(Z_x)=n\
\end{equation}
because $Z_x$ is of dimension $n$ at $y$.
The fiber~$T_x$ being a smooth $\hr x$-analytic space of pure dimension~$n$, 
equality (\ref{eq-diman-gtx})
implies that $(\supp G)_x$ contains the connected component of~$t$ in~$T_x$, and 
the 
annihilator of $\mathscr G_{T_x}$ is thus zero at $t$. Therefore
\begin{equation}\label{eq-dimkr-gtx}
\dim_{\mathrm {Krull}}\mathscr G_{T_x,t}=\dim_{\mathrm {Krull}}\mathscr O_{T_x,t}=\delta.
\end{equation}

As~$\pi\inv (t)=\{y\}$,
one has~$\mathscr G_{T_{x},t}=\mathscr F_{Y_{x},y}$ (Lemma \ref{lem-cohsheaf-finflat} (1)). It follows then from
\cite{ega41}, \Chp 0, \S  16.4.8 that 

\begin{equation}\label{eq-depth-F}
\mathrm{depth}_{\mathscr O_{T_{x},t}}\mathscr G_{T_{x},t}=\mathrm{depth}_{\mathscr O_{Y_{x},y}} \mathscr F_{Y_{x},y}.
\end{equation}

In view of (\ref{eq-dimkr-gtx})
this yields the equality 

\begin{equation}\label{eq-codepth-G}
\mathrm{codepth}_t \mathscr G_{T_x}=\mathrm{codepth}_y \mathscr F_{Y_x}=c.
\end{equation}
We now argue by induction on~$c$. 
Assume that~$c=0$. Then~$\mathscr G_{T_{x},t}$ is a finitely generated module of codepth~$0$ and of maximal Krull dimension over
the regular local ring~${\mathscr O}_{T_{x},t}$. It is thus free (\cite{ega41}, \Chp 0, 17.3.4). Let~$(f_{i})_{1\leq i\leq r}$ be a family of 
sections of~$\mathscr G$
over the affinoid space~$T$ such that~$(f_{i}(t))$ is a basis of~$\mathscr G_{\hr t}$; set~$\mathscr L=\mathscr O_{T}^{r}$ and consider the map~$\mathscr L\to \mathscr G$
that
sends~$(a_1,\ldots,a_r)$ to~$\sum a_i f_i$. By Nakayama's Lemma, this map is surjective at~$t$; moreover, its restriction to~$T_x$ is {\em bijective} at~$t$, because its stalk at~$t$ is a surjective map between free modules of the same finite rank over~${\mathscr O}_{T_x,t}$. Hence by suitably shrinking~$T$ (and all other data) we may assume that~$\mathscr L\to \mathscr G$ is surjective, and that its restriction to~$T_x$ is bijective. We get this way a~$\Gamma$-strict~$X$-d\'evissage of~${\mathscr F}$ at~$y$ in dimension~$n$. 

Suppose now that~$c>0$, and that the theorem has been proved in codepth~$<c$.
Choose a point $z\in T_x$
such that $z_x\al$ is the generic point
of the connected component of $T_x\al$ containing $t_x\al$.
Since the support of $\mathscr G_{T_x}$ contains the connected component of $t$, the vector space~$\mathscr G\al_{T_x\al,\kappa(z_x\al)}$ is of positive dimension; let us call it~$r$. Let~$(f_i)_{1\leq i\leq r}$ be a family of
sections of~$\mathscr G$
over~$T$ such that~$(f_i(z_x\al))_i$ is a basis of the~$\kappa(z_x\al)$-vector space~$\mathscr G\al_{T_x\al,z_x\al}$. Set~$\mathscr L=\mathscr O_T^r$ and consider the map~$\mathscr L\to \mathscr G$ that
sends~$(a_1,\ldots, a_r)$ to~$\sum a_if_i$. Since $\mathscr O_{T_x\al, z_x\al}=\kappa(z_x\al)$, the map $\mathscr L\al_{T_x\al,z_x\al}\to \mathscr G\al_{T_x\al, z_x\al}$ is an isomorphism; this implies that $z$ belongs to $\bij{\mathscr L_{T_x}}{\mathscr G_{T_x}}$. Being a Zariski-specialization of $z$, the point $t$
belongs to $\adht{ \bij{{\mathscr L}_{T_x}}{\mathscr G_{T_x}}}{T_x}$. 

The scheme~$T_x\al$ being regular, $\mathscr O_{T_x\al, t_x\al}$ is a domain whose fraction field is $\mathscr O_{T_x\al,z_x\al}$. As $\mathscr L$ is
a free $\mathscr O_T$-module, this implies that the map $\mathscr L\al_{T_x\al, t_x\al}\to \mathscr L\al_{T_x\al,z_x\al}$ is injective. Hence in the commutative diagram
$$\xymatrix{
{\mathscr L\al_{T_x\al, t_x\al} }\ar[r]\ar[d]&{\mathscr G\al_{T_x\al, t_x\al}}\ar[d]\\
{\mathscr L\al_{T_x\al, z_x\al} }\ar[r]&{\mathscr G\al_{T_x\al, z_x\al}}
}$$
the bottom horizontal arrow is an isomorphism and the left vertical arrow is injective. It follows that the top horizontal arrow 
$\mathscr L\al_{T_x\al, t_x\al}
\to \mathscr G\al_{T_x\al, t_x\al}$ is injective. The map $\mathscr L_{T_x,t} \to \mathscr G_{T_x,t}$ is
thus injective.
But it is not surjective. Indeed, if it were surjective it would be bijective
and the codepth of~$\mathscr G_{T_x,t}$
would then be equal to zero (the local ring ${\mathscr O}_{T_x,t}$ being regular), which would contradict the fact that ~$c>0$.

Set~$\mathscr P=\mathrm{Coker}\;({\mathscr L}\to \mathscr G)$. Since $\mathscr  L_{T_x,t} \to \mathscr G_{T_x,t}$
is not surjective, $t$ lies in~$\supp P_x$. 
Since $\supp P_{x}$ is obviously included in~$T_x\setminus\bij{\mathscr L_{|T_x}}{\mathscr G_{T_x}}$, we have
$\dim_t \mathscr P_{T_x}<n$. 
Hence there exists a global section $a$ on $T_x$ of the annihilator ideal of $\mathscr P_{T_x}$ whose zero-locus
contains no neighborhood of~$t$. Since the annihilator of~$\mathscr P_{T_x,t}$ contains the germ of $a$, 
it is non-zero, which
yields the inequality
\begin{equation}\label{eq-dimkr-ptx}
\dim_{\mathrm{Krull}} \mathscr P_{T_{x},t}<\delta.
\end{equation}

%We thus have~$\mathrm{depth}_{\mathscr O_{T_x,t}}\mathscr P_{T_{x},t}<\delta$;
%moreover, since~$r>0$ and since~$\mathscr O_{T_{x},t}$ is regular,~$\mathrm{depth}_{\mathscr O_{T_x,t}}\mathscr L_{T_{x},t}=\delta$.

By \Th 16.7 of \cite{matsumura1986}, the depth of any non-zero finitely generated $\mathscr O_{T_x,t}$-module $M$ is the 
smallest integer $i$ such that $\mathrm{Ext}_{\mathscr O_{T_x,t}}^i(\kappa(t_x), M)\neq 0$. Since $r>0$ and
$\mathscr O_{T_{x},t}$ is regular,~$\mathrm{depth}_{\mathscr O_{T_x,t}}\mathscr L_{T_{x},t}=\delta$; in particular, 
$\mathrm{Ext}_{\mathscr O_{T_x,t}}^i(\kappa(t_x),\mathscr L_{T_x})=0$ as soon as $i<\delta$. 
By considering the $\mathrm {Ext}_{\mathscr O_{T_x,t}}^\bullet(\kappa(t_x),\cdot)$
exact sequence associated with
$$0\to \mathscr L_{T_x,t} \to \mathscr G_{T_x,t}\to\mathscr P_{T_x,t}\to 0,$$
we deduce from the above that $\mathrm{Ext}^i_{\mathscr O_{T_x,t}}(\kappa(t_x),\mathscr G_{T_x,t})=\mathrm{Ext}^i_{\mathscr O_{T_x,t}}(\kappa(t_x),\mathscr P_{T_x,t})$
for every $i<\delta$. 
Now (\ref{eq-dimkr-ptx})
implies that $\mathrm{depth}_{\mathscr O_{T_x,t}}\mathscr P_{T_{x},t}<\delta$. 
Using again the caracterization of depth in terms of the $\mathrm{Ext}$ functors, it follows that

\begin{equation}\label{eq-depth-P}
\mathrm{depth}_{\mathscr O_{T_{x},t}}\mathscr P_{T_x,t}=\mathrm{depth}_{\mathscr O_{T_{x},t}}\mathscr G_{T_x,t}=\delta-c,
\end{equation}
where the second equality comes from (\ref{eq-codepth-G}). 
In view of (\ref{eq-dimkr-ptx}), we get the inequality
\begin{equation}\label{eq-codepth-hrec}
\mathrm{codepth}_{\mathscr O_{T_x,t}}\mathscr P_{T_x,t}<c.
\end{equation}
This allows
us to apply the induction hypothesis. It ensures that ${\mathscr P}$ admits a~$\Gamma$-strict~$X$-d\'evissage at~$t$,
in dimensions belonging to~$$I:=[\dim_t \mathscr P_{T_x}-\mathrm{codepth}_t\mathscr P_{T_x}\; ;\; \dim_t \mathscr P_{T_x}].$$ We are going to show
show that~$I\subset [n-c,n)$; by shrinking suitably~$V$,~$Z$, and~$T$, the d\'evissage
of~$\mathscr P$ together with~$V,T,\pi,\mathscr L,\mathscr P, \mathscr L\to \mathscr G=\pi_\ast\mathscr F$
will then provide a~$\Gamma$-strict d\'evissage of~$\mathscr F$ at~$y$ in dimensions belonging to~$[n-c\; ;\;n]$.

As~$\dim_t \mathscr P_{T_x}<n$, 
the interval~$I$ is strictly bounded above by~$n$. 
Let us now prove that it is bounded below by~$n-c$. One has

$$\begin{array}{clr}
&\dim_t \mathscr P_{T_x}-\mathrm{codepth}_t \mathscr P_{T_x}&\\
=&\dim_t \mathscr P_{T_x}-\dim_{\mathrm{Krull}} \mathscr O_{\supp P_x,t}+\mathrm{depth}_{\mathscr O_{T_x,t}}\mathscr P_{T_x,t}&\\
=&\mathrm{centdim}(\supp P_x,t)+\mathrm{depth}_{\mathscr O_{T_x,t}}\mathscr P_{T_{x},t}&\text{by Corollary \ref{cor-interp-centdim}}\\
=&\mathrm{centdim}(T_x,t)+\mathrm{depth}_{\mathscr O_{T_x,t}}\mathscr P_{T_x,t}&\text{by \ref{ss-basicprop-centdim}}\\
=&\mathrm{centdim}(Z_x,y)+\mathrm{depth}_{\mathscr O_{T_x,t}}\mathscr P_{T_{x},t}&\text{by \ref{ss-basicprop-centdim}}\\
=&\mathrm{centdim}(Z_x,y)+\mathrm{depth}_{\mathscr O_{T_x,t}}\mathscr G_{T_x,t}&\text{by (\ref{eq-depth-P})}\\
=&\mathrm{centdim}\;(Z_x,y)+\mathrm{depth}_{\mathscr O_{Y_x,y}}\mathscr F_{Y_x,y}&\text{by (\ref{eq-depth-F})}\\
=&\dim_y Z_x -\dim_{\mathrm{Krull}}{\mathscr O}_{Z_x,y}+\mathrm{depth}_{\mathscr O_{Y_x,y}}\mathscr F_{Y_x,y}&\text{by Corollary \ref{cor-interp-centdim}}\\
=&n-c.&
\end{array}
$$
\end{proof}

\section{Flatness can be checked naively
in the inner case}\markright{\thesection.~ FLATNESS IN THE INNER CASE}

Let~$Y\to X$ be a map between good~$k$-analytic spaces and let~$\mathscr F$ be a coherent module on~$Y$. 
We want to give some criteria (in terms of a d\'evissage at $y$) for~$\mathscr F$ to be~$X$-flat at a given point $y$ of $Y$, and to use them to show that in the boundaryless (and, more generally, overconvergent) case, naive~$X$-flatness at~$y$ is equivalent to~$X$-flatness at~$y$ (this means that it is preserved under
arbitrary good base change, including
ground field extension); this fact had already been proved by Berkovich, using a completely different method, in some unpublished work.

\begin{defi}\label{def-overconv}\index{coherent sheaf!overconvergent}\index{overconvergent coherent sheaf}
Let~$Y\to X$ be a map between $k$-analytic spaces and let
$\mathscr F$ be a coherent module on~$Y$. Let~$y$ be a point of~$Y$. We
shall say that~$\mathscr F$ is {\em~$X$-overconvergent at~$y$} if there exist an
analytic neighborhood~$W$ of~$y$ in 
$Y$, an $X$-isomorphism between~$W$ and an analytic domain of a {\em boundaryless}~$X$-space~$W'$,
and a coherent sheaf~$\mathscr G$ on~$W'$ such that~$\mathscr F_{W}\simeq \mathscr G_{W}$.
\end{defi}

\begin{rema}
If $y$ belongs to $\mathrm{Int}(Y/X)$, then~$\mathscr F$ is automatically $X$-overconvergent at~$y$. 
\end{rema}

\begin{rema}
We shall only use this notion when $Y$ and $X$ are good. In this case, $W$ and $W'$
can be chosen
to be affinoid (note that 
a boundaryless space over a good space is itself good {\em by definition}, see \cite{berkovich1993}, page 34). 
\end{rema}

\begin{theo}\label{thm-flat-naiveflat}
Let $Y\to X$ be a morphism between good $k$-analytic spaces, let $y$ be a point of $Y$, and let $\mathscr F$ be
a coherent sheaf at $Y$. 

\begin{enumerate}[1]

\item If $y\notin \supp F$, then $\mathscr F$ is $X$-flat at $y$.

\item Assume that $y$ belongs to $\supp F$ and let $(V,\{T_{i},\pi_{i},t_{i},{\mathscr L}_{i},{\mathscr P}_{i}\}_{i\in\{1,\ldots,r\}})$
be an $X$-d\'evissage of $\mathscr F$ at $y$ (such a d\'evissage always exists by Theorem \ref{thm-devissages}). The following
are equivalent. 

\begin{enumerate}[j]

\item The coherent sheaf $\mathscr F$ is~$X$-flat at~$y$. 

\item The arrow~$\mathscr L_1\to \pi_{1\ast}\mathscr F_V$ is injective at~$t_{1}$ and for every~$i\geq 2$, the arrow~$\mathscr L_i\to \pi_{i\ast}\mathscr P_{i-1}$ is injective at~$t_i$.

\end{enumerate}

\item Under the assumptions of \textnormal{(2)}, suppose moreover
that $\mathscr F$ is~$X$-overconvergent
at~$y$; \eg, $y\in \mathrm{Int}(Y/X)$. Assertions \textnormal{(i)} and \textnormal{(ii)}
are then equivalent to: 

\begin{enumerate}[j]\setcounter{enumii}{2}
\item The coherent sheaf $\mathscr F$ is naively~$X$-flat at~$y$.
\end{enumerate}
\end{enumerate}
\end{theo}

\begin{proof}
Assertion (1) is obvious (and was written down only for the sake of completeness).
By Lemma \ref{lem-cohsheaf-finflat}, 
naive $X$-flatness (\resp $X$-flatness) of~$\mathscr F$ at~$y$ is equivalent to that of~$\pi_{1\ast}\mathscr  F_V$ at~$t_1$ ; for the same reason, if~$i\leq r-1$, then
naive~$X$-flatness 
(\resp $X$-flatness) of~$\mathscr P_i$ at~$t_i$ is equivalent to that of~$\pi_{(i+1)\ast}\mathscr P_{i}$ at~$t_{i+1}$.
Hence the equivalence (i)$\iff$(iii), and the equivalence
(i)$\iff$(ii)$\iff$(iii) {\em when~$y\in \mathrm{Int}(Y/X)$},
follow from a repeated application of Proposition \ref{prop-pre-devisssage},
once one has remarked that since~$\mathscr P_{r}=0$, it is~$X$-flat at~$t_r$.

It remains to show that~(iii)$\Rightarrow$(i) under the assumption that~$\mathscr F$ is~$X$-overconvergent at~$y$. By shrinking
$Y$, we may assume that it is an affinoid domain of a (relatively) boundaryless $X$-analytic space $Y'$, and $\mathscr F=\mathscr G_Y$
for some coherent sheaf $\mathscr G$ on $Y'$. 
Assume that $\mathscr G_Y$ is
naively~$X$-flat at~$y$. This implies that $\mathscr G$ is naively $X$-flat at $y$
(\ref{flat-good-basics}). 
By the boundaryless case already established,~$\mathscr G$ is~$X$-flat at~$y$.
Therefore,~$\mathscr G_Y$ is~$X$-flat at~$y$ (\ref{ss-flat-andom}). 
\end{proof}

\begin{rema}\label{rem-devissage-flat-iso}
If properties (i) and (ii) are satisfied, it turns out that
the map~$\mathscr L_{r}\to \pi_{r\ast}\mathscr P_{r-1}$ (or~$\mathscr L_{1}\to \pi_{1\ast}\mathscr F_V$ when $r=1$) is {\em bijective at~$t_r$}, because it is injective by
(ii) and surjective since its cokernel~$\mathscr P_r$ is zero.
\end{rema}

Let us now give three easy (but important) consequences of Theorem \ref{thm-flat-naiveflat}.
The first one shows that checking flatness at a given point of a good analytic space does
not actually require to consider all possible base changes; the second one is an improvement 
of Theorem \ref{th-gagaflat-hard}; the third one explains how flatness
can be checked fiberwise in some cases (this is kind of an analytic analogue of 
\cite{sga1}, Exposé IV, \Cor 5.7, and our proof consists in reducing straightforwardly to the
latter result).

\begin{theo}\label{thm-check-flatness}
Let $Y\to X$ be a morphism of good $k$-analytic spaces, let $\mathscr F$
be a coherent sheaf on $Y$
and let $y$ be a point
of $Y$. Let $L$ be an analytic extension of $k$, and let $z$ be an $L$-rigid
point of $Y_L$ lying above $y$. If $\mathscr F_L$ is naively $X_L$-flat at $z$, then 
$\mathscr F$ is $X$-flat at $y$.
\end{theo}

\begin{proof}
Since $z$ is an $L$-rigid point, it belongs to $\mathrm{Int}(Y_L/X_L)$.
As $\mathscr F_L$ is naively $X_L$-flat
at $z$ by assumption, Theorem \ref{thm-flat-naiveflat}
ensures that $\mathscr F_L$ is $X_L$-flat at $z$; the coherent
sheaf $\mathscr F$ is then $X$-flat at $y$ by Proposition \ref{prop-flat-xl}. 
\end{proof}

\begin{theo}\label{thm-flat-improve}
Let~$Y\to X$ be a morphism between~$k$-affinoid spaces and
let~$Z$ be a closed analytic subspace of~$Y$ such that~$Z\to X$ is finite.
Let~$y$ be a point of $Z$ and let $\mathscr F$ be a coherent sheaf on~$Y$.
Assume that $\mathscr F\al$ is~$X\al$-flat at~$y\al$. The coherent sheaf $\mathscr F$
is then $X$-flat at $y$.
\end{theo}

\begin{proof}
Theorem \ref{th-gagaflat-hard} tells us that~$\mathscr F$ is
naively~$X$-flat at~$y$. As~$Z\to X$ is finite,
it is inner at $y$, hence $Y\to X$ is inner at~$y$ too
by \ref{ss-boundary-basics} (3); in view of Theorem \ref{thm-flat-naiveflat},~$\mathscr F$ is~$X$-flat at~$y$.
\end{proof}

\begin{theo}\label{th-flat-fiberwise}
Let $X$ be an analytic space, and let $Z\to Y$
be a morphism of $X$-analytic spaces. Let $z$
be a point of $Z$, and let $y$
and $x$ denote the images of $z$ on $Y$ and $X$, respectively; 
we assume that $Y$ is $X$-flat at $y$. 
Let $\mathscr F$ be a
coherent sheaf on $Z$.
The following are equivalent: 

\begin{enumerate}[i]
\item $\mathscr F$ is $Y$-flat at $z$; 
\item $\mathscr F$ is $X$-flat at $z$, and $\mathscr F_{Z_x}$ is $Y_x$-flat
at $z$.
\end{enumerate}
\end{theo}

\begin{proof}
By arguing G-locally on $X, Y$, and $Z$, we may and do
assume that all of them are good.
Proposition \ref{prop-flat-xl}
allows us to enlarge the ground field before proving the theorem; 
we thus can assume that $x,y$, and $z$ are rigid. In view of
Theorem \ref{thm-check-flatness}
and since
\[ \mathscr O_{Y_x,y}=\mathscr O_{Y,y}/\mathfrak m_x\mathscr O_{Y,y}\;\text{and}\;
\mathscr O_{Z_x,z}=\mathscr O_{Z,z}/\mathfrak m_x\mathscr O_{Z,z}\]
because $z$ is rigid, 
it suffices to prove that the following are equivalent:

\begin{enumerate}[i]
\setcounter{enumi}{2}
\item $\mathscr F_z$ is flat over $\mathscr O_{Y,y}$ ; 
\item $\mathscr F_z$ is  flat over $\mathscr O_{X,x}$
and $\mathscr F_z/\mathfrak m_x\mathscr F_z$
is flat over $\mathscr O_{Y,y}/\mathfrak m_x\mathscr O_{Y,y}$. 
\end{enumerate}

As
$Y$ is $X$-flat at $y$, it is naively $X$-flat at $y$, 
which means that $\mathscr O_{Y,y}$ is flat over $\mathscr O_{X,x}$. 
The equivalence (iii)$\iff$(iv)
then comes from a direct application 
of \Cor 5.9 of \cite{sga1}, Exposé IV
with
$A=\mathscr O_{X,x}$, $B=\mathscr O_{Y,y}$, $C=\mathscr O_{Z,z}$, and
$M=\mathscr F_z.$
\end{proof}

\section{The relative CM property}

In this section we shall use the CM property of a coherent sheaf at a given point of an
analytic space; this has to be understood in the sense of Lemma-Definition \ref{valid-at-concrete}. 
Note that being CM at a point amounts to being of codepth zero at it, in the sense of Definition
\ref{def-codepth-cohsheaf}. 

\begin{defi}\label{def-rel-cm}
Let~$Y\to X$ be a morphism of~$k$-analytic spaces, and
let~$\mathscr F$ be a coherent sheaf on~$Y$. Let~$y$ be a point of $Y$ and let~$x$ be its image on~$X$.
We say that~$\mathscr F$ is CM~{\em over $X$}
at~$y$ 
if~$\mathscr F_{Y_x}$ is CM at~$y$ and
$\mathscr F$ is~$X$-flat at~$y$. We say that~$Y$
is CM over~$X$ at~$y$ if $\mathscr O_Y$
is. We say ``CM over $X$" to
mean
``CM over $X$ at every point of $Y$". 

%If there is no ambiguity about the morphism~$Y\to X$ involved, 
%we shall sometimes write ``relatively CM"
%instead of ``CM~over~$X$". 
\end{defi}

\begin{rema}\label{rema-cm-andom}
Let $Y\to X$ be a morphism between $k$-analytic spaces and let $\mathscr F$
be a coherent sheaf on $Y$. Let $y$ be a point of $Y$, let $V$ be an analytic domain
of $Y$ containing $y$, and let $U$ be an analytic domain of $X$ containing the image of $V$. 
Then $\mathscr F$ is CM over $X$ at $y$ if and only if $\mathscr F_V$ is CM over $U$ at $y$:
this follows from the good behavior of flatness and of the CM property with respect to the
restriction to analytic domains (see \ref{ss-flat-andom} for flatness
and Remark \ref{rem-equiv-valid} for the CM property). 
\end{rema}

\begin{exem}\label{ex-cm-qf}
Let $Y\to X$ be a morphism of $k$-analytic spaces, let $y$ be a point of $Y$, and let $x$
be its image in $X$. Let $\mathscr F$
be a coherent sheaf on $Y$. Assume that $Y\to X$
is quasi-finite at $y$. We then have $\dim_y Y_x=0$, which implies that
$\mathscr F_{Y_x}$ is of codepth 0 at $y$ (\ref{ss-codepth-basics} (2)); therefore 
$\mathscr F$
is CM over $X$ at $y$ if and only if it is $X$-flat at $y$. 
\end{exem}

\begin{exem}\label{ex-cm-qsm}
Let $Y\to X$ be a morphism of $k$-analytic spaces. 
Let $y$ be a point of $Y$. If $Y$ is quasi-smooth over $X$ at $y$, then $\mathscr O_Y$
is CM at $y$: this follows from \ref{ss-codepth-basics} (3) together with Theorem \ref{thm-main-qsm}.
\end{exem}

Our purpose is now
to show by using d\'evissages that a coherent sheaf is relatively CM at a given point
if and only if it can be written around this point
as a ``combination"
of 
a relatively CM coherent sheaf as in Example \ref{ex-cm-qf}
(with $Y\to X$ even
\emph{finite} at $y$)
and of another one of the kind
described in Example \ref{ex-cm-qsm}.

\begin{lemm}\label{lem-local-cm}
Let~$Y\to T$ be a finite morphism and let~$T\to X$ be a quasi-smooth morphism. Let~$\mathscr F$ be a $T$-flat
coherent sheaf on~$Y$. The sheaf~$\mathscr F$ is
CM
over~$X$.
\end{lemm}

\begin{proof}
One can assume that~$Y,T$, and~$X$ are~$k$-affinoid. Since $T\to X$ is quasi-smooth, it is $X$-flat, whence the
$X$-flatness of~$\mathscr F$. Let~$x$ be a point of $X$. We are going to show that~$\mathscr F_{Y_x}$ is CM. Let~$y$ be a pre-image of $x$ in $Y$
and let~$t$ be its image on~$T_x$. The ring~$\mathscr O_{T_x,t}$ is regular by quasi-smoothness of~$T_x$,so
it is CM. The morphism~$\spec \mathscr O_{Y_x,y}\to \spec \mathscr O_{T_x,t}$ is finite, hence has zero-dimensional fibers; the finite module~$\mathscr F_{Y_x,y}$ is~$\mathscr O_{T_x,t}$-flat by~$T$-flatness of~$\mathscr F$. It follows then from \Prop 6.4.1 (ii)
of \cite{ega42} that~$\mathscr F_{Y_x,y}$ is CM. 
\end{proof}

\begin{theo}\label{thm-local-cm}
Let~$X$ be a good~$k$-analytic space, let~$Y$ be
a good~$X$-analytic space, and let~$y$ be a point of $Y$. Assume that
the germ $(Y,y)$ is~$\Gamma$-strict. 
Let~$n$ be an integer. Let~$\mathscr F$ be a coherent sheaf on~$Y$ such that~$y\in \supp F$
and let $U$ be the set of points of $Y$ at which $\mathscr F$
is CM over $X$. 
Assume that~$y\in U$ and $\dim_y \mathscr F_{Y_x}=n$. There exist: 

\begin{itemize}[label=$\bullet$]
\item  a~$\Gamma$-strict~$k$-affinoid neighborhood~$V$ of~$y$ in~$Y$
that is contained in $U$; 

\item a~$\Gamma$-strict~$k$-affinoid domain~$T$ of a smooth~$X$-space of pure relative dimension~$n$;
 
\item a finite~$X$-morphism~$\pi :\mathrm{Supp}(\mathscr F_V)\to T$  with respect to which~$\mathscr F_V$ is~$T$-flat.

\end{itemize}
\end{theo}

\begin{proof}
Let~$x$ be the image of~$y$ in~$X$. Since~$\mathscr F$
is
CM over~$X$
at~$y$, 
one has~$\mathrm{codepth}\mathscr F_{Y_x,y}=0$. 
By Theorem \ref{thm-devissages},
the coherent sheaf $\mathscr F$ admits a $\Gamma$-strict~$X$-d\'evissage~$(V,T,\pi,t,{\mathscr L},\mathscr P=
0)$ at~$y$ in dimension~$n$. As~$\mathscr F$ is~$X$-flat at~$y$,
Remark \ref{rem-devissage-flat-iso}
ensures
that~$\mathscr L\to \pi_\ast \mathscr F_V$ is bijective at~$t$. We can hence shrink the data so that~$\pi_\ast \mathscr F_V$ is a free~$\mathscr O_T$-module,
which implies $T$-flatness of $\mathscr F_V$ by Proposition \ref{prop-finite-flat}. 
And it follows from Lemma \ref{lem-local-cm}
that $V\subset U$. 
\end{proof}

\subsection{}\label{ss-cm-locus}
Let~$X$ be a~$k$-analytic space, let~$Y$ be an~$X$-analytic
space 
and let~$\mathscr F$ be a coherent sheaf on~$Y$. 
Let~$U$ be the set of points of $Y$ at which $\mathscr F$
is CM over $X$.

\begin{enumerate}[1]

\item It follows from Theorem \ref{thm-local-cm} that~$U$ is an open subset of~$Y$.
We shall see later that it is even {\em Zariski}-open (Theorem \ref{thm-constloc-main}).

\item Assume that $\mathscr F$ is $X$-flat, and let $x$ be a point of $X$. By our flatness assumption 
on $\mathscr F$, 
the intersection~$U\cap Y_x$ is the
CM~locus of~$\mathscr F_{Y_x}$. 
It is a Zariski-open subset of~$Y_x$ (Lemma \ref{locus-concrete} (2)), which is dense.
Indeed, to show
this one can assume that~$Y$ is affinoid. Now if~$\eta$ denotes
the generic point of an irreducible component of~$Y_x\al$, then~$\mathscr F\al_{Y_x\al,\eta}$ is CM because ~$\mathscr O_{Y_x\al,\eta}$ is artinian; this fact
together with
affinoid GAGA (Lemma \ref{gaga-concrete})
implies our claim. 

\end{enumerate}

\chapter{Quasi-finite multisections and images of maps}\label{IM}

A celebrated theorem by Raynaud asserts the following, in our language\footnote{
Raynaud's result is written in the rigid-analytic language, with the corresponding notion of
flatness; the consistency
with our notion will be established later; see 
Corollary \ref{cor-gagaflat-rigid}
}: if $\abs {k\gpm}\neq \{1\}$ and if $\phi \colon Y\to X$ is a flat
morphism
between 
strictly $k$-affinoid spaces, then $\phi(Y)$ is a strict analytic domain of $X$, 
\cf \cite{frg2},
\Cor  5.11.
In this chapter, we
slighlty generalize this result and prove the following (Theorem \ref{th-image-compact}):
if $\phi \colon Y\to X$ is map between $k$-affinoid spaces, and if $Y$ is the support of an $X$-flat coherent sheaf
$\mathscr F$,
then $\phi(Y)$ is an analytic domain of $X$, which is $\Gamma$-strict whenever $Y$ is $\Gamma$-strict
(as usual, $\Gamma$ denotes a subgroup of $\R_+\gpm$ which is non-trivial if $k$ is trivially valued). Moreover, 
our methods differ from Raynaud's (we replace the use of formal schemes by that of Temkin's graded reduction, 
and we do not perform any flattening), hence
we provide a new proof of Raynaud's theorem. 

Let us now roughly explain how we proceed. We first consider the case where $\phi$ is quasi-finite (Proposition
\ref{prop-image-cm}).
The coherent sheaf $\mathscr F$ is then CM over $X$, and Theorem
\ref{thm-local-cm}
enables us to reduce to the case where $\phi$ is a quasi-\'etale map. Then by arguing locally
and using Theorem \ref{thm-smallest-germ} (which ensures
the existence of a smallest analytic domain containing the image of a morphism 
of analytic germs, described through Temkin's reduction), we reduce to the situation
where $\phi$
is finite and \'etale, in which case $\phi(Y)$ is a union of connected components of $X$, and we are done. 

To handle the general case we first reduce, by performing
a ground field extension to $k_r$ for some suitable $k$-free
polyradius $r$ and using the corresponding Shilov section, to the case where $\abs{k\gpm}\neq \{1\}$ and
$Y$ is strict. We then prove (Theorem \ref{thm-multisections-global}) that under these assumptions, there exists a
strictly
$k$-affinoid space
$X'$, a quasi finite map $\psi \colon X'\to X$, and an $X$-morphism $\sigma\colon X'\to Y$
such that:

\begin{enumerate}[a]
\item $\mathscr F$ is CM over $X$ at every point of $\sigma(X')$; 

\item $\sigma^\ast\mathscr F$ is $X$-flat; 

\item $\psi(X')=\phi(Y)$. 
\end{enumerate}
Now assertion (c) together with the quasi-finite case already proven ensures that $\phi(Y)$ is a strictly $k$-analytic
domain of $X$. 

Let us say a few
words about the existence of a quasi-finite multisection
of $\phi$ satisfying (a), (b) and (c), which seems to us of independent interest. 
It is the analogue of a classical scheme-theoretic result (see \cite{ega44}, \Prop 19.2.9), but
it is more involved, because of boundary phenomena, and also because
in analytic geometry,
the Zariski topology of a fiber
is in general finer than the topology induced by the Zariski topology of the ambient
space. The core of our proof is a local construction (which can afterward easily be globalized, by compactness of $Y$); 
it is the object of an independent
theorem (Theorem \ref{thm-multisections-local}), whose proof uses the following ingredients:

\begin{itemize}[label=$\bullet$]

\item Once again, the
existence of a smallest analytic domain containing the image of a morphism 
of analytic germs and its description through Temkin's reduction
(Theorem \ref{thm-smallest-germ}). 
In fact when the source germ is strict (which is the case here, by strictness of $Y$),
one can get a more precise description, involving finitely many closed points
of some ``residue scheme" attached to the situation
(Theorem
\ref{theo-elim-sections}), and this is absolutely crucial for our purposes:  
these closed points
precisely indicate in which ``directions" one has to draw multisections if one wants
to be sure that they will cover the whole image of our germ. 

\item The local structure of coherent sheaves CM over the ground space (Theorem \ref{thm-local-cm}). 

\item The ``quasi-finite version of Raynaud's theorem", already proven (this is
the aforementioned Proposition \ref{prop-image-cm}). 

\item The fact that over a non-trivially valued
field, a smooth morphism
admits \'etale multisections locally on its image
(Corollary \ref{cor-qsm-sections}). 
\end{itemize}

We end the chapter by recording
the following
extra-results about the images of maps, which we deduce
from (our version of) Raynaud's theorem; other ingredients are the
coincidence of the
topological and analytic interiors for the inclusion of an analytic domain
(this is used for (1), and for deducing (2b) from (2a)), and our local description of 
morphisms of relative dimension $d$ (\Cor 4.7
of \cite{ducros2007}) for (2a). 

\begin{enumerate}[1]
\item Let $\phi \colon Y\to X$ be a boundaryless morphism. If $Y$ is the support of an $X$-flat sheaf, then 
$\phi$ is open
(Theorem \ref{th-image-compact-inner}; this had been proved by Berkovich in some unpublished work).

\item Let $n$ and $d$ be two non-negative integers, and let $\phi \colon Y\to X$ be a morphism between analytic spaces. 
Assume that $X$ is normal and purely $n$-dimensional, $\phi$ is of pure relative dimension $d$, and $Y$ is purely
$(n+d)$-dimensional. Then:

\begin{enumerate}[b]
\item If $Y$ and $X$ are affinoid, $\phi(Y)$ is a compact analytic domain of $X$, which is $\Gamma$-strict whenever
$Y$ is $\Gamma$-strict (Theorem \ref{th-image-equi}). 

\item If $\phi$ is boundaryless, it is open. 
\end{enumerate} 
\end{enumerate}

\section{Flat, quasi-finite multisections
of flat maps}\label{s-image-cm}\markright{FLAT, QUASI-FINITE MULTISECTIONS}

\begin{prop}\label{prop-image-cm}
Let~$Y$ be a~$\Gamma$-strict
quasi-compact
$k$-analytic space, let~$X$ be a separated~$k$-analytic space,
and let~$\phi\colon Y\to X$ be a morphism. Let~$\mathscr F$ be a coherent sheaf on~$Y$ which is
$X$-flat 
and whose support is quasi-finite
over~$X$.
The image~$\phi(\supp F)$ is a~$\Gamma$-strict compact analytic domain of~$X$.
\end{prop}

\begin{proof}
We can replace~$Y$ with the support of~$\mathscr F$; \ie, we can assume that
the support of $\mathscr F$ is equal to $Y$.
We first make some reductions, which are allowed because
we can argue G-locally on $Y$ since the latter is quasi-compact. 

\begin{enumerate}[a]

\item By Proposition
\ref{prop-image-subgstrict}, $\phi(Y)$ is contained in a compact $\Gamma$-strict analytic domain
of $X$; hence
we can assume that both $Y$ and $X$ are $\Gamma$-strict and affinoid.

\item Since $\mathscr F$
is $X$-flat with quasi-finite support, it 
is CM over $X$ (Example \ref{ex-cm-qf}).
By Theorem \ref{thm-local-cm}, 
we can thus assume that there exist a~$\Gamma$-strict~$k$-affinoid quasi-\'etale~$X$-space~$T$, and a factorization of~$Y\to X$ through a finite map~$\pi\colon Y\to T$ such that~$\pi_\ast  \mathscr F$ is a free~$\mathscr O_T$-module of positive rank. The latter condition implies that~$\pi(Y)=T$. Replacing~$Y$ with~$T$, we can assume that~$Y\to X$ is quasi-\'etale. 
\item By Theorem \ref{thm-qsm-embeddsm},
we can suppose that the quasi-\'etale map~$Y\to X$ can be written as a composition~$Y\hookrightarrow X'\to X$ where~$Y\hookrightarrow X'$ identifies~$Y$ with a~$\Gamma$-strict
affinoid domain of~$X'$, where~$X'$ is connected and where~$X'\to X$
factorizes through a finite \'etale map from~$X'$ to a connected~$\Gamma$-strict affinoid domain~$Z$ of~$X$. Let~$X''$ be
a connected finite Galois covering of~$Z$ dominating
$X'$. One can replace~$X$ by~$Z$ and~$Y$ by its preimage on~$X''$; the union of all Galois conjugates of~$Y$ is
then a~$\Gamma$-strict compact analytic domain of~$X''$ (possibly non-affinoid) whose image on~$X$ coincides with that of~$Y$. 

\end{enumerate}
We thus
now assume that~$X$ is a~$\Gamma$-strict connected~$k$-affinoid space and that~$Y$ is a (possibly non-affinoid)
Galois-invariant~$\Gamma$-strict compact analytic domain of a finite connected
Galois cover~$X''$ of~$X$.
Let~$y$
be a point of $Y$ and let~$x=\phi(y)$. Let $G$ be the set-theoretic stabilizer of $y$ inside $\mathrm{Gal}(X''/X)$ ; since $Y$ is Galois invariant, 
$G$ stabilizes the germ $(Y,y)$. 
By usual Galois theory, $\hr y$ is a Galois extension of $\hr x$ with group $G$. It follows then from
\ref{ss-gradval-res} that 
$\hrt y^\Gamma$ is a normal graded extension of $\hrt x^\Gamma$ (\ref{ss-graded-galois})
 and the natural map $G\to \mathrm{Gal}(\hrt y^\Gamma/\hrt x^\Gamma)$ is surjective. 
In view of \ref{ss-zr-normal}, 
the continuous map $\P_{\hrt y^\Gamma/\widetilde k^\Gamma}\to \P_{\hrt x^\Gamma/\widetilde k^\Gamma}$ identifies
topologically
$\P_{\hrt x^\Gamma/\widetilde k^\Gamma}$ with $\P_{\hrt y^\Gamma/\widetilde k^\Gamma}/ \mathrm{Gal}(\hrt y^\Gamma/\hrt x^\Gamma)$. 

Let~$(U,x)$ be the smallest analytic domain of
the germ $(X,x)$ through which
the map $(Y,y)\to (X,x)$ factorizes; it
exists by Theorem \ref{thm-smallest-germ}, which
also ensures that
$(U,x)$ is~$\Gamma$-strict
and its reduction~$\widetilde{(U,x)} ^\Gamma$ is the image of~$\widetilde{(Y,y)} ^\Gamma$ on~$\P_{\hrt x ^\Gamma/\widetilde k ^\Gamma}$. As~$G$ stabilizes $(Y,y)$ and surjects onto $\mathrm{Gal}(\hrt y^\Gamma/\hrt x^\Gamma)$,
the group $\mathrm{Gal}(\hrt y^\Gamma/\hrt x^\Gamma)$
stabilizes $\widetilde{(Y,y)}^\Gamma$; therefore the pre-image of  $\widetilde{(U,x)} ^\Gamma$ inside~$\P_{\hrt y ^\Gamma/\widetilde k ^\Gamma}$ is precisely~$\widetilde{(Y,y)} ^\Gamma$. As a consequence,
the map~$(Y,y)\to (U,x)$ is boundaryless by
the criterion \ref{ss-prop-gammared} (2). Being quasi-\'etale, it is then \'etale by Remark \ref{rem-qsm-goodness}.
Therefore, there exist a~$\Gamma$-strict compact analytic neighborhood~$V$ of~$x$ in~$U$ and a compact analytic neighborhood~$W$ of~$y$ in~$Y$ such that~$\phi$ induces a finite  \'etale map~$W\to V$ (note that $W$ is $\Gamma$-strict by \ref{ss-pullback-gstrict}, but we do not need
this). The image $\phi(W)$ is now
a finite union of connected components of~$V$; in particular, it is a~$\Gamma$-strict compact analytic domain of~$X$. This ends the the proof due to
the compactness of $Y$.
\end{proof}

We are now going to state and prove our results on the existence of flat, quasi-finite multisections
(for maps whose source space is the support of a coherent sheaf flat over the target). 

\begin{theo}[Existence of flat, quasi-finite multisections: the local case]\label{thm-multisections-local}
Assume that~$\abs{k\gpm}\neq \{1\}$, and let~$\phi\colon Y\to X$ be a morphism of~$k$-analytic spaces, with~$Y$
strict and~$X$ separated. 
Let $\mathscr F$ be an $X$-flat coherent sheaf on~$Y$; let $y$ be a point of $\supp F$
and let $x$ be its image in $X$. 
Denote by~$Z$ the 
set of points of $Y$ at which $\mathscr F$
is CM over $X$ (this is an open subset of $Y$ by \ref{ss-cm-locus})
and by~$(U,x)$ the smallest analytic domain of~$(X,x)$ through which~$(\supp F, y)\to (X,x)$
factorizes
(see Theorem \ref{thm-smallest-germ}).

There exist~$r\geq 1$, quasi-finite maps~$$\psi_1 \colon X_1\to X,\ldots,\psi_r \colon X_r\to X,~$$ and~$X$-morphisms~$$\sigma_1 \colon X_1\to Z\cap \supp F, \ldots, \sigma_r\colon X_r\to Z\cap \supp F$$ such that: 

\begin{enumerate}[1]

\item For every~$j$, the space~$X_j$ is compact and strictly~$k$-affinoid,
and the point~$x$ has a unique pre-image~$x_j$ on~$X_j$.  

\item For every~$j$, the coherent sheaf ~$\sigma_j^\ast  \mathscr F$ is
$X$-flat, and~$\psi_j(X_j)$ is thus a compact strictly~$k$-analytic domain of~$X$ (by Proposition \ref{prop-image-cm}). 

\item One has~$(U,x)=\bigcup (\psi_j(X_j),x)$. 

\end{enumerate}

Moreover:

\begin{enumerate}[A]

\item If~$Y\to X$ is quasi-smooth at~$y$ and if~$\mathscr F=\mathscr O_Y$, the~$\psi_j$'s can be chosen to be quasi-\'etale.  

\item If~$Y\to X$ is boundaryless at~$y$ (which implies that~$(U,x)=(X,x)$; see Example \ref{ex-smallest-boundaryless}) and if the germs~$(Y,y)$ and~$(X,x)$ are good, one can take~$r=1$, and~$\psi_1$ inner, hence finite,  at~$x_1$.
\end{enumerate}
\end{theo}

\begin{proof}
By replacing~$Y$ with~$\supp F$, we may assume that~$Y=\supp F$. 
We are first going to reduce all assertions to the case where both $Y$ and $X$ are strictly affinoid, 
by arguing locally or G-locally (hence we shall implicitly use the good behavior of flatness and of the ``relative"
`CM property with respect to restriction to analytic domains; see \ref{ss-flat-andom} 
and Remark \ref{rema-cm-andom}). 

Let us begin with assertion (B). So, we assume that $Y$ and~$X$ are good and $Y\to X$ is inner at~$y$. As noted in the statement  of the theorem, this
implies in view of  Example \ref{ex-smallest-boundaryless}
that~$(U,x)=(X,x)$; strictness of~$(Y,y)$ then implies that of~$(X,x)$ by Theorem \ref{thm-smallest-germ}.
We can
thus shrink~$Y$ and~$X$ so that both are strictly~$k$-affinoid. 

Let us now come to the other assertions (so, we do not assume anymore that
$Y$ and~$X$ are good nor
that $\phi$ is inner at~$y$). By replacing~$Y$ with a
strictly analytic compact neighborhood of~$y$,
one can assume that it is compact. Now, as~$X$ is separated,
by Proposition \ref{prop-image-subgstrict}
$\phi(Y)$ is contained
in a compact, strictly analytic domain~$X_0$ of~$X$; as~$(U,x)\subset (X_0,x)$ we can replace~$X$ with~$X_0$, and hence reduce to the case where~$X$ itself is strict.
The point~$x$ thus
has
a neighborhood in~$X$
that is a finite union of strictly affinoid domains containing~$x$; all assertions involved
are G-local on the germ~$(X,x)$, so we
can assume that~$X$
itself  is strictly~$k$-affinoid. The point~$y$ has a neighborhood in~$Y$ that is a finite union of strictly affinoid domains containing it; all assertions  involved being G-local on the germ~$(Y,y)$, one can assume that~$Y$ itself is strictly~$k$-affinoid.

\begin{enonce*}{Local convention}
As the proof will involve from now on only strictly~$k$-analytic spaces, it will be sufficient to consider \emph{non-graded} reductions; therefore,
in order to simplify notation, we shall write
\emph{for the rest pf the proof}~$\widetilde{(Y,y)},\widetilde k$, \etc, to denote~$\widetilde{(Y,y)}^1,\widetilde k^1$, \etc 
\end{enonce*}

\subsubsection{}
Let~$A$ (resp.~$B$) be the algebra of analytic functions on~$X$ (resp.~$Y$). Let
$A^\circ$  be the subring of~$A$ that consists
of functions whose spectral semi-norm is bounded above by~$1$,
let~$A^{\circ \circ}$ be the ideal of~$A^\circ$ that consists
of functions whose spectral semi-norm is {\em strictly} bounded
above by~$1$, and let~$\widetilde A$ be the quotient~$A^\circ/A^{\circ \circ}$; we define~$B^\circ, B^{\circ \circ}$ and $\widetilde B$ analogously; both~$\widetilde k$-algebras~$\widetilde A$ and~$\widetilde B$ are finitely generated
(\cite{bosch-g-r1984}, 6.3.4, \Cor 3). 
We denote by~$A'$ (\resp $B'$) the image of
the natural evaluation map~$\widetilde A\to \hrt x$
(\resp $\widetilde B \to \hrt y$); we denote by $B''$ the subring of $\hrt y$ generated by $\hrt x$
and $B'$. 
By Temkin's definition of the (non-graded) reduction of a strict analytic germ
(\cite{temkin2000}; see also
Remark \ref{rem-directdef-gred}), one has
\setcounter{equation}{0}
\[
\widetilde{(X,x)}=\P_{\hrt x/\widetilde k}\{A'\}\;\text{and}\;\widetilde{(Y,y)}=\P_{\hrt y/\widetilde k}\{B'\}.
\]
Let~$f_1,\ldots,f_n$ be elements of~$B^\circ$ whose images generate the~$\widetilde k$-algebra~$\widetilde B$. The~$\widetilde k$-algebra~$B'$
is then generated by~$\widetilde{f_1(y)}, \ldots, \widetilde{f_n(y)}$, whence the equality
\[
\widetilde{(Y,y)}=\P_{\hrt y/\widetilde k}\{\widetilde{f_1(y)},\ldots,\widetilde{f_n(y)}\}.
\]

\emph{The inner case.} If~$\phi$ is inner at~$y$, then~$\widetilde{(Y,y)}$ is equal to the pre-image of~$\widetilde{(X,x)}$ in~$\P_{\hrt y/\widetilde k}$; in other words,~$\P_{\hrt y/\widetilde k}\{B'\}=\P_{\hrt y/\widetilde k}\{A'\}$, which implies that~$B'$ is integral over~$A'$; in this case, $B''$
is a fortiori
algebraic over~$\hrt x$, hence is a field. 

\subsubsection{}\label{ss-closed-points}
By Theorem~\ref{theo-elim-sections}, there exist finitely many closed points~$\mathsf y_1,\ldots,\mathsf y_m$ of~$\spec B''$ such that

\[
\widetilde{(U,x)}=\bigcup p_{j}\left(\P_{\kappa(\mathsf y_j)/\widetilde k}\{\widetilde{f_1(y)}(\mathsf y_j),\ldots,\widetilde{f_n(y)}(\mathsf y_j)\}\right),
\]
where~$p_j$ denotes
the natural continuous map~$\P_{\kappa(\mathsf y_j)/\widetilde k}\to \P_{\hrt x/\widetilde k}$
for every $j$. Set~$\mathsf U_j=p_{j}\left(\P_{\kappa(\mathsf y_j)/\widetilde k
}\{\widetilde{f_1(y)}(\mathsf y_j),\ldots,\widetilde{f_n(y)}(\mathsf y_j)\}\right)\subset \P_{\hrt x/\widetilde k}$ for every $j$; by
Proposition~\ref{pro-elim-simple}, $\mathsf U_j$ is open and quasi-compact. For every $j$, choose a compact strictly analytic domain $U_j$ of~$X$ that contains~$x$ and satisfies the equalities~$\widetilde{(U_j,x)}=\mathsf U_j$. Since $\widetilde{(U,x)}$ is the union of the $\mathsf U_j$'s, we have the equality
\[
(U,x)=\bigcup_{j} (U_j,x).
\]

\emph{The inner case.}
If~$\phi$ is inner at~$y$, then as~$B''$ is a field,~$m=1$ and~$\mathsf y_1$ is the only point of~$\spec B''$. It follows that~$\P_{\kappa(\mathsf y_1)/\widetilde{k}}\{\widetilde{f_1(y)}(\mathsf y_1),\ldots,\widetilde{f_n(y)}(\mathsf y_1)\}$ is nothing but~$\P_{B''/\widetilde k}\{\widetilde{f_1(y)},\ldots,\widetilde{f_n(y)}\}$.
But since~$\P_{\hrt y/\widetilde k}\{\widetilde{f_1(y)},\ldots,\widetilde{f_n(y)}\}$ is, under our innerness assumption, the pre-image of~$\widetilde{(X,x)}$ in~$\P_{\hrt y/\widetilde k}$, the open subset~$\P_{B''/\widetilde k}\{\widetilde{f_1(y)},\ldots,\widetilde{f_n(y)}\}$ of~$\P_{B''/\widetilde k}$ is the pre-image of~$\widetilde{(X,x)}$ in~$\P_{B''/\widetilde k}$. One thus has

\[
\P_{\kappa(\mathsf y_1)/\widetilde k}\{\widetilde{f_1(y)}(\mathsf y_1),\ldots,\widetilde{f_n(y)}(\mathsf y_1)\}=p_1\inv (\widetilde{(X,x)}).
\]

\subsubsection{}\label{sss-chose-j}
We fix an integer~$j$ belonging to~$\{1,\ldots,m\}$. Let~$R$ be the subring of~$\mathscr O_{X,x}$
consisting of functions~$f$ such that~$\abs{f(x)}\leq 1$. For any~$i\in \{1,\ldots,n\}$, let~$P_i$ be a polynomial
in~$R[T_1,\ldots,T_i]$
that is monic in~$T_i$ and
is such that~$\widetilde{P_i}(\widetilde{f_1(y)}(\mathsf y_j),\ldots,\widetilde{f_{i-1}(y)}(\mathsf y_j),T)$
is the minimal polynomial of~$\widetilde{f_i(y)}(\mathsf y_j)$ over~$\hrt x[\widetilde{f_1(y)}(\mathsf y_j),\ldots,\widetilde{f_{i-1}(y)}(\mathsf y_j)]$ (by~$\widetilde{P_i}$ we denote of course the image of~$P_i$ under the natural map~$R[T_1,\ldots,T_i]\to\hrt x[T_1,\ldots,T_i]$). Let~$D$ be
a strictly affinoid neighborhood of~$x$ in~$X$ on which all the coefficients of the~$P_i$'s are defined. Let~$\Omega$
be the open subset of~$Y\times_X D$ defined as the
locus of
simultaneous validity of
the
inequalities~$$\abs{P_1(f_1)}<1,\abs{P_2(f_1,f_2)}<1,\ldots,\abs{P_n(f_1,\ldots,f_n)}<1.$$ 

Let us prove
by contradiction that $\Omega_x\neq \emptyset$. Suppose that~$\Omega_x=\emptyset$. Let~$I$ be the subset of $\{1,\ldots,n\}$ consisting of integers $i$
such that~$\abs{P_i(f_1(y),\ldots,f_i(y))}=1$ (a priori, this absolute value is \emph{at most} 1).
For every~$i\in I$, let~$Y_i$ be the affinoid domain of~$Y_x$ defined by the condition~$\abs{P_i(f_1,\ldots,f_i)}=1$. Under our assumption
that $\Omega_x=\emptyset$,
the union of the~$Y_i$'s for~$i\in I$ is a neighborhood of~$y$ in~$Y_x$. We thus have~$\widetilde{(Y_x,y)}=\bigcup\limits_{ i\in I} \widetilde{(Y_i,y)}$. Let us describe both terms of this equality.

\begin{itemize}[label=$\bullet$]

\item By \ref{ss-prop-gammared}, one has $\widetilde{(Y_x,y)}=\P_{\hrt y/\hrt x}\{\widetilde{f_1(y)},\ldots,\widetilde{f_n(y)}\}$.

\item  If~$i$ is any element of~$I$, then~$\widetilde{(Y_i,y)}$ is equal to~$$\P_{\hrt   y/\hrt  x}\{\widetilde{f_1(y)},\ldots,\widetilde{f_n(y)},\widetilde{P_i}(\widetilde{f_1(y)},\ldots,\widetilde{f_i(y)}), \widetilde{P_i}(\widetilde{f_1(y)},\ldots,\widetilde{f_i(y)})\inv \}.$$ 
\end{itemize}

There exists a valuation~$\mathsf v$ on~$\hrt y$ that is trivial on~$\hrt x$ and
whose ring~$\mathscr O_{\mathsf v}$ dominates~$\mathscr O_{\spec B'',\mathsf y_j}$.
As~$\widetilde{f_1(y)},\ldots, \widetilde {f_n(y)}$ belong to~$B''$,
they belong to~$\mathscr O_{\mathsf v}$; as~$\widetilde{P_i}(\widetilde{f_1(y)}(\mathsf y_j),\ldots,\widetilde{f_i(y)}(\mathsf y_j))=0$ for all~$i$, the element~$\widetilde{P_i}(\widetilde{f_1(y)},\ldots,\widetilde{f_i(y)})$ belongs to the maximal ideal of~$\mathscr O_{\mathsf v}$ for all~$i$. 
It now follows from the
above explicit descriptions of~$\widetilde{(Y_x,y)}$ and of the~$\widetilde{(Y_i,y)}$'s
that~$\mathsf v$ belongs to~$\widetilde{(Y_x,y)}$ but not to~$\bigcup\limits_{ i\in I} \widetilde{(Y_i,y)}$, contradiction. 

\subsubsection{}\label{sss-core-proof}
As~$\Omega_x\neq \emptyset$, it follows from \ref{ss-cm-locus}
(2)
that there exists a point~$\omega$ in~$\Omega_x$ lying on~$Z$.
By Theorem \ref{thm-local-cm}, there exists a strictly affinoid
neighborhood~$V$ of~$\omega$ in~$\Omega\cap Z$ such that~$V\to X$
admits a factorization~$V\to T\to X$, where~$T$ is a strictly affinoid  domain of a
smooth~$X$-space~$S$
and $V\to T$ is a finite map with respect to which~$\mathscr F_V$ is~$T$-flat. By
Theorem \ref{thm-qsm-embeddsm}, if~$\phi$ is quasi-smooth at~$y$ and $\mathscr F=\mathscr O_Y$, one can suppose that~$V=T$.

Let~$\varpi$ be the image of~$\omega$ in~$T$. By Corollary \ref{cor-finflat-open},
the image of~$\mathrm{Int}(V/Y)$ in~$T$ contains an open neighborhood
$W$ of~$\varpi$. As~$W_x$ is a non-empty strictly~$\hr x$-analytic space, it has an~$\hr x$-rigid point
(\ref{ss-analytic-nullst}), which automatically belongs to~$\mathrm{Int}(T_x/S_x)$.
This implies the existence of an open subset of~$S$ whose fiber at~$x$ is non-empty and is included in~$W_x$. 
Since $\abs{k\gpm}\neq \{1\}$, applying Corollary \ref{cor-qsm-sections}
to this open subset provides an \'etale~$X$-space~$X'$ and an~$X$-morphism~$X'\to S$ whose image intersects~$W_x$.

We fix a pre-image~$x'$ of~$x$ in~$X'$ whose image
in~$S$ belongs to~$W$ and is denoted by~$t$. We choose a pre-image~$v$ of~$t$ in
$\mathrm{Int}(V/Y)$. We denote by~$T'$ the analytic domain~$T\times_S X'$  of~$X'$, and by~$V'$ the fiber product~$V\times_TT'$.
We choose a point~$v'\in V'$ lying above both~$v$ and~$x'$. We then have the following commutative diagram of pointed spaces, 

$$\xymatrix {&(V,v)\ar[d]\ar@{^{(}->}[r]&(Z,v)\ar@/^1.5pc/[dd]\\(V',v')\ar[d]\ar[ru]&(T,t)\ar@{^{(}->}[r]&(S,t)\ar[d]\\(T',x')\ar[ru]\ar@{^{(}->}[r]&(X',x')\ar[r]\ar[ru]&(X,x)}$$
in which both squares are cartesian. Since~$X'\to X$ is \'etale, it is boundaryless. As it factorizes
through $X'\to S$, the latter
map is
boundaryless as well; thus~ $T'\to T$ and~$V'\to V$ are also boundaryless. 

As~$v\in \mathrm{Int}(V/Y)$, the germ~$(V,v)$ coincides with~$(Y,v)$; in other words,
we have~$\widetilde{(V,v)}=\P_{\hrt  v/\widetilde k}\{\widetilde{f_1(v)},\ldots,\widetilde{f_n(v)}\}$.
Since~$V'\to V$ is boundaryless, $\widetilde{(V',v')}=\P_{\hrt  {v'}/\widetilde k}\{\widetilde{f_1(v')},\ldots,\widetilde{f_n(v')}\}$
(one still writes~$f_i$ for the pull-back of~$f_i$ into the ring of functions on~$V'$); note that this implies
that $V'$ is good. 

By
the choice of~$V$, the point~$v$ belongs to~$\Omega_x$.
We hence have for every~$i$ the inequality~$\abs{P_i(f_1,\ldots,f_i)(v)}<1$. It implies
that~$\widetilde{P_i}(\widetilde{f_1(v')},\ldots,\widetilde{f_i(v')})=0$ for every~$i$. 
By the very definition of the~$P_i$'s (resting on $j$ as chosen at the beginning of \ref{sss-chose-j}), it follows
that there exists an~$\hrt x$-isomorphism between~$\kappa(\mathsf y_j)$
and~$\hrt x[\widetilde{f_1(v')},\ldots,\widetilde{f_n(v')}]$
that sends~$\widetilde{f_i(y)}(\mathsf y_j)$ to~$\widetilde{f_i(v')}$ for any~$i$. 
The image of~$\widetilde{(V',v')}$ inside~$\P_{\hrt  x/\widetilde k}$
therefore
coincides
with~$p_j\left(\P_{\kappa(\mathsf y_j)/\widetilde k
}\{\widetilde{f_1(y)}(\mathsf y_j),\ldots,\widetilde{f_n(y)}(\mathsf y_j)\}\right)=\mathsf U_j$.
As a consequence,~$(U_j,x)$ is the smallest analytic domain of~$(X,x)$ through which~$(V',v')\to (X,x)$
factorizes.

The morphism~$V'\to X$ is quasi-finite. The space~$T'$ is quasi-\'etale, and in particular
quasi-smooth, over~$X$, and $\mathscr F_{V'}$ is flat over~$T'$
because~$\mathscr F_V$ is flat over~$T$.
Moreover, if~$Y\to X$ is quasi-smooth and $\mathscr F=\mathscr O_Y$
then~$V'$ is quasi-\'etale over~$X$ (because in this situation~$V=W$).

\emph{The inner case}.
If~$\phi$ is inner at~$y$ we have seen at the end of \ref{ss-closed-points}
that~$j=1$ and that
$$\P_{\kappa(\mathsf y_1)/\widetilde k}\{\widetilde{f_1(y)}(\mathsf y_1),\ldots,\widetilde{f_n(y)}(\mathsf y_1)\}=p_1\inv (\widetilde{(X,x)}).$$ Therefore~$\widetilde{(V',v')}$ is the pre-image of~$\widetilde{(X,x)}$ in~$\P_{\hrt {v'}/\widetilde k}$, which exactly means that~$(V',v')\to (X,x)$ is boundaryless. 

\subsubsection{Conclusion}
As~$V'\to X$ is quasi-finite and $V'$ is good, 
there exists a strictly~$k$-affinoid
neighborhood~$X_j$ of~$v'$ in~$V'$ such that~$v'$ is the only pre-image of~$x$ inside~$X_j$.
To emphasize the dependance on~$j$, let us denote now by~$x_j$ the point~$v'$, by~$\psi_j$ the natural map~$X_j\to X$, and by~$\sigma_j$ the natural~$X$-map~$X_j\to Z$. 

The following follow from what was done in \ref{sss-core-proof}. 

\begin{itemize}[label=$\bullet$]

\item The morphism $\psi_j$ is quasi-finite. 

\item The coherent sheaf~$\sigma_j^\ast  \mathscr F$ 
is $X$-flat. 

\item The smallest analytic domain of~$(X,x)$ through which~$(X_j,x_j)$ factorizes is~$(U_j,x)$.

\item If~$\mathscr F=\mathscr O_Y$ and $\phi$ is quasi-smooth at~$y$
then~$\psi_j$ is quasi-\'etale. 

\item If~$\phi$ is inner at~$y$ then~$j=1$ and~$\psi_1$ is inner, hence finite, at~$x_1$. 

\end{itemize}

As the coherent sheaf~$\sigma^\ast  _j\mathscr F$ is
$X$-flat 
(and has support~$X_j$ because~$\mathscr F$ has support~$Y$), Proposition~\ref{prop-image-cm} ensures that~$\psi_j(X_j)$ is an analytic domain of~$X$.
We can shrink~$X_j$ so that~$\psi_j(X_j)\subset U_j$; thus $(\psi_j(X_j), x)=(U_j,x)$ by minimality of~$(U_j,x)$. 
Since~$(U,x)$ is the union of the~$(U_j,x)$'s (\ref{ss-closed-points}), the data~$(X_j,\psi_j,\sigma_j)_j$ satisfy the conclusions of the theorem.
\end{proof}

\begin{theo}\label{thm-multisections-global}
Assume that~$\abs{k\gpm}\neq \{1\}$.
Let~$Y$ be a
quasi-compact, strictly~$k$-analytic space and let~$X$ be a separated~$k$-analytic space. Let~$\phi: Y\to X$ be a morphism and let~$\mathscr F$
be an~$X$-flat coherent sheaf on~$Y$. Denote by~$Z$ the 
set of points of $Y$ at which $\mathscr F$
is CM over $X$ (this is an open subset of $Y$ by \ref{ss-cm-locus}).
There exist
a strictly~$k$-affinoid
space~$X'$, a
quasi-finite map~$\psi \colon X'\to X$, and an~$X$-morphism~$\sigma \colon X'\to Z\cap \supp F$ such that the following hold: 

\begin{enumerate}[1]
\item The coherent sheaf $\sigma^\ast  \mathscr F$ is 
$X$-flat (so $\psi(X')$ is a compact strictly analytic domain of~$X$ by Proposition \ref{prop-image-cm}). 

\item  The image $\phi(Y)$ is equal to $\psi(X')$.
\end{enumerate}

If moreover~$Y\to X$ is quasi-smooth and~$\mathscr F=\mathscr O_Y$, then~$\psi$ can be chosen to be quasi-\'etale.
\end{theo}
\begin{proof}

By replacing~$Y$ with~$\supp F$ we may assume that~$\supp F=Y$. Let~$y$ 
be a point of $Y$. Using the notation of Theorem \ref{thm-multisections-local}, and setting
\[X^y=\coprod X_j, \;\psi^y=\coprod \psi_j,\sigma^y=\coprod \sigma_j,
\]
one gets the existence of a strictly~$k$-affinoid space~$X^y$,
a morphism
$\psi^y: X^y\to X$ which is quasi-finite and
even quasi-étale if $\phi$ is quasi-smooth
and $\mathscr F=\mathscr O_Y$, and an~$X$-map~$\sigma^y : X^y\to Z$
such that the following are satisfied: 

\begin{enumerate}[a]

\item The coherent sheaf $(\sigma^y)^\ast  \mathscr F$ is
$X$-flat,  so $\psi^y(X^y)$
is a compact strictly~$k$-analytic domain of~$X$ by Proposition  \ref{prop-image-cm}.

\item The germ $(\psi^y(X^y),x)$ is equal to the smallest analytic domain $(U^y,x)$ of~$(X,x)$ through which~$(Y,y)\to (X,x)$
factorizes.

\item  As~$(Y,y)\to (X,x)$ factorizes through~$(\psi^y (X^y),x)$, there exists an 
analytic neighborhood~$V^y$ of~$y$ in~$Y$ such that~$\phi(V^y)\subset \psi^y(X^y)$. 

\end{enumerate}

By quasi-compactness of~$Y$,
there exist finitely many points~$y_1,\ldots, y_n$ on~$Y$ such that
$Y=\bigcup_iV^{y_i}$. Now set~$X'=\coprod X^{y_i}, \psi =\coprod \psi^{y_i}$, and~$\sigma=\coprod \sigma^{y_i}.$
By construction, (1) is satisfied
and~$\psi$ is quasi-\'etale if~$\mathscr F=\mathscr  O_Y$ and~$\phi$ is quasi-smooth; it thus remains to show (2).  
For every~$i$, one has~$\phi(V^{y_i})\subset \psi^{y_i}(X^{y_i})$. As $Y=\bigcup_iV^{y_i}$,
this implies that~$\phi(Y)\subset \psi(X')$; but the existence of~$\sigma$ provides the reverse inclusion, whence (2). 
\end{proof}

\begin{rema}
The strictness assumption
on $Y$ cannot
be dropped from the statement 
of Theorem \ref{thm-multisections-global}
(even if one does not require anymore
in the conclusion that
$X'$ be strict).
Indeed, let $Y$ be a compact $k$-analytic space.  The coherent 
sheaf $\mathscr O_Y$ is $k$-flat (Lemma \ref{lem-field-flat}); hence
we can apply Theorem 
\ref{thm-multisections-global} with $X=\mathscr M(k)$ and $\mathscr F=\mathscr O_Y$; and
even if one drops the strictness requirement
on $X'$, it
then simply states that if $Y$ is strict and non-empty, it has a rigid point. 
This is nothing but the analytic Nullstellensatz, and this does not hold in general if $Y$ is non-strict
(\eg, $Y=\mathscr M(k_r)$ for some non-empty $k$-free polyradius $r$). 

(Let us emphasize that we do not claim to have given a new proof of the
analytic
Nullstellensatz, because we used it at several places in the proof of Theorem \ref{thm-multisections-local}, 
either directly or indirectly; \eg, it
is used for the local existence of \'etale multisections
on the image of a smooth map.)

\end{rema}

\section{Images of maps}\label{s-images-maps}

We are now going to deduce general results on the images of maps
from Theorem \ref{thm-multisections-global}. Theorems \ref{th-image-compact}
and \ref{th-image-equi}
are stated under slightly
general assumptions
(for instance, the spaces involved are not assumed to be Hausdorff), but
for each of them the case of interest is
mainly that of affinoid source and target
(and in fact, the proof first reduces to this
case for both of them).

\begin{theo}\label{th-image-compact}
Let~$Y$ be a~$\Gamma$-strict~$k$-analytic space, let~$X$ be a~$k$-analytic space, and let~$\phi \colon
Y\to X$ be a
morphism. Assume that $\phi$ is
\emph{topologically}
proper (\ref{ss-remind-topo}; this amounts to requiring that $\phi\inv(V)$ is quasi-compact for every affinoid domain $V$ of $X$, see \ref{ss-topologies}; we emphasize that
we do not
assume that $\phi$ is topologically separated).
Let~$\mathscr F$ be a coherent sheaf on~$Y$ which 
is~$X$-flat. Assume that at least one of the two conditions below is satisfied: 

\begin{enumerate}[1]
\item The space $X$ is~$\Gamma$-strict. 
\item The space $X$ is separated and $Y$ admits a locally finite covering by $\Gamma$-strict affinoid domains (if $Y$ is Hausdorff, the latter condition 
is equivalent to paracompactness 
of $Y$; see \ref{ss-topologies} and Remark \ref{rem-topo-gstrict}). 
\end{enumerate}
Then~$\phi(\supp F)$ is a closed~$\Gamma$-strict analytic domain of~$X$.
\end{theo}

\begin{proof}
By replacing~$Y$ with~$\supp F$ we reduce to the case where~$\supp F=Y$. Since~$\phi$ is
topologically proper, $\phi(Y)$ is a closed subset of~$X$ and $\phi\inv(E)$ is quasi-compact for every quasi-compact
subset $E$ of $X$. We are now going to reduce to the case where both~$Y$ and~$X$ are~$\Gamma$-strict and~$k$-affinoid.

Let us first consider the case where (1) is fulfilled. One can check the result G-locally on~$X$, which allows to assume that~$X$ is~$\Gamma$-strict and~$k$-affinoid.
In this case,~$Y$ is quasi-compact, hence admits a finite covering by~$\Gamma$-strict, affinoid domains; one therefore immediately reduces to the case where~$Y$ is also~$\Gamma$-strict and affinoid.

Let us now consider the case where (2) is fulfilled. Choose a locally finite $\Gamma$-strict affinoid covering $(Y_i)$ of $Y$. It is sufficient to prove that~$\phi(Y_i)$ is a~$\Gamma$-strict compact analytic domain of~$X$ for any~$i$. Indeed, assume that it is the case. Then for every affinoid domain $V$ of $X$, the pre-image $\phi\inv(V)$ is quasi-compact, hence intersects only finitely many $Y_i$'s; this implies that $(\phi(Y_i))_i$
is a locally finite covering of~$\phi(Y)$ by~$\Gamma$-strict compact analytic domains of~$X$, which can be refined into a locally finite covering by $\Gamma$-strict affinoid domains. As~$X$ is separated, the intersection of two such domains will still be affinoid and $\Gamma$-strict; hence our covering is a $\Gamma$-strict affinoid atlas on $\phi(Y)$, and $\phi(Y)$ is a closed $\Gamma$-strict analytic domain of $X$.
We thus reduce
to the case where~$Y$ is compact. By Proposition~\ref{prop-image-subgstrict}, we can then assume~$X$ is compact and~$\Gamma$-strict, and even, since one can check the result G-locally on~$X$, that it is affinoid and~$\Gamma$-strict. And as~$Y$ admits a finite covering by~$\Gamma$-strict~$k$-affinoid domains, we eventually reduce to the case where~$Y$ is also~$\Gamma$-strict and~$k$-affinoid. 

So let us prove the theorem when both
spaces $Y$ and $X$ are affinoid and $\Gamma$-strict. Let~$r=(r_1,\ldots,r_n)$ be a~$k$-free polyradius such that the~$r_i$'s belong to~$\Gamma$, the valuation of $k_r$ is non-trivial, and $X_r$ and~$Y_r$ are strictly~$k_r$-affinoid.
Let~$\mathfrak s : X\to X_r$ be the Shilov section
(\ref{ss-shilov-section}).
By Theorem \ref{thm-multisections-global}, the image~$\phi_r(Y_r)$ is a compact strictly~$k_r$-analytic domain of~$X_r$. The subset~$\phi(Y)$ of~$X$ is nothing but~$\mathfrak{s}\inv (\phi_r(Y_r))$.
By the Gerritzen-Grauert theorem,~$\phi_r(Y_r)$ is a finite union of strictly~$k_r$-rational domains; it is then easily seen (using the explicit formula for~$\mathfrak s$, see the proof of Lemma 2.4 of \cite{ducros2003}) that~$\mathfrak{s}\inv (\phi_r(Y_r))$ is itself a finite union of rational domains whose definitions only involve elements of $\R_+\gpm$
that belong to~$\Gamma$; therefore,
$\phi(Y)$ is a compact~$\Gamma$-strict analytic domain of~$X$.
\end{proof}

\begin{theo}\label{th-image-equi}
Let~$n$ and~$d$ be two non-negative
integers, let~$Y$ be a~$\Gamma$-strict~$k$-analytic space, and let~$\phi$ be a topologically
proper morphism from~$Y$ to a normal~$k$-analytic space~$X$. Assume that $X$ is purely~$d$-dimensional, $Y$ is purely~$(n+d)$-dimensional,
and the fibers of $\phi$
are purely
$n$-dimensional. If~$X$ is~$\Gamma$-strict, or if~$X$ is separated and $Y$ admits a locally finite $\Gamma$-strict affinoid covering, then~$\phi(Y)$ is a~$\Gamma$-strict closed~$k$-analytic domain of~$X$.
\end{theo}

\begin{proof}
Exactly as at the beginning of the proof of Theorem~\ref{th-image-compact},
we reduce
to the case where both~$Y$ and~$X$ are~$\Gamma$-strict~$k$-affinoid spaces. By compactness, one can argue locally on~$Y$. Hence,
by combining \Cor 4.7 of \cite{ducros2007} with \ref{ss-pullback-gstrict}
(or Remark \ref{rem-gstrict-fini-transfer}),
we can assume that there exists a factorization~$Y\to T \to X$ where the map
$Y\to T$ is finite
and where~$T$ is a $\Gamma$-strict affinoid domain of a smooth $X$-space of pure relative dimension~$n$.

By flatness of quasi-smooth morphisms (Corollary \ref{cor-qsm-flat}) and Theorem \ref{th-image-compact}, the image of
any non-empty compact analytic domain $T'$ of $T$ in $X$ is a non-empty compact analytic
domain of $X$, hence is of dimension $d$.
Since the fibers of $T\to X$ are purely $n$-dimensional, it follows from \ref{ss-dim-sobafi} (3)
that $T'$ is of dimension $n+d$; therefore $T$ is
purely~$(n+d)$-dimensional. As~$Y\to T$ is finite and
$Y$ is purely $(n+d)$-dimensional, the image of
any irreducible component of $Y$ in $T$ is an irreducible Zariski-closed subset of $T$ of dimension~$n+d$
(again by \ref{ss-dim-sobafi}); \ie, 
this is an irreducible component of $T$. 
Therefore the image $Z$ of $Y$ in $T$ is a union of irreducible components of~$T$. Since~$X$ is normal and
$T\to X$ is quasi-smooth,~$T$ is normal by Proposition \ref{prop-tranferqsm-concrete}. Therefore~$Z$ is a union of connected components of~$T$, hence is a~$\Gamma$-strict affinoid domain of~$T$. By Theorem \ref{th-image-compact}  the image of~$Z$
in~$X$, which coincides with that of~$Y$, is a
compact~$\Gamma$-strict analytic domain of~$X$.
\end{proof}

\begin{theo}\label{th-image-compact-inner}
Let~$\phi\colon Y\to X$ be a morphism between~$k$-analytic spaces and let~$\mathscr F$ be a coherent sheaf on~$Y$ which is~$X$-flat. Let~$y$ be a point of~$\supp F$ at which~$\phi$ is inner, and let~$x$ be its image on~$X$. The image~$\phi(\supp F)$ is a neighborhood of~$x$.
\end{theo}

\begin{proof}
By shrinking $X$ around $x$ (and by shrinking $Y$ accordingly)
we may assume that $X$ is Hausdorff. By replacing
$Y$ with a compact analytic neighborhood of~$y$, we can further assume that~$Y$ is compact.
Then it
follows from Theorem \ref{th-image-compact}
that~$\phi(\supp F)$ is an analytic domain~$U$ of~$X$. As~$\phi$ is inner at~$y$,~$\phi_{|\supp F}$ is inner at~$y$ too. Therefore~$x$ belongs to~$\mathrm{Int}(U/X)$; \ie, to the
\emph{topological} interior of~$U$ in~$X$, whence the result.
\end{proof}

\begin{rema}
The openness of flat, boundaryless morphisms between good~$k$-analytic spaces has already been proved by Berkovich, in a slightly different way, in
unpublished
work. 
\end{rema}

\begin{theo}\label{th-image-equi-inner}
Let~$n$ and~$d$ be two
non-negative integers and let~$\phi
\colon Y\to X$ be a morphism between~$k$-analytic spaces. Assume that~$X$ is normal and
purely $d$-dimensional, $\phi$ is purely of relative dimension~$n$, and $Y$ is of pure dimension~$n+d$. Let~$y$ be a point of~$Y$ at which~$\phi$ is inner, and let~$x$ be its image on~$X$. The image~$\phi(Y)$ is a neighborhood of~$x$.
\end{theo}

\begin{proof}
By shrinking $X$ around $x$ (and by shrinking $Y$ accordingly)
we may assume that $X$ is Hausdorff. By replacing
$Y$ with a compact analytic neighborhood of~$y$, we can further assume that~$Y$ is compact.
Then it follows from
Theorem \ref{th-image-equi}
that~$\phi(Y)$ is an analytic domain~$U$ of~$X$. As~$\phi$ is inner at~$y$,~$x$ belongs to~$\mathrm{Int}(U/X)$; \ie, to the {\em topological} interior of~$U$ in~$X$, whence the result.
\end{proof}

\chapter{Constructible loci} \label{LOC}
This quite long chapter is devoted to the study of ``loci of validity". In order to describe them,
we need to fix some terminology. We shall say that a subset of an analytic space $X$
is \emph{constructible}
if it is a finite boolean combination of Zariski-open subsets of $X$ (Definition \ref{def-constr}
below). Though the properties of being Zariski-open or Zariski-closed are G-local, that of being constructible
is \emph{not}: we
sketch a counter-example (\ref{contrex-gcons}). But it is nonetheless ``almost" G-local: if $X$
is a \emph{finite-dimensional}
analytic space (which simply means that the dimensions of the irreducible components of $X$ are uniformly
bounded above), then every G-locally constructible subset of $X$ is actually constructible (Proposition
\ref{prop-cons-gloc}). 

We now aim to establish
that various loci of validity of relative properties
are locally constructible
(by the above, they will even be constructible as soon as the source space
is finite-dimensional), and sometimes Zariski-open when the property involved encapsulates some flatness
condition. 

For that purpose we develop in Section \ref{s-diagtr},
in a rather abstract categorical setting, 
a general strategy inspired by Kiehl's paper \cite{kiehl1967}
and by the technique of ``spreading out from the generic fiber"
in algebraic geometry -- the latter is not directly available here, but as explained in
the general Introduction
(see \ref{ss-obstacle3}), we bypass this obstacle
by using Theorem \ref{thm-localring-generic}
that applies to
local rings of generic fibers.
This strategy will enable us
to reduce various constructibility and Zariski-openness statements to simpler ones
(see the short introduction of Section \ref{s-diagtr}
for some
more specificity). One may skip those quite formal 
considerations and go directly to Section \ref{s-flatloc}, but one will then have to accept
repeatedly
arguments of the form 
``in view of this and that result of \ref{s-diagtr}, we may assume\ldots".

We do
not give further details here about the results of this chapter; we refer the reader to
the local introduction of each section, from \ref{s-diagtr}
to \ref{s-cons-maintheo}.

\section{Constructibility in analytic spaces}\label{s-construct}

\begin{defi}\label{def-constr}
Let~$X$ be an analytic space. 
We say that a subset~$E$ of~$X$ is \emph{constructible}\index{subset!constructible}\index{constructible subset}
if it can be written as a finite union~$\bigcup (U_i\cap F_i)$ where~$U_i$ (\resp~$F_i$) is
a Zariski-open (\resp a Zariski-closed) subset of~$X$ for every $i$. 

We say that a subset $E$ of $X$ is \emph{locally constructible} (\resp \emph{G-locally constructible})\index{subset!locally constructible}\index{subset!G-locally constructible}
if there exists an open covering (\resp a G-covering) $(X_i)$ of $X$ such that $E\cap X_i$ is\index{G-locally constructible subset}\index{locally constructible subset}
a constructible subset of $X_i$ for all $i$. 

\end{defi}

\begin{rema}
The set of constructible subsets of $X$ is a Boolean sub-algebra of $\mathscr P(X)$;
\ie, it is
stable under finite unions, finite intersections, and complements. 
We could also have defined it as the Boolean sub-algebra of $\mathscr P(X)$ generated by all Zariski-open subsets
(or all Zariski-closed subsets).

The set of locally constructible (\resp G-locally constructible) subsets of $X$ is also
a Boolean sub-algebra of $\mathscr P(X)$.

\end{rema}

\begin{rema}
In EGA, the Boolean algebra of 
constructible subsets of an arbitrary topological space $X$
is defined as the Boolean sub-algebra of $\mathscr P(X)$ 
generated by \emph{retrocompact}
open subsets; \ie, open subsets whose intersection with any
quasi-compact open subset
of $X$ is quasi-compact. If $X$ is a locally noetherian
topological space, then every open subset of $X$ is retrocompact.

Now if $X$ is a
\emph{quasi-compact}
analytic space, 
its Zariski topology is noetherian; therefore every Zariski-open subset of $X$
is retrocompact with respect to the Zariski topology, hence a subset of $X$
is constructible in our sense if and only if it is constructible in the sense of EGA
for the Zariski topology. We do not know whether this is the case
for an arbitrary analytic space, since it is unclear (at least to the author) which Zariski open subsets
of an arbitrary analytic space are retrocompact with respect to the Zariski topology.
\end{rema}

\begin{exem}
Let $\mathscr X$ be a scheme locally of finite type
over an affinoid algebra
$A$. 
If $E$ is a constructible (\resp locally constructible)
subset of $\mathscr X$, then
its pre-image $E\an$ in $\mathscr X\an$
is a
constructible (\resp locally constructible)
subset of $\mathscr X\an$; it is closed, \resp open, if
and only if $E$ is closed, \resp open (\cite{berkovich1993}, \Cor 2.6.6; this
is written for a constructible subset, but since the problem is local for the Zariski topology
on $\mathscr X$, one reduces immediately to this case). 
If $\mathscr X$ is proper over $A$
it follows from GAGA (\ref{ss-def-analytify})
that $E\mapsto E\an$ establishes a bijection between the set
of constructible subsets of $\mathscr X$ and that of constructible subsets of
$\mathscr X\an$, whose
converse bijection maps a constructible subset $F$ of $\mathscr X\an$
to its image $F\al$ in $\mathscr X$.
\end{exem}

\subsection{}
Let $X$ be analytic space, and let $Y$ be a Zariski-closed subset of $X$. 
Endow $Y$ with any structure of a closed analytic subspace of $X$. Then a subset $E$
of $Y$ is constructible (with respect to this given structure) if and only if it is constructible
as a subset of $X$. In particular, this does not depend on the chosen structure of $Y$, 
and we shall often
simply speak about constructible subsets of $Y$, without fixing 
any structure of closed analytic subspace on it. The same also holds for
locally constructible and G-locally constructible sets (for the latter, this rests on Remark \ref{rem-gerr-grau}).

\subsection{}
Let $Y\to X$ be a morphism of analytic spaces. If $E$ is a constructible
(\resp locally constructible, \resp G-locally constructible)
subset of $X$, it follows from the definition that the pre-image of $E$ in $Y$
is a constructible (\resp locally constructible, \resp G-locally constructible)
subset of $Y$. We shall often apply it for $Y$ an analytic domain 
of $X$, in which case the pre-image of $E$ in $Y$ is nothing but the intersection $E\cap Y$. 

\subsection{}\label{ss-constop-compact}
Let $X$ be a $k$-analytic space. Assume that
its Zariski-topology is noetherian; \eg, $X$
is quasi-compact. Since every irreducible
Zariski-closed subset $Y$ of $X$
has a Zariski-generic point (pick $y$ in $Y$
such that $d_k(y)=\dim Y$), it follows from \cite{ega31}, \Chp 0, \Cor 
9.2.4 that $X$ is quasi-compact for the constructible topology.
If $X$ is affinoid, this can also be deduced through
the assignement $E\mapsto E\an$
from the compactness of $X\al$
for the constructible topology. But be aware that the Hausdorff
property of the constructible topology of $X\al$ does
not transfer to the constructible topology of $X$ in general: 
indeed, two distinct points of $X$ lying over
the same point of $X\al$ belong to the same constructible subsets of $X$, hence
cannot be separated using such subsets. 

We are now going to show that some other basic
properties of the constructible subsets of a noetherian topological 
space (like \Prop 9.2.2 and \Prop 9.2.5 of \cite{ega31}, \Chp 0)
actually hold for the constructible subsets of an arbitrary analytic space. 
Our proofs are basically the same as those of the results in EGA alluded to, 
except that we replace noetherian induction by some arguments involving 
 dimension theory
and
the decomposition into irreducible components.

\begin{lemm}\label{lem-cons-dense}
Let $X$ be an analytic space and let $E$ be a constructible
subset of $X$. 

\begin{enumerate}[1]
\item The following are equivalent:
\begin{enumerate}[j]

\item $E$ contains a Zariski-dense open subset of $X$; 

\item $E$ is Zariski-dense in $X$. 
\end{enumerate}

\item The following are equivalent:
\begin{enumerate}[j]
\setcounter{enumii}{2}
\item $E$ is Zariski-closed in $X$; 
\item every irreducible Zariski-closed subset $Z$ of $X$ such that
$E\cap Z$ contains a non-empty Zariski-open subset of $Z$
is contained in $E$. 
\end{enumerate}
\item The following are equivalent:
\begin{enumerate}[j]
\setcounter{enumii}{4}
\item $E$ is Zariski-open in $X$; 
\item for every pair $(T,Z)$ of irreducible Zariski-closed subsets of 
$X$ with $T\subset Z$, if $T\cap E$ is Zariski-dense in $T$ then $Z\cap E$
is Zariski-dense in $Z$. 
\end{enumerate}
\end{enumerate}
\end{lemm}

\begin{proof}
Let us begin with (1).
We obviously have (i)$\Rightarrow$(ii).
Let us now assume that (ii) holds. Write $E=\bigcup_{i\in I} U_i\cap F_i$, where $I$ 
is a finite set and where $U_i$ is Zariski-open and $F_i$ Zariski-closed for all $i$.
Let $(X_j)$ be the family of irreducible
components of $X$. 
Fix $j$.
Any non-empty Zariski-open subset $U$ of $X_j$ contains a non-empty Zariski-open subset of $X$, namely
$U\cap \bigcup_{\ell \neq j}X_\ell$. Therefore the intersection $E\cap X_j$ is a Zariski-dense subset of $X_j$. This implies
that there exists some index $i$ such that $F_i\cap X_j=X_j$ and $U_i\cap X_j\neq \emptyset$.
Let $V_j$ be the intersection of $U_i\cap X_j$ with the complement
of $\bigcup_{\ell \neq j}X_\ell$. 
By construction, $\coprod_{j\in J}V_j$ is a Zariski-dense open subset of $X$ contained
in $E$, and (i) holds. 

Let us prove (2). The implication (iii)$\Rightarrow$(iv)
is obvious. Assume that (iv) holds. In order to prove that $E$ is Zariski-closed in $X$, it suffices
to prove that $E\cap Y$ is Zariski-closed for every irreducible component $Y$ of $X$ (because
the set of irreducible components is G-locally finite and the property of being Zariski-closed is G-local); 
we can thus assume that $X$ is irreducible, and we argue by induction on $\dim X$. There is nothing
to prove if $\dim X=0$; assume that $\dim X>0$ and 
the assertion is true in smaller dimensions. Write $E=\bigcup_{i\in I} U_i\cap F_i$, where $I$ 
is a finite set and where $U_i$ is Zariski-open; we can of course assume that $U_i\cap F_i\neq
\emptyset$ for all $i$. If there exists $i$ such that $F_i=X$ then $E$ contains the non-empty Zariski-open
subset $U_i$ of $X$, hence $E=X$ by assumption (iv). If $F_i\neq X$ for all $i$
then $F:=\bigcup_i F_i$ is a Zariski-closed subset of $X$
containing $E$ and of dimension $<\dim X$; then by
arguing componentwise on $F$ and using the induction hypothesis we see that 
$E$ is Zariski-closed in $F$ (hence in $X$). 

Let us prove (3). The implication (v)$\Rightarrow$(vi)
is obvious. Assume that (iv) holds, and set $F=X\setminus E$. Let $Z$ be an irreducible
subset of $X$ such that $F\cap Z$ contains a non-empty Zariski-open subset $U$ of $Z$. 
We are going to prove by contradiction that $F$ contains $Z$; this will ensure in view of (2) that $F$ is Zariski-closed,
and thus that $E$ is Zarsiki-open. 
Assume that there exists $z\in Z$ such that $z\notin F$. Then $z\in E$ and $E\cap \adhz{\{z\}}X$
is Zariski-dense in $\adhz{\{z\}}X$; by assumption (3) this implies that $E\cap Z$ is Zariski-dense
in $Z$, which contradicts the fact that $U\subset F$.
\end{proof}

\begin{rema}\label{rem-aff-topcons}
If $X$ is affinoid, Lemma \ref{lem-cons-dense}
can of course also be deduced through
the assignement $E\mapsto E\an$
from the corresponding statement on the constructible
subsets of $X\al$. 
\end{rema}

\begin{lemm}\label{lem-adh-const}
Let~$X$ be a
$k$-analytic space, let~$E$ be a constructible subset of~$X$,
and let~$V$ be an analytic domain of~$X$. Let~$L$ be an analytic extension of~$k$ and let~$E_L$ be
the pre-image of~$E$ on~$X_L$. 

\begin{enumerate}[1]

\item $\adht E X=\adhz E X$.

\item $\adhz{E\cap V}V=\adhz E X \cap V$. 

\item $\adht{E_L}{X_{L,\mathrm{Zar}}}=(\adhz E X)_L.$

\item $E$ is open (\resp closed) if and only if it is Zariski-open (\resp Zariski-closed). 
\end{enumerate}
\end{lemm}

\begin{proof}
Assertion (4) is an immediate consequence of (1). Let us now
prove (1), (2), and (3). All terms involved in these equalities commute
with finite unions; one can therefore assume that~$E=U\cap F$,
where~$U$ (\resp~$F$) is a Zariski-open (\resp
Zariski-closed) subset of~$X$; by replacing~$X$ with~$F$ one can assume that~$E$ is Zariski-open in~$X$.
The required equalities now follow from Lemma~\ref{lem-closure-zartop}, Corollary~\ref{cor-zarclosure-dom},
and Corollary~\ref{cor-zarclosure-xl}.
\end{proof}

\begin{coro}\label{cor-adh-const}
Let~$X$ be a~$k$-analytic space, let~$E$ be a subset of~$X$
and let $(X_i)$ be a G-covering
of~$X$ by
analytic domains such that
the intersection $E\cap X_i$ is a
constructible subset of~$X_i$ for every~$i$. Under those assumptions
$\adhz EX=\adht EX$ and $\adht E X \cap X_i=\adht{E\cap X_i}{X_i}$ for every~$i$.
\end{coro}

\begin{proof}
Let~$i$ and~$j$ be two indices. By Lemma \ref{lem-adh-const}
\[\adht{E\cap X_i}{X_i}\cap X_j=\adht{E\cap X_j}{X_j}\cap X_i=\adht{E\cap X_i\cap X_j}{(X_i\cap X_j)}.\]
Therefore if one sets~$F=\bigcup_i \adht{E\cap X_i}{X_i}$,
one has~$F\cap X_i=\adht{E\cap X_i}{X_i}$ for every~$i$. As a consequence,~$F$ is a Zariski-closed subset of~$X$;
it follows from
the definitions that if~$Z$ is a closed subset of~$X$ containing~$E$ then~$Z\supset F$. One has thus~$F=\adht E X=\adhz EX$. 
\end{proof}

%\begin{lemm}\label{lem-cons-dense}
%Let $X$ be an analytic space and let $E$ be a constructible
%subset of $X$. The following are equivalent:
%
%\begin{enumerate}[i]
%
%\item $E$ contains a Zariski-dense open subset of $X$; 
%
%\item $E$ is Zariski-dense in $X$. 
%
%
%\end{enumerate}
%
%\end{lemm}
%
%\begin{proof}
%We obviously have (i)$\Rightarrow$(ii).
%Let us now assume that (ii) holds. Write $E=\bigcup_{i\in I} U_i\cap F_i$, where $I$ 
%is a finite set and where $U_i$ is Zariski-open and $F_i$ Zariski-closed for all $i$.
%Let $(X_j)$ be the family of irreducible
%components of $X$. 
%
%Fix $j$.
%Any non-empty Zariski-open subset $U$ of $X_j$ contains a non-empty Zariski-open subset of $X$, namely
%$U\cap \bigcup_{\ell \neq j}X_\ell$. Therefore the intersection $E\cap X_j$ is a Zariski-dense subset of $X_j$. This implies
%that there exists some index $i$ such that $F_i\cap X_j=X_j$ and $U_i\cap X_j\neq \emptyset$.
%Let $V_j$ be the intersection of $U_i\cap X_j$ with the complement
%of $\bigcup_{\ell \neq j}X_\ell$. 
%
%By construction, $\coprod_{j\in J}V_j$ is a Zariski-dense open subset of $X$ contained
%in $E$. 
%\end{proof}

\begin{prop}\label{prop-cons-gloc}
Let~$X$ be an analytic space and let~$E$ be a subset of~$X$. The following are equivalent:

\begin{enumerate}[i]
\item $E$ is G-locally constructible. 

\item $E$ is locally constructible. 

\end{enumerate}

If moreover $X$ is finite-dimensional, those assertions
are also equivalent to

\begin{enumerate}[i]\setcounter{enumi}2
\item $E$ is constructible. 
\end{enumerate}

\end{prop}

\begin{proof}
It is clear that (iii)$\Rightarrow$(ii)$\Rightarrow$(i).

Assume that $X$ is finite-dimensional and satisfies (i); choose a G-covering $(X_i)$ of $X$ such that $E\cap X_i$
is a constructible subset of $X_i$ for all $i$, and let us prove
that
$X$ satisfies (iii). 
We argue by induction on
$\dim X\in \{-\infty\}\cup \N$. If $X=\emptyset$ there is nothing to prove. Assume that $X$ is non-empty,
and that the proposition is true in dimensions $<\dim X$. 

Set $Y=\adhz E X$
and $F=Y\setminus E$; for every $i$, we denote by $Y_i$ (\resp $E_i$, \resp $F_i$)
the intersection of $X_i$ with $Y$ (\resp $E$, \resp $F$). 

Fix $i$. It follows from Corollary~\ref{cor-adh-const}
that $E_i$ is a Zariski-dense constructible subset of $Y_i$ and 
$\adhz F Y\cap Y_i=\adht {F_i}{Y_{i,\mathrm{Zar}}}$. 
Since $E_i$ is a Zariski-dense constructible subset of
$Y_i$,  Lemma \ref{lem-cons-dense}
implies that $E_i$ contains a Zariski-dense open subset of $Y_i$; as a consequence, 
$\adht {F_i}{Y_{i,\mathrm{Zar}}}$ contains no irreducible components of $Y_i$. We thus have
$\dim \adht {F_i}{Y_{i,\mathrm{Zar}}}<\dim Y_i$ as soon as $Y_i\neq \emptyset$. 

By the above, $\dim (\adhz F Y\cap Y_i)<\dim Y_i$ for every $i$ such that $Y_i\neq \varnothing$. 
This implies that $\adhz F Y$ is of dimension $<\dim X$.
The induction hypothesis then ensures that~$F$ is a constructible subset of~$\adhz F Y$,
and therefore a constructible subset of~$X$. As a consequence,~$E=Y\setminus F$ is a constructible subset of~$X$, whence (iii). 

It remains to show that (i)$\Rightarrow$(ii). 
Assume that $E$ is G-locally constructible. Let $x$ be a point of $X$, and let
$V$ be a finite-dimensional
open neighborhood of $x$
in $X$ (\eg, $V$ is the topological interior of a compact analytic neighborhood of $x$).
Since $E$ is G-locally constructible, $E\cap V$ is a G-locally constructible
subset of $V$. Since (i)$\Rightarrow$(iii) for a finite-dimensional
ambient space, 
it follows that
$E\cap V$ is a constructible subset of $V$.
Hence (ii)
holds (by varying $x$). 
\end{proof}

\begin{rema}\label{cons-local-scheme}
If $\mathscr X$ is a finite-dimensional locally noetherian scheme, every
locally constructible subset of $\mathscr X$ is constructible; the proof is mutatis
mutandis the same as that of implication (i)$\Rightarrow$(iii)
in Proposition \ref{prop-cons-gloc} above. But note that
this is obvious
if $\mathscr X$ is of finite type (by transitivity of the Zariski topology in the
scheme-theoretic setting). 
\end{rema}

\begin{enonce}[remark]{Counter-example}\label{contrex-gcons}
In Proposition \ref{prop-cons-gloc}
above, the assumption that $X$ is finite-dimensional cannot be dropped. 
Indeed, let us denote by $X$ the closed unit disc
over $k$.
For every $n\in \Z$, let $j_n$
be the closed immersion $X^n\simeq X^n\times\{0\}\hookrightarrow X^{n+1}$. 
We define inductively a constructible subset $E_n$ of $X^n$ by the following conditions:

\begin{itemize}[label=$\bullet$]
\item $E_0=X^0=\{0\}$. 

\item $E_{n+1}=X^{n+1}\setminus j_n(E_n)$ for every $n$. 
\end{itemize}

By construction, the subset
$\coprod_n E_n$ of $\coprod_n X^n$
is locally constructible; but it is not constructible (exercise left to the reader). 
\end{enonce}

We end this section by proving
a kind of analytic Chevalley theorem for proper maps; the key point
will be Kiehl's theorem on the direct image of a coherent sheaves (\cf \ref{ss-kiehl-proper}). 

\begin{theo}[Proper Chevalley theorem]\label{thm-chevalley-proper}
Let $\phi \colon Y\to X$ be a proper morphism of $k$-analytic spaces and let $E$
be a locally constructible subset of $Y$. The image $\phi(E)$ is a locally constructible
subset of $X$.
\end{theo}

\begin{proof}
Let us first mention that since $\phi$ is proper, 
$\phi(T)$ is a Zariski-closed subset of $X$ for every Zariski-closed
subset $T$ of $Y$ (by Kiehl's theorem on the direct images of coherent
sheaves, \cf \ref{ss-kiehl-proper}); we shall use it repeatedly throughout the proof. 

By Proposition 
\ref{prop-cons-gloc},
the assertion is G-local on $X$; we can thus asume that it is affinoid.
The space $Y$
is then compact by topological properness and topological
separatedness, and $E$ is thus constructible by
Proposition \ref{prop-cons-gloc}. We can therefore
write $E=\bigcup_{i\in I} U_i\cap F_i$ where
$I$ is a finite set and $U_i$ (\resp $F_i$) is for every $i$ a Zariski-closed
(\resp Zariski-open) subset of $Y$; we can moreover assume that $F_i$ is
irreducible and $U_i\cap F_i \neq \emptyset$ for all $i$.

We are going to prove by noetherian induction on $X$ that $\phi(E)$ is a constructible
subset of $X$. We thus assume that
the intersection
of $\phi(E)$
with every proper Zariski-closed subset of $X$
is constructible,
and we shall prove that $\phi(E)$ is constructible. Let $(X_j)$ be the family of irreducible
components of $X$. 

If $X$ is not irreducible, then $X_j\subsetneq X$ for all $j$, and
$\phi(E)=\bigcup \phi(E)\cap X_j$ is thus a constructible subset of $X$. 
Assume now that $X$ is irreducible. For every $i$, the image $\phi(F_i)$ is an irreducible
Zariski-closed subset of $X$, and $\phi(E)\subset \bigcup_i \phi(F_i)$. As a consequence, 
if $ \bigcup_i \phi(F_i)\subsetneq X$, then $\phi(E)$ is constructible by
the induction hypothesis. It remains to consider the case where there exists $i$
such that $\phi(F_i)=X$; we write $F$ and $U$ instead of $F_i$ and $U_i$, respectively.
It suffices now to prove that there exists a proper Zariski-closed
subset $Z$ of $X$ such that $\phi(F\cap U)$ contains $X\setminus Z$. 
indeed, if this is the case, one will have
$\phi(E)=(X\setminus Z)\cup (\phi(E)\cap Z)$,
and since $\phi(E)\cap Z$ is contructible by the induction hypothesis, 
we shall be done. 

Set $n=\dim X$ and $m=\dim F$. Let $\xi$ be an Abhyankar point of $X$ 
(\ref{ss-abhyankar-points}). By Lemma \ref{lem-abhyankar-dimension}, 
the fiber $F_\xi$ is of pure dimension $n-m$. As it is non-empty (recall
that $\phi(F)=X$
by assumption),  there exists a point $y\in F_\xi$ such that $d_{\hr \xi}(y)=n-m$. We then have
$d_k(y)=n-m+d_k(\xi)=n$. The point $y$ is thus an Abhyankar point of $F$, hence is Zariski-dense
in $F$ (Remark \ref{lem-abhyankar-dimension}). Since the relative
dimension is upper semi-continuous for the Zariski topology (\cite{ducros2007}, \Th 4.9), 
the minimal relative dimension of $\phi|_F$ is equal to $n-m$, and the set of points 
$z$ of $F$ such that $\dim_z \phi>m-n$ is a proper Zariski-closed subset 
of $F$, whose image in
$X$ does not contain $\xi$ (because $F_\xi$ is of pure
dimension $n-m$), hence is a proper Zariski-closed subset of $X$. 

Let $G$ be the complement of $F\cap U$ in $F$. This is a proper Zariski-closed subset of $F$;
let $(G_\ell)_{\ell \in \Lambda}$
be the family (possibly empty if $G=\emptyset$)
of irreducible components of $G$, and let $\Lambda_0$ be the subset of $\Lambda$
consiting of indexes $\ell$ such that $\phi(G_\ell)=X$.

Let $\ell$ be an element of
$\Lambda_0$. By the same reasoning as above,
the minimal relative dimension of
$\phi|_{G_\ell}$ is equal to $\dim G_\ell-n$, the subset $H_\ell$
of $G_\ell$ consisting of points $z$ such that
$\dim_z \phi|_{G_\ell}>\dim G_\ell-n$ is a proper Zariski-closed
subset of $G_\ell$, and $\phi(H_\ell)$ is a proper Zariski-closed subset of $X$. 

Set $H=(\bigcup_{\ell \in \Lambda\setminus \Lambda_0}G_\ell)\cup \bigcup_{\ell \in \Lambda_0}H_\ell$. 
By construction, $H$ is a proper Zariski-closed subset of $F$, and $Z:=\phi(H)$ is a proper Zariski-closed subset of $X$. 
It suffices now to prove that $X\setminus Z\subset \phi(F\cap U)$.

Let $x$ be a point of $X$ that does not lie on $Z$.
The fiber $G_x$ does not intersect $\bigcup_{\ell \in \Lambda\setminus \Lambda_0}G_\ell$, nor
$\bigcup_{\ell \in \Lambda_0}H_\ell$. Since $\dim_y G_{\ell,\phi(y)}=\dim G_\ell-n<m-n$ for every
$\ell \in \Lambda_0$ and every $y\in G_\ell \setminus H_\ell$ (by definition of $H_\ell$), it follows that
$\dim G_x<m-n$. But  every fiber of $\phi|_F$ is non-empty and of dimension at least $m-n$; as a consequence, 
$G_x\subsetneq F_x$, which exactly means that $F_x$ intersects $F\cap U$; \ie, $x\in \phi(F\cap U)$. 
\end{proof}

\section{The diagonal trick}\label{s-diagtr}

We
are going to
describe here a general method which will be useful for
establishing the (local) constructibility or the Zariski-openness of some loci.
It is inspired by 
Kiehl's \cite{kiehl1967}, in which he proved the Zariski-openness of the flat
locus of a
\emph{complex-analytic} morphism.

Roughly speaking, it consists of the following. One wants to understand
the locus of validity of some relative property of a morphism $Y\to X$. 
If the formation of this locus commutes with
the
``tautological" base-change by $Y\to X$, one can replace
$Y\to X$ by the second projection $Y\times_X Y\to Y$
and only investigate
what happens on the diagonal. 
One thus
reduces to the case where $Y\to X$ has a section $\sigma$
and where it suffices to understand the intersection of our locus
of validity with $\sigma(X)$.  

Considering only points lying of $\sigma(X)$ will
bring two advantages, say, in the case
where both $Y$ and $X$ are affinoid (as can most of the time be assumed without loss of generality
by arguing G-locally):

\begin{itemize}[label=$\bullet$]

\item The scheme $\sigma(X)\al$ is of finite type over $X\al$ (it is even isomorphic to it!), which is
of course in general not the case for
$Y\al$; this finiteness condition plays a key role
in the study of the flatness
locus, through an intermediate theorem of Kiehl on morphisms between noetherian
schemes (see Theorem \ref{thm-kiehl-flat} below). 

\item Every point of $\sigma(X)$ belongs to $\mathrm{Int}(Y/X)$ by
\ref{ss-boundary-basics} (3). This
will allow
us to apply Theorem \ref{thm-localring-generic}
on the local rings of analytic fibers, and to only deal
with \emph{naive} flatness in view of Theorem
\ref{thm-flat-naiveflat}.

\end{itemize}

\subsection{Our general axiomatic setting}\label{ss-general-axioms-pfib}
We fix a subcategory~$\mathfrak C$ of 
the category of analytic spaces such that for every
object $X$ of $\mathfrak C$, 
the category of analytic spaces over $X$ (which contains all $X$-analytic spaces, but
also more generally all
$X_L$-analytic spaces for every complete extension $L$ of $k$)
embeds fully faithfully in~$\mathfrak  C$.
We still use the notations $\mathfrak T$, $\mathfrak L$ and
and $\mathfrak F$ of Section
\ref{s-category-framework}, and the notation~$\mathfrak{Coh}$
and~$\mathfrak{Coh}^{\mathfrak I}$ introduced in
Examples \ref{ex-fiber-coh}
and \ref{ex-fiber-diag}. 

We denote
by $\mathsf Q$ a property whose
validity at a given point
of an object $X$ of $\mathfrak C$
makes sense for every object of $\mathfrak F_X$, and which satisfies 
the following conditions:

\begin{enumerate}[1]

\item  For every~$X\in \mathfrak  C$,
every~$x\in X$, every analytic domain~$V$ of~$X$ containing~$x$, and every
object $D$
of $\mathfrak F_X$,
the object~$D$ satisfies~$\mathsf Q$ at~$x$ if and only if~$D_V$ satisfies~$\mathsf Q$ at~$x$.

\item
For every~$X\in \mathfrak C$, every~$x\in X$, every analytic extension~$L$ of~$k$, 
every point~$x'$ of~$X_L$ lying above~$x$, and
every~$D\in \mathfrak F_X$,
the object~$D$ satisfies~$\mathsf Q$ at~$x$ if and only if~$D_L$ satisfies~$\mathsf Q$ at~$x'$
(we insist that one only needs to check $\mathsf Q$ at \emph{one}
pre-image of $x$ on $X_L$, and not at all of them). 

\end{enumerate}

\begin{exem}\label{ex-category-trick1}
We can take for $\mathfrak C$ the category of all analytic spaces, and for $\mathfrak F$
any fibered category as in \ref{ss-fiber-category}. 
Let~$\mathsf P$ be a property making sense for any object of~$\mathfrak F_{\mathfrak L}$;
assume moreover that~$\mathsf P$ satisfies \hreg~(\ref{ss-list-hreg}). 
By Remark \ref{rem-equiv-valid}
and Proposition \ref{prop-behavior-extension}, we
can take for~$\mathsf Q$ the property of satisfying \emph{geometrically}
the property~$\mathsf P$ 
at a point of an analytic space (\ref{ss-geom-valid}); if~$\mathsf P$ satisfies
the stronger condition~\hci, 
we do not need 
to require the validity to be geometric (Remark \ref{rem-hci-geo}). 
\end{exem}

\begin{exem}\label{ex-category-trick2}
Let $Y$ be an analytic space
and let $E$ be a locally constructible subset
of~$Y$. We can take for $\mathfrak C$ the category of all
analytic spaces over $Y$, for $\mathfrak F$ the category $\mathfrak T$
itself, and for $\mathsf Q$ the property defined
as follows: 
if~$Z$ is an object of $\mathfrak C$ and if~$z\in Z$,
then~$Z$ satisfies~$\mathsf Q$ at~$z$ if~$z$ belongs to the Zariski-closure of the pre-image of~$E$
under the structure map~$Z\to Y$.
\end{exem}

\subsection{}\label{ss-def-functorS}
We also fix a functor $\mathscr S$
from $\mathfrak F_{\mathfrak C}$
to $\mathfrak{Coh}_{\mathfrak C}$
which is compatible with the fibered structures over $\mathfrak C$; \ie, for every arrow $p\colon Y\to X$
in
$\mathfrak C$ and
every object $D\in \mathfrak F_X$, one has a canonical isomorphism $\mathscr S(p^\ast D)\simeq p^\ast \mathscr S(D)$.

\begin{exem}\label{ex-functorS}
Let
$\mathfrak I$ be an interval of~$\Z$
(viewed as a category) and assume that~$\mathfrak F=\mathfrak{Coh}^{\mathfrak I}$
(and let $\mathfrak C$ and $\mathsf Q$ be arbitrary). For every $Y\in \mathfrak C$, 
the fiber category $\mathfrak F_Y$
is the category 
of diagrams of~$\mathscr O_Y$-linear maps
$\ldots \to \mathscr F_i\to \mathscr F_{i+1}\to \mathscr F_{i+2}\to\ldots$ where~$i$
runs through~$\mathfrak I$ and where the~$\mathscr F_i$'s are 
coherent sheaves on~$Y$. We can then take
for~$\mathscr S$ the functor sending such a diagram to the~$i$-th coherent sheaf involved (for given~$i\in \mathfrak I$), or
more generally to the direct sum
$\bigoplus_{i\in J}\mathscr F_i$ for some finite subset~$J$ of~$I$. 
\end{exem}

\begin{enonce}[remark]{Notation}\label{notation-uflat}
Let $Y\to X$ be a morphism of $k$-analytic spaces
and let $\mathscr E$ be
 a coherent sheaf on $Y$. The $X$-flat locus of $\mathscr E$ 
 will be denoted by $\uflat {\mathscr E}X$. 
 \end{enonce}

\subsection{Fiberwise validity}
Let~$\phi \colon Y\to X$ be a morphism of~$k$-analytic spaces, with
$Y\in\mathfrak C$. We shall
say that a given object~$D$ of $\mathfrak F_Y$ satisfies~$\mathsf Q$
\emph{fiberwise}
at some
point $y\in Y$ if~$D_{Y_{\phi(y)}}$ satisfies~$\mathsf Q$ at~$y$. The set of points of~$Y$ 
at which~$D$ 
satisfies~$\mathsf Q$ fiberwise
will be denoted by~$\pfb DX$; if~$\mathscr F$ is a coherent sheaf on~$Y$, we shall denote 
by~$\uflat {\mathscr F} X$
the~$X$-flatness locus of~$\mathscr F$.

\subsection{Three statements}\label{ss-claims}
The main purpose of this section is to 
develop general methods for proving the three following statements: 

\begin{itemize}

\item[$(\alpha)$]
For every morphism~$Y\to X$ between~$k$-analytic spaces with $Y\in \mathfrak C$
and every
object $ D\in \mathfrak F_Y$, 
the set~$\pfb D X$ is a locally constructible subset of~$Y$.  

\item[$(\beta)$] 
For every morphism~$Y\to X$ between~$k$-analytic spaces 
with $Y\in \mathfrak C$  and every
object $D\in \mathfrak F_Y$, 
the intersection~$\pfb D X\cap \uflat {\mathscr S(D)} X$ 
is a Zariski-open subset of~$Y$.  

\item[$(\gamma)$]
For every morphism~$Y\to X$ between~$k$-analytic spaces
and each coherent sheaf~$\mathscr F$ on~$Y$, the
set~$\uflat {\mathscr F}X$ is a Zariski-open subset of~$Y$. 
\end{itemize}

\begin{rema}\label{rem-gamma-beta}
Assertion~$(\gamma)$ is nothing but
assertion~$(\beta)$ 
when we take for $\mathfrak F$ 
the category $\mathfrak {Coh}$, for $\mathfrak C$ the category of all analytic spaces,~for $\mathscr S$ the identity functor, and
for $\mathsf Q$ the
property {\sc true}. 
\end{rema}

\begin{rema}\label{rem-S-sumSi}
Suppose that we are given a finite family~$(\mathscr S_i)$ of~$\mathfrak C$-functors from~$\mathfrak F_{\mathfrak C}$ 
to~$\mathfrak{Coh}_{\mathfrak C}$, and that~$\mathscr S=\bigoplus \mathscr S_i$. Then the subset of~$Y$
that is involved in~$(\beta)$ is nothing but the intersection of~$\pfb DX$
with~$\bigcap_i \uflat{\mathscr S_i(D)}X$.
\end{rema}

\begin{rema}
We emphasize that in statements $(\alpha)$ and $(\beta)$, 
the space $Y$ is assumed to be an object of $\mathfrak C$, but not
the space $X$. And we will actually apply the results of this section
in some cases where $X\notin \mathfrak C$, for instance while working
in the situation of Example \ref{ex-category-trick2}. 
\end{rema}

\subsection{Some auxiliary statements}\label{ss-auxiliary}
We shall need some (apparently) weaker versions of $(\alpha)$, $(\beta)$ and $(\gamma)$, which
we are now going to list. 

\begin{itemize}
\item [$(\alpha')$]
For every morphism~$Y\to X$ between~$k$-affinoid
spaces with $Y\in \mathfrak C$, every section~$\sigma : X\to Y$,
and every $D\in \mathfrak F_Y$, 
the preimage~$\sigma^{-1}(\pfb DX)$ is a constructible
subset of~$X$. 

\item [$(\alpha'')$] For every morphism~$Y\to X$ between~$k$-affinoid
spaces with
\emph{integral}
$X$ and with $Y\in \mathfrak C$, 
every section~$\sigma \colon X\to Y$, 
and every object~$D$ of $\mathfrak F_Y$, 
either~$\sigma\inv(\pfb XD)$ or~$X\setminus \sigma\inv(\pfb XD)$
contains a non-empty Zariski-open subset of~$X$; \ie, $ \sigma\inv(\pfb XD)$
either contains or is disjoint from a non-empty Zariski-open subset of $X$. 

\item[$(\beta^\flat)$]
For every morphism~$Y\to X$ between~$k$-affinoid
spaces 
with $Y\in \mathfrak C$  and every
object $D\in \mathfrak F_Y$ such that $\mathscr S(D)$ is $X$-flat, 
$\pfb D X$ 
is a Zariski-open subset of~$Y$.

\item[$(\beta')$]
For every morphism~$Y\to X$ between~$k$-affinoid
spaces with $Y\in \mathfrak C$, every section~$\sigma \colon X\to Y$,
and every~$D\in \mathfrak F_Y$,
the pre-image~
$$\sigma\inv\left(\pfb DX\cap \uflat {\mathscr S(D)}X\right)$$
is a Zariski-open subset of~$X$. 

\item [$(\beta '')$] For every 
\emph{flat}
morphism~$Y\to X$ between~$k$-affinoid
spaces with $Y\in \mathfrak C$, every section~$\sigma \colon X\to Y$
and every
object $D\in \mathfrak F_Y$,
the pre-image~
$$\sigma\inv\left(\pfb DX\cap \uflat {\mathscr S(D)}X\right)$$
is a Zariski-open subset of~$X$.

\item [$(\gamma ')$] For every 
morphism~$Y\to X$ between~$k$-affinoid
spaces, every section
$\sigma \colon X\to Y$,
and
every coherent sheaf~$\mathscr F$ on~$Y$, the
pre-image
$\sigma^{-1}(\uflat {\mathscr F}X)$ is a Zariski-open subset of~$X$. 

\end{itemize}
We have the following hierarchy between our statements:

\begin{itemize}[label=$\bullet$]
\item $(\alpha')$ is a particular case of $(\alpha)$, and $(\alpha'')$ is a particular case of $(\alpha')$; 

\item $(\beta^\flat)$ and $(\beta')$ are particular cases of $(\beta)$, and $(\beta'')$ is a
particular case of $(\beta')$;

\item  $(\gamma')$ is a particular case of $(\gamma)$. 

\end{itemize}
\begin{rema}
\label{com-beta-prime}
 Assertion~$(\gamma')$ is nothing but
assertion~$(\beta')$ 
when we take for $\mathfrak F$ 
the category $\mathfrak {Coh}$, for $\mathfrak C$ the category of all analytic spaces,~for $\mathscr S$ the identity functor, and
for $\mathsf Q$ the
property {\sc true}. 

\end{rema}

Our purpose is now to explain how various combinations
of the above
``auxiliary statements"
(with possibly some extra
assumptions) imply $(\alpha)$, $(\beta)$, or $(\gamma)$
(see
Lemmas \ref{lem-alpha-prime} -- \ref{alpha-gamma-beta}).

\begin{lemm}\label{lem-diag-trick}
Let $Y\to X$ be a morphism of $k$-affinoid spaces
with $Y\in \mathfrak C$. 
Let~$p_1$ and~$p_2$ be the two projections from~$Y\times_XY$ to~$Y$, and let
$\sigma$
be the diagonal immersion $Y\hookrightarrow Y\times_X Y$. Let $D$ be an object of $\mathfrak F_Y$,
let $y$ be a point of $Y$, and let $\mathscr F$ be a coherent sheaf on $Y$. 

\begin{enumerate}[1]
\item The following are equivalent:

\begin{enumerate}[j]
\item $D$ satisfies $\mathsf Q$ fiberwise at $y$; 

\item $p_1^\ast D$ satisfies $\mathsf Q$ fiberwise at $\sigma(y)$ with respect to $p_2$. 
\end{enumerate}

\item If moreover $Y$ is flat over $X$, the following are equivalent: 
\begin{enumerate}[j]\setcounter{enumii}{2}
\item $\mathscr F$ is $X$-flat at $y$; 

\item $p_1^\ast \mathscr F$ is $Y$-flat at $\sigma(y)$ with respect to $p_2$.
\end{enumerate}
\end{enumerate}
\end{lemm}

\begin{proof}
Let $x$ denote the image of $y$ on $X$. We have $p_2(\sigma(y))=y$, and $p_1$
induces an isomorphism $p_2\inv(y)\simeq Y_x\times_{\hr x}\hr y$, which sends $\sigma(y)$ to $y$.
The equivalence
(i)$\iff$(ii) thus follows from the good behavior of $\mathsf Q$ with respect to ground field extension, see axiom (2) in \ref{ss-general-axioms-pfib}. 

The implication (iii)$\Rightarrow$(iv) is true
without
the flatness assumption on $Y$ because flatness at a point is by definition preserved by any base change.
If $Y$ is flat over $X$,
the implication (iv)$\Rightarrow$(ii) follows from Proposition \ref{prop-flatbc}. 
\end{proof}

\begin{lemm}\label{lem-flatloc-nonempty}
Let $Y\to X$ be a morphism between $k$-affinoid spaces with $X$
integral, and let $\sigma\colon X\to Y$
be a section of $Y\to X$. 
Let $\mathscr F$ be any coherent sheaf on $Y$. 
The pre-image $\sigma\inv(\uflat {\mathscr F}X)$ contains each pre-image of the generic point of 
$X\al$, and is in particular non-empty. 
\end{lemm}

\begin{proof}
Let $x$ be a point of $X$ such that $x\al$ is the generic point of $X\al$. 
Since $\mathscr O_{X\al, x\al}$ is a field, $\mathscr F\al$ is $X\al$-flat at 
$\sigma(x)\al$. The closed immersion $\sigma$ makes $\sigma(X)$ a closed analytic
subspace of $Y$, which is finite over $X$ (the map $\sigma(X)\to X$ is even an isomorphism). 
It follows then from Theorem \ref{thm-flat-improve}
that $\mathscr F$ is $X$-flat at $\sigma(x)$. As a consequence, $x\in \sigma\inv(\uflat {\mathscr F}X)$. 
\end{proof}

\begin{lemm}\label{lem-alpha-prime}
Assume that $(\alpha')$ holds. Then $(\alpha)$ holds.
\end{lemm}

\begin{proof}
Let
$Y,X$ and $D$ be as in $(\alpha)$; since $(\alpha)$ is G-local, we can assume
that $Y$ and $X$ are affinoid. 
Let~$p_1$ and~$p_2$ be the two projections from~$Y\times_XY$ to~$Y$, and let
$\sigma$
be the diagonal immersion $Y\hookrightarrow Y\times_X Y$. 
By Lemma \ref{lem-diag-trick} (1), one has 
\[\pfb DX=\sigma\inv(\pfb {p_1^\ast D}Y),\]
where fiberwise validity on
the right hand side has to be understood with respect to $p_2$. By applying
$(\alpha')$ (which holds by assumption)
to $(Y\times_XY\stackrel{p_2}\longrightarrow Y, p_1^\ast D,\sigma)$, we see that
$\sigma\inv(\pfb {p_1^\ast D}Y)$ is a constructible subset of $Y$. Therefore
$ \pfb DX$
is a constructible subset of $Y$. 
\end{proof}

\begin{lemm}
\label{alphadash-alpha}
Assume that $(\alpha'')$ hold. Then $(\alpha)$ holds.
\end{lemm}

\begin{proof}
We shall prove that $(\alpha')$ holds, 
which will imply
that $(\alpha)$ holds in view of Lemma \ref{lem-alpha-prime}. 
Let $Y, X$ and $D$ be as in $(\alpha')$.
The affinoid space $X$ is quasi-compact for the constructible topology
(\ref{ss-constop-compact}). It thus suffices to prove
that for
every
point $x$ of $X$, there
exists a constructible subset of~$X$
containing~$x$
which is either included in~$\sigma\inv (\pfb DX)$
or in
its complement~$X\setminus \sigma\inv(\pfb DX)$. 

So, let~$x$ be a point of $X$. Since
we are interested in a fiberwise property, we may replace~$X$ with the reduced
Zariski closure of~$\{x\}$; hence we can assume that~$X$ is integral
and that~$x\al$ is the generic point of $X\al$.
As $(\alpha'')$ holds, either
$\sigma \inv(\pfb DX)$ or its complement contains a non-empty Zariski-open subset of $X$, 
which is a constructible subset of $X$ containing $x$.
\end{proof}

\begin{lemm}
\label{ss-beta-gamma-alphadash}
Assume that $(\beta)$ and $(\gamma)$ hold. Then $(\alpha)$ holds. 
\end{lemm}

\begin{proof}
We shall prove that $(\alpha'')$ holds, 
which will imply
that $(\alpha)$ holds in view of Lemma \ref{alphadash-alpha}. Let~$Y\to X$ be
as in~$(\alpha'')$. 
Since $(\beta)$ and~$(\gamma)$ are assumed to hold, the sets~$\pfb DX\cap \uflat{\mathscr S(D)}X$ and~$\uflat{\mathscr S(D)}X$
are Zariski-open subsets of $Y$. It follows from Lemma \ref{lem-flatloc-nonempty}
(applied to $\mathscr F=\mathscr S(D)$)
that
the Zariski-open subset $\sigma\inv(\uflat{\mathscr S(D)}X)$
of $X$
is non-empty.  
We now distinguish two cases: 

\begin{itemize}[label=$\bullet$]

\item If~$\sigma\inv(\pfb DX\cap \uflat{\mathscr S(D)}X)\neq \emptyset$, this is a non-empty Zariski-open
subset in~$X$ which is contained in~$\sigma\inv(\pfb DX)$, and we are done.

\item If
$\sigma\inv(\pfb DX\cap \uflat{\mathscr S(D)}X)=\emptyset$, 
then~$\sigma\inv(\uflat{\mathscr S(D)}X)$ is a non-empty Zariski-open subset of~$X$
which is contained in~$X\setminus \sigma\inv(\pfb DX)$, and we are done.
\end{itemize}
\end{proof}

\begin{lemm}\label{lem-betadash-beta}
We make the following assumptions:

\begin{enumerate}[a]
\item The fibered category $\mathfrak F$ is equal to $\mathfrak{Coh}^{\mathfrak I}$ for some
small category $\mathfrak I$
(this includes the cases~$\mathfrak F=\mathfrak T$
and~$\mathfrak F=\mathfrak{Coh}$;
\cf Remark \ref{rem-diag-includes}). 

\item The functor~$\mathscr S$
commutes (as a functor between fibered categories; \ie, incorporating
a natural compatibility with pullback isomorphisms)
with push-forwards by closed immersions. 

\item For
every~$D\in \mathfrak F_Y$ and every closed immersion~$\iota
\colon Y\hookrightarrow Z$
of~$X$-analytic spaces, $D$ satisfies fiberwise~$\mathsf Q$ at
a given point~$y\in Y$ if and only if~$\iota_\ast  D$ satisfies~$\mathsf Q$ 
fiberwise
at~$\iota(y)$. 

\end{enumerate}

Assume moreover that $(\beta'')$ holds. Then $(\beta)$ holds. 
\end{lemm}

\begin{proof}
Let $Y,X$, and $D$ be as in $(\beta)$; since $(\beta)$ is G-local, 
we can assume that $Y$ and $X$ are affinoid. 
The morphism~$Y\to X$ then factorizes through a closed
immersion~$\iota \colon Y\hookrightarrow \Delta\times_k X$ for
a suitable compact~$k$-polydisc~$\Delta$.
Set
$$U=\pfb {\Delta\times_k X}X\cap \uflat {\mathscr S(\iota_\ast D)}X.$$ 
By the assumptions of the lemma, $U\cap Y=\pfb D X\cap \uflat {\mathscr S(D)}X$; it is therefore
sufficient to prove that $U$ is a Zariski-open subset of $\Delta\times_k X$. Hence by replacing
$Y$ with $\Delta\times_k X$ and $D$ with $\iota_\ast D$, we reduce to the case where $Y$ is $X$-flat.

Let~$p_1$ and~$p_2$ be the two projections from~$Y\times_XY$ to~$Y$, and let
$\sigma$
be the diagonal immersion $Y\hookrightarrow Y\times_X Y$.
By Lemma \ref{lem-diag-trick} (note that flatness of $Y$ over $X$ is needed to apply statement (2) of \loccit), one has 
\[\pfb DX=\sigma\inv(\pfb {p_1^\ast D}Y)\]
and
\[\uflat {\mathscr S(D)}X=\sigma \inv(\uflat{p_1^\ast \mathscr S(D)}Y)=
\sigma \inv(\uflat{\mathscr S(p_1^\ast D)}Y),\]
where fiberwise validity and flatness over $Y$ are
understood to be
with respect to $p_2$. 
By applying
$(\beta'')$ (which holds by assumption)
to $(Y\times_XY\stackrel{p_2}\longrightarrow Y, p_1^\ast D,\sigma)$, we see that
\[\sigma\inv(\pfb {p_1^\ast D}Y\cap \uflat{\mathscr S(p_1^\ast D)}Y)\]
is a Zariski-open subset of $Y$ (note that we use once again the flatness of $Y$ over $X$, because this is one of the assumptions of $(\beta'')$); therefore $\pfb DX\cap \uflat {\mathscr S(D)}X$ is 
a Zariski-open subset of $Y$. 
\end{proof}

\begin{rema}\label{rem-redstat-fulfilled}
If we take for $\mathfrak F$ 
the category $\mathfrak {Coh}$, for $\mathfrak C$ the category of all analytic spaces,
for $\mathscr S$ the identity functor, and for $\mathsf Q$ the
property of being CM, or $S_m$ for some specified $m$, or of
a
given codepth, or more simply the property {\sc true}, 
then the assumptions (a), (b) and (c) of Lemma
\ref{lem-betadash-beta}
above are fulfilled. 
\end{rema}

\begin{lemm}\label{gamma-gammaprime}
Assume that $(\gamma')$ holds. Then $(\gamma)$ holds. 
\end{lemm}

\begin{proof}
Take for $\mathfrak F$ 
the category $\mathfrak {Coh}$, for $\mathfrak C$ the category of all analytic spaces,~for $\mathscr S$ the identity functor, and
for $\mathsf Q$ the
property {\sc true}. By Lemma
\ref{lem-betadash-beta}
and Remark \ref{rem-redstat-fulfilled}, we
have $(\beta'')\Rightarrow(\beta)$, and thus also $(\beta')\Rightarrow(\beta)$ since
$(\beta')$ is stronger that $(\beta'')$.
But in our context, $(\beta)=(\gamma)$
and $(\beta')=(\gamma')$, whence the claim. 
\end{proof}

\begin{lemm}
\label{gamma-betaprime-beta}
Assume that $(\gamma)$ and $(\beta')$ hold. Then $(\alpha)$
holds.
\end{lemm}

\begin{proof}
We shall prove that $(\beta)$ holds, 
which will imply
that $(\alpha)$ holds in view of Lemma \ref{ss-beta-gamma-alphadash}.
Let
$Y,X$, and $D$ be as in $(\beta)$; since $(\beta)$ is G-local, we can assume
that $Y$ and $X$ are affinoid. 
Let~$p_1$ and~$p_2$ be the two projections from~$Y\times_XY$ to~$Y$, and let
$\sigma$
be the diagonal immersion $Y\hookrightarrow Y\times_X Y$.
Set
\[E=\pfb {p_1^\ast D}Y\cap \uflat{p_1^\ast\mathscr S(D)}Y=\pfb {p_1^\ast D}Y\cap \uflat{\mathscr S(p_1^\ast D)}Y,\]
where fiberwise validity and flatness over $Y$ are understood
to be with respect to $p_2$. 
Since $\sigma(\uflat {\mathscr S (D)}X)) \subset \uflat {p_1^\ast \mathscr S(D))}Y$, we have
$$\pfb DX\cap \uflat {\mathscr S(D)}X=\sigma\inv (E)\cap \uflat{\mathscr S(D)}X,$$
due to Lemma \ref{lem-diag-trick}.
Since we assume that~$(\gamma)$
holds, the set~$\uflat{\mathscr S (D)}X$ is Zariski-open, and
it is thus sufficient to prove that~$\sigma\inv(E)$ is Zariski-open.
But the latter follows by applying $(\beta')$ (which holds by assumption)
to the list of data $(Y\times_XY\stackrel{p_2}\longrightarrow Y, p_1^\ast D,\sigma)$.
\end{proof}

\begin{lemm}
\label{alpha-gamma-beta}
Assume that $(\alpha)$, $(\beta^\flat)$, and $(\gamma)$ hold. Then $(\beta)$
holds.
\end{lemm}

\begin{proof}
Let $Y$, $X$, and $D$ be as in $(\beta)$; since $(\beta)$ is G-local, 
we can assume that $Y$ and $X$ are affinoid. 
Since $(\alpha)$ and $(\gamma)$ hold,~$\pfb DX \cap \uflat{\mathscr S(D)} X$ is a constructible subset of $Y$. 
In order to prove~$(\beta)$, it thus suffices
to prove that
the intersection
$\pfb DX\cap \uflat{\mathscr S(D)}X$
is open
(Lemma~\ref{lem-adh-const} (4)). Let~$y$
be a point of~$\pfb DX\cap \uflat{\mathscr S(D)}X$. Since~$y\in \uflat{\mathscr S(D)}X$, 
it follows from~$(\gamma)$ that there exists an affinoid
neighborhood~$V$ of~$y$ in~$Y$ which is included in~$\uflat{\mathscr S(D)}X$. 
Since we assume that $(\beta^\flat)$ holds, the set $\pfb {D_V}X=\pfb DX\cap V$ is a Zariski-open subset of $V$; 
in particular, it contains a neighborhood of $y$. 
\end{proof}

\section{The flat locus}\label{s-flatloc}

The main theorem of this section (Theorem \ref{thm-flatloc-zaropen})
says the following: if $Y\to X$ is a morphism of $k$-analytic spaces
and $\mathscr F$ is a coherent sheaf on $Y$, the $X$-flat locus of $\mathscr F$
is a Zariski-open subset of $Y$; our proof follows Kiehl's strategy developed
in \cite{kiehl1967}. 

Zariski-openness of the flat locus has the following consequences, in view of the Nullstellensatz: if $\abs {k\gpm}\neq \{1\}$
and if $Y$ and $X$ are strict, then $\mathscr F$ is $X$-flat if and only if it is $X$-flat at every rigid point of $Y$; and
$\mathscr F$ is $X$-flat at every point lying over a given point $x$
of $X$ if and only if it is $X$-flat at every rigid point of $Y_x$. 

We use the first consequence
together with Theorem \ref{thm-flat-improve}
to get the compatibility between our notion
of flatness and that of rigid flatness (Corollary \ref{cor-gagaflat-rigid}) as well 
as that of formal flatness (Corollary \ref{cor-flat-formal}).

We use the second one
to prove that flatness
holds automatically over any Abhyankar point of the target space -- provided the latter
is reduced (Theorem \ref{thm-flat-abhyankar}). Let us make this precise. 
Once having reduced to the strict
good case in a standard way, we only have to check 
flatness at rigid points of the fiber under investigation (by the above). We then use a ``deboundarization" result
(Lemma \ref{lem-deboundarize}) for such a rigid point, which rests on the notion of smallest analytic domain
containing the image of a morphism of analytic germs (Theorem \ref{thm-smallest-germ});
this enables us to reduce to the inner case, 
for
which it suffices to check naive flatness (Theorem \ref{thm-flat-naiveflat}). But the latter holds for free
because the local ring of an Abhyankar point of a good analytic space is artinian
(Example \ref{ex-centdim-recap}), hence a field whenever the space is reduced.

%\begin{enonce}[remark]{Notation}\label{notation-uflat}
%If $Y\to X$ is a morphism of schemes (\resp of $k$-analytic spaces) and if $\mathscr E$ is a coherent sheaf on $Y$, we denote
%by $\uflat {\mathscr E}X$ the set of points of $Y$ at which $\mathscr E$ is $X$-flat.
%\end{enonce}

For the reader's convenience, we state and prove 
the following theorem
of Kiehl, which will be crucial for our description of the flat locus. 

\begin{theo}[Kiehl, \cite{kiehl1967} Satz 1]\label{thm-kiehl-flat}
Let $Y\to X$ be a morphism of noetherian schemes, let $\mathscr E$ be a
coherent sheaf on $Y$, and let $Z$ be a closed subscheme of $Y$
of finite type over $X$. The intersection
\[Z\cap \uflat{\mathscr E}X\]
is a Zariski-open subset of $Z$.

\end{theo}

\begin{proof}
If $Y$ is itself of finite type over $X$ and if $Z=Y$, this is \cite{ega43}, \Th 11.1.1.
For the general case, one can follow \emph{mutatis mutandis}
the proof of \loccit, except that classical ``generic flatness" (\cite{ega42}, Lemme 6.9.2) has to be replaced
with the following
stronger statement: \emph{let $A\to B$ be a morphism
of noetherian rings with $A$ a \emph{domain}, let $M$ be a finite $B$-module, and let $J$ be an ideal of $B$
such that the $A$-algebra $B/J$ is finitely generated;
there exists $a\neq 0$ in $A$ such that $M_{\mathfrak p}$ is $A$-flat
for every prime ideal $\mathfrak p$ of $B$ with $J\subset \mathfrak p$ and $a\notin \mathfrak p$.}

Let us prove this claim. The ring 
$C:=\bigoplus_n J^n/J^{n+1}$ is finitely generated over $B/J$ (since $J$
is finitely generated by noetherianity of $B$), hence it is also
finitely generated over $A$. By classical generic flatness, there exists $a\neq 0$ in $A$ such that $(C\otimes_B M)_a$ is $A$-flat. 
Let $\mathfrak p$ be a prime ideal of $B$ with $J\subset \mathfrak p$ and $a\notin \mathfrak p$. Since $a\notin \mathfrak p$, 
the $B$-module $(C\otimes_B M)_{\mathfrak p}$ is a localization of $(C\otimes_B M)_a$, hence is $A$-flat. This means
that $ J^nM_{\mathfrak p}/J^{n+1}M_{\mathfrak p}$ is $A$-flat for every $n$. Therefore $M_{\mathfrak p}/J^nM_{\mathfrak p}$ 
is $A$-flat for every $n$. Since $J\subset \mathfrak p$ this implies by \cite{ega31}, Chapitre 0, \Prop 10.2.6 that
the $A$-module $M_{\mathfrak p}$ is flat. 
\end{proof}

\begin{theo}\label{thm-flatloc-zaropen}
Let~$Y\to X$ be a morphism of~$k$-analytic spaces and let~$\mathscr E$ be a coherent sheaf on~$Y$. The subset
$\uflat{\mathscr E}Y$
of $Y$
is Zariski-open.
\end{theo}

\begin{proof}
We want to prove assertion~$(\gamma)$ of~\ref{ss-claims}. By Lemma \ref{gamma-gammaprime}, 
it is sufficient to prove assertion~$(\gamma')$ of~\ref{ss-auxiliary}. We thus may assume that~$Y$ and~$X$
are~$k$-affinoid and $Y\to X$ has a section~$\sigma$, and it is then
sufficient to prove that~$\sigma^{-1}(\uflat{\mathscr E}X)$
is a Zariski-open subset of~$X$; or, what amounts to the same, that~$Z\cap \uflat{\mathscr E}X$
is a Zariski-open subset of~$Z$, where~$Z$ is the closed analytic subspace of~$Y$ defined by
the closed immersion~$\sigma$. 

By Theorem \ref{thm-flat-improve}, the intersection~$\uflat{\mathscr E}X\cap Z$ is nothing but the pre-image
of~$\uflat{\mathscr E\al}{X\al}\cap Z\al$ in $Z$. 
Both~$Y\al$ and~$X\al$ are noetherian schemes, and~$Z\al$ is a Zariski-closed subscheme of~$Y\al$
that is of finite type over~$X\al$. By Kiehl's theorem stated above
(Theorem \ref{thm-kiehl-flat}), $\uflat{\mathscr E\al}{X\al}\cap Z\al$ is a Zariski-open subset of~$Z\al$;
as a consequence, $\uflat{\mathscr E}X\cap Z$ is a Zariski open subset of $Z$. 
\end{proof}

Due to this theorem we recover the fact that in the rigid setting, global algebraic flatness implies global analytic flatness:

\begin{coro}\label{cor-gagaflat-rigid}
Let~$Y\to X$ be a morphism of \emph{strict}
$k$-affinoid spaces. Let
$\mathscr E$ be a coherent sheaf on~$Y$. Then~$\mathscr E$ is~$X$-flat if and only if~$\mathscr E(Y)$ is a flat $\mathscr O_X(X)$-module.
\end{coro}

\begin{proof}
If~$\mathscr E$ is~$X$-flat, it is in particular naively~$X$-flat, and Lemma
\ref{lem-gagaflat-easy}
then ensures that~$\mathscr E(Y)$ is flat over $\mathscr O_X(X)$.
Conversely, assume that $\mathscr E(Y)$ is flat over $\mathscr O_X(X)$.
Then by Theorem \ref{thm-flat-improve}, $\mathscr E$ is~$X$-flat at any rigid point of~$Y$.
Since $\uflat{\mathscr E}X$ is a Zariski-open subset of $Y$ by Theorem \ref{thm-flatloc-zaropen},  
it follows from the analytic (\resp algebraic) Nullstellensatz if~$\abs{k\gpm}\neq \{1\}$ (\resp~$\abs{k\gpm}= \{1\}$), that $\uflat {\mathscr E}X=Y$. 
\end{proof}

\begin{rema}
Corollary~\ref{cor-gagaflat-rigid}
above is false in general without any strictness assumption; for a counter-example, 
see~\ref{ss-counterex-algflat}. 
\end{rema}

\begin{coro}\label{cor-flat-formal}
Let~$\mathfrak Y\to \mathfrak X$ be a
morphism between topologically finitely presented~$\mathrm{Spf}\;k^\circ$-formal schemes,
and let~$\mathscr E$ be a coherent sheaf on~$\mathfrak Y$ which is~$\mathfrak X$-flat.
The associated coherent sheaf~$\mathscr E_\eta$ on~$\mathfrak Y_\eta$ is~$\mathfrak X_\eta$-flat.
\end{coro}

\begin{proof}
We can assume
that both~$\mathfrak Y$ and~$\mathfrak X$ are affine formal schemes.  By assumption,~$\mathscr E(\mathfrak Y)$ is
a flat~$\mathscr O_{\mathfrak X}(\mathfrak X)$-module; therefore,~$\mathscr E_\eta(\mathfrak Y_\eta)=\mathscr E(\mathfrak Y)\otimes_{k^\circ}k$ is
flat over $\mathscr O_{\mathfrak X_\eta}(\mathfrak X_\eta)=\mathscr O_{\mathfrak X}(\mathfrak X)\otimes_{k^\circ}k$.
In view of the preceding corollary, this implies that~$\mathscr E_\eta$ is~$\mathfrak X_\eta$-flat.
\end{proof}

Our purpose is now to prove that flatness holds
automatically over Abhyankar points
(\ref{ss-abhyankar-points}) of reduced spaces.
This is
Theorem \ref{thm-flat-abhyankar}
below, which witnesses the fact
that the best analytic analogue of a scheme-theoretic generic fiber 
is a fiber over an Abhyankar point of a reduced space.

Theorem  \ref{thm-flat-abhyankar}
rests  on the following ``deboundarization" lemma, which
is of independent interest. 

\begin{lemm}\label{lem-deboundarize}
Let~$Y\to X$ be a morphism of
separated~$\Gamma$-strict~$k$-analytic spaces, let~$y$
be a point of $Y$
and let~$x$ be its image in~$X$. If~$\hr y$ is finite over~$\hr x$, there exists
a~$\Gamma$-strict $k$-analytic space~$X'$, a 
quasi-\'etale
morphism~$X'\to X$, and a pre-image~$y'$ of~$y$ on~$Y':=Y\times_X X'$
such that~$Y'\to X'$ is inner at~$y'$.
\end{lemm}

\begin{proof}
Let~$L$ be the separable closure
of~$\hr x$ inside~$\hr y $, and let~$Z$ be a 
finite \'etale
cover of an analytic neighborhood of~$x$
in~$X$ such 
that~$x$ has one pre-image~$x'$ on~$Z$ with~$\hr x'\simeq L$
(such a $Z$ exists by \cite{berkovich1993}, \Th 3.4.1). 

There is a canonical point~$y'$ of~$T:=Y\times_X Z$ lying above~$y$ 
and~$x'$ 
such that~$\hr {y'}$ is a finite radicial 
extension of~$\hr {x'}$, which implies
that $\hrt {y'}$ is a finite radicial
extension of $\hrt {x'}$ (\ref{ss-gradval-res}).

Let~$(X',x')$ be the smallest analytic domain of~$(Z,x')$ through
which~$(T,y')$
factorizes
(Theorem \ref{thm-smallest-germ}); it is~$\Gamma$-strict
by \loccit, 
and the quasi-compact open
subset~$\widetilde{(X',x')}^\Gamma$ of~$\P_{\hrt{x'}^\Gamma/\widetilde k^\Gamma}$
is the image of~$\widetilde{(T,y')}^\Gamma$ 
under~$\P_{\hrt{y'}^\Gamma/\widetilde k^\Gamma}\to \P_{\hrt{x'}^\Gamma/\widetilde k^\Gamma}$,
which is a homeomorphism 
because~$\hrt{y'}^\Gamma$ is radicial
over~$\hrt{x'}^\Gamma$ (Remark \ref{rem-zr-radicial}).
Hence 
the inverse image
of~$\widetilde{(X',x')}^\Gamma$ in~$\P_{\hrt{y'}^\Gamma/\widetilde k^\Gamma}$
is exactly~$\widetilde{(T,y')}^\Gamma$. This 
implies by
the criterion \ref{ss-prop-gammared} (2)
that
the map $(T,y')\to (X',x')$
is boundaryless, and the~$X$-analytic space~$X'$ satisfies the required 
conditions.
\end{proof}

\begin{theo}\label{thm-flat-abhyankar}
Let~$Y\to X$ be a morphism of~$k$-analytic spaces, with~$X$ reduced.
Let~$x$ be an Abhyankar point of $X$ (\ref{ss-abhyankar-points}),
and let~$\mathscr E$
be a coherent sheaf on~$Y$. The sheaf~$\mathscr E$ is~$X$-flat above~$x$.
\end{theo}

\begin{proof}
We set $n=d_k(x)=\dim_xX$.
We may assume that~$Y$ and~$X$ are~$k$-affinoid. 
Let~$r$ be a~$k$-free polyradius such that~$\abs {k_r\gpm}\neq \{1\}$
and
such that both~$Y_r$ and~$X_r$ are strictly~$k_r$-affinoid. 
The~$k_r$-analytic space~$X_r$
is~$n$-dimensional and is reduced (as follows from the fact that~$k_r$ is
analytically separable over~$k$, or in a more
elementary way from
Lemma \ref{lem-xxr-scheme}
applied with $\mathscr X=\spec A$). If~$\mathfrak s : X\to X_r$
denotes the Shilov section, one has~$d_{k_r}(\mathfrak s(x))=d_k(x)=n$ (\ref{ss-dk-shilov}). 
Hence, due to Proposition \ref{prop-flat-xl},
we can replace
the field~$k$ by~$k_r$, the spaces~$X$ and~$Y$ by~$X_r$ and~$Y_r$,
the sheaf~$\mathscr E$ by $\mathscr E_{Y_r}$, 
and
the point~$x$ by~$\mathfrak s(x)$; \ie, we can assume that~$\abs {k\gpm}\neq \{1\}$ and that~$Y$ and~$X$
are strictly~$k$-affinoid.

Theorem
\ref{thm-flatloc-zaropen}
ensures that $\uflat {\mathscr E}X$ is a Zariski-open subset of $Y$. 
We want to prove that it contains $Y_x$. 
Since $Y_x$ is strictly affinoid over the non-trivially valued field~$\hr x$,
it is sufficient to prove that it contains every rigid point of~$Y_x$. 

Let~$y$ be a rigid point of~$Y_x$. By 
Lemma~\ref{lem-deboundarize}
above, there exists a strictly~$k$-analytic space~$X'$, a quasi-\'etale map~$X'\to X$,
and a pre-image~$y'$ of~$y$ in~$Y':=Y\times_X X'$ such that~$Y'\to X'$ is
inner at~$y'$; we can assume (by shrinking it if necessary) that~$X'$ is strictly~$k$-affinoid, 
in which case so is~$Y'$. 

Since~$X'$ is quasi-\'etale over the reduced analytic space $X$, 
it is reduced by Proposition \ref{prop-tranferqsm-concrete}; 
moreover, it is of dimension~$\leq n$ by \ref{ss-dim-sobafi}. Since $\hr {x'}$
is finite over $\hr x$, we have $d_k(x')=n$ (hence $\dim_{x'} X'=n$). It follows from Example
\ref{ex-centdim-abhyankar}
that $\mathrm{centdim}(X',x')=n$, which implies in view of Corollary 
\ref{cor-interp-centdim}
that~$\mathscr O_{X',x'}$ is artinian, hence a field because it is also reduced. 

Therefore~$\mathscr E_{Y',y'}$ is flat over~$\mathscr O_{X',x'}$. As
$(Y',y')\to (X',x')$
is boundaryless, this implies
in view of Theorem \ref{thm-flat-naiveflat}
that $\mathscr E_{Y'}$
is~$X'$-flat 
at~$y'$. Since~$X'\to X$ is quasi-\'etale, it is flat
(Corollary \ref{cor-qsm-flat}). Therefore~$\mathscr E$
is~$X$-flat at~$y$ by Proposition \ref{prop-flatbc}; 
otherwise said, $y\in \uflat {\mathscr E}X$.
\end{proof}

\section{The fiberwise closure of a locally constructible set}\label{s-fib-closure}

In Section \ref{s-cons-maintheo}, we shall reduce our study of the locus of fiberwise validity of the $S_n$ property 
(\resp geometric $R_n$ property) to that of fiberwise validity of the CM property
(\resp the quasi-smooth) property, through a result that expresses the former in terms
of the latter and of some codimension considerations (Lemma \ref{lem-before-main}). 

We thus
have to investigate, being given
a morphism $\phi \colon Y\to X$ of
$k$-analytic spaces
and a locally constructible subset $E$ of $Y$ such that $E_x$ is a Zariski-closed
subset of $Y_x$ for every $x\in X$, 
the \emph{fiberwise
codimension function} $y\mapsto \codim_y (E_{\phi(y)}, Y_{\phi(y)})$.
For that purpose we prove the following results in this section, with $\phi \colon Y\to X$
as above. 

\begin{enumerate}[1]

\item Let $E$ be a locally constructible subset of $Y$, and set $\adht E\phi=\bigcup_{x\in X}\adht {E_x} {Y_x}$.
Then $\adht E\phi$ is a locally constructible subset of $Y$ (Theorem \ref{thm-fiber-closure}). 

\item For every non-negative integer $d$, the set of points $y\in \adht E\phi$ such that $\dim_y (\adht E\phi)_{\phi(y)}=d$
is locally constructible (Theorem \ref{thm-fiber-closure}). 

\item Let $F$ be a locally constructible subset of $Y$ contained in $E$. For every non-negative integer $d$,
the set of points $y\in \adht E\phi$ such that $\codim_y (\adht F\phi,\adht E\phi)_{\phi(y)}=d$
is locally constructible (Proposition \ref{prop-relat-codim}). 
\end{enumerate}
Note that the situation adressed by these statements might look slightly more general than the one we
have described as
motivation,
which involves only the fiberwise codimension in $Y$ of a locally constructible subset which is fiberwise
Zariski-closed,
but in fact even for proving (3) in this particular
case, we do not see how to avoid using (1) and (2) in their full generality
(and we think moreover that they are of independent interest). 

Let us mention that our results are quite analogous to the ones proved in \cite{ega43} 9.5
but our proofs are more involved, because the Zariski topology of a fiber is in general
strictly finer than the one induced by the global Zariski topology. The key tool to bypass
this problem is Theorem \ref{thm-localring-generic}
about the local rings of generic fibers.

\subsection{}\label{def-ebar-phi}
Let~$\phi : Y\to X$ be a morphism of~$k$-analytic spaces, and let~$E$ be a
locally constructible subset of~$Y$. We set
\[\adht E \phi=\bigcup_{x\in X}\adht{E\cap Y_x}{Y_x}=\bigcup_{x\in X}\adht{E\cap Y_x}{Y_{x,\mathrm {Zar}}}\]
(the second equality comes from Corollary \ref{cor-adh-const}). \label{IN-Ephi}
We say that $\adht E \phi$ is the
\emph{fiberwise closure of~$E$ with respect to~$\phi$ (or over $X$)}. 
We obviously have~$\adht E \phi\subset \adht E Y$. 

Let~$y$
be a point of $\adht E \phi$ and let~$x$ be its image in~$X$. The
dimension of the Zariski-closed subset~$\adht E \phi _x$ of~$Y_x$ 
at~$y$ is called the \emph
{relative dimension of~$\adht E\phi$ at~$y$.}

We shall prove
below that $\adht E\phi$ is locally constructible (Theorem \ref{thm-fiber-closure}).
For that purpose, we need the following technical proposition. 

\begin{prop}\label{prop-fiber-closure}
Let $\phi \colon Y\to X$ be a morphism between $k$-affinoid spaces, 
let $Z$ be a
closed analytic subspace of $Y$ with 
dense
(Zariski-open) complement $U$, and let $t$ be a point of $Z$. 
Assume that $t\in \mathrm{Int}(Y/X)$ and both $Y$ and
$Z$ are $X$-flat at $t$. The point $t$ then
belongs to $\adht U\phi$.
\end{prop}

\begin{proof}
Let~$s$ be the image of $t$ in $X$. 
We can perform any base change
consisting
of replacing~$X$ with an affinoid neighborhood of~$s$
without modifying
our assumptions
(as far as density of $U$ in $Y$ is concerned, this
rests on Corollary  \ref{cor-zarclosure-dom}). 
We can
therefore assume that $\mathfrak m_s$ is generated by an ideal of 
$\mathscr O_X(X)$, and we denote by $S$ the corresponding
closed analytic subspace of $X$; we set 
$T=Y\times_X S$. 
The local ring $\mathscr O_{S,s}=\mathscr O_{X,s}/\mathfrak m_s\mathscr O_{X,s}$ 
is then a field, and $\mathscr O_{S\al, s\al}$
is thus also a field by \ref{ss-gaga-domain}. 
We use the conventions of \ref{ss-conventions-gaga} 
(3); be aware
that $t_s\al:=(t_s)\al$ is a point of the scheme
$T_s\al:=(T_s)\al$, while $t\al$ is a point of $T\al$
lying in the fiber $T\al_{s\al}:=(T\al)_{s\al}$ of $T\al$
over the point $s\al$ of $S\al$. We are now going to prove the two
following statements. 

\begin{enumerate}[1]
\item The point
$t\al$ belongs to the closure of~$U\al_{s\al}$ in~$Y\al_{s\al}$.

\item The scheme $Y_s\al$ is flat
over $Y\al_{s\al}$ at $t_s\al$.
\end{enumerate}

Let us prove (1) by contradiction. Assume
that it does not hold.  The
(support of) the closed subscheme $Z\al_{s\al}$ of $Y\al_{s\al}$ is
then a neighborhood of $t\al$
in $Y\al_{s\al}$, and
we thus have
\begin{equation}\label{d-krull-notadh}
\dim_{\mathrm{Krull}}\mathscr O_{Y\al_{s\all}, t\al}=
\dim_{\mathrm{Krull}}\mathscr O_{Z\al_{s\all}, t\al}.\end{equation}

By assumption, the affinoid spaces $Y$ and~$Z$ are~$X$-flat at~$t$; therefore
the schemes $Y\al$ and~$Z\al$ are
$X\al$-flat
at~$t\al$ (Lemma \ref{lem-gagaflat-easy}),
whence the equations

\begin{eqnarray}
\dim_{\mathrm{Krull}}\mathscr O_{Y\al,t\al}&=&\dim_{\mathrm{Krull}}\mathscr O_{X\al,s\al}
+\dim_{\mathrm{Krull}}\mathscr O_{Y\al_{s\all}, t\al}\label{d-krull1}\\
\dim_{\mathrm{Krull}}\mathscr O_{Z\al,t\al}&=&\dim_{\mathrm{Krull}}\mathscr O_{X\al,s\al}
+\dim_{\mathrm{Krull}}\mathscr O_{Z\al_{s\all}, t\al}\label{d-krull2}
\end{eqnarray}
which together with (\ref{d-krull-notadh})
yield
the equality
$\dim_{\mathrm{Krull}}\mathscr O_{Y\al,t\al}=\dim_{\mathrm{Krull}}\mathscr O_{Z\al,t\al}$. 
This implies that (the support of) $Z\al$ contains at
least one irreducible component of~$Y\al$ going through~$t\al$
and contradicts
the assumption that~$U\al=Y\al\setminus Z\al$ is a dense open subset of~$Y\al$, 
whence (1).

Let us now prove (2). The local ring $\mathscr O_{S,s}$ is a field, 
and $t$ lies in $\mathrm{Int}(T/S)$ since it lies in $\mathrm{Int}(Y/X)$
by assumption; therefore Theorem
\ref{thm-localring-generic}
implies that $\mathscr O_{Y_s,t}=\mathscr O_{T_s,t}$ if flat over $\mathscr O_{T,t}$. 
The latter is flat over $\mathscr O_{T\al, t\al}$ (by \ref{ss-gaga-alg}), and since
$\mathscr O_{S\al, s\al}$ is a field we have
\[\mathscr O_{T\al, t\al}=\mathscr O_{T\al_{s\all},t\al}=\mathscr O_{Y\al_{s\all}, t\al}.\]
As a consequence, $\mathscr O_{Y_s,t}$ is flat over $\mathscr O_{Y\al_{s\all}, t\al}$. 
On the other hand, $\mathscr O_{Y_s,t}$ is flat over $\mathscr O_{Y_s\al, t\al}$
(again by \ref{ss-gaga-alg}). The vertical and horizontal arrows of
the commutative diagram 
\[\xymatrix{
{\mathscr O_{Y_s,t}}&{\mathscr O_{Y\al_s,t_s\al}}\ar[l]\\
{\mathscr O_{Y\al_{s\all}, t\al}}\ar[u]\ar[ru]&}
\]
are thus flat; hence
$\mathscr O_{Y\al_{s\all}, t\al}\to \mathscr O_{Y\al_s,t_s\al}$
is flat too and (2) is proven. 

Due to (1), there exists $\omega\in U\al_{s\al}$ which
specializes to~$t\al$.
This implies in view of (2)
that there exists a point
$\omega'$ on~$Y\al_s$
lying above~$\omega$ and specializing to~$t_s\al$. 
Since $\omega'$ lies above $\omega$, it belongs to $U\al_s$; we thus have shown
that~$t_s\al\in \adht{U\al_s}{Y\al_s}$; 
but this
means that~$t\in \adht{U_s}{Y_{s,\mathrm{Zar}}}=\adht {U_s}{Y_s}$
or, in other words, that~$t$ belongs to $\adht U \phi$.
\end{proof}

\begin{theo}\label{thm-fiber-closure}
Let~$\phi \colon Y\to X$ be a morphism of~$k$-analytic spaces and let~$E$ be a locally
constructible subset of~$Y$. Let~$N$ be a subset of $\N$. The subset
of $\adht E\phi$ consisting of points at which the relative dimension of $\adht E\phi$
belongs to $N$ is a locally constructible subset of~$Y$. In particular,
$\adht E\phi$ is locally constructible (take $N=\N$). 
\end{theo}

\begin{proof}
We first consider the case where $N=\N$; \ie, we
first prove
that $\adht E \phi$ is locally constructible.

By arguing locally on $Y$
(which is possible due to Corollary \ref{cor-zarclosure-dom}), 
we can assume that
$E$ is constructible. 
Write $E=\bigcup U_i\cap F_i$, where~$(U_i)$, \resp ~$(F_i)$, 
is a finite family of Zariski-open, \resp Zariski-closed, subsets of~$Y$. 
Since $\adht {E} \phi=\bigcup_i \adht{U_i\cap F_i}\phi$,
it suffices to treat the case where~$E = U \cap F$ , with~$U$, \resp $F$, 
a Zariski-open, \resp Zariski-closed, subset of~$Y$. 
By replacing~$Y$ with~$F$ (equipped with any structure of a closed analytic subspace; \eg, its
reduced structure)
we can assume that~$E=U$. 

We are now going to apply some of the general results in \ref{s-diagtr},
which involve two categories $\mathfrak F$ and $\mathfrak C$ 
and a property $\mathsf Q$
as in \ref{ss-general-axioms-pfib}.
We take $\mathfrak F, \mathfrak C$ 
and $\mathsf Q$
as in Example \ref{ex-category-trick2}
(with $E=U$). 
What we want to prove is assertion $(\alpha)$ of
\ref{ss-claims}. By Lemma \ref{alphadash-alpha},  it is
sufficient to prove assertion $(\alpha'')$ of 
\ref{ss-auxiliary}. We thus reduce to the following situation: $Y$ and $X$ are affinoid, $X$ is integral,
and $\phi$ admits a section $\sigma$; 
and we have to prove that there exists a non-empty Zariski-open subset of $X$ which is
either contained in $\sigma\inv(\adht U\phi)$ or
disjoint
from it. 
This amounts to proving that there exists a non-empty Zariski-open subset of $\sigma(X)$ 
which is either contained in $\adht U\phi$ or disjoint of it. 
We choose $x\in X$ such that $x\al$ is the generic point of $X\al$
and we distinguish three cases. 

If $\sigma(x)\in U$, then $U\cap \sigma(X)$ is a non-empty Zariski-open subset of $\sigma(X)$
contained in $\adht U\phi$, and we are done. 
If $\sigma(x)\notin \adht UY$, then $\sigma(X)\setminus \adht UY$ is a non-empty
Zariski-open subset of~$\sigma(X)$ 
disjoint from $\adht UY$ (hence disjoint from $\adht U\phi$)
and we are done. 

We thus can assume that $\sigma(x)$ lies on $\adht U Y$
but not on $U$.
Under this assumption,~$\sigma$ goes 
through~$\adht UY_{\mathrm{red}}$; therefore by replacing~$Y$ with~$\adht UY_{\mathrm{red}}$, we
can assume that~$U$ is Zariski-dense in~$Y$. Let~$Z=(Y\setminus U)_{\mathrm{red}}$;
this is closed analytic subspace
of~$Y$ through which $\sigma$ factorizes. 
According to Theorem \ref{thm-flatloc-zaropen}, 
the~$X$-flat
loci $\uflat YX$ and $\uflat ZX$
are Zariski-open subsets of~$Y$
and~$Z$ respectively, and both of them
contain~$\sigma(x)$
by Lemma \ref{lem-flatloc-nonempty}. Therefore
$\uflat YX\cap \uflat ZX\cap \sigma(X)$
is a Zariski-open subset of $\sigma(X)$ which contains $\sigma(x)$,
and it suffices to prove that it is contained in $\adht U\phi$. But this follows
from Proposition \ref{prop-fiber-closure}
since $\sigma(X)\subset \mathrm{Int}(Y/X)$ by \ref{ss-boundary-basics} (3); this ends the 
proof when $N=\N$.

Let us consider now the case
of an arbitrary subset $N$ of $\N$.
We go back to the general assumptions of the theorem.
By arguing locally 
on $Y$ (which is possible due to Corollary \ref{cor-zarclosure-dom}),
we may assume that $Y$ is finite-dimensional and 
$E$ is constructible. Since the dimension of any fiber of $\phi$
is bounded by $\dim Y$, 
the subset of $\adht E\phi$ we are interested in 
is a finite Boolean combination
of sets of the form
\[\adht E\phi_{\geq \delta}:=\{y\in \adht E\phi| \dim_y \adht E\phi\geq \delta\}.\]
It thus suffices to prove that $\adht E\phi_{\geq \delta}$
is constructible for every $\delta \in \N$. 

The constructible subset $E$ can be written $\bigcup U_i\cap F_i$, where~$(U_i)$, \resp
$(F_i)$, 
is a finite family of Zariski-open, \resp Zariski-closed, subsets of~$Y$; 
the set $\adht E\phi_{\geq \delta}$ is then equal to $\bigcup_i \adht{U_i\cap F_i}\phi_{\geq \delta}$.
It is thus suffices to prove the theorem for each of the $U_i\cap F_i$'s; \ie, 
we may assume that $E=U\cap F$ with $U$, \resp $F$, a Zariski-open, \resp Zariski-closed, 
subset of $Y$. By replacing $Y$ with $F$ (equipped with any structure of a closed
analytic subspace of $Y$; \eg, its reduced structure), we reduce to the case where $E=U$. 

For every integer $\delta\geq 0$, let us denote by $G_{\geq \delta}$
be the subset of $Y$
consisting of points at which $\phi$ is of dimension $\geq \delta$; since $y\mapsto \dim_y \phi$
is upper semi-continuous for the Zariski topology (\cite{ducros2007}, \Th 4.9), $G_{\geq \delta}$ 
is Zariski-closed, and
$G_{\geq \delta} \cap U$ is thus constructible.

Let $x$ be a point of $X$
and let $(Y_i)$ be the family of irreducible components of $Y_x$. For every $i$,
we denote by $d_i$ the dimension of $Y_i$. Let $I$ be the set of indices $i$ such that $Y_i$ intersects $U$. 
Fix $\delta\in \N$. 
We have the following equalities:
\begin{eqnarray}\setcounter{equation}{1}
\adht U\phi_x&=&\bigcup_{i\in I} Y_i\\
(\adht U\phi_{\geq \delta})_x&=&\bigcup_{i\in I,d_i\geq \delta} Y_i\\
(G_{\geq \delta}\cap U)_x&=&\bigcup_{i\in I, d_i\geq \delta}(Y_i\cap U)\\
\adht {G_{\geq \delta}\cap U}\phi_x&=&\bigcup_{i\in I,d_i\geq \delta}Y_i
\end{eqnarray}
(note that (a) and (d) rest on density of $Y_i\cap U$ in $Y_i$ for every $i\in I$).  
We deduce from (b) and (d) (which hold for every $x\in X$) that
\begin{eqnarray*}
\adht U\phi_{\geq \delta}&=&\adht {G_{\geq \delta}\cap U}\phi.
\end{eqnarray*}
By the case $N=\N$ already proven, $\adht {G_{\geq \delta}\cap U}\phi$ is constructible, whence
the constructibility of $\adht U\phi_{\geq \delta}$. 
\end{proof}

We now come to our original motivation, namely the fiberwise codimension function.

\begin{prop}\label{prop-relat-codim}
Let~$\phi \colon Y\to X$ be a morphism of~$k$-analytic spaces, and let~$E$ and~$F$ be two
locally
constructible subsets of~$Y$ such that~$F\subset E$
(recall that by
Proposition \ref{prop-cons-gloc}, a locally constructible subset of $Y$
is constructible as soon as $Y$ is finite-dimensional). Let
$N$ be a subset of $\N\cup\{+\infty\}$. 
The set
\[\{y\in \adht E\phi|\codim_y (\adht F\phi_{\phi(y)},\adht E\phi_{\phi(y)})\in N\}\] is a locally constructible subset of~$Y$.
\end{prop}

\begin{rema}
We recall that a  a locally constructible subset of $Y$
is constructible as soon as $Y$ is finite-dimensional, by Proposition \ref{prop-cons-gloc}. 
We also recall that for $y$ a point of  $\adht E\phi$, the codimension 
$\codim_y (\adht F\phi_{\phi(y)},\adht E\phi_{\phi(y)})$ is equal to $+\infty$
if and only if $y\notin \adht F\phi$
(\ref{ss-codim-def}). 
\end{rema}

\begin{proof}[Proof of Proposition \ref{prop-relat-codim}]
By arguing locally (which is possible 
in view of \ref{ss-codim-basics} (1)), we can assume that~$Y$
is finite dimensional. 
Let~$D\in \N$ be
such that all fibers of~$Y\to X$ are of dimension bounded
by~$D$ (if~$Y\neq \emptyset$ one may take~$D=\dim Y$).
For every $(n,m)$ in $\{0,\ldots,D\}^2$, we denote
by $\Lambda_{n,m}$
the set of~$y\in \adht F\phi$ for which there exist 
an irreducible component~$T$ of~$\adht F\phi_{\phi(y)}$ 
of dimension~$m$
and an irreducible component~$Z$ of~$\adht E\phi_{\phi(y)}$ of dimension~$n$ 
with~$y\in T\subset Z$ (note that $\Lambda_{n,m}=\emptyset$
if $m>n$).

Fix $(n,m)\in \{0,\ldots, D\}^2$; we are going to give an alternative description of $\Lambda_{n,m}$
from which we shall deduce that it is constructible. For that purpose, let us denote by $G$
(\resp $H$) the subset of $\adht  E\phi$ (\resp $\adht F \phi$) consisting of points
at which $\adht E\phi$ (\resp $\adht F \phi$) is of
relative dimension $n$ (\resp $m$)
over 
$X$.
By Theorem
\ref{thm-fiber-closure}
and Proposition \ref{prop-cons-gloc}, the subsets $H$ and $G$ of $Y$ are constructible, and
thus so 
are $\adht G \phi$ and $\adht H \phi$. 

Let $x$ be a point of $X$. By definition, 
\[G_x=\{y\in \adht E\phi_x|\dim_y \adht E\phi_x=n\}\]
and $\adht G\phi_x$ is the closure 
of $G_x$ inside $Y_x$. As a consequence, $\adht G\phi_x$ is the union of all $n$-dimensional irreducible components of $\adht E\phi_x$. Analogously, 
$\adht H \phi_x$ is the union of all $m$-dimensional irreducible components of $\adht F\phi_x$. Now let $T$ be an $m$-dimensional irreducible component of
$\adht F\phi_x$. By the above, $T$ is contained in an $n$-dimensional irreducible component of  $\adht E\phi_x$ if and only if it is contained
in $\adht G \phi_x$ or, what amounts to the same, if and only if it is an irreducible component of $\adht F\phi_x\cap \adht G\phi_x$; conversely, any $m$-dimensional
irreducible component of $\adht F\phi_x\cap \adht G\phi_x$ is contained in $ \adht G\phi_x$ and is an irreducible
component of $\adht F\phi_x$
by a dimension argument. Therefore $\Lambda_{n,m,x}$ is the union of all $m$-dimensional irreducible components of $\adht F\phi_x\cap \adht G\phi_x$, whence
we get
the equality  
\[\Lambda_{n,m,x}=
\{y\in \adht F\phi_x\cap \adht G\phi_x|\dim_y \adht F\phi_x\cap \adht G\phi_x=m\}\]
(because $\adht F\phi_x\cap \adht G\phi_x$ is contained in $\adht F\phi_x$ and is
thus of dimension $\leq m$).

We have thus proved that $\Lambda_{n,m}$ is the subset of
$\adht G\phi \cap \adht H\phi$ consisting of points at which $\adht G\phi \cap \adht H\phi$
is of relative dimension $m$; using again Theorem \ref{thm-fiber-closure}, we get the constructibility of $\Lambda_{n,m}$.

For any~$y\in \adht E\phi$, let us denote by~$\daleth(y)$ the set of~pairs $(n,m)\in \{0,\ldots,D\}^2$ such that~$y\in \Lambda_{(n,m)}$. 
If~$P$ is a subset of ~$ \{0,\ldots,D\}^2$,  then the set of~$y\in \adht{E}\phi$ such that~$\daleth(y)=P$
is a boolean combination of some of the~$\Lambda_{(n,m)}$'s and is therefore a constructible subset of~$Y$. 
Now the function from~$\adht{E}\phi$ to~$\N\cup\{+\infty\}$ 
that sends~$y$ to~$\codim_y  (\adht F\phi_{\phi(y)}, \adht E\phi_{\phi(y)})$
is constant on every fiber of~$\daleth$
(note that the set of points at which it takes the value $+\infty$ is
exactly $\daleth\inv(\emptyset)$), whence the proposition.
\end{proof}

\section{The fiberwise exactness locus}\label{s-fibexact} 

Let $Y\to X$ be a morphism of $k$-analytic spaces, and let $D=(\mathscr G\to \mathscr F\to \mathscr E)$
be a sequence of coherent sheaves on $Y$. In this section we prove the following. 

\begin{enumerate}[1]

\item The set of points
of $Y$ at which $D$ is fiberwise a complex, or fiberwise exact, or more generally
at which it is
a complex with fiberwise homology at $\mathscr F$
having a given fiber rank, is locally constructible (Theorem \ref{thm-fiberwise-exact}). 

\item The set of points
of $Y$ at which $D$ is fiberwise exact and $\mathscr E$ is $X$-flat is Zariski-open (Theorem
\ref{thm-exactflat-locus}). 

\item The set $U$ of points
of $Y$ at which $D$ is fiberwise exact and $\mathscr E$ and $\mathscr F$ are $X$-flat is Zariski-open,  
and $D_U$ is exact (Proposition \ref{prop-exactflat-locus}).

\item The set of points of $Y$ at which $\mathscr E$ is fiberwise locally free
is locally constructible (Theorem \ref{thm-locfree-locus} (1)).

\item The set $U$ of points of $Y$ at which $\mathscr E$ is fiberwise locally free and $X$-flat
is locally constructible, and $\mathscr E_U$ is locally free (Theorem \ref{thm-locfree-locus} (2)).
\end{enumerate}
Those are more or less
analogues of classical results in scheme theory (with often
slightly different proofs): see \cite{ega43}, \Prop 9.4.4 for (1); 
 \Prop 12.3.3 for (2) and (3); and  \Prop 9.4.7 for (4). But note that this analogy (and the intrinsic interest of these
statements) was not our only motivation: we shall also need (1) and (2) for studying the fiberwise codepth and the fiberwise Gorenstein
property, which we will investigate by homological
methods. 

\begin{enonce}[remark]{Notation}
If~$Y$ is an analytic space and if~$D=(\mathscr G\to \mathscr F\to \mathscr E)$ is a sequence
of morphisms of coherent sheaves on~$Y$,
we shall denote by~$\mathscr Z(D)$ (\resp $\mathscr B(D)$) the kernel (\resp the image) of
$\mathscr F\to \mathscr E$ (\resp ~$\mathscr G\to \mathscr F$);
note that both~$\mathscr Z(D)$ and~$\mathscr B(D)$ are coherent subsheaves of~$\mathscr F$;
their internal sheaf-theoretic sum inside~$\mathscr F$ will be denoted by~$\mathscr C(D)$. If~$y\in Y$, then~$D$
is a complex at~$y$ if and only if~$(\mathscr C(D)/\mathscr Z(D))_{\hr y}=0$, and it is exact at~$y$
if and only if
\[(\mathscr C(D)/\mathscr Z(D))_{\hr y}=0\;\text{and}\;(\mathscr C(D)/\mathscr B(D))_{\hr y}=0.\]
\end{enonce}

\begin{rema}\label{rem-cdz-flat}
The functors that send
$D$
to $\mathscr C(D), \mathscr B(D)$,
and $\mathscr Z(D)$ commute with flat base change and with ground
field extension: this is an immediate consequence of
Proposition \ref{prop-flat-univexact} (1). 
\end{rema}

\begin{theo}\label{thm-fiberwise-exact}
Let~$\phi \colon Y\to X$ be a morphism of~$k$-analytic spaces, 
let $N$ and~$M$ be two subsets of $\N$ and
let~$D=(\mathscr G\to \mathscr F\to \mathscr E)$
be a sequence of morphisms of coherent sheaves on~$Y$.

\begin{enumerate}[1]
\item The set~$C$ of points~$y\in Y$ such that
\[\rk y (\mathscr C(D_{Y_{\phi(y)}})/\mathscr Z(D_{Y_{\phi(y)}}))\in N\]
and
\[\rk y (\mathscr C(D_{Y_{\phi(y)}})/\mathscr B(D_{Y_{\phi(y)}}))\in M\]
is a locally constructible subset of~$Y$. 

\item The set of points~$y\in Y$ such that $D$ is fiberwise a complex at $y$ 
(\ie, $D_{Y_{\phi(y)}}$ is a complex at $y$) is a locally constructible subset of $Y$.

\item The set of points~$y\in Y$ such that $D$ is fiberwise exact at $y$
(\ie, $D_{Y_{\phi(y)}}$ is exact  at $y$) is a locally constructible subset of $Y$.

\end{enumerate}
\end{theo}

\begin{proof}
We first note that (2) is nothing but (1) for $N=\{0\}$ and $M=\N$, and that (3)
is nothing but (1) for $N=M=\{0\}$. It suffices thus to prove (1). We introduce the following
notation: if $Z$ is an analytic space and  if $\Delta$ is a 3-term sequence of coherent sheaves
on $Z$, 
then for any point $z$ of $Z$ we set $\lambda(\Delta,z)=\rk z(\mathscr C(\Delta)/\mathscr B(\Delta))$ 
and $\mu(\Delta,z)=\rk z(\mathscr C(\Delta)/\mathscr Z(\Delta))$.

We are now going to apply some of the general results described in \ref{s-diagtr},
which involve two categories $\mathfrak F$ and $\mathfrak C$ 
and a property $\mathsf Q$
as described in \ref{ss-general-axioms-pfib}. We take for $\mathfrak F$ the category 
of 3-term sequences of coherent sheaves (\ie, $\mathfrak{Coh}^{\{0,1,2\}}$ with the notation of 
Example \ref{ex-fiber-diag}), for $\mathfrak C$ the category of all analytic spaces, 
and for $\mathsf Q$ the property defined as follows: if $Z$ is an analytic space and 
if $\Delta$ is an object of $\mathfrak F_Z$, then $\Delta$ satisfies $\mathsf Q$ at a
point $z$ of $Z$ if $\lambda(\Delta,z)\in  N$ and
$\mu(\Delta, z) \in M$ (axioms (1) and (2)
of \ref{ss-general-axioms-pfib}
are fulfilled due to Remark \ref{rem-cdz-flat}). 

What we want to prove is assertion $(\alpha)$ of
\ref{ss-claims}. By Lemma \ref{alphadash-alpha},  it is
sufficient to prove assertion $(\alpha'')$ of 
\ref{ss-auxiliary}. We thus reduce to the following situation: $Y$ and $X$ are affinoid, $X$ is integral,
and $\phi$ admits a section $\sigma$; 
and we have to prove that there exists a non-empty Zariski-open subset of $X$ which is
either contained in $\sigma\inv(C)$ or
disjoint
from
it. Let us introduce some notation. We denote by~$\mathscr K$ the kernel of~$\mathscr G\to \mathscr F$
and by~$\mathscr F'$ its cokernel; we denote by~$\mathscr I$ the
image of~$\mathscr F\to \mathscr E$ and by~$\mathscr E'$ its cokernel;
we denote by~$\mathscr L$ the cokernel of~$\mathscr C(D)\to \mathscr F$. We have the following exact sequences:

\[0\to \mathscr K\to \mathscr G\to \mathscr B(D)\to 0, \;0\to \mathscr B(D)
\to \mathscr F\to \mathscr F'\to 0,$$~$$0\to \mathscr Z(D)\to \mathscr F\to \mathscr I\to 0, \;0\to \mathscr I\to \mathscr E\to \mathscr E'\to 0,$$~$$0\to \mathscr C(D)
\to \mathscr F\to \mathscr L\to 0\;,\;\mathscr B(D)\oplus\mathscr Z(D)\to \mathscr C(D)\to 0.\]

Let~$U$ be the~$X$-flat
locus of~$\mathscr B(D)\oplus  \mathscr F'\oplus \mathscr I\oplus \mathscr E'\oplus \mathscr L$. By
Theorem \ref{thm-flatloc-zaropen}, $U$ is a Zariski-open subset of $Y$; and 
by~Lemma \ref{lem-flatloc-nonempty}, 
the Zariski-open subset~$\sigma\inv(U)$ of~$X$ is non-empty. 
As $\mathscr F'_U$ and $\mathscr B(D)_U$ are $X$-flat, 
so is $\mathscr F_U$ and hence $\mathscr C_D(U)$ since $\mathscr L_U$
is $X$-flat too (Lemma \ref{lem-flat-extension}).

Let $x$ be a point of $X$. 
By construction, the restriction to $U$
of each of the above exact sequences
has an $X$-flat
right
term, hence remains exact after arbitray base change
$X'\to X$ by Proposition \ref{prop-flat-univexact}. In particular, they remain exact
after restriction to the fiber $U_x$; therefore 
the coherent sheaves~$\mathscr B(D_{U_x}),\mathscr Z(D_{U_x})$ and~$\mathscr C(D_{U_x})$ are respectively naturally isomorphic to
$\mathscr B(D)_{U_x}, \mathscr Z(D)_{U_x}$, and~$\mathscr C(D)_{U_x}$, whence the equalities
\begin{eqnarray}\setcounter{equation}{1}
\lambda(D,y)&=&\lambda(D_{U_{\phi(y)}}, y)\\
\mu(D,y)&=&\mu(D_{U_{\phi(y)}}, y)
\end{eqnarray}
for all $y\in U$. 

Let~$E$ be the set of points~$y\in Y$ such that $\lambda(D,y)$ and  $\mu(D,y)$
belong respectively to $N$ and $M$. Since the pointwise rank function of a given coherent sheaf on the affinoid
space $Y$ takes
only finitely many values, $E$ is a constructible subset of $Y$. Moreover we have 
in view of (a) and (b)
the equality $C\cap U=E\cap U$.

The pre-image~$\sigma\inv(E)$
is a constructible subset of the integral 
affinoid space~$X$.
Therefore, there exists a non-empty Zariski-open subset~$V$ 
of~$X$ such that~$V\subset \sigma\inv(E)$
or~$V\subset X\setminus \sigma\inv(E)$. Now~$W:=V\cap \sigma\inv(U)$ is
a non-empty Zariski-open subset
of~$X$. 
From the equality~$C\cap U=E\cap U$
it follows that $W\subset \sigma\inv(C)$ if~$V\subset \sigma\inv(E)$
and that $W\subset X\setminus \sigma\inv(C)$ otherwise. 
\end{proof}

\begin{theo}\label{thm-exactflat-locus}
Let~$\phi \colon
Y\to X$ be a morphism between~$k$-analytic spaces and 
let~$D=(\mathscr G\to \mathscr F\to \mathscr E$) be a complex
of coherent sheaves on~$Y$.
The set~$A$ of points of~$Y$ at which~$\mathscr E$ is~$X$-flat and 
$D$ is fiberwise exact is a Zariski-open subset of~$Y$.
\end{theo}

\begin{proof}
We are going to
apply some of the general results described in \ref{s-diagtr},
which involve two categories $\mathfrak F$ and $\mathfrak C$ 
and a property $\mathsf Q$
as described in \ref{ss-general-axioms-pfib}, and a functor
$\mathscr S$ as described in \ref{ss-def-functorS}. 
We take for $\mathfrak F$ the category 
of 3-term
complexes of coherent sheaves, for $\mathfrak C$ the category of all analytic spaces, 
and for $\mathsf Q$ the exactness
property. We take for $\mathscr S$ the functor that sends
a 3-term complex to its right
term. We want to prove assertion $(\beta)$ 
of \ref{ss-claims}. Since assertion $(\gamma)$ of \loccit~holds
(this is Theorem \ref{thm-flatloc-zaropen}), Lemma \ref{gamma-betaprime-beta}
ensures that it suffices to prove assertion $(\beta')$ of \ref{ss-auxiliary}; \ie, 
we can assume that $Y$ and $X$ are affinoid
and  $Y\to X$ has a section~$\sigma$, 
and we only need to prove that~$\sigma\inv(A)$ is a Zariski-open subset of~$X$.

Let $A_0$ be the subset of $Y$ consisting
of points at which $D$ is fiberwise exact. By Theorem \ref{thm-fiberwise-exact} (3),   $A_0$
is a constructible subset of $Y$ and since $Y$ is affinoid (hence finite-dimensional),
$A_0$ is constructible by Proposition \ref{prop-cons-gloc}. 
By Theorem \ref{thm-flatloc-zaropen},
$\uflat {\mathscr E}X$
is Zariski-open subset of $Y$. As a consequence, $A=A_0\cap \uflat{\mathscr E}X$
is a constructible subset of $Y$,
and $\sigma\inv(A)$ is thus
a constructible subset of $X$. 
%In order to prove that $\sigma\inv(A)$ is Zariski-open, 
%it suffices to prove that for every pair $(T,Z)$ of irreducible Zariski-closed
%subsets of $X$ with $T\subset Z$
%and $\sigma\inv(A)\cap T$ Zariski-dense in $T$,  the intersection $\sigma\inv(A)\cap Z$
%is Zariski-dense in $Z$: indeed, this is nothing but a rephrasing of the
%fact that a constructible subset of $X\al$ is Zariski-open if (and only if) it is stable under generization.  
%(Recall that 
%a constructible subset of an irreducible analytic space is Zariski-dense
%if and only if it contains a non-empty \emph{Zariski-open}
%subset of this space; this comes from
%Lemma \ref{lem-cons-dense}
%but in the affinoid case -- the only one we will
%need here -- this can also be deduced from the corresponding scheme-theoretic statement). 

Let $(T,Z)$ be a pair of irreducible Zariski-closed
subsets of $X$ with $T\subset Z$ and $T\cap \sigma\inv(A)$ Zariski-dense in $T$; 
we are going to prove
that $Z\cap \sigma\inv (A)$ is Zariski-dense in $Z$, which will yield
to the Zariski-openness of $A$
(Lemma \ref{lem-cons-dense} (3); see also Remark \ref{rem-aff-topcons}).

By \ref{ss-codim-basics} (2) there exists a chain
$T_0=T\subset T_1\ldots \subset T_{\codim(T,Z)}=Z$ where every $T_i$ 
is an irreducible Zariski-closed subset of $X$ and where we have
$\codim(T_i,T_{i+1})=1$
for all $i<\codim(T,Z)$. Hence we reduce by induction on
codimension 
to the case where $\codim(T,Z)=1$. 
Note that since $\sigma\inv(A)\cap T$ is a Zariski-dense
constructible subset of $T$, it contains a non-empty Zariski-open subset of $T$
(Lemma \ref{lem-cons-dense} (1);  see also Remark \ref{rem-aff-topcons}). 
Let $\widetilde Z$ be the normalization of $Z_{\mathrm{red}}$.  
Choose an irreducible component $\widetilde T$ of  $\widetilde Z \times_Z T$
that dominates $T$ (this exists by surjectivity of the scheme-theoretic normalization map 
$\widetilde Z\al\to Z\al$).
Set $\widetilde Y=Y\times_X \widetilde Z$, let $\widetilde A_0$ and $\widetilde A$ 
be the pre-images of $A_0$ and $A$ on $\widetilde Y$, and let $\widetilde \sigma \colon \widetilde Z\to \widetilde Y$
denote the map induced by $\sigma$. Let $B$ be the set of points of $\widetilde Y$ at which 
$\mathscr E_{\widetilde Y}$ is $\widetilde Z$-flat and at which $D_{\widetilde Y}$ is fiberwise exact;
since the formation of the 
fiberwise exactness locus commutes to arbitrary base change, we
have $\widetilde A\subset \widetilde A_0\cap \uflat{\mathscr E_{\widetilde Y}}{\widetilde Z}=B.$
As the intersection
$\sigma\inv(A)\cap T$ contains a non-empty Zariski-open subset of $T$,
its pre-image $\widetilde \sigma \inv(\widetilde A)\cap \widetilde T$ under the finite dominant map 
$\widetilde T\to T$
is Zariski-dense in $\widetilde T$; hence
$\widetilde \sigma \inv(B)\cap \widetilde T$ is a fortiori Zariski-dense in $\widetilde T$, and
it suffices to prove that this implies 
Zariski-density of
the intersection $\widetilde \sigma \inv(B)\cap \widetilde Z$
in $\widetilde Z$. 
Indeed, assume that the latter holds. Then $\widetilde \sigma\inv(\widetilde A_0)\cap \widetilde Z$ is a
Zariski-dense subset of $\widetilde Z$;  since $\widetilde Z\to Z$
is birational, 
it follows that $\sigma \inv(A_0)\cap Z$ is Zariski
dense in $Z$. On the other hand, $\sigma\inv(\uflat {\mathscr E}X)\cap T$
contains $\sigma\inv(A)\cap T$, hence is Zariski-dense in $T$ and in particular non-empty; since  $\sigma\inv(\uflat {\mathscr E}X)$
is a Zariski-open subset of $X$, the intersection
$\sigma\inv(\uflat {\mathscr E}X)\cap Z$ is a non-empty Zariski open subset of
the irreducible space $Z$, 
so
$\sigma\inv(A)\cap Z=\sigma \inv(A_0)\cap \sigma\inv(\uflat{\mathscr E}X)\cap Z$
is Zariski-dense in $Z$. 

As a consequence, by performing base change from $X$ to $\widetilde Z$ and by considering
the pair $(\widetilde T, \widetilde Z)$ we reduce to the case where
$Z=X$ and $X$ is integral and normal. Set $n=\dim X$; we then have $\dim T=n-1$. 

The set $T\cap \sigma\inv(A)$
contains
a non-empty Zariski-open subset of $T$,
so
there exists $t\in T$ such that $d_k(t)=n-1$. 
By Example \ref{ex-centdim-recap}, both local rings
$\mathscr O_{X,t}$ and $\mathscr O_{X\al,t\al}$ are of dimension 1, 
hence are discrete valuation rings as $X$ is normal. 
Let $\varpi\in \mathscr O_X(X)$ be a function that generates $\mathfrak m_{t\al}$; \ie, 
it is a uniformizing parameter of $\mathscr O_{X\al, t\al}$. 
As
$\mathfrak m_{t\al} \mathscr O_{X,t}=\mathfrak m_t$ by Example \ref{ex-centdim-recap},
$\varpi$
is a uniformizing parameter of $\mathscr O_{X,t}$ as well. 

Let us endow $T$ with its reduced structure, and let us choose a point 
$x$ in $X$ with $d_k(x)=n$.  By Remark \ref{rem-abh-point}, 
the point $x$ is Zariski-dense in $X$ and the point $t$ is Zariski-dense in $T$; \ie, 
$x\al$ is the generic point of $X\al$ and $t\al$ is the generic point of $T\al$. 
By Example \ref{ex-centdim-recap}, both local rings
$\mathscr O_{X,x}$ and $\mathscr O_{T,t}$ are
artinian, hence are fields because $X$ and $T$ are reduced. 

In order to prove that $\sigma\inv (A)$
is Zariski-dense, it suffices 
to prove that $x\in \sigma\inv (A)$. 
Since $t$ belongs to $\sigma\inv(A)$, it belongs to $\sigma\inv(\uflat {\mathscr E}X)$; 
the latter is thus
a
non-empty Zariski-open subset of $X$, hence
it contains $x$.
It remains only
to show that 
$x$ belongs to $\sigma\inv(A_0)$; \ie, $D$ is fiberwise exact at $\sigma(x)$. 
We shall need the following flatness assertions:

\begin{enumerate}[a]

\item The local ring~$\mathscr O_{Y_x,\sigma(x)}$ is flat
over~$\mathscr O_{Y\al,\sigma(x)\al}$. 

\item The local ring~$\mathscr O_{Y_t,\sigma(t)}$ is flat
over~$\mathscr O_{Y\al,\sigma(t)\al}/\varpi\mathscr O_{Y\al,\sigma(t)\al}$.
\end{enumerate}

Let us first prove (a). 
The point $\sigma(x)$ belongs to~$\mathrm{Int}(Y/X)$ by
\ref{ss-boundary-basics} (3)
and~$\mathscr O_{X,x}$ 
is a field; therefore~$\mathscr O_{Y_x,\sigma(x)}$ is flat
over~$\mathscr O_{Y,\sigma(x)}$ by Theorem \ref{thm-localring-generic}. 
On the other hand, 
$\mathscr O_{Y,\sigma(x)}$
is flat over~$\mathscr O_{Y\al, \sigma(x)\al}$ (\ref{ss-gaga-alg}), so
$\mathscr O_{Y_x,\sigma(x)}$ is flat
over~$\mathscr O_{Y\al,\sigma(x)\al}$.

Let us now prove (b). 
The point $\sigma(t)$ belongs to $\mathrm{Int}(Y/X)$ by
\ref{ss-boundary-basics} (3), hence also to $\mathrm{Int}(Y\times_X T/T)$.
Since the local ring
$\mathscr O_{T,t}$ is a field,
Theorem \ref{thm-localring-generic}
ensures that
$\mathscr O_{(Y\times_X T)_t,\sigma(t)}=\mathscr O_{Y_t, \sigma(t)}$ is flat
over~$\mathscr O_{Y\times_XT, \sigma(t)}$. The latter
is itself flat over~$\mathscr O_{ Y\al\times_{X\al}T\al, \sigma(t)\al}$, which is nothing
but~$\mathscr O_{Y\al,\sigma(t)\al}/\varpi\mathscr O_{Y\al,\sigma(t)\al}$
because since $T\al$ is reduced, it is defined by the equation~$\varpi=0$ around
its generic point $t\al$.
Hence~$\mathscr O_{Y_t,\sigma(t)}$ is flat
over $\mathscr O_{Y\al,\sigma(t)\al}/\varpi\mathscr O_{Y\al,\sigma(t)\al}$.

By assumption, $\sigma(t)\in A$. Therefore,~$D$
is fiberwise
exact at~$\sigma(t)$ and $\mathscr E$
is $X$-flat at $\sigma(t)$.  Fiberwise exactness of $D$ at $\sigma(t)$ 
means
that $D_{Y_t,\sigma(t)}$ is exact. This implies
in view of (b) that 
\[D\al_{Y\al,\sigma(t)\al}\otimes_{\mathscr O_{Y\al, \sigma(t)\al}}
\mathscr O_{Y\al, \sigma(t)\al}/\varpi\mathscr O_{Y\al, \sigma(t)\al}\]
is exact too. 

Since $\mathscr E$ is $X$-flat at $\sigma(t)$, Lemma \ref{lem-gagaflat-easy}
ensures that $\mathscr E_{Y\al,\sigma(t)\al}$ is
flat over~$\mathscr O_{X\al, t\al}$. Applying
Lemma 12.3.3.1 of \cite{ega43}
with~$B=\mathscr O_{Y\al,\sigma(t)\al}$, $t=\varpi$, 
and~$M,N$, and~$P$ respectively equal to 
$\mathscr G_{Y\al, \sigma(t)\al}, \mathscr F_{Y\al, \sigma(t)\al}$,
and $\mathscr E_{Y\al, \sigma(t)\al}$ (which is possible since
$\mathscr E_{Y\al, \sigma(t)\al}$ has no
non-zero $\varpi$-torsion because
it is flat over $\mathscr O_{X\al, t\al}$), we 
get the exactness of the complex~$D_{Y\al, \sigma(t)\al}$. 
Moreover the local ring~$\mathscr O_{Y\al, \sigma(x)\al}$
is a localization of~$\mathscr O_{Y\al, \sigma(t)\al}$ because~$x\al$ is a generalization
of $t\al$.  It follows that
the complex~$D_{Y\al, \sigma(x)\al}$ is exact. This implies
in view of (a) that 
$D_{Y_x,\sigma(x)}$ is exact.
\end{proof}

\begin{prop}\label{prop-exactflat-locus}
Let~$\phi\colon Y\to X$ be a morphism of~$k$-analytic spaces and let~$D=(\mathscr G\to \mathscr F\to \mathscr E)$ be a complex
of coherent sheaves on~$Y$.
The set~$A$ of points of~$Y$ at which~$\mathscr E$ and~$\mathscr F$
are~$X$-flat and at which~$D$ is fiberwise exact is a Zariski-open subset of~$Y$ on
which~$D$ is exact.
\end{prop}

\begin{proof}
It follows from Theorem \ref{thm-exactflat-locus}
above and from
Zariski-openness of the~$X$-flat locus of~$\mathscr F$ 
(Theorem \ref{thm-flatloc-zaropen})
that~$A$ is Zariski-open. It remains to prove that~$D$
is exact on~$A$. 

We can assume that~$Y$ and~$X$ are affinoid. Let~$y$ be a point of $A$
and let~$x$ be its image in~$X$. We shall prove that~$D_{Y,y}$ 
is exact. This can be proved
after enlarging the ground field, which allows us
to assume
that $x$ is rigid. 
The local ring~$\mathscr O_{Y_x,y}$ is
then equal to~$\mathscr O_{Y,y}/\mathfrak m_x\mathscr O_{X,x}$.

By our assumptions, 
the complex
$D_{Y,y}=(\mathscr G_{Y,y}\to \mathscr F_{Y,y}\to \mathscr E_{Y,y})$
enjoys the following properties: 

\begin{itemize}[label=$\bullet$]

\item The~$\mathscr O_{X,x}$-modules~$\mathscr F _{Y,y}$
and~$\mathscr E_{Y,y}$ are flat. 

\item The complex~$D_{Y,y}\otimes_{\mathscr O_{Y,y}}\mathscr O_{Y,y}/\mathfrak m_x\mathscr O_{Y,y}$
is exact. 
\end{itemize}

A repeated application of Lemma 12.3.3.5 of \cite{ega43}
ensures that
the complex~$D_{Y,y}\otimes_{\mathscr O_{Y,y}}/\mathfrak m_x^n\mathscr O_{Y,y}$ is exact for every~$n>0$. 
As the complex~$D_{Y,y}$ only involves finitely
generated modules over the noetherian ring~$\mathscr O_{Y,y}$
and as $\H(D_{Y,y}\otimes_{ \mathscr O_{Y,y}}
\mathscr O_{Y,y}/\mathfrak m_x^n\mathscr O_{Y,y})=0$
for all $n>0$, 
it follows
from \cite{ega32}, 7.4.7.2 
that
\[\H(D_{Y,y})\otimes_{\mathscr O_{Y,y}}\lim_\leftarrow \mathscr O_{Y,y}/\mathfrak m_x^n\mathscr O_{Y,y}=0.\]

Being finitely generated over~$\mathscr O_{Y,y}$, the module~$\H(D_{Y,y})$
is separated for the~$\mathfrak m_y$-adic topology
and a fortiori
for the $\mathfrak m_x$-adic topology, so it is itself zero, which ends the proof.
\end{proof}

\begin{theo}\label{thm-locfree-locus}
Let~$\phi \colon
Y\to X$ be a morphism of~$k$-analytic spaces and let $\mathscr E$ be a coherent sheaf on~$Y$. 

\begin{enumerate}[1]

\item The set~$A$ of points of~ $Y$
at which~$\mathscr E$ is fiberwise locally free is a constructible subset of~$Y$.

\item The set~$B$ of points of~ $Y$
at which~$\mathscr E$ is~$X$-flat and fiberwise locally free is a Zariski-open subset of~$Y$ over which~$\mathscr E$ is locally free.
\end{enumerate}
\end{theo}
\begin{proof}
We may assume that~$Y$ and~$X$ are affinoid. Let~$y$
be a point of $Y$. We set 
$r=\rk y (\mathscr E)$
and we choose~$r$ global sections~$f_1,\ldots,f_r$ of~$\mathscr E$ on~$Y$ such that~$(f_1(y),\ldots,f_r(y))$ is a basis of~$\mathscr E_{\hr y}$; let~$\ell$ be the morphism~$\mathscr O_Y^r\to \mathscr E$ 
that is induced by the~$f_i$'s. Let~$E$ be the set of~$z\in Y$ such that~$\ell_{Y_{\phi(z)}}$ is an isomorphism at~$z$. By assertion (3) of Theorem \ref{thm-fiberwise-exact}, 
$E$ is a constructible subset of~$Y$. 

If $y\in A$, the
constructible subset $E$ of $Y$ contains~$y$, and it is included in~$A$ by definition. 
Now assume that~$y\notin A$.
Let~$F$ be the set of points~$z$ of~$Y$ such that $\rk z(\mathscr E)=r$.  It is a constructible subset of~$Y$ which contains~$y$. Let~$V$ 
be the complement of~$\mathrm{Supp}(\mathrm {Coker}\;\ell)$; it is a Zariski-open subset which contains~$y$. Let~$G$ be the set of~$z\in Y$ at which~$\ell_{Y_{\phi(z)}}$ is
\emph{not}
an isomorphism. By assertion \ref{thm-fiberwise-exact} (3) $G$ is a constructible subset of~$Y$, and it contains~$y$
since $y\notin A$. 
Now  consider $z\in F\cap V\cap G$. The rank $\rk z(\mathscr E)$ is equal to~$r$ and~$\ell_{Y_{\phi(z)}}$ is surjective at~$z$; if~$\mathscr E_{Y_{\phi(z)}}$ were locally free at~$z$, then~$\ell_{Y_{\phi(z)}}$ would be an isomorphism at~$z$, which would contradict the assumption that~$z\in G$. Therefore, 
the set~$F\cap V\cap G$ is a constructible subset of~$Y$ which contains~$y$ and is included in~$Y\setminus A$. 
Assertion (1) then follows by
quasi-compactness of~$Y$ for the constructible topology
(\ref{ss-constop-compact}).

Assume that~$y\in B$. Let~$U$ be the~$X$-flat locus of~$\mathscr E$; it is a Zariski-open subset of~$Y$ by
Theorem \ref{thm-flatloc-zaropen}. It follows
from Lemma~\ref{lem-iso-fiber} that~$\ell$ is an isomorphism at every point of~$U\cap E$. This shows that the constructible set~$U\cap E$ is open, hence Zariski-open by \cite{berkovich1993}, \Cor 2.6.6,
and that~$\mathscr E$ is locally free on~$U\cap E$.
Hence~$U\cap E$ is a Zariski-open subset of~$Y$ containing~$y$, contained in~$B$, and on which~$\mathscr E$ 
is locally free; this proves~(2).
\end{proof}

\section{Regular sequences}\label{s-reg-seq}

Motivated by the study of the fiberwise CI property, we introduce in this section the 
notion of a \emph{regular sequence} with respect to a coherent sheaf (but we shall in fact mainly use it 
for the structure sheaf). The main result is Proposition \ref{prop-regular-seq}
below, which investigates the locus of ``fiberwise regularity" of a given sequence. 
Its proof is quite easy modulo our
preceding
work on fiberwise homology
(Theorem \ref{thm-fiberwise-exact} and Proposition \ref{prop-exactflat-locus}),
and some basic results on flatness (Proposition
\ref{prop-flat-univexact} and Lemma \ref{lem-flattrick-reg}).

\begin{defi}\label{def-regular-seq}\index{regular sequence}
Let~$X$ be a~$k$-analytic space, let $\mathscr F$ be a coherent sheaf on $X$,
and let~$(g_1,\ldots, g_n)$ be a family of analytic functions
on~$X$. For every~$i$, 
we set~$\mathscr F_i=\mathscr O_X/(g_1,\ldots, g_{i-1})$,
and 
we denote by~$h_i$ the endomorphism~$a\mapsto g_i a$ of~$\mathscr F_i$. 
Let~$x$ be a point of~$X$ at which all the~$f_i$'s vanish. The sequence~$(g_1,\ldots, g_n)$
is
called \emph{$\mathscr F$-regular}
at~$x$ if every~$h_i$ is injective at~$x$. We shall say 
often
say 
\emph{regular}
instead of $\mathscr O_X$-regular.
\end{defi}

\subsection{Basic properties}
We keep the notation of Definition \ref{def-regular-seq}. If~$V$ is an analytic domain
of~$X$ containing~$x$, then~$(g_i)$ is regular at~$x$ if and only if~$(g_{i|V})$ is regular at~$x$. 
If~$V$ is good, this is equivalent to requiring that~$(f_i)$ is a regular sequence of the module
$\mathscr F_{V,x}$ over the noetherian 
local ring~$\mathscr O_{V,x}$ ; since the latter property is invariant under any
permutation of the~$g_i$'s, 
the property of~$(g_1,\ldots, g_n)$
being regular at~$x$ is invariant under any
permutation 
of the~$g_i$'s. 

\subsection{}\label{ss-generation-i}
Let~$Y\to X$ be a morphism between~$k$-analytic spaces, let $\mathscr F$
be a coherent sheaf on $Y$. Let~$y$ be a point of $Y$
and let~$x$ be its image in~$X$. Let~$(f_1,\ldots, f_n)$ be a family of analytic functions on~$Y$
vanishing at~$y$. 
We shall say that~$(f_i)$ is \emph{fiberwise $\mathscr F$-regular} at~$y$ (or fiberwise
regular at $y$ if $\mathscr F=\mathscr O_Y$)
if the family~$(f_{i|Y_x})$
of analytic functions on the~$\hr x$-analytic space~$Y_x$ is $\mathscr F_{Y_x}$-regular at~$y$. 

Let~$\mathscr I$ be a coherent sheaf of ideals on~$Y$ 
containing the~$f_i$'s. It gives rise to two 
a priori different
objects on the fiber~$Y_x$: its pull-back~$\mathscr I_{Y_x}$ on~$Y_x$
\emph{as an abstract coherent sheaf}, and a coherent sheaf of \emph{ideals}
on~$Y_x$, namely the image of the natural map~$\mathscr I_{Y_x}\to \mathscr O_{Y_x}$,
which is not injective in general. 

We say that the~$f_i$'s generate~$\mathscr I$ fiberwise at~$y$
\emph{as a coherent sheaf}
if the~$f_{i|Y_x}$'s generate~$\mathscr I_{|Y_x}$ at~$y$. 
This is equivalent to the fact that the~$f_i$'s generate~$\mathscr I$ at~$y$ 
(consider the morphism
$(a_i)\mapsto \sum a_if_i$ 
from~$\mathscr O_Y ^n$ to  $\mathscr I$ and apply \ref{ss-surj-nakayama}). 

We say that   the~$f_i$'s generate~$\mathscr I$ fiberwise at~$y$
\emph
{as an ideal sheaf}
if the~$f_{i|Y_x}$'s generate the image
of~$\mathscr I_{Y_x}\to \mathscr O_{Y_x}$ at~$y$. Of course, if the~$f_i$'s 
generate~$\mathscr I$ fiberwise at~$y$ as a coherent sheaf 
they generate it fiberwise at~$y$
as an ideal sheaf, but the converse is not true in general. 

\begin{prop}\label{prop-regular-seq}
Let~$Y\to X$ be a morphism
of~$k$-analytic spaces, and
let~$\mathscr I$ be a coherent sheaf of ideals
on~$Y$. Let~$(f_1,\ldots, f_n)$
be a family of analytic functions on $Y$ that belong
to $\mathscr I$, and let $Z$ denote
the closed analytic subspace of $Y$ defined by $\mathscr I$.

\begin{enumerate}[1]
\item The set of points
of~$Z$ at which~$(f_i)$ is fiberwise regular and
generates $\mathscr I$ fiberwise
as an ideal sheaf
is a locally
constructible subset of~$Z$. 

\item Let $y$ be a point of $Z$ at which $Y$ is~$X$-flat.
The following are equivalent:

\begin{enumerate}[j]
\item
$(f_i)$ generates $\mathscr I$
fiberwise
as an ideal sheaf
at~$y$, is fiberwise regular at~$y$,
and~$Z$ is~$X$-flat at~$y$.

\item $(f_i)$ generates~$\mathscr I$ at~$y$ and is fiberwise regular at~$y$. 

\end{enumerate}
\item
If~$Y$ is~$X$-flat, 
the locus of validity in $Z$
of the equivalent properties \textnormal{(i)} and \textnormal{(ii)}
of \textnormal{(2)}
is a Zariski-open subset
of~$Z$. 
\end{enumerate}

\end{prop}

\begin{proof}
For every~$i$, we
denote by~$\mathscr F_i$ the coherent sheaf~$\mathscr O_Y/(f_1,\ldots,f_{i-1})$
and by~$h_i$ the endomorphism~$a\mapsto f_ia$ of~$\mathscr F_i$. 

Let us prove (1). Denote by $D$
the complex
\[\xymatrix{\mathscr O_Y^n\ar[rrr]^(0.5){(a_i)\mapsto \sum a_i f_i}&&&\mathscr O_Y\ar[r]&\mathscr O_Y/\mathscr I\ar[r]&0}\]
of coherent sheaves on $Y$. 
Let $y$ be a point of $Z$. The family $(f_1,\ldots, f_n)$ generates $\mathscr I$
fiberwise as an ideal sheaf
at~$y$
if and only if~$D$ is fiberwise exact at~$y$, and~$(f_1,\ldots,f_n)$
is fiberwise regular at~$y$ if and only if every~$h_i$ is fiberwise injective at~$y$. 
Therefore assertion (1) comes from Theorem \ref{thm-fiberwise-exact}.

Let us now prove (2). Suppose that (i) holds. 
Since~$Z$ is~$X$-flat at~$y$, the coherent
sheaf~$\mathscr O_Y/\mathscr I$ is~$X$-flat at~$y$. It follows then
from Proposition \ref{prop-flat-univexact} (2)
that the exact sequence~$$0\to \mathscr I\to \mathscr O_Y\to \mathscr O_Y/\mathscr I\to 0$$ 
is fiberwise exact at~$y$; hence the natural map
$\mathscr I_{Y_{\phi(y)}}\to \mathscr O_{Y_{\phi(y)}}$ is injective 
at~$y$. Therefore the~$f_i$'s 
generate~$\mathscr I_{|Y_{\phi(y)}}$ at~$y$; 
this implies that
they generate~$\mathscr I$ at~$y$ (\ref{ss-generation-i}), so~(ii)
holds. Suppose conversely that~(ii) holds. 
Since the~$f_i$'s 
generate~$\mathscr I$ at~$y$, they generate it fiberwise as an ideal sheaf at~$y$. 
Since all~$h_i$'s
are fiberwise injective at~$y$ and s$\mathscr O_Y$ is~$X$-flat
at $y$, 
a repeated application of Lemma~\ref{lem-flattrick-reg}
ensures that~$\mathscr O_Y/(f_1,\ldots, f_n)$ is~$X$-flat at~$y$. But since
the~$f_i$'s generate~$\mathscr I$ at~$y$, this means that~$Z$ is~$X$-flat at~$y$, 
so (i) holds. 

Let us now prove (3). We denote by~$\mathscr K$ the cokernel of~$$\xymatrix{\mathscr O_Y^n\ar[rrr]^(0.5){(a_i)\mapsto \sum a_i f_i}&&&\mathscr I}.$$
The locus of validity of (ii)  on~$Z$ is the set
of points of~$Z$ at which the following
properties are satisfied: 

\begin{enumerate}[a]

\item $\mathscr K$ is zero; 

\item every~$h_i$ is fiberwise injective. 
\end{enumerate}

Since~$\mathscr O_Y$ is~$X$-flat, a repeated application Lemma~\ref{lem-flattrick-reg}
shows that at every point of $Z$ at which~(b) is true, 
the coherent sheaves~$\mathscr F_i$ are~$X$-flat. Therefore~(b) is equivalent to

\begin{enumerate}[a]\setcounter{enumi}{2}
\item every~$\mathscr F_i$ is~$X$-flat and every~$h_i$ is fiberwise injective.
\end{enumerate}

The locus of validity of (a) on~$Z$ is Zariski-open, and that of~(c) is also Zariski-open
by Proposition \ref{prop-exactflat-locus}. Therefore the set of points of~$Z$ at which the equivalent
conditions~(i) and (ii) are fulfilled is Zariski-open.
\end{proof}

\section{The main theorem}\label{s-cons-maintheo}

This section is essentially
devoted to the (lengthy) proof of Theorem \ref{thm-constloc-main}. 
The latter establishes the local constructibility of the loci of fiberwise validity of the usual algebraic properties, as well as Zariski-openness
results under extra
assumptions (always involving flatness), thus providing analytic counterparts of classical scheme-theoretic results
(which we \emph{do not use} in the proof); \cf the following statements of
\cite{ega43}:

\begin{itemize}[label=$\bullet$]
\item \Prop 9.9.2 (iv) (v) (vii) (viii) (ix);
\item \Prop 9.9.4 (i) (ii) (iii) (iv);
\item  \Th 12.1.1 (v) (vi) (vii); 
\item \Th 12.1.6 (i) (ii) (iii)
(iv);
\item \Cor 12.1.7. 
\end{itemize}

We then prove two additional theorems. The first one is Theorem \ref{thm-constloc-alg}, 
which essentially asserts that some of the constructible loci exhibited at various parts
of this memoir
(Theorem \ref{thm-fiber-closure}, Proposition \ref{prop-relat-codim}
Theorem \ref{thm-constloc-main}) are algebraizable as soon as the source
space $Y$
and all
data living on it are algebraizable and the map
$Y\to X$ under investigation is not too ``widely analytic".
The proof is quite easy and consists in reducing
to Theorem \ref{thm-constloc-main} by using GAGA and the 
extension of coherent sheaves (for dense open immersions
in scheme theory). 

The second one is Theorem \ref{thm-locus-target}, which roughly speaking
turns all local constructibility
or Zariski-openness assertions of Theorem \ref{thm-constloc-main}
into local constructibility
or  Zariski-openness assertions \emph{on the target}
when the map involved is proper. It rests on Kiehl's theorem on 
the direct images of coherent sheaves by a proper map
(which ensures that a proper map
is closed for the Zariski topologies involved; see \ref{ss-kiehl-proper}) and on our
``proper Chevalley theorem" (Theorem \ref{thm-chevalley-proper}).

\begin{lemm}\label{lem-before-main}
Let~$X$ be a
$k$-analytic space, let
$Y$ be its
non-regular locus and let~$\mathscr E$ be a coherent sheaf on~$X$. 
Let $x$ be a point of $X$ and let $m$
be a non-negative integer. 

\begin{enumerate}[1]

\item For any~$n$, the subset~$U_n$
of $X$ consisting of points at which codepth of~$\mathscr E$ at~$x$ is bounded
above by~$n$ is a Zariski-open subset of~$X$.

\item The coherent sheaf
$\mathscr E$ is~$S_m$ at~$x$ if and only if
\[\codim_x (X\setminus U_n, \supp E)>n+m\]
for every $n$. 

\item The space~$X$ is~$R_m$ at~$x$ if and only if~$\codim_x (Y,X) >m$.
\end{enumerate}
\end{lemm}

\begin{proof}
We can assume that~$X$ is affinoid. By GAGA principles
(Lemma \ref{gaga-concrete} and \ref{ss-codim-basics} (2))
we reduce to the corresponding scheme-theoretic statements on $X\al$. 

Now~(1) comes from the fact that $X\al$ is isomorphic to a closed
subscheme of a regular scheme
(see \ref{ss-embed-regular} (1))
and from a theorem by Auslander; \cf \cite{ega42}, \Prop 6.11.2 (i). 
Assertion (2)
comes from from \Prop 5.7.4 (i) of \cite{ega42}, and assertion (3)
from the definition of an $R_m$-scheme. 
\end{proof}

\begin{theo}\label{thm-constloc-main}
Let~$X$ be a~$k$-analytic space, let~$Y$
be an~$X$-analytic space, let~$\mathscr E$ be a coherent sheaf on~$Y$, let~$n$ and~$d$ be two
non-negative integers. 
Let us consider the following subsets of~$Y$ (fiberwise notions are understood
with respect to~$X$). 

\begin{itemize}[label=$\bullet$]

\item The set~$A_n$ of points at which~$\mathscr E$ is fiberwise of codepth~$n$.

\item The set~$A'_n$ of points at which~$\mathscr E$ is $X$-flat
and fiberwise of codepth~$\leq n$. 

\item The set $A_\infty$ of points  at which~$\mathscr E$ is fiberwise~CM. 

\item The set $A'_\infty$ of points  at which~$\mathscr E$ is $X$-flat and
fiberwise~CM; \ie, CM over $X$ in the sense of Definition \ref{def-rel-cm}.

\item The set~$B_n$ of points  at which~$\mathscr E$ is fiberwise~$S_n$. 

\item The set~$B'_n$ of points at which~$\mathscr E$ is $X$-flat
and fiberwise $S_n$.

\item The set~$C$ of points at which~$Y$ is fiberwise Gorenstein. 

\item The set~$C'$ of points at which~$Y$ is $X$-flat
and fiberwise Gorenstein. 

\item The set~$D$ of points at which~$Y$ is fiberwise CI. 

\item The set~$D'$ of points at which~$Y$ is $X$-flat and
fiberwise CI. 

\item The set~$E_n$ of points at which~$Y$ is fiberwise
geometrically~$R_n$.

\item The set~$E'_n$ of points at which~$Y$ is $X$-flat and
fiberwise geometrically~$R_n$

\item The set~$E_\infty$ of points at which~$Y$ is fiberwise
quasi-smooth. 

\item The set $E'_{\infty,d}$ of points at which~$Y$ is 
quasi-smooth of relative
dimension~$d$ over~$X$. 

\item The set $E'_\infty$ of points at which $Y$ is quasi-smooth
over $X$. 

\item The set $\Delta$ of points at which $Y$ is fiberwise geometrically
reduced.

\item The set $\Delta'$ of points at which $Y$ is $X$-flat and
fiberwise geometrically reduced.

\item The set $\Theta$ of points at which $Y$ is fiberwise geometrically
normal.
\item The set $\Theta'$ of points at which $Y$ is $X$-flat
and fiberwise geometrically
normal. 
\end{itemize} 

The sets~$A_n$, $A_\infty$, $B_n$, $C$, $D$, $E_n$, $E_\infty$, $\Delta$, and $\Theta$
are locally
constructible (hence constructible if
$Y$ is finite dimensional, by Proposition \ref{prop-cons-gloc}),
and the sets~$A'_n$, $A'_\infty$, $C'$, $D'$, $E'_{\infty,d}$, and $E'_\infty$
are Zariski-open. If~$\supp E$ is purely of 
relative dimension~$d$ over $X$, then~$B'_n$ is Zariski-open too. If $Y$
is purely of relative dimension $d$ over $X$, then $E'_n, \Delta'$, and $\Theta'$ 
are Zariski-open too. 
\end{theo}

\begin{proof} By Remark \ref{rem-norm-r1s2}, the G-local constructibility of
$\Delta$ will follow from that of $B_1$ and $E_0$, and the G-local
constructibility of
$\Theta$ will follow from that of $B_2$ and $E_1$; analogously, the Zariski-openness 
of $\Delta'$ and $\Theta'$ in the equidimensional case
will respectively follow from that of $B'_1$ and $E'_0$, and from that 
of $B'_2$ and $E'_1$.

\subsubsection{Study of~$A_n, A'_n, C,$
and $C'$: first reductions}\label{sss-ac-first}
Let us denote by~$A_{\leq n}$
the set of points of~$Y$ at which~$\mathscr E$
is fiberwise of codepth~$\leq n$. Instead of proving
directly
the constructibility of~$A_n$, 
we shall prove that of~$A_{\leq n}$; this will
clearly imply 
the former one because~$A_n=A_{\leq n}\setminus A_{\leq n-1}$. 

We are going to
apply some of the general results in \ref{s-diagtr},
which involve two categories $\mathfrak F$ and $\mathfrak C$ 
and a property $\mathsf Q$
as in \ref{ss-general-axioms-pfib}, and a functor
$\mathscr S$ as in \ref{ss-def-functorS}. 
For the study of $A_{\leq n}$ and $A'_n$, we take
$\mathfrak F$
to be the category 
$\mathfrak {Coh}$, $\mathfrak C$ to be the category of all analytic spaces, 
and for $\mathsf Q$ to be the property of being of codepth $\leq n$, and $\mathscr S$
to be the 
identity functor. For the study of $C$ and $C'$,
we take
$\mathfrak F$ to be the category $\mathfrak T$,  
$\mathfrak C$ to be the category of all analytic spaces, 
$\mathsf Q$ to be the property of being of Gorenstein, and $\mathscr S$ to be the 
functor that sends any analytic space to its structure sheaf.

In both cases, we are interested in assertions $(\alpha)$ and $(\beta)$
of \ref{ss-claims}. Since assertion $(\gamma)$ of \loccit~holds
(this is Theorem \ref{thm-flatloc-zaropen}), Lemma
\ref{gamma-betaprime-beta}
ensures that it suffices to prove assertion $(\beta')$ of \ref{ss-auxiliary}. 
And as far as $A_{\leq n}$ and $A'_n$
are concerned, it even suffices to prove assertion
$(\beta^\flat)$ of \loccit, in view of Lemma 
\ref{lem-betadash-beta} and Remark \ref{rem-redstat-fulfilled}.

More explicitly, those reductions 
mean that we
can assume that $Y$ and~$X$ are affinoid and the morphism $Y\to X$ has a 
section~$\sigma$, and we have to prove the following: 

\begin{enumerate}[1]

\item If~$Y$ is $X$-flat then~$\sigma\inv(A'_n)$
is a Zariski-open subset of~$X$. 

\item The subset $\sigma\inv (C')$ of~$X$ is Zariski-open
(without any flatness assumption on $Y$). 

\end{enumerate}

\subsubsection{Some
homological computations}\label{sss-homolog-comput}
We fix a non-negative
integer $d$ such that
the relative dimension of~$\supp  E$ 
over~$X$ is bounded above
by~$d$  (if~$\supp E\neq \emptyset$ one can take~$d=\dim \supp E$). Choose a resolution
\[\mathscr F_{d+1}\to \mathscr F_d\to\ldots\to \mathscr F_1\to \mathscr F_0\to \sigma_\ast\mathscr O_X\to 0,\]
where the~$\mathscr F_i$'s are free~$\mathscr O_Y$-modules of finite rank, and let $\mathsf F$
denote the complex 
\[\mathscr F_{d+1}\to \mathscr F_d\to\ldots\to \mathscr F_1\to \mathscr F_0\to \mathscr F_{-1}=0.\]

Let $x$ be a point of $X$ such $Y$ is~$X$-flat at~$\sigma(x)$. For any $i\in\{0,\ldots,d\}$ 
we have natural isomorphisms 

\begin{eqnarray}\setcounter{equation}{1}
\mathrm{Ext}^i_{\mathscr O_{Y_x,\sigma(x)}}(\hr x,\mathscr E_{Y_x,\sigma(x)})
&\simeq &\H_i\left(\mathrm{Hom}(\mathsf F_{Y_x,\sigma(x)},\mathscr E_{Y_x,\sigma(x)})\right)\\
&\simeq& \mathscr H_i\left(\mathscr H om(\mathsf F_{Y_x}, \mathscr E_{Y_x})\right)_{\sigma(x)}\\
&\simeq & \mathscr H_i\left(\mathscr H om(\mathsf  F, \mathscr E)_{Y_x}\right)_{\sigma(x)}
\end{eqnarray}

Indeed, since $Y$ is $X$-flat at $\sigma(x)$, so is $\mathscr F_i$ for every $i$ (because it is free
over $\mathscr O_Y$): morevoer, $\sigma_\ast X$ is $X$-flat (everywhere). Proposition \ref{prop-flat-univexact} then
ensures that 
\[\mathscr F_{d+1,Y_x,\sigma(x)}\to \mathscr F_{d,Y_x,\sigma(x)}\to\ldots\to \mathscr F_{1,Y_x,\sigma(x)}\to (\sigma_\ast \mathscr O_X)_{Y_x,\sigma(x)}\to 0\]
is exact. Isomorphism (a)
now follows since $\mathscr F_{i,Y_x,\sigma(x)}$ is for every $i$
a free $\mathscr O_{Y_x,\sigma(x)}$-module (because $\mathscr F_i$ is free over $\mathscr O_Y$) and
the $\mathscr O_{Y_x,\sigma(x)}$-module $(\sigma_\ast \mathscr O_X)_{Y_x,\sigma(x)}$ is nothing
but the residue field $\hr x$ of $\mathscr O_{Y_x,\sigma(x)}$. Isomorphism (b) comes from 
the coherence of the sheaves involved, 
and isomorphism (c) is due to the fact that $\mathscr Hom(\mathscr F_i,\mathscr E)_{Y_x}=
\mathscr Hom(\mathscr F_{i,Y_x},\mathscr E_{Y_x})$ for every $i$, again by freeness of $\mathscr F_i$. 

We thus have
\begin{eqnarray}
\mathrm{Ext}^i_{\mathscr O_{Y_x,\sigma(x)}}(\hr x,\mathscr E_{Y_x,\sigma(x)})&\simeq&
\mathscr H_i\left(\mathscr H om(\mathsf  F, \mathscr E)_{Y_x}\right)_{\sigma(x)}\\
&=&\mathscr H_i\left(\mathscr H om(\mathsf  F, \mathscr E)_{Y_x}\right)_{\hr x}
\end{eqnarray}
(equality (e) comes from the fact 
that the $\mathscr O_{Y_x,\sigma(x)}$-module $\mathscr H_i\left(\mathscr H om(\mathsf  F, \mathscr E)_{Y_x}\right)_{\sigma(x)}$
is in fact an $\hr x$-vector space by (d)). 
In view of the description of depth through $\mathrm{Ext}$ functors (\cf \Th
16.7 of \cite{matsumura1986}), we deduce from (d) that
\begin{equation}
\mathrm{depth}_{\mathscr O_{Y_x,\sigma(x)}}\mathscr E_{Y_x,\sigma(x)}=
\inf \left\{i|\;\mathscr H_i\left(\mathscr H om(\mathsf  F, \mathscr E)_{Y_x}\right)_{\sigma(x)}\neq 0\right\}
\end{equation}
(note that this is true even if $\mathscr E_{Y_x,\sigma(x)}=0$ because in this case 
its depth is equal to $+\infty$). 

\subsubsection{}\label{sss-open-lambda}
Let $a$ be an element of $\in \{0,\ldots, d+1\}$, and
let $\mathsf F^{\leq a}$ be
the truncated complex 
$\mathscr F_{a}\to \mathscr F_{a-1}\to\ldots\to \mathscr F_1\to \mathscr F_0\to \mathscr F_{-1}=0$. We
denote
by $\Lambda_a$ the subset of $Y$ consisting of points at which $\mathscr Hom(\mathsf F^{\leq a},\mathscr E)$
is fiberwise exact and $\mathscr E$ is $X$-flat. 

For every $i\in\{-1,\ldots, d+1\}$, the coherent sheaf $\mathscr F_i$ is free over $\mathscr O_Y$, hence
$\mathscr Hom(\mathscr F_i,\mathscr E)$ is a direct sum of finitely many copies of $\mathscr E$.
Moreover, since $\mathscr F_0$
surjects onto $\sigma_\ast \mathscr O_X$, it is of positive rank
as soon as $X\neq \varnothing$, which is always the case when $Y\neq \varnothing$. 
We have therefore the equality
\[\uflat {\mathscr E}X=\uflat {\mathscr Hom(\mathscr F_0,\mathscr E)}X=\bigcap_{-1\leq i\leq a}\uflat  {\mathscr Hom(\mathscr F_i,\mathscr E)}X.\]

We can therefore
describe $\Lambda_a$ as the subset of $Y$ consisting of points at which
the complex $\mathscr Hom(\mathsf F^{\leq a},\mathscr E)$ is fiberwise exact
and at which every sheaf involved in this complex is $X$-flat. 
It follows then from Theorem \ref{thm-exactflat-locus}
that $\Lambda_a$ is a Zariski-open subset of $Y$.

\subsubsection{Proof of assertion \textnormal{(1)} of \ref{sss-ac-first}}\label{sss-conclu-a}
We assume that~$Y$ is~$X$-flat, and we are going to prove that~$\sigma\inv(A'_n)$
is a Zariski-open subset of~$X$. Let us begin with a remark. If $y$ is a point of $Y$ and
if $x$ denotes its image in $X$, we define the fiberwise
\emph{depth}
of $\mathscr E$ at $y$ as the depth of the $\mathscr O_{Y_x,y}$-module
$\mathscr E_{Y_x,y}$. If $y\notin \supp E$, this fiberwise depth is equal to $+\infty$; 
if $y\in \supp E$, it 
is bounded by $\dim_{\mathrm{Krull}} \mathscr O_{\supp E_x,y}$, hence by $d$ in view
of Corollary \ref{cor-interp-centdim}. 

Let~$m$ be an integer such that $0\leq m\leq d$. We denote by $F_m$
the subset of $Y$ consisting of points 
at which
$\mathscr E$ is $X$-flat and of fiberwise depth $\geq m$.  
Let us first show that~$\sigma\inv(F_m)$ is a Zariski-open
subset of~$X$. Let $x$ be a point of $X$.
By equation (f) of \ref{sss-homolog-comput},
$\sigma(x)\in F_m$ if and only if the two following conditions are fulfilled:

\begin{enumerate}[i]
\item The complex $\mathscr Hom(\mathsf F^{\leq m+1},\mathscr E)$
is fiberwise exact at $\sigma(x)$. 
\item The coherent sheaf $\mathscr E$ is $X$-flat at $\sigma(x)$. 
\end{enumerate}
 Otherwise said, $\sigma\inv(F_m)=\sigma\inv(\Lambda_{m+1})$.
Since $\Lambda_{m+1}$
is a Zariski-open subset of $Y$ (\ref{sss-open-lambda}),  $\sigma\inv(F_m)$ is Zariski-open subset of $X$, 
as desired. 

For
every~$\delta\in \N$, we
denote by~$H_\delta$, \resp $H_{\leq \delta}$, 
the set of points of~$\supp E$ at which the relative dimension of~$\supp E$ over~$X$ is equal
to~$\delta$ (\resp $\leq \delta$). By Zariski-semi-continuity of the relative dimension
(\cite{ducros2007}, \Th 4.9), the sets $H_{\leq \delta}$ and $H_{\delta}$
are respectively Zariski-open and constructible in $\supp E$. 

Let~$x\in \sigma\inv(\supp E)$. Let~$\delta$
and $m$ be respectively the relative dimension and the fiberwise codepth 
of~$\supp E$ at~$\sigma(x)$.
Since $\hr{\sigma(x)}=\hr x$,
the Krull dimension of~$\mathscr O_{\supp E_x,\sigma(x)}$ is equal to~$\delta$
(Corollary \ref{cor-interp-centdim}); the fiberwise codepth of~$\mathscr E$
at~$\sigma(x)$ is thus equal to~$\delta-m$. 

As a consequence, 
\[\sigma\inv(A'_n)=\left(X\setminus \sigma\inv (\supp E)\right)\bigcup \bigcap_{\delta \leq d, m\leq d, \delta-m\leq n} \sigma\inv(F_m\cap H_{\leq \delta}).\]
We have seen that $\sigma \inv(F_m)$ is a Zariski-open subset of $X$
for every $m\leq d$, 
and that $H_{\leq \delta}$
is a Zariski-open subset of $\supp E$
for every $\delta$. Therefore
$\bigcap_{\delta \leq d, m\leq d, \delta-m\leq n} \sigma\inv(F_m\cap H_{\leq \delta})$ is
a Zariski-open subset of $\sigma\inv(\supp E)$,
so $\sigma\inv(A'_n)$ is a Zariski-open subset of $X$. 
We have thus proved assertion (1) of \ref{sss-ac-first}.

\subsubsection{Proof of assertion \textnormal{(2)} of \ref{sss-ac-first}}
The affinoid space $Y$ is no longer assumed to be $X$-flat, 
and we are going to prove that~$\sigma\inv(C')$ is a Zariski-open
subset of~$X$. We will use the results and notation of~\ref{sss-ac-first}, \ref{sss-homolog-comput}
\ref{sss-open-lambda}, and \ref{sss-conclu-a} for $\mathscr E=\mathscr O_Y$; so $d, H_\delta$, $H_{\leq \delta}$, and $F_m$ refer from now on
to the fiberwise dimension 
and depth of $Y$ over $X$, and $\Lambda_a$ involves the complex $\mathscr Hom(\mathsf F^{\leq a},\mathscr O_Y)$.

We fix an integer~$\delta\leq d$, and we are going to
show that~$\sigma \inv(C'\cap H_\delta)$ 
is a Zariski-open subset of~$X$, which will be sufficient since $X=\bigcup_{\delta \leq d}\sigma \inv(H_\delta)$. 
Let $U$ be the subset of $X$ consisting
of points $x$ such that $\sigma(x)\in H_{\leq \delta}$
and
the $\hr x$-vector space~$\mathrm{Ext}^i_{\mathscr O_{Y_x,\sigma(x)}}(\hr x, \mathscr O_{Y_x,\sigma(x)})=
\mathscr H_i\left(\mathscr H om(\mathsf  F, \mathscr E)_{Y_x}\right)_{\sigma(x)}$
(see equation (d) of \ref{sss-homolog-comput})
has rank~$0$ for~$i<\delta$ and rank~$1$ for~$i=\delta$. 

Let $x$ be a point of $U$. The depth of $\mathscr O_{Y_x,\sigma(x)}$ is equal to~$\delta$ by equation (f)
of \ref{sss-homolog-comput}; therefore the Krull dimension of that ring is at least~$\delta$. Since 
$\dim_{\sigma(x)}Y_x\leq \delta$ by definition of $U$, it follows from 
Corollary  \ref{cor-interp-centdim}
that
\[\dim_{\mathrm{Krull}}\mathscr O_{Y_x,\sigma(x)}=\dim_{\sigma(x)}Y_x=\delta.\] 

 We thus
see that~$U\cap \sigma\inv(\uflat YX)$ coincides with~$\sigma \inv(C'\cap H_\delta)$; it remains to show that~$U\cap \sigma\inv(\uflat YX)$
is a Zariski-open subset of~$X$. 
By definition of $U$
and equation (f) of \ref{sss-homolog-comput}
we have the equality 
\[U\cap \sigma\inv(\uflat YX)=\sigma \inv(H_{\leq \delta}\cap \uflat YX \cap V\cap V'),\] where: 

\begin{itemize}[label=$\bullet$]

\item $V$ is the subset of $Y$ consisting of points at which
the complex~\[\mathscr Hom(\mathsf F^{\leq \delta+1},\mathscr O_Y)\] is fiberwise exact
in degrees $<\delta$. 

\item $V'$  is the subset of $Y$ consisting of points $y$ such that 
\[\rk y(\mathscr H_\delta(\mathscr Hom(\mathsf F^{\leq \delta+1}_{Y_x},\mathscr E_{Y_x})))=1,\]
with
$x$ the image of $y$ in $X$. 
\end{itemize}
But we also have 
\[U\cap \sigma\inv(\uflat YX)=\sigma \inv(H_{\leq \delta}\cap \uflat YX \cap V\cap V'')\]
where~$V''$ is  the set of points $y$ of~$Y$ such that
\[\rk y(\mathscr H_\delta(\mathscr Hom(\mathsf F^{\leq \delta+1}_{Y_x},\mathscr E_{Y_x})))\leq 1,\]
for $x$ the image of $y$ in $X$.  

 Indeed, let $x$ be a point of $\sigma \inv(H_{\leq \delta}\cap \uflat YX \cap V\cap V'')$. 
The Krull dimension of $\mathscr O_{Y_x,\sigma(x)}$ is bounded above by $\delta$, 
so its depth is too. Now since $\sigma(x)\in V$, it follows from equation (e)
and (f) of \ref{sss-homolog-comput}
that the depth of $\mathscr O_{Y_x,\sigma(x)}$ is $\geq \delta$. It is therefore
equal to $\delta$, which implies by \ref{sss-homolog-comput}
that 
\[\mathscr H_\delta(\mathscr Hom(\mathsf F^{\leq \delta+1}_{Y_x},\mathscr E_{Y_x}))_{\hr x}\neq 0.\]
As a consequence, $\rk {\sigma(x)}(\mathscr H_\delta(\mathscr Hom(\mathsf F^{\leq \delta+1}_{Y_x},\mathscr E_{Y_x})))=1$
and  the point $x$ belongs to $\sigma \inv(H_{\leq \delta}\cap \uflat YX \cap V\cap V')$, whence our claim. 

We know that $H_{\leq \delta}$ is a Zariski-open subset of~$Y$. 
The intersection $V\cap \uflat Y X$ is nothing but the subset $\Lambda_\delta$ of $Y$, which is Zariski-open 
(\ref{sss-open-lambda}).
It follows from Proposition \ref{prop-flat-univexact} (4)
(and from the fact that $\uflat YX\subset \bigcap_i \uflat{\mathscr F_i}X$, see \ref{sss-open-lambda})
that $V''\cap \uflat YX$ is the set of points $y$ of~$\uflat YX$ 
such that $\rk y (\mathscr H_\delta(\mathscr Hom(\mathsf F, \mathscr O_Y)))\leq 1$. 
By Zariski-upper semi-continuity of the  pointwise rank of a coherent sheaf
(and Theorem \ref{thm-flatloc-zaropen}),
$W''\cap \uflat YX$ is a Zariski-open subset of~$Y$. 
This implies
that~$H_{\leq \delta}\cap \uflat YX\cap W\cap W''$ is a Zariski-open subset of~$Y$; therefore
\[\sigma\inv(C'\cap H_\delta)=U\cap \sigma\inv(\uflat YX)=\sigma \inv(H_{\leq \delta}\cap \uflat YX \cap V\cap V'')\]
is a Zariski-open subset of $X$, and assertion (2) of \ref{sss-ac-first}
is proved. 

\subsubsection{Study of $A_\infty$ and $A'_\infty$}
By arguing G-locally on $Y$ we may assume that
it is good and that
there is an integer $\delta$ such that $\dim Y_x$ is bounded by $\delta$
for every $x\in X$; we then also have $\dim_{\mathrm{Krull}}\mathscr O_{Y_x,y}\leq \delta$
for every $y\in Y_x$ (Corollary \ref{cor-interp-centdim}). We thus see that $A_m=Y$
for $m>\delta$. As a consequence, 
\[A_\infty=\bigcap_{m\leq \delta}A_m, \;A'_\infty=\bigcap_{m\leq \delta}A'_m.\]
We have already proven that the $A_m$'s are locally constructible, and that the $A'_m$
are Zariski-open.  
Therefore $A_\infty$ is locally constructible, and
$A'_\infty$ is Zariski-open.

\subsubsection{Study of~$B_n$ and~$B'_n$} 
We begin with the local constructibility of $B_n$. 
By arguing
locally, we may assume that~$Y$
is finite-dimensional. Under this assumption, the subset $A_m$ of~$Y$ is 
empty for all but finitely many
$m$, and any G-locally constructible
subset of $Y$ is constructible by Proposition \ref{prop-cons-gloc}.
It follows then from the
(G-local) constructibility of $A_m$ for every $m$
(which has already been proved), 
Lemma~\ref{lem-before-main} (2),
and
the (G-local) constructibility of the fiberwise
codimension function (Proposition \ref{prop-relat-codim}) that~$B_n$
is constructible. 

Let us assume now that $\supp E$ is
of pure
relative dimension~$d$ over $X$, and let us prove that $B'_n$ is open. 
We are going to
apply some of the general results in \ref{s-diagtr},
which involve two categories $\mathfrak F$ and $\mathfrak C$ 
and a property $\mathsf Q$
as in \ref{ss-general-axioms-pfib}, and a functor
$\mathscr S$ as in \ref{ss-def-functorS}. 
We take
$\mathfrak F$ to be the category 
$\mathfrak {Coh}$, $\mathfrak C$ to be the category of all analytic spaces, 
$\mathsf Q$ to be the property of being of $S_n$, and $\mathscr S$ to be the 
identity functor. We are interested in assertion $(\beta)$
of \ref{ss-claims}. Since assertions $(\alpha)$ and $(\gamma)$ of \loccit
already hold ($(\alpha)$ is the G-local constructibility of $B_n$ we have just established, 
and $(\gamma)$ is Theorem \ref{thm-flatloc-zaropen}), Lemma \ref{alpha-gamma-beta}
ensures that it is sufficient to prove
assertion $(\beta^\flat)$ of \ref{ss-auxiliary}; \ie, 
we may assume that $Y$ and $X$ are affinoid and $\mathscr E$ is $X$-flat
(the proof of Lemma \ref{alpha-gamma-beta}
only involves arguing G-locally, hence it does not
modify the dimension of~$\supp E\to X$). 

Let $m$ be a non-negative integer. We have already proved that $A'_m$
is Zariski-open; hence $Y\setminus A'_m$ is a Zariski-closed subset of $Y$
which is contained in $\supp E$. Let
$J_m$ be
the set of points of~$\supp E$ at which the fiberwise codimension
of~$Y\setminus A'_m$ in~$\supp E$ 
is~$>n+m$. By Zariski-upper semi-continuity of the
dimension of a morphism (\cite{ducros2007}, \Th 4.9),
and because~$\supp E \to X$ is of pure dimension~$d$, the subset~$J_m$ 
of~$Y$ is Zariski-open \emph{in the Zariski-closed subset~$\supp E$}; note that it
is equal to the whole
of $\supp E$ for large enough~$m$. 
By Lemma \ref{lem-before-main} (2) and
$X$-flatness of $\mathscr E$, 
the
set $B'_n$ is equal to
\[(Y\setminus \supp E)\bigcup \;\bigcap_m J_m\]
and
hence is Zariski-open in $Y$. 

\subsubsection{Study of~$D$ and~$D'$} 
We may assume that~$Y$ and~$X$ are
affinoid. Under this assumption, 
there exists an 
affinoid~$X$-space $X'$
which is flat with regular fibers, 
such that~$Y$ can be identified (over~$X$)
with a closed analytic subspace
of~$X'$;  one can take for 
example $X'$ to be the product of~$X$ and a suitable compact polydisc, but
we want to emphasize that what follows holds for an arbitrary such~$X'$. 
Let~$\mathscr I$ be the ideal sheaf on~$X'$ that corresponds to~$Y$; 
we fix a system of global sections~$f_1,\ldots, f_n$ generating~$\mathscr I$. 

For every~$J\subset \{1,\ldots, n\}$, let~$ P_J$ be the set of points
of~$Y$ at which~$(f_i)_{i\in J}$ generates~$\mathscr  I$ fiberwise
\emph
{as an ideal} (\ref{ss-generation-i})
and at which~$(f_i)_{i\in J}$ is fiberwise regular (Definition \ref{def-regular-seq}).
It follows from
Proposition \ref{prop-regular-seq}
(1) that~$P_J$ is a constructible subset of~$Y$. 

Let us prove that $D$ is the union of the $P_J$'s, hence is constructible. 
Let $y$ be a point of $\bigcup P_J$ and let $x$ be its image in $X$. 
By assumption, there exists $J$ such that $(f_i)_{i\in J}$ is a regular sequence
of the regular local ring $\mathscr O_{X'_x, y}$ and $\mathscr O_{Y_x,y}=\mathscr O_{X'_x,y}/(f_i)_{i\in J}$. 
Then $\mathscr O_{Y_x,y}$ is CI by definition, hence $y\in D$. 
Assume conversely that $y$ lies on $D$. Let
$J$ be a minimal subset of $\{1,\ldots, n\}$ such that the $f_i$'s
for $i\in J$ generate $\mathscr I$ fiberwise at $y$; \ie, the $f_i$'s
for $i\in J$ generate the kernel of $\mathscr O_{X'_x,y}\to \mathscr O_{Y_x,y}$.
Since $y\in D$, the local ring
$\mathscr O_{Y_x,y}$
is CI. Since 
$\mathscr O_{X'_x,y}$ is regular, the family 
$(f_i)_{i\in J}$ is a regular sequence of the local ring $\mathscr O_{X'_x,y}$, 
\cf \cite{ matsumura1986}, 17.4(iii), (1) $\Leftrightarrow$ (3); 
\ie,
$(f_i)_{i\in J}$ is fiberwise regular at $y$ and $y\in P_J$. 

Now $D'=D\cap \uflat YX$. By the above, 
$D'$ is the union of the~$\uflat YX\cap P_J$'s
for~$J\subset\{1,\ldots, n\}$. Since~$X'\to X$ is
flat,
$\uflat YX \cap P_J$ is a Zariski-open subset of~$Y$ for every~$J$ by
Proposition \ref{prop-regular-seq}. Therefore~$D'$
is a Zariski-open subset of~$Y$. 

\subsubsection{Study of~$E_n, E'_n, E_\infty, E'_\infty$ and~$E'_{\infty,d}$} 
For every~$\delta$, 
let~$E_{\infty, \delta}$ be the set of points of~$Y$ at which~$Y$ 
is fiberwise quasi-smooth of dimension~$\delta$. This is the set of points
at which~$Y\to X$ is of dimension~$\delta$
and $\Omega_{Y/X}$ has fiber rank~$\delta$;  
it is therefore locally
constructible by Theorem 4.9
of~\cite{ducros2007}. Now~$E_\infty$ is
the union
of all~$E_{\infty,\delta}$'s; since any compact analytic domain of~$Y$ intersects only finitely many~$E_{\infty, \delta}$'s, we see
that~$E_\infty$ is locally constructible.

The set~$E'_{\infty,d}$ is the intersection of two subsets of~$Y$: 

\begin{itemize}[label=$\bullet$]
\item the set~$E_{\infty,d}$, which we have just seen is constructible; 

\item the~$X$-flatness locus of~$Y$, which is Zariski-open by Theorem \ref{thm-flatloc-zaropen}.
\end{itemize}
Therefore~$E'_{\infty,d}$ is locally
constructible. By Corollary \ref{coro-qsm-open}, it is
also open; it is therefore Zariski-open by Lemma~\ref{lem-adh-const} (4).
It follows that
$E'_{\infty}=\bigcup_{\delta} E'_{\infty,\delta}$ is Zariski-open.

For any $x\in X$, the intersection $E_\infty \cap Y_x$ is nothing but the set $E'_\infty$ understood
with respect to the map $Y_x\to \mathscr M(\hr x)$: it is thus Zariski-open
in $Y_x$ (this was already known; 
see \cite{ducros2009}, \Prop 6.6). 
By Lemma~\ref{lem-before-main} (3), 
one
can describe~$E_n$ as the set of points of~$Y$ at which the fiberwise
codimension of~$Y\setminus E_\infty$ in~$Y$ is~$>n$. By Proposition
\ref{prop-relat-codim}, this is a constructible subset of~$Y$. 

Let us assume now that $Y$ is
of pure
relative dimension~$d$ over $X$, and let us prove that $E'_n$ is open. 
We are going to
apply some of the general results in \ref{s-diagtr},
which involve two categories $\mathfrak F$ and $\mathfrak C$ 
and a property $\mathsf Q$
as in \ref{ss-general-axioms-pfib}, and a functor
$\mathscr S$ as in \ref{ss-def-functorS}. 
We take
$\mathfrak F$ to be the category 
$\mathfrak T$, $\mathfrak C$ to be the category of all analytic spaces,
$\mathsf Q$ to be  the property of being geometrically $S_n$, and $\mathscr S$ the 
functor that sends any analytic space to its structure sheaf.
We are interested in assertion $(\beta)$
of \ref{ss-claims}. Since assertions $(\alpha)$ and $(\gamma)$ of \loccit
already hold ($(\alpha)$ is the G-local constructibility of $E_n$ we have just established, 
and $(\gamma)$ is Theorem \ref{thm-flatloc-zaropen}), Lemma \ref{alpha-gamma-beta}
ensures that it is sufficient to prove
assertion $(\beta^\flat)$ of \ref{ss-auxiliary}; \ie, 
we may assume that $Y$ and $X$ are affinoid and $Y$ is $X$-flat
(the proof of Lemma \ref{alpha-gamma-beta}
only involves arguing G-locally, hence it
does not
modify the dimension of~$Y\to X$). 

Now as~$Y\to X$ is flat, the 
constructible subset~$E_\infty$ is equal to $E'_\infty$, hence 
is Zariski-open by the above. 
By Lemma~\ref{lem-before-main} (3), one 
can describe~$E'_n$ as the set of points of~$Y$
at which the fiberwise
codimension of~$Y\setminus E_\infty$ in~$Y$ is~$>n$. 
The relative dimension of~$Y\to X$ is equal to~$d$ everywhere, and
the relative
dimension of~$Y\setminus E_\infty\to X$ is Zariski-upper semi-continuous by
\Th 4.9 of \cite{ducros2007}
(we use the fact that~$Y\setminus E_\infty$ is a Zariski-closed
subset of~$Y$). 
It follows that~$E'_n$ is a Zariski-open
subset of~$Y$. 
\end{proof}

\begin{rema}\label{rem-notopti}
We have proved the Zariski-openness of $B'_n$, $E'_n$, $\Delta'$, 
and $\Theta'$ only under the assumption that $Y\to X$ is equidimensional.
We think that this is not optimal. For instance, in scheme theory, the analogues
of $\Delta'$ and $\Theta'$ (for a finitely presented morphism between noetherian schemes)
are Zariski-open without any assumption on the dimension (\cite{ega43}, \Th 1.2.4 (iv) and (v)), 
and so is the analogue of $B'_n$ when $\mathscr E=\mathscr O_Y$ (\Th 2.2.6 (i) of \opcit). 
Note nevertheless that there is a counter-example to the scheme-theoretic version of Zariski-openness
of $E'_1$ when one allows non-equidimensional fibers (\opcit, \Rem12.1.8 (ii)), which can be turned into a counter-example
to Zariski-openness of $E'_1$ in our setting by endowing the ground field with the trivial absolute value and
applying GAGA principles. 

In order to get stronger statements, we should probably investigate
carefully the variation of 
the number
of (local) geometric irreducible or embedded components on which a point lies \emph{it its fiber}, which in turn
woud require for  the ``spreading out" process
a deeper understanding of the links between the local rings of a generic fiber and those
of the source space, far beyond Theorem  \ref{thm-localring-generic}.
\end{rema}

We are now going to show that most constructible loci considered in this
memoir
tend to be algebraizable when one starts from algebraic data.

\begin{theo}\label{thm-constloc-alg}
Let $V$ be a $k$-affioid space, ant let $\mathscr Y$ be a $V\al$-scheme
locally of finite type; set $Y=\mathscr Y\an$. Let $\phi \colon Y\to X$ be a morphism of $k$-analytic spaces
such that (at least) one of the following assertions holds:

\begin{enumerate}[A]
\item The morphism $\phi$ is the composition of the structure map $Y\to V$ and an arbitrary morphism 
$V\to X$.  

\item There exists a $k$-affinoid space $U$, a morphism $V\to U$, and a scheme $\mathscr X$ locally of finite
type over $U\al$ such that $X=\mathscr X\an$ and $\phi$ is induced by 
an $U\al$-map $\mathscr Y\to \mathscr X$.

\end{enumerate}
Let $\mathscr E$ be a coherent
sheaf on $Y$ and let $E$ and $F$ be two locally constructible subsets of $Y$ with $F\subset E$. 
Assume moreover that $\mathscr E$ arises from an algebraic coherent sheaf 
on $\mathscr Y$, 
and both $E$ and $F$ arise from locally
constructible subsets 
of $\mathscr Y$; \ie, $E\al$ and $F\al$ are locally constructible, $E=(E\al)\an$ and $F=(F\al)\an$. 
Let $n$ and $d$ be two integers. 

\begin{enumerate}[1]
\item The following subsets of $Y$ are of the form $P\an$, for $P$ a locally constructible subset of $\mathscr Y$
(note that such a $P$ is constructible as soon as $\mathscr Y$ is of finite type,
or more generally finite-dimensional, see Remark \ref{cons-local-scheme}). 

\begin{enumerate}[b]
\item The set of points $y$ such that $\dim_y \phi=d$. 

\item The subset of $\adht E\phi$ consisting of points at which the fiberwise dimension of $\adht E\phi$ 
over $X$ belongs to a given subset of $\N$; in particular, $\adht E\phi$ iself. 

\item The subset of $\adht E\phi$ consisting of points at which the fiberwise codimension of $\adht F\phi$ 
in $\adht E\phi$ belongs to a given subset of $\N\cup\{+\infty\}$. 

\item The sets $A_n$ , $A_\infty$, $B_n$, $C$,$ D$, $E_n$, $E_\infty$, $\Delta$, and $\Theta$
of Theorem \ref{thm-constloc-main}. 

\end{enumerate}

\item The following sets are of the form $P\an$, for $P$ a Zariski-open subset of $\mathscr Y$:

\begin{enumerate}[b]
\item The set of points $y$ such that $\dim_y \phi\leq d$. 

\item The set of points at which $\mathscr E$ is $X$-flat.

\item The sets $A'_n$ , $A'_\infty$, $C'$, $D'$, $E'_{\infty,d}$, and $E'_\infty$
of Theorem \ref{thm-constloc-main}. 

\item The sets $E'_n$, $\Delta'$, and $\Theta'$ of Theorem \ref{thm-constloc-main}
if $\phi$ is of pure relative dimension $d$. 

\item The set $B'_n$ of Theorem \ref{thm-constloc-main} if $\supp E\to X$
is of pure relative dimension $d$.

\end{enumerate}
\end{enumerate}
\end{theo}

\begin{proof} 
In case (A), the theorem is 
local on $\mathscr Y$, and we can thus assume that the latter admits
a proper compactification $\overline{\mathscr Y}$ over $V\al$. 

In case (B), the theorem is local on $\mathscr Y$ and $\mathscr X$. We can thus assume that $\mathscr X$
admits a proper compactification $\overline {\mathscr X}$ over $U\al$, and 
that $\mathscr Y\to \overline{\mathscr X}\times_{U\al}V\al$ admits a proper compactification 
$\overline{\mathscr Y}$. 

Hence in both cases we can assume that $\mathscr Y$ admits a proper compactification $\overline{\mathscr Y}$
over $V\al$, and that 
$\phi$ extends to a morphism $\overline{\mathscr Y}\an\to X$ (for achieving this in case (B) one first needs to replace
$\mathscr X$ with $\overline{\mathscr X}$); note 
that $E\al$ and $F\al$ are now constructible since $\mathscr Y$ is of finite type 
over $V\al$. 
By \cite{ega1}, \Cor 9.4.8,
the coherent sheaf on $\mathscr Y$ from which 
$\mathscr E$ arises can be extended to
a coherent sheaf on $\overline{\mathscr Y}$; hence $\mathscr E$
can be extended to a coherent sheaf $\overline{\mathscr E}$ on $\overline{\mathscr Y}\an$. Moreover, 
$E$ and $F$ remain constructible inside
$\overline{\mathscr Y}\an$ (in contrast with non-transitivity of the analytic Zariski-topology in general)
because
$E\al$ and $F\al$ are constructible in $\overline{\mathscr Y}$. 

By replacing the scheme $\mathscr Y$ with $\overline{\mathscr Y}$
and the coherent sheaf 
$\mathscr E$ with $\overline{\mathscr E}$ (and still
working with $E$ and $F$), we reduce to the
case where $\mathscr Y$ is proper, \emph{except possibly for \textnormal{(2d)} and 
 \textnormal{(2e)}}, because the relative equidimensionality property
they require might be lost. Now, modulo GAGA:

\begin{itemize}[label=$\bullet$]
\item Cases (1a) and (2a) come from Zariski-upper semi-continuity of the
relative dimension (\cite{ducros2007}, \Th 4.9), case (1b) from Theorem \ref{thm-flatloc-zaropen}, 
and cases (1c) and (2c) from Theorem \ref{thm-constloc-main}.

\item Cases (1c) and (1d) come respectively from Theorem \ref{thm-fiber-closure} and Proposition \ref{prop-relat-codim}.

\end{itemize}

It remains to consider cases (2d) and (2e) (the scheme $\mathscr Y$ is no longer
assumed to be proper). We already know that the sets considered in (1d) and (2b)
are of the form $P\an$ for $P$ a locally
constructible
subset of $\mathscr Y$; it follows that the sets considered in (2d)
and (2e) are also of this form. Moreover,
they are Zariski-open
by Theorem \ref{thm-constloc-main}; 
in particular, they are open in $\mathscr Y\an$. 
But if $P$ is a locally constructible
subset of $\mathscr Y$ such that $P\an$ is open in $\mathscr Y\an$, then $P$
is open in $\mathscr Y$
by \cite{berkovich1993}, \Cor 2.6.6
(which is stated for the case of a constructible subset, but 
extends to the locally constructible case
by arguing locally), hence we are
done. 
\end{proof}

We end this section by showing that if one starts from a proper map, 
one can get
some local constructibility and Zariski-openness results \emph{on the target.}

\begin{theo}\label{thm-locus-target}

Let~$X$ be a~$k$-analytic space
and let~$Y$ be a \emph{proper}
$X$-analytic space. 
Let~$\mathscr E$ be a coherent sheaf on~$Y$
and let~$n$ and~$d$ be two
non-negative integers. For every subset $\Pi$ of $Y$, 
we denote by $\brac \Pi$
the set of points $x$ of $X$ such that $Y_x\subset \Pi$. 
We use the notation of Theorem \ref{thm-constloc-main}. 

\begin{enumerate}[1]

\item The sets~$\brac {A_n}$, $\brac{A_\infty}$, $\brac{B_n}$
$\brac C$, $\brac D$, $\brac{E_n}$, $\brac{E_\infty}$, $\brac \Delta$
and $\brac \Theta$
are locally
constructible subsets of $X$ (hence are constructible if $X$ is finite-dimensional, see Proposition \ref{prop-cons-gloc}). 

\item The sets~$\brac{A'_n}$, $\brac {A'_\infty}$, $\brac{C'}$, $\brac {D'}$,
$\brac{E'_{\infty,d}}, \brac{E'_\infty}$ and $\brac{\uflat{\mathscr E} X}$
are Zariski-open
subsets of $X$.
The subset $\Omega$ of $X$ consisting of points
$x$ such that $\dim Y_x\leq d$ is Zariski-open. 

\item If~$\supp E$ is purely of 
relative dimension~$d$ over $X$, then the set $\brac{B'_n}$ is a Zariski-open 
subset of $X$. If $Y$
is purely of relative dimension $d$
over $X$, then $\brac{E'_n}, \brac{\Delta'}$ and $\brac{\Theta'}$ 
are Zariski-open subsets of $X$. 
\end{enumerate}
\end{theo}

\begin{proof}
Since $Y$ is proper over $X$, 
the image in $X$ of any Zariski-closed subset (\resp locally constructible subset)
of $Y$
is a Zariski-closed subset of $X$ (\resp a locally constructible
subset of $X$): this follows from \ref{ss-kiehl-proper} (\resp Theorem \ref{thm-chevalley-proper}). 
Hence the theorem 
is a straightforward consequence of Theorem \ref{thm-constloc-main}, 
Theorem \ref{thm-flatloc-zaropen}
as far as  $\brac{\uflat{\mathscr E} X}$
is concerned, and \Th 4.9 of \cite{ducros2007}
as far as $\Omega$ is concerned.
\end{proof}

\chapter{Algebraic properties: target, fibers and source}\label{c-TARG}

In commutative algebra and algebraic geometry,
many familiar properties satisfy a principle of the following kind: if one is given a flat map,
and if the property under investigation holds on the target and fiberwise, then it holds on the source. 
This has proved of fundamental importance, and the goal of this chapter is to get similar results 
in analytic geometry. 

Our general strategy is not the same as in \Chp \ref{LOC}
\footnote{We will nevertheless \emph{use} two results from \Chp \ref{LOC}, namely 
the Zariski-openness of the relative flat locus (Theorem \ref{thm-flatloc-zaropen})
and our ``deboundarization" statement (Lemma \ref{lem-deboundarize}).}.
In the latter, 
we established analogues of various algebraic results by writing direct analytic proofs,
without using these results. 
Here we shall start from a property satisfying some principle of
the kind described above in the algebraic setting
(we use again a rather abstract, axiomatic presentation as we did in 
\ref{s-category-framework}-\ref{s-appendix-analytic}: 
see Section \ref{s-fib-valid}
below), and deduce from this that it satisfies the analogous
principle in analytic geometry; this is Theorem \ref{thm-base-fibers},
see also its ``concrete" version in Theorem \ref{thm-base-fibers-concrete}. 

Section \ref{s-weak-fact}
is devoted to some preparatory work, and ends with a result
that will be crucial for the proof of Theorem \ref{thm-base-fibers},
but also has independent interest (Theorem \ref{thm-weak-chevalley}). Let
us quickly explain what it consists of. Let $\mathscr Y\to \mathscr X$
be a morphism of integral schemes of finite type over a field, and let $d$ be the 
generic dimension of the fibers. It is then obvious that $\mathscr Y\to \mathscr X$
admits a factorization 
$\mathscr Y\to \mathscr Z\hookrightarrow \mathscr X$, with $\mathscr Z$ an integral
closed subscheme of $\mathscr X$ of dimension $\dim \mathscr Y-d$, and $\mathscr Y\to \mathscr Z$
dominant: simply take for $\mathscr Z$ the reduced closure of the image of the generic point
of $\mathscr Y$ in $\mathscr X$. This is of course very useful for various purposes: induction on dimension 
of the ground scheme, reduction to the dominant situation in order to use genericity arguments, \etc

Unfortunately, there is no such theorem in analytic geometry,
but Theorem \ref{thm-weak-chevalley}
provides a weak (and nonetheless useful) substitute for it. Let~$\phi \colon Y\to X$
be a morphism between~$k$-analytic spaces
and let~$d$ and~$\delta$ be two non-negative
integers. 
Assume that~$Y$ is of dimension~$d$
(so $Y\neq \emptyset$)
and that~$\phi$ is purely
of relative dimension~$\delta$. 
Then Theorem \ref{thm-weak-chevalley} asserts that
there exists a non-empty affinoid domain~$V$ of~$Y$, an affinoid domain~$U$ of~$X$,
and a Zariski-closed subset~$S$ of~$U$ of pure dimension~$d-\delta$
such that~$\phi(V)\subset S$. One thus gets a nice factorization, but only after restriction to some affinoid domain of the source. 
We shall not say here anything about its proof, except that one proceeds in a perhaps unsual way: in some sense, one reduces to the \emph{maximally non-strict
case}.

\section{A weak ``dominant factorization theorem"}\label{s-weak-fact}

\begin{enonce}[remark]{Notation}
If~$X$ is a~$k$-analytic space and if~$x$
is a point of $X$, we
denote
by~$\rho_k(x)$ the dimension of the~$\Q$-vector space
$\abs{\hr x\gpm}^\Q/\abs{k\gpm}^\Q$. 
\end{enonce}
 
\begin{rema}\label{rem-dk-rhok}
Let $X$ be a $k$-analytic space of finite dimension
$d\geq 0$ and let $x$ be a point of $X$. Since $d_k(x)$
is the sum of $\rho_k(x)$ and of the transcendence degree
of $\hrt x^1$ over $\widetilde k^1$, we have
$\rho_k(x)\leq d_k(x)$, with equality if and only if 
$\hrt x^1$ is algebraic over $\widetilde k^1$.
As $d_k(x)\leq d$, we see that $\rho_k(x)=d$ if and only if
$d_k(x)=d$ and $\hrt x^1$ is algebraic over 
 $\widetilde k^1$.
\end{rema}

\begin{lemm}\label{lem-dim-ratrk}
Let~$X$ be a
$k$-analytic space of
finite dimension~$d\geq 0$.
Assume that
the~$\Q$-vector space~$\R_+\gpm/\abs{k\gpm}^\Q$
is infinite-dimensional.
There exists a point~$x\in X$
such that~$\rho_k(x)=d.$
\end{lemm}

\begin{proof}
The problem we are dealing with is insensitive to nilpotents,
hence we can assume that~$X$ is reduced. The space $X$ has a
$d$-dimensional 
affinoid domain, and the latter has a $d$-dimensional
irreducible component. We can thus can assume that $X$ is affinoid and
integral. Let us argue now by induction on $d$.

If~$d=0$ the space~$X$ consists of one
rigid point and we are done. Let us assume now that~$d>0$ and that
the lemma has been proven for smaller dimension. We distinguish two cases. 

Let us first assume that the Krull dimension of~$X\al$ is zero.
Since~$X$ is reduced, this means that $\mathscr O_X(X)$ is a
field. The generalized affinoid \emph{Nullstellensatz}
(\cite{ducros2007}, \Th 2.7) then tells that (at least) one of the following
assertions holds: 

\begin{enumerate}[a]
\item There exists a $k$-free polyradius $r$ such that $X=\mathscr M(k_r)$. 

\item The absolute value of $k$ is trivial, and there exist two real numbers
$r$ and $s$ with $0<r\leq s<1$ such that $X$ is the compact annulus
$\{x\in \A^{1,\mathrm{an}}_k, r\leq \abs{T(x)}\leq s\}$. 
\end{enumerate}

If (a) holds then $X$ consists of one point $x$ and $d=d_k(x)=\rho_k(x)$,
hence we are done. If (b) holds then $d=1$ and we may take $x=\eta_r$. 

Let us now assume that the Krull dimension of~$X\al$ is positive.
Under this assumption,~$X\al$ admits a non-empty proper
Zariski-closed subset.
As a consequence, there exists an analytic function~$f$ on~$X$ whose
zero-locus is non-empty, and not the whole space~$X$. The image~$\abs{f}(X)$ is then
a compact interval~$[0,R]$ with~$R>0$; by our assumption on the value group
$\abs{k\gpm}$, 
this interval contains an element~$r$ that
does not belong to~${\abs k\gpm}^\Q$. 

Let~$V$ be the affinoid domain of~$X$ defined by the equation~$\abs f=r$. By the choice
of~$r$, the domain~$V$ is non-empty (hence~$d$-dimensional), and
its~$k$-analytic structure factorizes through a~$k_r$-analytic structure
induced by~$f$. For any $x\in V$ one has~$d_k(x)=d_{k_r}(x)+1$,
which implies that the~$k_r$-analytic dimension of~$V$ is~$d-1$
Since $\R_+\gpm/\abs{k_r\gpm}^\Q$ is still infinite dimensional
(because $\abs{k_r\gpm}^\Q/\abs{k\gpm}^\Q=r^\Q$), we may apply
the induction hypothesis. It asserts that there exists 
a point~$x$
in~$V$ such that~$\rho_{k_r}(x)=d-1$, and one has
$\rho_k(x)=\rho_{k_r}(x)+1=d$.
\end{proof}

\begin{lemm}\label{lem-im-rigid}
Let $d$ be a non-negative integer, and let
$\phi \colon
Y\to X$ be a morphism of pure
relative dimension~$d$ between $k$-analytic spaces. 
Assume that $Y$ is quasi-compact and $d$-dimensional.
The image~$\phi(Y)$ consists of finitely
many rigid points.
\end{lemm}

\begin{proof}
Let~$V$ be a non-empty affinoid domain of~$Y$, and let~$V'$
be an irreducible component of~$V$ (say, with its reduced structure).
We shall prove that~$\phi(V')$ 
consists of a single rigid point, which will yield the required result
by quasi-compactness of~$Y$. 
Let~$r$ be a~$k$-free polyradius such
that~$\abs{k_r\gpm}\neq \{1\}$ and such that
$V'_r$ is strictly~$k_r$-affinoid; note that~$V'_r$ remains
irreducible (Lemma \ref{lem-xxr-scheme} applied with $\mathscr X=\spec A$, or \cite{ducros2007}, Lemma 1.3). 

Let~$V''$ be the union of all irreducible components of~$V$
that are not equal to~$V'$. The open subset~$V'_r\setminus V''_r$
of~$V'_ r$ is non-empty and strictly~$k_r$-analytic, hence has a~$k_r$-rigid
point~$y$ (\ref{ss-analytic-nullst}); let~$x$ be its image on~$X_r$. Since~$x$ is a rigid point, 
$\phi_{|V_r}^{-1}(x)$ is a Zariski-closed subset of~$V_r$, which is purely~$d$-dimensional
by assumption and contains~$y$; since~$y\notin V''_r$
and $\dim V'\leq d$, this forces~$\phi_{|V_r}^{-1}(x)$ to contain
the whole component~$V'_r$
(and the dimension of~$V'$ to be equal to~$d$: as a by-product of our proof, $Y$ is purely~$d$-dimensional). 
In other words, $\phi(V'_r)=\{x\}$.

Let~$t$ be the image of~$x$ on~$X$. One has~$\phi(V')=\{t\}$; it remains to show that~$t$ is rigid. 
For any pre-image~$z$ of~$t$ on~$X_r$, the fiber of~$V'_r$ over~$z$
is naturally isomorphic
to~$V'\times_{\hr t}{\hr z}$; in particular, it is non-empty, which shows that~$z=x$; in other words, 
the set of pre-images of~$t$ on~$X_r$ is the singleton~$\{x\}$. 
Since this set can be identified with~$\mathscr M(\hr t_r)$, 
we see that~$r$ is~$\hr t$-free and $\hr x=\hr t_r$
(\ref{ss-intro-etar}). As
$\hr x$
is a finite Banach $k_r$-algebra, $\hr t$ is a finite Banach~$k$-algebra
by \Prop 2.1.11
of~\cite{berkovich1990}, which means 
that~$t$ is rigid. 
\end{proof}

\begin{theo}\label{thm-weak-chevalley}
Let~$\phi \colon Y\to X$ be a morphism between~$k$-analytic spaces
and let~$d$ and~$\delta$ be two non-negative
integers. 
Assume that~$Y$ is of dimension~$d$
(so $Y\neq \emptyset$)
and~$\phi$ is purely
of relative dimension~$\delta$ (so $\delta\leq d$). 
There exists a non-empty affinoid domain~$V$ of~$Y$, an affinoid domain~$U$ of~$X$,
and a Zariski-closed subset~$S$ of~$U$ of pure dimension~$d-\delta$
such that~$\phi(V)\subset S$. (Note that in the case where $\delta=d$, this
is an immediate consequence of Lemma \ref{lem-im-rigid}). 
\end{theo}

\begin{rema}
Even if $\abs{k\gpm}\neq \{1\}$
and $Y$ and $X$ are strictly $k$-analytic, the affinoid
domain built by the proof of the theorem
is not strict in general: if the analytic field
$k$ is topologically of countable
type over its prime complete field
(\eg, $k=\Q_p, \C_p$ or $\F_p(\!(t)\!)$)
and if $d>\delta$, the domain $V$ is $k_r$ analytic
for some $k$-free
polyradius $r=(r_1,\ldots, r_{d-\delta})$, hence non-strict. 
\end{rema}

\begin{proof}[Proof of the theorem \ref{thm-weak-chevalley}]
Choose a $d$-dimensional irreducible component
$Y'$ of $Y$,  and then choose a non-empty affinoid domain $Y''$ inside the open subset of 
$Y'$ consisting 
of points that do not belong to another component; by construction, $Y''$ is purely $d$-dimensional. 
Let~$X'$
be an affinoid domain of~$X$ intersecting~$\phi(Y'')$. By replacing~$Y$ with a non-empty affinoid
domain of~$\phi^{-1}(X')\cap Y''$ and~$X$ with~$X''$, we can assume that~$X$ and~$Y$ are affinoid,
and that~$Y$ is of pure dimension~$d$.

The definition of~$\phi \colon Y \to X$
only involves countably many parameters; therefore
there exists a complete subfield~$k_0$ of~$k$ topologically
of countable type
over its prime complete field, and a morphism~$Y_0\to X_0$ between~$k_0$-affinoid spaces
such that~$Y\to X$ is deduced from~$Y_0\to X_0$ by ground field extension to~$k$. 
Since the assertion we are interested in is ``stable under
ground field extension'', we can 
replace $k$ with~$k_0$, and~$Y\to X$ with~$Y_0\to X_0$,
and hence assume
that~$k$ is topologically of countable type over its prime complete field.

Since $k$ is of countable type over its prime complete field, 
$\abs{k\gpm}$ is countable
and $\R_+\gpm/\abs{k\gpm}^\Q$
is thus infinite-dimensional. It follows then from 
Lemma \ref{lem-dim-ratrk}
that there exists~$y\in Y$ with~$\rho_k(y)=d$; set~$x=\phi(y)$. 
Since $Y$ is $d$-dimensional, Remark \ref{rem-dk-rhok}
ensures that $d_k(y)=d$ and that $\hrt y^1$ is
algebraic over $\widetilde k^1$. Since $\phi$ is purely of relative dimension $\delta$,
it follows from \ref{ss-dim-sobafi} (4) that $d_k(x)=d-\delta$. As $\hrt x^1\subset \hrt y^1$,
$\hrt x^1$
is also algebraic over $\widetilde k$; we thus have $\rho_k(x)=d_k(x)=d-\delta$, 
again by Remark \ref{rem-dk-rhok}.

Since~$X$ is affinoid, the group~$\abs{\hr x\gpm}$ is generated by elements
of the form~$\abs{f(x)}$ where~$f$ is an analytic function on~$X$ which does not
vanish at~$x$. 
Therefore there exist analytic functions~$f_1,\ldots, f_{d-\delta}$ on~$X$
which do not vanish at~$x$ and such that~$(\abs{f_i(x)})_{1\leq i\leq d-\delta}$
is a basis of $\abs{\hr x\gpm}^\Q/\abs{k\gpm}^\Q$; we set~$r_i=\abs{f_i(x)}$ for every~$i$, 
and~$r=(r_1,\ldots, r_{d-\delta})$. Note that~$r$ consists by construction of positive real numbers
which are multiplicatively $\Q$-linearly independent modulo $\abs{k\gpm}$; otherwise said, 
$r$ is $k$-free. 

Let~$U$ be the affinoid domain of~$X$ defined by the equations~$\{\abs{f_i}=r_i\}_{1\leq i\leq \delta}$
and let~$V$ be its pre-image in~$Y$; note that $y\in V$ by construction, 
so $V\neq \emptyset$.
As $V$ is a non-empty affinoid domain of~$Y$, it is~$d$-dimensional
and compact; by construction, $y\in V$. Now~$U$ and~$V$ inherit through 
the~$f_i$'s a compatible $k_r$-analytic structure. For every~$z\in V$ one has
$d_k(z)=d_{k_r}(z)+d-\delta$; therefore $V$ is of~$k_r$-analytic dimension~$d-(d-\delta)=\delta$. Since
the $k_r$-analytic map
$V\to U$ is purely of relative
dimension~$\delta$, Lemma~\ref{lem-im-rigid}
above ensures that~$S:=\phi(V)$ consists of finitely many $k_r$-rigid points of~$U$; 
in particular, $S$ is Zariski-closed in~$U$ and~$d_k(t)=d-\delta$ for every~$t\in S$;
therefore~$S$ is of pure~$k$-analytic dimension~$d-\delta$.
\end{proof}

\section{The axiomatic setting}\label{s-fib-valid}

\subsection{}\label{ss-assumptions-bfb}
We use the general categorical setting
and the notation introduced in~\ref{s-category-framework}, especially
$\mathfrak T, \mathfrak F$, and $\mathfrak L$ (Definition \ref{def-category-framework}, 
\ref{ss-fiber-category}, \ref{ss-mathfrak-l}), 
but we assume that $\mathfrak F=\mathfrak T$
or $\mathfrak F=\mathfrak{Coh}$ (Examples \ref{ex-fiber-trivial}
and \ref{ex-fiber-coh}; see also Convention \ref{ss-conf-mathfrakl}). 

We fix a property~$\mathsf P$ as in~\ref{s-alg-properties}. We assume
that $\mathsf P$ satisfies conditions \hreg, \field and~
\open~of \ref{ss-list-hreg}
(see Examples \ref{ex-hreg-rings} and \ref{ex-hreg-modules}
for properties of interest that satisfy those conditions). If
$X$ is an object of $\mathfrak T$ and if $D$ is an object of $\mathfrak F_X$, it 
makes sense to say that $D$ satisfies $\mathsf P$ at a given point of $X$, or
that $D$ satisfies $\mathsf P$
(see Remark \ref{rem-valid-localization}
and Lemma-Definition \ref{lem-equiv-valid}).
Let $\phi \colon Y\to X$ be a morphism between $k$-analytic
spaces or between schemes belonging to $\mathfrak T$. We shall say that an object $D$
of $\mathfrak F_Y$ satisfies $\mathsf P$ fiberwise at a point $y$ of $Y$ if $D_{Y_{\phi(y)}}$
satisfies $\mathsf P$ at $y$, and that it satisfies $\mathsf P$ fiberwise if satisfies $\mathsf P$ fiberwise
at every point of $Y$. 

\begin{lemm}\label{lem-basfib-local}
Let $Y\to X$ be a morphism of~$k$-affinoid spaces, 
let $y$ be a point of $Y$ and let 
$x$ be its image in $X$. Let 
$\mathfrak p$ be a prime ideal of
$\mathscr O_{X,x}$, and let $Z$ be an integral
closed analytic subspace of $X$ inducing
the quotient map $\mathscr O_{X,x}\to \mathscr O_{X,x}/\mathfrak p$. 
Let $D$ be an object of~$\mathfrak F_Y$. Assume that there exists a non-empty
Zariski open subset $Z'$ of $Z$ such that $D_{Y\times_XZ'}$ satisfies $\mathsf P$. Then 
$D_y$ satisfies $\mathsf P$ fiberwise over $\mathfrak p$ with respect to
$\spec \mathscr O_{Y,y}\to \spec \mathscr O_{X,x}$. 
\end{lemm}

\begin{proof}
Set $T=Y\times_XZ$.
By definition of $Z$, the local 
ring $\mathscr O_{Z,x}$ is
the domain $\mathscr O_{X,x}/\mathfrak p$, and $\mathscr O_{T,y}$
is equal to $\mathscr O_{Y,y}/\mathfrak p\mathscr O_{Y,y}$; hence the fiber
of $\spec \mathscr O_{Y,y}\to \spec \mathscr O_{X,x}$ over $\mathfrak p$
is equal to the generic
fiber of
$\spec \mathscr O_{T,y}\to \spec \mathscr O_{Z,x}$.
We denote by $\xi$ the generic point of $Z\al$; 
note that the generic point of $\spec \mathscr O_{Z,x}$ lies over $\xi$
(by flatness). 

By assumption, $D_{T\times_ZZ'}=D_{Y\times_XZ'}$ satisfies $\mathsf P$. 
Since $\mathsf P$ satisfies \hreg~and
the map
$T\to T\al$
is flat as a morphism of locally ringed spaces, $D\al_{(T\times_ZZ')\al}$
satisfies $\mathsf P$ as well; as $(T\times_ZZ')\al$
contains by definition of $Z'$ the generic fiber of $T\al\to Z\al$, the object
$D\al_{T\al}$ satisfies $\mathsf P$ at every point of $T\al$ lying above $\xi$. 
Since $\mathsf P$ satisfies \hreg
~and
the map $\spec \mathscr O_{T,y}\to \spec \mathscr O_{T\al,y\al}$ is regular,
it follows that $D_{T,y}$ satisfies 
$\mathsf P$ at every point of $\spec \mathscr O_{T,y}$ lying above $\xi$. In particular, 
$D_{T,y}$
satisfies $\mathsf P$ at every point of the generic fiber of 
$\spec \mathscr O_{T,y}\to \spec \mathscr O_{Z,x}$; 
but at such a point, the
validity of $\mathsf P$ is equivalent to \emph{fiberwise}
validity
with respect to $\spec \mathscr O_{T,y}\to \spec \mathscr O_{Z,x}$.
The generic
fiber of the map
$\spec \mathscr O_{T,y}\to \spec \mathscr O_{Z,x}$
is
equal to the fiber of $\spec \mathscr O_{Y,y}\to \spec \mathscr O_{X,x}$ over $\mathfrak p$, 
so
$D_y$ satisfies $\mathsf P$ fiberwise above $\mathfrak p$
with respect to $\spec \mathscr O_{Y,y}\to \spec \mathscr O_{X,x}$. 
\end{proof}

\begin{lemm}\label{lem-basfib-shilov}
Let $Y\to X$ be a morphism 
of $k$-analytic spaces and let $r$
be a $k$-free polyradius. Let us denote by 
$\mathfrak s$ both Shilov sections $X\to X_r$ and~$Y\to Y_r$
(\ref{ss-shilov-section}). 
Let $y$ be a point of $Y$, let $x$ be its image in $X$,
and let $y'$ be any pre-image of $y$ in $Y_r$ lying above
$\mathfrak s(x)$; \eg, $y'=\mathfrak s(y)$. Let $D$ be an object of $\mathfrak F_Y$.

\begin{enumerate}[1]
\item The following are equivalent:

\begin{enumerate}[j]
\item $D_{Y_r}$ satisfies~$\mathsf P$ at~$y'$. 

\item $D$ satisfies $\mathsf P$ at $y$. 
\end{enumerate}

\item The following are equivalent:

\begin{enumerate}[j]\setcounter{enumii}{2}
\item $D_{Y_r}$ satisfies $\mathsf P$ fiberwise at~$y'$.

\item $D$ satisfies $\mathsf P$ fiberwise at $y$. 
\end{enumerate}

\end{enumerate}
\end{lemm}

\begin{proof}
Since $k_r$ is analytically separable over $k$
(Example \ref{ex-etar-ansep}),  (i)$\iff$(ii)
follows from Proposition \ref{prop-behavior-extension} and the fact that $\mathsf P$ satisfies
\hreg.

Let us prove that (iii)$\iff$(iv). The fiber $(Y_r)_{\mathfrak s(x)}$ is equal 
to $Y_x\times_{\hr x}\hr {\mathfrak s(x)}$. By definition of the Shilov section, 
$\mathfrak s(x)$ is the point $\eta_{\hr x,r}$
of the fiber 
$\mathscr M(\hr x_r)$; this implies 
that $\hr {\mathfrak s(x)}$ is
an analytically separable extension of $\hr x$
(Example \ref{ex-etar-ansep}, applied
to
the field $\hr x$). 
Therefore it follows again from 
Proposition \ref{prop-behavior-extension}
and the fact that $\mathsf P$ satisfies
\hreg~ that
$D_{Y_x}$ satisfies $\mathsf P$ at $y$ if and only if $D_{(Y_r)_{\mathfrak s(x)}}$
satisfies $\mathsf P$ at $y'$. Otherwise said, (iii)$\iff$(iv). 
\end{proof}

\subsection{New technical conditions}\label{ss-hweak}
Let us 
introduce now some technical conditions
that make sense for $\mathsf P$. Each of these involves a
local morphism $A\to B$
between local noetherian rings, and ``fiberwise" will be understood
with respect to $\spec B\to \spec A$; the maximal ideal of $A$
will be denoted by $\mathfrak m_A$.

\subsubsection*{Condition \hwk, when $\mathfrak F=\mathfrak T$}
For every \emph{flat} morphism~$A\to B$  of $\mathfrak L$, 
 the following implications hold: 
\begin{itemize}[label=$\bullet$] 

\item If~$B$ satisfies~$\mathsf P$, then~$A$ satisfies~$\mathsf P$. 

\item If~$A$ satisfies~$\mathsf P$, then~$B$ satisfies~$\mathsf P$ if
it satisfies $\mathsf P$ fiberwise. 
\end{itemize}

\subsubsection*{Condition \hstr, when $\mathfrak F=\mathfrak T$}
For every \emph{flat} morphism~$A\to B$ of $\mathfrak L$, the following implications hold: 
\begin{itemize}[label=$\bullet$] 

\item If~$B$ satisfies~$\mathsf P$, then~$A$ satisfies~$\mathsf P$. 

\item If~$A$ satisfies~$\mathsf P$, then~$B$ satisfies~$\mathsf P$ if
it satisfies $\mathsf P$ fiberwise at the closed point of $\spec B$; \ie, if $B/\mathfrak m_AB$
satisfies $\mathsf P$. 
\end{itemize}

\subsubsection*{Condition \hwkk,  when $\mathfrak F=\mathfrak {Coh}$}
For every morphism~$A\to B$ of $\mathfrak L$, every finitely 
generated $A$-module $M$ and every
non-zero and
\emph{$A$-flat} finitely generated $B$-module $N$
the following implications hold: 
\begin{itemize}[label=$\bullet$] 

\item If~$N\otimes_A M$ satisfies~$\mathsf P$, then~$M$ satisfies~$\mathsf P$. 

\item If~$M$ satisfies~$\mathsf P$, then~$N\otimes_A M$ satisfies~$\mathsf P$ if
$N$ satisfies $\mathsf P$ fiberwise. 
\end{itemize}

\subsubsection*{Condition \hstrr, when $\mathfrak F=\mathfrak {Coh}$}
For every morphism~$A\to B$ of $\mathfrak L$, every finitely
generated $A$-module $M$ and every non-zero and
\emph{$A$-flat} finitely generated $B$-module $N$
the following implications hold:

\begin{itemize}
[label=$\bullet$]

\item If~$N\otimes_A M$ satisfies~$\mathsf P$, then~$M$ satisfies~$\mathsf P$. 

\item If~$M$ satisfies~$\mathsf P$, then~$N\otimes_A M$ satisfies~$\mathsf P$ if
$N$ satisfies $\mathsf P$ fiberwise at the closed point of $\spec B$; \ie, 
if the $B/\mathfrak m_AB$-module $N/\mathfrak m_A N$
satisfies $\mathsf P$. 
\end{itemize}

\begin{exem}\label{ex-t-strong}
We assume that $\mathfrak F=\mathfrak T$. 
The following properties satisfy \hstr: being regular (\cite{ega42}, \Prop 6.5.1);  being CI
(\cite{avramov1975}); 
being Gorenstein (\cite{matsumura1986}, \Th 23.4). 
The property of being $R_m$ for some specified $m$ satisfies \hwk  (\cite{ega42}, \Prop 6.5.3). 
\end{exem}

\begin{exem}\label{ex-tprime-strong}
We assume that $\mathfrak F=\mathfrak {Coh}$. The property of being CM satisfies \hstrr
~by \cite{ega42}, \Cor 6.3.5. 
The property of being $S_m$ for some specified $m$ satisfies \hwkk
by \cite{ega42}, \Prop 6.4.1. 
\end{exem}

\begin{rema}
All properties mentioned in Example \ref{ex-t-strong}
and\ref{ex-tprime-strong}
also satisfy \hreg, \field, and
\open; see Examples \ref{ex-hreg-rings} and \ref{ex-hreg-modules}
for precise references. 
\end{rema}

\section{The main theorem}

\begin{theo}\label{thm-base-fibers}
Let $\mathfrak F$ and $\mathsf P$ be as in
\ref{ss-assumptions-bfb}. 
Let~$Y\to X$ be a morphism of~$k$-analytic spaces.
Let~$y$ be a point of $Y$ and let~$x$ be its image in~$X$; in what follows, ``fiberwise"
will 
always be relative to the morphism~$Y\to X$. 

\begin{enumerate}[1]
\item Assume that~$\mathfrak F=\mathfrak T$ and that~$\mathsf P$
satisfies~\hwk. 

\begin{enumerate}[b]

\item If~$Y$ satisfies~$\mathsf P$ at~$y$
and is~$X$-flat at~$y$,
then~$X$ satisfies~$\mathsf P$ at~$x$. 

\item If~$Y$ satisfies $\mathsf P$ fiberwise everywhere
and is $X$-flat 
at $y$, and $X$ satisfies~$\mathsf P$
at~$x$, then~$Y$ satisfies~$\mathsf P$ at~$y$. 

\item If~$Y$ is~$X$-flat at~$y$ and satisfies $\mathsf P$ fiberwise at~$y$,
 the space $X$ satisfies~$\mathsf P$
at~$x$, 
and~$\mathsf P$
satisfies \hstr, then~$Y$ satisfies~$\mathsf P$ at~$y$. 
\end{enumerate}

\item Assume now that~$\mathfrak F=\mathfrak{Coh}$, 
let~$\mathscr E$ be a coherent sheaf on~$X$, and let~$\mathscr F$
be a coherent sheaf on~$Y$. Assume that~$\mathsf P$ satisfies \hwkk. 

\begin{enumerate}[b]
\item If~$\mathscr F\boxtimes \mathscr E$ satisfies~$\mathsf P$ at~$y$
and  is~$X$-flat at~$y$
then~$\mathscr E$ satisfies~$\mathsf P$ at~$x$.

\item If~$\mathscr F$ is~$X$-flat at~$y$ and satisfies~$\mathsf P$ fiberwise everywhere, 
and $\mathscr E$ satisfies~$\mathsf P$
at~$x$,
then~$\mathscr F\boxtimes \mathscr E$ satisfies~$\mathsf P$ at~$y$. 

\item If~$\mathscr F$ is~$X$-flat at~$y$ and satisfies~$\mathsf P$ fiberwise at~$y$,
the coherent sheaf $\mathscr E$ satisfies~$\mathsf P$
at~$x$,
and $\mathsf P$
satisfies \hstrr, then~$\mathscr F\boxtimes \mathscr E$ satisfies~$\mathsf P$ at~$y$. 
\end{enumerate}
\end{enumerate}

\end{theo}

\begin{proof}
We are going to prove~(2). The proof of~(1)
is \emph{mutatis mutandis}
exactly the same (one simply has
to replace
everywhere~$\mathscr E$
with~$X$ or $\mathscr O_X$, replace
$\mathscr F\boxtimes \mathscr E$ and~$\mathscr F$
with~$Y$ and $\mathscr O_Y$, and
properties~\hwk~and \hstr~with \hwkk~and \hstrr, respectively). 

We can assume that~$Y$ and~$X$ are~$k$-affinoid. 
Assertion (2a) then follows from the first of the two properties
in \hwk. We
now focus on assertions
(2b) and (2c). 

\subsubsection{First reductions}\label{sss-basefib-reduc}
We assume that one of the following conditions 
is satisfied: 

\begin{enumerate}[A]
\item  $\mathscr E$ satisfies~$\mathsf P$ at~$x$, the coherent sheaf~$\mathscr F$ is~$X$-flat at~$y$,
and $\mathscr F$ satisfies~$\mathsf P$ fiberwise; 

\item  $\mathscr E$ satisfies~$\mathsf P$ at~$x$, 
the coherent sheaf~$\mathscr F$ is~$X$-flat at~$y$
and satisfies~$\mathsf P$ fiberwise at~$y$, and~$\mathsf P$ satisfies \hstr. 

\end{enumerate}
Our goal is to show that~$\mathscr F\boxtimes \mathscr E$ satisfies~$\mathsf P$  at~$y$. 
Since $\uflat {\mathscr F}X$ is (Zariski)-open by Theorem \ref{thm-flatloc-zaropen}
and~$\mathsf P$ satisfies~\open~by assumption, we can shrink $Y$
so that 
(B) 
can be replaced with: 

\begin{itemize}
\item [$(\mathrm B')$]  $\mathscr E$ satisfies~$\mathsf P$ at~$x$, 
the coherent sheaf~$\mathscr F$ is~$X$-flat everywhere and 
satisfies~$\mathsf P$ fiberwise {\em at every point of~$Y_x$}, and~$\mathsf P$ satisfies \hstr.
\end{itemize}

\subsubsection{Beginning of the proof under assumption \textnormal{(A)}}
We assume that
(A) holds, and we want to prove that~$\mathscr F\boxtimes \mathscr E$ satisfies~$\mathsf P$
at~$y$; we argue by induction on~$\dim_x X$. 

Assume that $\dim_x X=0$. 
This means that~$x$
is an isolated rigid point of~$X$
(\ref{ss-x-dimzero}); note that~$\mathscr O_{X,x}$ is then artinian. 
As~$\mathsf P$ satisfies
\hwk~
and $\mathscr E$ satisfies~$\mathsf P$ at~$x$, it
suffices to prove that~$\mathscr F_y/\mathfrak m_x\mathscr F_y$ 
satisfies~$\mathsf P$. But since~$x$ is rigid,
$\mathscr F_y/\mathfrak m_x\mathscr F_y=\mathscr F_{Y_x,y}$,
which satisfies~$\mathsf P$ by assumption~(i); hence we are done. 

Assume now
that~$\dim_x X>0$, and that the assertion holds
for strictly smaller local dimensions. 
As~$\mathsf P$ satisfies \hwk~
and
$\mathscr E$ satisfies~$\mathsf P$ at~$x$, it suffices to prove that
$\mathscr F_{Y,y}$ satisfies~$\mathsf P$ fiberwise 
with respect to~$\spec \mathscr O_{Y,y}\to \spec \mathscr O_{X,x}$.

Let us fix a prime ideal~$\mathfrak p$ of~$\mathscr O_{X,x}$, and 
prove 
that~$\mathscr F_y$ satisfies~$\mathsf P$ fiberwise 
above~$\mathfrak p$ with respect to
$\spec \mathscr O_{Y,y}\to \spec \mathscr O_{X,x}$.
One can shrink~$X$ (and accordingly~$Y$) 
so that~$\mathfrak p$
is induced by an integral closed analytic subspace~$Z$ of~$X$
(\ref{ss-gaga-domain}). 
By Lemma \ref{lem-basfib-local}, it suffices
to exhibit a non-empty Zariski-open subset $Z'$ of $Z$
satisfying the following:

\begin{itemize}
\item[$(\star)$] The coherent sheaf $\mathscr F_{Y\times_X Z'}$ satisfies~$\mathsf P$.
\end{itemize}  

Since $Z$ is integral, Lemma \ref{lem-loci-abstract}
ensures that the subset $Z'$ of 
$Z$ consisting of points
at which $Z$ is normal and
$\mathscr O_Z$ satisfies
$\mathsf P$ is Zariski-open and non-empty. We shall prove that
$Z'$ satisfies ($\star$). We argue by contradiction, so
\emph{we assume from now on that $Z'$ does \emph{not} satisfy $(\star)$}. 

\subsubsection{}
Let $D$
be the subset of $Y\times_XZ$
consisting of points at which 
at which~$\mathscr F$ does not satisfy~$\mathsf P$. Since
$\mathsf P$ satisfies \open, $D$ 
is a Zariski-closed subset of~$Y\times_X Z$. Since we assume that $Z'$
does not satisfy $(\star)$, the intersection $D\cap (Y\times_X Z')$
is non-empty. 
In particular, 
there exists an irreducible component~$D_0$
of~$D$
whose intersection~$D_0'$
with~$Y\times_X Z'$  is non-empty, hence
is a dense Zariski-open subset of~$D_0$. 
Let $D''_0$ be the intersection of $D'_0$ with the set of points
of $D_0$, that do not lie on any other component of $D$; this is
a dense Zariski-open subset of $D_0$ which
is Zariski-open in $D$ as well. 
Let~$d$ be the dimension of~$D_0$,
and let~$\delta$ be the infimum of the
relative
dimension of the morphism~$D''_0\to Z'$
(seen as a function from $D''_0$ to $\N$). 
The set of points
of~$D''_0$
at which the relative
dimension 
of~$D''_0\to Z'$ is precisely
equal to~$\delta$ is a non-empty
Zariski-open subset of~$D''_0$
(\cite{ducros2007}, \Th 4.9). We can therefore find an affinoid domain~$Y_0$ of~$Y\times_X Z'$
enjoying the following properties: 

\begin{itemize}[label=$\bullet$]
\item $\emptyset \neq Y_0\cap D\subset D''_0$ (in particular, $Y_0\cap D$ is purely of
dimension~$d$); 

\item  $Y_0\cap D\to Z'$ is purely of relative dimension $\delta$. 
\end{itemize}
By Theorem \ref{thm-weak-chevalley}, 
there exists a non-empty affinoid domain
$D_1$ of~$Y_0\cap D$ (endowed, say, with its reduced structure), 
an affinoid domain~$V$ of~$Z'$, and a 
Zariski-closed subset~$S$ of~$V$ of pure dimension $d-\delta$ which 
contains the image of~$D_1$. 
By construction, 
$D_1\subset D\cap (Y_0\times_{Z'}V)$
and~$\dim D_1=d$. 
Moreover, we can shrink $D_1$ so that $D_1=T\cap D$ for some affinoid
domain~$T$ of~$Y_0\times_{Z'}V$ (Remark \ref{rem-gerr-grau}). By replacing
$V$ with any of its connected components that intersect 
the image of $D_1$ and $T$
by the pre-image of this component
in $Y_0\times_{Z'}V$, we
get a triple $(T,V,S)$ satisfying the following: 

\begin{itemize}[label=$\bullet$]
\item $V$ is a connected affinoid domain of $Z'$; 
\item $S$ is a purely $(d-\delta)$-dimensional Zariski-closed subset of $V$; 
\item $T$ is an affinoid domain of $Y\times_X V$ such that $T\cap D$ is non-empty and purely $d$-dimensional, 
$T\cap D\to V$ is purely of relative dimension $\delta$, and the image of $T\cap D$
in $V$ is contained in
$S$. 
\end{itemize}

Note that under these assumptions $S$ is non-empty (it contains the image of $T\cap D$), hence $V$ is non-empty.
And since $Z$ 
is irreducible, it is purely of dimension $\dim Z$, whence the equality $\dim V=\dim Z$.
Moreover, 
$V$ is contained $Z'$; the latter is normal, and $\mathscr O_Z$ satisfies $\mathsf P$ on $Z'$. As a
consequence $V$ is normal (hence integral, because it is non-empty and connected) and $\mathscr O_V$
satisfies $\mathsf P$. 

\subsubsection{Proof of the inequality $d-\delta<\dim Z$}\label{sss-d-delta}
The next step
consists in proving that $d-\delta<\dim Z$, which will be crucial
for applying the induction hypothesis.  Note that
since $Z$ contains the non-empty $(d-\delta)$-dimensional
Zariski-closed subset $S$ of $V$, we have $d-\delta\leq \dim Z$. 
We are going to
show by contradiction that the strict inequality holds; so we suppose that $d-\delta=\dim Z$ $(=\dim V)$. 
Choose $t$ in $T\cap D$ such that~$d_k(t)=d$, and let~$z$ be the image of~$t$ in~$V$. 
We then have $d_k(z)=d-\delta=\dim V$ by \ref{ss-dim-sobafi}. 

Let~$r$ be a~$k$-free polyradius such that~$\abs{k_r\gpm}\neq\{1\}$
and such that~$T_r$ and~$V_r$ are strictly~$k_r$-affinoid.
Let~$\mathfrak s$ denote both Shilov sections  $X\to X_r$
and $Y\to Y_r$, and set~$\mathscr F_r=\mathscr F_{Y_r}$. 
The field~$k_r$ 
is analytically separable over~$k$. This implies, 
as $V$ is normal as an affinoid domain of~$Z'$, that~$V_r$ is normal. This
also implies, since~$\mathsf P$ satisfies~\hreg, that~$D_r\cap T_r$ is 
precisely the set of points of~$T_r$ at which~$\mathscr F_r$ does not
satisfy~$\mathsf P$. Moreover, it follows from~\ref{ss-dk-shilov}
that
\[d_{k_r}(\mathfrak s(z))=d_k(z)=\dim_k V =\dim_{k_r} V_r,\]
and from~Lemma \ref{lem-basfib-shilov}
that~$\mathscr F_r$ satisfies~$\mathsf P$ fiberwise at every point
lying over~$\mathfrak s(z)$
(because we are currently working under assumption (A); \ie, $\mathscr F$ satisfies $\mathsf P$ fiberwise). 

Since $V_r$ and $T_r$ are strict
over the non-trivially valued field~$k_r$, the fiber of~$T_r\cap D_r$ above~$\mathfrak s(z)$
(which is non-empty because it contains~$\mathfrak s(t)$) has an~$\hr{\mathfrak s(z)}$-rigid point~$s$. 
By Lemma~\ref{lem-deboundarize}, there exists a
quasi-\'etale morphism
$V'\to V_r$, and a pre-image $s'$ of
$s$ in 
$T':=V'\times_{V_r}T_r$
such that~$T'\to V'$ is inner at~$s'$;
we can shrink~$V'$ so that it 
is~$k_r$-affinoid. Let~$z'$ be the image of~$s'$ on~$V'$. 
Let us write down some consequences of the
quasi-\'etaleness of~$V'\to V_r$. 

\begin{enumerate}[i]
\item The space~$V'$ is normal by Proposition \ref{prop-tranferqsm-concrete}, 
and
\[d_{k_r}(z')=d_{k_r}(\mathfrak s(z))
=\dim V_r=\dim V'; \]
this implies that~$\mathscr O_{V',z'}$ is artinian (Corollary \ref{cor-interp-centdim}),
hence a field by normality (reducedness
would be sufficient). Since~$T'\to V'$ is inner at~$s'$, Theorem \ref{thm-localring-generic}
then
ensures that~$\mathscr O_{T'_{z'}, s'}$ is flat over~$\mathscr O_{T',s'}$. 

\item The mophism~$\spec \mathscr O_{T'_{z'}, s'}\to \spec \mathscr O_{T_z, s}$
is regular (Theorem \ref{thm-qsm-schemereg}). 

\item
The morphism~$\spec \mathscr O_{T',s'}\to \spec \mathscr O_{T,s}$ is
flat (Corollary \ref{cor-qsm-flat}; it is in fact
even regular
by Theorem \ref{thm-qsm-schemereg}). 
\end{enumerate}
The property $\mathsf P$
satisfies \hreg, and the coherent sheaf~$\mathscr F_r$ satisfies
$\mathsf P$ fiberwise at~$s$. This implies the following, by using
successively (ii), (i), and (iii): 

\begin{itemize}[label=$\bullet$]

\item $\mathscr F_{T'}$ 
satisfies~$\mathsf P$ fiberwise at~$s'$; 

\item $\mathscr F_{T'}$ 
satisfies~$\mathsf P$ at~$s'$; 

\item $\mathscr F_r$
satisfies~$\mathsf P$ at $s$. 

\end{itemize}
But $s$ lies on $D_r$, which is the
set of points of~$T_r$ at which~$\mathscr F_r$ does not
satisfy~$\mathsf P$, a contradiction. Hence our assumption 
that $\dim Z=d-\delta$ was wrong,
so $d-\delta<\dim Z$, as announced. 
\subsubsection{End of the proof 
under assumption \textnormal{(A)}}
We choose
a point $\tau$ in the non-empty space $T\cap D$, 
and we denote by $\sigma$ the image of~$\tau$ in~$S$. Since $\mathscr O_V$
satisfies $\mathsf P$,
the module $\mathscr O_{V,\sigma}$ satisfies~$\mathsf P$. We shall now 
prove that for every
prime ideal~$\mathfrak q$ of~$\mathscr O_{V,\sigma}$, the module~$\mathscr F_{T,\tau}$ satisfies~$\mathsf P$
fiberwise above~$\mathfrak q$; 
since~$\mathsf P$ satisfies \hwk, this will
imply that~$\mathscr F$ satisfies~$\mathsf P$ at~$\tau$, hence contradict the
fact that~$\tau$
lies on~$D$ -- and end the proof when (A) holds.

We fix $\mathfrak q$ as above, and we shrink~$V$ so that
the quotient map $\mathscr O_{V,\sigma}\to \mathscr O_{V,\sigma}/\mathfrak q$
is induced by an integral closed analytic subspace~$W$ of~$V$
(\ref{ss-gaga-domain}).
By Lemma \ref{lem-basfib-local}, 
it
suffices
to prove the existence of a non-empty Zariski-open subset~$W'$ of~$W$
such that~$\mathscr F_{T\times_V W'}$ satisfies~$\mathsf P$.

Let us first assume that $W=V$. 
By definition of $D$, the coherent
sheaf $\mathscr F_{T}$ satisfies~$\mathsf P$
at any point of~$T$ outside~$D$, and in particular
above the Zariski-open subset
$V\setminus S$. One has~$\dim S=d-\delta$, and 
$d-\delta<\dim Z=\dim V$
(\ref{sss-d-delta}). Therefore $S\subsetneq W$ and we may take
$W'=W\setminus S$. 

Assume now that $W\neq V$. Since $V$ is irreducible one has then
\[\dim W<\dim V=\dim Z\leq \dim_x X\]
(the latter inequality comes from the fact that $Z$ contains $x$ by definition).  
Since~$W\al$ is integral, 
it follows from Lemma  \ref{lem-loci-abstract}
that there exists a non-empty Zariski 
open subset~$W'$ of~$W$
such that the coherent sheaf~$\mathscr O_{W'}$
satisfies~$\mathsf P$ everywhere. Now~$\mathscr F_{T\times_V W'}$
satisfies $\mathsf P$ fiberwise (because so does $\mathscr F$)
and~$\dim W'=\dim W<\dim_x X$.
By the induction hypothesis, $\mathscr F_{T\times_V W'}$ satisfies~$\mathsf P$,
and we are done. 

\subsubsection{Proof under assumption $(\mathrm B')$}\label{sss-proof-bprime}
We assume that $(\mathrm B')$ holds, and
we shall prove that~$\mathscr F\boxtimes \mathscr E$ satisfies~$\mathsf P$
at every point of $Y_x$ (note that
the original point $y$ does not play any role in $(\mathrm B')$). 

Let~$r$ be a~$k$-free polyradius such that~$\abs{k\gpm}\neq\{1\}$ and
such that~$Y_r$ and~$X_r$ are strictly~$k_r$-affinoid; let~$\mathfrak s$
be the Shilov section of~$X_r\to X$. 
Due to Lemma \ref{lem-basfib-shilov},
we can replace~$k$ with~$k_r$, the spaces $Y$ and $X$ with~$Y_r$ and $X_r$ respectively, the point~$x$ with~$\mathfrak s(x)$,
and the sheaves~$\mathscr F$ and~$\mathscr E$ with $\mathscr F_{Y_r}$ and~$\mathscr E_{X_r}$
respectively,
to reduce
to
the case
where~$\abs{k\gpm}\neq\{1\}$ and $Y$ and~$X$ are strict. 

The set of points of~$Y$ at which~$\mathscr F\boxtimes \mathscr E$ satisfies~$\mathsf P$ is a Zariski-open subset of~$Y$
because~$\mathsf P$ satisfies~$\mathsf O$. It is therefore sufficient to prove that
the coherent sheaf~$\mathscr F\boxtimes \mathscr E$ satisfies~$\mathsf P$
at every \emph{rigid}
point of the strictly~$\hr x$-analytic space~$Y_x$. So, let $\omega$ be any~$\hr x$-rigid
point of $Y_x$. 
Lemma~\ref{lem-deboundarize}
ensures that there exists a strictly~$k$-analytic space~$X'$, a quasi-\'etale
morphism~$X'\to X$ and a point~$\omega'$ on~$Y':=Y\times_X X'$ lying above~$y$
and such that~$Y'\to X'$ is inner at~$\omega'$; we can  
assume that~$X'$
is affinoid. Since~$\mathsf P$ satisfies~\hreg, 
it follows from from Proposition~\ref{prop-transferqsm-general}
that we can, by pulling-back all 
the data to~$X'$, 
reduce to the case where~$Y\to X$ is inner at~$\omega$. 

Since~$\mathsf P$ satisfies~\hstrr~and
$\mathscr E$ satisfies~$\mathsf P$ at~$x$, it
suffices to prove that
the quotient $\mathscr F_{Y,\omega}/\mathfrak m_x\mathscr F_{Y,\omega}$ satisfies~$\mathsf P$. 
We can shrink~$X$ so that $\mathfrak m_x$
is induced by a closed analytic subspace~$Z$ of~$X$. One has~$\mathscr O_{Z,x}=\mathscr O_{X,x}/\mathfrak m_x$; as a consequence,
$\mathscr O_{Z,x}$ is a field. One also has~$\mathscr O_{Y,\omega}/\mathfrak m_x\mathscr O_{Y,\omega}=\mathscr O_{Y\times_XZ, \omega}$. The morphism~$Y\times_XZ \to Z$
is boundaryless at~$\omega$, and~$\mathscr O_{Z,x}$ is a field. It follows from Theorem \ref{thm-localring-generic}
that~$\mathscr O_{Y_x,\omega}$ is flat over~$\mathscr O_{Y\times_XZ,\omega}$. 
The module~$\mathscr F_{Y_x, \omega}$ satisfies~$\mathsf P$ by assumption~$(\mathrm B')$, and~$\mathsf P$ satisfies~\hstr. Therefore
$\mathscr F_{Y\times_X Z,\omega}=\mathscr F_{Y,\omega}/\mathfrak m_x\mathscr F_{Y,\omega}$ satisfies~$\mathsf P$.
\end{proof}

\begin{rema}
In the proof of the theorem, we needed twice (in \ref{sss-d-delta}
and \ref{sss-proof-bprime}) to ``spread out" a property from the generic fiber, 
which we achieved by using Theorem \ref{thm-localring-generic}.
But the latter only holds at an \emph{inner}
point. Therefore we first needed to reduce to the inner case, by applying
Lemma \ref{lem-deboundarize}, which requires to start from a \emph{rigid}
point of the relevant
fiber. This need of enough rigid points was the reason 
why we had twice to reduce to the strict case (by extending
scalars to some $k_r$).

\end{rema}
Theorem \ref{thm-base-fibers}
above deals with a general property satisfying some axioms, because we wanted to give
a unified proof, 
and to emphasize the assumptions we actually need. But of course, the properties of interest are those
mentioned in Examples \ref{ex-t-strong}
and \ref{ex-tprime-strong}.  For that reason, we are
now
going to write 
(a particular case of) this theorem 
with \emph {explicit}
properties involved. Note that assertion (1a) below 
is part of Lemma \ref{lem-flat-desc2}; assertion (3)
is a simple application of (1) and (2).

\begin{theo}[Concrete version of Theorem \ref{thm-base-fibers}]\label{thm-base-fibers-concrete}
Let~$Y\to X$ be a morphism of~$k$-analytic spaces. Let~$\mathscr E$ be a coherent
sheaf on~$X$, and let~$\mathscr F$ be a coherent sheaf on~$Y$. 
Let~$y$ be a point of $Y$ and let~$x$ be its image in~$X$; in what follows, ``fiberwise" will 
always be relative to the morphism~$Y\to X$. Let~$m$ be a non-negative integer. 

\begin{enumerate}[1]
\item Properties of the ambient spaces. 

\begin{enumerate}[b]

\item If~$Y$ is regular (\resp $R_m$, \resp CI, \resp Gorenstein) at~$y$
and $Y$ is~$X$-flat at~$y$, 
then~$X$ is regular (\resp $R_m$, \resp CI,  \resp Gorenstein) at~$x$.

\item If~$Y$ is~$X$-flat at~$y$ and is fiberwise~$R_m$ everywhere, and $X$ is~$R_m$
at~$x$,
then~$Y$ is~$R_m$ at~$y$. 

\item If~$Y$ is~$X$-flat at~$y$ and is fiberwise regular (\resp Gorenstein, \resp CI)
at $y$, and $X$ is regular (\resp Gorenstein, \resp CI)
at~$x$,  then~$Y$ is regular (\resp Gorenstein, \resp CI)
at~$y$. 

\end{enumerate}

\item Properties of coherent sheaves. 

\begin{enumerate}[b]

\item If~$\mathscr F\boxtimes \mathscr E$ is CM (\resp $S_m$)
at~$y$ and $\mathscr F$ is~$X$-flat at~$y$, then
$\mathscr E$ is CM (\resp $S_m$) at~$x$. 

\item If~$\mathscr F$ is~$X$-flat at~$y$ and is fiberwise $S_m$
everywhere, and $\mathscr E$ is~$S_m$ at $x$, then 
$\mathscr F\boxtimes \mathscr E$ is~$S_m$ at~$y$. 

\item If~$\mathscr F$ is~$X$-flat at~$y$ and fiberwise CM at $y$, and $\mathscr E$ is CM
at~$x$, then $\mathscr F\boxtimes \mathscr E$ is CM at~$y$. 

\end{enumerate}

\item Reducedness and normality. 

\begin{enumerate}[b]

\item If~$Y$ is reduced (\resp normal) at~$y$
and $Y$ is~$X$-flat at~$y$, 
then~$X$ is reduced (\resp normal) at~$x$.

\item If~$Y$ is~$X$-flat at~$y$ and is fiberwise reduced (\resp fiberwise normal)
everywhere, and $X$ is~reduced (\resp normal)
at~$x$,
then~$Y$ is reduced (\resp normal)
at $y$. 
\end{enumerate}

\end{enumerate}
\end{theo}
%Passage à l'appendice
\appendix

\chapter{Graded commutative algebra}\label{c-graded}

The purpose of this appendix is 
to introduce {\em graded} commutative algebra, after Temkin
\cite{temkin2004}. Most classical notions of classical commutative algebra have graded counterparts, and the usual theorems often remain {\em mutatis mutandis} true in the graded context,
with similar proofs; one only has essentially to add the word ``graded" or ``homogeneous" at suitable places. We shall therefore
give almost no proofs. The justifications are left to the reader, who can also fruitfully read \cite{temkin2004}.

\section{Basic definitions}\label{s-grad-alg}\index{graded!ring}\index{ring!graded}

\begin{defi}\label{def-graded-ring}
In this memoir, a {\em graded ring} will always be 
an {\em~$\R\gpm_+$-graded ring}; \ie, a ring $A$ equipped with a
decomposition $A=\bigoplus_{r\in \R_+\gpm} A^r$ as an abelian group, satisfying the condition
$A^r\cdot A^s\subset A^{rs}$ 
for all $r,s$ (the notation relative to the graduation is then {\em multiplicative}).
\end{defi}

\subsection{}Let $A$
be a graded ring. It follows from the definition that $1\in A^1$, and that $A^1$ is a usual ring. 
For any $r>0$, the summand $A^r$ is called the set of {\em homogeneous elements
of degree $r$}; any homogeneous nonzero element has thus a well-defined degree, but $0$ is homogeneous of degree $r$ for all positive $r$. 

Note that any ring can be considered as a trivially graded ring; \ie,
a graded ring in which any element is homogeneous of degree~$1$.

\subsection{}\label{ss-grad-agamma}
If~$A$ is a graded ring and if~$\Gamma$ is a subgroup of~$\R\gpm_+$, we shall denote by~$A^\Gamma$ the graded subring~$\bigoplus_{\gamma \in \Gamma} A^\gamma$ of~$A$. 

\subsection{}
A morphism of graded rings $f\colon A\to B$ is a usual ring homomorphism $f$ from $A$ to $B$ such that $f(A^r)\subset B^r$ for every $r>0$. 

\subsection{}An ideal $I$ of a graded ring $A$ is called
{\em homogeneous}
if $I=\bigoplus_r I_r$ (as an abelian group) or, what amounts to the same, 
if $I$ admits a set of generators consisting of homogeneous elements. 
If $I$ is such an ideal, the quotient $A/I$ inherits a graduation such that $A\to A/I$ is a
morphism of graded rings.

\subsection{}A {\em graded domain}
is a graded ring whose underlying commutative ring is a domain. A graded ring $A$ is a
graded domain if and only if $A\neq \{0\}$ and
the product of two {\em homogeneous}
nonzero elements of $A$ is always nonzero.

\subsection{}\label{sss-grad-field}\index{graded!field}\index{field!graded}
A {\em graded field} is a graded ring in which every {\em homogeneous}
nonzero element is invertible. If $K$ is a graded field, $K^1$ is a field
(whose characteristic will be also called the characteristic of $K$), and 
the set of degrees of homogeneous nonzero elements
of $K$ is a subgroup of $\R_+\gpm$; we shall denote it by $\mathfrak D(K)$.

The reader should be aware that the commutative ring underlying a graded field is {\em not}
a field in general. For example, let $K$ be a field and let $r$ be an element of $\R\gpm\setminus\{1\}$.
Let $L$ be the graded
ring whose underlying commutative ring is $K[T,T\inv]$ and whose graduation is such 
that $L^{r^i}=KT^i$ for every $i$ (and $L^s=\{0\}$ if $s\notin r^{\Z}$). Then any homogeneous nonzero
element of $L$ is invertible, hence 
$L$ is a graded field; but the ring
$L$ is not a field. 

\subsection{}
Let $I$ be a homogeneous
ideal of a graded ring $A$. We shall
say that $I$ is {\em prime}
if $A/I$ is a graded domain; this is the case if and only if $I\neq A$ and $ab\in I\Rightarrow a\in I\;{\rm or}\;b\in I$
for every pair $(a,b)$ of homogeneous elements of~$A$. 
We shall say that $I$ is {\em maximal}
if $A/I$ is a graded field; this is the case if and only if $I\neq A$ and $I$ is maximal among all proper
homogeneous ideals of $A$; note that every maximal homogeneous ideal is prime.
By Zorn's lemma any proper homogeneous ideal of $A$ is contained in a maximal 
homogeneous ideal; in particular if $A\neq \{0\}$ it contains a maximal homogeneous ideal.

\subsection{}
If $A$ is a graded ring and if $S$ is a multiplicative set of homogeneous elements of $A$, there is a well-defined
{\em graded
localization} $S\inv A$. For $r>0$, every element of $(S\inv A)^r$ can be written as a fraction $\frac a b$ with $a\in A^{rs}$ and $b\in S^s$ for some $s>0$; moreover two fractions $\frac a b$
and $\frac cd$ are equal if and only if there exists $e\in S$ such that $e(ad-bc)=0$. The localization $S\inv A$ comes
with a natural morphism of graded rings $A\to S\inv A$, which sends any homogeneous element $a$ to $\frac a 1$. 

For instance, any graded domain $A$ has a {\em graded field of fractions} $\mathrm{Frac}\;A$, given by the above construction 
for $S$ consisting of all non-zero homogeneous elements of $A$.

\subsection{}\label{sss-notation-ktr}
The above construction of a graded field whose underlying ring is not a field
(\ref{sss-grad-field})
can be generalized as follows. 
Let us start from any graded field $K$, let 
$r=(r_1, \ldots, r_n)$ be a family of positive real numbers and let $T=(T_1, \ldots, T_n)$
be a family of indeterminates. We denote by $K[T_1/r_1, \ldots, T_n/r_n]$, 
or by $K[T/r]$ for short, the graded domain whose underlying domain is $K[T]$, and whose 
graduation is such that $K[T/r]^s$ is equal to $\bigoplus_{I\in \N^n} K^{sr^{-I}}T^I$ for every positive $s$
(in other words, this is the
only graduation extending that of $K$ and such that every $T_i$ is homogeneous of degree $r_i$); if $r_i=1$ 
for some $i$ we shall write $T_i$ instead of $T_i/1$.  

We denote by $K(T/r)$ the graded field of fractions of the graded domain $K[T/r]$. 
If $r_1, \ldots,r_n$ 
are linearly independent as elements of the $\Q$-vector space $\R\gpm_+/\mathfrak D(K)^{\Q}$, the 
commutative ring underlying $K(T/r)$ is $K[T, T\inv]$. Thus if $n=1$, 
if $K$ is trivially graded and if $r\neq 1$ we recover the aforementioned
example.

\section{Graded linear algebra}\label{s-graded-linear}

We fix a graded ring $A$. 

\begin{defi}\label{def-grad-mod}\index{graded!module}\index{module!graded}
A {\em graded $A$-module}
is a (usual) $A$-module $M$ equipped with a decomposition $M=\bigoplus_{r>0}M^r$ as an abelian group, 
such that $A^rM^s\subset M^{rs}$ for all $r,s$.
\end{defi}

\subsection{}Let 
$M$ and $N$ be two graded $A$-module and let $f\colon N\to M$
be an $A$-linear map. For $r$ a positive
real number, the map $f\colon M\to N$
is called {\em homogeneous of degree $r$}
if $f(M^s)\subset N^{rs}$ for every $s>0$. The
map $f$ is called \emph{homogeneous}
(without a mention of the degree)
if it is homogeneous of degree $r$ for some $r>0$.

We denote by
$A\text{-}\mathsf{Mod}_{\rm g}$ the category whose
objects are graded $A$-modules and whose arrows are homogeneous
$A$-linear maps of degree $1$. 

\subsection{}
Let $M$, $N$ and $P$ be graded $A$-modules. 
An $A$-bilinear map $b\colon M\times N\to P$
is called {\em homogeneous}
if for every $r>0$, every $m\in M^r$ and every $n\in N^r$, the two maps
$b(m,\cdot)\colon N\to P$ and $b(\cdot,n)\colon M\to P$ are
homogeneous of degree $r$. 

\subsection{}
A submodule $N$ of a graded $A$-module $M$ is called a {\em graded}
submodule of $M$
if $N=\bigoplus_{r>0} N\cap M^r$ as an abelian group
or, what amounts to the same, if 
$N$ can be generated by a set of homogeneous elements. If $N$ is a graded submodule of 
$A$, it inherits a natural graduation with $N^r=N\cap M^r$ for every $r>0$, and the inclusion 
$N\hookrightarrow M$ is homogeneous of degree 1. The quotient $M/N$ also inherits a natural 
graduation for which
$(M/N)^r=M^r/N^r$ for every $r>0$ (as an abelian group); it makes $M/N$
a graded $A$-module, 
and the quotient map $M\to M/N$ is then homogeneous 
of degree $1$. 
Note that we may view $A$ as a graded $A$-module; its graded submodules are precisely its
homogeneous ideals. 

If $f\colon N\to M$ is a homogeneous $A$-linear map bewteen two graded $A$-modules, 
its image is a graded submodule of $M$, and its kernel is a graded submodule of $N$. 
Therefore $f$ is injective if and only if $(f(n)=0)\Rightarrow (n=0)$ for every 
{\em homogeneous}
element $n$ of $N$. 

If $M$ is a graded $A$-module, for every $r>0$ we define $M(r)$ as
the graded $A$-module whose underlying $A$-module is $M$ and whose
graduation is such that $M(r)^s=M^{rs}$ for every $s>0$. 

\subsection{}
If $(M_i)_{i\in I}$ is any family of graded $A$-modules, 
the usual direct sum $\bigoplus_i M_i$ inherits a graduation, 
with $(\bigoplus_i M_i)^r=\bigoplus M_i^r$ for all $r$; the
graded module $\bigoplus_i M_i$ is the disjoint sum of the
$M_i$'s in the category $A\text{-}\mathsf{Mod}_{\rm g}$.  
Combining this construction with that of quotients, we see that
$A\text{-}\mathsf{Mod}_{\rm g}$ admits arbitrary colimits.

\subsection{}
Let $M$ be a graded $A$-module
and let $(m_i)_{i\in I}$ be a family
of homogeneous elements of $M$, say $m_i\in M^{r_i}$ for every $i$. This family gives rise
to a homogeneous $A$-linear map of degree $1$ from $\bigoplus_i A(r_i^{-1})$
to $M$, that sends any homogeneous element $(a_i)$ to $\sum a_i m_i$. 
The family $(m_i)$ is said to be free
(\resp generating, \resp a basis) if this map is injective
(\resp surjective, \resp bijective). 

\subsection{}Let $K$ be a graded field. Graded $K$-modules will be rather called\index{graded!vector space}\index{vector space!graded}
graded $K$-vector spaces. If $E$ is such a space, it has a basis; moreover, 
all bases of $E$ have the same cardinality, which is called the {\em dimension}
of $E$. If $E'$ is any graded subspace of $E$, there exists a graded subspace $E''$ of $E$ such that
$E=E'\oplus E''$.

\begin{defi}
\label{def-grad-tens}\index{graded!tensor product}\index{tensor product!graded}
Let $A$ be a graded ring, and let $M$ and $N$ be two graded $A$-modules. 
The covariant functor from $A\text{-}\mathsf{Mod}_{\rm g}$ that sends
$P$ to the set of homogeneous $A$-bilinear maps from $M\times N$ to $P$
is representable. The object that represents it is called the {\em graded tensor product}
of $M$ and $N$ over $A$. 
\end{defi}

\subsection{}\label{ss-desc-gradotimes} If $M$ and $N$ are two graded $A$-modules, the $A$-module underlying their
graded tensor product over $A$ is the usual tensor product $M\otimes_A N$, and its summand
of homogeneous elements
of given degree $r$ is generated {\em as an abelian group}
by all elements of the form $m\otimes n$ with $m\in M^s$ and $n\in N^{r\inv s}$ for some $s$. 

Hence we shall also denote by $M\otimes_A N$ the graded tensor product of $M$ and $N$ over $A$. 

\subsection{}For every graded $A$-module $M$, the functor $M\otimes$ from $A\text{-}\mathsf{Mod}_{\rm g}$ to itself
commutes with colimits. 

\subsection{}\label{ss-grad-flat}
A graded $A$-module $M$ is called {\em flat}
if $M\otimes\colon A\text{-}\mathsf{Mod}_{\rm g}\to A\text{-}\mathsf{Mod}_{\rm g}$
preserves injections. Any graded vector space over
a graded field $K$ is flat over $K$. 

Let $M$ be a graded $A$-module. Since $M\otimes$ commutes with colimits, 
$M$ is flat if and only if $M\otimes$ preserves injections $N\hookrightarrow N'$
in $A\text{-}\mathsf{Mod}_{\rm g}$ such that $N'$ is generated by $N$ and a {\em finite}
family $(e_1,\ldots, e_m)$
of homogeneous elements. By induction on $m$, this
is the case if and only if $M\otimes$
preserves injections as above 
with $m=1$. Let $N\hookrightarrow N'$ be such an injection. The quotient $N'/N$ can then be generated by a single
homogeneous element, hence is isomorphic to $(A/I)(r)$ for some homogeneous ideal $I$ of $A$
and some positive real number~$r$. By the same piece of diagram-chasing as in the classical case, $M\otimes_A N\to M\otimes_A N'$ is
injective as soon as $M\otimes_A I\to M$ is so.
Using again commutation with colimits, we eventually see
that $M$ is flat if and only if $M\otimes_A J\to M$ is injective for every finitely
generated homogeneous ideal $J$ of $A$.

\section{Graded algebras and graded extensions}

\begin{defi}\label{def-grad-alg}
A {\em graded $A$-algebra}
is a graded ring $B$ endowed with a morphism of graded rings from $A\to B$; a morphism
of graded $A$-algebras is a morphism of graded rings commuting with the structure maps from $A$.
\end{defi}
 
\subsection{}Any graded $A$-algebra inherits a structure of a graded $A$-module. 
In particular if $B$ and $C$ are two graded $A$-algebras, the graded tensor product
$B\otimes_AC$ makes sense. Since $B$
and $C$ are usual $A$-algebras, $B\otimes_AC$ admits a natural
structure of a ring; the latter together with the graduation of $B\otimes_AC$
makes $B\otimes_AC$ a graded ring, which is the amalgamated sum of $B$ and $C$
along $A$ in the category of graded rings. 

\subsection{}\label{ss-graded-algebraic}
Let $K$ be a graded field and let $L$ be a graded $K$-algebra. If $L$ is nonzero then 
the structure morphism $K\to L$ is injective. We shall call $L$ a {\em graded extension}
of $K$ if $L$ is a graded field; let us assume from now on that this is the case.

Let $(x_i)_{i\in I}$ be a family of
homogeneous elements
of $L$, say  $x_i\in L^{r_i}$ for every $i$. Evaluating polynomials at $(x_i)$ yields a morphism
of graded $K$-algebras from $K[T/r]$ to $L$ (here $T=(T_i)_{i\in I}$, and we extend straightforwardly the definition 
of $K[T/r]$ given in~\ref{sss-notation-ktr}
to the case of an arbitrary family of indeterminates). The elements $x_i$ are said to be {\em algebraically
independent} 
over $K$ if this morphism is injective. A maximal family of
homogeneous elements of $L$ that are
algebraically independent over $K$
is called a transcendance basis of $L$ over $K$; there exist transcendance bases
of $L$ over $K$, all of which have the same cardinality; the latter is called
the {\em transcendance degree}
of $L$ over $K$. 

If $x$ is a homogeneous element  of $L$ of degree $r>0$,
the singleton family $\{x\}$ is algebraically independent over $K$ if $P(x)\neq 0$
for all nonzero homogeneous elements $P$ of $K[T/r]$ (where $T$ is now a single indeterminate).
If this is the case we say
that $x$ is \emph{transcendental} over $K$. The element $x$ is said to be
\emph{algebraic} over $K$ if it is not transcendental over $K$.
If $x$ is algebraic over $K$, the ideal $I$ generated by the homogeneous elements
$P\in K[T/r]$ such that $P(x)=0$ is principal (like any homogeneous ideal of $K[T/r]$); 
if $x\neq 0$, the \emph{minimal polynomial}
of $x$ is the unitary homogeneous generator of $I$ (the condition $x\neq 0$ ensures
that $r$ is well-defined; if $x=0$,  its minimal polynomial should certainly be $T$, but it could be
seen as belonging to $K[T/r]$ for
any $r$, without canonical choice).

A homogeneous element $x$ of $L$ is algebraic over $K$
if and only if the graded $K$-algebra generated by $x$ is finite-dimensional (and hence a graded
field). As a consequence, the homogeneous
elements of $L$ that are
algebraic over $K$ are the homogeneous elements of a graded subfield of $L$, which is called
the \emph{algebraic closure}
of $K$ inside $L$. 

A graded extension $L$ of $K$ is said to be algebraic if all its homogeneous
elements are algebraic over $K$ or, otherwise said, if its transcendance degree is zero. 

\subsection{}\label{ss-criterion-alg}\index{graded!extension}\index{extension!graded}
Let $K\hookrightarrow L$ be a graded extension of graded fields.
Let $r$ be a positive real number, let $x$ be a non-zero element of $L^r$ algebraic over $K$, 
and let $P\in K[T/r]$ be its minimal polynomial. If $n$ denotes the monomial degree of $P$, the
constant coefficient of $P$ (which is non-zero) has degree $r^n$, hence $r^n\in \mathfrak D(K)$. 
Morevoer every non-zero coefficient is of degree $r^i$ for some $i$; as a consequence, $x$ is algebraic 
over $K^{r^\Z}$.

Hence if $L$ is algebraic over $K$, then
$\mathfrak D(L)/\mathfrak D(K)$ is torsion and $L^\Gamma$ is 
algebraic over $F^\Gamma$ for every subgroup $\Gamma$ of $\R_+\gpm$; in particular, $L^1$ is algebraic
over $K^1$.
 Conversely, assume that 
$L^1$ is 
algebraic over $F^1$ and $\mathfrak D(L)/\mathfrak D(K)$ is torsion. Let $x$ be a 
homogeneous element of $L$. Since  $\mathfrak D(L)/\mathfrak D(K)$ is torsion, there exists
a non-zero homogeneous element $a$ of $K$ and $n\in \N$ such that $ax^n\in L^1$; hence $ax^n$
is algebraic over $K^1$, and $x$ is algebraic over $K$. 

\subsection{}\label{ss-graded-equalities}
Let $K\hookrightarrow L$ be a graded extension of graded fields. 
Let $(r_i)$ be a system of representatives of
the quotient $\mathfrak D(L)/\mathfrak D(K)$. For every $i$, 
let $x_i$ be a non-zero element of $L^{r_i}$; let $(y_j)$ be a basis of the $K^1$-vector space $L^1$. 
The family $(x_iy_j)_{i,j}$ is then a basis of the graded $K$-vector space $L$ (the verification 
is straightforward and left to the reader). In particular, 

\setcounter{equation}{0}
\begin{equation}
[L:K]=[L^1:K^1]\cdot[\mathfrak  D(L):\mathfrak D(K)]
\end{equation}
(this is an equality of cardinal numbers, possibly infinite). For every subgroup $\Gamma$
of $\R_+\gpm$ we thus have
\begin{equation}
[L^\Gamma:K^\Gamma]=[L^1:K^1]\cdot[\mathfrak  D(L^\Gamma):\mathfrak D(K^\Gamma)]\leq
[L^1:K^1]\cdot[\mathfrak  D(L):\mathfrak D(K)]=[L:K].
\end{equation}

Let $(s_\ell)$ be a family of positive real numbers lifting a basis of the $\Q$-vector space
$\mathfrak D(L)^\Q/\mathfrak D(K)^\Q$. For every $\ell$, let $z_\ell$ be a non-zero element of $L^{s_\ell}$; let
$(t_\lambda)$ be a transcendence basis of $L^1$ over $K^1$. The family $(t_\lambda)\coprod(z_\ell)$ is then 
a transcendance basis of $L$ over $K$.
Indeed, a straightforward computation (left to the reader) show that it consists of algebraically independent elements. Denote
by $F$ the graded subfield of $L$ generated by $K$ and the family $(t_\lambda)\coprod(z_\ell)$. By construction, 
$L^1$ is algebraic over $F^1$ and $\mathfrak D(L)/\mathfrak D(F)$ is torsion, whence our claim. In particular, 

\begin{equation}
\text{tr.\@~deg}(L/K)=\text{tr.\@~deg}(L/K)+\dim_{\Q}D(L)^\Q/\mathfrak D(K)^\Q
\end{equation}
(this is an equality of cardinal numbers, possibly infinite). 

\subsection{}\label{ss-graded-galois}
The classical theory of algebraic extensions admits a graded counterpart, for which
a general reference
is the first section of \cite{ducros2013b}. Let us simply recall here some basic definitions and properties.
The proofs
(most of which are outlined in \cite{ducros2013b})
essentially consist of a 
 ``graded transcription" of Bourbaki's approach; see \cite{Bourbaki1981}, \Chp V and
especially \S 7 about separable extensions, \S9 about normal extensions (they are called ``quasi-galoisiennes"
by Bourbaki), and \S 10 about Galois extensions. 

We fix an algebraic graded extension $K\hookrightarrow L$. 
We denote by $\mathrm{Gal}(L/K)$ the Galois group of $L$ over $K$; \ie, the
group of $K$-automorphisms
of the graded field $L$. If $[L:K]$ is finite, then $\mathrm{Gal}(L/K)$ 
is finite of cardinality $\leq [L:K]$; in general, 
it has a natural topology making it a profinite group (by considering the system
of all finite graded subextensions of $L$ stable under $\mathrm{Gal}(L/K)$). 

We say that the graded extension
$K\hookrightarrow L$ is \emph{radicial}
if for
every homogeneous element $a$ of $L$ there exist $n\geq 0$ such that $a^{p^n}\in L_{\mathrm{sep}}$, where
$p$ is the characteristic exponent of $L$; \ie, $p$ is equal to the characteristic of
$L$ if the latter is positive, and to $1$ otherwise. 
If $L$ is radicial over $K$, we have
$\mathrm{Gal}(L/K)=\{\mathrm{Id}_L\}$ (because $a\mapsto a^p$ is an endomorphism 
of the graded ring $L$, hence is injective since $L$ is a graded field).

Let $r$ be a positive real number. An element $x$ of $L^r$ is
called \emph{separable} over $K$
if $x=0$ or if its minimal polynomial $P$ satisfies the condition $P'(x)\neq 0$. 
There exists a graded subfield $L_{\mathrm{sep}}$ of $L$ such that
$L_{\mathrm{sep}}^r$ is for every $r>0$ the set of separable elements of $L^r$.
The graded subfield
$L_{\mathrm{sep}}$
is called the \emph{separable closure}
of $K$ inside $L$, and $L$ is called a separable extension of $K$ if $L_{\mathrm{sep}}=L$
(this is automatically
the case if \car $L=0$). The graded extension $L_{\mathrm{sep}}\hookrightarrow L$ is \emph{radicial}. 

The graded extension $K\hookrightarrow L$ is called \emph{normal}
if for every $r>0$ and every element
$x$ of $L^r\setminus \{0\}$, 
the minimal polynomial $P$ of $x$ splits in $L$; \ie, $P$ can be written $\prod_{1\leq i\leq n} (T-x_i)$ 
with $x_i\in L^r$ for all $i$.

The graded extension $K\hookrightarrow L$
is called \emph{Galois}
if it is both separable and normal. This is the case if and only if
the fixed graded field of $\mathrm{Gal}(L/K)$ is equal to $K$. 
If moreover $[L:K]$ is finite, $L$ is Galois over $K$ if
and only if the cardinality of $\mathrm{Gal}(L/K)$ is equal to $[L:K]$. 

The graded extension $K\hookrightarrow L$
is normal if and only if the 
fixed graded field of $\mathrm{Gal}(L/K)$ is radicial over $K$. 

In fact, for the graded extension $K\hookrightarrow L$ to be normal (\resp Galois), it suffices
that there exists a subgroup $G$ of $\mathrm{Gal}(L/K)$
whose fixed graded field is radicial over $K$ (\resp is equal to $K$); and if
this is the case, $G$ is dense in $\mathrm{Gal}(L/K)$.  

Let $\Gamma$ be a subgroup of $\R_+\gpm$. The graded subfield $L^\Gamma$
of $K^\Gamma$ is stable under $\mathrm{Gal}(L/K)$, whence
we get a restriction 
morphism $\mathrm{Gal}(L/K)\to \mathrm{Gal}(L^\Gamma/K^\Gamma)$. Let $H$
denote the image of this map (this is a closed subgroup of $\mathrm{Gal}(L^\Gamma/K^\Gamma)$),
and let $F$ denote the fixed graded field of $\mathrm{Gal}(L/K)$. 
The fixed graded field of $H$ is equal to $F\cap K^\Gamma=F^\Gamma$. Since $F^\Gamma$
is radical over $K^\Gamma$ (\resp equal to $K^\Gamma$) if $F$ is radicial over $K$ (\resp equal to $K$), 
we see that if $L$ is normal (\resp Galois) over $K$
then $L^\Gamma$ is normal
(\resp Galois) over $K^\Gamma$ and 
$\mathrm{Gal}(L^\Gamma/K^\Gamma)=G$.

\section{Graded valuations}\label{s-grad-val}\index{graded!valuation}\index{valuation!graded}

\begin{defi}\label{def-grad-val}
Let~$K$ be a graded field and let~$\Gamma$ be an ordered
abelian group with multiplicative notation. A $\Gamma$-valued \emph{graded valuation} on~$K$ is a map~$\abs {\cdot}$ defined on the set of homogeneous elements of~$K$ with values in~$\Gamma\cup\{0\}$ which satisfies the following conditions: 
\begin{enumerate}
\item $\abs 1=1, \abs 0 =0$, and~$\abs{ab}=\abs a \cdot \abs b$ for every pair~$(a,b)$ of homogeneous elements
of $K$;
\item for every pair~$(a,b)$ of homogeneous elements of $K$ {\em of the same degree} we have~$\abs {a+b}\leq \max(\abs a, \abs b)$.

\end{enumerate}
\end{defi}

\subsection{} If we do not need to focus on the group~$\Gamma$, or if the latter is clear from the context, we shall simply talk about a {\em graded valuation} on~$K$; if~$K$ is a field (viewed as a trivially graded field), a graded valuation on~$K$ is nothing but a classical Krull valuation on~$K$. Two graded valuations
$$\abs{\cdot}\colon \coprod K^r\to \Gamma_0\;{\rm and}\;\abs{\cdot}'\colon \coprod K^r\to \Gamma'_0$$ are called {\em equivalent}
if there exist an ordered abelian group $\Gamma''$, two increasing embeddings $i\colon \Gamma''\hookrightarrow \Gamma$
and $j\colon \Gamma''\hookrightarrow \Gamma'$, and a $\Gamma''$-valuation $\abs{\cdot}''$ on $K$ such that $\abs{\cdot}=i\circ \abs{\cdot}''$
and $\abs{\cdot}'=j\circ \abs{\cdot}''$.

\subsection{}\label{ss-valgrad-ring} For any graded  valuation $\abs{\cdot} $ on~$K$, the set~$\{\abs a\}_{a \in \bigcup_r K^r\setminus \{0\}}$ is an ordered group which is called the {\em value group}
of $\abs{\cdot}$. And the set~$$\bigoplus_r \{\lambda  \in K^r,\;\abs \lambda \leq 1\}$$ is a graded subring of~$K$ which is called the {\em graded ring of~$\abs{\cdot}$.} It is a local graded ring; \ie, it has a unique
maximal homogeneous ideal, namely~$\bigoplus_r\{\lambda\in K^r,\;\abs \lambda <1\}$.
The residue graded field of this local graded ring is called the {\em residue graded field} of~$\abs{\cdot}$.
Two graded valuations on~$K$ are equivalent if and only if they have the same graded ring. 

\subsection{} A graded subring~$A$ of~$K$ is a graded valuation ring of $K$ (\ie, the graded ring of 
some graded valuation on $K$) if and only if for every non-zero homogeneous element~$\lambda$ of~$K$, one has~$\lambda\in A$
or~$\lambda\inv \in A$; or, what amounts to the same, if and only if~$A$ is a graded local subring of~$K$ which is maximal for the domination relation
(the latter is defined similarly as for usual local rings). 
Hence by Zorn's lemma, every graded local subring of~$K$ is dominated by a graded valuation ring of $K$. It follows that any graded valuation 
on a graded subfield of $K$ extends to $K$ (with possibly larger value group). 

\subsection{}Let~$\abs{\cdot}$ be a graded valuation on~$K$, let~$A$ be its graded ring and let~$k$ be its graded residue field. If~$ \abs{\cdot}'$ is a graded valuation on~$k$, the pre-image of the graded ring of~$\abs {\cdot}'$ inside~$A$ is a graded valuation ring of~$K$. The corresponding graded valuation is called the  {\em composition} of~$\abs{\cdot}$ and~$\abs{\cdot}'$ and has the same residue graded field as~$\abs {\cdot}'$.

\subsection{}\label{ss-notation-gauss} Let~$S=(S_i)_{1\leq i\leq n}$
be a family of indeterminates
and let~$s=(s_i)_{1\leq i\leq n}$ be a family positive real numbers.
If~$\abs{\cdot}$ is any graded valuation on~$K$, we shall
denote by~$\abs{\cdot}_{\mathrm{Gau\ss}}$ the valuation on~$K(S/s)$ that sends any homogeneous element~$\sum_{I\in \N^n} a_I S^I$ to~$\max_I \abs a_I$. It has the
same value group as $\abs{\cdot}$ and is characterized by the following properties:

\begin{itemize}
\item$\abs{\cdot}_{\mathrm{Gau\ss}}$ is an extension of $K$ to $K(S/s)$ such that $\abs {S_i}_{\mathrm{Gau\ss}}=1$ for every $i$. 

\item the images of the $S_i$'s in the residue graded field of $\abs{\cdot}_{\mathrm{Gau\ss}}$ are algebraically independent over the residue graded field of $\abs{\cdot}$. 

\end{itemize}

\subsection{Graded reduction}\label{ss-gradred-general}
Let $A$ be a ring equipped with a sub-multiplicative semi-norm $\norm {\cdot}$; \ie, $\norm{\cdot}$ is a map
from $A$ to $\R_+$
such that $\norm 0=0$, $\norm 1\leq 1$\index{graded!reduction}\index{reduction!graded}
and 
$$\norm {-a}=\norm a, \;\norm{a+b}\leq \max(\norm a, \norm b), \;\;\text{and}\;\;\norm{ab}\leq \norm a \cdot \norm b$$
for every $(a,b)\in A^2$; one then has $\norm 1=1$ unless $\norm \cdot=0$. We shall denote by $\widetilde A$
the residue {\em graded} ring of $A$
in the sense of Temkin \cite{temkin2004}, which is by definition equal to
$$\bigoplus_{r>0} \{x\in A, \norm x\leq r\}/\{x\in A, \norm x<r\}.$$ Note that $\widetilde A^1$ is the usual residue ring of $A$. If~$a$ is any element of~$A$ and if~$r$ is a positive real number such that~$\norm a\leq r$, we
shall denote by~$\widetilde a^r$ the image of~$a$ in~$\widetilde A^r$. If~$\norm a \neq 0$ we shall
write~$\widetilde a$ instead of~$\widetilde a^{\norm a}$; if~$\norm a=0$ we set~$\widetilde a=0$.

\subsection{}\label{ss-gradred-field}
Let $k$ be a field equipped with a valuation $\abs {\cdot}\colon k \to \R_+$. The previous construction provides a graded residue ring
$\widetilde k$, which is easily seen to be a graded field. The field $\widetilde k^1$ is the residue field of the valuation $\abs{\cdot}$ in the classical sense, and the group of 
degrees $\mathfrak D(\widetilde k)$ is equal to $\abs {k\gpm}$. Hence $\widetilde k$ encodes information on both 
the residue field and the value group of $\abs{\cdot}$ (for other manifestations of this phenomenon, see
\ref{ss-index-graded}
below).

Note that if~$\abs{\cdot}$ is the {\em trivial} valuation 
(\ie, $\abs x=1$ for all $x\in k\gpm$) then~$\widetilde k=\widetilde k^1$; if not, it does not seem that~$\widetilde k$
can be interestingly interpreted as a residue graded field in the sense of \ref{ss-valgrad-ring}.

\begin{exem}
Let $p$ be a prime number. There
is an isomorphism 
of graded algebras over $\F_p=\widetilde{\Q_p}^1$
from $\F_p(T/\abs p)$
to $\widetilde{\Q_p}$, which sends $T$ to $\widetilde p$
(see \ref{sss-notation-ktr} for the meaning of the notation
$\F_p(T/\abs p)$). 

\end{exem}

\begin{exem}\label{exem-gradres-etar}
Let $k$ be a field endowed with a valuation $\abs{\cdot} \colon k\to \R_+$. 
Let $(r_1,\ldots, r_n)$ be a family of positive real numbers, and let $T=(T_1,\ldots, T_n)$ 
be a family of indeterminates. The formula
$\sum a_I T^I\mapsto \max \abs a_I\cdot r^I$ defines a real-valued valuation on $k(T)$. 
It follows from its very definition that there exists a (unique) $\widetilde k$-isomorphism
of graded fields 
$$\widetilde k\left(\frac \tau r\right)\simeq \widetilde {k(T)}$$ sending $\tau_i$ to $\widetilde{T_i}$ for every $i$, 
where $\tau=(\tau_1,\ldots, \tau_n)$ is a family of indeterminates. 
\end{exem}

\subsection{Graded interpretation of classical invariants}\label{ss-index-graded}
Let $k\hookrightarrow L$ be an isometric extension of real-valued fields.
Classical valuation theory assigns four invariants
to such an extension, which are (possibly infinite) cardinal numbers: 

\begin{itemize}[label=$\bullet$]
\item The ramification index $e$, which is the cardinality of $\abs {L\gpm}/\abs {k\gpm}$. 

\item  The inertia index $f$, which is the dimension of the $\widetilde k^1$-vector space $\widetilde L^1$. 

\item  The rational rank $r$, which is the dimension of the $\Q$-vector space $\abs {L\gpm}^\Q/\abs {k\gpm}^\Q$. 

\item  The residue transcendence degree $d$, which is the transcendence degree of $\widetilde L^1$ over $\widetilde k^1$. 
\end{itemize}

The product $ef$ and the sum $r+d$ admit natural interpretations in terms
of graded reduction:

\begin{enumerate}[1]
\item The  product $ef$ is the dimension of the $\widetilde k$-graded vector space $\widetilde L$, by equality
(a) of \ref{ss-graded-equalities}.

\item The sum $r+d$ is the transcendence degree of the graded extension $\widetilde k\hookrightarrow \widetilde L$
by equality (c)
of \ref{ss-graded-equalities}.
Note that it is always bounded by the usual transcendence degree of $L$ over $k$ (this is the so-called
\emph{Abhyankar inequality}, \cf \cite{bourbaki1985}, Chapitre VI, \S 10, \no 3, \Cor 1). 
\end{enumerate}

\subsection{Graded reduction of algebraic extensions}\label{ss-gradval-res}
Let $k\hookrightarrow L$ be an algebraic isometric extension of
real-valued (non-graded) fields, and let $\Gamma$ be a subgroup of $\R_+\gpm$.
The graded field $\widetilde L$ is algebraic over
$\widetilde k$, as we see by reducing to the case where $L$ is finite
over $k$, in which case $\widetilde L$ is finite over $\widetilde k$ by \ref{ss-index-graded} (1). 
This implies that $\widetilde L^\Gamma$ is algebraic over $\widetilde k^\Gamma$ (\ref{ss-criterion-alg}). 

Assume first that $L$ is radicial over $k$. In this case, $\widetilde L$ is radicial
over $\widetilde k$ (which implies that $\widetilde L^\Gamma$ is radicial over $\widetilde k^\Gamma$). 
Indeed, if \car $k=0$ then $L=k$ and $\widetilde L=\widetilde k$. And if \car $k=p>0$ then for every element $a$ 
of $L$, there exists $n\in \N$ such that $a^{p^n}\in k$; hence for every homogeneous element $\alpha$ 
of $\widetilde L$, there exists $n\in \N$ such that $\alpha^{p^n}\in \widetilde k$.

Assume now that $L$ is henselian (\eg, $L$ is complete)
and Galois over $k$. In this case, the valuation of $L$ is preserved by the Galois action, 
and $\mathrm{Gal}(L/k)$ acts therefore in a natural way on $\widetilde L$.
By \Prop 2.11
of \cite{ducros2013b}, the graded field $\widetilde L$ is normal over $\widetilde k$ 
and
the natural map $\mathrm{Gal}(L/k)\to \mathrm{Gal}(\widetilde L/\widetilde k)$ is surjective. 

By \ref{ss-graded-galois}, this implies that
the graded field
$\widetilde L^\Gamma$ is normal over $\widetilde k^\Gamma$ 
and
the natural map $\mathrm{Gal}(L/k)\to \mathrm{Gal}(\widetilde L^\Gamma/\widetilde k^\Gamma)$ is surjective.

\backmatter
\chapter*{Index of notation}

{\small
\noindent
{\em General conventions.} From the beginning of chapter \ref{c-gstrict}
till the end of the memoir (except the appendix), 
$k$ denotes an analytic field without any specific assumption (Convention before the introduction of Chapter \ref{c-gstrict}). 
In Chapter \ref{c-gstrict} and parts of Chapters \ref{GR}, \ref{DEV}
and \ref{IM}, $\Gamma$ denotes a subgroup of $\R_+\gpm$
such that $\Gamma \cdot \abs{k\gpm}\neq \{1\}$; \ie, $\Gamma\neq \{1\}$ whenever $k$
is trivially valued. 

\vspace{1cm}
\noindent
{\bf Topology}

\medskip
$\adht EX\z$ (\ref{ss-remind-topo})

\bigskip
\noindent
{\bf Scheme theory}

\medskip
$\mathscr O_{X,x}\z \mathfrak m_x\z \kappa(x) \z \mathfrak F_x \z \mathfrak F_{\kappa(x)} \z Y_x\z
\mathscr F_Y\z$ (\ref{ss-conv-scheme}) 

$\uflat {\mathscr E}X\z$ (\ref{notation-uflat})

\bigskip
\noindent
{\bf Analytic fields and affinoid algebras}

\medskip
$\widetilde k^r\z \widetilde k^1 \z \widetilde x^r \z \widetilde x \z$ (\ref{ss-conv-tilde})

$d_k(L)\z$ (\ref{sss-analytic-extension})

$A_L\z$ (\ref{ss-space-nofield})

$k_r\z$ (\ref{ss-intro-etar}) 

$A_r\z$ (\ref{ss-shilov-section})

$\dim_k A\z$ (\ref{ss-dim-affalg})

\bigskip
\noindent
{\bf Analytic spaces}

\medskip
$X_L\z$ (\ref{ss-space-nofield})

$\hr x\z$ (\ref{ss-hrx})

$Y_x \z$ (\ref{ss-def-xanalytic})

$\eta_{k,r}\z\eta_r\z$ (\ref{ss-intro-etar})

$X_r\z$ (\ref{ss-shilov-section})

$\mathrm {Int}(Y/X)\z \mathrm{Int}(X)\z \partial(Y/X)
\z \partial(X)\z$ (\ref{ss-boundary})

$\mathscr O_X\z$ (\ref{conv-ox})

$\mathscr F_Y\z \mathscr F_L\z \mathscr F_r\z$ (\ref{ss-cohsheaf-notations})

$\mathscr F\boxtimes \mathscr G\z$ (\ref{ss-fboxg})

$\mathscr F_x\z \mathscr O_{X,x}\z \mathfrak m_x\z\kappa(x)\z \mathscr F_{\kappa(x)}\z$ (\ref{ss-stalks-fibers})

$x_V\z \kappa(x_V)\z\mathscr F_{\kappa(x_V)}\z t_x\z$ (\ref{rem-x-xv})

$\adhz EX\z$ (\ref{ss-zariski-topology})

$Y_{\mathrm{red}}\z$ (\ref{ss-red-structure})

$\dim_k X\z$ (\ref{ss-dim-arbitrary-space})

$d_k(x)\z$ (\ref{ss-dimx-dkx})

$\dim_{k,x} X\z$ (\ref{ss-def-dimloc})

$\dim X\z \dim_x X\z \dim_y \phi\z$ (\ref{ss-dim-nofield})

$\codim (Y,X)\z \codim_x (Y,X)\z$ (\ref{ss-codim-def})

$\mathscr X\an\z X\al\z x\al \z F\al \z E\an \z \mathscr F\an \z \mathscr G\al\z$ (\ref{ss-def-analytify})

$x_V\al \z T_x\al \z t_x\al \z T\al_{x\al}\z$ (\ref{ss-conventions-gaga})

$\mathscr F_{\hr x}\z \rk x(\mathscr F)\z$ (\ref{defi-fhx-rkx})

$\supp F\z\dim \mathscr F\z\dim_x \mathscr F$ (\ref{ss-interpret-support})

$\Omega_{Y/X}\z$ (\ref{ss-remind-omega})

$\bij {\mathscr F}{\mathscr G}\z$ (\ref{ss-bij-locus})

$\mathrm{codepth}_y \;\mathscr F\z$ (\ref{def-codepth-cohsheaf})

$\uflat {\mathscr F}X\z$ (\ref{notation-uflat})

$\adht E \phi\z$ (\ref{def-ebar-phi})

\bigskip
\noindent
{\bf Analytic germs and Temkin's theory}

\medskip
$(X,x)\z$ (\ref{ss-recall-germs})

$ \mathrm{centdim (X,x)}\z$ (\ref{def-centdim})

$\P_{L/K}\z \P_{L/K}\{S\}\z$ (\ref{ss-def-plkzr})

$\mathscr S_{L/K}\z$ (\ref{ss-def-slk})

$\mathsf X^\Gamma\z$ (\ref{ss-dstrict-ff})

$\widetilde{(X,x)}\z$ (\ref{ss-defredgrad-general})

$\widetilde{(X,x)}^\Gamma\z$ (\ref{ss-graded-reduction})

\bigskip
\noindent
{\bf Abstract formalisation of algebraic properties}

\medskip
$\mathfrak T \z$ (\ref{def-category-framework})

$\mathfrak F\z \mathfrak F_X\z \mathfrak F_A \z
D_Y \z D_A\z D_x \z D_L \z D\an\z$ (\ref{ss-fiber-category})

$\mathfrak L\z \mathfrak F_{\mathfrak L}\z$ (\ref{ss-mathfrak-l})

$\mathfrak{Coh}\z$ (\ref{ex-fiber-coh})

$\mathfrak {Coh}^{\mathfrak I}\z$  (\ref{ex-fiber-diag})

$\mathsf P\z$ (\ref{ss-alg-properties})

\gen~~~~~(\ref{rem-valid-localization})

\hreg~~\hci~~\hv~~\field~~\open~~~~~(\ref{ss-list-hreg})

$\mathfrak C\z \mathsf Q\z$ (\ref{ss-general-axioms-pfib})

$\mathscr S\z$ (\ref{ss-def-functorS})

\hwk~~~ \hstr~~~\hwkk~~~\hstrr~~~~(\ref{ss-hweak})

\bigskip
\noindent
{\bf Graded commutative algebra}

\medskip
$A^r\z$ (\ref{def-graded-ring})

$A^\Gamma\z$ (\ref{ss-grad-agamma})

$\mathfrak D(K)\z$ (\ref{sss-grad-field})

$K[T/r]\z K(T/r)\z$ (\ref{sss-notation-ktr})

$M\otimes_A N\z$ (\ref{ss-desc-gradotimes})

$\abs \cdot_{\mathrm{Gau\ss}}\z$ (\ref{ss-notation-gauss})
}

\printindex
\bibliographystyle{smfalpha}
\bibliography{aducros.bib}

\end{document}